\documentclass[12pt,twoside,reqno,a4paper]{amsart}
\usepackage[latin1]{inputenc} 
 \usepackage[english]{babel}
\usepackage{amsthm}
\usepackage{amssymb, verbatim}
\usepackage{amsmath} 
\usepackage{eucal}
\usepackage{mathrsfs} 
\usepackage{hyperref}
\usepackage[margin=3.5cm]{geometry}

\usepackage{changepage}

 \topmargin = - 30 pt 
\usepackage[all]{xy}
\usepackage{amsmath,amsthm,amssymb,amsopn,amscd,tikz-cd}

\newcommand{\unit}{I}

\newcommand     {\A}    {{\mathcal A}}

\renewcommand   {\S}    {{\mathcal S}}

\newtheorem{thm}{Theorem}[section]

\newtheorem{lemma}[thm]{Lemma}
\newtheorem{prop}[thm]{Proposition}
\newtheorem{cor}[thm]{Corollary}

\theoremstyle{definition}

\newtheorem{defn}[thm]{Definition}

\newtheorem{ex}[thm]{Example}

\newtheorem{rem}[thm]{Remark}

\newtheorem{problem}{Problem}

\numberwithin{equation}{section}

\begin{document}

\title[Unitary coboundary weak quasi-Hopf algebras]{Unitary coboundary weak quasi-Hopf algebras  in  the   Kazhdan-Lusztig-Finkelberg    theorem}

\author[S.~Ciamprone]{Sergio Ciamprone}
\author[M.V. Giannone]{Marco Valerio Giannone}
 \author[C.~Pinzari]{Claudia Pinzari }
 \email{
 sergiociamprone@gmail.com\\ 
marcov89@gmail.com \\ pinzari@mat.uniroma1.it}
\address{Dipartimento di Matematica, Sapienza Universit\`a  di
Roma\\ P.le Aldo Moro, 5 -- 00185 Rome, Italy}\maketitle

\centerline\small{\it Dedicated to the memory of Sergio Doplicher and   John E. Roberts.}

\begin{abstract}
We construct, for every complex simple Lie algebra $\mathfrak{g}$ and every
positive integer level $k$, a unitary modular tensor category structure on the
module category of the affine vertex operator algebra $V_{\mathfrak{g}_k}$,
equivalent to the quantum group fusion category
$\mathcal{C}(\mathfrak{g},q,\ell)$, within the Doplicher--Roberts quantum group
duality program for the Wess--Zumino--Witten model. The construction reaches
positive integer levels that Finkelberg's argument does not cover, namely $E_6$
and $E_7$ at $k=1$ and $E_8$ at $k=2$, and it treats simply-laced and
non-simply-laced types uniformly, whereas the Kazhdan--Lusztig construction was
carried out for the simply-laced case.
 
We establish that a weak quasi-fibre structure on a functor preserving a
coboundary symmetry in the sense of Drinfeld extends a symmetric functor under
permutation symmetry. Using Wenzl's functor associated with the unitary quantum
group fusion category, we construct such a structure compatible with the
associativity morphisms, yielding unitary coboundary weak Hopf $C^*$-algebras
for all Lie types and positive integer levels, a new class axiomatically
introduced here. Via a specialized Drinfeld twist for our algebras and the Wenzl
de-quantization curve, this structure is transported onto the Zhu algebra, which
becomes a unitary coboundary weak quasi-Hopf $C^*$-algebra with a
$3$-coboundary associator. Rigidity is obtained from the quantum group side,
without the Knizhnik--Zamolodchikov equations or the Verlinde formula.
 
We further analyze the type $A$ case via classification methods based on
Kazhdan--Wenzl theory and our weak Hopf algebra framework, which gives insight
into the determination of associativity from the braiding in the general case.
Our methods settle several further problems: Galindo's question on the
uniqueness of unitary tensor structures, Kirillov's conjecture on the positivity
of a Hermitian form on the braided tensor category of modules of an affine Lie
algebra introduced by Beilinson, Feigin and Mazur, the quantum group structure
on the Zhu algebra sought by Frenkel and Zhu, and the
Kazhdan--Lusztig--Finkelberg equivalence for all Lie types and all positive
integer levels without passing through the negative shifted levels, as raised by
Huang. We also compare our braided tensor structure with that of Huang and
Lepowsky.
\end{abstract}
 
  \tableofcontents

  \section{Introduction}\label{1}

\subsection{Historical  context of the KLF equivalence}\label{1.1}

Let $\mathfrak{g}$ be a f.d. complex simple Lie algebra and $k$ a positive integer, 
defined as the physical level. Let $d$ denote the ratio of the squared lengths 
of the long to the short roots ($d=1$ for simply laced types, $d=2$ for types 
$B, C, F$, and $d=3$ for $G_2$). We define the shifted level as 
$\ell = d(k + \check{h})$, where $\check{h}$ is the dual Coxeter number of $\mathfrak{g}$.

  The current algebra of what is now known  as the Wess--Zumino-Witten (WZW) model was shown by Witten \cite{Witten_WZW}
     to satisfy the commutation relations of affine Kac-Moody algebras at positive integer levels.     Knizhnik and Zamolodchikov showed that    the $n$-point correlation functions of primary fields of the WZW model  satisfy
  a system of partial differential equations, the Knizhnik--Zamolodchikov (KZ) equations   \cite{KZ}. 

Pioneering work by Kohno \cite{Kohno2} established that the braid group representations with respect to the vector representation of a non-exceptional Lie algebra, arising from the monodromy of the  KZ equations in his framework, are equivalent to those originating from quantum groups. Kohno initially constructed these integrable connections on the full tensor product of   finite-dimensional representations of
the Lie algebra $\mathfrak{g}$ under a generic parameter $\lambda$. In the particular case of $\mathfrak{sl}_n$, he found the associated Hecke algebra representations.

While he noted that the KZ equations  appearing in conformal field theory were related to his framework when the parameter specializes to $\lambda = 1/\ell$, he restricted his explicit implementation of the affine framework at a positive integer level
 strictly to the case of the 
affine Kac--Moody algebra $\widehat{\mathfrak{sl}}_2$, drawing upon the operator formalism of Tsuchiya and Kanie \cite{Tsuchiya_Kanie}.
By utilizing Wenzl's classification of unitarizable Hecke algebra representations at roots of unity, Kohno established the unitarizability of these geometric braid group representations for both the vector and higher irreducible representations of $\mathfrak{sl}_2$.

Crucially, this framework relied fundamentally on the multiplicity-free nature of the $\mathfrak{sl}_2$ fusion rules and the direct path-algebra descriptions of the corresponding Temperley--Lieb algebra. For higher-rank Lie algebras, the emergence of non-trivial multiplicities and higher-dimensional intertwining spaces introduces deep combinatorial obstructions.

Inspired by these results, Drinfeld \cite{Drinfeld_quasi_hopf, Drinfeld_galois} introduced a deep braided tensor category structure on the category of finite-dimensional $\mathfrak{g}$-modules by equipping it with an associator explicitly derived from the KZ   equations. He proved the equivalence at a categorical level between his category---known as the Drinfeld category---and the rigid braided tensor category of modules over the formal quantum group $U_h(\mathfrak{g})$, marking a foundational breakthrough in representation theory.
Drinfeld obtained rigidity of the Drinfeld category by constructing an antipode on the associated quasi-bialgebra, following the Drinfeld--Kohno theorem, via an operation that he called a \textit{twist} of a quantum group  \cite{Drinfeld_quasi_hopf}.

 Subsequently, inspired by Drindeld's results, Kazhdan and Lusztig \cite{KLseries} extended Kohno's   braid group equivalence to the affine setting at negative rational shifted levels, constructing a   braided tensor category structure on certain non-semisimple module categories $\mathcal{O}_{-\ell}$ over the corresponding affine Kac-Moody algebras (specifically for simply laced types). They proved a   categorical equivalence between $\mathcal{O}_{-\ell}$ and the representation categories of corresponding quantum groups at  the root of unity $q=e^{\frac{i\pi}{\ell}}$. A main result of their work is a direct proof of rigidity of  $\mathcal{O}_{-\ell}$.

For the positive integer levels,     Finkelberg
    followed the algebraic geometric approach developed by Beilinson, Feigin, and Mazur \cite{BFM}
for the semisimple braided tensor category  $\tilde{\mathcal O}_{\ell}$ associated to an affine Lie algebra,  and constructed a functor  to a semisimple subquotient fusion   category $\tilde{\mathcal O}_{-\ell}$  of a Kazhdan-Lusztig category at  the opposite shifted level  \cite{Finkelberg},  and thus, by the Kazhdan-Lusztig equivalence, to a semisimple subquotient fusion category ${\mathcal C}({\mathfrak g}, q, \ell)$ of the corresponding quantum group at a root of unity.

     Frenkel and Zhu established a vertex operator algebra structure   on a certain representation $L_{k, 0}$ of the
  affine Kac-Moody---and analogously for the Virasoro algebra---in such a way that the representations of $L_{k, 0}$ give rise to the WZW model \cite{Frenkel_Zhu}.   We denote by $V_{{\mathfrak g}_k}$ this vertex operator algebra.

 Following an analytic approach based on differential equations, Huang and Lepowsky constructed a  braided tensor category structure on the module category of a vertex operator algebra in a very general setting,     including the affine cases, for which they relied on the KZ equations \cite{Huang_LepowskiI}--\cite{HuangIV}. Huang  gave a rigorous proof of the Verlinde formula  and  the fusion rules  \cite{Verlinde} and showed its relevance to
 prove rigidity and modularity  \cite{Huang1}, \cite{Huang2}. 
   
To  complete the proof of the tensor equivalence of his functor $\tilde{\mathcal O}_{\ell}\to\tilde{\mathcal O}_{-\ell}$,   Finkelberg        needed the fusion rules and the Verlinde formula as well \cite{Finkelberg_erratum}, established in his setting
    by Teleman and Faltings  \cite{Faltings}, \cite{Teleman}, guaranteeing that the $K$-rings of the positive and negative shifted level categories are based-isomorphic. He also relied on rigidity of Kazhdan-Lusztig category ${\mathcal O}_{-\ell}$ and of
    the   category of finite-dimensional representations of the simple Lie algebra ${\mathfrak g}$. He  also derived rigidity of   $\tilde{\mathcal O}_{\ell}$ from his tensor equivalence.
    
Summarizing, the existing proof of the tensor equivalence
$\tilde{\mathcal O}_\ell\to{\mathcal C}({\mathfrak g}, q, \ell)$ relies on
results spanning distinct algebraic and geometric settings across both positive
integer and negative shifted levels.
It requires separate quotient constructions to yield the semisimple subquotients
where rigidity defines the necessary tensor ideals, and uses rigidity to obtain
tensoriality of the equivalence.
Furthermore, because the Kazhdan--Lusztig construction does not apply to certain
negative rational shifted levels, Finkelberg's equivalence, and consequently the
proof of rigidity of $\tilde{\mathcal O}_{\ell}$, do not cover some exceptional
cases---namely, $E_6$ and $E_7$ at level $k=1$, and $E_8$ at levels $k=1,2$.
Additionally, it is not easy to find literature explicitly extending the
Kazhdan--Lusztig rigid braided tensor category construction to non-simply laced
Lie algebras, due to the greater complexities of the KZ equations.
Regarding the non-simply laced types,
it is noted in Sect. 2.6 of \cite{Finkelberg} and Sect. 2, part (f) of
\cite{Finkelberg_erratum} that due to a uniform lower bound constraint of $\ell$
needed to guarantee rigidity of ${\mathcal O}_{-\ell}$, the rigidity of
$\tilde{\mathcal{O}}_{\ell}$ and the tensor equivalence property of the
Finkelberg functor may be unproved for low levels $k\leq2$, requiring further
case-by-case study.
 
Taken together, the existing proof of the equivalence is distributed over
\cite{KLseries}, \cite{Finkelberg}, \cite{Finkelberg_erratum}, \cite{Faltings}
and \cite{Teleman}. Kazhdan and Lusztig carried out their construction for the
simply-laced types, while Finkelberg's functor is stated for the other types as
well.
 
Consequently, the need to describe a direct connection between affine Lie
algebra module categories at positive integer levels and quantum groups emerged,
see e.g. Gannon's observation in Sect. 6.2.3 of \cite{Gannon_book}.
Notably, in \cite{Huang2018}, Huang explained in detail the history of
Finkelberg's theorem and the need of finding a direct proof---one that bypasses
the obstructions arising from the negative shifted levels and includes the
exceptional cases (see Problem \ref{problem4}).
The interest of a direct construction extends beyond the case at hand. There are
conjectural equivalences between braided tensor categories of modules for other
vertex operator algebras and categories of modules for quantum groups for which
no Kazhdan--Lusztig category is available; in those cases the route through
negative shifted levels  is unavailable. Huang has pointed out to us that this is his main reason of interest in the problem. A direct construction of Finkelberg's equivalence is therefore a natural starting point for those cases.
 
In this spirit, we start with the quantum group fusion category. Wenzl showed in
\cite{Wenzl} that this category can be constructed intrinsically for all Lie
types and all positive integer levels, treating the simply-laced and
non-simply-laced cases uniformly, for the roots of unity of the kind
$q=e^{i\pi/\ell}$ starting with a fundamental (generating) representation $V$ of
$U_q({\mathfrak g})$, see Theorem 5.4 in \cite{CP} for an explicit proof.
 
Starting from Wenzl's realization of ${\mathcal C}({\mathfrak g}, q, \ell)$ and
adopting a direct, operator-algebraic framework, we give a direct proof of the
equivalence with a category of modules of the corresponding affine vertex
operator algebra, in a single self-contained text. In this way we circumvent the
Verlinde formula dependency of the equivalence with the categories of modules of
affine vertex operator algebras (or affine Lie algebras), the passage to the
negative integer levels, and we derive rigidity of these categories from the
equivalence, thereby also rendering boundary cases uniformly accessible. We
describe our approach in the next subsection, state the problems solved in this
paper in subsection \ref{1.3}, and describe the main results in section
\ref{5++}. In section \ref{History} we describe the organization of the rest of
the paper.

\subsection{Overview of our approach.}\label{1.2}  In the operator algebraic approach to quantum field theory (AQFT) introduced by Haag and Kastler \cite{Haag_Kastler}, the subsequent   work of Doplicher, Haag, and Roberts established   the rigidity of the associated symmetric tensor $C^*$-categories of superselection sectors   from physically motivated principles, yielding Dirac's particle-antiparticle duality as a   structural consequence \cite{DHR}, \cite{DHR2}, \cite{DR_qft}, \cite{Haag}.

For high-dimensional theories,
 a deep result by Doplicher and Roberts established that the category of localizable charges is equivalent to the representation category of a unique compact group  having the role of a global gauge group   \cite{DR1}.  
 The Doplicher--Roberts compact group associated with the category describes rigidity via the passage to the conjugate representation of the group. 
 A closely related result was independently obtained by Deligne   for affine group schemes over a field \cite{Deligne}.

 In low-dimensional theories, one encounters unitary rigid braided tensor $C^*$-categories, where the braided symmetry is given by unitary representations of the braid group acting on tensor powers of an object \cite{FRS, FRS2}.

 Parallel to this, as recalled in Subsect. \ref{1.1},, Jones index theory led to the construction of rigid tensor $C^*$-categories from inclusions of von Neumann algebras with trivial center (subfactors) and yielded unitary representations of the braid group satisfying Temperley--Lieb relations; conversely, inclusions of subfactors arise from such relations \cite{Jones_subfactors}.

 Shortly after Vaughan Jones' discovery of subfactor theory and the Jones 
polynomial, Wenzl connected the representation theory of Hecke algebras 
$H_n(q)$ of type $A$ at roots of unity to inclusions of finite-index 
$\mathrm{II}_1$ subfactors. This breakthrough provided a foundational framework 
for the unitary structures linking subfactor theory, quantum groups, and 
two-dimensional conformal field theory \cite{Wenzl_thesis, We1}. Wenzl 
explicitly constructed and classified  finite-dimensional, 
irreducible representations of the non-semisimple algebras $H_n(q)$ that are 
unitarizable, thereby constructing new hyperfinite $\mathrm{II}_1$ subfactors 
that generalize the Temperley-Lieb-Jones constructions. As observed in \cite{We1}, 
via quantum Schur-Weyl duality, these representations of Hecke algebras can 
naturally be viewed as arising from the endomorphism algebras of tensor powers 
of the fundamental representation of $U_q(\mathfrak{sl}_N)$ \cite{Jimbo}.

Drinfeld   introduced   {\it quantum groups} at the  ICM in Berkeley in 1986 \cite{Drinfeld_qg}.  Quantum groups provide rigidity via the antipode.
 The representation theory of quantum groups at roots of unity  was developed   by Lusztig \cite{Lusztig, Lusztig1990GeoDed}, Reshetikhin, Turaev \cite{RT}.

  We are naturally led to ask whether rigidity of categories on the conformal field theory side, and in particular
  for $V_{{\mathfrak g}_k}$, can be explained by an equivalence with a category of representations of a semisimple quantum gauge group in the sense of the Doplicher-Roberts program.
  We therefore look for a natural construction of a quantum version of the Doplicher--Roberts compact group for braided categories. A major difference with Drinfeld setting is that, in our setting, the involved braided tensor categories, and especially the quantum group fusion category, are not given from the start with a fibre functor to the category of vector spaces.
 
   Mack and Schomerus pioneered the   question of bridging the Doplicher--Roberts program on global quantum gauge groups and field algebras  in conformal field theory  with Drinfeld's work on the Drinfeld-Kohno theorem mentioned in Subsect. \ref{1.1}. Specifically, they introduced the {\it weak quasi-Hopf algebras} that extend Drinfeld's original structures, implementing them for the $\mathfrak{sl}_2$-case as a quantum analogue of Doplicher--Roberts compact groups  \cite{MS1, Schomerus}. 
  
Their early construction of both the weak quasi-Hopf algebras and the analogue of the Doplicher--Roberts field algebra encountered   complications, arising primarily from the non-triviality of the associator and a lack of uniqueness or naturality beyond the ${\mathfrak sl}_2$ case.

Wenzl extended the construction of subfactors from quantum groups at roots of unity
  and constructed unitary braided tensor structures on the quantum  group fusion category   uniformly for all Lie types,     in  \cite{We2}, \cite{Wenzl}.


Mack and Schomerus program was not   developed further until recently, when two of us found a natural construction of Mack-Schomerus weak quasi-Hopf algebras associated to  the quantum group fusion category of $U_q({\mathfrak sl}_N)$ for all positive integer levels \cite{CP}.

Our first aim in this paper is to put the examples of \cite{CP} into a new theory of {\it weak Hopf algebras}, and then extend 
the construction to all other Lie types. 

Specifically, we adopt an operator-algebraic approach where unitarity and   $C^*$-structures with unitary representations of the braid group play a primary role. 

Our approach heavily relies on the foundational $C^*$-algebraic techniques for quantum group fusion categories   developed by Wenzl.  
We consider Wenzl linear functor $$W: {\mathcal C}({\mathfrak g}, q, \ell)\to{\rm Hilb}$$ and construct on it a {\it weak tensor structure} in the   sense of Definition \ref{wtf}.

We achieve   new examples  $A_W({\mathfrak g}, q, \ell)$ in this paper for all Lie types and levels, with different methods as compared to the type $A$ case treated in \cite{CP}, without using the fusion rules, see Theorem \ref{main_wh}.
For ${\mathfrak g}={\mathfrak sl}_N$ we recover our earlier examples.
We prove that these algebras satisfy properties of
{\it compatible unitary coboundary weak Hopf algebras} in a sense that we introduce in Definition \ref{strongly_Hermitian_ribbon_wqh}.

The   class of \emph{unitary coboundary weak (quasi-)Hopf algebras} 
(Definition~\ref{Hermitian_ribbon_wqh}) admits a categorical characterization 
that significantly simplifies their axiomatic description; see 
Theorem~\ref{TK_unitary_ribbon}. In the specific case of a symmetric tensor 
$C^*$-category endowed with a unitary fibre functor to the category of Hilbert 
spaces, this characterization reduces precisely to the requirement that the 
functor is symmetric. This constraint is required and, in fact, 
plays a fundamental role in the   Doplicher--Roberts duality theory for compact groups,  in determining the fibre functor uniquely.


 In the braided setting, the categorical characterization of a unitary 
coboundary weak quasi-Hopf algebra requires finding a weak quasi-tensor 
structure on the given functor  implemented by {\it coisometric} 
fibre maps $F_{\rho, \sigma}$ associated with the tensor product bifunctor of the 
underlying braided category. It is required that the structure maps respect the coboundary symmetries 
of both the source and target categories. Indeed, the unique coboundary symmetry 
on the category of Hilbert spaces is the standard permutation symmetry. Since the 
coboundary symmetry of a braided tensor category is obtained by multiplying the 
braiding by a $2$-coboundary derived from a square root of the ribbon element, 
this characterization forces a rigid constraint: via the fibre functor, the 
braided symmetry of the category is uniquely determined as the permutation 
symmetry modified by this square root $2$-coboundary.


 Formulas of this nature   indeed govern the explicit structure of the braided 
symmetry in categories arising from conformal field theory, such as loop group 
conformal nets or affine vertex operator algebras (see 
Remarks~\ref{DR_symmtric_functor}, \ref{HL_braiding}, and references therein). 
This operator algebraically realizes the categorical balancing relation where the double braiding equals the $2$-coboundary defined by the twist. On an irreducible component of {\it   the minimum energy subspace} of the tensor product, the braiding evaluates precisely to the square root of the Virasoro phases multiplied by the permutation symmetry, see Theorem 4.1 in \cite{Huang2}.

While the braided 
symmetry in these settings is traditionally derived from the monodromy of the 
Knizhnik--Zamolodchikov equations, our framework provides an independent 
method for constructing this symmetry entirely via unitary coboundary weak 
quasi-Hopf algebras with coisometric structure maps. In our view, this 
construction beautifully realizes an insightful perspective of Vaughan Jones on 
the profound relevance of \emph{manifest unitarity} in conformal field 
theory (cf.~\cite{pinzari_survey}).



In the quantum group setting, the action of the $R$-matrix on tensor product spaces does not take an analogous simple diagonal form. Consequently, the associated coboundary matrix does not act trivially. Since this matrix induces the unitary structure on the fusion spaces \cite{Wenzl}, this unitary structure effectively measures the deviation of the $R$-matrix from acting strictly via the square roots of the ribbon phases--a purely quantum effect. We explicitly quantify this phenomenon. Specifically, we show that the unitary modular quantum group fusion category $\mathcal{C}(\mathfrak{g}, q, \ell)$, which is described by the unitary representations of  
our unitary coboundary weak Hopf algebra $A_W(\mathfrak{g}, q, \ell)$ by Theorem \ref{main_wh}, is fundamentally incompatible with the existence of coisometric fusion projections due to the weakly trivial associator of $A_W(\mathfrak{g}, q, \ell)$ and the non-integrality of the quantum dimensions of  $\mathcal{C}(\mathfrak{g}, q, \ell)$ (see Theorem \ref{unitarity_obstruction}). Thus, extracting coisometric fusion projections via polar decomposition inevitably twists the weak Hopf algebra into a genuine weak quasi-Hopf algebra, necessarily equipped with a 3-coboundary associator.

By applying an analytic Drinfeld twist method and Wenzl's de-quantization curve to the weak Hopf algebras 
$A_W({\mathfrak g}, q, \ell)$, we obtain a self-contained construction of a 
unitary rigid braided tensor structure on the module category of the
 affine 
vertex operator algebra $V_{{\mathfrak g}_k}$ 
(see 
Corollary~\ref{cor_Zhu_as_a_compatible_unitary_wqh}). 

A pivotal feature of 
our unitary structure on ${\rm Rep}(V_{{\mathfrak g}_k})$ and  
analytic twist is that it equips  Zhu's {\it minimum energy functor}
with a canonical weak quasi-tensor structure $(F_0, G_0)$, in a way that $F_0$ coincides with the restriction of Huang-Lepowsky tensor product structure map, and moreover the twist forces   $F_0$   to be 
coisometric and $G_0$ to be the isometric Hilbert space adjoint of $F_0$,
a property we   establish in Theorem~\ref{Kirillov}.
Comparing with Huang-Lepowsky ribbon braided tensor structure, at this stage we have established   an equivalence of the two   categories, the ribbon structures and the tensor product bifunctors for the affine vertex operator alegbras at positive integer levels.

In searching the Drinfeld twist  we have been    inspired   by   the categorical framework of Drinfeld-Kohno theorem established by Neshveyev and Tuset  for compact quantum groups \cite{NT_KL}. However, our twist seems different from their twist. Our twist may be regarded as a specialization at the root of unity to the fusion algebra  $A_W({\mathfrak g}, q, \ell)$ of the original Drinfeld twist in the formal setting, arising from a square root of  Drinfeld coboundary matrix \cite{Drinfeld_quasi_hopf}, that we consider in an operator algebraic sense.

We prove all the properties of the affine vertex operator algebra fusion category at a positive integer level using the methods of quantum symmetry and the Drinfeld twist. Specifically, the ribbon-braided tensor structure, fusion rules, and rigidity of the unitary quantum group fusion category ${\mathcal C}({\mathfrak g}, q, \ell)$ induce a corresponding structure on the weak Hopf $C^*$-algebra $A_W({\mathfrak g}, q, \ell)$.
From this position, we transfer the structure to the Zhu algebra via our Drinfeld twist and the  de-quantization curve used in \cite{Wenzl} to prove unitarity of the fusion category.

Wenzl's curve  gives  an algebra isomorphism between $A_W({\mathfrak g}, q, \ell)$ and the Zhu algebra $A(V_{{\mathfrak g}_k})$.
The twist transports the rest of the structure: coproduct, (weak trivial) associator, antipode, $R$-matrix, unitary structure, making the Zhu algebra into a unitary coboundary weak quasi-Hopf $C^*$-algebra with 3-coboundary associator.
From this position, we lift the structure to the module category ${\rm Rep}(V_{{\mathfrak g}_k})$ of the vertex operator algebra via Zhu's linear equivalence and its inverse.
 In this way, we obtain the structure of a unitary rigid ribbon-braided tensor category on ${\rm Rep}(V_{{\mathfrak g}_k})$ braided tensor equivalent to the quantum group fusion category.

We deduce the rigidity of our VOA fusion category ${\rm Rep}(V_{{\mathfrak g}_k})$ directly from our tensor equivalence, driven by the antipode of $A_W({\mathfrak g}, q, \ell)$ and our specialized Drinfeld twist.

 
  This approach allows us to
  bypass the need to dualize through negative levels, provide an explicit, self-contained proof of the  Kazhdan-Luszting-Finkelberg theorem that   includes the exceptional boundary cases, and simultaneously settle the   problem of Galindo on uniqueness of the unitary structure in a tensor category,  as well as those of Frenkel-Zhu and Doplicher-Roberts for the WZW model within the framework of vertex operator algebras. Along our path, we solve Kirillov's conjecture
 in the algebraic-geometric framework of Beilinson, Feigin, and Mazur utilized by Finkelberg, using general Tannakian methods for weak quasi-Hopf algebras, Finkelberg's equivalence theorem,  and Wenzl's unitarity of ${\mathcal C}({\mathfrak g}, q, \ell)$.

We  complete the  comparison between our unitary ribbon braided structure with the Huang--Lepowsky structure. We find a complete identification of the full categorical data (ribbon structure, tensor product bifunctor, braiding, and associativity morphisms) for the classical Lie types and $G_2$. For all  Lie types, we obtain identifications between the two structures for the braiding and associativity morphisms on many objects involving tensor products with the fundamental representations, see Theorem \ref{Finkelberg_HL}. To keep the focus on the core ideas developed in this paper, the detailed proofs of Theorems \ref{claim0} and \ref{claim1}, as well as their applicability to the stated Lie types, are deferred to our companion paper \cite{On_a_problem_posed_by_Huang}. In this way we   obtain an  independent construction of the Huang-Lepowsky braided tensor structure
on ${\rm Rep}(V_{{\mathfrak g}_k})$ which does not depend on the Knizhnik-Zamolodchikov equations, and a
proof of rigidity of this category
that does not depend on the Verlinde formula or fusion rules.

Finally, in \cite{pinzari_survey}, we provide a conceptual and expository historical overview of the main ideas of our proof of the KLF theorem and its applications to the unitarizability of categories arising from conformal field theory, and we outline directions for future research.

The remainder of the paper is organized as follows. In the subsequent Subsect. \ref{1.3}, we we   state the open Problems \ref{problem1}, \ref{problem2}, \ref{problem_kirillov},  \ref{problem3}, \ref{problem4}, and \ref{problem5}   which are of primary interest in our paper (with Problem \ref{problem1} serving as our original motivation). In Sect. \ref{5++}, we state our main results. In Sect. \ref{History}, we  provide an outline of the rest of the paper. In doing so, we discuss the utility of the KLF theorem (or other instances of equivalences coming, e.g., from pointed fusion categories) to transport unitary structures from the quantum group fusion category to the affine vertex operator algebra fusion category.

  \smallskip

  \subsection{Statements of the problems under consideration in this paper
  }\label{1.3}

As outlined in the previous subsections, our uniform operator-algebraic framework is connected to several   long-standing problems across algebraic quantum field theory (AQFT), conformal field theory, and tensor category theory. For clarity, we   state these problems here.

As recalled in Subsect. \ref{1.2}, the   Doplicher--Roberts theorem characterizes the category of high-dimensional localizable charges by  the representation category of a unique compact group. However, in low dimensions, the statistics dimension may be non-integer, meaning the category cannot correspond to a standard compact group. Because the braid group naturally appears in low-dimensional quantum field theories, Doplicher and Roberts posed the following problem, which served as the original motivation for this paper:

 \begin{problem}\label{problem1} (Sect. 7 in  \cite{DR})
 {\it Can we extend the duality theory for compact groups \cite{DR1} Theorem 6.1 
  to more general objects, by replacing strict symmetric monoidal C*-categories by
strict braided monoidal C*-categories? What is the full class of compact group-like
objects which arise in this way?}
 \end{problem}  

Parallel to this, in the study of conformal field theory (CFT), Moore and Seiberg conjectured analytic properties of quantum fields and discovered polynomial equations that later evolved into the notion of modular tensor categories \cite{Moore-Seiberg1, Moore-Seiberg2, RT2, RT, Turaev_modular_categories}. Noting that in various examples of CFT, the braid matrix is closely related to the $R$-matrix of specific quantum groups, they posed the following:

\begin{problem}\label{problem2} (Sect. 9 in  \cite{Moore-Seiberg2}) {\it Understand chiral algebras of conformal field theory as generalization of quantum groups. }
  \end{problem}

  Regarding construction of  unitary structures on   modular fusion categories,
  Kirillov defined   hermitian forms on the quantum group fusion category ${\mathcal C}({\mathfrak g}, q, \ell)$ in \cite{Kirillov} and on 
  the Beilinson-Feigin-Mazur
   category    
   $\tilde{\mathcal O}_\ell$ at positive integer levels in \cite{Kirillov3}, and conjectured that they are unitary. The quantum group case was settled by Wenzl in \cite{Wenzl}, while the affine Lie algebra case, to the authors' knowledge,  remained open.

  \begin{problem}\label{problem_kirillov} (Remark 10.7 in \cite{Kirillov3})
  {\it  The
natural conjecture, parallel to the quantum group case, is that so defined inner
product on the spaces of morphisms is positive definite. So far, we have no proof
of it except for sl2 case where it can be checked by direct calculation.}
  \end{problem}

Observing the deep connections between the categorical data of CFTs and Drinfeld's quasi-Hopf algebras \cite{Drinfeld_quasi_hopf}, Frenkel and Zhu asked for a rigorous algebraic bridge in the framework of vertex operator algebras:

  \begin{problem}\label{problem3} (\cite{Frenkel_Zhu}) {\it  A complete conceptual explanation of this phenomenon is needed. In the setting of vertex operator algebras,  associate a general construction of quasi-triangular quasi-Hopf algebra with a vertex operator algebra with the same tensor category of representations.}
  \end{problem}
  
  Further related conjectures may be found in \cite{GG}, \cite{Montgomery}.

Addressing the tensor category structure directly, Huang and Lepowsky constructed a general tensor product theory for module categories of VOAs \cite{Huang_LepowskiI, Huang_LepowskiII, Huang_LepowskiIII, HuangIV}, establishing, in particular, the modular tensor structure for affine VOAs at positive integer levels \cite{Huang(modularity), Huang1, Huang2}. To reconcile this intrinsic VOA approach with the historical algebraic-geometric proof of the Kazhdan--Lusztig--Finkelberg equivalence, Huang formulated the following challenge:

\begin{problem}\label{problem4} (Problem 4.4 in \cite{Huang2018})
  {\it Find a direct construction of this   equivalence without using the equivalence given by Kazhdan-Lusztig so that this equivalence covers all the cases, including the important ${\mathfrak g}=E_8$ and $k=2$ case.} \medskip
\end{problem}

Finally, our operator-algebraic methods based on weak quasi-bialgebras require establishing unique $C^*$-structures on tensor categories. This directly addresses a fundamental question raised by Galindo \cite{Gal}:

 \begin{problem}\label{problem5} {\it May a fusion category   admit more than
a unitary structure making it into a unitary tensor category?}
\end{problem}

In the following Section 2, we state our main results, demonstrating how the construction of weak quasi-Hopf $C^*$-algebras and the transport of unitary coboundary structures via Drinfeld twists provides a direct, unified solution to these five problems.

\section{Main results}\label{5++}

 In this section we state  some of our main results that we prove in this paper regarding Problems  \ref{problem1}, \ref{problem2}, \ref{problem_kirillov}, \ref{problem3}, \ref{problem4}, \ref{problem5}.
 
\subsection{Unitarizability of   fusion categories: The  Kirillov and Galindo conjectures} \bigskip

Let ${\mathfrak g}$ be a complex simple Lie algebra and let $U_q({\mathfrak g})$ be the Drinfeld-Jimbo quantum group
specialized  at the root of unity  $q=e^{i\pi/\ell}$, with $\ell$ a positive integer divisible by $d$,   the ratio of the square   lengths of the long and the short roots,
and   $\ell/d$ larger than    the dual Coxeter number   ${h}^\vee$  of ${\mathfrak g}$.
 Let ${\mathcal C}({\mathfrak g}, q, \ell)$ denote the associated   fusion category.
 Let us define the level $k$ associated to $\ell$ by $$\ell= d(h^\vee+k).$$

   Using general Tannakian methods for weak quasi-Hopf algebras, we develop a main unitarizability criterion for braided tensor categories, Theorem \ref{unitarizability}. A main consequence    is the solution of Kirillov's Problem \ref{problem_kirillov} given in Theorem \ref{BFM}. The proof is 
based on the general theorem to   the original Finkelberg's theorem \cite{Finkelberg}, \cite{Finkelberg_erratum} and Wenzl's unitarity  result of ${\mathcal C}({\mathfrak g}, q, \ell)$ \cite{Wenzl}.

   \begin{thm}\label{BFM2}
Under the same assumptions as in Finkelberg's theorem \ref{Finkelberg},  $\tilde{\mathcal O}_\ell$ becomes a unitary modular fusion category with Kirillov's Hermitian structure.
   \end{thm}

Another consequence of Theorem \ref{unitarizability} is the first positive answer to Galindo's question Problem \ref{problem5} for a wide class of tensor categories 
with possibly infinitely many simple objects, the proof is  in Sect. \ref{14}, see Theorem \ref{Galindo2}.  Another  proof   has been given by Reutter in \cite{Reutter} with different methods. 

\begin{thm}\label{Galindo}
Let ${\mathcal C}_1$ and ${\mathcal C}_2$ be tensor equivalent $C^*$-tensor categories endowed with an integral weak dimension function (e.g. they are finite semisimple tensor categories). Then ${\mathcal C}_1$ and ${\mathcal C}_2$ are also unitarily tensor equivalent.
\end{thm}

These theorems are interesting especially taking into consideration the fact that the abstract theorem \ref{unitarizability} on which they rely uses only the choice of an integral weak dimension function on the affine Lie algebra fusion category. 
This category naturally admits such a   dimension function given by representations of the simple Lie algebra $\mathfrak{g}$ and no inner structure is needed in detail, other than the two cited milestone theorems.

  Theorem \ref{BFM2} is,  in other   words, a strong motivation
 to refine our Tannakian methods compatibly with the inner structure of the categories. This leads us to seek a direct proof of
 the equivalence between quantum group fusion categories and CFT fusion categories.
 To get a deeper understanding of our unitarizability criterion for affine Lie algebra fusion categories, we need to examine their inner structures and utilize their naturally arising integral dimension functions.      \smallskip
 
\subsection{Lie algebra and level set up}\label{constraints}  Our framework relies on the setup in \cite{Wenzl}, which applies to all simple Lie algebras, whether simply laced or non-simply laced.  The fundamental representation $V$ of ${\mathfrak g}$ is defined in \cite{Wenzl} for each Lie type.
 We also denote by $V$ the  corresponding quantized counterpart of $U_q({\mathfrak g})$, $q=e^{i\pi/\ell}$, and
the associated object of  ${\mathcal C}({\mathfrak g}, q, \ell)$.

The minimal level for which the fundamental representation belongs to the   Weyl alcove is $k=1$ for all types $\mathfrak{g} \neq E_8$, and $k=2$ for $\mathfrak{g}=E_8$. Accordingly, we assume throughout this paper that $k \geq 1$ if $\mathfrak{g} \neq E_8$ and $k \geq 2$ if $\mathfrak{g} = E_8$.   \medskip
 
 \begin{rem}  The case $\mathfrak{g}=E_8$ at level $k=1$ yields a trivial
 pointed braided category, and it is the only case.
  \end{rem}
 
 \begin{rem} 
 The sole exception to the constraints in \ref{constraints} is part (c) of Theorem \ref{Finkelberg_HL}, 
 which establishes   the equivalence of  our structure with the Huang--Lepowsky braided tensor structure.

Our strategy to prove part (c) of Theorem \ref{Finkelberg_HL}   can be viewed as a wide generalization of the classification of braided tensor structures  for pointed categories case (examples of such categories are described    in Subsect. \ref{22.1} for  the type $A$   level $k=1$ case). To this aim, we develop a    classification of braided tensor structures in a semisimple pre-tensor category endowed with a generating object, see   Theorems \ref{claim0} and \ref{claim1}.
 \end{rem}

 By the main result of \cite{Wenzl}, ${\mathcal C}({\mathfrak g}, q, \ell)$ is a unitary ribbon braided rigid tensor category
 with natural unitary representation of the braid group.

\subsection{Construction of the weak Hopf algebras $A_W({\mathfrak g}, q, \ell)$ (global quantum gauge groups)} 
The weak Hopf algebras $A_W({\mathfrak g}, q, \ell)$  summarized in this subsection will be the solution of Doplicher-Roberts Problem \ref{problem1} for the braided tensor categories associated to the WZW model by Corollary \ref{cor_Zhu_as_a_compatible_unitary_wqh} and Theorem \ref{Finkelberg_HL}.

Starting from   \cite{Wenzl}, 
in Sect. \ref{unitary_structure_of_fusion_category} we recall the construction  of a
 strict   braided tensor category  $\tilde{\mathcal G}_q$ with objects the fusion tensor powers of   $V$, and we recall an equivalence
from  $\tilde{\mathcal G}_q$   to   ${\mathcal C}({\mathfrak g}, q, \ell)$ with an embedding full functor.
 Wenzl identified a
 natural linear functor $$W:  \tilde{\mathcal G}_q\to{\rm Hilb}.$$   
 
In Sect. \ref{6} we introduce the new definition of   {\it weak Hopf algebra} and   Theorem \ref{TK_algebraic} gives a categorical characterization by weak tensor functors in the semisimple case.

Moreover, motivated by the hermitian structure of $U_q({\mathfrak g})$ 
at $q=e^{i\pi/\ell}$, the specialization of Drinfeld coboundary matrix $\overline{R}$ of the generic $U_h({\mathfrak g})$ \cite{Drinfeld_quasi_hopf} and its explicit relation with the hermitian structure  of $U_q({\mathfrak g})$ 
 detailed by Wenzl in \cite{Wenzl}, we introduce the definition of {\it compatible unitary coboundary weak quasi-Hopf algebra} in Def. \ref{strongly_Hermitian_ribbon_wqh}, which replicates the   properties as for the hermitian structure of the non-semisimple $U_q({\mathfrak g})$  in a unitary setting.  
  We give a new construction extended to all Lie types and levels of such algebras $A_W({\mathfrak g}, q, \ell)$ associated to Wenzl's functor $W$.

\begin{thm}\label{main_wh} 
  Wenzl functor  $W$   admits a natural weak tensor structure 
making $A_W({\mathfrak g}, q, \ell):={\rm Nat}(W)$ into a unitary coboundary weak  Hopf $C^*$-algebra with   
$\Omega$-involution and ribbon tensor structure compatible with ${\mathcal C}({\mathfrak g}, q, \ell)$, and with antipode of Kac type.  
There is  a natural epimorphism of $^*$-algebras
  $$\pi: U_q({\mathfrak g})\to A_W({\mathfrak g}, q, \ell)$$
   that has support the simple representations of $U_{q}({\mathfrak g})$ in the Weyl alcove $\Lambda^+(q)$
and   that satisfies
\begin{equation}\label{relation_between_the_two_coproducts}P\pi\otimes\pi(\Delta^U(a))=\Delta(\pi(a))=\pi\otimes\pi(\Delta^U(a))P, \quad\quad P=\Delta(I).\end{equation}
The unitary coboundary matrix $\overline{R}^{A_W({\mathfrak g}, q, \ell)}$  of $A_W({\mathfrak g}, q, \ell)$ 
arises naturally from the (non-semisimple) Hermitian
  coboundary matrix $\overline{R}^U$ for 
  $U_q({\mathfrak g})$
   through Tannakian reconstruction via
\begin{equation}\label{the_formula_for_the_weak_coboundary}\overline{R}^{A_W({\mathfrak g}, q, \ell)}={\pi}\otimes{\pi}(\overline{R}^U)\Delta(I),
\end{equation}
  where $\overline{R}^U$ has been constructed in Theorem \ref{U_q_as_a_Hermitian_ribbon_h}.  
 This is a $2$-cocycle for $A_W({\mathfrak g}, q, \ell)$ 
as defined in Def.
\ref{def_2-cocycle}.  
A similar relation holds for the corresponding $R$-matrices.   \end{thm}

 The definition of weak tensor functor is given in \ref{wtf}. Axioms of {\it unitary coboundary weak (quasi-)Hopf algebras}, their categorical characterization,
  their main properties are introduced and studied in
 Sects. \ref{18}, \ref{DR}, \ref{19}. These are special kind of $\Omega$-involutive weak quasi-Hopf algebras in the sense of Sect.
 \ref{8}, with $\Omega$ explicitly associated to the $R$-matrix and the ribbon structure. The notion of antipode of Kac type extending a notion used in the setting of compact quantum groups \cite{CQGRC} to the $\Omega$-involutive case,  is given  in Def.
 \ref{Kac_type_definition}. Theorem \ref{main_wh} will be proved in Sect. \ref{20}.
 
 Establishing the stated  properties for $A_W({\mathfrak g}, q, \ell)$ is not immediately linked to the properties
 of $U_q({\mathfrak g})$ because the ideal of this algebra defining $A_W({\mathfrak g}, q, \ell)$ is not a co-ideal.
 
\subsection{Classification of type $A$ braided fusion categories}
 The following result is our first application of the weak Hopf algebras $A_W({\mathfrak sl}_N, q, \ell)$ (first
 constructed in \cite{CP}) to the study  of Problems \ref{problem3}--\ref{problem4}.
 It gives an equivalence   between fusion categories  in   type $A$ with different methods, by classifying ribbon fusion categories with the same fusion rules
 as ${\mathcal C}({\mathfrak sl}_N, q, \ell)$ by their ribbon structure.
The Grothendieck ring of ${\mathcal C}({\mathfrak sl}_N, q, \ell)$ is denoted by $R_{N, \ell}$.
 
To establish equivalences between the fusion categories of vertex operator algebras (which are not inherently unitary) and those of quantum groups or conformal nets, we replace the   unitarity assumption with the weaker condition of pseudounitarity,
 see Sect. \ref{KW}.

 \begin{thm}\label{classification_type_A}
  Let ${\mathcal C}$ and ${\mathcal C}'$ be  pseudo-unitary  ribbon fusion categories with  ribbon structures $\theta$ and $\theta'$,
  assumed positive for $N$ even, and with        based Grothendieck rings isomorphic to $R_{N, \ell}$ with $N+1\leq\ell<\infty$.
  Let   
  $f:{\rm Gr}({\mathcal C})\to{\rm Gr}({\mathcal C}')$ be a based ring  isomorphism such that  
   for each irreducible $\rho\in{\mathcal C}$, $\theta'_{\rho'}=\theta_\rho$ 
 where $\rho'$ is an irreducible in ${\mathcal C}'$ in the class of $f[\rho]$.
 Then there is an equivalence of ribbon braided tensor categories ${\mathcal F}: {\mathcal C}\to{\mathcal C}'$ inducing $f$. 
 If the categories are unitary, ${\mathcal F}$ may be chosen unitary.
  \end{thm}

   The previous result is 
  based on Kazhdan-Wenzl theory \cite{KW}.   The main conceptual argument of  proof consists of Theorem \ref{classification_type_A} in
extending   an analogous result by
 Neshveyev and Yamashita \cite{NY_twisting} in the setting of compact quantum groups to our weak Hopf algebras.
 This result is a first important step in our study of determining braided tensor structures of categories
 related to conformal field theory, in that is gives an insight of how the braided symmetry determines the associator, a feature
 that will emerge again in our next results with  more  direct methods with the uniqueness Theorems \ref{claim0} and \ref{claim1}
 for braided tensor structures in a semisimple category with a tensor product bifunctor, a generating object admitting a positive integer valued weak dimension function.

 \begin{rem} Theorem \ref{corollary_of_KW} applies to establish in particular equivalence of the categories ${\mathcal C}({\mathfrak sl}_N, q, \ell)$
 and ${\rm Rep}(V_{({{\mathfrak sl}_N})_k})$ giving an alternative proof of  KLF theorem  in the type $A$ case. The proof gives information on the structure of the equivalence.
 The theorem   applies also to other  
 approaches to braided tensor structures in categories motivated by conformal field theory, provided
 they have the prescribed fusion rules and ribbon structure.  \end{rem}

 Kazhdan-Wenzl theory classifies all the ribbon fusion categories with prescribed Gro\-then\-dieck rings of type $A$ up to a twist of the associativity morphisms. These twists have been classified
 by Kazhdan and Wenzl  in \cite{KW}, see also Sect. \ref{KW} for more details.
 It follows that the categories in the statement correspond to representation categories of the weak Hopf algebra
 $A_W({\mathfrak sl}_N, q, \ell)$ up to a twist. Then as anticipated by extending   an argument by Neshveyev and Yamashita to our weak Hopf algebra
we prove that the twist of the associativity morphisms is trivialized by the presence of the braided symmetry. It is  remarkable
 that their arguments in the setting of compact quantum groups extend to weak Hopf algebras, because our algebras are
 not co-associative, strictly speaking. This contributes to the idea that their lack of co-associativity
 is limited. We take the opportunity to note that  similar phenomena also occurred in the development of the theory of weak Hopf algebras
 Sects. \ref{6}, \ref{11},
 and 
 in the study of amenability in the setting of weak quasi-Hopf algebras and weak Hopf algebras, see Sect. \ref{13}.
 Moreover, the relevance of weak Hopf algebras that first emerged in the proof of Theorem \ref{classification_type_A}, it emerged again
 in the direct proof of Theorem \ref{Zhu_as_a_compatible_unitary_wqh}.
 
 Moreover, Neshveyev and Yamashita proved in the same paper that the categories with twisted associator in the case where $q$ is real,
 correspond to compact quantum groups, and classified these groups, as deformation of the defining compact quantum group. An analogue of this result for the fusion categories
 ${\mathcal C}({\mathfrak sl}_N, q, \ell)$ and our weak Hopf algebras  $A_W({\mathfrak sl}_N, q, \ell)$ has been studied in 
 \cite{Giannone_thesis}.
 \smallskip
 
\subsection{Vertex operator algebras and the KLF equivalence} Let  
 ${\rm Rep}(V_{{\mathfrak g}_k})$ be the module category of the   affine vertex operator algebra $V_{{\mathfrak g}_k}$ at positive integer level $k$. Let $A(V_{{\mathfrak g}_k})$ be the Zhu algebra.
 
 To connect with the vertex operator algebra fusion category, we develop an analogue of Drinfeld-Kohno Theorem, see Theorem \ref{Drinfeld_Kohno}, which leaves the structure of compatible unitary cobundary weak quasi-Hopf algebra invariant.
 We apply this theorem to our $A_W({\mathfrak g}, q, \ell)$ to connect to the Zhu algebra $A(V_{{\mathfrak g}_k})$ via Wenzl continuous de-quantization curve described in \cite{Wenzl}.

The following result solves Frenkel-Zhu Problem \ref{problem3} for affine vertex operator algebras at positive integer levels.
This problem has roots  in the question by Moore and Seiberg in conformal field theory, see Problem \ref{problem2}.

     \begin{thm}\label{Zhu_as_a_compatible_unitary_wqh} 
     With the same notation as above, with $q=e^{i\pi/\ell}$ with $\ell=d(h^\vee+k)$, we have:
 \begin{itemize}
\item[(a)]
the Zhu algebra $A(V_{{\mathfrak g}_k})$ admits a canonical structure of compatible unitary coboundary weak quasi-Hopf
$C^*$-algebra. The $C^*$ structure of $A(V_{{\mathfrak g}_k})$, regarded as a quotient of
$U({\mathfrak g})$, is induced by the classical compact real form of ${\mathfrak g}$.
The weak quasi-bialgebra structure is induced by
 a weak quasi-tensor  structure $(F_0, G_0)$ for Zhu forgetful $^*$-functor
$$ Z: {\rm Rep}(A(V_{{\mathfrak g}_k}))\to{\rm Hilb}$$
obtained via a natural isomorphism with Wenzl functor $W$ and a Drinfeld twist. 

More in detail, 
the triple $(Z, F_0, G_0)$ enriched over the category of Hilbert spaces,  is
obtained by transferring the untwisted structure of $A_W({\mathfrak g}, q, \ell)$ via Drinfeld-Kohno theorem \ref{Drinfeld_Kohno} and Wenzl continuous path on the arc ${\mathbb T}_{q, 1}$ connecting $q$ to $1$ clockwise through Kashiwara-Lusztig specialized canonical bases (Cor. \ref{corollary_of_positivity}), that is there is an isomorphism $\phi$ and twist $T$ such that
$$A(V_{{\mathfrak g}_k})\simeq^\phi (A_W({\mathfrak g}, q, \ell))_T.$$ 
\item[(b)]
The structure in (a) satisfies the strong unitarity property $F_0^*=G_0$ on pairs $(V_\lambda, V)$ and $(V, V_\lambda)$where $V_\lambda$ is a simple representation and $V$ is the fundamental representation 
(cf. part b) of Theorem \ref{Drinfeld_Kohno}).

\end{itemize}

 \end{thm}

We have used the same notation $V$ and $V_\lambda$ for 
    the  
  representations of the vertex operator algebra $V_{{\mathfrak g}_k}$ and of the Zhu algebra  $A(V_{{\mathfrak g}_k})$  corresponding to the irreducible objects  of the quantum group 
  fusion   categories. Parts (a) and (b) of Theorem \ref{Zhu_as_a_compatible_unitary_wqh} will be proved in 
  Subsect. \ref{30.2}.
  \medskip

  The following result is a consequence of   Theorem \ref{Zhu_as_a_compatible_unitary_wqh} and a combination of our previous  Tannakian constructions and  of the general unitarizability and tensor category constructions Sect \ref{12} applied to Zhu's linear equivalence.

\begin{cor}\label{cor_Zhu_as_a_compatible_unitary_wqh} For any complex simple Lie algebra ${\mathfrak g}$ and all the positive integer levels $k$ (with $k\geq2$ for ${\mathfrak g}=E_8$),
the linear category ${\rm Rep}(V_{{\mathfrak g}_k})$ becomes a   unitary modular tensor category ${\rm Rep}_{\rm T}(V_{{\mathfrak g}_k})$ with the structure induced by Zhu's linear equivalence 
$$Z: {\rm Rep}(V_{{\mathfrak g}_k})\to {\rm Rep}_{\rm T}(A({V_{{\mathfrak g}_k}})),$$ where  ${\rm Rep}_{\rm T}(A({V_{{\mathfrak g}_k}}))$ is endowed with the unitary modular fusion category structure
induced by the   structure of $A(V_{{\mathfrak g}_k})$
   constructed   in    Theorem \ref{Zhu_as_a_compatible_unitary_wqh} (a) via the Drinfeld twist $T$ and Wenzl's continuous curve.
   Therefore these constructions give  unitary ribbon braided tensor equivalences
   \begin{equation}\label{3.2}{\rm Rep}_{\rm T}(V_{{\mathfrak g}_k}) \xrightarrow{Z} {\rm Rep}_{\rm T}(A({V_{{\mathfrak g}_k}}))\xrightarrow{\ref{Zhu_as_a_compatible_unitary_wqh},  (a)} {\rm Rep}(A_W({\mathfrak g}, q, \ell)\xrightarrow{{\rm TK\  ribbon \ equiv}}{\mathcal C}({\mathfrak g}, q, \ell).
   \end{equation}
\end{cor}

  The middle equivalence in (\ref{3.2})  is application of the explicit symmetry provided by our analogue of Drinfeld-Kohno theorem
\ref{Drinfeld_Kohno} described in Theorem \ref{Zhu_as_a_compatible_unitary_wqh} (a). 
  The twist $T$ described by that theorem and the categorical counterpart 
  given by the middle equivalence above  is given by a Drinfeld twist explicitly defined by the action of the
  $R$-matrix, more precisely by the braided and ribbon structure in the two settings, that allow to see   in the vertex operator algebra setting
    the  same structure as that in the quantum group setting, via the naturally isomorphic algebras $A_W({\mathfrak g}, q, \ell)$ and $A(V_{{\mathfrak g}_k})$.
    \medskip

  \begin{rem}
  Before proceeding further we explicitly notice that      Corollary \ref{cor_Zhu_as_a_compatible_unitary_wqh} gives a self-contained structure
  of unitary rigid  ribbon braided tensor category on the module category ${\rm Rep}_{\rm T}(V_{{\mathfrak g}_k})$ of the affine vertex operator algebra $V_{{\mathfrak g}_k}$ for all Lie types and levels, equivalence with the quantum group fusion category by direct construction.


In particular,   we have considered both simply laced and non-simply laced Lie types in a uniform way. Furthermore, we do not use the Knizhnik-Zamolodchikov equations.
Our proof of rigidity of ${\rm Rep}_{\rm T}(V_{{\mathfrak g}_k})$  does not depend on the Verlinde formula.  
 \end{rem}
 
 \subsection{Comparing our braided tensor structure of   ${\rm Rep}_{\rm T}(V_{{\mathfrak g}_k})$ with Huang-Lepowsky braided tensor structure ${\rm Rep}_{\rm HL}(V_{{\mathfrak g}_k})$}

Let us next endow the linear category ${\rm Rep}(V_{{\mathfrak g}_k})$  of modules of the vertex operator algebra
with the ribbon braided tensor structure  ${\rm Rep}_{\rm HL}(V_{{\mathfrak g}_k})$ introduced by Huang and Lepowsky \cite{HL_tensor_products_of_modules}--\cite{Huang2}.
The following result, compares the structures of ${\rm Rep}_{\rm HL}(V_{{\mathfrak g}_k})$  and ${\rm Rep}_{\rm T}(V_{{\mathfrak g}_k})$, is our solution to Huang's Problem \ref{problem4}.

\begin{thm}\label{Finkelberg_HL}  
We have that:
\begin{itemize}

\item[(a)] 
The representation category ${\rm Rep}(A(V_{{\mathfrak g}_k}))$  
admits a natural structure of rigid
ribbon braided tensor category ${\rm Rep}_{\rm HL}(A(V_{{\mathfrak g}_k}))$  
obtained by transporting  Huang-Lepowsky structure to ${\rm Rep}(A(V_{{\mathfrak g}_k}))$.
In this way,
Zhu's linear equivalence 
$$Z: {\rm Rep}_{\rm HL}(V_{{\mathfrak g}_k}) \to{\rm Rep}_{\rm HL}(A(V_{{\mathfrak g}_k}))$$ 
admits a natural structure of
ribbon braided tensor equivalence.

\item[(b)] 
Let $A(V_{{\mathfrak g}_k})$ be endowed with the twisted unitary compatible coboundary weak quasi-Hopf structure following  Theorem \ref{Zhu_as_a_compatible_unitary_wqh}   (a), and let, as before,
  ${\rm Rep}_{\rm T}(A(V_{{\mathfrak g}_k}))$ be endowed with the corresponding 
 unitary rigid ribbon braided tensor structure. 
  Then 
$$Z: {\rm Rep}_{\rm HL}(V_{{\mathfrak g}_k}) \to{\rm Rep}_{\rm T}(A(V_{{\mathfrak g}_k}))$$
  preserves the ribbon for all objects.

Moreover ${\rm Rep}_{\rm T}(A(V_{{\mathfrak g}_k}))$ and ${\rm Rep}_{\rm HL}(A(V_{{\mathfrak g}_k}))$ have the same 
tensor product bifunctor.   
The equivalence structure maps of  $Z$ as in (a) w.r.t.  the  modified structure for the target category
satisfy the braided tensor equivalence equations for   
the braid morphisms for pairs $(V_\lambda, V)$ and $(V_\lambda, V)$, and the associativity morphisms for triples
$$(V_\lambda, V, V), \quad (V, V_\lambda, V), \quad (V, V, V_\lambda)$$ with $V_\lambda$ an arbitrary irreducible object.

\item[(c)]  If  ${\mathfrak g}$
 is of one of the Lie types
   $A$, $B$, $C$, $D$, $G_2$     then Zhu's equivalence in (b) is a ribbon braided tensor equivalence.
It follows that the composition 
\begin{equation}\label{31.2}{\rm Rep}_{\rm HL}(V_{{\mathfrak g}_k}) \xrightarrow{Z}{\rm Rep}_{\rm T}(A(V_{{\mathfrak g}_k}))\xrightarrow{\ref{Zhu_as_a_compatible_unitary_wqh},  (a)}{\rm Rep}(A_W({\mathfrak g}, q, \ell))\xrightarrow{{\rm TK\  ribbon \ equiv}}{\mathcal C}({\mathfrak g}, q, \ell)\end{equation}
is a ribbon braided  tensor equivalence by application of the indicated theorems.

\end{itemize}

\end{thm}  

 Theorem \ref{Finkelberg_HL} establishes the equivalence with Huang-Lepowsky 
structure for simply laced and non simply laced Lie types.
By  bypassing the analytic monodromy of the KZ equations entirely and the Verlinde formula, we establish an alternative construction of the Huang--Lepowsky braided tensor structure. 

Our proof   starts with the  the work of Wenzl \cite{Wenzl} alongside   the description of the {\it primary fields} by Wassermann \cite{Wassermann} for the minimum energy functor  of affine Lie algebras  in the type $A$ case.
 This leads   to the  description by Frenkel and Zhu of the
fusion rules for the affine vertex operator algebras \cite{Frenkel_Zhu}. The comparison of the two tensor product befunctprs results   in the tensor product Theorem \ref{tensor_product}. Connection between the work by Frenkel and Zhu with the basic tensor product bifunctor  
by  Huang and Lepowsky is given in Theorem \ref{FZ_condition}. 
The  comparison   led us to apply the description of the primary fields and also of the braided symmetry studied
by Toledano Laredo in \cite{Toledano_laredo}. For the connection with the vertex operator algebra setting, see Remark \ref{HL_braiding}.
Part (a)  
is proved in Sect. \ref{32}.  
The corresponding abstract result is explicitly stated in Sect. 15, see also   Sect. 7 of
\cite{On_a_problem_posed_by_Huang}.
The rest of part (b) summarizes the comparison on the associators described in Theorems  \ref{30.4} and  \ref{HL_is_of_CFT_type} and of the braiding morphisms Theorems \ref{braiding_Zhu_algebra},
\ref{Toledano_laredo},  Remark \ref{HL_braiding}.


  In  part (c), one needs to extend the coincidence of the braiding and associativity morphisms from the special cases considered in part (b) to all such morphisms. Establishing this result for a given Lie algebra allows one to pass from local to global coincidence between our analytic tensor structure and the   Huang-Lepowsky structure. Because the corresponding quantum group fusion category is known to be rigid, this global coincidence provides a new, direct proof of rigidity for the Huang-Lepowsky braided tensor category.

To prove part (c), we establish our main abstract uniqueness theorem, Theorem \ref{claim1} (a generalization of Theorem \ref{claim0}) for braided tensor structures within a category equipped with a common tensor product bifunctor and a multiplicative generating object (which, in our application, is the fundamental representation). This theorem operates via cohomological reduction: by exploiting the pentagon and hexagon equations, we show that the global braided tensor structure is completely and uniquely determined by its values on a  subclass of special    objects, directly treatable in the application. The hypotheses of the theorem are specifically formulated to match the explicit data available in our setting. Consequently, our final equivalence between our quantum group based structure and Huang-Lepowsky structure, holds for all classical Lie types and $G_2$, driven by state-of-the-art results concerning the generalized Schur-Weyl dualities between the fundamental representation and the braid group within the quantum group fusion category. Ultimately, this approach yields an independent, algebraic derivation of the Huang-Lepowsky braided tensor structure that entirely bypasses the historical analytic difficulties of the KZ equations and the Verlinde formula.

Building on the results established in part (b), in this paper we   verify   the hypotheses of Theorem \ref{claim1} concerning direct comparison between the braiding and associativity morphisms on the side of the quantum group fusion category for its original structure and for  the Huang-Lepowsky structure transferred to the Zhu algebra, via our Drinfeld twist. The further verification regarding  the duality property of the braid group (i.e. the braid group generates the morphism spaces between tensor powers of the fundamental representation) is performed
on the quantum group side on the algebra $A_W({\mathfrak g}, q, \ell)$. This
 requires separate verification across the various Lie types. To accomplish this,   we compare the original ribbon braided tensor category structure of the quantum group fusion category with the untwisted braided tensor category structure coming from the Zhu algebra endowed with Huang-Lepowsky structure. We do this untwist operation applying the same twist as before on the module category of the Zhu algebra equipped with the Huang-Lepowsky braided tensor structure from part (a). At this stage, the final   assumption of Theorem \ref{claim1} is the requirement that the centralizer algebras of the truncated tensor powers of the generating representation $V$ in the quantum group fusion category
 ${\mathcal C}({\mathfrak g}, q, \ell)$ are generated by the representation of the braid group. For the classical Lie types and $G_2$, this verification is detailed in Section 11 of \cite{On_a_problem_posed_by_Huang}. It relies on well-known extensions of quantum Schur-Weyl duality for types $A$ and $C$, alongside more recent results concerning the generators and relations of the centralizer algebras for the vector representation in types $B$, $D$, and $G_2$.  
 The detailed proof of Theorems \ref{claim0} and \ref{claim1} and verification of the duality between the fundamental representation and the braid group is deferred to \cite{On_a_problem_posed_by_Huang}.

\begin{rem} Theorem 
  \ref{TheoremUnitaryBraidRepAffine} gives a
unitarizability result for  ${\rm Rep}_{\rm HL}(V_{{\mathfrak g}_k})$.  The proof is based on our approach to  KLF Theorem \ref{Finkelberg_HL} for Huang-Lepowsky ribbon braided tensor structure, proved in later sections. This unitarizability result   is further specified  by   Corollary \ref{cor_Zhu_as_a_compatible_unitary_wqh} combined with 
 Theorem \ref{Finkelberg_HL}  as follows.

Specifically,  Corollary \ref{cor_Zhu_as_a_compatible_unitary_wqh} is obtained as an application of the general unitarizability and tensor category construction
criterion given in Theorem \ref{transportability} for linear $C^*$-categories. To apply this theorem we choose
${\mathcal C}^+$ to be    the linear $C^*$-category   
of unitary representations of $V_{{\mathfrak g}_k}$ regarded as a unitary vertex operator algebra (see  \cite{GuiI} and references therein for the notion of unitary representations of unitary vertex operator algebras), 
 $A=A(V_{{\mathfrak g}_k})$,  the unitary compatible weak quasi-Hopf algebra structure of the Zhu algebra 
 given by
Theorem \ref{Zhu_as_a_compatible_unitary_wqh}, and  
${\mathcal C}={\rm Rep}(A(V_{{\mathfrak g}_k}))$ regarded as endowed with the   ribbon braided tensor structure induced
by $A(V_{{\mathfrak g}_k})$. The linear $C^*$-structure of 
${\rm Rep}(A(V_{{\mathfrak g}_k}))$ is compatible with that of ${\mathcal C}^+$ via Zhu's linear equivalence ${\mathcal E}^+=Z$ and its inverse ${\mathcal S^+}$    obtained by Zhu \cite{Zhu}, since it corresponds to the classical compact real form
of ${\mathfrak g}$.
Theorem \ref{transportability} then applies and makes
 ${\mathcal C}^+$ into a unitary  tensor category    with tensor product bifunctor described in the proof of \ref{transportability} via the pair 
 (${\mathcal E}^+$, ${\mathcal S}^+$). The unitary tensor category structure obtained in this way
 defines ${\rm Rep}_{\rm T}(V_{{\mathfrak g}_k})$ and    is    also   rigid and unitary ribbon   braided as so is ${\rm Rep}(A(V_{{\mathfrak g}_k}))$.
 
 Thanks to the ribbon braided tensor equivalence between ${\rm Rep}_{\rm T}(V_{{\mathfrak g}_k})$ and ${\rm Rep}_{\rm HL}(V_{{\mathfrak g}_k})$ by  (c) of Theorem \ref{Finkelberg_HL} for the specified Lie types, the unitarizability then holds also for Huang-Lepowsky ribbon braided tensor category ${\rm Rep}_{\rm HL}(V_{{\mathfrak g}_k})$.

Applicative aspects of hermitian forms of vertex operator algebras and connections with hermitian forms on the Zhu algebra are discussed in Sect. \ref{VOAnets2}, see in particular the general criterion given in Theorem
  \ref{TheoremUnitaryRepVOA}. More examples of unitarization of module categories of   vertex operator algebras fulfilling
  the needed assumptions are discussed in conclusion of the same section.
\end{rem}

   \medskip

 \section{Outline of the paper}\label{History}
 
 The first part of the present work is a presentation of  the general theory of tensor categories and weak quasi-Hopf algebras, and general Tannakian results both in the algebraic and $C^*$ setting, with or without braiding and ribbon structures. These topics are the content of Sects. \ref{2}--\ref{5}, \ref{6}--\ref{11}.

 Section \ref{6} is dedicated to developing the notion    of  {\it weak Hopf algebra} in a new sense,
as the Hopf algebra counterpart of Mack-Schomerus weak quasi-Hopf algebras.
Our definition differs from the usual notion of weak Hopf algebra, because the coproduct is not strictly coassociative,  but the defining
fiber functor is weakly monoidal and the associator explicitly depends only on the weakness of the coproduct.
This includes our previous examples
  $A_W({\mathfrak sl}_N, q, \ell)$ of \cite{CP} and the new examples $A_W({\mathfrak g}, q, \ell)$ constructed later in   this paper,  in a general theoretical context.

In   Sects.   \ref{2}, \ref{3} we review tensor categories and their functors. 
Because CFT  fusion categories are   unitary, and are equipped with unitary representations of the braid group, any   reconstruction of quantum gauge group
must be equipped with a compatible $C^*$-structure, which requires considering a new and general notion of an $\Omega$-involution.

Sects \ref{4}, \ref{5}, \ref{6}--\ref{10}
are devoted to  Tannakian duality theorems, both in the algebraic and unitary case.
We introduce  a general  notion of
    $\Omega$-involutions for weak Hopf or weak quasi-Hopf algebras extending and developing the case of quasi-Hopf algebras
considered by Gould and Lekatsas \cite{quasi_star}.

The notion  of $\Omega$-involution describes
the unitary structure of our first examples in \cite{CP} in the type $A$ case, with a specific positive
matrix $\Omega$  explicitly associated to the ribbon structure and the braiding of   $A_W({\mathfrak sl}_N, q, \ell)$,
describing positivity of the inner product of tensor product of representations, and
it induces
 unitary representations of the braid group.
   This special kind of $\Omega$-involution is a main interest in later sections of our paper.
 
  Sect. \ref{8} is dedicated to the   development of a theory of weak (quasi-)Hopf $C^*$-algebras with general $\Omega$-involutions,  which 
 turns out to unify the theory of the discrete duals of of the compact quantum groups, where $\Omega$ is trivial and the coproduct is unital, with the unitary
  structure of the quantum group fusion category  ${\mathcal C}({\mathfrak g}, q, \ell)$ studied by Wenzl, associated to
    the hermitian form of $U_q({\mathfrak g})$ at roots of unity  \cite{Wenzl}. 
    
$\Omega$-involutions arise in  Tannakian constructions  of quasi-Hopf or weak quasi-Hopf algebras  from tensor $C^*$-categories
when the structure maps are not   unitary. This property is inherited by tensor categories equivalent to tensor $C^*$-categories with such maps.

     Sects. \ref{11}, \ref{12}, \ref{16} study unitarizability of braided tensor categories with methods of the associated
    weak quasi-Hopf algebras.

The main result of Sect. \ref{12} is Theorem \ref{unitarizability}, giving a method for unitarizing braided tensor categories  knowing that the category is tensor equivalent to a unitary braided tensor category.

 Section  \ref{13} regards   amenability, its extends to weak quasi-Hopf algebras known results from several sources.   A main result establishes that
 when the categorical dimensions are not integer and the unitary category corresponds to the representation category of a weak Hopf algebra then necessarily the corresponding $\Omega$-involution can not be trivial. Thus in cases where one has a trivial $\Omega$-involution then the corresponding algebra can not be a weak Hopf algebra, see Corollary \ref{unitarity_obstruction}.
 
   A main application of the abstract unitarizability results of Sect. \ref{12}  is the first  solution
 to the problem posed by Galindo on uniqueness of the unitary structure, Theorem \ref{Galindo} in Section \ref{14}.
 
 We analyse
      different natural integral weak dimension functions in Sect.\ref{15}.
  We apply    Tannakian results to the unitary quantum group fusion categories   in   Sect. \ref{74}. These results are formulated fixing the  choice of an
  integral weak (i.e. submultiplicative)   dimension functions.
 
  The second main application of the unitarizability criterion \ref{unitarizability}, is   Theorem \ref{BFM}, giving the  unitarization of the Beilinson-Feigin-Mazur category $\tilde{\mathcal O}_\ell$ for all the cases where it is known that Finkelberg's equivalence with ${\mathcal C}({\mathfrak g}, q, \ell)$ holds \cite{Finkelberg}, \cite{Finkelberg_erratum}. This theorem relies on   Wenzl'unitarity result of ${\mathcal C}({\mathfrak g}, q, \ell)$ in \cite{Wenzl}.

  The proof of the unitarizability of  $\tilde{\mathcal O}_\ell$, Theorem \ref{BFM}, depends on the Finkelberg equivalence, and this  was our original   motivation to look for a direct proof. We do this in the rest of the paper.
 
  We discuss our   unitarization results of module categories 
${\rm Rep}(V_{{\mathfrak g}_k})$
   of affine vertex operator algebras at positive integer levels $k$ as a first result on the subject, see     Sect. \ref{VOAnets2}.
The preceding       Sect. \ref{VOAnets} is dedicated to preliminaries on vertex operator algebras and their modules.   This result depends on our main result Theorem \ref{Finkelberg_HL}.

          Section \ref{VOAnets3}   deals with weak quasi-Hopf algebras  obtained via Tannakian duality and conformal nets, along the same  lines.

In the rest of the work we develop direct, canonical constructions of weak quasi-Hopf algebras naturally associated to the
     inner fusion structure of the categories via their coproduct.

To approach our study   of a direct proof of  KLF equivalence theorem, we first study the type $A$ case, from a classification perspective
thanks to Kazhdan-Wenzl theory \cite{KW}.

We do this in Sect. \ref{KW} that can be read independently of the rest of the paper, and the result is both interesting in its own and quite enlightening for the rest of the paper because Kazhdan-Wenzl methods are abstract, and thus may be applied to concrete
examples beyond quantum groups, including the case of vertex operator algebras.
 The main result is stated as Theorem \ref{classification_type_A}
that gives a classification of these fusion categories   in terms of the representation ring and the ribbon structure.  

Thanks to the braided symmetry, the emerging categories from the classification can only correspond to 
to  ${\mathcal C}({\mathfrak sl}_N, q, \ell)$.
In particular, the proof shows how any associativity morphism when realized 
over ${\mathcal C}({\mathfrak sl}_N, q, \ell)$ becomes trivial thanks to the constraint given by the braided symmetry. This result extends a remarkable result by 
Neshveyev and Yamashita in the setting of compact quantum groups, Prop. 4.4 in \cite{NY_twisting}.
In our case, we    apply their arguments   to the weak Hopf algebras $A_W({\mathfrak sl}_N, q, \ell)$ of \cite{CP}.
The result highlights how the presence of the braided symmetry determines  the associator and trivializes it for the quantum group fusion category, setting up the framework for the arbitrary Lie types handled at the very end of the paper with a mixture of direct comparison, cohomological reduction, and quantum Schur-Weyl duality.

Our study for the general Lie types beyond type $A$ is   starts from Sects. \ref{73} and occupies the rest of the paper
  (with the exception of Sect. \ref{KW}  as said).

 To make our approach  to unitarization of module categories of affine vertex operator algebras useful,   we  first need to understand
   KLF equivalence 
   in the setting of Huang-Lepowsky tensor product structure.
   Our main result is Theorem \ref{Finkelberg_HL} and its proof is a main aim of this paper.

 Sects.  \ref{17}, \ref{18}, \ref{DR}, \ref{19} develop central categorical and algebraic notions of {\it unitary coboundary weak quasi-Hopf algebras} and their {\it compatible} and {\it hermitian}   counterparts. While hermiticity is a weakening of unitarity to the case where the $^*$-involution is not a $C^*$-algebra involution, and is useful to include $U_q({\mathfrak g})$ as the main example,
 the mentioned {\it compatibility} property is a strengthening to the notion of hermitian/unitary coboundary weak quasi-Hopf algebra,
 to the case where, just like  $U_q({\mathfrak g})$, the coproduct anticocommutes with the $^*$-involution.  The weaker case admits a nice Tannakian characterization established in Sect. \ref{DR} extending to the coboundary symmetry the symmetry property of functors between symmetric tensor categories, a property envisaged for a long time in our previous work.  
 
 Sect. \ref{20} is devoted to the construction of canonical unitary compatible coboundary  weak Hopf $C^*$-algebras $A_W({\mathfrak g}, q, \ell)$ for all
 the Lie types associated to ${\mathcal C}({\mathfrak g}, q, \ell)$, which, as $^*$-algebras, are semisimple quotients of $U_q({\mathfrak g})$, by an ideal which is not a coideal. Their coproduct  naturally arises from the fusion tensor product of 
 ${\mathcal C}({\mathfrak g}, q, \ell)$.
 We study their naturally arising $\Omega$-involution
associated to their coboundary matrices $\overline{R}$, extending   those considered in \cite{CP}, \cite{Ciamprone} for the type $A$ case with different methods.

Similarly to the type $A$ case, the equivalence between ${\mathcal C}({\mathfrak g}, q, \ell)$ and the corresponding ${\rm Rep}(V_{{\mathfrak g}_k})$ relies on these algebras $A_W({\mathfrak g}, q, \ell)$, but in these cases we adopt  a direct comparison. 
Similarly to the type $A$ case, the interplay between the braiding
and the associativity morphisms are relevant. But differently from the type $A$ case, we develop an abstract cohomological reduction for    the pentagon and hexagon equations of a braided tensor category with. a generating object, reminding the classification of braided pointed fusion categories \cite{EGNO}.
We shall be able   to explain structural properties of our unitarization of  ${\rm Rep}(V_{{\mathfrak g}_k})$ and its
  deep symmetry with the quantum group fusion category  by the end of the paper.

In Sect.  \ref{21} we use these canonical unitary coboundary weak Hopf $C^*$-algebras  $A_W({\mathfrak g}, q, \ell)$   
 to transport, via the construction of a Drinfeld twist of   $A_W({\mathfrak g}, q, \ell)$, the unitary structure to ${\rm Rep}(V_{{\mathfrak g}_k})$. To describe this unitary structure, we pass through   the Zhu algebra $A(V_{{\mathfrak g}_k})$, which becomes a unitary coboundary weak quasi-Hopf algebra with 3-coboundary associator, by our analytic analogue of the Drinfeld-Kohno theorem \ref{Drinfeld_Kohno}.
Up to this point we have constructed, via transport from the Zhu algebra and Zhu's linear equivalence and its inverse \cite{Zhu}, the structure of a rigid unitary modular category on ${\rm Rep}(V_{{\mathfrak g}_k})$
for all the Lie types and positive integer levels, with   fusion multiplicities transported from the quantum group fusion category via  $A_W({\mathfrak g}, q, \ell)$  and the Drinfeld twist.
These constitutes  parts (a), (b) of Theorem \ref{Zhu_as_a_compatible_unitary_wqh}.
We are left to   compare our ribbon braided structure  with Huang-Lepowsly structure.

Thus we pass to introduce Huang-Lepowsky tensor structure in generality, 
following  historical presentation starting from   the setting of affine Lie algebras, where one can see the first connections between
{\it primary fields} 
with the Verlinde fusion structure of corresponding quantum groups at roots of unity. This is the content of Sects. \ref{22}--\ref{33}
which also contains   the proof of parts (a), (b)  Theorem \ref{Finkelberg_HL}.
    
     Historically, the comparison between the fusion tensor product of quantum groups and of affine Lie algebras
has roots {\it primary fields}, see  e.g. \cite{pinzari_survey}, or Appendix A in \cite{On_a_problem_posed_by_Huang}.  This 
connection emerges in the works by
Tsuchiya and Kanie, Kohno, Wassermann and relation with the work of Wenzl \cite{Tsuchiya_Kanie}, \cite{Kohno2}, \cite{Wenzl}, \cite{Wassermann}. 

The notion of primary field and its relation to quanto groups plays a central role    in our proof.
Sect.  \ref{22} introduces these fundamental notion of primary field in the setting of loop groups or affine Lie algebras. It gives a central Theorem \ref{tensor_product} based on the work by   Fenkel and Zhu \cite{Frenkel_Zhu} and our construction of weak quasi-Hopf algebra on the Zhu algebra.
 This is the first connection between our truncated tensor product bifunctor $\boxtimes$ and what will later correspond to Huang-Lepowsky bifunctor $\boxtimes_{\rm HL}$, via a known connection between intertwining operators in the setting of vertex operator algebras and primary fields.
 This section also contains Corollary  \ref{Toledano_laredo} that is a reduction of the final result to the associator, using work by Toledano-Laredo \cite{Toledano_laredo}. 
 
The methods that we have used
so far (mainly our version of Drinfeld-Kohno theorem \ref{Drinfeld_Kohno}) to prove parts  (a) and (b) of Theorem \ref{Zhu_as_a_compatible_unitary_wqh} and a direct comparison to prove parts (a), (b)  Theorem \ref{Finkelberg_HL},
show
coincidence of the two tensor product bifunctors, of the ribbon structure, and the braiding  
for certain special pairs   of variables where one 
of them is the fundamental representation $V$. We also know that special
associativity morphisms coincide when two out of the three variables are the fundamental representation, by direct verification.
We call these special associativity morphisms {\it of CFT type}, and we develop an abstract cohomology theory in Sect. \ref{5+} to determine the remaining braiding and associativity morphisms   uniquely using the pentagon and hexagon equations. The main results are Theorems \ref{claim0} and \ref{claim1}, whose proofs are deferred to \cite{On_a_problem_posed_by_Huang}.

 This general comological reduction applies to the case where
the category has  a generating object in a given linear category with a tensor product bifunctor,
with braiding morphisms that generate the centralizer algebras of tensor powers of the generating object. This applies to the quantum group fusion categories for the classical Lie types $ABCD$ and $G_2$, by results known in the literature on generalization of Schur-Weyl duality, see \cite{Chari_Pressley} for the classical vector representation, and due to   Wenzl and his collaborators for the considered fundamental representations.
This application is discussed  in Sect. 11 in \cite{On_a_problem_posed_by_Huang} as well.

In summary, our approach systematically bridges quantum groups   with the analytical realities of conformal field theory. By substituting strict coassociativity with our weak quasi-Hopf framework, we bypass the historical reliance on the analytic monodromy of the Knizhnik--Zamolodchikov equations, passage to negative shifted levels, use of the Verlinde formula to establish rigidity, use of rigidity to establish braided tensor equivalence. We also expand our results to all positive integer levels and to all Lie types for our UMFC structure on ${\rm Rep}_{\rm T}(V_{{\mathfrak g}_k})$ and to the classical Lie types and $G_2$ for the Huang-Lepowsky structure. This also allows us to uniformly establish the uniqueness of unitary structures across tensor categories 
(Theorem \ref{Galindo}) and settle the long-standing unitarization problem for the Beilinson--Feigin--Mazur category (Theorem \ref{BFM}), together with the other problems mentioned in Subsect. \ref{1.3} for the WZW model. When classification is unavailable from the outset, we utilize the weaker condition of pseudounitarity alongside Kazhdan--Wenzl theory to fully classify Type $A$ ribbon fusion categories and complete the operator-algebraic proof of the Kazhdan--Lusztig--Finkelberg equivalence.

Perspectives for future research have been discussed in \cite{pinzari_survey}. Among them, the extension of the theory of quantum homogeneous spaces \cite{PRergodic} and induced representations \cite{PRinduction} of Woronowicz compact quantum groups \cite{Wor_su2, Wor, Wor_compact_quantum_groups} to unitary coboundary weak Hopf $C^*$-algebras, to study extensions of the WZW model in the setting of vertex operator algebras and conformal nets.
 \smallskip

\section{Preliminaries on tensor categories and their  functors, generating object}\label{2}

In this section we recall the   basic terminology concerning  tensor categories and unitary tensor categories. Our main references are \cite{EGNO, Mueger3} and \cite{CQGRC} respectively.
 We also give  the main definitions of certain functors between these categories.
 The most familiar notion is that of   tensor functor but we  need suitable weak generalizations, known in the literature as {\it quasi-tensor functors} and more importantly for us their weak versions, the {\it weak quasi-tensor functors}.
  We also introduce  a new notion, that of     {\it weak tensor functor} between tensor categories as  a slight generalisation of   notions already considered in the literature. We shall describe a cohomological interpretation in the setting of weak quasi-Hopf algebras later on. Finally, we introduce a notion of unitarity for weak quasi-tensor functors between unitary tensor categories and discuss a unitarization procedure for general weak quasi-tensor functors which
  will be fruitful later on.

  All categories in this paper will be essentially small, thus  they will  admit a small skeleton. 
     The morphism space from an object $\rho$ to $\sigma$ 
  is denoted by $(\rho, \sigma)$.
By a {\it linear  category} ${\mathcal C}$ we mean a category  whose morphism spaces are 
 complex vector spaces and such   that composition is bilinear.

The notion of {\it   semisimple category} is central  in this paper, we briefly recall the definition   directing our attention to linear categories,     we refer  the reader to Ch. 1 in \cite{EGNO} for details.

  A   {\it linear additive} category ${\mathcal C}$ is a linear category with a zero object $0$, that is
  $(0, 0)=0$, and direct sums,
  that is for any pair of objects $\rho$, $\sigma\in{\mathcal C}$ there is an object $\tau\in{\mathcal C}$ and morphisms
  $S\in(\rho, \tau)$,   $T\in(\sigma, \tau)$, $S'\in(\tau, \rho)$, $T'\in(\tau, \sigma)$ such that $S'S=1$, $T'T=1$,
  $SS'+TT'=1$. The object $\tau$ is defined up to isomorphism and denoted $\rho\oplus\sigma$. A  {\it  linear abelian} category is a
linear additive category with extra structure.
The central additional notion is that of kernel and symmetrically of cokernel of a morphism. For a morphism  $A\in (\rho, \sigma)$  
  the kernel ${\rm Ker}(A)$  
  is   an object $k$ and a  morphism $K\in(k, \rho)$ such that $AK=0$, and universal with this property. 
  Kernels and cokernels    are assumed to exist for every morphism, among other things. 
    A {\it subobject} of an object $\rho$
  is an object $\sigma$ together with a morphism $S\in(\sigma, \rho)$ with ${\rm Ker}(S)=0$. 
  An object $\rho$ is called
   {\it simple}, or {\it irreducible}, if $\rho\neq0$ and the only subobjects are $0$ and $\rho$.
         
        It follows from  Schur's Lemma,  see e.g. Lemma 1.5.2 in \cite{EGNO}  and Prop. 5.4.5 in \cite{Cohn} that
in a linear abelian category  
   with   finite dimensional morphism spaces,
when $\rho$ and $\sigma$ are simple, $(\rho, \sigma)$ is either the trivial vector space or it is formed by scalar multiples of a unique isomorphism, it follows that 
   $(\rho, \rho)={\mathbb C}1$.
     In our paper, all our categories will have finite dimensional morphism spaces.

   A {\it semisimple} category is a linear abelian category  such that every nonzero object is a finite direct sum of simple objects, the decomposition is unique up to isomorphism.

A splitting idempotent, or a summand, of an object $\rho$ is an object $\sigma$, an idempotent $E\in(\rho, \rho)$ together with morphisms $S\in(\sigma, \rho)$, $S'\in(\rho, \sigma)$ such that  $S'S=1$, $SS'=E$.  In particular, $\sigma$ is a subobject of $\rho$. For example, a direct sum $\rho\oplus\sigma$ as previously defined has $\rho$ and $\sigma$ as summands defined by complementary idempotents.
In a semisimple category every idempotent splits, thus every subobject is a summand.

  The next notion is that of {\it tensor category}. We  follow Sect. 1.2 in \cite{Mueger3}, and the notion of monoidal category of Ch. 2 in \cite{EGNO} except for 
we assume the linear structure.

\begin{defn}\label{pre_tensor_category}
By a {\it  pre-tensor category}  we mean a linear category ${\mathcal C}$
  endowed with a tensor product operation $\otimes$, which is a bilinear  bifunctor 
${\mathcal C}\times {\mathcal C}\to{\mathcal C}$, a distinguished tensor unit object $\iota$
 satisfying the {\it unit axioms}, that is the functors $\rho\to\rho\otimes \iota$ and $\rho\to\iota\otimes\rho$ are autoequivalences of ${\mathcal C}$.
 \end{defn}
 
 \begin{defn}\label{tensor_category} A {\it tensor category} is a pre-tensor category
 endowed with
natural isomorphisms
 $\alpha_{\rho,\sigma,\tau}: (\rho\otimes\sigma)\otimes\tau\to\rho\otimes(\sigma\otimes\tau)$. 

The associativity morphisms $\alpha_{\rho, \sigma, \tau}$   are required to satisfy the {\it pentagon equation}
\begin{equation}\label{pentagon_equation}
\xymatrix@C=1em{
  ((\rho\otimes \sigma)\otimes \tau)\otimes \upsilon\ar[d]_{\alpha}\ar[r]^{\alpha\otimes {1}}
  &(\rho\otimes(\sigma\otimes \tau))\otimes \upsilon\ar[r]^{\alpha}
  &\rho\otimes((\sigma\otimes \tau)\otimes \upsilon)\ar[d]^{{1}\otimes\alpha}\\
  (\rho\otimes \sigma)\otimes(\tau\otimes \upsilon)\ar[rr]_{\alpha}&&\rho\otimes(\sigma\otimes (\tau\otimes \upsilon))
}
\end{equation}

  \end{defn}
 By Sect. 2.9 in \cite{EGNO}
 one can identify $\rho\otimes\iota$ and $\iota\otimes\rho$ by a simple passage which
uses only the unit isomorphisms, in this way   $\iota$ becomes  strict, meaning that   
 $\iota\otimes\rho=\rho\otimes\iota=\rho$ for every object and $1_\iota\otimes T=T\otimes 1_\iota=T$ for every morphism $T$.
 To simplify our discussion,    we shall assume that   $\iota$ is  strict in our abstract results, and we shall tacitly
 use this passage  in   our applications where it is not natural to work with a strict unit, e.g. Sect.  \ref{VOAnets}, \ref{VOAnets2}.
 Moreover, we assume the normalization condition $\alpha_{\rho, \iota, \tau}=1$.

 Alternatively, one may adopt the definition of monoidal units originally given by Saavedra, following the treatment in \cite{Kock}. This viewpoint has been adopted in
  \cite{Giannone_thesis}. (See also
 Remark \ref{closely_related} for a brief description of the aims of the thesis and   more closely related papers in the setting of Kazhdan-Wenzl theory,   Sect. \ref{KW}.) 
  
  Here we just point out an interesting and quite handy feature: in terms of   Saavedra notion, a tensor functor $\mathcal{F}:\mathcal{C}\to\mathcal{C}'$ is compatible with units exactly if it sends units into units, or equivalently if it sends any one unit, say $\iota$ into a unit. Moreover, given any other unit $\iota'$ of $\mathcal{C}'$, the compatibility is expressed by a unique isomorphism between $\iota'$ and $\mathcal{F}(\iota)$.

We shall only deal with  tensor    categories for which the tensor unit satisfies    $(\iota,\iota)={\mathbb C}1$.
The category is called strict if the tensor unit is strict and associativity morphisms are identity.
With abuse of language, we shall also regard the category ${\rm        Vec}$ of finite dimensional vector spaces as strict.
  Whenever convenient, for a given semisimple category ${\mathcal C}$,
we fix a set ${\rm Irr}({\mathcal C})$ of simple objects in ${\mathcal C}$ such that
every simple object of ${\mathcal C}$ is isomorphic to exactly
one element of ${\rm Irr}({\mathcal C})$. This can be done because we
are assuming that
${\mathcal C}$ is essentially small.
A semisimple tensor category  with
    finitely many inequivalent irreducible objects  will be  called {\it finite semisimple}.
  If ${\mathcal C}$ is in addition rigid,   it is a (complex) {\it  fusion} category \cite{ENO}.\bigskip
  
  The following definition will play an important role in our paper.
    \begin{defn}\label{Def_generating_object}
  Let ${\mathcal C}$ be a semisimple pre-tensor category. An object $\rho\in{\mathcal C}$ is called {\it generating} if  
   every simple object of ${\mathcal C}$ is a summand of some tensor power of $\rho$ defined by some 
   parenthesisation.

  \end{defn}

A tensor power of $\rho$ is an object defined as the result of iterative applications of the bifunctor $\otimes$
to pairs with entries $\iota$ or $\rho$. 
We define the order of   a tensor power of $\rho$ in a natural way.
  Different    parenthesisations of the same order of two tensor powers of $\rho$ are equivalent in a tensor category,
  thus it suffices to verify the generating property on a preferred choice of    parenthesisation in this case, and is independent of the choice of the associativity morphisms.
  
  Functors   between linear categories 
   are ${\mathbb C}$-linear maps between morphism spaces.  
   
\begin{defn}   A linear  functor ${\mathcal F}:{\mathcal C}\to{\mathcal C}'$ between linear categories is called a  {\it linear equivalence}
if there is a linear functor, called a {\it quasi-inverse}, ${\mathcal G}:{\mathcal C}'\to{\mathcal C}$ such that ${\mathcal F}{\mathcal G}$ and ${\mathcal G}{\mathcal F}$ are naturally isomorphic to the identity functors of ${\mathcal C}'$ and ${\mathcal C}$ respectively.
\end{defn}

\begin{rem}\label{equivalent_def_of_equivalence}
It is well known (Theorem 1  in IV.4 of \cite{MacLane}) that a linear functor ${\mathcal F}: {\mathcal C}\to{\mathcal C}'$ is a linear equivalence if and only if it is full and faithful (i.e. bijective between the morphism spaces) and essentially surjective (every object of ${\mathcal C}'$ is isomorphic to one in the image of 
${\mathcal F}$.) We shall use these definitions interchangeably. When ${\mathcal C}$ and ${\mathcal C}'$ are semisimple,
this is equivalent to the property that  $\{{\mathcal F}(\rho), \  \rho\in{\rm Irr}({\mathcal C})\}$
is a complete set of pairwise non-isomorphic simple objects in ${\mathcal C}'$.
\end{rem}

The following notion of {\it weak quasi-tensor functor}  for a tensor category was   introduced by H\" aring-Oldenburg in  \cite{HO}
in connection with the study of duality for weak quasi-Hopf algebras. We note that the definition  ignores the associativity structure of ${\mathcal C}$ and ${\mathcal C}'$ and thus 
is meanginful in the 
setting of pre-tensor categories.

\begin{defn}\label{wqtf} Let ${\mathcal C}$ and ${\mathcal C'}$ be pre-tensor categories. A  
{\it weak quasi-tensor} functor is defined by a
${\mathbb C}$-linear functor  ${\mathcal F}:{\mathcal C}\to{\mathcal C}'$
satisfying  ${\mathcal F}(\iota)=\iota$  together with two morphisms
$F_{\rho, \sigma}: {\mathcal  F}(\rho)\otimes {\mathcal F}(\sigma)\to {\mathcal F}({\rho\otimes \sigma})$ and $G_{\rho, \sigma}: {\mathcal F}({\rho\otimes \sigma})\to {\mathcal F}(\rho)\otimes {\mathcal F}(\sigma)$ satisfying
\begin{equation}\label{normalized_iota}
F_{\iota, \rho}=F_{\rho,\iota}=1_{{\mathcal F}(\rho)},\quad\quad G_{\iota, \rho}=G_{\rho,\iota}=1_{{\mathcal F}(\rho)},\end{equation}
\begin{equation}\label{right_inverse}
F_{\rho, \sigma}\circ G_{\rho, \sigma}=1_{{\mathcal F}({\rho\otimes \sigma})}\end{equation}
\begin{equation}\label{naturality}
F_{\rho', \sigma'}\circ {\mathcal F}(S)\otimes{\mathcal F}(T)={\mathcal F}(S\otimes T)\circ F_{\rho,\sigma},\quad
{\mathcal F}(S)\otimes{\mathcal F}(T) \circ G_{\rho,\sigma}=G_{\rho',\sigma'}\circ {\mathcal F}(S\otimes T)
\end{equation}
for  objects $\rho$, $\sigma$, $\rho'$, $\sigma'\in{\mathcal C}$ and morphisms $S:\rho\to\rho'$, $T:\sigma\to\sigma'$.
\end{defn}
Property (\ref{naturality})   expresses naturality of $F$ and $G$ in $\rho$ and $\sigma$, while the right inverse condition (\ref{right_inverse}) implies that
\begin{equation}\label{Definition_of_P}
P_{\rho,\sigma}=G_{\rho,\sigma}\circ F_{\rho, \sigma}: {\mathcal F}(\rho)\otimes {\mathcal F}(\sigma)\to 
{\mathcal F}(\rho)\otimes {\mathcal F}(\sigma)\end{equation} is an idempotent satisfying 
$$F_{\rho,\sigma}\circ P_{\rho, \sigma}=F_{\rho,\sigma}, \quad\quad P_{\rho,\sigma}G_{\rho, \sigma}=G_{\rho, \sigma}.$$
If $P_{\rho,\sigma}=1_{{\mathcal F}(\rho)\otimes {\mathcal F}(\sigma)}$
for all $\rho$, $\sigma$
  (i.e. all $F_{\rho,\sigma}$ are isomorphisms), we recover the notion of quasi-tensor functor  of \cite{Drinfeld_quasi_hopf, Majid4}.
   
 \begin{defn}\label{monoidal_transformation} Let ${\mathcal F}$, ${\mathcal F}': {\mathcal C}\to{\mathcal C'}$ be
 two weak quasi-tensor functors between pre-tensor categories defined by  
  $(F_{\rho,\sigma}, G_{\rho,\sigma})$,  $(F'_{\rho,\sigma}, G'_{\rho,\sigma})$, respectively.
 A natural transformation $\eta:{\mathcal F}\to{\mathcal F}'$   is called {\it monoidal} if $\eta_\iota=1_\iota$ and if
 $$F'_{\rho,\sigma}\circ\eta_{\rho}\otimes\eta_\sigma=\eta_{\rho\otimes\sigma}\circ F_{\rho,\sigma},\quad G'_{\rho,\sigma}\circ\eta_{\rho\otimes\sigma}=\eta_\rho\otimes\eta_\sigma\circ G_{\rho,\sigma}.$$
 \end{defn}

 The following definition is motivated by the requirement of   compatibility between the functor and the associativity morphisms.

\begin{defn}\label{wtf} Let ${\mathcal C}$ and ${\mathcal C}'$ be tensor categories with associativity morphisms $\alpha$ and $\alpha'$ respectively.
A {\it weak tensor} functor  is a weak quasi-tensor functor ${\mathcal F}:{\mathcal C}\to{\mathcal C}'$ for which the associated natural transformations $F_{\rho,\sigma}$, $G_{\rho, \sigma}$ satisfy for all objects $\rho$, $\sigma$, $\tau$,
\begin{equation}
{\mathcal F}(\alpha_{\rho,\sigma,\tau})=F_{\rho, \sigma\otimes\tau}\circ 1_{{\mathcal F}(\rho)}\otimes F_{\sigma,\tau}\circ \alpha'_{{\mathcal F}(\rho), 
{\mathcal F}(\sigma), {\mathcal F}(\tau)}\circ G_{\rho, \sigma}\otimes 1_{{\mathcal F}(\tau)}\circ G_{\rho\otimes\sigma,\tau}\label{wt1}
\end{equation}
\begin{equation}{\mathcal F}(\alpha_{\rho,\sigma,\tau}^{-1})=  F_{\rho\otimes\sigma,\tau}\circ  F_{\rho,\sigma}\otimes 1_{{\mathcal F}(\tau)}\circ 
{\alpha'}^{-1}_{{\mathcal F}(\rho), 
{\mathcal F}(\sigma), {\mathcal F}(\tau)}\circ 1_{{\mathcal F}(\rho)}\otimes G_{\sigma, \tau} \circ G_{\rho, \sigma\otimes\tau}.\label{wt2}\end{equation}

\end{defn}

\begin{defn}\label{tf}
In the case that  all $F_{\rho, \sigma}$ are isomorphisms then $G_{\rho, \sigma}=F_{\rho, \sigma}^{-1}$ thus only one of the   equations (\ref{wt1}) and (\ref{wt2}) suffices and we   recover 
the   notion of a {\it tensor functor}  \cite{EGNO, Kassel, Mueger3, CQGRC}. 

\end{defn}

 \begin{rem} \label{sufficient_conditions_for_associativity_morphisms} a) In Sect. \ref{20} we shall construct a weak tensor structure $(F, G)$ on a certain faithful functor ${\mathcal F}: {\mathcal C}({\mathfrak g},   
 q, \ell)\to{\rm Hilb}$ on 
  the strict fusion tensor category ${\mathcal C}({\mathfrak g},   
 q, \ell)$ associated to quantum groups at roots of unity to ${\rm Hilb}$, the category of finite dimensional Hilbert spaces.
 To this aim, the notion of {\it negligible morphism} in the category of tilting modules for a suitable quantum group at root of unity $U_q({\mathfrak g})$ plays a central role.

b) In the paper, see e.g. Sect. \ref{5+},   \ref{21}, \ref{32}, \ref{33}, we shall also consider a   notion weaker than that of weak tensor functor, in that
we shall need to work with functors for which (\ref{wt1}) or (\ref{wt2}) hold for a subcollection ${\mathcal V}$ of triples
of objects, which is sufficiently large in a way that will be made precise.

c) The following seems an interesting insight gained from   the notion of weak tensor functor, and more specifically
from the utility of having the two equations (\ref{wt1}), (\ref{wt2}), which were first noticed  in the course of the construction
of the examples of weak  Hopf algebras in \cite{CP}, and then again in the course of the construction of the generalizations
of these examples to all Lie types in Sect. \ref{20}. These properties were shown to follow
 from the properties of negligible morphisms in the representation category of $U_q({\mathfrak g})$ at root of unity, but have not been not used   in \cite{CP} or for the most part of  this paper, except for the
 following possibly non trivial application.   
 
If we have a faithful functor ${\mathcal F}:{\mathcal C}\to{\rm Vec}$ on a semisimple pre-tensor category
${\mathcal C}$ and a weak quasi-tensor structure $(F, G)$ for ${\mathcal F}$
 and an invertible  natural transformation $$\alpha_{\rho, \sigma, \tau}:  
 (\rho\otimes\sigma)\otimes\tau\to\rho\otimes(\sigma\otimes\tau)$$    
making ${\mathcal C}$ into a tensor category and we  have an understanding in a specific situation of validity of one preferred equation   among the two  (\ref{wt1}), (\ref{wt2}) on a {\it fixed
triple} $(\rho, \sigma, \tau)$, and if moreover $\beta(\rho, \sigma, \tau)$ is a fixed invertible
morphism associated to the same fixed triple $(\rho, \sigma, \tau)$ and acting a morphism of ${\mathcal C}$ between the same objects as $\alpha_{\rho, \sigma, \tau}$,
 $$\beta_{\rho, \sigma, \tau}: 
 (\rho\otimes\sigma)\otimes\tau\to\rho\otimes(\sigma\otimes\tau)$$    
and satisfies the same property as $\alpha$ among (\ref{wt1}) or (\ref{wt2}),
then $\alpha_{\rho, \sigma, \tau}=\beta_{\rho, \sigma, \tau}$ on that fixed triple by faithfulness of ${\mathcal F}$ and uniqueness of the inverse.
If moreover one finds a collection ${\mathcal V}$ of   triples $(\rho, \sigma, \tau)$ for which this property is known to hold 
on each element of ${\mathcal V}$,  which is sufficiently large to determine $\alpha$ uniquely
on all triples of objects of ${\mathcal C}$ (for example ${\mathcal V}$ is as   in Sect. \ref{5+}, Corollary \ref{determination_of_associator_with_generating_object}, Theorem \ref{claim0}) then $\beta$ admits a unique extension to associativity morphisms
of ${\mathcal C}$ satisfying the axioms of a tensor category, and the extension is given by $\alpha$.
 In our main application, we wish to achieve this setting with ${\mathcal F}(\alpha)$ and   $(F, G)$   defined by a Drinfeld twist construction 
from a situation as in a) and moreover   ${\mathcal F}(\beta)$ with the same $(F, G)$, arise from conformal field theory ,
see Sect. \ref{33}. 
\end{rem}
\begin{defn}\label{tensor_equivalence2}
A {\it tensor equivalence}  between tensor categories ${\mathcal C}$ and ${\mathcal C}'$ is a tensor functor ${\mathcal E}: {\mathcal C}\to{\mathcal C}'$ which is an equivalence of linear categories.\end{defn}
It is known that a quasi-inverse ${\mathcal G}:{\mathcal C'}\to{\mathcal C}$ may be chosen tensorial  
and the natural transformations $1_{{\mathcal C}'}\to{\mathcal F}{\mathcal G}$, $1_{\mathcal C}\to{\mathcal G}{\mathcal F}$ monoidal see Remark 2.4.10 in \cite{EGNO}. In particular, ${\mathcal G}$ is a tensor equivalence as well.

In general,   we are making no assumption on compatibility of $\alpha'$ with the two subobjects of $({\mathcal F}(\rho)\otimes {\mathcal F}(\sigma))\otimes 
{\mathcal F}(\tau)$ and ${\mathcal F}(\rho)\otimes({\mathcal F}(\sigma)\otimes {\mathcal F}(\tau))$ corresponding respectively to the right invertible maps $F_{\rho\otimes\sigma, \tau}\circ F_{\rho, \sigma}\otimes 1_{{\mathcal F}(\tau)}$ and $F_{\rho, \sigma\otimes\tau}\circ 1_{{\mathcal F}(\rho)}\otimes F_{\sigma,\tau}$.

A weak quasi-tensor functor monoidally isomorphic to a weak tensor functor is itself weak tensor.

The notion of weak  (quasi) tensoriality for a functor  applies to contravariant functors ${\mathcal C}\to{\mathcal C}'$ as well, but in this case the defining natural transformations are required to act as $F_{\rho,\sigma}: {\mathcal F}(\rho)\otimes {\mathcal F}(\sigma)\to
 {\mathcal F}(\sigma\otimes \rho)$,  $G_{\rho,\sigma}:  {\mathcal F}(\sigma\otimes \rho)\to {\mathcal F}(\rho)\otimes {\mathcal F}(\sigma)$ and the diagrams
(\ref{wt1}) and (\ref{wt2}) have to be appropriately modified.
  Equivalently, such functors
 may be regarded as covariant (quasi) tensor functors after replacing ${\mathcal C'}$ with the opposite category
 $({\mathcal C}')^{{\rm        op}}$ that is the category with same objects and morphisms, but opposed morphisms and  reversed   tensor products.

 We shall also consider categories with involutions and involution preserving functors. We shall follow \cite{DR1} and 
\cite{CQGRC}. These structures   will not be needed until Sect. \ref{9}.
 
\begin{defn} 
A {\it $^*$-category} is a linear category ${\mathcal C}$ endowed with an antilinear, contravariant, involutive functor $^*: {\mathcal C}\to{\mathcal C}$ acting trivially on objects. A {\it tensor $^*$-category} is a tensor category equipped with the structure of a $^*$-category satisfying $(S\otimes T)^*=S^*\otimes T^*$ for any pair of morphisms $S$, $T\in{\mathcal C}$.
The associativity morphisms are assumed unitary, $\alpha_{\rho,\sigma,\tau}^*=\alpha_{\rho, \sigma,\tau}^{-1}$.\end{defn}

\begin{defn}
A {\it $C^*$-category} is a $^*$-category where morphism spaces are Banach spaces such that the   norm satisfies $\|S\circ T\|\leq\|S\|\|T\|$ and $\|T^*\circ T\|=\|T\|^2$ for every pair of morphisms $S$, $T$ and $S^*S$ is positive
 (i.e. has positive spectrum) in the algebra $(\rho,\rho)$ for every morphism $S \in (\rho,\sigma)$.
Finally, a {\it tensor $C^*$-category} is a tensor $^*$-category which is also   a $C^*$-category with respect to the given $^*$-involution. 
\end{defn}

The positivity condition is equivalent to the existence of $S'\in(\rho,\rho)$ such that $S^*S=S'^*S'$.
 It follows in particular that $(\rho, \rho)$ is a $C^*$-algebra for any object $\rho$.
In a $C^*$-category,  two isomorphic objects $\rho$, $\sigma$ are called {\it  unitarily isomorphic}  if there is a unitary $U\in(\rho, \sigma)$, that is $U^*U=1$, $UU^*=1$.   An {\it orthogonal summand} of $\rho$ is a summand 
defined by 
a selfadjoint idempotent $E\in(\rho, \rho)$ which is the range of an isometry (there is $S\in(\sigma, \rho)$ such that $S^*S=1$ and $SS^*=E$). An {\it orthogonal direct sum} $\rho\oplus\sigma$ is defined by isometries $S_1\in(\rho, \rho\oplus\sigma)$, $S_2\in (\sigma, \rho\oplus\sigma)$ such that
 $S_1S_1^*+S_2S_2^*=1$. 
 
It follows from the positivity of   $T^*T$ that  a left  invertible morphism
$T\in(\sigma, \rho)$ admits polar decomposition in ${\mathcal C}$.  Thus $S=T(T^*T)^{-1/2}\in(\sigma, \rho)$ 
is an isometry.    In particular,  two isomorphic objects $\rho$, $\sigma$ are also unitarily isomorphic.  

It also follows that a summand or a direct sum is isomorphic to an orthogonal one.  Indeed,
 by Prop. 4.6.2 in \cite{Blackadar}
  every idempotent in a unital $C^*$-algebra  is similar to a selfadjoint idempotent.
Thus  a summand  $\sigma$ of $\rho$ up to isomorphism corresponds to a selfadjoint  idempotent in $E\in(\rho, \rho)$, and  it follows that
polar decomposition of the corresponding morphism $S\in(\sigma, \rho)$ gives the
needed isometry. Similarly, the defining complementary idempotents  of a direct sum $\rho\oplus\sigma$ may be assumed  selfadjoint and it follows that the direct sum is orthogonal.

In particular, 
  a semisimple $C^*$-category
has orthogonal summands and direct sums.
    It is also easy to see   that the positivity condition of $T^*T$ follows from the other properties of a $C^*$-category and  existence  of orthogonal direct sums, cf. Ch. 2 in \cite{CQGRC}.  

\begin{defn} A $^*$-functor ${\mathcal F}:{\mathcal C}\to{\mathcal C}'$ between $^*$-categories is a linear functor 
satisfying  ${\mathcal F}(T^*)={\mathcal F}(T)^*$ for all morphisms $T\in{\mathcal C}$.
If ${\mathcal C}$ and ${\mathcal C}'$ are tensor $^*$-categories, a $^*$-functor endowed with a 
weak quasi tensor  structure
 will be called a 
 {\it  weak quasi tensor $^*$-functor}.
\end{defn}

Let ${\mathcal F}$ be a weak quasi tensor $^*$-functor defined by   $F_{\rho, \sigma}$, $G_{\rho, \sigma}$. 
Then  the adjoint pair $F'_{\rho, \sigma}=G^*_{\rho, \sigma}$, $G'_{\rho, \sigma}=
F^*_{\rho, \sigma}$ defines another weak quasi tensor structure on ${\mathcal F}$.

\begin{defn} A {\it $^*$-equivalence} 
 between   $^*$-categories ${\mathcal C}$ and ${\mathcal C}'$  is an equivalence compatible with the $^*$-structure, that is a 
$^*$-functor ${\mathcal E}:{\mathcal C}\to{\mathcal C}'$ admitting a quasi-inverse  ${\mathcal E}': {\mathcal C'}\to{\mathcal C}$ which is  a $^*$-functor with
 natural unitary transformations $\eta: 1\to {\mathcal E}{\mathcal E}'$ and $\eta': 1\to{\mathcal E}'{\mathcal E}$.
If ${\mathcal C}$ and ${\mathcal C}'$ are tensor $^*$-categories, ${\mathcal E}$ is a {\it tensor $^*$-equivalence} if   ${\mathcal E}$ and ${\mathcal E}'$ are tensor $^*$-functors.
 \end{defn}

We   note the following   $C^*$-version of  the characterisation of equivalences between categories of Remark
\ref{equivalent_def_of_equivalence}.

 \begin{prop}\label{MacLane}
 Let ${\mathcal F}: {\mathcal C}\to{\mathcal C'}$ be a $^*$-functor between $C^*$--categories. Then ${\mathcal F}$ is a $^*$-equivalence  if and only if it is a  $^*$-functor which is an equivalence of linear categories.
 If ${\mathcal C}$ and ${\mathcal C}'$ are tensor $C^*$-categories then ${\mathcal F}$ is a tensor $^*$--equivalence
 if and only if it is a $^*$-functor and a tensor equivalence.
 
 \end{prop}
 
 \begin{proof} We start with the definition of a linear equivalence as a full, faithful and essentially surjective functor ${\mathcal F}$, as in Remark \ref{equivalent_def_of_equivalence}.   Theorem IV.4.1 \cite{MacLane} constructs a 
   linear functor ${\mathcal G}: {\mathcal C}'\to{\mathcal C}$ and invertible natural transformations $\eta: 1\to {\mathcal F}{\mathcal G}$ and $\eta': 1\to{\mathcal G}{\mathcal F}$. We are thus left to show that we can always choose $\eta$ and $\eta'$ unitary and ${\mathcal G}$ a $^*$-functor. To this aim,  it is not difficult to adapt   the
 proof of that theorem  to the needed framework as follows. The isomorphisms $\eta_c$ defined there, corresponding to our $\eta$, may be chosen unitary passing to polar decomposition available with the $C^*$-structure of ${\mathcal C'}$. This implies that the quasi-inverse equivalence constructed
 there and denoted $T$,  in turn corresponding to ${\mathcal G}$, satisfies that $\eta: 1\to{\mathcal F}{\mathcal G}$ is a unitary natural transformation.
 This fact, together with the fact that ${\mathcal F}$ is a faithful $^*$-functor,  implies that   ${\mathcal G}$  is linear and  $^*$-preserving on morphism spaces. If $\eta'': 1\to{\mathcal G}{\mathcal F}$ is any invertible natural transformation, one of which is found in the same theorem,  then the unitary part in the polar decomposition $\eta'$ of $\eta''$ will be a unitary natural transformation between the same functors  thanks to the $^*$-preserving properties of the involved functors.
The last statement follows from the fact that when ${\mathcal C}$ and ${\mathcal C}$' are tensor $C^*$-categories
then we already know that we may construct a  tensorial quasi-inverse ${\mathcal G}$ and then we apply the first part of the proof.
 \end{proof}
 
 \begin{rem}\label{equivalent_def_of_equivalence_for_*_categories}
 We note that   a faithful and essentially surjective $^*$-functor between $^*$-categories ${\mathcal F}: {\mathcal C}\to{\mathcal C}'$ does not necessarily admit a quasi-inverse $^*$-functor. 
 An example is given by the immersion of the category ${\rm Hilb}$ of finite dimensional Hilbert spaces into the category ${\rm Herm}$ of finite dimensional Hermitian spaces.  This category will be introduced and studied starting with Sect. \ref{9}. For the subclass of semisimple $^*$-categories  we have the following useful criterion analogous
in analogy to Remark \ref{equivalent_def_of_equivalence}.
 Let ${\rm Irr}^u({\mathcal C})$ be a   set of pairwise unitarily inequivalent
 simple objects in ${\mathcal C}$ such that
every other simple object  is unitarily isomorphic to  
one element of ${\rm Irr}^u({\mathcal C})$.
A faithful   $^*$-functor between $^*$-categories ${\mathcal F}: {\mathcal C}\to{\mathcal C}'$ is a $^*$-equivalence if
and only if  the set of objects ${\mathcal F}(\rho)$ with $\rho\in{\rm Irr}^u({\mathcal C})$
is a complete set of pairwise unitarily inequivalent simple objects in ${\mathcal C}'$.
 
 \end{rem}

In the theory of $C^*$-tensor categories, or more generally of tensor $^*$-categories, we have the following notion of unitarity for a tensor functor and a tensor equivalence, see  \cite{CQGRC}.

 \begin{defn} Let ${\mathcal C}$ and ${\mathcal C}'$ be tensor $^*$-categories.
 A {\it unitary tensor functor} $({\mathcal F}, F, G=F^{-1})$, is a tensor $^*$-functor such that $F$ is unitary.
 A  {\it unitary tensor equivalence} is a tensor $^*$-equivalence which is unitary as a tensor $^*$-functor and with a 
 unitary quasi-inverse.
 \end{defn}
 
 Unitary tensor functors from $C^*$-tensor categories to ${\rm Hilb}$ arise as
forgetful functors of compact quantum groups   see e.g. \cite{CQGRC}.
As fusion categories do not in general admit tensor functors to ${\rm Vec}$, but always admit
weak quasi-tensor functors, we introduce a 
notion of unitarity in the following more general setting.

 We next begin to discuss a problem that has  relevence in how paper, that is how
to associate to a given weak quasi-tensor structure
$(F, G)$ another one that has in some sense a     more trivial unitary structure.
Historically, the first condition considered in the literature is $G=F^*$ and $G$ unitary see e.g. \cite{CQGRC}, or more generally
isometry \cite{HO}.

\begin{defn}\label{unitarity}  Let ${\mathcal C}$ and ${\mathcal C}'$ be tensor $^*$-categories.  A  {\it unitary weak quasi tensor  functor} is a weak quasi tensor $^*$-functor ${\mathcal F}:{\mathcal C}\to{\mathcal C}'$ 
defined by $(F, G)$ such that $F^*$ and $G$ are isometries. A {\it strongly unitary weak quasi tensor} functor we have
in addition $F^*=G$.
\end{defn}

 For quasi-tensor $^*$-functors we recover  the usual notion of unitarity $F_{\rho, \sigma}^*=F_{\rho, \sigma}^{-1}$. 

\begin{rem}
The   definition of unitarity may equivalently be formulated by   the properties
$$F_{\rho, \sigma}^*F_{\rho, \sigma}=P^*_{\rho, \sigma}P_{\rho, \sigma}, \quad
 G_{\rho, \sigma}G_{\rho, \sigma}^*=P_{\rho, \sigma}P^*_{\rho, \sigma}$$
where $P_{\rho,\sigma}$ is the idempotent defined in (\ref{Definition_of_P}). 
\end{rem}

 In general, if $(F, G)$ is unitary then we may have two new strongly unitary structures $(F, F^*)$ and $(G^*, G)$
 arising from $(F, G)$. However, in the $C^*$-case all these structures coincide. More precisely, we note the following
 simple result.

\begin{prop}\label{strongly_unitary_prop} Let ${\mathcal C}$ be a tensor $^*$-category, ${\mathcal C}'$ a tensor $C^*$-category
and $(F, G)$ a   weak quasi-tensor structure for a $^*$-functor ${\mathcal F}:{\mathcal C}\to{\mathcal C}'$.
Let $\rho$, $\sigma\in{\mathcal C}$ be a pair of objects. If $F_{\rho, \sigma}^*$  and
$G_{\rho, \sigma}$ are isometries then  $F_{\rho, \sigma}^*=G_{\rho, \sigma}$.
In particular, any unitary weak quasi-tensor structure is automatically strongly unitary.
\end{prop} 

\begin{proof}
We have that $F_{\rho, \sigma}G_{\rho, \sigma}=1_{{\mathcal F}(\rho\otimes\sigma)}=G_{\rho, \sigma}^*G_{\rho, \sigma}=
F_{\rho, \sigma}F_{\rho, \sigma}^*$. It follows that
$$G_{\rho, \sigma}^*(1_{{\mathcal F}(\rho)\otimes{\mathcal F}(\sigma)}-F_{\rho, \sigma}^*F_{\rho, \sigma})G_{\rho, \sigma}=G_{\rho, \sigma}^*G_{\rho, \sigma}-(F_{\rho, \sigma}G_{\rho, \sigma})^*(F_{\rho, \sigma}G_{\rho, \sigma})=0.$$
The $C^*$-property of ${\mathcal C}'$ implies
$(1-F_{\rho, \sigma}^*F_{\rho, \sigma})G_{\rho, \sigma}=0$  thus  $G_{\rho, \sigma}=F_{\rho, \sigma}^*$.
\end{proof}

To construct unitary   weak quasi-tensor structures from a given weak quasi-tensor structure, structure it is natural to try with polar decomposition.

We consider a 
  weak quasi tensor $^*$-functor $({\mathcal F}, F, G):{\mathcal C}\to{\mathcal C}'$ between $C^*$-tensor categories and 
 we describe a 
 a unitarization of the weak quasitensor structure $(F, G)$.
We set $$\Omega_{\rho, \sigma}:=F_{\rho, \sigma}^*\circ F_{\rho, \sigma}\in ({\mathcal F}(\rho)\otimes{\mathcal F}(\sigma), 
 {\mathcal F}(\rho)\otimes{\mathcal F}(\sigma)).$$ Note that $\Omega_{\rho, \sigma}$ is partially invertible (in the sense of Def. \ref{partial_inverse})
 with partial inverse  
  $$\Omega_{\rho, \sigma}^{-1}:=G_{\rho, \sigma}\circ G_{\rho, \sigma}^*\in ({\mathcal F}(\rho)\otimes{\mathcal F}(\sigma), 
 {\mathcal F}(\rho)\otimes{\mathcal F}(\sigma))$$ satisfying  
  $\Omega_{\rho, \sigma}^{-1}\Omega_{\rho, \sigma}=P_{\rho, \sigma}$ and 
  $\Omega_{\rho, \sigma}\Omega_{\rho, \sigma}^{-1}=P_{\rho, \sigma}^*$. Since they  are both positive, we may
take   the respective square roots  $\Omega_{\rho, \sigma}^{1/2}$ and $(\Omega_{\rho, \sigma}^{-1})^{1/2}$.

It can easily be shown  that these operators are partially invertible between the same idempotents as $\Omega_{\rho, \sigma}$ and $\Omega_{\rho, \sigma}^{-1}$ respectively.
If we in addition know that $({\Omega_{\rho, \sigma}^{-1}})^{1/2}$ is a left inverse
of  $\Omega_{\rho, \sigma}^{1/2}$, that is 
\begin{equation}\label{left_inverse_condition} ({\Omega_{\rho, \sigma}^{-1}})^{1/2}\Omega_{\rho, \sigma}^{1/2}=P_{\rho, \sigma}
\end{equation}  then we shall just write 
${\Omega_{\rho, \sigma}^{-1/2}}$ for $({\Omega_{\rho, \sigma}^{-1}})^{1/2}$.
We have
\begin{equation}\label{definition_polar_decomposition} F =S^* \circ \Omega^{1/2}, \quad\quad G=\Omega^{-1/2}\circ T,\end{equation} where $S$ and $T$ are isometries as $G$ is a right inverse of $F$.

\begin{prop}\label{polar_decomposition}
Let ${\mathcal F}: {\mathcal C}\to{\mathcal C}'$ be a weak quasi-tensor $^*$-functor between tensor $C^*$-categories defined by   $(F, G)$ such that $({\Omega_{\rho, \sigma}^{-1}})^{1/2}\Omega_{\rho, \sigma}^{1/2}=P_{\rho, \sigma}$
(e.g. $P=1$).
 Then 
\begin{itemize}
\item[{\rm        a)}]
the pair $(F', G')$, where $$F'=F\Omega^{-1/2}=S^*\Omega^{1/2}\Omega^{-1/2}, \quad\quad G'=\Omega^{1/2}G=\Omega^{1/2}\Omega^{-1/2}T,$$  is a unitary weak quasi-tensor structure for    ${\mathcal F}$, and therefore
strongly unitary, $F'=G'^*$, 
\item[{\rm        b)}] In particular, if $({\mathcal F}, F, G)$ is quasi tensor then $F'=S^*$, $G'=S=T$ is always
well defined and is a unitary   quasi-tensor structure,
\item[{\rm        c)}]
if ${\mathcal F}$ is  full and if $(F, G)$ is a tensor structure  then    $(S^*, S)$     is a unitary tensor structure for ${\mathcal F}$.
\end{itemize}

 \end{prop}

\begin{proof}
a)   It follows from $^*$-invariance of ${\mathcal F}$ that $\Omega_{\rho, \sigma}$ 
  is natural in $\rho$, $\sigma$, and from continuous functional calculus that  $\Omega_{\rho, \sigma}^{1/2}$ and $\Omega_{\rho, \sigma}^{-1/2}$ are natural as well, hence the same holds for $F'$ and $G'$. We have $F'G'=FPG=1$, so $(F', G')$ is a weak quasi tensor structure. The associated idempotent is given by
  $P':=G'F'=\Omega^{1/2}P\Omega^{-1/2}=\Omega^{1/2}\Omega^{-1/2}$. Furthermore 
  $F'F'^*=F(\Omega^{-1})^{1/2}(\Omega^{-1})^{1/2}F^*=FGG^*F^*=1$, $G'^*G'=G^*\Omega^{1/2}\Omega^{1/2}G=G^*F^*FG=1$,
   thus $(F' , G')$ is unitary, and by Prop. \ref{strongly_unitary_prop} also strongly unitary.
  c) In this case $F$, $G$ are invertible and   $G=F^{-1}$, thus   $P=1$,  $S$, $T$ are unitary and $S^*T=1$.
  d) Since
  $(G_{\rho, \sigma}^*\circ G_{\rho, \sigma})^{1/2}$ is a positive invertible element in the $C^*$-algebra
  $({\mathcal F}(\rho\otimes\sigma), {\mathcal F}(\rho\otimes\sigma))$ and ${\mathcal F}$ is full,
  we may  write 
  $(G_{\rho, \sigma}^*\circ G_{\rho, \sigma})^{1/2}={\mathcal F}(A_{\rho, \sigma})$ 
with $A_{\rho, \sigma}\in( \rho\otimes\sigma, \rho\otimes\sigma)$ positive, and  $G_{\rho, \sigma}=S_{\rho, \sigma}\circ{\mathcal F}(A_{\rho, \sigma})$ with $S$ unitary. It follows that ${\mathcal F}(1_\rho\otimes A_{\sigma,\tau})$
is positive by $^*$-invariance of ${\mathcal F}$ and also invertible by naturality of $G$.
 Furthermore, $$1_{{\mathcal F}(\rho)}\otimes G_{\sigma, \tau}\circ G_{\rho, \sigma\otimes\tau}=1_{{\mathcal F}(\rho)}\otimes S_{\sigma, \tau}\circ {\mathcal F}(1_{\rho})\otimes{\mathcal F}(A_{\sigma,\tau})\circ G_{\rho, \sigma\otimes\tau}=$$
  $$1_{{\mathcal F}(\rho)}\otimes S_{\sigma, \tau}\circ G_{\rho, \sigma\otimes\tau}\circ{\mathcal F}(1_\rho\otimes A_{\sigma, \tau})=1_{{\mathcal F}(\rho)}\otimes S_{\sigma, \tau}\circ S_{\rho, \sigma\otimes\tau}\circ B_{\rho, \sigma,\tau},$$
  where
  $B_{\rho, \sigma, \tau}:={\mathcal F}(A_{\rho, \sigma\otimes\tau})\circ {\mathcal F}(1_\rho\otimes A_{\sigma, \tau}).$
  A similar computation starting with the same element  but relying now on naturality of $S$ in place of $G$, see a), leads to conclude that 
  ${\mathcal F}(A_{\rho, \sigma\otimes\tau})$ and ${\mathcal F}(1_\rho\otimes A_{\sigma, \tau})$ commute, and this implies that   
  $B_{\rho, \sigma, \tau}$
 is   positive,   besides   invertible. In a similar way
 $G_{\rho, \sigma}\otimes 1_{{\mathcal F}(\tau)}\circ G_{\rho\otimes\sigma, \tau}=S_{\rho, \sigma}\otimes 1_{{\mathcal F}(\tau)}\circ S_{\rho\otimes\sigma, \tau}\circ C_{\rho, \sigma, \tau}$ for some other positive invertible morphism 
 $C_{\rho, \sigma, \tau}$. Inserting  these relations into    the    tensoriality diagram  
   $1_{{\mathcal F}(\rho)}\otimes G_{\sigma, \tau}\circ G_{\rho, \sigma\otimes\tau}\circ {\mathcal F}(\alpha)= \alpha' \circ G_{\rho, \sigma}\otimes 1_{{\mathcal F}(\tau)}\circ G_{\rho\otimes\sigma,\tau}$ gives another tensoriality diagram satisfied by
  $S$ in place of $G$ by unitarity of the associativity morphisms and uniqueness of polar decomposition.
  \end{proof}
 
 \begin{defn}\label{unitarization_of_a_functor}
Let $({\mathcal F}, F, G):{\mathcal C}\to{\mathcal C}'$ be a weak quasi-tensor $^*$-functor between tensor $C^*$-categories satisfying the left inverse property (\ref{left_inverse_condition}). Then the same functor ${\mathcal F}$ together  with the new unitary weak quasi-tensor structure $(F', G')$   defined in part a) of
Prop. \ref{polar_decomposition} will be called the {\it unitarization} of $({\mathcal F}, F, G)$.
\end{defn}

  \begin{rem} 
We would   like to warn the reader that   it is not clear to us whether   
(\ref{left_inverse_condition}) holds  in our main   late applications as in Sect.
\ref{19} and following. It follows that it is unclear whether the
polar decomposition construction of Prop. \ref{polar_decomposition} can be used. We shall need to develop a modification of the unitarization construction for a functor in Sect. \ref{19}. On the other hand, the unitarization of a functor will be fruitful for us in case
of full domains ($P=1$), see Sect. \ref{14}, where we shall discuss uniqueness of unitary structures in tensor categories.
\end{rem}

Part c) shows that in the important case of tensor $^*$-equivalence the unitarization gives a unitary tensor equivalence.  
We have the following consequence.

  \begin{cor}
 Two tensor $^*$-equivalent tensor $C^*$-categories are also unitarily tensor equivalent.
\end{cor}

 \begin{rem}
Note that we do not have  a statement about unitarization of a weak tensor $^*$-functor.
On this subject we shall see that the notion of unitary weak tensor $^*$-functor  is 
too strong   for unitary fusion categories of interest for us. Specifically, a unitary weak tensor $^*$-functor to the category of Hilbert spaces is automatically   tensor  for  large classes of semisimple unitary tensor tensor categories and the category   necessarily has an integer-valued dimension function, we refer to Corollary \ref{unitarity_obstruction} for details.
 It follows that the unitarization   of a weak  tensor $^*$-functor in general is only a unitary weak quasi-tensor $^*$-functor. 
  In Sect. \ref{20} we shall construct examples of weak tensor $^*$-functors associated to unitary fusion categories of quantum groups at roots of unity, and 
 part a)  of Prop. \ref{polar_decomposition}  will turn out   useful. 
 \end{rem}

\section{Rigidity,  weak tensor functors, braided symmetry, ribbon and modular category}\label{3}

In this brief section we recall the notion of rigidity, braided and ribbon tensor category and we show a simple result
that weak tensor functors are always compatible with rigidity.  
 
\begin{defn}
 Let  ${\mathcal C}$  be a tensor category with associativity morphisms $\alpha_{\rho, \sigma, \tau}\in ((\rho\otimes\sigma)\otimes\tau, \rho\otimes(\sigma\otimes\tau))$. An object $ {\rho}^\vee$ is a {\it right dual} of $\rho$ if there are morphisms  $d\in({\rho}^\vee\otimes \rho,\iota)$ and  $b\in(\iota, \rho\otimes {\rho}^\vee)$ satisfying the right duality equations
\begin{equation}\label{first_right}
 1_\rho\otimes d\circ\alpha_{\rho, {\rho}^\vee,\rho}\circ b\otimes 1_\rho=1_\rho,\end{equation}
\begin{equation}\label{second_right}
{d}\otimes 1_{{\rho}^\vee}\circ \alpha^{-1}_{{\rho}^\vee,\rho,{\rho}^\vee}\circ 1_{{\rho}^\vee}\otimes b=1_{{\rho}^\vee}.\end{equation} 
A {\it left dual} object ${}^{\vee}\rho$ is defined by morphisms $b'\in(\iota, {}^{\vee}{\rho}\otimes \rho)$, $d'\in(\rho\otimes {}^{\vee}{\rho},\iota)$ satisfying the  left duality equations
\begin{equation}\label{first_left} d'\otimes 1_\rho\circ\alpha^{-1}_{\rho, {}^\vee\rho,\rho}\circ 1_\rho\otimes b'=1_\rho,\end{equation}
\begin{equation}\label{second_left}1_{{}^\vee\rho}\otimes d'\circ\alpha_{ {}^\vee\rho,\rho, {}^\vee\rho}\circ b'\otimes 1_{{}^\vee\rho}=1_{{}^\vee\rho}.\end{equation}
A tensor category is called {\it rigid} if every object has  left and right duals.
\end{defn}

The following   facts are well known:  another right dual $(\tilde{\rho}, \tilde{b}, \tilde{d})$ is isomorphic to $\rho^\vee$, the isomorphism is 
\begin{equation}\label{xi}
\xi:=\tilde{d}\otimes 1_{\rho^\vee}\circ 1_{\tilde{\rho}}\otimes b: \tilde{\rho}\to\rho^\vee,\end{equation} and similarly for left duals.
 If $\rho$ and $\sigma$ have right   duals $\rho^\vee$ and $\sigma^\vee$, then so does $\rho\otimes\sigma$ , and it is given by $\sigma^\vee\otimes\rho^\vee$ via the morphisms  
$d_{\rho\otimes\sigma}=d_\sigma\circ1_{{\sigma}^\vee}\otimes(d_\rho\otimes 1_\sigma)\circ\alpha\in(({\sigma}^\vee\otimes{\rho}^\vee)\otimes(\rho\otimes\sigma),\iota)$, $b_{\rho\otimes\sigma}=\alpha'\circ1_\rho\otimes(b_\sigma\otimes 1_{{\rho}^\vee})\circ b_\rho\in(\iota, (\rho\otimes\sigma)\otimes({\sigma}^\vee\otimes{\rho}^\vee))$, where $\alpha$ and $\alpha'$ are suitable associativity morphisms.
  
\begin{defn} A {\it right duality} is defined by the choice
 of a right dual $(\rho^\vee, b_\rho, d_\rho)$ for each object $\rho$ such that $\iota^\vee=\iota$ with $b_\iota=d_\iota=1_\iota$.  
 A  {\it left duality}   is   defined in a  similar way. 
   \end{defn}
  
Every right duality gives rise to a  contravariant functor $D:{\mathcal C}\to{\mathcal C}$ acting as 
\begin{equation}\label{duality_functor}\rho\to\rho^\vee, \quad T\in(\rho,\sigma)\to T^\vee:=d_\sigma\otimes 1_{\rho^\vee}\circ1_{\sigma^\vee}\otimes T\otimes 1_{\rho^\vee}\circ 1_{\sigma^\vee}\otimes b_\rho\in(\sigma^\vee, \rho^\vee),\end{equation}  called the {\it right duality functor}, which turns out tensorial.
A different right duality structure  $(\tilde{\rho}, \tilde{b}_\rho, \tilde{d}_\rho)$ gives a corresponding  duality functor
$\tilde{D}$ related to $D$ via   the isomorphisms $\xi_\rho:\tilde{\rho}\to\rho^\vee$ defined in (\ref{xi}), which is    a natural   monoidal isomorphism  $\xi: \tilde{D}\to D$.

 Right and left dualities naturally arise in  representation categories of Hopf algebras and their generalisations, where canonical choices   are induced by the antipode, we shall discuss this in detail in Sect. \ref{5}.
 A well-behaved choice of right and left dualities lead to the notion of {\it spherical category}. In a spherical category a theory of categorical dimension can be developed. 
 By a theorem of Deligne \cite{Yetter},  see also Sect. \ref{17},
when the category is braided there is a correspondence between spherical structures and ribbon structures for the braided symmetry.

\begin{defn}\label{compatibility_with_duality}  Let ${\mathcal C}$ be a tensor category with right duality $(\rho^\vee, b_\rho, d_\rho)$. A natural isomorphism $\eta\in(1, 1)$ of the identity functor of ${\mathcal C}$ is called {\it compatible with duality} if $$\eta_{\rho^\vee}=(\eta_\rho)^\vee.$$
 \end{defn}
 
 If $\rho$ is simple then $\eta_\rho$ is a nonzero scalar multiple of $1_\rho$. 
 It is easy to  see that the property of being compatible 
with duality for an isomorphism $\eta\in(1,1)$ does not depend on the choice of the right duality.

 We   recall the definition of braided symmetry,   ribbon (or premodular) and modular category category.

\begin{defn}\label{braided_symmetry_tensor_category}
A {\it braided symmetry} for ${\mathcal C}$  is a natural isomorphism $c(\rho, \sigma)\in(\rho\otimes\sigma, \sigma\otimes\rho)$ satisfying
the normalization property
\begin{equation}\label{normalization_symmetry} c(\rho,\iota)=c(\iota,\rho)=1_\rho\end{equation} and such that the following two {\it hexagonal diagrams} commute
\begin{equation}\label{braided_symmetry1} \begin{CD}
 (\rho\otimes \sigma)\otimes \tau @>{{\alpha}}>>  \rho\otimes(\sigma\otimes \tau)@>{c}>> (\sigma\otimes\tau)\otimes \rho \\
@V{c\otimes 1}VV @. @VV{\alpha}V \\
 (\sigma\otimes\rho)\otimes\tau@>{\alpha}>>  \sigma\otimes(\rho\otimes\tau)@>{1\otimes c}>> \sigma\otimes(\tau\otimes\rho) \end{CD}\end{equation}
 \begin{equation}\label{braided_symmetry2} \begin{CD}
 (\rho\otimes \sigma)\otimes \tau @>{{c}}>>  \tau\otimes(\rho\otimes \sigma)@>{\alpha^{-1}}>> (\tau\otimes\rho)\otimes\sigma \\
@A{\alpha^{-1}}AA @. @AA{c\otimes 1}A \\
 \rho\otimes(\sigma\otimes\tau)@>{1\otimes c}>>  \rho\otimes(\tau\otimes\sigma)@>{\alpha^{-1}}>> (\rho\otimes\tau)\otimes\sigma \end{CD}\end{equation}
 
 A weak quasi-tensor functor between braided tensor categories $({\mathcal F}, F, G): ({\mathcal C}, c)\to({\mathcal D}, d)$  is called {\it braided} if the structure maps are compatible with the braided symmetries, meaning
 $${\mathcal F}(c(\rho, \sigma))=F_{\sigma, \rho} d({\mathcal F}(\rho), {\mathcal F}(\sigma)) G_{\rho, \sigma},$$
  $${\mathcal F}(c(\rho, \sigma)^{-1})=F_{\rho, \sigma} d({\mathcal F}(\rho), {\mathcal F}(\sigma))^{-1} G_{\sigma, \rho}.$$
  For a tensor equivalence, the natural transformations $F$, $G$, are one another inverse isomorphisms, consequently only one equation suffices.
 \end{defn}

\begin{prop}\label{braided_symmetry_with_generating_object}
Let ${\mathcal C}$ be a semisimple tensor category with associativity morphisms 
$\alpha$. Let $c(\tau, \sigma)\in(\tau\otimes\sigma, \sigma\otimes\tau)$ be natural isomorphisms satisfying (\ref{normalization_symmetry}) and (\ref{braided_symmetry1}) (or (\ref{braided_symmetry2} resp.) Let $\rho$ be a generating object of ${\mathcal C}$ as
in Def. \ref{Def_generating_object}. Then $c$ is uniquely determined by the morphisms
$c(\rho_\lambda, \rho)$ ($c(\rho, \rho_\lambda)$ resp.) with $\rho_\lambda$ varying in a complete set of simple objects of ${\mathcal C}$.

In particular, if $c$ is a braided symmetry then it is determined by either class of morphisms of the form
$c(\rho_\lambda, \rho)$ or of the form $c(\rho, \rho_\lambda)$.

\end{prop}

\begin{proof} By  (\ref{braided_symmetry1}), if
  $\rho$  is replaced by a fixed tensor power $\rho^r$  of the given generating object,
  $\sigma=\rho$, $\tau=\rho^t$ with $t\geq0$ and $\rho^t$ a tensor power of the generating object with parentheses on the right, it follows by induction that $c(\rho^r, \rho^s)$ for $r, s>0$  is determined by $c(\rho^r, \rho)$. Thus by naturality and semisimplicity, $c$ is determined
  by all the $c(\rho_\lambda, \rho)$. A similar reasoning holds if $c$ satisfies  (\ref{braided_symmetry2}) and shows that $c$ is determined by the collection of all $c(\rho, \rho_\lambda)$. The last statement follows from these two cases.

\end{proof}

\begin{defn}\label{ribbon_category} Let ${\mathcal C}$ be a rigid   tensor category with braided symmetry $c$.
A {\it ribbon structure} is a natural isomorphism $v\in(1,1)$ such that $c(\sigma, \rho)\circ c(\rho, \sigma)=v_\rho\otimes v_\sigma\circ v_{\rho\otimes\sigma}^{-1}$ and compatible with some right duality.
\end{defn} 

In a ribbon category one has a notion of ${\mathbb C}$-valued categorical trace ${\rm Tr}_\rho(T)$ for all morphisms
$T\in(\rho, \rho)$ see Sect. \ref{17} for the definition. For any pair of objects $\rho$, $\sigma$, consider the matrix
$$S_{\rho, \sigma}={\rm Tr}_{\rho\otimes \sigma}(c(\sigma, \rho)c(\rho, \sigma)).$$
The number $S_{\rho, \sigma}$ depends on $X$ and $Y$ only up to isomorphism.
When ${\mathcal C}$ is a fusion ribbon category, let $\{\rho_i\}$ be a complete set of irreducible objects.

\begin{defn} A ribbon fusion category is called {\it modular} if the matrix $S=(S_{\rho_i, \rho_j})$ is  invertible.
\end{defn}

Let $T=(T_{\rho_i, \rho_j})$ be the diagonal matrix defined by $T_{\rho_i, \rho_j}=\delta_{i, j} v_{\rho_i}^{-1}$.
  If ${\mathcal C}$ is modular then 
  \begin{equation}
S\to s=\left(\begin{array}{cc} 0 & 1\\ -1 & 0 \end{array}\right)\quad\quad 
T\to t= \left(\begin{array}{cc} 1 & 1\\ 0 & 1 \end{array}\right)
\end{equation}
defines a projective representation of ${\rm SL}(2, {\mathbb Z})$. Recall that ${\rm SL}(2, {\mathbb Z})$.
is generated by $s$, $t$ with relation $s^2=-1$, $(st)^3=1$.

Let ${\mathbb Z}_2({\mathcal C})$ be {\it symmetric center} of ${\mathcal C}$, defined in \cite{Br2}, \cite{Mueger_Galois}
as the full subcategory of ${\mathcal C}$ with objects the set of $\rho\in{\mathcal C}$ for which $c(\rho, \sigma)c(\sigma, \rho)=1_{\sigma\otimes\rho}$ for all $\sigma\in{\mathcal C}$.  
We shall need the following characterization of modularity, ${\mathcal C}$ is modular if and only if 
  ${\mathbb Z}_2({\mathcal C})$ is trivial.
  We refer to \cite{Mueger3} for more complete information on modular categories, see also \cite{EGNO}.

\bigskip

  Unitary braided symmetries are central notions for this paper see e.g. Sect. \ref{10}, \ref{18}, \ref{20}, \ref{KW}, \ref{VOAnets}, \ref{VOAnets2}.
In Sect. \ref{17}   we shall extend Deligne theorem
to a class of symmetries more general than braided symmetries which play a central role
  in the study of unitary structures in this paper in Sect. \ref{18}, \ref{19}, \ref{20}.   
Ribbon structure and categorical dimension are used in our applications, the classification result of ${\mathfrak sl}_{N, \ell}$-type categories in Sect. \ref{KW}.

Prop. \ref{braided_symmetry_with_generating_object}
 will be  useful  to describe the braiding  of the representation category of
the Zhu algebra associated to an affine  vertex operator algebra at positive integer with respect to a ribbon braided tensor structure that 
we shall construct, see  Theorem \ref{braiding_Zhu_algebra}.  See also Sect. \ref{23},  for a statement of the main result
and Sect. \ref{22} for further comparison of our
 braiding with previously known braided tensor structures in the setting of loop groups.

\begin{prop}\label{weak_duality}
Let $({\mathcal F}, F, G): {\mathcal C}\to{\mathcal C}'$ be a weak tensor functor between tensor categories.  If $\rho^\vee$ is a right dual of $\rho$ defined by $d\in({\rho}^\vee\otimes \rho,\iota)$ and  $b\in(\iota, \rho\otimes {\rho}^\vee)$ 
then
${\mathcal F}(\rho^\vee)$ is a right dual of ${\mathcal F}(\rho)$ defined by
$d_1={\mathcal F}(d)\circ F_{\rho^\vee, \rho}$ and $b_1=G_{\rho, \rho^\vee}\circ {\mathcal F}(b)$,
similarly for  left duals.
\end{prop}

\begin{proof} 
We only show that $d_1$ and $b_1$
solve
 (\ref{first_right})
 for ${\mathcal F}(\rho)$. We have
$$  1_{{\mathcal F}(\rho)}\otimes d_1\circ\alpha'_{{\mathcal F}(\rho),  {\mathcal F}({\rho}^\vee),{\mathcal F}(\rho)}\circ b_1\otimes 1_{{\mathcal F}(\rho)}=$$
$$F_{\rho,\iota}\circ 1_{{\mathcal F}(\rho)}\otimes d_1\circ\alpha'_{{\mathcal F}(\rho),  {\mathcal F}({\rho}^\vee),{\mathcal F}(\rho)}\circ b_1\otimes 1_{{\mathcal F}(\rho)}\circ G_{\iota,\rho}= $$
$$ F_{\rho,\iota}\circ 1_{{\mathcal F}(\rho)}\otimes {\mathcal F}(d)\circ 1_{{\mathcal F}(\rho)}\otimes F_{\rho^\vee, \rho}\circ\alpha'_{{\mathcal F}(\rho),  {\mathcal F}({\rho}^\vee),{\mathcal F}(\rho)}\circ G_{\rho, \rho^\vee}\otimes 1_{{\mathcal F}(\rho)}\circ {\mathcal F}(b)\otimes 1_{{\mathcal F}(\rho)}\circ G_{\iota,\rho}=$$
$$
   {\mathcal F}(1_\rho \otimes d)\circ F_{\rho, \rho^\vee\otimes\rho}    \circ 1_{{\mathcal F}(\rho)}\otimes F_{\rho^\vee, \rho}\circ\alpha'_{{\mathcal F}(\rho),  {\mathcal F}({\rho}^\vee),{\mathcal F}(\rho)}\circ G_{\rho, \rho^\vee}\otimes 1_{{\mathcal F}(\rho)}\circ G_{\rho\otimes\rho^\vee, \rho}\circ{\mathcal F}(b\otimes 1_ \rho)=
$$
$$
{\mathcal F}(1_\rho \otimes d)\circ   {\mathcal F}(\alpha_{\rho, \rho^\vee, \rho})\circ{\mathcal F}(b\otimes 1_ \rho)=1_{{\mathcal F}(\rho)}.
$$
 
\end{proof}

\begin{cor}\label{weak_dim}
Let ${\mathcal C}$ be a rigid tensor category and ${\mathcal F}: {\mathcal C}\to{\rm        Vec}$ be a weak tensor functor.
Then ${\rm        dim}({\mathcal F}(\rho))={\rm        dim}({\mathcal F}(\rho^\vee))={\rm        dim}({\mathcal F}({}^\vee\rho))$ for every object $\rho$.
\end{cor}

If a tensor category is rigid, left and right duals  need not be isomorphic. It is easy to see that this is the case if and only if $\rho\simeq\rho^{\vee\vee}$ and,  
 following M\"uger, we call ${\rho}^\vee$ a {\it two-sided dual} of $\rho$.  We shall say that ${\mathcal C}$ has two-sided duals if every object has a two-sided dual. For example,   
 duals are two-sided if ${\mathcal C}$ is a semisimple tensor category, see e.g. Prop. 2.1 in  \cite{ENO},  
 a tensor category with a coboundary, e.g. a braided symmetry, by Prop. \ref{left_duality_coboundary}, or
  a tensor $^*$-category  \cite{LR}.
  In the last case,   a solution $d$ and $b$ of the right duality equations gives one of the left duality equations via  
  $\rho^{\vee}:={}^{\vee}\rho $,
      $b'=d^*$ and $d'=b^*$. This  dual  is also  called a {\it conjugate} of $\rho$ and denoted $\overline{\rho}$. The duality equations are written in terms of  $r:=d^*$ and $\overline{r}:=b$, and referred to as the {\it conjugate equations}:
\begin{equation}\label{conjugates_def}\overline{r}^*\otimes 1_\rho\circ\alpha^{-1}_{\rho, \overline{\rho},\rho}\circ 1_\rho\otimes r= 1_\rho,\quad\quad
  {r}^*\otimes 1_{\overline{\rho}}\circ \alpha^{-1}_{\overline{\rho}, \rho, \overline{\rho}}\circ 1_{\overline{\rho}}\otimes \overline{r}= 1_{\overline{\rho}}.\end{equation}
Let ${\mathcal C}$ be a tensor $C^*$-category. 
The {\it intrinsic dimension} of $\rho$ is defined as $d(\rho)=\inf\|r\|\overline{r}\|$
over all solutions of the conjugate equations  for $\rho$  \cite{LR}.

\begin{cor}\label{intrinsic_dimension_bound}
Let ${\mathcal C}$ and ${\mathcal C}'$ be tensor $C^*$-categories and ${\mathcal F}: {\mathcal C}\to{\mathcal C}'$ a   weak tensor $^*$-functor defined by $(F, G)$. If $\rho\in{\mathcal C}$ has a conjugate then 
$d({\mathcal F}(\rho))\leq \|F_{\overline{\rho}, \rho}\|\|G_{\rho,\overline{\rho}}\|d(\rho)$.
\end{cor}

\begin{proof}
Let $b$, $d$ solve the right duality equations for $\rho$ and consider the associated solution $b_1$, $d_1$ for ${\mathcal F}(\rho)$ as in Prop.
\ref{weak_duality}, so $r_1=d_1^*$, $\overline{r}_1=b_1$ solves the conjugate equations for the same object. We have   $r_1^*r_1\leq\|F_{\overline{\rho}, \rho}\|^2{\mathcal F}(r^*r)$ 
  so $\|r_1\|\leq \|F_{\overline{\rho}, \rho}\| \|r\|$ by the $C^*$-property. Similarly $\|\overline{r}_1\|\leq \|G_{\rho, \overline{\rho}}\|\|\overline{r}\|$
  and the conclusion follows.\end{proof}
  
In particular if ${\mathcal F}$ is a unitary weak tensor functor we have
$d({\mathcal F}(\rho))\leq \|d(\rho)\|$, and if ${\mathcal F}$ is in turn unitary tensor we recover a well known upper bound
in representation theory of compact quantum groups of the vector space dimension of a representation by the quantum dimension.
More precisely, this case corresponds to   ${\mathcal C}$ the representation category of the compact quantum group,
   ${\mathcal C}'={\rm Hilb}$ and ${\mathcal F}$ the forgetful functor,  see   Cor. 2.2.20 in \cite{CQGRC}.
   
  As already remarked before Def. \ref{unitarity}, we shall see that by Prop. \ref{intrinsic_dimension_bound}
    together with the results  in Sect. \ref{13} and more specifically
    Cor. \ref{unitarity_obstruction},  in ${\mathcal C}$ and ${\mathcal C}'$ are rigid
    $C^*$-tensor categories and ${\mathcal C}$ is amenable then every unitary weak tensor functor ${\mathcal F}: {\mathcal C}\to{\mathcal C}'$ preserves the intrinsic dimensions. In particular,
    non-integrality of the intrinsic dimension    is
    an obstruction to the concurrence of both unitarity and weak tensoriality   for a weak quasitensor structure
$({\mathcal F}, F, G)$ to ${\rm Hilb}$. In the  non-weak case this   result was shown in \cite{LR}, see also Cor. 2.7.9
 in \cite{CQGRC}
and references therein.
   Examples of non-unitary weak tensor structures or unitary weak quasitensor structures arising from fusion categories associated to quantum groups at roots of unity and conformal field theory will be discussed in Sect. \ref{19}, \ref{KW}, \ref{VOAnets2}.

\section{Weak quasi-Hopf algebras}\label{4}
In \cite{Drinfeld_quasi_hopf} Drinfeld introduced the notion of quasi-Hopf algebra as an extension of that of 
  Hopf algebra to 
the case where the coproduct is not coassociative. Quasi-Hopf algebras are more flexible than Hopf algebras in that they admit a so called twist operation.

Quasi-Hopf algebras play an important role in the proof of the Drinfeld-Kohno theorem on the connection between  conformal field theory and quantum groups \cite{Drinfeld_quasi_hopf}, see also \cite{NT_KL}.
However, quasi-Hopf algebras are not sufficiently general to describe fusion categories from CFT.  This follows from Frobenius-Perron theorem, according to which a fusion category 
 ${\mathcal C}$ admits a unique positive   dimension function, it is the Frobenius-Perron dimension function, 
 $\rho\in{\rm Irr}({\mathcal C})\to {\rm FPdim}(\rho)$, see Sect. 5 in  \cite{EGNO}, see also Sect. \ref{13}, \ref{KW}.
This implies that ${\mathcal C}$ is tensor equivalent to ${\rm Rep}(A)$ for a quasi-Hopf algebra $A$ if and only if ${\rm FPdim}$ takes values in ${\mathbb N}$, in this case $A$ is unique up to twist deformation. However the integrality condition is not satisfied already for the fusion category associated the Ising model, which may be realised by an affine vertex operator algebra over ${\mathfrak sl}_2$ at level $2$ \cite{MS}. 

In the early 90s Mack and Schomerus \cite{MS} suggested to give up the request that the coproduct is unital. This leads to the   notion  of weak quasi-Hopf algebra, that is the main subject of this section and
 plays a  central role in this paper.
  As we shall see, 
Drinfeld notion of twist deformation extends in a natural way to weak quasi-Hopf algebras.

 \medskip

\begin{defn}\label{partial_inverse}
Let $B$ be an algebra, and consider the linear category with objects idempotents of $B$ and morphism spaces 
between two
idempotents $p$, $q\in B$ defined by
$$(p, q):=qBp=\{T\in B: qT=T=Tp\}.$$ Given an element $T\in(p, q)$, we shall refer to $D(T):=p$ and $R(T):=q$ as the {\it domain} and {\it range} of $T$.
We shall call $T$  {\it partially invertible} if it is invertible as  a morphism of that category. In other words,     if there is an element $T^{-1}\in (q,p)$  satisfying
\begin{equation}\label{partial_inverse2} T^{-1}T=p,\quad TT^{-1}=q.\end{equation}
 Clearly $T^{-1}$ is unique in $(q, p)$. We shall refer to   $T^{-1}$ as the partial inverse, or simply the inverse of $T$. 
 \end{defn}
 
 In most of our applications, $p$ is given. Assume that   we have $T$ and $T^{-1}$ such that $T^{-1}$ is a partial left inverse of
 $T$ in the sense of the first equation (\ref{partial_inverse2}), then we have a unique range $q=TT^{-1}$ such that $T$ is partially invertible.

\begin{defn}\label{wqh}  A {\it weak quasi-bialgebra} $A$  
 is defined by the following data
 
\noindent {\bf a)}  {\it algebra}: a    complex, associative  algebra $A$ with unit $I$, \smallskip

 \noindent {\bf b)}  {\it coproduct}: a  possibly non-unital homomorphism $\Delta:A\rightarrow A\otimes A$   \smallskip

  \noindent {\bf c)} {\it counit}:  a homomorphism $\varepsilon:A\rightarrow\mathbb{C}$ satisfying
\begin{equation}
(\varepsilon\otimes 1)\circ\Delta=1=(1\otimes\varepsilon)\circ\Delta, \label{eqn:intro1}
\end{equation}

\smallskip

  \noindent {\bf d)} {\it associator}:  a partially invertible element
    $\Phi\in A\otimes A\otimes A$
  with 
   \begin{equation}
D(\Phi)=\Delta\otimes1(\Delta(I)),\quad\quad R(\Phi)=1\otimes\Delta(\Delta(I)), \label{eqn:intro2} 
\end{equation}
   \begin{equation}
 \Phi\Delta\otimes1(\Delta(a))=1\otimes\Delta(\Delta(a))\Phi{\rm        , }\quad\quad\quad a\in A \label{eqn:intro4},
   \end{equation}
    \begin{equation}  
(1\otimes1\otimes\Delta(\Phi))(\Delta\otimes1\otimes1(\Phi))=(I\otimes\Phi)(1\otimes\Delta\otimes1(\Phi))(\Phi\otimes I) \label{eqn:intro6},   \end{equation}
   \begin{equation} 1\otimes\varepsilon\otimes1(\Phi)=\Delta(I) \label{eqn:intro7}.
\end{equation}
 
 \end{defn}

   The relations 
$\varepsilon\otimes 1\otimes1(\Phi)=\Delta(I)=1\otimes1\otimes\varepsilon(\Phi)$ hold automatically, a result extending a known result for quasi-bialgebras algebras.   For example, the first follows  from  the fact   that the domain and range of $\varepsilon\otimes 1\otimes1(\Phi)$ is $\Delta(I)$, and then   as in the quasi-bialgebra case \cite{Drinfeld_quasi_hopf}, one evaluates     $\varepsilon\otimes\varepsilon\otimes1\otimes1$ on (\ref{eqn:intro6}) and takes into account
properties (\ref{eqn:intro1}), (\ref{eqn:intro2}),   (\ref{eqn:intro7}).  
Equation (\ref{eqn:intro6}) is called {\it pentagon equation} in the setting of weak quasi-bialgebras, and plays an important role in connection with construction 
of tensor categories.

We shall find it useful to introduce a weaker terminology. 

\begin{defn}\label{pre-associator}
If we have a quadruple $(A, \Delta, \varepsilon, \Phi)$ such that $(A, \Delta, \varepsilon)$ satisfies axioms a), b), c)
and $\Phi\in A\otimes A\otimes A$ satisfies (\ref{eqn:intro2}),  (\ref{eqn:intro7})  then $\Phi$ is called a {\it pre-associator}.
\end{defn}

We shall see that the notion of pre-associator is useful as in many interesting cases coming from categories from CFT,
where a pre-associator can be introduced using only data of the coproduct.  If we disregard $\Phi$ we shall see that
${\rm Rep}(A)$ is a pre-tensor category. But if one can show that $\Phi$  satisfies the remaining 
axioms  of an associator (partial invertibility, the intertwining relation (\ref{eqn:intro4}) and the pentagon equation
(\ref{eqn:intro6}))
  then $(A, \Delta, \varepsilon, \Phi)$ is  a weak quasi-bialgebra and   ${\rm Rep}(A)$ is a tensor category that by construction depends
  up to tensor equivalence only on $(A, \Delta, \varepsilon)$. We shall develop this in Sect. \ref{5} and we see an application in Sect. \ref{21}, \ref{22}, \ref{23}.

   \begin{defn}\label{antipode}
A  {\it weak quasi-Hopf algebra} is a weak quasi-bialgebra with an {\it antipode}: an antiautomorphism $S$ of $A$ together with elements $\alpha$, $\beta\in A$ for which
\begin{align}
& S(a_{(1)})\alpha a_{(2)}  =\varepsilon(a)\alpha,\quad a_{(1)}\beta S(a_{(2)}) =\varepsilon(a)\beta, \quad\quad\quad  a\in A \label{eqn:antip1}\\
 &  x\beta S(y)\alpha z=I= S(x')\alpha y'\beta S(z'), \label{eqn:antip2}
\end{align}
where  $m: A\otimes A\to A$ is the multiplication map and we use the notation
$\Phi=x\otimes y \otimes z$, $\Phi^{-1}= x'\otimes y'\otimes z' $.

 \end{defn}
 
  If $\Delta$ is unital, the definition of weak quasi-Hopf algebra   reduces to that of  quasi-Hopf algebra introduced by Drinfeld in  \cite{Drinfeld_quasi_hopf}. 
  The following example provides the simplest family of  quasi-Hopf algebras. 

 \begin{ex}\label{pointed_wqh}
 Let $G$ be a finite group. The   algebra $\text{Fun}_\omega(G)$ of complex valued functions on $G$ is a commutative quasi-bialgebra with
   coproduct $\Delta(f)(g, h)=f(gh)$,   counit $\varepsilon(f)=f(e)$, associator given  by a  normalized  $3$-cocycle $\omega: G^{3}\to {\mathbb T}$.
     If $\omega$ is trivial we recover the usual   Hopf algebra $\text{Fun}(G)$.
 If $\omega$ is a $3$-cocycle and  $\omega_F(g, h, k)=F(h, k) F(g, hk)\omega(g, h, k) F^{-1}(gh, k) F^{-1}(g, h)$ is a cohomologous $3$-cocycle via a  normalized $2$-cochain $F$ then  $\text{Fun}_{\omega_F}(G)=(\text{Fun}(G)_\omega)_F$. It follows that the twist isomorphism class of  
 $\text{Fun}_\omega(G)$ is determined by the class of $\omega$ in $H^3(G, {\mathbb T})$.
 An antipode is given by
  $S(f)(g)=f(g^{-1})$, $\alpha(g)=\omega(g, g^{-1}, g)^{-1}$, $\beta(g)=1$. 
  (Note that the $3$-cocycle relation for $\omega$ yields the equality 
   $\omega(g, g^{-1}, g)=\omega(g^{-1}, g, g^{-1})^{-1}$, which is useful to verify the antipode axioms.)
\end{ex}
  
   \begin{defn}
 An antipode $(S, \alpha, \beta)$ will be called {\it strong} if   $\alpha=\beta=I$.  
\end{defn}

 \begin{rem} An antiautomorphism $S$ of $A$ can be a strong antipode only if it satisfies the following compatibility conditions with the associator,
\begin{equation}\label{associator_strong_antipode}
 x S(y) z=I,\quad\quad   S(x') y' S(z')=I.
\end{equation}
For example, when $A$ is a bialgebra, that is $\Phi=I\otimes I\otimes I$, then the above  equations obviously hold 
and the notion of a strong antipode reduces to the usual notion of antipode of a Hopf algebra.
More generally, in the weak case  we shall see that   equations (\ref{associator_strong_antipode}) are satisfied by the associator of a weak Hopf algebra, see Sect. \ref{6}. 
  \end{rem}

  \begin{defn}\label{definition_of_twist}  Let $A$ be a weak quasi-bialgebra with coproduct $\Delta$ and counit $\varepsilon$.
  \begin{itemize}
 \item[{\rm        a)}]
A {\it twist} is a pair of elements       $T, T^{-1}\in A\otimes A$  such that $T^{-1}$ is a partial left inverse of $T$, that is $T^{-1}T=\Delta(I)$ and such that $\varepsilon\otimes 1(T)=1\otimes \varepsilon(T)=I$. 
  \item[{\rm        b)}]  
A {\it trivial twist} of $A$ is a twist of the form $E=P\Delta(I)$ where   $P\in A\otimes A$ is an idempotent,
  $E^{-1}=\Delta(I)P$, $EE^{-1}=P$.
 \end{itemize}

 \end{defn}
 
If $P$ is a trivial twist then    $P=\Delta_P(I)$. In particular, in
the framework of quasi-bialgebras the only trivial twist is the identity, 
and this motivates our terminology. Trivial twists may informally be thought as the necessary adjustment between two
weak bialgebra structures that that would   be coinciding except for the value the coproducts take on the identity.
Trivial twists will arise in the study of unitary structures in Sect. \ref{7} and unitary ribbon structures in Sect. \ref{18},
\ref{19},
\ref{20}.

\begin{prop}\label{twisted_wqh}
A twist $T$ of a weak quasi bialgebra $A$  gives rise to another weak quasi-bialgebra, denoted $A_T$, with the same  algebra structure and counit as $A$   but 
coproduct and associator given by
$$\Delta_T(a)=T\Delta(a)T^{-1}$$
\begin{equation}\label{twisted_associator}\Phi_T= I\otimes T 1\otimes\Delta(T)\Phi\Delta\otimes1(T^{-1})T^{-1}\otimes I.\end{equation}
If $A$ has   antipode $(S, \alpha, \beta)$, then $A_T$ has antipode $(S, \alpha_T, \beta_T)$ where 
 \begin{equation}\label{twisted_antipode}
  \alpha_T= S(f')\alpha g', \quad \beta_T= f\beta S(g), 
\end{equation} and  $T= f\otimes g$,  $T^{-1}=f'\otimes g'$.

\end{prop}

\begin{proof}
Verification of the axioms can be done as in the unital case,   \cite{Kassel}, with slight modifications due to  non triviality of domain idempotents.
\end{proof}

  \medskip

In the last part of the section  we extend to weak quasi-Hopf algebras  properties of   antipodes   of quasi-Hopf algebras \cite{Drinfeld_quasi_hopf}.
 
  \begin{prop}\label{unique_antipode} Let $A$ be a weak quasi-Hopf algebra with antipode $(S,\alpha,\beta)$. Then for every invertible $u\in A$,
    the triple   $(\overline{S},\overline{\alpha},\overline{\beta})$ defined by
       \begin{equation}
\overline{S}(a)=uS(a)u^{-1},\label{intertwining_antipode} 
\end{equation}\begin{equation}
\overline{\alpha}=  u\alpha,\quad \quad\overline{\beta}= \beta u^{-1}\label{eqn:genantprop}
\end{equation}
 is another antipode of $A$. Conversely, any antipode is of this form with $u\in A$ uniquely determined by (\ref{intertwining_antipode}) and one of the equations in (\ref{eqn:genantprop}).    \end{prop}
 
\begin{proof} From (\ref{eqn:intro6}) it follows that
$$
(1\otimes1\otimes\Delta(\Phi))(\Delta\otimes1\otimes1(\Phi))(\Phi^{-1}\otimes I)= $$
$$(I\otimes\Phi)(1\otimes\Delta\otimes1(\Phi))((1\otimes\Delta(\Delta(I)))\otimes I)=
(I\otimes\Phi)(1\otimes\Delta\otimes1(\Phi)).$$
We may extend  the proof of the   quasi-Hopf case, i.e. Prop. 1.1   of \cite{Drinfeld_quasi_hopf},
  to the weak case. 
\end{proof}

  Notice that $u$ and $u^{-1}$ can be derived from (\ref{eqn:genantprop}) if one of the antipodes is strong. 

 \begin{cor}\label{strong_antipode_twist} Let $A$ be a weak quasi-Hopf algebra and $(S, \alpha, \beta)$ an antipode. Then
 \begin{itemize}
\item[{\rm        a)}]  $A$ admits a strong antipode if and only if   $\alpha$ and $\beta$ are invertible and $\beta=\alpha^{-1}$.
 In particular a strong antipode is unique and given by ${\rm ad}(\alpha^{-1})S$.
  \item[{\rm   b)}]  If $A$ admits a strong antipode $S$ then the same holds for a twisted algebra $A_T$ if and only  if  
\begin{equation}\label{strong_antipode_twist_eq}
 m\circ S\otimes1(T^{-1})=(m\circ 1\otimes S(T))^{-1}.
 \end{equation}
  \end{itemize}
\end{cor}

\begin{proof}
The proof follows from (\ref{twisted_antipode}) and Prop. \ref{unique_antipode}.
\end{proof}

By \cite{Drinfeld_quasi_hopf}, p. 1424, when $\Phi=I$, thus $A$ is a bialgebra, and $(S, \alpha, \beta)$ is an antipode
then $\beta=\alpha^{-1}$, thus we may always assume that the antipode is strong. We shall see that this property extends
to any weak bialgebra with an antipode of a weak quasi-Hopf algebra, see Prop. \ref{weak_strong_antipode}.
 We illustrate these notions for the  quasi-Hopf algebras defined in Example \ref{pointed_wqh}.

   \begin{ex} 
  It follows from Cor. \ref{strong_antipode_twist} that $A={\rm Fun}_\omega(G)$ has a strong antipode if and only if 
  $\omega(g, g^{-1}, g)=1$ for all $g$.
   For example,  when $G={\mathbb Z}_N$, 
   each complex  $N$-th root of unity $w$ induces  the $3$-cocycle $\omega_w(a, b, c)=w^{\gamma(a, b)c}$, with 
 $\gamma(a, b)=\lfloor\frac{a+b}{N}\rfloor-\lfloor\frac{a}{N}\rfloor-\lfloor\frac{b}{N}\rfloor$, where $\lfloor\lambda\rfloor$ is the greatest integer not exceeding  $\lambda$.
 Furthermore this association gives an isomorphism of the group  of $N$-th roots of unity with
 $H^3({\mathbb Z}_N, {\mathbb T})$. If $h$ is the natural generator of ${\mathbb Z}_N$, $\omega(h, h^{-1}, h)=w$. It follows that $\text{Fun}_{\omega_w}({\mathbb Z}_N)\in{\mathcal H}$   if and only if  $w=1$.  
 Quite interestingly, elements of ${\rm Fun}_\omega(G)\in{\mathcal H}'$ can be twisted to  elements of ${\mathcal H}$ which are not   Hopf algebras, but this 
  can happen only if a certain obstruction of the associator vanishes. More in detail,  $F$ is a twist such that 
  $({\rm Fun}_\omega(G))_F\in{\mathcal H}$ if and only if $\beta_F=\alpha_F^{-1}$
  which amounts to solve the equation    \begin{equation}\label{obstruction}
  F(g^{-1}, g)\omega(g, g^{-1}, g)=F(g, g^{-1})\end{equation}
  When there are elements $g\in G$ with $g^2=e$ and such that $\omega(g, g, g)\neq 1$ then clearly the equation has no solution. For example, for $G={\mathbb Z}_2$,   $\omega_{-1}(h, h, h)=-1$. 
  Note that this is a general property, $\lambda_g:=\omega(g, g, g)=\pm 1$ when $g^2=e$, and it is not difficult to see 
  that the property that 
  $\lambda$ take the value $-1$ on some involutive element $g$ is the only obstruction  to solve equation (\ref{obstruction}) for a normalized twist $F$. 
     For example the obstruction vanishes   if $G$ has odd order.
     We shall come back to 
 $3$-cocycles on   ${\mathbb Z}_N$ in Sect. \ref{KW}, cf. (\ref{root}).   
   \end{ex}

 Drinfeld  showed that the antipode of a quasi-Hopf algebra  satisfies a twisted anticomultiplicativity property with the coproduct which extends  the usual (i.e. untwisted) anticomultiplicativity in the framework of Hopf algebras. We in turn extend this to weak quasi-Hopf algebras. Since our arguments are a straightforward generalisation of 
 \cite{Drinfeld_quasi_hopf}, we shall only briefly  sketch the needed modifications.  
Set
\begin{equation}\label{gamma_delta}
\gamma=V((\unit\otimes\Phi^{-1})(1\otimes1\otimes\Delta(\Phi))), \quad\quad
\delta=V'((\Delta\otimes1\otimes1(\Phi))(\Phi^{-1}\otimes\unit)) 
\end{equation}
where $V,V':A^{\otimes 4}\rightarrow A^{\otimes 2}$ are defined by $V(a\otimes b\otimes c\otimes d)=S(b)\alpha c\otimes S(a)\alpha d$ and $V'(a\otimes b\otimes c\otimes d)=a\beta S(d)\otimes b\beta S(c)$.

\begin{prop}
\label{anticom} Let $A$ be a weak quasi-Hopf algebra. Then the  new weak quasi-Hopf algebra with same algebra structure and counit   but coproduct $S\otimes S\circ\Delta^{{\rm        op}}\circ S^{-1}$ and associator $S\otimes S\otimes S(\Phi_{321})$
is a twist of $A$ by a unique 
   partially invertible element $f\in A\otimes A$ such that 
    \begin{equation}
\gamma =f\cdot\Delta(\alpha),\quad\quad 
\delta =\Delta(\beta)\cdot f^{-1}.
\end{equation} Explicitly,  $D(f)=\Delta(I)$, $R(f)=S\otimes S\circ\Delta^{{\rm        \rm op}}(I)$,
\begin{equation} f\Delta(S(a))f^{-1}=S\otimes S(\Delta^{{\rm        \rm op}}(a)),  \quad \quad \label{ef}\quad 
S\otimes S\otimes S(\Phi_{321})=\Phi_f.\end{equation}
 We have
$f= S\otimes S(\Delta^{{\rm        \rm op}}(p))\gamma\Delta(q\beta S(r)) $ and $f^{-1}=   \Delta(S(p)\alpha q){\delta}S\otimes S(\Delta^{{\rm        \rm op}}(r))$.
In particular, if the antipode is strong then $f=\gamma$, $f^{-1}=\delta$.

\end{prop}
\begin{proof}
The proof of the first relation in (\ref{ef}) follows from   the following two lemmas, in turn  extending Lemmas 1 and 2 of \cite{Drinfeld_quasi_hopf} to weak quasi-Hopf algebras. More precisely, thanks to Lemma \ref{lemma_1.7} we may apply    lemma \ref{lemma} to $B=A\otimes A$, $p=\Delta(I)$, $q=S\otimes S(\Delta^{{\rm        \rm op}}(I))$, $f=\Delta$, $g=\Delta\circ S$, $\rho=\Delta(\alpha)$, $\sigma=\Delta(\beta)$, $\overline{g}=S\otimes S\circ\Delta^{{\rm        \rm op}}$, $\overline{\rho}=\gamma$, $\overline{\sigma}=\delta$. We omit the proof of the second relation of (\ref{ef}).
\end{proof}

\begin{lemma}\label{lemma_1.7}
We have: \\

\begin{itemize}
 \item[{\rm        a)}]  \begin{equation}
\gamma=V((\Phi\otimes\unit)(\Delta\otimes1\otimes1(\Phi^{-1}))),\quad 
\delta=V'((1\otimes1\otimes\Delta(\Phi^{-1}))(\unit\otimes\Phi)),
\end{equation}
\item[{\rm        b)}] 
for $a\in A$, 
\begin{equation}
(S\otimes S(\Delta^{{\rm        \rm op}}(a_{(1)})))\gamma\Delta(a_{(2)})=\varepsilon(a)\gamma \\
\Delta(a_{(1)})\delta(S\otimes S(\Delta^{{\rm        \rm op}}(a_{(2)})))=\varepsilon(a)\delta
\end{equation} 
\item[{\rm        c)}] 
\begin{equation}
 \Delta(x)\delta(S\otimes S(\Delta^{{\rm        \rm op}}(y)))\gamma\Delta(z)=\Delta(I)=\end{equation}
\begin{equation}\Delta(\unit) 
 (S\otimes S(\Delta^{{\rm        \rm op}}(p)))\gamma\Delta(q)\delta(S\otimes S(\Delta^{{\rm        \rm op}}(r)))
\end{equation}
\end{itemize}
\end{lemma}
\begin{proof}
a)
By the cocycle property \eqref{eqn:intro6} we can write
\[
\gamma=V(\unit\otimes(\Delta\otimes1(\Delta(\unit)))(1\otimes\Delta\otimes1(\Phi))(\Phi\otimes\unit)(\Delta\otimes1\otimes1(\Phi^{-1}))).
\]
\noindent By the defining antipode property \eqref{eqn:antip1} we have, for $T\in A^{\otimes 4}$,
\[
V(a\otimes\Delta(b)\otimes c\cdot T)=\varepsilon(b)V(a\otimes I\otimes I\otimes c\cdot T)=V(1\otimes 1\otimes\varepsilon\otimes 1(a\otimes b\otimes c)_{134}T).
\]
It suffices to choose $a\otimes b\otimes c=I\otimes\Delta(I)\Phi$ and $T=(\Phi\otimes\unit)(\Delta\otimes1\otimes1(\Phi^{-1})).$
 The identity involving $\delta$ can be proved in a similar way. 
 The proof of  b) and c) is a straightforward generalisation  of the case of quasi-Hopf algebras. We refrain from giving details, and we refer the interested reader to \cite{Drinfeld_quasi_hopf}.

\end{proof}

\begin{lemma}\label{lemma}
Let $B$ be a algebra, $p$ an idempotent in $B$, $f:A\rightarrow B$ a homomorphism and $g:A\rightarrow B$ an anti-homomorphism with $f(I)=g(I)=p$, and $\rho,\sigma\in pBp$ such that:
\begin{equation}
\label{eqn:equaz1}
 g(a_{(1)})\rho f(a_{(2)})=\varepsilon(a)\rho, \quad\quad
f(a_{(1)})\sigma g(a_{(2)})=\varepsilon(a)\sigma
\end{equation}
where $a\in A$. Moreover,
\begin{equation}
 f(x)\sigma g(y)\rho f(z)=p \label{eqn:equaz2},\quad\quad 
  g(p)\rho f(q)\sigma g(r)=p \end{equation}
 In addition, we have an idempotent $q\in B$, $\overline{\rho},\overline{\sigma}\in qBp$ and an anti-homomorphism $\overline{g}:A\rightarrow B$ with $\overline{g}(I)=q$ also satisfying \eqref{eqn:equaz1} - \eqref{eqn:equaz2} (in \eqref{eqn:equaz2} $q$ replaces $p$).
Then there exists a unique partially invertible element $F\in B$  with $D(F)=p$, $R(F)=q$, such that
\begin{equation}\label{F_relations}
 F\rho =\overline{\rho},\quad\quad
\overline{\sigma}F=\sigma \end{equation}
\begin{equation}
\overline{g}(a)=Fg(a)F^{-1}.
\end{equation}
We have
\begin{equation}\label{F_explicit}
 F= \overline{g}(p)\overline{\rho}f(q)\sigma g(r), \quad\quad
F^{-1}=\sum_i g(p)\rho f(q){\overline{\sigma}}{\,}\overline{g}(r).
\end{equation}
\end{lemma}
\begin{proof}
We first show uniqueness. Let $F$ be partially invertible with the stated domain and range and satisfying  
(\ref{F_relations}).  
 Inserting the explicit form of $p$ and $q$ given in 
(\ref{eqn:equaz2})   in the equalities $F=Fp$ and $F^{-1}=qF^{-1}$, respectively,  and   taking into account  the mentioned relations (\ref{F_relations}), 
gives  formulas   (\ref{F_explicit}).

We apply the map $W:A^{\otimes 3}\rightarrow B$, $W(b\otimes c\otimes d)=\overline{g}(b)\overline{\rho}f(c)\sigma g(d)$, respectively to $(\Delta\otimes1(\Delta(a)))\Phi^{-1}$ and $\Phi^{-1}(1\otimes\Delta(\Delta(a)))$ and obtain, if $F$ is defined as in (\ref{F_explicit}), $Fg(a)=\overline{g}(a)F$. Similarly, applying the map $X:A^{\otimes 4}\rightarrow B$, $X(b\otimes c\otimes d\otimes e)=\overline{g}(b)\overline{\rho}f(c)\sigma g(d)\rho f(e)$, to the equality:
\[
\begin{split}
(1\otimes1\otimes\Delta(\Phi))(\Delta\otimes1\otimes1(\Phi))(\Phi^{-1}\otimes\unit)= \\
=(\unit\otimes\Phi)(1\otimes\Delta\otimes1(\Phi))(1\otimes\Delta(\Delta(\unit))\otimes\unit)
\end{split}
\]
gives $F\rho=\overline{\rho}$.
  The relations $FF^{-1}=q$, $F^{-1}F=p$ follow again from 
(\ref{F_relations}). \end{proof}

We next show that   a strictly  coassociative  coproduct with trivial associator in the of a weak case, quasi-Hopf algebra is not compatible with non-unitality of the coproduct.

\begin{prop}\label{Hopf_algebra}
Let $A$ be a weak quasi-Hopf algebra with coassociative coproduct and associator $\Phi=\Delta\otimes 1\circ\Delta(I)=\Phi^{-1}$.
Then $A$ is a Hopf algebra.
\end{prop}

\begin{proof}
It is easy to see that $\Phi$ is an associator and that the elements $\alpha$ and $\beta$ defining an antipode are invertible, hence $A$ admits a strong antipode, say $S$. We are left to show that $\Delta(I)=I\otimes I$. The element  $\gamma$  defined by relation (\ref{gamma_delta})
turns out to be $I$ thanks to coassociativity of $\Delta$. Hence $S$ satisfies the untwisted anticomultiplicative relation
$\Delta\circ S=S\otimes S\circ\Delta^{{\rm        op}}$ by the previous proposition. We use the notation $\Delta(x)=x_1\otimes x_2$ and $\Delta(I)=a\otimes b$ and compute
 
$$
\Delta(I)=\Delta(I)\varepsilon(a)b\otimes I=\Delta(\varepsilon(a)I)b\otimes I= \Delta(a_1S(a_2)) b\otimes I= 
$$
$$
\Delta(aS(b_1))b_2\otimes I=a_1S(b_{1,2})b_2\otimes a_{2}S(b_{1,1})= a_{1}S(b_{2,1)})b_{2,2}\otimes a_{2}S(b_{1})= $$
$$ a_{1}\varepsilon(b_{2})\otimes a_{2}S(b_{1})= a_{1}\otimes a_{2}S(b_{1}\varepsilon(b_{2}))= a_{1}\otimes a_{2}S(b)=$$
$$a\otimes b_{1}S(b_{2})= a\otimes\varepsilon(b)I=a\varepsilon(b)\otimes I=I\otimes I.
$$

\end{proof}

In conclusion of the section we introduce a class of most interest in this paper, those for which the underlying algebra  is isomorphic to
a   direct sum of full matrix algebras. Although   we are mostly interested in finite dimensional algebras,
in the following definition we allow infinite dimensionality.  The  {\it direct sum} of full matrix algebras
$$A=\bigoplus_r M_{n_r}({\mathbb C}),$$ is the algebra with elements of the form $(a_r)$ with $a_r\in M_{n_r}$, and only finitely many of them are nonzero. The identity of $M_{n_r}$ is a minimal central projection of $A$ and will be denoted by $e_r$. Similarly, the {\it direct product} $$M(A)=\prod_{r}  M_{n_r}({\mathbb C})$$ is the algebra of elements $(a_r)$ of the same form but with no further restriction on the entries.  There is no distinction between $A$ and $M(A)$ precisely when the index set is finite, which amounts to say that $A$ is unital.

\begin{defn}\label{discrete}  
An algebra $A$ is called {\it discrete} if it is isomorphic to a direct sum of full matrix algebras. A 
{\it discrete weak quasi bialgebra (Hopf algebra)}  is a discrete algebra endowed with coproduct, counit and associator where the axioms of a   weak quasi bialgebra  are modified as follows.
A coproduct $\Delta: A\to M(A\otimes A)$ is     assumed to take values 
in $M(A\otimes A)=\prod_{r,s} M_{n_r}\otimes M_{n_s}$. For fixed integers $r$, $s$, the sum $\sum_j\Delta(e_j)e_r\otimes e_r$ is well defined as only finitely many entries are nonzero. Then
the coproduct $\Delta$ extends   to a map $M(A)\to M(A\otimes A)$ via the formula $\Delta(a)e_r\otimes e_s=\sum_j\Delta(a_j)e_r\otimes e_s$
for  $a=(a_j)$,    and  the extension is a homomorphism. In particular, $\Delta(I)$ is a well defined idempotent of $M(A\otimes A)$. Similarly, $\Delta\otimes 1$ and $1\otimes\Delta$
extend to $M(A\otimes A)$.
The associator $\Phi$, counit $\varepsilon$ (and the antipode $(S, \alpha, \beta)$  in the Hopf case) are  defined as in the unital case,  except that $\Phi$, $\alpha$, $\beta$  may lie in the corresponding multiplier algebras.
\end{defn}

Most of the results of this section hold for discrete weak quasi  bialgebras (Hopf algebras). In Sect. \ref{8}  we shall   introduce involutive and $C^*$-versions.
As we shall see in later sections, such a class is useful to study semisimple tensor categories.  
We also note that   Van Daele developed a theory for the {\it multiplier Hopf algebras},
a class of algebras more general
than the discrete Hopf algebras
\cite{Van_Daele_An_algebraic_framework_for_group_duality}. An analogous generalization
from the theory of weak quasi-Hopf algebras goes beyond the aim of this paper.

\section{Tannaka-Krein duality and 
 integral weak dimension functions (wdf)}\label{5}

The problem of constructing weak quasi-Hopf algebras
from an abstract fusion category was introduced in \cite{MS, Schomerus} and developed in \cite{HO}. Their   motivation was that  the framework of   quasi-Hopf algebras is an important notion for conformal field theory
but too restrictive for many related fusion categories  as they
 may not  admit integral valued dimension functions. Their central idea
consists  in a weakened notion of a dimension function taking integral values whose existence can easily be proven for
all fusion categories and still allows  Tannakian  recostruction theorems. In this section we review and   expand these results far beyond fusion categories.

  In the first part of this section we describe how
weak quasi-Hopf algebras   lead to rigid tensor categories. We then   discuss Tannaka-Krein duality results for semisimple rigid tensor categories.
We shall then see that every fusion category may be described by a weak quasi-Hopf algebra associated to an 
integral   weak dimension function on the Grothendieck ring of the category. Moreover, we shall extend this result
far beyond the class of fusion categories.

Our description
originates from the work in \cite{HO} and will be   fruitful later on, for different purposes. For example, the weak quasi-Hopf algebra representation provided by an
   integral weak dimension function provides a cohomological
insight into the category that will be further investigated in the paper. Moreover,     weak dimension functions  
will play a central role in our study of unitary structures in fusion categories of affine vertex operator algebras.
Furthermore, we shall describe examples of   algebras naturally associated to certain fusion categories for which
 the integral dimensions arising
 from their representations   satisfy the weak dimension   property, see  Sect. \ref{20} and \ref{VOAnets2}.

Let $A$ be a complex unital algebra. By a representation of   $A$ we mean a unital left action of $A$, $\rho:A\to{\mathcal L}(V)$ on a  finite dimensional complex vector space $V$.  It is customary to pass to the language of (left) $A$-modules, 
 dropping  reference to $\rho$. We shall conform to this notation when no confusion arises.
 The representation category ${\rm        Rep}(A)$  is the category with objects representations of $A$ and morphisms  
 between two objects the subspace $(\rho, {\rho'})$ of  ${\mathcal L}(V_\rho, V_{\rho'})$ consisting of all $A$-linear maps. 
 The {\it forgetful functor} is the  functor $${\mathcal F}: {\rm        Rep}(A)\to{\rm        Vec}$$      associating a representation with its vector space, and acting trivially on morphisms.

 If $A$ admits the structure of a weak quasi-bialgebra $(\varepsilon, \Delta, \Phi)$ then the   counit 
 $\varepsilon$   is  a $1$-dimensional representation. We may form the   tensor product representation  $\rho\underline{\otimes}\rho'$ which is the representation acting on   the   subspace  
 $$V_{\rho\underline{\otimes}\rho'}:=\Delta(I)V_\rho\otimes V_{\rho'}$$ 
 of the tensor product vector space $V_\rho\otimes V_{\rho'}$ with left action induced by the coproduct: 
 $$\rho\underline{\otimes}\rho':=\rho\otimes\rho'\circ\Delta.$$ 
 Given two morphisms $S\in(\rho, \sigma)$, $S'\in(\rho', \sigma')$, the tensor product map    $S\otimes S'\in{\mathcal L}(V_\rho\otimes V_{\rho'}, V_\sigma\otimes V_{\sigma'})$ commutes with the action of 
  $\Delta(I)$, thus takes  $V_{\rho\underline{\otimes}\rho'}$ to  $V_{\sigma\underline{\otimes}\sigma'}$. 
  The restriction $S\underline{\otimes} T$  to $V_{\rho\underline{\otimes}\rho'}$ is a morphism in $(\rho\underline{\otimes}\rho', \sigma\underline{\otimes}\sigma')$.
  Given representations $\rho$, $\sigma$, $\tau$, $(\rho\underline{\otimes}\sigma)\underline{\otimes}\tau$ and $\rho\underline{\otimes}(\sigma\underline{\otimes}\tau)$ act respectively on the ranges of $\Delta\otimes 1\circ\Delta(I)$ and $1\otimes\Delta\circ\Delta(I)$.
  The  {\bf restriction of the action of $\Phi$} to the space of $(\rho\underline{\otimes}\sigma)\underline{\otimes}\underline{\tau}$  is an isomorphism    $\alpha_{\rho, \sigma, \tau}: (\rho\underline{\otimes} \sigma)
\underline{\otimes}\tau\to \rho\underline{\otimes}(\sigma\underline{\otimes} \tau)$. 
In this way   ${\rm        Rep}(A)$ becomes a tensor category with unit object the counit of $A$.

 \begin{prop}
 The forgetful functor ${\mathcal F}: {\rm        Rep}(A)\to{\rm        Vec}$ of a weak quasi-bialgebra $A$  is weak quasi-tensor   with $F_{\rho,\sigma}=\Delta(I)$ and $G_{\rho,\sigma}$ the inclusion map.
\end{prop}

  We   give a categorical interpretation of the notion of twist of a weak quasi-Hopf algebra, extending   
  properties known for quasi-Hopf algebras.
 Let $A$ be a unital discrete algebra endowed with two weak quasi-bialgebra structures $(A, \varepsilon, \Delta, \Phi)$, $(A, \varepsilon, \Delta', \Phi')$.  We may correspondingly form two tensor categories ${\rm        Rep}(A)$, ${\rm        Rep}'(A)$ 
 and the functor ${\mathcal E}: {\rm        Rep}(A)\to{\rm        Rep}'(A)$ acting identically on objects and morphisms. This functor fixes the tensor units, it is
  full,  faithful on morphisms and essentially surjective, and hence ${\mathcal E}$  is an equivalence of linear categories.
Furthermore, the two forgetful functors
 ${\mathcal F}:{\rm        Rep}(A)\to{\rm        Vec}$, ${\mathcal F}':{\rm        Rep}'(A)\to{\rm        Vec}$ satisfy the property that  ${\mathcal F}'{\mathcal E}={\mathcal F}$  just as linear  functors. We would like to  make ${\mathcal E}$   into an equivalence of tensor categories.

\begin{prop}\label{class} Let the discrete unital algebra $A$ be endowed with two weak quasi-bialgebra structures
$A=(A, \varepsilon, \Delta, \Phi)$ and $A'=(A, \varepsilon, \Delta', \Phi')$. Then there is a bijective correspondence
between tensor structures on  the identity linear equivalence
${\mathcal E}: {\rm Rep}(A)\to{\rm Rep}'(A)$ and  twists $F\in M(A\otimes A)$ such that 
$A'=A_F$ as weak quasi-bialgebras. Given $F$, the tensor structure $E_{\rho,\sigma}: {\mathcal E}_{\rho}\otimes {\mathcal E}_\sigma\to {\mathcal E}_{\rho\otimes \sigma}$  is given by the action of $F^{-1}$.
   \end{prop}

\begin{proof} The proof is a straightforward extension of the case of quasi-bialgebras, for which we refer the reader to Prop. 2.1 in \cite{NT_KL}. We briefly comment on how to construct the twist from the tensor structure.
Given a  tensor structure $E_{\rho, \sigma}$ on  ${\mathcal E}:{\rm        Rep}(A)\to:{\rm        Rep}'(A)$ we consider the unique   elements $F^{-1}, F\in M(A\otimes A)$ having components $E_{\rho,\sigma}$, and ${E_{\rho,\sigma}}^{-1}$ respectively in the representation $\rho\otimes\sigma$ of $A\otimes A$.
Then $\rho\otimes\sigma(F^{-1}F)=E_{\rho,\sigma}\circ E^{-1}_{\rho, \sigma}=1_{F(\rho\otimes\sigma)}=\rho\otimes\sigma(\Delta(I))$, hence $F^{-1}F=\Delta(I)$. 
The relation $\varepsilon\otimes 1(F)=I=1\otimes\varepsilon(F)$ can be checked in a similar way, hence $F$ is a twist. The relations $\Delta'=\Delta_F$ and $\Phi'=\Phi_F$ correspond respectively to the   intertwining relations   $E_{\rho, \sigma}\in ({\mathcal E}_\rho\otimes {\mathcal E}_\sigma, {\mathcal E}_{\rho\otimes\sigma})$ and tensoriality property. \end{proof}

Extending the terminology of \cite{CQGRC} to non-coassociative Hopf algebras, a twist $V\in A\otimes A$ is called
{\it invariant} if $\Delta_V=\Delta$ and $\Phi_V=\Phi$. 
For example, if $v\in A$ is central invertible then $\Delta(v) v^{-1}\otimes v^{-1}$ is an invariant twist.
By the previous proposition, invariant twists induce   tensor autoequivalence structures on the identity functor ${\rm        Rep}(A)\to{\rm        Rep}(A)$ and they are all of this form in the discrete case.
  
More generally,  if $A$ is discrete, given $(A,  \varepsilon, \Delta, \Phi)$ and $(A,  \varepsilon, \Delta', \Phi')$, the  weak quasi-tensor structures
on  $E:{\rm        Rep}(A)\to{\rm        Rep}'(A)$ correspond to  the twists $F\in A\otimes A$ such that $\Delta'=\Delta_F$.
Given such a structure, the composite functor ${\mathcal F}'{\mathcal E}$ becomes a weak quasi-tensor with the composed structure.  Since ${\mathcal F}={\mathcal F}'{\mathcal E}$ as functors, this also induces  a new  weak quasi-tensor structure on ${\mathcal F}$. Of course, this is given by the action of $F^{-1}$, with $F$ the twist corresponding to ${\mathcal E}$, so the induced structure on ${\mathcal F}$ determines that of ${\mathcal E}$.   
Thus the construction of  a tensor structure on ${\mathcal E}$   can be regarded as that of a   weak quasi-tensor structure
of the forgetful functor ${\mathcal F}:{\rm        Rep}(A)\to{\rm        Vec}$ defined by a twist $F\in A\otimes A$ solving  $(A, \varepsilon, \Delta', \Phi')=A_F$.

 Two weak quasi-tensor structures on ${\mathcal F}$ are monoidally isomorphic if and only if the corresponding twists  $F_1$ and $F_2$ are related by an invertible $u\in M(A)$ such that $F_2=u\otimes u F_1\Delta(u^{-1})$. This corresponds to
 say that $A_{F_1}$ and $A_{F_2}$ are isomorphic as weak quasi-bialgebras.
For example, the weak quasi-tensor structures on ${\mathcal F}$ monoidally isomorphic to the original one
correspond to twists of the form $u\otimes u\Delta(u^{-1})$, where
 $u\in A$ is an invertible element. These twists are called $2$-coboundaries.
 The monoidal isomorphism $\eta_\rho$ acts as  $\rho(u)$ on $V_\rho$.
\medskip

  Rigidity in  ${\rm        Rep}(A)$   is described similarly to  quasi-Hopf algebras.
  
\begin{defn}\label{duals_wqh} Let $\rho$ be a representation of a weak quasi-Hopf algebra. The contragredient representations $\rho^c$
and ${}^c\rho$ are the 
 the representations of $A$   acting on the dual space $V_\rho'$ respectively as   
  $$\langle\rho^c(a)\phi, \xi\rangle=\langle\phi, \rho(S(a))\xi\rangle,\quad\quad \langle{}^c\rho(a)\phi, \xi\rangle=\langle\phi, \rho(S^{-1}(a))\xi\rangle.$$\end{defn}

  \begin{prop}\label{rigidity}   If $A$ is a weak quasi-Hopf algebra the category ${\rm        Rep}(A)$ is rigid.  Right and left duals of an object $\rho$ are respectively given  by   
  $$\rho^\vee=\rho^c,\quad\quad {}^\vee\rho={}^c\rho.$$
    Solutions of 	the right and left duality equations are respectively given by 
    $$d_\rho(\phi\otimes\xi)=\phi(\alpha\xi)\quad\quad b_\rho=\sum_i \beta e_i\otimes e^i,$$ and 
$$b'_\rho=\sum_i e^i\otimes S^{-1}(\beta) e_i,\quad\quad d'_\rho(\xi\otimes\phi)=\phi(S^{-1}(\alpha)\xi)$$ where $(e_i)$ and $(e^i)$ is a dual pair of bases.
  
  \end{prop}
  
 Thus  
  ${\rm        Rep}(A)$ is rigid and by the above proposition, an antipode of $A$ induces    right and left duality structures, $(b_\rho, d_\rho)$ and $(b'_\rho, d'_\rho)$, respectively, and consequently a
  (say, right) duality functor $c: \rho\to\rho^c$ acting as transposition of $\alpha T\beta$ 
    on a morphism $T$. 
  By Prop. \ref{anticom}    the collection of operators 
$f_{\sigma, \rho}:=\Sigma\sigma^c\otimes\rho^c(S^{-1}\otimes S^{-1}(f_{21}))$ is an invertible  natural transformation $\sigma^c\otimes\rho^c\to(\rho\otimes\sigma)^c$ making $c$ into a    contravariant  tensor functor.
 We compute the   natural transformation associated to $c^2$. We canonically identify the double dual space $V_{\rho}^{''}$ of a representation with $V_\rho$, so $\rho^{cc}$ identifies with $\rho\circ S^2$. 
Reading   (\ref{ef}) as an intertwining relation $f: \Delta\to S\otimes S\circ\Delta^{{\rm        op}}\circ S^{-1}$, it  implies that $S\otimes S(f_{21}^{-1}): S\otimes S\circ\Delta^{{\rm        op}}\circ S^{-1}\to S^2\otimes S^2\circ\Delta\circ S^{-2}$, hence we can form the composite which intertwines
$$S\otimes S(f_{21}^{-1}) f: \Delta\to  S^2\otimes S^2\circ\Delta\circ S^{-2}.$$
 This implies that 
$\rho\otimes\sigma(f^{-1}S\otimes S(f_{21}))$ can be regarded as an intertwiner  $\rho^{cc}\otimes\sigma^{cc}\to(\rho\otimes\sigma)^{cc}$, and this is
  the natural transformation    of   $c^2$.
  
 Note that  left and right duals of the same object of   ${\rm        Rep}(A)$   are equivalent whenever $S^2$ is an inner
  automorphism of $A$ and a converse holds if $A$ is discrete, that is $S^2$ is induced by an invertible in $M(A)$.
  For example  if $A$ is not assumed discrete, $S^2$ is inner whenever
   $A$ has an $\Omega$-involution in the sense of the Sect. \ref{8}
  commuting with    $S$, by Cor. \ref{Kac}, or for the   class  weak  Hopf algebras introduced in Sect. \ref{6} with a quasitriangular structure, by Prop.
  \ref{inner_antipode}.

  If $S^2$ is inner,
 any invertible $x\in A$ such that $S^2(a)=xax^{-1}$ induces an invertible natural transformation $ \eta: 1\to c^2$, where $\eta_\rho$ is defined 
 by the action of $\rho(x)$, but  to construct a pivotal structure we need a monoidal natural transformation.
 
 \begin{defn}\label{pivotal_wqh}
 A {\it pivotal} weak quasi-Hopf algebra is a pair $(A, \omega)$ with $A$ a weak quasi-Hopf algebra and $\omega\in A$  an invertible
element, called the {\it pivot}, such that  $S^2(a)=\omega a\omega^{-1}$ for all $a\in A$ and 
 $f^{-1}S\otimes S(f_{21})=\Delta(\omega)\omega^{-1}\otimes \omega^{-1}$. 
 \end{defn}
 
 The pivot is not unique but determined up to multiplication by an invertible central   element $z$ satisfying $\Delta(z)=z\otimes z$.
 In Sect. \ref{7} we shall  see that if $A$ is a ribbon weak quasi-Hopf algebra, then there is a canonically associated $\omega$ such that $\eta$ becomes a monoidal.   Note that since the identity functor is tensorial, we may use this property to derive tensoriality of $c^2$ more easily for such class of algebras. Indeed, a quasi-tensor functor which is monoidally  isomorphic to a tensor functor must be tensorial as well. This endows  ${\rm  Rep}(A)$ with the structure of a pivotal tensor category. But  more is true:  ${\rm        Rep}(A)$ becomes a spherical category  in the sense of \cite{BW}, 
 a result extending to the weak case, results known
 for ribbon Hopf algebras. Thus, there is a well-behaved theory of dimension in ${\rm        Rep}(A)$,  see Sect. \ref{13}.

 The following Tannakian reconstruction results are due to  \cite{HO} and
extend to the weak case an earlier result of Majid for discrete quasi-Hopf algebras \cite{Majid4}.  For 
a review for discrete Hopf algebras, see \cite{MRT}.
The starting point is an abstract   semisimple category equipped with a fibre functor ${\mathcal F}:{\mathcal C}\to{\rm        Vec}$.
We let ${\rm  Nat}_0({\mathcal F})$ denote the discrete algebra of natural transformations of ${\mathcal F}$ to itself
with finite support.

 \begin{thm}\label{TK_algebraic_quasi} Let ${\mathcal C}$ be a semisimple   category and
   ${\mathcal F}:{\mathcal C}\to{\rm        Vec}$   a faithful   functor. Then   
    \begin{itemize}
\item[(a)] $A={\rm        Nat}_0({\mathcal F})$   is  a discrete     algebra  and there is a linear equivalence ${\mathcal E}:{\mathcal C}\to{\rm        Rep}(A)$   which, after composition with the forgetful functor   ${\mathcal F}_A: {\rm        Rep}(A)\to{\rm        Vec}$, is   isomorphic to ${\mathcal F}$.
 Up to isomorphism, $A$ is   determined by the last      property among discrete algebras. 
 
 \item[(b)]  If  ${\mathcal C}$ is tensorial and ${\mathcal F}$ has a weak quasi-tensor structure then $A$ is a weak quasi-bialgebra, ${\mathcal E}$ is a tensor equivalence,  the isomorphism ${\mathcal F}_A{\mathcal E}\simeq{\mathcal F}$ is monoidal and $A$ is determined among discrete weak quasi-bialgebras.
  \end{itemize}
 
Let $({\mathcal C}$, ${\mathcal F})$ satisfy the same assumptions as in {\rm (b)}.

 \begin{itemize}

 \item[(c)]    If  ${\mathcal C}$ is braided  then $A$ is a quasitriangular
 weak quasi-bialgebra    and ${\mathcal E}$ is braided.
 
  \item[(d)] If
${\mathcal C}$ is rigid and  ${\rm dim}({\mathcal F}(\rho))={\rm dim}({\mathcal F}(\rho^\vee))$
 then a solution of the right duality equations induces an antipode on $A$   making it into a weak quasi-Hopf algebra.
 
 \item[(e)] If
${\mathcal C}$ satisfies (d) and is ribbon then $A$ is a ribbon weak quasi-Hopf algebra.
\end{itemize}
\end{thm}

\begin{proof}   We briefly discuss
  a few  aspects that we shall need.
(a) A natural transformation $\eta\in{\rm        Nat}_0({\mathcal F})=A$ is determined by the values it takes on a complete set of simple objects $\{\rho_i\}_i$, and this gives an algebra isomorphism  of $A\simeq \bigoplus_i{\mathcal L}(V_i)$, with $V_i={\mathcal F}(\rho_i)$, so $A$ is discrete. (b)
As before, $\alpha_{\rho,\sigma,\tau}:(\rho\otimes\sigma)\otimes\tau\to\rho\otimes(\sigma\otimes\tau)$ denote the associativity morphisms of ${\mathcal C}$ and $F_{\rho, \sigma}$ and 
$G_{\rho, \sigma}$   the natural transformations defining the quasi-tensor structure of ${\mathcal F}$.
Counit, coproduct,  and associator    of $A$ are respectively defined as follows. We identify $A\otimes A$ with natural transformations on two variables $\zeta_{\rho,\sigma}: {\mathcal F}(\rho)\otimes {\mathcal F}(\sigma)\to
{\mathcal F}(\rho)\otimes {\mathcal F}(\sigma)$, and similarly for $A^{\otimes 3}$. We set: $\varepsilon(\eta)=\eta_\iota,$
\begin{equation}\label{coproduct}
\Delta(\eta)_{\rho, \sigma}=G_{\rho,\sigma}\circ \eta_{\rho\otimes\sigma}\circ F_{\rho,\sigma},\end{equation}
\begin{equation}\label{associativity}\Phi_{\rho,\sigma, \tau}=1_{{\mathcal F}(\rho)}\otimes G_{\sigma, \tau}\circ G_{\rho, \sigma\otimes\tau}\circ{\mathcal F}(\alpha_{\rho, \sigma,\tau})\circ F_{\rho\otimes\sigma,\tau}\circ F_{\rho, \sigma}\otimes 1_{{\mathcal F}(\tau)}.
\end{equation}
It follows that
$$\Phi^{-1}_{\rho, \sigma, \tau}=G_{\rho, \sigma}\otimes 1_{{\mathcal F}(\tau)}\circ G_{\rho\otimes\sigma, \tau}\circ {\mathcal F}(\alpha^{-1}_{\rho,\sigma,\tau})\circ F_{\rho, \sigma\otimes\tau}\circ 1_{{\mathcal F}(\rho)}\otimes F_{\sigma,\tau}.$$
The axioms can be checked with routine computations. In comparison with the quasi-tensor setting where the natural transformations are invertible, the relations $F_{\rho, \sigma}\circ G_{\rho, \sigma}=1_{{\mathcal F}(\rho\otimes\sigma)}$ is used here to show partial invertibility of $\Phi$.
The tensor equivalence ${\mathcal E}$ is   ${\mathcal F}$ regarded as a functor with values in ${\rm        Rep}(A)$ and tensor structure obtained  by restricting that of ${\mathcal F}$. (c), (e) The notion of braided or ribbon tensor category is recalled in Sect. \ref{3}, Definitions \ref{braided_symmetry_tensor_category} and \ref{balancing_ribbon_tensor_category} respectively.
Quasitriangular   and ribbon structures for weak quasi-bialgebras are given in Sect. \ref{7}, Definition \ref{quasi_triangular_structure} and \ref{balanced_ribbon_wqh}.
 If $c(\rho, \sigma)$ is a braided symmetry in ${\mathcal C}$, and $\Sigma(V, W)$ is the permutation symmetry of ${\rm        Vec}$,  then the element $R\in M(A\otimes A)$ defined by
$\Sigma({\mathcal F}(\rho), {\mathcal F}(\sigma))\circ R_{\rho,\sigma} =G_{\sigma, \rho}\circ{\mathcal F}(c(\rho, \sigma))\circ F_{\rho, \sigma}$ makes $A$ quasitriangular.  When ${\mathcal C}$ has a ribbon structure $v_\rho$
then $A$ has a ribbon structure defined by the ribbon element $v\in M(A)$, where $v$ is the natural transformation ${\mathcal F}(v_\rho)$.
(d) A weak quasi-Hopf algebra antipode
$(S, \alpha, \beta)$
is constructed as follows. For $\rho\in{\rm Irr}({\mathcal C})$, we fix  linear isomorphisms from the dual vector spaces $U_\rho:  {\mathcal F}(\rho)'\to {\mathcal F}(\rho^\vee)$, and extend $U$ to a natural transformation from the functor $\rho\to {\mathcal F}(\rho)'$ to
the functor $\rho\to {\mathcal F}(\rho^\vee)$. We set $S(\eta)_\rho=U_\rho^t\eta_{\rho^\vee}^t{U_\rho^t}^{-1}$, where $L^t:
W'\to V'$ is the transposed of the linear map $L: V\to W$, and $\alpha$, $\beta$ are determined by ${\mathcal F}(d_\rho)\circ F_{\rho^\vee, \rho}\circ U_\rho\otimes 1(f\otimes\xi)=f(\alpha_\rho\xi)$,
$1\otimes U_\rho^{-1}\circ G_{\rho, \rho^\vee}\circ{\mathcal F}(b_\rho)=\sum_i \beta_\rho e_i\otimes e^i$, for 
$\rho\in{\rm Irr}({\mathcal C})$, 
$f\in {\mathcal F}(\rho)'$, 
$\xi\in{\mathcal F}(\rho)$, $e_i\in{\mathcal F}(\rho)$ a linear basis and $e^i\in{\mathcal F}(\rho)'$ the dual basis.
We refer to Lemma 12 in \cite{HO} or to Prop. 2.5 in \cite{NT_KL} for the verification  of the antipode axioms.

The equivalence ${\mathcal E}$ takes an object $\rho\in{\mathcal C}$ to the representation ${\mathcal E}(\rho):\eta\to\eta_\rho$
of ${\rm Nat}_0({\mathcal F})$ on the vector space ${\mathcal F}(\rho)$. The natural transformation $E_{\rho, \sigma}$ making ${\mathcal E}$ into a tensor equivalence is the restriction of $F_{\rho, \sigma}$ to $\Delta(I)F(\rho)\otimes F(\sigma)\to{\mathcal F}(\rho\otimes\sigma)$, with inverse $E_{\rho, \sigma}^{-1}$ given by $G$ considered
as a map ${\mathcal F}(\rho\otimes\sigma)\to \Delta(I)F(\rho)\otimes F(\sigma)$.

\end{proof}

  \begin{rem}\label{faithfulness_and_antipode_reconstruction} a) By semisimplicity of ${\mathcal C}$, faithfulness of ${\mathcal F}$ is equivalent to requiring that ${\mathcal F}(\rho)\neq 0$ for all   simple objects $\rho$.  In particular, ${\mathcal F}$ is always faithful on the morphism spaces $(\rho, \sigma)$ where both $\rho$ and  $\sigma$ are $\neq0$.  
  b) The   requirement of dimension equality in (d) is automatic in the case where ${\mathcal C}$   has finitely many inequivalent simple objects,  (i.e. is a fusion category), see \cite{NT_KL} for a discussion and references, and also where ${\mathcal F}$ is a weak tensor functor, by Cor. \ref{weak_dim}.  c) When we start with a given semisimple weak quasi-Hopf algebra $A$ then  Tannakian reconstruction of Theorem \ref{TK_algebraic_quasi} 
  applied to 
   the forgetful functor ${\mathcal F}:{\rm Rep}(A)\to{\rm Vec}$ with the natural weak quasi-tensor structure
   provides a discrete weak quasi-Hopf structure on ${\rm Nat}_0({\mathcal F})$ which corresponds to 
   the original   structure of $A$
    under the   natural inclusion of $A$ with ${\rm Nat}_0({\mathcal F})$.
   Note that the construction of an antipode of ${\rm Nat}_0({\mathcal F})$ as in the proof of Theorem 
   \ref{TK_algebraic_quasi}
   depends on the choice of
     a right duality $(\rho^\vee, b_\rho, d_\rho)$ of ${\rm Rep}(A)$ and the natural transformation $U$. 
     In particular, by Prop. \ref{rigidity} 
a  given antipode $(S, \alpha, \beta)$  of $A$    
  corresponds to the antipode of  ${\rm Nat}_0({\mathcal F})$
  defined   by the canonical right duality associated to $\rho^\vee=\rho^c$ as in Prop. \ref{rigidity} and to the
  identity natural transformation $U$
   (note that this is an admissible choice as   the functor $\rho\to {\mathcal F}(\rho^\vee)$ coincides
  with $\rho\to{\mathcal F}(\rho)'$).
  d) In general, the algebras ${\rm Nat}_0({\mathcal F})$ and ${\rm Nat}({\mathcal F})$ of general natural transformations
  of ${\mathcal F}$ to itself may have different representation categories, see \cite{Guido_Tuset}.
  However,   regarding ${\rm Nat}({\mathcal F})=M({\rm Nat}_0({\mathcal F}))$ as a topological algebra
  with the strict topology defined by ${\rm Nat}_0({\mathcal F})$ the     category of  nondegenerate representations 
  of ${\rm Nat}_0({\mathcal F})$ coincides  with the
  full subcategory of strictly continuous 
  representations of ${\rm Nat}({\mathcal F})$. We shall touch on again the relevance of   
 the   Tannakian algebra ${\rm Nat}({\mathcal F})$  as a topological algebra for  the forgetful functor associated
  to   $U_q({\mathfrak g})$   for the construction of the $R$-matrix, see  Sects.
  \ref{18}, \ref{19}, \ref{20}.
   \end{rem}

We next introduce the notion
positive weak dimension function.

\begin{defn}\label{weak_dimension_function} Let ${\mathcal C}$ be a  semisimple tensor category.
  A {\it positive  weak dimension function} is a  positively valued function $D $    defined  on a complete set 
  ${\rm Irr}({\mathcal C})$ of irreducible objects   and
satisfying  $D(\iota)=1$, 
and \begin{equation}\label{weak_dim_funct_inequality} \sum_{\tau\in{\rm        Irr}({\mathcal C})} D(\tau){\rm        dim}(\tau, \rho\otimes\sigma)\leq   D(\rho)D(\sigma).    
\end{equation} 
When ${\mathcal C}$ is rigid a   weak dimension function satisfying  $D(\rho)=D(\rho^\vee)=D({}^\vee\rho)$, for all  
$\rho$,  is called {\it symmetric}. 

\end{defn}

If the inequality is always an equality we recover the notion of positive   dimension function.
 We     tacitly extend a  weak dimension function to all the objects of ${\mathcal C}$ via additivity and isomorphism invariance, and  
  (\ref{weak_dim_funct_inequality}) reads as $$D(\rho\otimes\sigma)\leq D(\rho)D(\sigma)$$ for every pair of objects 
  $\rho$ and $\sigma$.
 A weak dimension function $D$ for ${\mathcal C}$ may be   regarded  as  a (unital, additive, and submultiplicative) function on the Grothendieck ring ${\rm Gr}({\mathcal C})$,  and ${\rm Irr}({\mathcal C})$ as a ${\mathbb Z}$-basis.

   For a large part of this paper, we shall    consider     weak dimension functions taking positive integral values.
Furthermore, when the categories have duals, we shall also assume the symmetry condition.
However, in   Sect. \ref{13} and \ref{KW}
 we shall also consider   dimension functions  for a different purpose,    which may not be positive   or integral, but the context should lead to no confusion.

If $A$ is a weak quasi bialgebra and ${\mathcal F}:{\rm Rep}(A)\to{\rm Vec}$  is   the forgetful functor of $A$ then
$D(\rho)={\rm dim}({\mathcal F}(\rho))$ is an integral weak dimension function. 
It follows that a semisimple (rigid) tensor category ${\mathcal C}$  equivalent to the representation category of a weak quasi-bialgebra (quasi-Hopf algebra) admits an integral (symmetric) weak dimension function.
 The following result shows that under suitable conditions existence of an integral weak dimension on ${\mathcal C}$ function is also a sufficient
to represent ${\mathcal C}$ in this way.

\begin{thm}\label{propweakdim}
Let ${\mathcal C}$ be a  semisimple linear category. 
 \begin{itemize}
\item[(a)]
The assignment ${\mathcal F}\to D$, $D(\rho):={\rm dim}(F(\rho))$,
is a bijective correspondence between faithful functors 
${\mathcal F}:{\mathcal C} \to{\rm        Vec}$ up to natural isomorphism and    functions $D:{\rm Irr}({\mathcal C})\to{\mathbb N}$.
\item[(b)]  If ${\mathcal C}$ is   tensorial then the functor   ${\mathcal F}$ admits  a weak quasi-tensor structure  if and only if  $D$ is an integral  weak dimension function. Furthermore, quasi-tensor structures correspond to genuine dimension functions.
\item[(c)]
 The   weak quasi bialgebra   structures on $A={\rm Nat}_0({\mathcal F})$ associated to the various  weak quasi-tensor structures on ${\mathcal F}$ of dimension $D$ as in Theorem \ref{TK_algebraic_quasi} are pairwise   twist isomorphic.
 \end{itemize}

\end{thm}

\begin{proof} (a) Obviously naturally isomorphic functors are associated to the same   function $D:{\rm Irr}({\mathcal C})\to{\mathbb N}$. Conversely, given   $D$, choosing, for $\rho\in{\rm        Irr}({\mathcal C})$,
 a vector space ${\mathcal F}(\rho)$ with ${\rm dim}({\mathcal F}(\rho))=D(\rho)$ gives rise to
 a faithful functor ${\mathcal F}: {\mathcal C}\to{\rm Vec}$, determined 
  up to natural isomorphism.  (b)
If ${\mathcal F}:{\mathcal C}\to {\rm        Vec}$ admits a (weak) quasi-tensor structure then $D(\rho):={\rm        dim}(F(\rho))$ is a (weak) dimension function. 
For the converse, since
by assumption, $ {\rm        dim}({\mathcal F}(\rho)\otimes {\mathcal F}(\sigma))\geq {\rm        dim}({\mathcal F}(\rho\otimes\sigma))$
for all $\rho$, $\sigma\in{\rm        Irr}({\mathcal C})$, we may pick   epimorphisms $F_{\rho, \sigma}: {\mathcal F}(\rho)\otimes{\mathcal F}(\sigma)\to{\mathcal F}(\rho\otimes\sigma)$ and
monomorphisms $G_{\rho,\sigma}: {\mathcal F}(\rho\otimes\sigma)\to{\mathcal F}(\rho)\otimes{\mathcal F}(\sigma)$ 
satisfying $F_{\rho, \sigma}\circ G_{\rho,\sigma}=1$ and acting identically if either $\rho$ or $\sigma$ is the tensor unit.
We  extend these maps to all the objects $\mu$, $\nu$ using complete reducibility:   choose
 $\alpha_\rho^i\in(\rho, \mu)$, $\beta_\rho^i\in(\mu,\rho)$ with $\beta_\rho^j\alpha_\rho^i=\delta_{i, j}1_\rho$, $\sum_{i, \rho}\alpha_\rho^i\beta_\rho^i=1_\mu$, and similarly for   $\gamma_\sigma^j\in(\sigma, \nu)$, $\delta_\sigma^i\in(\nu,\sigma)$.
 Set $F_{\mu, \nu}=\sum{\mathcal F}(\alpha_\rho^i\otimes\gamma_\sigma^j)\circ F_{\rho,\sigma}\circ{\mathcal F}(\beta_\rho^i)\otimes{\mathcal F}(\delta_\sigma^j)$. It is easy to see that naturality holds, that is $F_{\mu', \nu'}\circ{\mathcal F}(S)\otimes{\mathcal F}(T)={\mathcal F}(S\otimes T)\circ F_{\mu, \nu}$. Naturality also shows that $F_{\mu, \nu}$ is independent of the choice of the morphisms involved in the decompositions.
 We similarly obtain a natural transformation 
  $G_{\mu, \nu}$ and it is easy to see that $F_{\mu, \nu}\circ G_{\mu, \nu}=1$. 
We thus have a weak quasi-tensor structure, which is  quasi-tensor   if $D$ is a dimension function. (c)
If $(F, G)$, $(F', G')$ define two weak quasi-tensor structures on ${\mathcal F}$ then we know from 
Theorem \ref{TK_algebraic_quasi} and its proof that the   coproduct associated to the latter is   defined by 
$\Delta'(\eta)_{\rho, \sigma}=G'_{\rho,\sigma}\circ \eta_{\rho\otimes\sigma}\circ F'_{\rho,\sigma}$, and similarly for $\Delta$. We may then write $\Delta'(\eta)=G'F\Delta(\eta) GF'$ since $FG=1$. Setting $T=G'F$ and $T^{-1}=GF'$ we see that these natural transformations may be regarded as   elements of $A\otimes A$ and that $T^{-1}T=GF=\Delta(I)$, $TT^{-1}=G'F'=\Delta'(I)$.  A similar computation shows that the corresponding associators are related by the corresponding twist relation.
  \end{proof}

It follows from Remark \ref{rational} and Theorem \ref{propweakdim} that any   finite semisimple (fusion) category is tensor equivalent to that of a weak quasi   bialgebra (Hopf algebra), and a tensor equivalence 
corresponds to a twist isomorphism between two associated such algebras.
     
 \begin{cor} Let ${\mathcal C}$ and ${\mathcal C}'$ be   semisimple  tensor categories endowed with 
 integral weak dimension functions $D$ and $D'$ respectively compatible with a linear equivalence 
${\mathcal E}: {\mathcal C}\to{\mathcal C}'$. Then ${\mathcal E}$ admits the structure of a tensor equivalence if and only if the corresponding weak quasi-bialgebras  
are   isomorphic up to twist.

\end{cor}

\begin{proof}
If the categories are tensor equivalent then   we apply Th. \ref{propweakdim} and Th. \ref{TK_algebraic_quasi}.  
Conversely,  let
${\mathcal F}:{\mathcal C}\to {\rm Vec} $ and ${\mathcal F}':{\mathcal C}'\to  {\rm Vec}$ be weak quasi-tensor functors of dimensions
 $D$ and $D'$ and associated  weak quasi-bialgebras $A$ and $A'$ 
respectively. Then ${\mathcal F}'{\mathcal E}$ and ${\mathcal F}: {\mathcal C} \to  {\rm Vec}$ have the same dimension $D$, so they are isomorphic by Th. \ref{propweakdim} (a). It follows that ${\mathcal F}'{\mathcal E}$ 
admits 
a weak quasi-tensor structure with weak quasi-bialgebra isomorphic to $A$, thus there is 
a tensor equivalence ${\mathcal E}_1:{\mathcal C}\to{\rm Rep}(A)$ and a
 monoidal isomorphism  ${\mathcal F}'{\mathcal E}\simeq {\mathcal F}_A{\mathcal E}_1 $
with   ${\mathcal F}_A: {\rm Rep}(A)\to{\rm Vec}$ the forgetful functor. On the other hand, we similarly have a   monoidal isomorphism of ${\mathcal F}'\simeq {\mathcal F}_{A'}{\mathcal E}_2$ with
${\mathcal E}_2: {\mathcal C}'\to{\rm Rep}(A')$ a tensor equivalence and  ${\mathcal F}_{A'}: {\rm Rep}(A')\to{\rm Vec}$ the forgtful functor.
Since $A$ is isomorphic to a twist of $A'$, there is a tensor equivalence   ${\mathcal E}_3: {\rm Rep}(A')\to
{\rm Rep}(A)$ and an isomorphism  ${\mathcal F}_A{\mathcal E}_3\simeq{\mathcal F}_{A'}$ by Prop. \ref{class}.
We have an isomorphism of functors ${\mathcal F}_A{\mathcal E}_1\simeq{\mathcal F}_A{\mathcal E}_3{\mathcal E}_2{\mathcal E}$ and since ${\mathcal E}_1$ admits the structure of a tensor equivalence, the same holds for  
${\mathcal E}_3{\mathcal E}_2{\mathcal E}$.   Let ${\mathcal E}'_2$ and ${\mathcal E}'_3$ be quasi-inverse tensor equivalences of ${\mathcal E}_2$ and ${\mathcal E}_3$ respectively. Then ${\mathcal E}'_2{\mathcal E}'_3{\mathcal E}_3{\mathcal E}_2{\mathcal E}$ is a tensor equivalence naturally isomorphic to ${\mathcal E}$ as a linear equivalence, thus ${\mathcal E}$
admits the structure of a tensor equivalence.
\end{proof}

In Sect. \ref{VOAnets2} we shall use weak quasi-Hopf algebras associated to tensor equivalent fusion categories to gain      insight into the study of unitarizability of fusion categories and this will find fruitful applications   to CFT.
    We formulate a   simple   criterion that will eventually   be useful to  construct ribbon tensor equivalences, see Sect. \ref{KW}.

If a weak quasi bialgebra $A'$ is obtained from another such   bialgebra $A$   by replacing the associator of the latter with a new one but leaving the rest of the  structure unchanged, then ${\rm Rep}(A)$ and ${\rm Rep}(A')$ have isomorphic Grothendieck rings. The following   proposition, inspired by a  similar statement in 
\cite{NY_towards} for   Hopf algebras, shows that at an abstract level an isomorphism of Grothendieck rings of fusion categories can always be visualized in this way.

\begin{prop}\label{deforming_Phi}
Let ${\mathcal C}$ and ${\mathcal C'}$ be semisimple  tensor categories and let 
$f:  {\rm Gr}({\mathcal C})
\to {\rm Gr}({\mathcal C}')$ be an isomorphism between their Grothendieck rings.
Let $(A, \Delta,  \Phi')$ be a weak quasi bialgebra corresponding to an integral weak dimension function $D'$ on ${\mathcal C}'$. Then there is an   associator
$\Phi$ for $A$ defining a new weak quasi bialgebra $(A, \Delta, \Phi)$ which corresponds to ${\mathcal C}$ with respect to  $D=D'f$. In particular, if ${\mathcal C}'$ is a finite semisimple tensor category then ${\mathcal C}$   is tensor equivalent to one with the same category and tensor product structure as ${\mathcal C}'$ but possibly different associativity morphisms.
\end{prop}

\begin{proof}  Consider a complete set ${\rm Irr}({\mathcal C}')$ of irreducible objects of ${\mathcal C}'$.
Let ${\mathcal F}':{\mathcal C}'\to{\rm Vec}$ be a weak quasi-tensor functor corresponding to $D'$ and defining $(A, \Delta, \Phi')$. 
By   Theorem  \ref{propweakdim} a weak quasi-tensor structure on ${\mathcal F}'$ is determined by the choice, for $\rho$, $\sigma\in{\rm Irr}({\mathcal  C}')$,  of
(normalized) epimorphisms 
$F_{\rho, \sigma}: {\mathcal F}'(\rho)\otimes{\mathcal F}'(\sigma)\to{\mathcal F}'(\rho\otimes\sigma)$ and
monomorphisms $G_{\rho,\sigma}: {\mathcal F}'(\rho\otimes\sigma)\to{\mathcal F}'(\rho)\otimes{\mathcal F}'(\sigma)$ 
satisfying $F_{\rho, \sigma}\circ G_{\rho,\sigma}=1$. These maps are in turn specified
by the choice of linear maps maps ${G}^{\tau, i}_{\rho, \sigma}: {\mathcal F}'(\tau)\to {\mathcal F}'(\rho)\otimes{\mathcal F}'(\sigma)$, ${F}^{\tau, j}_{\rho, \sigma}: {\mathcal F}'(\rho)\otimes{\mathcal F}'(\sigma)\to{\mathcal F}'(\tau)$ for $\tau\in{\rm Irr}({\mathcal C}')$, 
via $\sum_{\tau, i}{G}^{\tau, i}_{\rho, \sigma}{\mathcal F}'(T^\tau_i)=:{G}_{\rho,\sigma}$ and $\sum_{\tau, i}{\mathcal F}'(S^\tau_i){F}^{\tau, i}_{\rho, \sigma}=:{F}_{\rho, \sigma}$,
where $S^\tau_i\in(\tau, \rho\otimes\sigma)$, $T^\tau_i\in(\rho\otimes\sigma, \tau)$ satisfy $T^\tau_jS^\tau_i=\delta_{i,j}$, $\sum_{\tau, i}S^\tau_iT^\tau_i=1$,
in turn subject to ${F}^{\tau, i}_{\rho, \sigma}{G}^{\upsilon, j}_{\rho, \sigma}=\delta_{\tau, \upsilon}\delta_{i,j}$. 
 Writing $A={\rm Nat}_0({\mathcal F}')$, the coproduct formula of $A$ given in  (\ref{coproduct}) can be written as
  $\Delta(\eta)_{\rho, \sigma}=\sum_{\tau, i}G^{\tau, i}_{\rho, \sigma}\eta_\tau F^{\tau, i}_{\rho, \sigma}$
 by  naturality of 
$\eta$.

Note that  we may establish a bijective correspondence
$\rho\in{\rm Irr}({\mathcal C})\to\rho'\in {\rm Irr}({\mathcal C}')$ and linear isomorphisms $(\tau, \rho\otimes\sigma)\to(\tau', \rho'\otimes\sigma')$. We then set ${\mathcal F}(\rho):={\mathcal F}'(\rho')$, extend ${\mathcal F}$ to a faithful functor
${\mathcal F}:{\mathcal C}\to{\rm Vec}$, and consider the weak quasi-tensor structure of ${\mathcal F}$ defined by the same maps ${F}^{\tau, i}_{\rho, \sigma}$, ${G}^{\tau, i}_{\rho, \sigma}$ under the correspondence $\rho\to\rho'$. It follows that
the corresponding weak quasi bialgebras may be chosen with the same algebra and coproduct structures.
\end{proof}

 \begin{ex}\label{pointed}  Let $G$ be a finite group.
Consider the category ${\mathcal C}={\rm Vec}_G$ of finite dimensional $G$-graded vector spaces with tensor product defined by convolution and trivial associativity morphisms. 
The representation  ring is
  ${\mathbb Z}G$. The constant   function $D=1$   is a dimension function, 
 giving rise to the commutative  bialgebra $C(G)$ of complex functions $f$ on $G$ with  usual coproduct $\Delta(f)(g, h)=f(gh)$.
   Prop. \ref{deforming_Phi} reduces to the known classification of tensor categories with this representation ring.
 Indeed, in this special case it shows that any such category is tensor equivalent to  some ${\rm Vec}^\omega_G$,      obtained from  ${\rm Vec}_G$ with a new associativity morphism given by  
  a normalised ${\mathbb C}^\times$-valued $3$-cocycle $\omega$.
  It corresponds to the quasi-bialgebra $C_\omega(G)$ coinciding with $C(G)$ except for the
    associator, which is given by $\omega$.  
    Since ${\rm Vec}^\omega_G$ is a pointed fusion category, $D=1$ is the only dimension function on  ${\mathbb Z}G$.
    Thus $C_\omega(G)$ is, up to twist, the only   quasi-bialgebra
    that can be associated to ${\rm Vec}^\omega_G$. Twist isomorphism 
  corresponds to cohomologous   cocycles. It follows that the fusion categories
    ${\rm Vec}^\omega_G$ are parameterised by $H^3(G, {\mathbb C}^\times)$.
It also follows that
 ${\rm Vec}^\omega_G$ admits a   faithful   tensor functor to ${\rm Vec}$  if and only if $\omega$ is cohomologically trivial. 
 
   For example, the category ${\rm Vec}^\omega_{{\mathbb Z}_2}$, with $\omega$ the non trivial element of $H^3({\mathbb Z}_2, {\mathbb C}^\times)$, arises from the representation theory of the affine vertex operator algebra associated to ${\mathfrak sl}_2$ at level $1$, a topic that will be discussed in more detail in Sections
   \ref{KW},  \ref{VOAnets}, \ref{VOAnets2}. We shall   come back to this in more detail and generality later on. We shall see that 
 this category also admits 
  a weak tensor functor to ${\rm Vec}$ with weak dimension function   $D(\rho)=2$, and  $\rho$ the unique non trivial irreducible object,
  cf. Example \ref{pointed_case}.
    \end{ex}

  The following   result will be useful to construct   a tensor structure on a given linear equivalence between semisimple tensor categories.
  
  \begin{prop}\label{integration}
  Let ${\mathcal C}$ and ${\mathcal C}'$ be semisimple tensor categories, ${\mathcal G}:{\mathcal C}\to{\mathcal C}'$
  a tensor equivalence and ${\mathcal F}: {\mathcal C}\to{\mathcal C}'$ a linear equivalence. If ${\mathcal F}$ and ${\mathcal G}$ induce the same isomorphism between the corresponding Grothendieck  rings then ${\mathcal F}$
  can be made into a tensor equivalence.
  \end{prop}
  
  \begin{proof} By assumption, for every simple object $\rho\in{\mathcal C}$, ${\mathcal F}(\rho)$ and ${\mathcal G}(\rho)$ are equivalent simple objects in ${\mathcal C}'$, and any simple object of ${\mathcal C}'$ is equivalent to one of them.
  It follows that ${\mathcal F}$ and ${\mathcal G}$ are related by an invertible  natural   transformation $\eta$, and therefore
  ${\mathcal F}$ may be endowed with a unique weak quasi-tensor structure making  $\eta$ monoidal. It also follows that this is a tensor structure for ${\mathcal F}$ since so is the quasi-tensor structure of ${\mathcal G}$.

  \end{proof}
  
  In Sect. \ref{5+} we shall study methods to construct tensor equivalences between tensor categories
  motivated by fusion categories of quantum groups at roots of unity, vertex operator algebras and conformal nets.

  \section{Two abstract uniqueness results of braided tensor structures}\label{5+}

This section aims to discuss at an  abstract level 
braided tensor equivalences between quantum group tensor categories, and tensor categories arising from conformal field theories, in the setting of vertex operator algebras and conformal nets.

The main results of this section are Theorems \ref{claim0} and the more general Theorem \ref{claim1} that apply to show part (c) of Theorem
\ref{Finkelberg_HL}.

In particular, we discuss here  the general theory concerning the full braided tensorial part of   Kazhdan-Lusztig-Finkelberg equivalence. 
With this we   understand the extension of the braided tensor equivalence equations from certain specific
classes of objects for which verification can be performed directly, where a generating object $V$
takes all but one variables in specific coordinates in the braiding and associativity morphisms, to all the objects.
We next motivate the notion of ${\mathcal V}$-pre-associator of CFT type
that we introduce in this section, by explaining how we shall apply it.

We regard   construction of weak quasi-Hopf algebras from fusion categories
of quantum groups  as a generalization of Drinfeld twist method for Drinfeld category, with the aim of leading, in this analogy,
to a direct proof a Kazhdan-Lusztig-Finkelberg theorem.
In our case,   our weak Hopf algebras are regarded as playing the role of Drinfeld-Jimbo quantum groups and the Zhu algebra that
of Drinfeld quasi-Hopf algebra.

In this analogy, we shall develop our Drinfeld-Kohno theorem \ref{Drinfeld_Kohno}.  In analogy to the original Drinfeld-Kohno theorem, our Drinfeld-Kohno
theorem leads to the construction of   a twisted braided symmetry    in the setting of affine vertex operator algebras from fusion categories of quantum groups at roots of unity.

The main difficulty in constructing tensor structures in CFT, is the construction of associativity morphisms. In the original
Drinfeld-Kohno theorem, the associativity morphisms have been constructed by Drinfeld based on the use of the KZ differential equations.
In our analogy, they have  been constructed by Huang and Lepowsky in their tensor product theory.

In our case,
 the Zhu algebra needs the construction of a weak quasi-Hopf algebra compatible with the already existing  braided
tensor structure constructed by  Huang and Lepowky,

Our Drinfeld-Kohno theorem \ref{Drinfeld_Kohno} compares the braiding morphisms, but does not compare the two associativity morphisms,
one obtained as application of our Drinfeld-Kohno theorem and the other associated to Huang-Lepowsky theory.
We do this in this section in abstract form. To this aim,
the main related notion is that of ${\mathcal V}$-pre-associator of CFT-type.
This is a function defined only on special triples of objects that belong to a given family ${\mathcal V}$, with values in the triple tensor product of a semisimple algebra with a (non unital) coproduct.

We  state  two  uniqueness theorems, Theorem \ref{claim0} and Theorem \ref{claim1} on  the associativity morphisms  in presence
of two braided symmetries respectively, which extend the ${\mathcal V}$-pre-associator. The proofs of these theorems are given in \cite{On_a_problem_posed_by_Huang}.  We assume that both the associativity morphisms and the braided symmetry coincide on a subfamily of objects, which is not enough to generate all the objects additively. In the application, reaching this equality is easier with respect to general objects, because of common structural properties
on fusion, when one has corresponding generating representations in the two settings,  and on the braided symmetry, which in turn can be reached by our Drinfeld-Kohno theorem. We shall derive that
these unique associativity morphisms  and braided symmetries coincide with those arising from Huang-Lepowsky theory when
the quantum group structure is transferred to the vertex operator algebra side via our Drinfeld-Kohno theorem and Wenzl quantization continuous curve.

The starting   observation  is an  analogy between the coproduct and associator of a weak  bialgebra arising from a weak-tensor
functor of a   tensor category  to ${\rm Vec}$ and the form taken by the tensor tensor product and associativity morphisms associated to a vertex operator algebra  by Huang and Lepowsky in their papers.  
  
  To make this observation useful, one needs to compare projections that enter into the construction of the coproducts and associativity morphisms in the two cases, associated to triples of arbitrary irreducible representations. This can be done directly for {\it special tensor products} because of similar structural properties, but globally multiplicities make it difficult to select  common projections.
  
Our main tool to this aim is the use of a generating representation useful to reduce the comparison problem
  to a substantially smaller collection of projections, which also benefits of the existence of well defined projections onto
  irreducible components, following the work of Wenzl in \cite{Wenzl}, which plays a key role in our paper.

 Following Wenzl in the setting of quantum groups at roots of unity, we
emphasize the use of a generating object $V$ in the tensor category.  For the classical Lie algebras, $V$   is the defining (vector) representation of the quantum groups
for the Lie types A and C,   the spin representation for the Lie types B and the sum of the two spin representations for the Lie type D.
For the Lie type $G_2$ is the 7-dimensional representation. The fundamental representations for the $E$ and $F$ types are
described in \cite{Wenzl}.

The  braided symmetry is defined by two equations that
 closely link it with the   associativity morphisms in a braided tensor category, and are called hexagonal diagrams.
 
The first immediate consequence of braided symmetry is Proposition \ref{braided_symmetry_with_generating_object}, which shows that the associativity morphisms and the knowledge
of the braided symmetry when one variable is restricted to the generating object and the other is free,
determines the braided symmetry uniquely.

Our Drinfeld-Kohno theorem applied to fusion categories of quantum groups at roots of unity allows to construct   the   braided symmetry which coincides with the braided symmetry
  known  in the setting of  loop groups, or affine vertex operator algebras at a positive integer level.

Given the difficulty of constructing associativity morphisms in the setting of affine vertex operator algebras, one would like to reverse 
 proposition  \ref{braided_symmetry_with_generating_object} and study the following question:  determine a restricted family ${\mathcal V}$
 of triples of objects that as variables   determine the associativity morphisms   on all triples of objects in a unique way.
In the application, such a family reduces the comparison of the two associativity morphisms to this   restricted family, for which
 the comparison is simpler by the properties of the specific generating representations in the Lie types, that have been described in 
 \cite{Wenzl}.

The original observation together with above question motivate
our   definition of pre-associator of CFT-type and the weaker definition of ${\mathcal V}$-pre-associator of CFT-type. 
By a pre-associator we understand a function defined on triples of representations of a discrete algebra with coproduct, which does
not necessarily satisfy the pentagon equation. If it does, then the algebra is a weak bialgebra, with associator analogous to
the trivial case of Hopf algebras.

The case of representation categories of vertex operator algebras at positive integer level, motivates the weaker definition
of ${\mathcal V}$-pre-associator of CFT-type, because we do not know the weak Hopf property in the vertex operator algebra case,
but we show that we have such a pre-associator on special triples of representations and we wish to
have a uniqueness result on the possible associators that extend the restriction of the pre-associator to
${\mathcal V}$.

The main  abstract uniqueness result on associativity morphisms with the same braided symmetry  is Theorem \ref{claim0}.

Our approach to the associativity morphism part is centered on showing vanishing
 of certain cohomological obstructions in the associator of related  weak quasi-bialgebras that we associate, see Theorem \ref{Zhu_from_qg_if_of_CFT_type}. For a more detailed description of our strategy   we refer the reader to the previous   Sect. \ref{5++}.

 \subsection{Definition of ${\mathcal V}$-pre-associator of CFT-type for a  generating object $V$}
 \bigskip

We start with the same scenario and notation as in the basic Tannakian   Theorem \ref{TK_algebraic_quasi}, and 
  remark on the form of the Tannakian coproduct and  associator of the algebra of $A={\rm Nat}_0(\mathcal F)$ induced
  by a faithful weak quasi-tensor functor $({\mathcal F}, F, G)$  $ ({\mathcal C}, \otimes, \iota, \alpha) \to {\rm Vec}$
  that are given in general by  formulas (\ref{coproduct}) and  (\ref{associativity}) respectively.

  \begin{rem}\label{strictness2} Note that the coproduct of $A={\rm Nat}_0(\mathcal F)$ 
  depends only on the given weak quasi-tensor structure
   $(F, G)$ on ${\mathcal F}$ by the formula (\ref{coproduct}),   and the associator formula
   (\ref{associativity})
  depends on    $(F, G)$ but also on the image ${\mathcal F}(\alpha)$ of the associativity morphisms of 
   ${\mathcal C}$.

  We are interested in cases where the dependence of the Tannakian associator of $A$ 
  on ${\mathcal F}(\alpha)$   can be  reduced to
  the dependance on $(F, G)$ only.
  This may happen for several reasons. 
  
The simplest class of examples is the case where $({\mathcal C}, \otimes, \iota, \alpha=1)$ is a strict semisimple tensor category and
  ${\mathcal F}:{\mathcal C}\to{\rm Vec}$ a linear faithful functor. Then    given a weak quasi-tensor structure $(F, G)$ for ${\mathcal F}$, the associator of $A={\rm Nat}_0({\mathcal F})$
   induced by Tannakian duality is
   given by
    \begin{equation}\label{strict_case1_equation}\Phi_{\rho,\sigma, \tau}=1_{{\mathcal F}(\rho)}\otimes G_{\sigma, \tau}\circ G_{\rho, \sigma\otimes\tau}\circ F_{\rho\otimes\sigma,\tau}\circ F_{\rho, \sigma}\otimes 1_{{\mathcal F}(\tau)}.
        \end{equation}
  
  Another class of examples  is that  where $({\mathcal F}, F, G)$ is a weak tensor functor. Then the image
   ${\mathcal F}(\alpha)$ of the associativity morphisms depends only on $(F, G)$, by the formulas (\ref{wt1}), (\ref{wt2}).  
  As a consequence the Tannakian associator  (\ref{associativity}) of $A$  also explicitly depends  only
  on $(F, G)$ in this case.

  Finally, another case that is the main interest of this section,
    is that where Tannakian associator (\ref{associativity}) of $A$ is uniquely determined by its restriction to  a proper subcollection
   of triples of representations of $A$, that we denote by ${\mathcal V}\subset{\rm Ob}({\mathcal C})^{\times 3}$,
   restriction depending only on
 $(F, G)$ as in the previous case of weak tensor functors. 
 
If   the associator of $A={\rm Nat}_0({\mathcal F})$ can be uniquely determined by the weak quasi-tensor structure  $(F, G)$ only
 as in   the examples above considered rather than the full knowledge of
the associativity morphisms 
 of ${\mathcal C}$,   then given another $({\mathcal F}', F', G'): {\mathcal C}'\to{\rm Vec}$ with the same property,
  the most difficult part
 of the  verification of the tensor equivalence property
of an equivalence  between two
  tensor categories ${\mathcal E}: {\mathcal C} \to {\mathcal C}'$ (the associator-preserving property, equations (\ref{wt1}), (\ref{wt2}), Def. \ref{tensor_equivalence2}, \ref{tf}) 
   reduces to a comparison of two  weak quasi-tensor structures $(F, G)$  and $(F', G')$.

 \begin{rem}\label{Tannakian_equivalence_for_pre-tensor_categories}  Recall that pre-tensor  categories   admit a tensor product and a unit
object, but they do not have associativity morphisms in the definition, see Def. \ref{pre_tensor_category}.
The reason why we have introduced them is that we would like to make them into tensor categories.
Conversely, we would also like
to construct    a weak quasi-Hopf algebra associator from a tensor category with a given linear equivalence
${\mathcal E}$
and natural transformation $E$.
 possibly 
different than those given by the Tannakian theorem \ref{TK_algebraic_quasi}.
These matters
are   discussed in  this section and also   in    
 Sect.\ref{12} at an abstract level,   and   applications to CFT in Sect. \ref{21}, \ref{23}, \ref{22}. 
 
 \end{rem}

  Let $({\mathcal C}, \otimes, \iota)$ be a semisimple linear pre-tensor category (as in Def. \ref{pre_tensor_category}) 
and $({\mathcal F}, F, G):{\mathcal C}\to{\rm Vec}$  a faithful weak quasi-tensor functor.  Then the discrete algebra
$A={\rm Nat}_0({\mathcal F})$ has a coproduct $\Delta$
 defined in (\ref{coproduct}), a counit $\varepsilon$ given by the representation corresponding to $\iota$, and in this way
    ${\rm Rep}(A)$ is a pre-tensor category.  
      In this case,
Tannakian duality gives a quasi-equivalence $({\mathcal E}, E):{\mathcal C}\to{\rm Rep}(A)$ of pre-tensor categories
defined as in the case where ${\mathcal C}$ is a tensor category, by an immediate generalization of Theorem
\ref{TK_algebraic_quasi}.  
\end{rem}

\begin{rem}\label{construction_of_associativity_morphism_from_CFT_type_property}

Let $({\mathcal C}, \otimes, \iota)$ be a semisimple pre-tensor category, $({\mathcal F}, F, G):{\mathcal C}\to{\rm Vec}$
a faithful weak quasi-tensor structure.
Usually in Tannakian reconstruction, one starts from associativity morphisms in ${\mathcal C}$ to endow ${\rm Rep}(A)$ 
with associativity morphisms via $({\mathcal F}, F, G)$,  where $A={\rm Nat}_0({\mathcal F})$ has   associated 
coproduct  $\Delta_{F, G}$.
But the construction also works in the other direction.
If $(A, \Delta_{F, G})$ has an associator $\Phi$ such that $(A, \Delta_{F, G}, \Phi)$ is a weak quasi-bialgebra then Tannakian quasi equivalence
$({\mathcal E}, E):  {\mathcal C}\to  {\rm Rep}(A)  $
pulls back $\Phi$ uniquely to associativity morphisms $\alpha$ in ${\mathcal C}$ in such a way that 
$({\mathcal C}, \otimes, \iota, \alpha)$
is a tensor category, and $({\mathcal E}, E)$ becomes a tensor equivalence.
Applying again usual Tannakian construction to $({\mathcal C}, \otimes, \iota, \alpha)$ we obtain  $(A, \Delta_{F, G}, \Phi_{\rm TK})$ 
  and moreover
$\Phi_{\rm TK}=\Phi$.

To verify this,  let $({\mathcal E}, E): {\mathcal C}\to{\rm Rep}(A)$ be the Tannakian quasi-equivalence, cf.
  Remark \ref{Tannakian_equivalence_for_pre-tensor_categories}.
We define the associativity morphisms in ${\mathcal C}$ requiring the associator preserving property (\ref{wt1})  to
  $({\mathcal E}, E, E^{-1})$. Thus
  $\Phi$, ${\rm Rep}(A)$, $({\mathcal E}, E, E^{-1})$  replace   $\alpha'$, ${\mathcal C}'$, $({\mathcal F}, F, G)$ 
  respectively at the r.h.s. of  (\ref{wt1}). 
We get a morphism in ${\rm Rep}(A)$ and we know that Tannakian equivalence is a full functor, thus 
this morphism is in the image of ${\mathcal E}$. Let $\alpha$ be the preimage of this morphisms in ${\mathcal C}$.
  We have ${\mathcal E}(\alpha)=E_{1,2}\Phi E_{2,1}^{-1}$ by definition.
  Then $\alpha$ satisfies naturality, the normalization condition and the pentagon equation, since this holds for the r.h.s
  and ${\mathcal F}$ is faithful.
  Thus   ${\mathcal C}$ becomes a tensor category with associativity morphisms $\alpha$.

    On the other hand, ${\mathcal E}$ acts as
  ${\mathcal F}$, $E$ as $F$ and $E^{-1}$ as $G$. Thus 
   ${\mathcal F}(\alpha)=
 F_{1,2}\Phi G_{2,1}$. The Tannakian associator of
 $A$ corresponding to $\alpha$ is $$\Phi_{\rm TK}=G_{1,2}{\mathcal F}(\alpha)F_{2,1}=G_{1,2}F_{1,2}\Phi G_{2,1}F_{2,1}=\Phi.$$ 
\end{rem}

The following definition, and its weaker version, Def. \ref{CFT_type_associator},
are the starting point  to construct tensor equivalences between fusion categories arising from quantum groups and different formulations of conformal field theories, vertex operator algebras and conformal nets, as explained in the second and third class of examples, respectively
in Remark \ref{strictness2}.

       \begin{defn}\label{CFT_type} ({\it  pre-associator of CFT-type})  Let $({\mathcal C}, \otimes, \iota)$ be a semisimple linear pre-tensor category (as in Def. \ref{pre_tensor_category}) 
and $({\mathcal F}, F, G):{\mathcal C}\to{\rm Vec}$  a faithful weak quasi-tensor functor.  Consider the discrete algebra
$A={\rm Nat}({\mathcal F})$ with coproduct $\Delta$ defined by $(F, G)$ as in (\ref{coproduct}).
We set 
  \begin{equation}\label{associativity2}
  (\Phi_{F, G})_{\rho,\sigma, \tau}=(G_{1,2}F_{1,2})_{\rho, \sigma, \tau}\circ(G_{2,1}F_{2,1})_{\rho, \sigma,\tau}
         \end{equation}
               where
$$(F_{1,2})_{\rho, \sigma,\tau}=F_{\rho, \sigma\otimes\tau}\circ  1_{{\mathcal F}(\rho)}\otimes F_{\sigma, \tau}, \quad\quad (G_{1,2})_{\rho, \sigma, \tau} =1_{{\mathcal F}(\rho)}\otimes G_{\sigma, \tau}\circ G_{\rho, \sigma\otimes\tau}, $$
$$(F_{2,1})_{\rho, \sigma,\tau}=F_{\rho\otimes\sigma,\tau}\circ F_{\rho, \sigma}\otimes 1_{{\mathcal F}(\tau)},\quad\quad (G_{2,1})_{\rho, \sigma, \tau} =G_{\rho, \sigma}\otimes 1_{{\mathcal F}(\tau)} \circ G_{\rho\otimes\sigma, \tau}.$$
Then $\Phi_{F, G}$ is a pre-associator for $A$ as in Def. \ref{pre-associator}. The  pre-associator  $\Phi_{F, G}$ is called {\it of CFT-type.}
    \end{defn}

For simplicity of notation we are dropping  the associator of ${\rm Vec}$, in the middle of (\ref{associativity2}) and 
(\ref{associativity3}).

 When the pre-associator $\Phi_{F, G}$ defined by a weak quasi-tensor functor $({\mathcal F}, F, G):{\mathcal C}\to{\rm Vec}$ as in (\ref{associativity2}) satisfies also the dropped axioms (that is in addition $\Phi_{F, G}$ is partially invertible    and   satisfies
  (\ref{eqn:intro4}), (\ref{eqn:intro6}))
then $(A, \Delta, \Phi_{F, G}, \varepsilon)$ becomes a discrete weak quasi-bialgebra as in Def. \ref{wqh}.

\begin{rem}\label{restriction_of_CFT_type}
By definition, the associativity morphisms $\alpha_{\rho, \sigma, \tau}$ in ${\rm Rep}(A)$ induced by an associator (or a pre-associator) $\Phi$ of a weak quasi-bialgebra $A$, is given by the {\it restriction} of $\Phi$ to the space $V_{(\rho\underline{\otimes}\sigma)\underline{\otimes}\tau}$ of the representation $(\rho\underline{\otimes}\sigma)\underline{\otimes}\tau$ of $A$, with range the space $V_{\rho\underline{\otimes}(\sigma\underline{\otimes}\tau)}$ of  $\rho\underline{\otimes}(\sigma\underline{\otimes}\tau)$. In particular, a
pre-associator of CFT-type $\Phi_{F, G}$ induces associativity morphisms $\alpha_{\rho, \sigma, \tau}$ acting between the spaces of  {\it a given triple} $(\rho, \sigma, \tau)$ of representations just as composition 
$$\alpha_{\rho, \sigma, \tau}={(F_{1, 2})}_{\rho, \sigma, \tau}{(G_{2, 1})}_{\rho, \sigma, \tau}:$$
of two  inclusion maps from
 {\it truncated left-parenthesized tensor products spaces }  to  {\it  full tensor product spaces}
 $${(G_{2, 1})}_{\rho, \sigma, \tau}:
V_{(\rho\underline{\otimes}\sigma)\underline{\otimes}\tau}\to V_{(\rho\underline{\otimes}\sigma) {\otimes}\tau}\to V_{(\rho {\otimes}\sigma) {\otimes}\tau}$$
  with two more projection maps from {\it  full tensor product spaces} to  {\it    right-parenthesized tensor product spaces}
 $${(F_{1, 2})}_{\rho, \sigma, \tau}:$$
 $$
  V_{\rho {\otimes}(\sigma {\otimes}\tau)}\to V_{\rho {\otimes}(\sigma\underline{\otimes}\tau)}\to V_{\rho\underline{\otimes}(\sigma\underline{\otimes}\tau)}.$$
  Similarly, the inverse associativity morphisms
  $$\alpha^{-1}_{\rho, \sigma,\tau}=(F_{2, 1})_{\rho, \sigma, \tau}(G_{1, 2})_{\rho, \sigma, \tau}$$
  act as inclusion maps 
  from 
   {\it   truncated right-parenthesized  tensor product spaces} to   {\it   full tensor product spaces}
 $${(G_{1, 2})}_{\rho, \sigma, \tau}:
V_{\rho\underline{\otimes}(\sigma\underline{\otimes}\tau)}\to V_{\rho{\otimes}(\sigma \underline{\otimes}\tau)}\to V_{\rho {\otimes}(\sigma{\otimes}\tau)}$$
  with two more projection maps from {\it  full tensor product spaces}  to {\it  truncated left-parenthesized }  tensor product spaces
 $${(F_{2, 1})}_{\rho, \sigma, \tau}:$$
 $$
  V_{(\rho {\otimes}\sigma) {\otimes}\tau}\to V_{(\rho \underline{\otimes}\sigma){\otimes}\tau}\to V_{(\rho\underline{\otimes}\sigma)\underline{\otimes}\tau}.$$

\end{rem}

Our approach to compare associativity morphisms in categories of modules of vertex operator algebras with those arising
from quantum group fusion categories or the strict tensor  $C^*$-categories arising from conformal nets categories originated 
from the following question,   briefly mentioned in the introduction.

\begin{rem}\label{how_to_define_F_G}  If  we have two semisimple tensor categories ${\mathcal C}$ and ${\mathcal C'}$ arising from two different settings, and two natural functors ${\mathcal F}: {\mathcal C} \to {\rm Vec}$, ${\mathcal F}': {\mathcal C}' \to {\rm Vec}$,   can we well define
 in a natural way   weak quasi-tensor structures $(F, G)$ and $(F', G')$ in the two cases for all pairs of objects?
 
 For example, if in one case for ${\mathcal F}: {\mathcal C} \to {\rm Vec}$ we succeed to   obtain $(F, G)$ well defined for all pairs of objects
such that the corresponding pre-associator of CFT-type equals the Tannakian associator  of the associated weak quasi-bialgebra, then this will be a weak  Hopf algebra.
 
 If in another   case the algebra associated to ${\mathcal F}': {\mathcal C}' \to {\rm Vec}$   is naturally endowed with an associator  that corresponds to the tensor product and associativity morphisms
 of ${\mathcal C}'$, by some   specific method available in that setting, can we define in a natural way a weak quasi-tensor structure
 $(F', G')$ for ${\mathcal F}'$   that describes the associator as a corresponding
 CFT-type associator?  The answer will be negative if we know of obstacles to obtain a weak  Hopf algebra in that setting.
  Thus there will not be an everywhere well-defined weak quasi-tensor structure $(F', G')$ that corresponds to a CFT-type associator
  for all triples.
  
In the specific setting arising from  vertex operator algebras,
  the tensor product of two irreducible representations is defined starting with
  the fusion rules defined in that setting,
  roughly speaking  a decomposition into irreducibles. By the impressive work by Huang and Lepowsky, the associator
     depends on the fusion rules
  only, thus we have a similar situation for
 the associativity morphisms. 
   
  The Zhu algebra is an invaluable notion to study representation theory in this setting. In the semisimple case,
  there is a natural equivalence of linear categories from representations of the vertex operator algebra to representations of the Zhu algebra, given by Zhu's functor, that corresponds to ${\mathcal F}'$ in our abstraction.
  
   In the affine case at a positive integer level, 
  the dimension of a tensor product representations reported
  to the corresponding Zhu algebra is lower that the product of the two vector space dimensions of the respresentations of the Zhu algebra. The fusion rules describe space of intertwiners of the associated Lie algebra from a finite dimensional
  representation to a tensor product of  other two,
  but this does not immediately define a   subrepresentation of the tensor product  in general, because of the  different multiplicities occurring in the decomposition of the fusion tensor product and of the Lie algebra.
  References and more explanations may be found in the second  part of this paper.
  
  One thus meets the problem of constructing a natural inclusion of this fusion tensor product into the full tensor product of the two respresentations of the Lie algebra, covariant for the action of the simple Lie algebra.
  
  This problem can be solved in a coherent way passing to the quantum group at roots of unity, where to our knowledge the analogous problem was first
  described and studied by Wenzl for some special tensor products of representations  \cite{Wenzl}.
    
  Once one is able to embed  tensor products of specific pairs  of representation of the Lie algebra that correspond to representations of the vertex operator algebra with the correct fusion rules,
  then one can also embed the fusion tensor product of suitable specific triples of representations in  the corresponding full
   tensor product of representations of the Lie algebra, starting with a preferred    parenthesization which indicates the embedding following the previous step for pairs of specific representations.
Then we change the parenthesization using the associator in the setting of affine vertex operator algebras based on the fusion rules, and then again embed
   into a full tensor product of Lie algebra representation using again the construction from the quantum group.
   In this way, the associator arising from the setting of vertex operator algebras, for that specific triple
   acts as a CFT-type associator.
  
 More in detail,  the case of the fusion category of quantum groups at roots of unity is similar, and first motivated our observations above. The corresponding discussion
  may be found in \cite{Wenzl}, where a solution is found for tensor products of representations of the form $V_\lambda\otimes V$ or $V\otimes V_\lambda$, with $V_\lambda$ arbitrary and irreducible
  and $V$ a specific generating representation which provides multiplicity free decompositions (except for a few cases including $E_8$,  handled separately).

  In our work we have extended the projections corresponding to the fusion rules of these special tensor products
  to arbitrary pairs of representations, by extending Wenzl projections to a weak tensor functor from the fusion category of quantum groups at roots of unity to ${\rm Vec}$. This was possible thanks to the quantum group $U_q({\mathfrak g})$ and most importantly to its non-semisimple
  structure, which is of great help to define the full weak tensor structure $(F, G)$ of Wenzl functor ${\mathcal F}$.

Thus such maps $(F,G)$ or  $(F', G')$ initially   well defined as naturally associated to the tensor products
of ${\mathcal C}$ and ${\mathcal C}' $ only  for a restricted class of pairs of representations, and covariant
   with respect to the action of the quantum group at roots of unity or a simple Lie algebra, in the quantum group case can be extended to
  a full weak tensor structure for ${\mathcal F}$. Then   an impressive idea of Drinfeld can be applied to this case, and give a Drinfeld twist
  that describes the  weak quasi-tensor structure for Zhu functor ${\mathcal F}'$ in such a way that
  the corresponding pre-tensor structure on the Zhu algebra identifies with that arising from the tensor product theory by Huang and Lepowski in the setting of vertex operator algebras by application of the Drinfeld twist. At this point we apply the above argument on the associator on special triples of 
  representations.
  
  Then the question becomes that
of asking whether a restricted class of triples
  of representations 
will suffice to identify the  rest of the structure completely, the braiding and the associativity morphisms.

Due to the fact that we can apply an analogue of the original Drinfeld twist method (a Drinfeld-Kohno theorem that we develop in our setting), we shall obtain a positive answer
for the braiding by construction, on the special pairs of representations, where explicit computations indeed identify the braiding for such tensor products.

 A positive answer to our uniqueness question on braiding and associativity morphisms  will also be helpful
to relate and compare this construction, and especially the associativity morphisms, which are the most complicated part of the structure, with the associativity morphisms arising from quantum groups at roots of unity,   provided we relate the pair $(F, G)$ to $(F', G')$ for those special pairs of representations.

To study the question, a first step
will be that of    identifying pairs $(\rho, \sigma)$ of irreducible
representations for which $(F, G)$ and $(F', G')$ are well defined and can be identified.

As said, our  approach to the study of these questions in   our paper is to follow the work by Wenzl   in the setting of quantum groups at roots of unity and   then apply to the setting of both fusion categories of quantum groups at roots of unity 
and that of CFT, mainly categories of modules of vertex operator algebras.

As said, Wenzl approach starts  by  fixing one of the variables as being a generating representation $V$ (which can be reducible
in the application in certain cases, for example the sum of the two spinor representations in the type $D$ fusion categories associated to quantum groups or affine vertex operator algebras, but this is  not cause of difficulty)  and the other  an arbitrary irreducible representation. 

Our reasoning on the difficulty of well defining inclusion and
projection maps $(F, G)$    is  
suggested by the situation arising in the setting of quantum groups at roots of unity and described in \cite{Wenzl}. In that case as already mentioned, an everywhere well defined weak  tensor structure, and therefore a weak  Hopf algebra, is possible thanks to non-semisimplicity of $U_q({\mathfrak g})$ at roots of unity, see  \cite{CP} for the type $A$ case and Sect. \ref{19} for all Lie types.

In this section we do abstract work to apply later to these weak  Hopf algebras arisig from the fusion category of quantum groups at roots of unity.
By the generating property of $V$, one may equivalently work with
all tensor powers $V^r$  of the generating representation in place of an arbitrary irreducible representation,
  which have the virtue of allowing an analysis on a minimal family that determines the braiding and the associativity morphisms.
  
  It will be important to take into consideration further work    in the literature that we shall see,
concerning the study of the intertwining spaces in representation theory of quantum groups.
This is what in type $A$ case is known as Schur-Weyl duality, that is a property that describes braiding morphisms as generating
for tensor powers of the fundamental representation $V$. 
We shall need to use this property for all Lie types for which it is known to hold,
in the case of quantum groups.

In the setting of affine vertex operator algebras at a positive integer level, a weak tensor structure will not be possible in general, due to obstacles arising from
the study of the relations between amenability and unitary structure in our setting, see Sect. \ref{13}.
We shall use the methods developed in this subsection to determine the braided symmetry and associativity morphisms completely.

Our results on having found a weak  Hopf algebra on one side and a weak quasi-Hopf algebra
on the other side related by an isomorphism and a specific Drinfeld twist given by a square root of a coboundary matrix,
 is completely analogous to the case of the original Drinfeld-Kohno theorem, especially for the form of the twist, for the category of representations of a quantum group and Drinfeld category \cite{Drinfeld_cocommutative}, \cite{Drinfeld_quasi_hopf}. This theorem first motivated Mack and Schomerus \cite{MS} to introduce weak quasi-Hopf algebras, and also our approach to  KLF theorem. Except for   our case seems more direct than the original Drinfeld-Kohno theorem, in that the twist determines the associator 
 of the module category of affine vertex operator algebras completely.
\end{rem}\medskip

   \begin{prop}\label{wh_as_CFT_type}
Let $({\mathcal C}, \otimes, \iota, \alpha)$ be a semisimple tensor category, ${\mathcal F}:{\mathcal C}\to{\rm Vec}$ a linear faithful functor
 and $(F, G)$  a weak quasi-tensor structure for ${\mathcal F}$.  Then the associator $\Phi_{\rm TK}$ defined in (\ref{associativity}) for $A={\rm Nat}_0({\mathcal F})$     by the Tannakian Theorem \ref{TK_algebraic_quasi}  
coincides with the pre-associator $\Phi_{F, G}$ of CFT-type (\ref{associativity2})  with inverse given by 
 \begin{equation}\label{associativity3}
 ( \Phi_{F, G})_{\rho,\sigma, \tau}^{-1}=(G_{2,1}F_{2,1})_{\rho, \sigma, \tau}\circ(G_{1,2}F_{1,2})_{\rho, \sigma,\tau}.
         \end{equation}
 if and only if $(F, G)$ is a weak tensor structure.
 In particular  $\Phi_{F, G}$  is an associator.
 
  \end{prop}
 
 \begin{proof} 
   The condition $FG=1$ easily shows that the associator  of ${\rm Nat}_0({\mathcal F})$ 
 defined by the Tannakian Theorem \ref{TK_algebraic_quasi}, see (\ref{associativity}), coincides with (\ref{associativity2})
with inverse given by (\ref{associativity3})
if and only if ${\mathcal F}(\alpha)$ is given  by the relation (\ref{wt1}) with inverse satisfying
 (\ref{wt2}), and this is the definition of a weak tensor structure.
 The last statement follows from the fact that  the associator defined by (\ref{associativity}) automatically satisfies the pentagon equation.

    \end{proof}

The algebras defined by weak tensor structures as in 
Prop. \ref{wh_as_CFT_type},   are called weak  bialgebras, (or weak  Hopf algebras depending on the existence of an antipode). In Sect. \ref{6} we shall see their basic properties. We shall construct examples associated to quantum groups at roots of unity later on, extending the result of \cite{CP} to all Lie types.
 
Thus Def. \ref{CFT_type} of pre-associator of CFT-type is useful to provide a tensor category only when the given weak quasi-tensor structure $(F, G)$
 is a weak tensor structure. We next  consider a   definition  enriched with a collection of triples
 ${\mathcal V}$ of objects of ${\mathcal C}$  that is more flexible, and still provides an associator, possibly for genuine weak quasi-bialgebras,  and thus still leads to a tensor category.

 \begin{defn}\label{CFT_type_associator} 
 ({\it  ${\mathcal V}$-pre-associator of CFT-type})  Let $({\mathcal C}, \otimes, \iota)$ be a semisimple linear pre-tensor category (as in Def. \ref{pre_tensor_category}),
  $({\mathcal F}, F, G):{\mathcal C}\to{\rm Vec}$  a faithful weak quasi-tensor functor. 
  Consider the discrete algebra
$A={\rm Nat}({\mathcal F})$ with coproduct $\Delta$ defined by $(F, G)$ as in (\ref{coproduct}) and corresponding pre-associator of CFT-type $\Phi_{F, G}$ as in Def. \ref{CFT_type}. Given a collection ${\mathcal V}$ in ${\rm Ob}({\mathcal C})^{\times 3}$
we say that the pre-associator $\Phi_{F, G}$ is a ${\mathcal V}$-pre-associator of CFT-type
if there is a     $\Phi\in A\otimes A\otimes A$   
 such that 
\begin{itemize} 
 \item[(a)] 
$\Phi$ is an associator for  $\Delta$ 
(thus $(A, \Delta, \Phi)$ is a discrete weak quasi-bialgebra) 
 \item[(b)]    $\Phi$ restricts to    $\Phi_{F, G}$  
 on ${\mathcal V}$.
 
 \end{itemize}

 \end{defn}

If ${\mathcal V}={\rm Ob}({\mathcal C})^{\times 3}$ then the definition gives $\Phi_{F, G}=\Phi$,
 and  $({\mathcal F}, F, G)$ is a weak tensor functor.   
 
 \subsection{A first uniqueness result of associativity morphisms  for  semisimple braided tensor categories 
 with a
generating object satisfying braid group duality with completely fixed braiding}    

Let  $({\mathcal C}, \otimes, \iota)$ be a semisimple pre-tensor category with a faithful weak quasi-tensor functor
 $({\mathcal F}, F, G)$ to ${\rm Vec}$.
 If $\Phi$ is an   associator     satisfying Def. \ref{CFT_type_associator}
 for a given pre-associator
 $\Phi_{F, G}$ on a given collection  ${\mathcal V}$,
 then $(A={\rm Nat}_0({\mathcal F}), \Delta, \varepsilon, \Phi)$ is a discrete weak quasi-bialgebra and therefore ${\rm Rep}(A)$ is a tensor category. By Remark \ref{construction_of_associativity_morphism_from_CFT_type_property},
 $({\mathcal C}, \otimes, \iota)$ can be made uniquely into a tensor category $({\mathcal C}, \otimes, \iota, \alpha)$ under the requirement that Tannakian quasi-equivalence ${\mathcal E}: {\mathcal C}\to{\rm Rep}(A)$ become a tensor equivalence.

  In how many ways  can we form tensor categories in this way with the same pre-associator $\Phi_{F, G}$
  on ${\mathcal V}$?
  We next describe a uniqueness result for such an extension associator $\Phi$ of the restriction of the pre-associator $\Phi_{F, G}$ to ${\mathcal V}$ that will turn out important for our applications.

\begin{defn}      Let $V^r$ be a  tensor power    of an object $V$ corresponding to some fixed parenthesization. We call $r$ the order of  $V^r$,  independently of   the parenthesization defining $V^r$.
\end{defn}

\begin{ex}\label{Examples_of_different_tensor_powers_of_the_same_order_obtained_from_braiding}
Let $V^r$ be a  tensor power    of $V$ with order $r>1$. Thus 
  $V^r$  is the tensor product of two further tensor powers of $V$ of smaller orders:
 $V^r=V^{r_1}\otimes V^{r_2}$ with $r_i\geq1$. Then $V^{r_2}\otimes V^{r_1}$ is an example of another tensor power of $V$ of the same order as the original $V^r$.

If  $r_1>1$, we can apply the above procedure to 
$V^{r_1}$. So writing $V^{r_1}=V^{s_1}\otimes V^{s_2}$ with $s_1$, $s_2\geq 1$ then 
$(V^{s_2}\otimes V^{s_1})\otimes V^{r_2}$ in another tensor power of $V$ of the same order as $V$.

If $r_2>1$ and $V^{r_2}=V^{t_1}\otimes V^{t_2}$ then $V^{r_1}\otimes (V^{t_2}\otimes V^{t_1})$ is another   example of tensor power of $V$ of order the same as that of the original $V^r$.

We may inductively further decompose every tensor power of $V^k$ of order $k>1$, decompose it into a tensor product of two further tensor powers of $V$, $V^k=V^{k_1}\otimes V^{k_2}$ with $k_i\geq 1$
and then reverse the order of these factors in the tensor product $V^{k_2}\otimes V^{k_1}$ in place of $V^k$.

\end{ex}

 \begin{defn} \label{braid_group_generating_property} Let 
  $({\mathcal C}, \otimes, \iota)$ be a semisimple linear pre-tensor category (Def. \ref{pre_tensor_category})
 and let $$c(\rho, \sigma): \rho\otimes\sigma\to\sigma\otimes\rho$$ be an invertible natural transformation
 such that $c(\iota, \sigma)$ and $c(\rho, \iota)$ are identity morphisms for all objects $\rho$ and $\sigma$.
 
 We shall say that an object $V$ of ${\mathcal C}$ satisfies {\it the   duality property with respect to $c$} if  
 for any positive integer $r\geq1$ the morphism space $(V^r, V^r)$ between any two   tensor powers of $V$  of the kind of Example
 \ref{Examples_of_different_tensor_powers_of_the_same_order_obtained_from_braiding}
   is linearly generated by finite compositions of  the component morphisms of $c$, their tensor products with identity morphisms and tensor products of morphisms in $(V, V)$ with identity morphisms.
   
   We shall say that two tensor powers  of $V$ of the same order of the kind described  in \ref{Examples_of_different_tensor_powers_of_the_same_order_obtained_from_braiding} belong to the same orbit under the action of $c$.

 \end{defn}   

   Definition \ref{braid_group_generating_property} is independent of the associativity morphisms. But in the case where
 $({\mathcal C}, \otimes, \iota, c)$
  has associativity morphisms $\alpha$ making it into a   tensor category
   $({\mathcal C}, \otimes, \iota, \alpha)$  then
 any   tensor power of $V$ of the same order as that of a given $V^r$, is isomorphic to $V^r$ via an isomorphism given by composition of associativity morphisms. So if $V$ satisfies the   duality property  with respect to $c$  then   the structure
 of morphism spaces between any two arbitrary tensor powers of $V$ of the same order is known.
 Any such morphism
 $T\in (V^r, V^r)$ is of the form
 $$T=A_1BA_2$$
 with $A_1$ and $A_2$ suitable associativity morphisms and $B$
  a linear combination of a composition of   braiding morphisms, morphisms in $(V, V)$ and their compositions 
 with identity morphisms. 
 If the range of $A_2$ and the domain of $A_1$ are two fixed objects in the same orbit under the action of $c$, then $A_1$ and $A_2$ are
 uniquely determined by McLane coherence, and hence so is $B$.

 \begin{defn} \label{braid_group_generating_property2}
 Let $({\mathcal C}, \otimes, \iota, c)$ be a pre-tensor category with invertible natural transformation
 $c(\rho, \sigma): \rho\otimes\sigma\to\sigma\otimes\rho$. 
 Let $V\in {\mathcal C}$ be an object satisfying the duality property with respect to $c$.
 Let $\alpha$ be associativity morphisms in ${\mathcal C}$ making $({\mathcal C}, \otimes, \iota, \alpha, c)$
 into a braided tensor category.
 {\it We shall say that $V$ satisfies the braid group duality property. }
 
 Two tensor powers 
 of $V$ of the same order of the kind described in Examples
 \ref{Examples_of_different_tensor_powers_of_the_same_order_obtained_from_braiding} will be referred to as {\it in the same orbit under the action of the braid group.}
 \end{defn}

Examples of generating objects $V$ satisfying the braid group duality property, or of
tensor powers of   $V$ in the same orbit under the action of the braid group
   with respect to given associativity morphisms    appear 
in  Theorem \ref{claim0}.

We next describe an important   case of non-trivial collection of triples ${\mathcal V}$,
that we consider for applications in this paper.
Let ${\mathcal V}$ be the collection of triples of objects
\begin{equation}\label{definition_of_V}
 {\mathcal V}=\{(V_\lambda, V, V), (V, V_\lambda, V), (V, V, V_\lambda), \quad V_\lambda\in{\rm Irr}({\mathcal C})\}.\end{equation}

    \begin{thm}\label{claim0}
 Let $({\mathcal C}, \otimes, \iota)$ be a semisimple pre-tensor category with  a generating object $V$ and admitting a faithful weak quasi-tensor functor $({\mathcal F}, F, G): {\mathcal C}\to{\rm Vec}$ into the category of finite dimensional vector spaces.
 
  Let  $c(\rho, \sigma): \rho\otimes\sigma\to\sigma\otimes\rho$ be a normalized invertible natural transformation and let
  $V$ satisfy the   duality property with respect to $c$.

 Let $\alpha$ and $\beta$ be two associativity morphisms for $({\mathcal C}, \otimes, \iota)$ such that $({\mathcal C}, \otimes, \iota, \alpha, c)$ and $({\mathcal C}, \otimes, \iota, \beta, c)$ are braided tensor categories.
 
Let ${\mathcal V}$ be defined as in (\ref{definition_of_V}).

Assume that     
 \begin{equation}\label{assumption_on_associators2}
 \alpha=\beta \quad\text{on }
 {\mathcal V}.\end{equation}  
  Then  $\alpha=\beta$ everywhere.  
 \end{thm}

 The proof of this theorem relies on a detailed combinatorial analysis of the pentagon and hexagon coherence axioms across the truncated tensor powers. To preserve the focus of the present   framework, the proof of this theorem as well as that of its generalization Theorem \ref{claim1}  are are presented in full   detail in Section 10 of our companion paper  \cite{On_a_problem_posed_by_Huang}.

     \begin{rem}\label{reduction_of_assumptions} It follows from the proof of Theorem \ref{claim0} given in Theorems 10.4 and 10.6 of \cite{On_a_problem_posed_by_Huang} that
     we may weaken the starting assumption on $\alpha$ and $\beta$ in Theorem    \ref{claim0} without modifying the rest
by    only
 assuming  that they coincide on one of the  smaller families of triples of the form 
 $${\mathcal V}_l:=\{(V_\lambda, V, V), (V, V_\lambda, V)\},$$ 
 or  
 $${\mathcal V}_r:=\{(V, V_\lambda, V), (V, V, V_\lambda)\}$$
 and derive coincidence on the other as a consequence along the proof  as follows.
         
The only point in the proof where we need  to know coincidence of $\alpha$ and $\beta$ on both $(V_\lambda, V, V)$ and $(V, V, V_\lambda)$ with $V_\lambda$ arbitrary is at the beginning on the inductive assumption on $g$.
    The argument of paragraph l) of the proof for $\ell=1$ leads to the claim. This paragraph
  does not use the assumption that $\alpha$ and $\beta$ coincide
   on both $(V^u, V, V)$ and   $(V, V, V^u)$ for all $u>1$.  
      The part where we verify that the $T$-maps are defined in the same way in the intersection objects
   works using again the pentagon equation in the same way. But we do not need to use any inductive assumption on $f$
   to this aim in this case, but only the fact that on triples $(V, V, V)$ and $(V, V^2, V)$ we know that the two associators coincide by assumption. 
   We may anticipate     paragraph l) in the case
   $\ell=1$ to the inductive assumption on $g$ at the beginning of the proof, to reduce this assumption to coincidence
   on only one among ${\mathcal V}_l$ or ${\mathcal V}_r$.

       In other   cases there may be other structural reasons that imply     coincidence of $\alpha$ and $\beta$ on one of ${\mathcal V}_l$ or ${\mathcal V}_r$ knowing coincidence on the other.
    For example, we shall see that this is always the case if we have a coboundary symmetry, which we do have in the application of the fusion category of quantum groups at roots of unity for all Lie types, see Sect. \ref{17}.

  \end{rem}

By Remark \ref{reduction_of_assumptions}  one of the assumptions on coincidence 
   of $\alpha$ and $\beta$ on   triples with the arbitrary term on the left  $(V_\lambda, V, V)$  or on the right $(V, V, V_\lambda)$, is redundant in the presence of the rest of the assumptions.

    The following corollary   follows  from  Theorem \ref{claim0}.

\begin{cor}\label{determination_of_associator_with_generating_object}
 Let $({\mathcal C}, \otimes, \iota)$ be a semisimple linear pre-tensor category (Def. \ref{pre_tensor_category})
 with a generating object $V$. 
    Let  $$c(\rho, \sigma): \rho\otimes\sigma\to\sigma\otimes\rho, \quad\quad \rho, \sigma\in{\mathcal C},$$ be a normalized invertible natural transformation and let
  $V$ satisfy the   duality property with respect to $c$. Let
  $({\mathcal F}, F, G):{\mathcal C}\to{\rm Vec}$  a faithful weak quasi-tensor functor. 
Then an   associator $\Phi$ of $A={\rm Nat}_0({\mathcal F})$ making $c$ into a braided symmetry
and extending the pre-associator $\Phi_{F, G}$ on ${\mathcal V}$ is unique.

\end{cor}

 \begin{rem}\label{on_assumption_on_associators}({\it    On the method of verification 
 of the assumptions (\ref{assumption_on_associators2}) of Theorem \ref{claim0} in the application}.)
In our application we shall verify  the assumptions on the associativity morphisms   of Corollary 
\ref{determination_of_associator_with_generating_object} to Theorem \ref{claim0}, 
see Theorem \ref{Zhu_from_qg_if_of_CFT_type} and Sect. \ref{33}.
 
 More precisely we shall show property
  (\ref{assumption_on_associators2}) by showing that
   two specific  associativity morphisms $\alpha$, $\beta$ arising respectively from fusion categories associated to
   quantum groups at roots of unity and from the tensor product theory of Huang and Lepowsky
   for affine vertex operator algebras at a positive integer level, satisfy 
 \begin{equation}\label{verification} \alpha_{V_\lambda, V, V}=\beta_{V_\lambda, V, V}, \quad \alpha_{V, V_\lambda,  V}=\beta_{V, V_\lambda, V}, \quad
 \alpha^{-1}_{V, V, V_\lambda}=\beta^{-1}_{V, V, V_\lambda}.\end{equation}
 
 These relations will be obtained by verifying the equations of the weak tensor functor property
 (\ref{wt1}) for the first two identities and equation (\ref{wt2}) for
 the last identity on the special triples. More precisely, we shall verify that the two associators (or their inverse according to
 the triples) both act as described
 in Remark  \ref{restriction_of_CFT_type} for each of triple in ${\mathcal V}$ for the same inclusion and projection maps
 $F$ and $G$.
 
 Note that this verification for the corresponding weak quasi-Hopf algebra $A$ does not necessarily imply that we have a 
 weak tensor functor  (or a weak  Hopf algebra) as we are considering the property only on a restricted class of triples.
As already anticipated, and we shall see in detail in the second part of the paper,  our verification of  (\ref{verification}) depends on the form
taken by
 the associativity
morphisms on the associative discrete algebras associated to the two categories in question, the fusion category of quantum groups at roots of unity and  the tensor category of affine vertex operator algebras at a positive integer level.
They both take the form of ${\mathcal V}$-preassociators of CFT type for the same weak tensor structure $(F, G)$ on the same functor
${\mathcal F}$.

Depending on whether we are reducing the verification on the side of the quantum group fusion category, or on the side of vertex operator algebra fusion category, the corresponding $A$ will be a weak  Hopf algebra in the first case, or the Zhu algebra
with a weak quasi-Hopf structure   in the vertex operator algebra case.

Indeed, for the case of quantum group at roots of unity, we shall construct weak  Hopf algebras for all Lie types in this paper
and we shall transfer on it also the braided tensor structure from Huang-Lepowsky tensor product theory via an isomorphism
and a the construction of a Drinfeld twist.

 For the case of affine vertex operator algebras,
 the construction of braided symmetry and associativity morphisms is due to the
highly non-trivial work by Huang-Lepowsky on the construction
 of vertex tensor categories, and we may transfer our braided tensor structure
of the weak  Hopf algebras to the Zhu algebra  with a procedure inverse
  of the previous case. 
  
  In either setting one decides to work, either on the quantum group side or the vertex operator algebra side,
  we shall reduce our work to verify the assumptions of Theorem \ref{claim0} for the occurring weak (quasi)-Hopf algebra 
  $A$.
(See also   part c) of Remark \ref{sufficient_conditions_for_associativity_morphisms}, where we comment
on the fact that these equations imply that $\beta$ extends uniquely to associativity morphisms satisfying the pentagon equation.)
 \end{rem}

 \begin{rem} \label{on_assumption_on_braided_symmetry}  ({\it On the method of verification of the   assumption of a common braided symmetry of Theorem \ref{claim0} in the application.})
 On the other hand,
the assumption that the braided symmetry is the same for the two associators as required in Theorem
\ref{claim0}, does not meet precisely our approach in the application, where
  we shall limit ourselves to verify only that the two braided
symmetries arising from the setting of quantum groups  at roots of unity and that
of affine vertex operator algebras at positive integer level, coincide on tensor products of the form $V_\lambda\otimes V$ and $V\otimes V_\lambda$ on associated weak quasi-Hopf algebras.
While more direct work in the setting of vertex operator algebras might be possible to reach the assumption of Theorem \ref{claim0} on coincidence of
 the two full  braided symmetries in the application, we shall refrain from doing this in our applicative work, in the attempt to
emphasize general properties in common for the two setting, quantum groups at roots of unity and vertex operator algebras.

To reach coincidence of the two full  braided symmetries and two corresponding associativity morphisms simultaneously from knowledge of part of them, in the next subsection we
state  a second general uniqueness result, Theorem \ref{claim1} in the  setting of braided tensor categories, which will suffice for our application, for reasons similar to those anticipated in Remark \ref{on_assumption_on_associators}. The proof of this theorem is deferred to Theorem 10.6 in \cite{On_a_problem_posed_by_Huang} as well.
\end{rem}

   \begin{rem}  
 In  \cite{KW} Kazhdan and Wenzl classify fusion categories with Grothendieck semiring isomorphic to the semiring associated to  ${\mathfrak sl}_N$ or to the fusion categories
 of the associated quantum group at suitable roots of unity (also known as with Verlinde fusion rules), as equivalent to those arising from quantum group $U_q({\mathfrak sl}_N)$ with
 $q$ generic or a suitable root of unity respectively with associator modified by scalar factor. We have recalled their result in Theorem \ref{Kazhdan_Wenzl}.    In Theorem \ref{classification_type_A}  we have shown how to make this
   scalar trivial, and therefore how to fix  the associator to be the usual trivial associator of  fusion categories of quantum groups at roots of unity, by fixing  the   ribbon structure on the Grothendieck semiring as a datum.    
      
   There is a   similarity between   Theorem \ref{classification_type_A}  and Theorem \ref{claim0}, in the role of the braiding to determine the associator and also in
   the occurrence of the   semisimple quotients $\tilde{A}=\tilde{A}_{N, M}$ and the invertible central elements $U=U_{N, M}\in \tilde{A}_{N, M}$. These central elements appearing in our proof eventually define a central element in $A$, and in this sense they
    remind us of   the scalar factor appearing in 
   Kazhdan-Wenzl theory  in the Lie type $A$ of quantum group fusion categories.
   
   Both Theorem \ref{classification_type_A} and Theorem \ref{claim0} use discrete weak quasi-Hopf algebras in a crucial way, and more than this, in Theorem
    \ref{classification_type_A} we need weak Hopf algebras.
    
    On the other hand, for Lie types different from $A$, Tuba and Wenzl \cite{TW} have studied a similar question as in  \cite{KW}, although with some
    differences in the starting assumptions. Recently, Yamashita has announced results closely related to 
 Kazhdan-Wenzl theory and Tuba-Wenzl theory for the classical Lie types $B$, $C$, $D$ in a joint work with Grossman and Neshveyev. 
     
    On the other hand, we have weak discrete Hopf algebras available also for all Lie types beyond A.
   One might thus look for more connections between the two approaches. For example, try to see whether our central element $U$ of the weak quasi-Hopf algebra $A$ which differentiates the two associators $\alpha$ and $\beta$ by a scalar factor  appearing in the proof of Theorem \ref{claim0} (see Theorem 10.4 in \cite{On_a_problem_posed_by_Huang}) is related to the scalar factors appearing in these theories or try
   to draw the conclusions of
   the proof of Theorem \ref{claim0} along  lines similar to Sect. \ref{KW} for all Lie types.
    \end{rem}

\subsection{  A second uniqueness result of braided tensor structure  for  semisimple braided tensor categories with a
generating object satisfying braid group duality and braiding fixed on special 
pairs of objects}

As above remarked,  Theorem \ref{claim0} gives a uniqueness result for associativity morphisms of a braided tensor category 
 $(\alpha, c)$ among associativity morphisms $(\beta, c)$
all admitting the same braided symmetry $c$ from the start.
We next give   a uniqueness theorem for   braided tensor structures   $(\alpha, c)$ among pairs
$(\beta, d)$, where we let
not only  the associativity morphisms but also the braided symmetry vary, assuming coincidence of a part of their components
that we shall verify in our application.
This result will prove that the braided symmetries coincide simultaneously with proving that the associativity morphisms
coincide.

The  uniqueness result of this subsection Theorem \ref{claim1}   
starts with  weaker assumptions on the braided symmetries
assuming some more knowledge on the associativity morphisms. This formulation will find verification in our applicative work with methods already
 summarized in Remarks \ref{on_assumption_on_associators} and \ref{on_assumption_on_braided_symmetry}

To describe what we shall need to know more on the associativity morphisms 
we start with special kind of associativity morphisms $\alpha'$ and $\beta'$ on {\it $n$ variables} arising from $\alpha$ and $\beta$ respectively.
 
\begin{defn}

Let us consider associativity morphisms for $n\geq 4$:
$$\alpha'_{W_1, W_2, W_3, \dots, W_n}: ((W_1\otimes W_2)\otimes W_3)\otimes\dots )\otimes W_n \to W_1\otimes (\dots \otimes(W_{n-2}\otimes(W_{n-1}\otimes W_n))$$
that pass from left-parenthesized   to the right-parenthesized tensor products of the   objects $W_i$.
Let us also consider their inverses
$$(\alpha'_{W_1, W_2, W_3, \dots, W_n})^{-1}: 
W_1\otimes (\dots \otimes(W_{n-2}\otimes(W_{n-1}\otimes W_n))\to ((W_1\otimes W_2)\otimes W_3)\otimes\dots )\otimes W_n$$
that pass from right-parenthesized   to the left-parenthesized tensor products of the   objects $W_i$.
We refer to such associativity morphisms  $\alpha'_{W_1, W_2, W_3, \dots, W_n}$   on $n$ variables and their inverses
$(\alpha'_{W_1, W_2, W_3, \dots, W_n})^{-1}$
 as {\it extremal}.
 \end{defn}

Extremal associativity morphisms $\alpha'_{W_1, W_2, W_3, \dots, W_n}$ on $n$ variables are   explicitly determined by compositions of  
$n-2$ associativity morphisms $\alpha\circ\dots\circ\alpha$ on  three variables. 

Let us consider the   example   $n=4$, which is instructive for the following result.
In this case, 
we have two ways of defining extremal associativity morphisms,
following  the two possible paths of
the pentagon equation
(\ref{pentagon_equation}) connecting the right-parenthesized   with the left-parenthesized tensor products of the four objects $W_i$.
We may either
 compose $\alpha\circ\alpha$, the     left vertical with bottom horizontal associativity morphisms,   or compose the other three (the two upper horizontal with
the right vertical) associativity morphisms in the same diagram, and we obtain the desired extremal associativity morphisms.

Given an  object $V$ of a semisimple pre-tensor category ${\mathcal C}$, let us consider the collection
${\mathcal V}'$ of quadruples of objects where three coordinates equal $V$, and the fourth is free:

\begin{equation}
{\mathcal V}':=\{(V_\lambda, V, V, V), (V, V_\lambda, V, V), (V, V, V_\lambda, V), (V, V, V, V_\lambda), \quad\quad
V_\lambda\in{\rm Irr}({\mathcal C})\}.
\end{equation}

We also consider the subcollection ${\mathcal V}''\subset{\mathcal V}'$ where the arbitrary object $V_\lambda$
takes only the first or last coordinate

\begin{equation}{\mathcal V}'':=\{(V_\lambda, V, V, V), \quad (V, V, V, V_\lambda), \quad\quad
V_\lambda\in{\rm Irr}({\mathcal C})\}\end{equation}

   \begin{prop}\label{characterization_of_extremality_on_V'}
   Let $({\mathcal C}, \otimes, \iota)$ be a semisimple pre-tensor category with  an object $V$.
   Let $\alpha$ and $\beta$ be associativity morphisms making $({\mathcal C}, \otimes, \iota)$ into a tensor category such that
  $\alpha=\beta$ on ${\mathcal V}$. Then the following properties are equivalent:
   \begin{itemize} 
   \item [(a)]  $\alpha=\beta$ in addition on triples of the form
   $\{(V_\lambda, V^2, V),  (V, V^2, V_\lambda),  \quad V_\lambda\in{\rm Irr}({\mathcal C})\}$,
   \item [(b)]   extremal associativity morphisms satisfy
 $\alpha'=\beta'$ on  ${\mathcal V}'$,
 \item [(c)]   extremal associativity morphisms satisfy
 $\alpha'=\beta'$ on  ${\mathcal V}''$.

      \end{itemize} 
    \end{prop}
    
    \begin{proof}
    Let us define extremal associativity morphisms $\alpha'$ of four variables using the composition of the three associativity morphisms $\alpha$
    in three variables defined in the upper path of the pentagon equation (\ref{pentagon_equation}).
    Then property (b) becomes
    $$1\otimes\alpha\circ\alpha\circ\alpha\otimes1=1\otimes\beta\circ\beta\circ\beta\otimes1 \quad \text{on } {\mathcal V}'.$$
    The first and last factors of this equation in $\alpha$ on the left hand side coincide with the corresponding factors in $\beta$ at the right hand side, by the assumption 
      $\alpha=\beta$ on ${\mathcal V}$. It follows that (b) is equivalent to
      $$\alpha_{\rho, (\sigma\otimes\tau), \nu}=\beta_{\rho, (\sigma\otimes\tau), \nu},  \quad (\rho, \sigma, \tau, \nu)\in {\mathcal V}'.$$
      On quadruples of  ${\mathcal V}'$ for which the arbitrary object $V_\lambda$ takes the second or third coordinate
      this equation 
      holds  by semisimplicity
      of ${\mathcal C}$ and the assumption $\alpha=\beta$ on ${\mathcal V}$.
      On quadruples of ${\mathcal V}'$ for which the arbitrary object $V_\lambda$ takes the first or last coordinate, this equation becomes
      property (a).
      This argument also shows that   (b) is equivalent to (c) under our assumption.
    \end{proof}

    \begin{thm}\label{claim1}
 Let $({\mathcal C}, \otimes, \iota)$ be a semisimple pre-tensor category with  a generating object $V$ and admitting a faithful weak quasi-tensor functor $({\mathcal F}, F, G): {\mathcal C}\to{\rm Vec}$ into the category of finite dimensional vector spaces.
 
 Let  $c(\rho, \sigma): \rho\otimes\sigma\to\sigma\otimes\rho$ be a normalized invertible natural transformation and let
  $V$ satisfy the   duality property with respect to $c$.
  
   Let  $d(\rho, \sigma): \rho\otimes\sigma\to\sigma\otimes\rho$ be another normalized invertible natural transformation.
   
    Let $\alpha$ and $\beta$ be two associativity morphisms for $({\mathcal C}, \otimes, \iota)$ such that $({\mathcal C}, \otimes, \iota, \alpha, c)$ and $({\mathcal C}, \otimes, \iota, \beta, d)$ are braided tensor categories.

Let ${\mathcal V}$ be defined as in (\ref{definition_of_V}).

Assume that for all 
$\lambda \in{\rm Irr}({\mathcal C})$: 

\begin{itemize} 

\item[(a)] $c(V_\lambda, V)=d(V_\lambda, V)$ and $c(V, V_\lambda)=d(V, V_\lambda)$ .

\item[(b)]
 \begin{equation}  
 \alpha=\beta \quad\text{on }
 {\mathcal V}.\end{equation}

 \item[(c)] Extremal associativity morphisms  $\alpha'$ and $\beta'$  coincide on   $1+k$-tuples of the form
 $(V^r, V, V, \dots, V)$ with $V$ repeated a number of times $k\leq r$, for all $r\geq2$,

  \item[(d)] Extremal associativity morphisms  $(\alpha')^{-1}$ and $(\beta')^{-1}$  coincide on all $h+1$-tuples of the form
 $(V, V, \dots, V, V^s)$ with $V$ repeated a number of times $h\leq s$, for all $s\geq2$.

\end{itemize}

  Then $\alpha=\beta$ and $c=d$ everywhere.
 
 \end{thm}
 
 As for Theorem \ref{claim0}, the proof of Theorem \ref{claim1} is given as Theorem 10.6 in \cite{On_a_problem_posed_by_Huang}.
 
 \begin{rem}\label{on_assumption_on_claim1}({\it    On the method of verification 
 of the assumptions    of Theorem \ref{claim1} in the application}.)
 In Remarks \ref{on_assumption_on_associators}  we have summarized
 how we shall verify property (\ref{assumption_on_associators2})     of Theorem \ref{claim0} in the application. This is the same property
as   (b)  of Theorem  \ref{claim1}. In Remark \ref{on_assumption_on_braided_symmetry}  we have commented on the usefulness of our assumption (a) in Theorem  \ref{claim1} for the application, and we shall verify this in Sect. \ref{21}.
 Regarding verification of properties (c) and (d) of the same theorem, we shall show that the two associativity morphisms
 $\alpha$ and $\beta$ arising respectively from quantum groups and vertex operator algebras at positive integer level
 satisfy  
 $$\alpha'_{V_\lambda, V, \dots, V}=\beta'_{V_\lambda, V, \dots, V}, \quad {\alpha'}^{-1}_{V, \dots, V, V_\lambda}={\beta'}^{-1}_{V, \dots, V, V_\lambda}, \quad
 V_\lambda\in{\rm Irr}({\mathcal C}).$$
 To this aim, we shall use arguments   summarized in Remarks \ref{restriction_of_CFT_type}, \ref{how_to_define_F_G} and  \ref{on_assumption_on_associators},
  for $n$-variable  extremal associativity morphisms on special $n$-tuples of representations.\end{rem}
 \bigskip

\subsection{Strictification}
    
     We conclude the section with a general discussion on the construction of tensor equivalences between tensor categories focusing on the case where ${\mathcal V}$-pre-associators of CFT-type are available. The central result
     for our purposes is Cor. \ref{tensor_equivalence_CFT_strict}.
    
      We next describe the functorial dependence of the weak quasi-bialgebra $A={\rm Nat}_0({\mathcal F})$  on the inducing category ${\mathcal C}$.
    
  \begin{prop}
  Let $({\mathcal E}, E, E'): {\mathcal C}\to {\mathcal D}$ be a  weak quasi-tensor functor between semisimple tensor categories  and let ${\mathcal F}:{\mathcal D}\to{\rm Vec}$ be a weak quasitensor functor with structure   $(F_2, G_2)$. Then ${\mathcal G}={\mathcal F}{\mathcal E}: {\mathcal C}\to{\rm Vec}$ is a weak quasi tensor functor with
  structure  
  \begin{equation}\label{composed_functor} (F_1)_{\rho, \sigma}:={\mathcal F}(E_{\rho, \sigma})(F_2)_{{\mathcal E}(\rho), {\mathcal E}(\sigma)},\quad\quad (G_1)_{\rho, \sigma}:=(G_2)_{{\mathcal E}(\rho), {\mathcal E}(\sigma)}{\mathcal F}(({E'}_{\rho, \sigma})).\end{equation}
  Moreover  a tensor equivalence $({\mathcal E}, E, E'=E^{-1})$  induces an isomorphism $\phi: A_{\mathcal F}\to A_{\mathcal G}$ between the weak quasi bialgebras
  associated to    $({\mathcal F}, F_2, G_2)$ and $({\mathcal G}, F_1, G_1)$  via Tannakian duality as follows.
 For $\eta\in{\rm Nat}_0({{\mathcal F}})$, $\phi(\eta):\rho\in{\mathcal {\mathcal C}}\to
 \eta_{{\mathcal E}(\rho)}\in ({\mathcal G}(\rho), {\mathcal G}(\rho))$.  
  \end{prop}

  We next reverse the question, and ask how to upgrade a 
 {linear equivalence between semisimple tensor categories }
  \begin{equation}   {\mathcal E}: {\mathcal C}\to {\mathcal D}  
  \end{equation}
to a tensor equivalence.

\begin{ex}\label{Example_strict} Let ${\mathcal E}: {\mathcal C}\to{\mathcal D}$ be a linear equivalence 
between semisimple tensor categories with ${\mathcal D}$ strict.
The simplest solution, that is 
$E_{\rho, \sigma}=1$ for all objects $\rho$, $\sigma$,  is possible only if ${\mathcal E}$ is multiplicative on objects,
${\mathcal E}(\rho)\otimes{\mathcal E}(\sigma)={\mathcal E}(\rho\otimes\sigma)$.   
By  (\ref{wt1}), (\ref{wt2}), $({\mathcal E}, 1)$ is a tensor equivalence if and only if ${\mathcal C}$ is strict as well.

  \end{ex}

    In the applications to CFT in the setting of vertex operator algebras, tensor categories are not strict.
   If we try to use Mac Lane strictification theorem then we add many undesired objects to the categories, thus  Example \ref{Example_strict}
    may    turn too restrictive.
We look for general solutions that may be of help for one of our
 main questions.
   On one hand, conformal nets   give rise to strict tensor categories, as tensor product is described by composition
   of localized endomorphisms in Doplicher-Haag-Roberts theory, thus we are reduced to study the comparison
   between a strict tensor category and a non-strict tensor category. Moreover we relax the condition $E=1$ to pairs of simple objects.
  We    write down the equations that must hold
in a way that
  the   natural transformation $E =1$ holds only on pairs of simple objects extends to a tensor equivalence.
  This holds if and only if   the following equation
  (\ref{associativity_with_a_strict_category}) is satisfied by the associativity morphisma. Then point (d) shows that this obstruction may be vanished by a
  pre-associator of CFT-type, and this is the central theme of the following result and Cor.
  \ref{tensor_equivalence_CFT_strict}.
  
  In this paper, we are interested in the following cases. In particular, case (a) although very simple,
  will be applied in applications to vertex operator algebras, and in particular to transfer tensor structure
  of a module category of a vertex operator algebra to the Zhu algebra via Zhu functor. An abstract form will be studied in Sect. \ref{5++}, with results in Sect. \ref{12} and then used in Sect. \ref{23}, see Theorem \ref{Finkelberg_HL}, Sects. \ref{Interlude}, \ref{32}, \ref{33}.

   \begin{thm}\label{strict_equivalence} Let ${\mathcal E}: ({\mathcal C},\otimes, \iota)\to({\mathcal D}, \otimes, \iota)$ be a linear equivalence between semisimple   pre-tensor categories  such that ${\mathcal E}(\rho_i\otimes\rho_j)={\mathcal E}(\rho_i)\otimes{\mathcal E}(\rho_j)$
    with $\{\rho_i\}$ a complete set of    simple objects in ${\mathcal C}$. Let $E_{\rho, \sigma}: {\mathcal E}(\rho)\otimes{\mathcal E}(\sigma)\to{\mathcal E}(\rho\otimes\sigma)$ be the (invertible) natural transformation extending $1_{{\mathcal E}(\rho_i)\otimes{\mathcal E}(\rho_j)}$ by naturality and let 
   ${\mathcal F}:{\mathcal D}\to{\rm Vec}$ be a faithful weak quasi-tensor functor
 with structure $(F_2, G_2)$. Set
 $$(F_1)_{\rho_i, \rho_j}:=(F_2)_{{\mathcal E}(\rho_i), {\mathcal E}(\rho_j)}, \quad (G_1)_{\rho_i, \rho_j}:=(G_2)_{{\mathcal E}(\rho_i), {\mathcal E}(\rho_j)}.$$
 Then  
 \begin{itemize} 
 \item [(a)] The simplest solution $E_{\rho, \sigma}=1$ for all objects is a tensor equivalence if and only if 
 ${\mathcal E}(\alpha)=\beta$;
 \item[(b)] formulas (\ref{composed_functor}) extend uniquely
 $(F_1)_{\rho_i, \rho_j}, (G_1)_{\rho_i, \rho_j}$   to a weak quasi tensor structure for  
 ${\mathcal G}={\mathcal F}{\mathcal E}$; 
 \item[(c)] ${\mathcal E}$ induces an isomorphism of  algebras $\phi: {\rm Nat}_0({\mathcal F}) \to{\rm Nat}_0({\mathcal G})$ via $\phi(\eta)_\rho:=\eta_{{\mathcal E}(\rho)}$ 
 that intertwines the coproducts induced by $(F_2, G_2)$ and $(F_1, G_1)$ via Tannakian duality respectively.
  \item[(d)]
  Assume that $({\mathcal C},\otimes, \iota, \alpha)$ and $({\mathcal D},\otimes, \iota, \beta)$ are tensor categories and that ${\mathcal D}$ is strict ($\beta=1$). Then
   $\phi$ is an isomorphism of weak quasi-bialgebras associated via Tannakian duality if and only if  $({\mathcal E}, E)$ is a tensor equivalence
 and this holds by definition if and only if
 \begin{equation}\label{associativity_with_a_strict_category}
 {\mathcal E}(\alpha_{\rho_i, \rho_j, \rho_k})=E_{\rho_i, \rho_j\otimes\rho_k}(E_{{\rho_i\otimes\rho_j}, \rho_k})^{-1}.
 \end{equation} 
   \item[(e)]
 $\phi$ is automatically an isomorphism between the   pre-associators 
 $\Phi_{F_2, G_2}$   for ${\rm Nat}_0({\mathcal F})$,
 and  $\Phi_{F_1, G_1}$   for ${\rm Nat}_0({\mathcal G})$ of CFT-type, that is 
  $\phi\otimes\phi\otimes\phi(\Phi_{F_2, G_2})=\Phi_{F_1, G_1}$.
    
  \end{itemize}

   \end{thm}

   \begin{proof} (a) is an immediate consequence of the definition of tensor equivalence.
   (b) By  tensor multiplicativity
   of ${\mathcal E}$ on simple objects,  
   $(F_1)_{\rho_i, \rho_j}:=(F_2)_{{\mathcal E}(\rho_i), {\mathcal E}(\rho_j)}$ acts between ${\mathcal G}(\rho_i)\otimes{\mathcal G}(\rho_j)\to{\mathcal G}(\rho_i\otimes\rho_j)$
   and $(G_1)_{\rho_i, \rho_j}:=(G_2)_{{\mathcal E}(\rho_i), {\mathcal E}(\rho_j)}$ 
   acts between ${\mathcal G}(\rho_i\otimes\rho_j)\to {\mathcal G}(\rho_i)\otimes{\mathcal G}(\rho_j)$
   and   they satisfy  $(F_1)_{\rho_i, \rho_j}(G_1)_{\rho_i, \rho_j}= 1_{{\mathcal G}(\rho_i\otimes\rho_j)}$. Thus $F_1$, $G_1$ extend by naturality to   a weak quasi-tensor structure for ${\mathcal G}$. By construction,
   (\ref{composed_functor}) holds. In this way   
    ${\mathcal G}$   becomes a weak quasitensor functor.
   (c) Let ${\mathcal G}={\mathcal F}{\mathcal E}: {\mathcal C}\to{\mathcal D}\to{\rm Vec}$ be the composed functor.
   For $ \eta\in{\rm Nat}_0({\mathcal F})$, $\phi(\eta):\rho\in{\mathcal C}\to\eta_{{\mathcal E}(\rho)}\in ({\mathcal G}(\rho), {\mathcal G}(\rho))$ is a natural transformation of ${\mathcal G}$. Thus $\phi: {\rm Nat}_0({\mathcal F})\to{\rm Nat}_0({\mathcal G})$ is an algebra isomorphism.
     We regard $A={\rm Nat}_0({\mathcal F})$ and $B={\rm Nat}_0({\mathcal G})$ with the weak bialgebra structures
      $(\Delta_A, \Phi_A)$, $(\Delta_B, \Phi_B)$ defined by $(F_2, G_2)$ and $(F_1, G_1)$ respectively as in (\ref{coproduct}), (\ref{associativity}).
      By duality, ${\mathcal C}$ identifies with ${\rm Rep}(B)$ and ${\mathcal D}$ with ${\rm Rep}(A)$.
       Moreover  the intertwining relation  between coproducts $\Delta_B\phi=\phi\otimes\phi\Delta_A$ easily follows from tensor multiplicativity of ${\mathcal E}$ on tensor products of simple objects.
       and $\phi\otimes\phi\otimes\phi(\Phi_A)=\Phi_B$ on pairs or triples of simple objects of ${\mathcal C}$. 
       (d) If ${\mathcal D}$ is strict
       the associator of $A$ is uniquely determined by   $(F_2, G_2)$
      as in (\ref{strict_case1_equation}).  We have
       $$\phi\otimes\phi\otimes\phi(\Phi_A)_{\rho_i, \rho_j, \rho_k}=(\Phi_A)_{{\mathcal E}(\rho_i), {\mathcal E}(\rho_j), {\mathcal E}(\rho_k)}=$$
   $$ 1_{{\mathcal F}({\mathcal E}(\rho_i))}\otimes (G_2)_{{\mathcal E}(\rho_j), {\mathcal E}(\rho_k)}\circ (G_2)_{{\mathcal E}(\rho_i), {\mathcal E}(\rho_j)\otimes{\mathcal E}(\rho_k)}\circ (F_2)_{{\mathcal E}(\rho_i)\otimes{\mathcal E}(\rho_j),{\mathcal E}(\rho_k)}\circ (F_2)_{{\mathcal E}(\rho_i), {\mathcal E}(\rho_j)}\otimes 1_{{\mathcal F}({\mathcal E}(\rho_k))}=
    $$
    $$1_{{\mathcal G}(\rho_i)}\otimes (G_1)_{\rho_j, \rho_k}\circ (G_1)_{\rho_i, \rho_j\otimes \rho_k}\circ{\mathcal F}(E_{\rho_i, \rho_j\otimes\rho_k})\circ {\mathcal F}(E_{\rho_i\otimes\rho_j, \rho_k}^{-1})(F_1)_{ \rho_i\otimes \rho_j, \rho_k}\circ (F_1)_{\rho_i, \rho_j}\otimes 1_{{\mathcal G}(\rho_k)}$$
    and equals
    $$\Phi_B=1_{{\mathcal G}(\rho_i)}\otimes (G_1)_{\rho_j, \rho_k}\circ (G_1)_{\rho_i, \rho_j\otimes \rho_k}\circ
    {\mathcal F}({\mathcal E}(\alpha_{\rho_i, \rho_j, \rho_k}))\circ
    (F_1)_{ \rho_i\otimes \rho_j, \rho_k}\circ (F_1)_{\rho_i, \rho_j}\otimes 1_{{\mathcal G}(\rho_k)}$$
    if and only if $({\mathcal E}, E)$ is a tensor equivalence and if and only if (\ref{associativity_with_a_strict_category}) holds.
    (e) We follow the notation of Def. 
       \ref{CFT_type}. The previous computations show that for any triple of simple objects $\rho_i$, $\rho_j$, $\rho_k$,
    $$((G_2)_{1,2})_{{\mathcal E}(\rho_i), {\mathcal E}(\rho_j), {\mathcal E}(\rho_k)}=((G_1)_{1,2})_{\rho_i, \rho_j, \rho_k}\circ {\mathcal F}(E_{\rho_i, \rho_j\otimes \rho_k})$$
    $$((F_2)_{1,2})_{{\mathcal E}(\rho_i), {\mathcal E}(\rho_j), {\mathcal E}(\rho_k)}={\mathcal F}(E_{\rho_i, \rho_j\otimes \rho_k}^{-1})\circ((F_1)_{1,2})_{\rho_i, \rho_j, \rho_k} $$
    $$((G_2)_{2,1})_{{\mathcal E}(\rho_i), {\mathcal E}(\rho_j), {\mathcal E}(\rho_k)}=((G_1)_{2,1})_{\rho_i, \rho_j, \rho_k}\circ {\mathcal F}(E_{\rho_i\otimes\rho_j,  \rho_k})$$
    $$((F_2)_{2,1})_{{\mathcal E}(\rho_i), {\mathcal E}(\rho_j), {\mathcal E}(\rho_k)}={\mathcal F}(E_{\rho_i\otimes \rho_j,  \rho_k}^{-1})\circ((F_1)_{2,1})_{\rho_i, \rho_j, \rho_k}.$$
    Thus the compositions 
    $((G_2F_2)_{1,2})_{{\mathcal E}(\rho_i), {\mathcal E}(\rho_j), {\mathcal E}(\rho_k)}$, 
  $((G_2F_2)_{2,1})_{{\mathcal E}(\rho_i), {\mathcal E}(\rho_j), {\mathcal E}(\rho_k)}$
  eliminate the factors coming from the natural transformation $E$, that is
  $$((G_2F_2)_{1,2})_{{\mathcal E}(\rho_i), {\mathcal E}(\rho_j), {\mathcal E}(\rho_k)}=((G_1F_1)_{1,2})_{\rho_i,
   \rho_j,  \rho_k},$$, 
  $$((G_2F_2)_{2,1})_{{\mathcal E}(\rho_i), {\mathcal E}(\rho_j), {\mathcal E}(\rho_k)}=((G_1F_1)_{2,1})_{\rho_i,
   \rho_j,  \rho_k}.$$
 It follows that
   $$\phi\otimes\phi\otimes\phi(\Phi_{F_2, G_2})=\Phi_{F_1, G_1}.$$

         \end{proof}

         \begin{cor}\label{tensor_equivalence_CFT_strict} Let $({\mathcal C}, \otimes, \iota)$
         and $({\mathcal D}, \otimes, \iota)$ be  semisimple pre-tensor categories  and $$({\mathcal F}, F_2, G_2):{\mathcal D}\to {\rm Vec}$$ a weak quasi-tensor functor. Let ${\mathcal E}:{\mathcal C}\to{\mathcal D}$ be a linear equivalence such that $${\mathcal E}(\rho_i\otimes\rho_j)={\mathcal E}(\rho_i)\otimes{\mathcal E}(\rho_j)$$
  on a complete set $\{\rho_i\}$ of simple objects of ${\mathcal C}$. Consider the   weak quasitensor structure
  $(F_1, G_1)$ on ${\mathcal G}={\mathcal F}{\mathcal E}$ as in Theorem \ref{strict_equivalence}.
   Assume that  
    $$\Phi_{F_1, G_1}, \quad \Phi_{F_2, G_2} \quad \text{are} \quad {\mathcal V}_1\text{-}\quad \text{and}\quad  {\mathcal V}_2\text{-}\text{pre-associators of CFT-type}$$
     $$\text{on}\quad A={\rm Nat}_0({\mathcal G}), \quad  B={\rm Nat}_0({\mathcal F})\quad \text{resp.},$$
 with ${\mathcal E}({\mathcal V}_1)={\mathcal V}_2$. Let us regard $A$ and $B$ as a weak quasi-bialgebras with their defining associators, see Def. \ref{CFT_type_associator}.
Then
  ${\mathcal E}$ induces a tensor equivalence ${\rm Rep}(A)\to{\rm Rep}(B)$. Equipping ${\mathcal C}$ and ${\mathcal D}$ with tensor category structures as in Prop. \ref{construction_of_associativity_morphism_from_CFT_type_property} gives a tensor equivalence
  ${\mathcal C}\to{\mathcal D}$.
 \end{cor}
         
         \begin{proof}
         This follows from Theorem \ref{strict_equivalence} (b), (c), (d), Prop.  \ref{construction_of_associativity_morphism_from_CFT_type_property}.

         \end{proof}

    In the applications  we would like to   compare conformal field theories representation categories between themselves
    and  to quantum group fusion categories. 
    
    We take   ${\mathcal C}$  arising from
  vertex operator algebras and ${\mathcal G}$ is Zhu's functor. We shall see that $F_1$ arises from tensor product
  theory of Huang and Lepowsky, and $G_1$ from the unitary structure, and we shall study this in Sect. \ref{12}, \ref{21}, \ref{22}, \ref{23}.
  
We study two    classes of examples of categories ${\mathcal D}$ that satisfies the assumptions of the Cor. \ref{tensor_equivalence_CFT_strict}. The first  
is motivated by the theory of conformal nets.

The second arises from quantum groups, we are in the setting of
Kazhdan-Lusztig-Finkelberg theorem. In this case ${\mathcal F}$ is Wenzl functor, $F_2$, $G_2$ is a weak tensor structure on this functor that we describe in Sect. \ref{20}, \ref{21}, \ref{22}, \ref{23} related to   Wenzl fusion tensor product. In this case as we shall see  $F_2^*\neq G_2$.
  On the other hand, by a suitable twist $T$, we shall try to make the twisted structure $F_2^T$, $G_2^T$ satisfy $(F_2^T)^*=G_2^T$ at least locally, in a suitable sense that we shall explain later on.
  \bigskip

   \begin{prop}\label{CFT_type_condition_strict_case}
    Let $({\mathcal C}, \otimes, \iota, \alpha=1)$ be a semisimple strict unitary tensor
      category, $({\mathcal F}^s,  F^s, G^s):{\mathcal C}\to {\rm Hilb}^s$ a weak quasitensor $^*$-functor   to a strictified
      category of Hilbert spaces,
        such that
        for any pair of objects $\rho$, $\sigma\in{\mathcal C}$, ${\mathcal F}^s(\rho\otimes\sigma)$ is a subspace
       of ${\mathcal F}^s(\rho)\otimes{\mathcal F}^s(\sigma)$, $G^s_{\rho, \sigma}: {\mathcal F}^s(\rho\otimes\sigma)\to
        {\mathcal F}^s(\rho)\otimes{\mathcal F}^s(\sigma)$
       is the inclusion map and $F^s_{\rho, \sigma}:   
        {\mathcal F}^s(\rho)\otimes{\mathcal F}^s(\sigma)\to {\mathcal F}^s(\rho\otimes\sigma) $ is
       the orthogonal projection   with respect to the usual  inner product 
  of ${\mathcal F}^s(\rho)\otimes{\mathcal F}^s(\sigma)$. Then 
  \begin{equation}\label{CFT_type_strict_case1}
  F^s_{\rho\otimes\sigma, \tau}\circ F^s_{\rho, \sigma}\otimes 1_{{\mathcal F}^s(\tau)}=F^s_{\rho, \sigma\otimes\tau}\circ 1_{{\mathcal F}^s(\rho)}\otimes F^s_{\sigma, \tau},
  \end{equation}
   \begin{equation} \label{CFT_type_strict_case2}    
  G^s_{\rho, \sigma}\otimes 1_{{\mathcal F}^s(\tau)}\circ 
   G^s_{\rho\otimes\sigma, \tau}=1_{{\mathcal F}^s(\rho)}\otimes G^s_{\sigma, 
    \tau}\circ G^s_{\rho, \sigma\otimes\tau}.\end{equation}  
    In particular $({\mathcal F}^s, F^s, G^s)$ is a weak tensor functor.
       Let $({\mathcal F}, F, G)$ be the composition of $({\mathcal F}^s,  F^s, G^s)$ with a tensor equivalence
       ${\rm Hilb}^s\to{\rm Hilb}$. Then   $A={\rm Nat}_0({\mathcal F})$ endowed with
     the structure induced by the Tannakian theorem \ref{TK_algebraic_quasi} has associator of CFT-type and   the associator on ${\mathcal C}$ 
     defined
   following Prop.  \ref{construction_of_associativity_morphism_from_CFT_type_property} 
  coincide with the original trivial associativity morphisms $1$.
      
      \end{prop}
      
     \begin{proof}
    Both sides of (\ref{CFT_type_strict_case2}) are identity maps on the same subspace  of a Hilbert space by the strictness properties, thus they
       coincide. Equation (\ref{CFT_type_strict_case1}) follows 
         from (\ref{CFT_type_strict_case2}) taking the adjoint.
      In the notation of Def. \ref{CFT_type} applied to a strict category of vector spaces, $F^s_{1,2} =
  F^s_{2,1}$, $G^s_{1, 2}=  G^s_{2,1}$.
     It follows that $F^s_{1,2} \circ G^s_{2,1}=F^s_{2,1}\circ G^s_{2,1}=1$, thus this is a weak tensor structure on ${\mathcal F}^s$ to ${\rm Hilb}^s$.
The composed structure $(F, G)$ is a weak tensor structure for the composed functor ${\mathcal F}$ to ${\rm Hilb}$, thus $(F)_{1,2}\alpha^{\rm Hilb}G_{2,1}={\mathcal F}(1)=(F)_{2,1}(\alpha^{\rm Hilb})^{-1}G_{1,2}$. It follows   that the associator (\ref{associativity}) given by Tannakian equivalence to $({\mathcal F}, F, G)$ is of CFT-type.    Since Tannakian duality is a tensor equivalence, the last statement follows 
     from Prop.  \ref{construction_of_associativity_morphism_from_CFT_type_property}.
     
   \end{proof}

   \begin{rem}
   If the subspaces ${\mathcal F}^s(\rho\otimes\sigma)$ 
   coincide with ${\mathcal F}^s(\rho)\otimes{\mathcal F}^s(\sigma)$
   then the assumptions of Prop. \ref{CFT_type_condition_strict_case} reduce to the requirement that
   ${\mathcal F}$ is a strict tensor functor. When
   ${\mathcal C}$ is a unitary symmetric strict tensor category then ${\mathcal C}$ admits a unique symmetric strict 
   tensor functor to ${\rm Hilb}^s$ by the Doplicher-Roberts duality theorem \cite{DR1}, the original proof is based
   on the  strictification of ${\rm Hilb}$ induced by the Cuntz algebras \cite{Cuntz}.
   
   \end{rem}

 We next discuss a strictification of tensor categories with a generating object, see Def.\ref{Def_generating_object}.
      The class of tensor categories with a generating object   may be seen as motivated by geometry.
   For example, if ${\mathcal C}$ is the category of unitary representations of a compact group $G$ on finite dimensional Hilbert spaces then ${\mathcal C}$ admits a generating object if and only if it is a Lie group.

 We shall use natural  generating objects for affine VOAs and quantum group fusion categories following the work of Wenzl. More generally,
 any fusion category admits a generating object given by the   sum of the simple objects.  
   We wish to explain in some detail  a strictification of a  semisimple tensor category with a generating 
   object. In this case, the set of objects of the strictified category may be kept under control.

Let ${\mathcal C}$ be a semisimple tensor category and $\rho\in{\mathcal C}$   an  object.   Then we may pass to the full subcategory ${\mathcal C}_\rho$ of ${\mathcal C}$ with objects   tensor powers of $\rho$ with different paranthesisations. 
If $\rho$ is a generating object,   ${\mathcal C}_\rho$ is tensor equivalent to ${\mathcal C}$ with the inclusion map, but is not strict yet. We may strictify ${\mathcal C}_\rho$ taking objects under control.
   That is, let ${\mathcal C}_\rho^{s}$ be the full subcategory of ${\mathcal C}_\rho$ with
   objects $\iota$, $\rho$, $\rho\otimes\rho$, $(\rho\otimes\rho)\otimes\rho$ and so on, that is we put parentheses on the left, and we denote by $\rho_n$ the $n$-th tensor power of $\rho$ defined in this way.
   We regard ${\mathcal C}_\rho^s$ as a linear category and define a new  tensor structure on ${\mathcal C}_\rho^{s}$ by $\rho_n\underline{\otimes}\rho_m:=\rho_{n+m}$.
   For $S\in(\rho_n, \rho_{n'})$, $T\in(\rho_m, \rho_{m'})$
   we set  
   $$S\underline{\otimes} T: \rho_{n+m}\xrightarrow{\alpha}\rho_n\otimes\rho_m\xrightarrow{S\otimes T}\rho_{n'}\otimes\rho_{m'}\xrightarrow{\alpha}\rho_{n'+m'}$$
   where $\alpha$ denote the unique morphisms that can by obtained as compositions of associativity morphisms 
   of ${\mathcal C}$ with identity isomorphisms, by Mac Lane coherence theorem.
   
   We next discuss the braiding.
 Let $c$ be a braided symmetry for ${\mathcal C}$ and define
$$\overline{c}(\rho_n, \rho_m):\rho_n\underline{\otimes}\rho_m=\rho_{n+m}\xrightarrow{\alpha}\rho_n\otimes\rho_m\xrightarrow{c(\rho_n, \rho_m)}\rho_m\otimes\rho_n\xrightarrow{\alpha}\rho_m\underline{\otimes}\rho_n=\rho_{m+n}.$$

Let us see the ribbon structure. Let $v\in(1,1)_{\mathcal C}$ be a ribbon structure for the braided tensor category ${\mathcal C}$
and regard  $\overline{v}:=v\in(1,1)_{{\mathcal C}_\rho^s}$ as a natural transformation of
the identity functor of   ${\mathcal C}_\rho^s$.

\begin{thm}\label{strictification} Let ${\mathcal C}$ be a semisimple tensor category with generating object $\rho$.
The category ${\mathcal C}_\rho^s$ is a strict tensor category. Moreover $({\mathcal H}, H): {\mathcal C}_\rho^s\to{\mathcal C}$ is a tensor equivalence, where
${\mathcal H}: {\mathcal C}_\rho^s\to{\mathcal C}_\rho$ is the inclusion 
 and
$${H}_{\rho_n, \rho_m}:=\rho_n\otimes\rho_m\to\rho_{n+m}$$  the unique associativity morphism given by Mac Lane coherence.
If $c$ is a braided symmetry for ${\mathcal C}$   then $\overline{c}$ is a braided symmetry for ${\mathcal C}_\rho^s$;   if ${\mathcal C}$ is rigid then ${\mathcal C}_\rho^s$ is rigid; if    $v$ is a ribbon structure for ${\mathcal C}$ then $\overline{v}$ is a ribbon structure for ${\mathcal C}_\rho^s$. In these cases, $({\mathcal H}, H)$ becomes a braided or ribbon tensor equivalence, accordingly.
\end{thm}

\begin{proof}
One may verify the following relations on objects and morphisms, using Mac Lane coherence theorem,
$$1_{\rho_n}\underline{\otimes} 1_{\rho_{m}}=1_{\rho_{n+m}},$$
$$S'\underline{\otimes} T'\circ S\underline{\otimes}T=(S'S)\underline{\otimes} (T'T),$$
$$(\rho_n\underline{\otimes} \rho_m)\underline{\otimes}\rho_p=\rho_{n+m+p}=\rho_n\underline{\otimes}(\rho_m\underline{\otimes}\rho_p),$$
$$\overline{\alpha}_{\rho_n, \rho_m, \rho_p}:=1_{\rho_{n+m+p}}$$ is natural,
$$S\underline{\otimes}(T\underline{\otimes} U)=(S\underline{\otimes}T)\underline{\otimes} U$$
$$\iota\text{ is a strict tensor unit for }   {\mathcal C}_\rho^s.$$
We omit the   computations.
 Thus the abstract completion of ${\mathcal C}_\rho^s$ with subobjects and direct sums is a strict tensor category.
Routine computations show that ${\mathcal C}_\rho^s$ becomes a braided tensor category in this way with $\overline{c}$, moreover the remaining statements may be easily verified using Mac Lane coherence again.

\end{proof}

\begin{ex}
Let $G$ be a finite group, $\omega$ a normalized ${\mathbb C}^{\times}$-valued cocycle and consider the pointed fusion category ${\rm Vec}^\omega_G$ discussed in more detail in Example \ref{pointed}. In this case, if we 
work with a strict realization of ${\rm Vec}$ but $\omega$ 
does not arise from a $2$-cocycle then the category is not strict, but it is skeletal,
that is there is a unique object in each isomorphism class \cite{EGNO}. When we apply the strictification to the full subcategory
generated by the tensor powers of the direct sum of the simple objects we have a strict and skeletal category, that 
does not have subobjects. Adding subobjects following Karoubi completion gives a strict tensor category that is not
skeletal in general, in agreement with remark 2.8.7 in \cite{EGNO}.

\end{ex}

\section{Weak  Hopf algebras}\label{6}
  Hopf algebras are characterised among quasi-Hopf algebras by the property
of having  trivial associator \cite{Drinfeld_quasi_hopf}. This characterization  gives   insight into the cohomological interpretation
of quasi-Hopf algebras, in that it leads to  the   notion of a $3$-coboundary associator.
 In this section we develop   a weak analogue of  the notion of Hopf algebra  among weak quasi-Hopf algebras. The corresponding special subclass   will be termed 
   weak  Hopf algebras. We shall see that
    there is no strictly coassociative weak example. We shall construct    examples later on.

 \begin{defn}\label{3-coboundary}
 Let $A$ be a weak quasi bialgebra with associator $\Phi$ and coproduct $\Delta$.
 We shall call $\Phi$ a $3$-coboundary associator if there is a twist $F\in A\otimes A$ such that
 \begin{equation}\label{coho1}\Phi=1\otimes \Delta(F^{-1})I\otimes F^{-1}F\otimes I\Delta\otimes 1(F),\end{equation}
\begin{equation}\label{coho2}\Phi^{-1}=\Delta\otimes 1(F^{-1})F^{-1}\otimes I I\otimes F  1\otimes\Delta(F).\end{equation}
 \end{defn}

If $A$ is a quasi bialgebra and $F$ is an invertible twist then only one equation suffices among (\ref{coho1}) and ({\ref{coho2}), and Def. \ref{3-coboundary} reduces to the corresponding notion of a $3$-coboundary associator. 
We next introduce weak  Hopf algebras.

Let   $A$ be an algebra with a   coproduct $\Delta$ and  a counit $\varepsilon$. 
To shorten some formulas, we set:
$$P=\Delta(I),$$
$$P_3=\Delta\otimes 1(P),\quad\quad\quad Q_3=1\otimes\Delta(P),$$
$$P_4=\Delta\otimes 1\otimes 1(P_3),\quad\quad\quad Q_4=1\otimes1\otimes\Delta(Q_3)$$

Assume that the coproduct satisfies  the following intertwining relations, espressing coassociativity in  a weak sense. For $a\in A$,
\begin{align}
&Q_3\Delta\otimes 1\circ\Delta(a)  = 1\otimes\Delta\circ\Delta(a)P_3,\label{int1} \\
& P_3 1\otimes\Delta\circ\Delta(a)=\Delta\otimes 1\circ\Delta(a) Q_3\label{int2}.
\end{align}

\begin{prop}
The element
 $\Phi:=Q_3P_3$
  satisfies  Def.  \ref{wqh} d), with partial inverse $\Phi^{-1}=P_3Q_3$ if and only if
\begin{align}
&P_3Q_3P_3=P_3, \quad Q_3P_3Q_3=Q_3,\label{basic0}\\
&Q_41\otimes\Delta\otimes1(I\otimes PP\otimes I)P_4=
Q_4\Delta\otimes\Delta(P)  P_4. \label{cocycle}
\end{align}
\end{prop}

\begin{proof} Relations (\ref{basic0}) correspond obviously to (\ref{eqn:intro2}), and  (\ref{int1}) to (\ref{eqn:intro4}). We explicit the cocycle condition (\ref{eqn:intro6}). We have $I\otimes P_3=1\otimes\Delta\otimes 1(I\otimes P)$ and $I\otimes PQ_3=Q_3$, and similarly  $Q\otimes IP_3=P_3$. This implies, taking into account
(\ref{int1}) and (\ref{int2}),
$$I\otimes\Phi1\otimes\Delta\otimes 1(\Phi)\Phi\otimes I=I\otimes Q_3P_3 1\otimes\Delta\otimes 1(Q_3P_3) Q_3P_3\otimes I=$$
$$I\otimes Q_3 1\otimes\Delta\otimes 1(Q_3P_3)P_3\otimes I=I\otimes Q_3 1\otimes\Delta\otimes 1(Q_3)1\otimes\Delta\otimes 1(P_3)P_3\otimes I=$$
$$I\otimes Q_3 1\otimes\Delta\otimes 1(1\otimes\Delta(P)1\otimes\Delta\otimes 1(\Delta\otimes1(P))P_3\otimes I=Q_4I\otimes P_3Q_3\otimes I P_4=$$
$$Q_41\otimes\Delta\otimes1(I\otimes PP\otimes I)P_4.$$
On the other hand,
$$1\otimes 1\otimes \Delta(\Phi)\Delta\otimes 1\otimes1(\Phi)=1\otimes1\otimes\Delta(Q_3P_3)\Delta\otimes1\otimes1(Q_3P_3)=$$
$$Q_4\Delta\otimes\Delta(P)P_4.$$
Finally, the normalisation condition relation (\ref{eqn:intro7}) is an immediate consequence of the counit axioms (\ref{eqn:intro1}).
\end{proof}

\begin{rem}\label{cocycle_variant} The cocycle relation (\ref{cocycle}) can   alternatively be written as 
$$Q_41\otimes\Delta\otimes1(Q_3P_3)P_4=
Q_4\Delta\otimes\Delta(P)  P_4.$$
Indeed, the computations in the last proof show that
$$Q_41\otimes\Delta\otimes1(I\otimes P)=I\otimes Q_31\otimes\Delta\otimes 1(Q_3),$$
(and a similar identity involving $P_4$ and $P_3$) hence multiplying on the left by $Q_4$, this term also equals
$Q_41\otimes\Delta\otimes 1(Q_3)$.
\end{rem}

\begin{defn}\label{wbialgebra} An algebra $A$ with coproduct $\Delta$ and counit $\varepsilon$ for which the projections $P$, $P_j$, $Q_j$, $j=3,4$, satisfy the requirements of the previous proposition is a weak quasi-bialgebra with associator $\Phi=Q_3P_3$ and will
  be called a be called a {\it weak  bialgebra}.

\end{defn}

 \begin{prop}\label{weak_strong_antipode}
If a weak  bialgebra $A$ admits an antipode $(S, \alpha,\beta)$ in the sense of weak quasi-Hopf algebras then $\alpha$, $\beta$ are invertible and $\beta=\alpha^{-1}$. Hence ${\rm        ad}(\alpha^{-1})S$  is the unique strong antipode of $A$.
\end{prop}

\begin{proof} A computation shows that if (\ref{eqn:antip1}) holds for $(S,\alpha,\beta)$ where $S$ is an antiautomorphism of $A$, then equations  (\ref{eqn:antip2}) for the associator $\Phi=Q_3P_3$ reduce to $\beta\alpha=I$ and $\alpha\beta=I$. The last statement follows from  Prop. \ref{strong_antipode_twist} a)
\end{proof}

\begin{defn}\label{wha} A weak  bialgebra with a  (unique) strong antipode,  will be called a {\it weak  Hopf algebra}.\end{defn}

\begin{rem}  The first examples of   weak   quasi-Hopf algebras  appeared in the physics literature, in the work by Mack and Schomerus \cite{MS}, who were motivated by the need of constructing a quantum analogue of a global gauge group for certain models of algebraic quantum field theories in low dimensions. They started with  a nonsemisimple category of representations of  $U_q({\mathfrak sl}_2)$   at roots of unity and indicated how to construct a such an algebra \cite{MS1, MS}. In a previous work \cite{CP}, Mack-Schomerus construction was studied in detail in the more general case of  $U_q({\mathfrak sl}_N)$, and it was shown    that these are indeed weak  Hopf algebras in the sense of this section.\end{rem}

    We next state, without proof, a few simple properties of weak  Hopf algebras (and in fact already of weak quasi-Hopf algebras)  useful 
    to construct new examples from given ones.
    
     \begin{prop}\label{constructions} Let $A$ be a weak  Hopf algebra.
     \begin{itemize}
\item[(a)] (tensor products)
    If   $B$ is another weak  Hopf algebra then the   natural weak quasi-Hopf structure on the tensor product algebra  $A\otimes B$ is  a weak  Hopf algebra structure.
    \item[(b)] (subalgebras) let   $C$ be a unital subalgebra of $A$ which is invariant under coproduct and antipode. Then $C$ is a weak  Hopf algebra with the restricted structure and there is a natural inclusion of rigid tensor categories ${\rm Rep}(A)\to{\rm Rep}(B)$.
        \item[(c)] (quotients) If $D$ is a weak Hopf algebra related to $A$ via an algebra epimorphism $A\to D$ compatible with coproduct and antipode then there is an inclusion ${\rm Rep}(C)\to{\rm Rep}(A)$ as a full rigid tensor subcategory.
        \end{itemize}
      
    \end{prop}

 \begin{prop}
 Let $A$ and $B$ be weak  Hopf algebras, and let $\alpha:A \to B$ an algebra isomorphism which intertwines the corresponding coproducts and antipodes. Then $\alpha$ is automatically an isomorphism of weak quasi-Hopf algebras.
 \end{prop}

Semisimple bialgebras are described via Tannaka-Krein duality by semisimple tensor categories endowed with a tensor functor to ${\rm Vec}$. This characterization extends  to weak bialgebras, and is based on the simple observation that they have a weak tensor furgetful functor.

\begin{thm}\label{TK_algebraic} Let ${\mathcal C}$  be a semisimple  (rigid) tensor category with finite dimensional morphism spaces and ${\mathcal F}:{\mathcal C}\to{\rm        Vec}$ a faithful weak quasi-tensor functor (taking an object and a dual to spaces with the same dimension).
Then $A={\rm        Nat}_0({\mathcal F})$ is a weak  bialgebra (weak Hopf algebra) if and only if ${\mathcal F}$ is a weak tensor functor. 
 \end{thm}

\begin{proof} Let $A={\rm        Nat}_0({\mathcal F})$  be a weak  bialgebra. The forgetful functor of   $A$ is weak tensor and this implies that the same holds for  ${\mathcal F}$ since it is monoidally isomorphic to the composition of a tensor equivalence with the forgetful  functor. Conversely, if ${\mathcal F}$ is weak  tensor then   the associator $\Phi$ of $A$   and its inverse $\Phi^{-1}$ are derived from (\ref{wt1}) and (\ref{wt2}), and a computation shows that $\Phi=1\otimes\Delta(\Delta(I))\Delta\otimes 1(\Delta(I))$, $\Phi^{-1}=\Delta\otimes 1(\Delta(I)) 1\otimes\Delta(\Delta(I))$, that is $A$ is a weak  bialgebra. For the last assertion note that the equality requirement on the dimensions of an object and a dual are automatically satisfied in our case, thanks to Cor. \ref{weak_dim}. Hence Theorem \ref{TK_algebraic_quasi} guaranties that $A$ has an antipode.

\end{proof}

  It follows that the constructions of Prop. \ref{constructions} have a   description in terms of pairs of abstract tensor categories endowed with a weak tensor functor. In particular, the following will turn out  useful to construct new weak Hopf algebras from given examples, see Sect. \ref{20}.
  
  \begin{cor}\label{quotients}
     Let ${\mathcal C}$ be a fusion category endowed with a weak tensor functor to ${\rm Vec}$.    Under Tannaka-Krein correspondence,  full fusion subcategories ${\mathcal D}\subset {\mathcal C}$ endowed with the restricted functor correspond to  quotient  weak Hopf algebras of $A={\rm        Nat}_0({\mathcal F})$.  
     
  \end{cor}

 The class of weak  Hopf is   not invariant under general twists, but we next see that it is so under   a suitable subclass of twists, that play the role of $2$-cocycles in our framework. 
  
 \begin{defn}\label{def_2-cocycle} Let $A$ be a weak  bialgebra.
A twist $F\in A\otimes A$  is called a {\it   $2$-cocycle} of $A$ if   it satisfies the following equations,
\begin{equation}\label{cocycle1}1\otimes \Delta(F^{-1})I\otimes F^{-1}F\otimes I\Delta\otimes 1(F)=Q_3P_3,\end{equation}
\begin{equation}\label{cocycle2}\Delta\otimes 1(F^{-1})F^{-1}\otimes I I\otimes F  1\otimes\Delta(F)=P_3Q_3.\end{equation}
\end{defn}
 
  Note that $P_3$ and  $P_3^F:=\Delta_F\otimes 1(FF^{-1})$ are respectively domain and range for
$F\otimes I\Delta\otimes 1(F)$,   and the partial  inverse of this element is   $\Delta\otimes 1(F^{-1})F^{-1}\otimes 1$, and similarly for $I\otimes F1\otimes\Delta(F)$.
The $2$-cocycle equations can equivalently be written in the following form  $$Q_3^FF\otimes I\Delta\otimes 1(F)=I\otimes F1\otimes\Delta(F)P_3,$$
$$P_3^F I\otimes F1\otimes\Delta(F)=F\otimes I\Delta\otimes 1(F)Q_3,$$
with $Q_3^F:=1\otimes\Delta_F(FF^{-1})$,
as well as in a form
which emphasises a categorical feature,  
$$\Phi_FF\otimes I\Delta\otimes 1(F)=I\otimes F1\otimes\Delta(F)\Phi,$$
$$\Phi_F^{-1} I\otimes F1\otimes\Delta(F)=F\otimes I\Delta\otimes 1(F)\Phi^{-1}.$$
This last form also shows that the notion of a $2$-cocycle has an extension to weak quasi-Hopf algebras which in turn extends the corresponding notion for quasi-Hopf algebras, see, e.g.,  \cite{Kassel}.

\begin{prop}\label{2-cocycle} Let $A$ be a weak quasi-bialgebra with coproduct $\Delta$ and associator $\Phi$,  and let $F\in A\otimes A$ be a twist. Then   $A_F$  is a weak  bialgebra if and only if $\Phi$ is the $3$ coboundary associator defined by $F$ as in (\ref{coho1}), (\ref{coho2}).
In particular, if $A$ is a weak  bialgebra,   $A_F$ is a weak  bialgebra as well if and only if
$F$ is a $2$-cocycle. \end{prop}

\begin{proof} 
We already know that $A_F$ is a weak quasi-bialgebra with    coproduct $\Delta_F(a)=F\Delta(a)F^{-1}$ and  associator $\Phi_F=I\otimes F 1\otimes \Delta(F)\Phi\Delta\otimes 1(F^{-1})F^{-1}\otimes I$. We have $\Phi_F^{-1}=F\otimes I\Delta\otimes 1(F)\Phi^{-1}1\otimes\Delta(F^{-1})I\otimes F^{-1}$. Hence for $A_F$ to be a   weak  bialgebra it suffices that    the associator and its inverse satisfy $\Phi_F=Q_3^FP_3^F$, $\Phi_F^{-1}=P_3^FQ_3^F$. A simple computation shows that these equations are   equivalent to the   equations in the statement. If in particular $A$ is a weak  bialgebra as well, these equations reduce to the $2$-cocycle equations (\ref{cocycle1}), (\ref{cocycle2}).
\end{proof}

\begin{rem} If $A$ is a quasi-Hopf algebra, equations (\ref{cocycle1}) and (\ref{cocycle2}) are precisely the cohomological equations  which characterise a  cohomologically trivial associator. Quite interestingly, these equations are meaningful for weak quasi-Hopf algebras   with the weak counterparts of associator and twist, with  no extra requirement on  $F$. The previous proposition  shows that   weak  bialgebras arise naturally when one tries to solve  them for a given associator $\Phi$ of a weak quasi-bialgebra $A$. 
This gives a cohomological motivation for regarding the associator of a weak  Hopf algebra as trivial. 
\end{rem}

The following corollary extends to weak  Hopf algebras a property known for
 Hopf algebras, see, e.g., \cite{Turaev}.

\begin{cor}
Let $A$ be  a weak  Hopf algebra and $F\in A\otimes A$ a $2$-cocycle.
Then the element $u_F=m\circ S\otimes 1(F^{-1})$ is invertible and $u_F^{-1}=m\circ 1\otimes S(F)$.
\end{cor}

\begin{proof}  
The twisted weak quasi-bialgebra $A_F$ is a weak  bialgebra thanks to Prop. \ref{2-cocycle}. If $S$ is the strong antipode of $A$ then $A_F$ has weak quasi-Hopf algebra antipode    $(S, \alpha_F, \beta_F)$ where   $\alpha_F=m\circ S\otimes 1(F^{-1})$, $\beta_F=m\circ 1\otimes S(F)$, by (\ref{twisted_antipode}). Hence we can apply Prop. \ref{weak_strong_antipode} to $A_F$ and deduce that $\alpha_F$ and $\beta_F$ are inverses of one another.
\end{proof}

\begin{prop}
If $F$ is a $2$-cocycle of $A$ and $G$ is a $2$-cocycle of $A_F$ then $GF$ is a $2$-cocycle of $A$.
\end{prop}

We introduce  two examples of $2$-cocycles that will be useful.
 
\begin{prop}\label{2-cocycle_deformation}
Let $v\in A$ be an invertible element with $\varepsilon(v)=1$ and $F\in A\otimes A$ a $2$-cocycle, then $F_v:=v\otimes vF\Delta(v^{-1})$ is a $2$-cocycle as well. 
\end{prop}

\begin{proof} Obviously $F_v^{-1}=\Delta(v)F^{-1}v^{-1}\otimes v^{-1}$.
A  computation shows that   the left hand side of (\ref{cocycle1}) equals  
$$1\otimes\Delta\circ\Delta(v)1\otimes\Delta(F^{-1})I\otimes F^{-1}F\otimes I\Delta\otimes 1(F)\Delta\otimes1\circ\Delta(v^{-1})=$$
$$1\otimes\Delta\circ\Delta(v)Q_3P_3\Delta\otimes1\circ\Delta(v^{-1})=$$
$$1\otimes\Delta\circ\Delta(v)1\otimes\Delta\circ\Delta(v^{-1})Q_3P_3=Q_3P_3.$$ Relation (\ref{cocycle2}) for $F_v$ can be proved in  a similar way.
  
\end{proof}

\begin{prop}\label{trivial_cocycle}
Let $E\in A\otimes A$ be an idempotent satisfying 
$$\varepsilon\otimes 1(E)=1\otimes\varepsilon(E)=I,$$
$$EP_2E=E,\quad\quad P_2EP_2=P_2.$$
 Then $F=EP_2 $ defines a trivial twist with $D(F)=P_2$, $R(F)=E$ and $F^{-1}=P_2E.$ It is a   $2$-cocycle
 if and only if the following additional  relations hold,
 $$Q_3 1\otimes\Delta(E)I\otimes E E\otimes  1\Delta\otimes 1(E)P_3=Q_3P_3,$$
$$P_3\Delta\otimes I(E) E\otimes I I\otimes E 1\otimes\Delta(E) Q_3=P_3Q_3.$$
\end{prop}

We omit the proof as it follows from a simple computation.

\section{Quasitriangular and    ribbon structures}\label{7}

The notion of quasitriangular  Hopf algebra was introduced by Drinfeld in \cite{Drinfeld_qg} and extended to the quasi-Hopf algebra case in
\cite{Drinfeld_quasi_hopf}.
In this section we introduce  and study  quasitriangular structure for   weak quasi-Hopf algebras.   We shall then restrict 
to weak quasi-Hopf algebras with a strong antipode   and introduce 
 the notion of ribbon structure in this case. In particular,   we   develop the basic   results for this special subclass.
For some results for which computational difficulties would arise, we further restrict to 
 the   special subclass of weak  Hopf algebras. In this case, we are able to present arguments   extending   the corresponding results for Hopf algebras. We conclude the section 
 explaining how  later on we shall extend all  the results of this section concerning weak  Hopf algebras to weak quasi-Hopf algebras with a strong antipode. This extension  will be useful for the forthcoming developments of the paper of Sect. \ref{18}  and for our applications
 of Sect. \ref{20}, \ref{KW}.

 \medskip

With any weak quasi-bialgebra  $A$,   we associate the {\it opposite algebra} $A^{{\rm        op}}$
 with data given by
 \begin{equation}
 {\varepsilon}^{{\rm        op}}=\varepsilon,\quad\quad  {\Delta}^{{\rm        op}}(a):=\sigma\circ\Delta(a),\quad\quad  {\Phi}^{{\rm        op}}:=\Phi_{321}^{-1},\label{opp}
 \end{equation}
 where $\sigma$ is the transposition automorphism of $A\otimes A$ and $\Phi_{321}^{-1}$ understood in a partial sense. Note that $A^{{\rm        op}}$ is a weak  bialgebra if so is $A$.

\begin{defn}\label{quasi_triangular_structure}
A quasitriangular structure on $A$, also referred to $R$-matrix axioms, is defined by a partially  invertible element $R\in A\otimes A$,  ($R\in M(A\otimes A)$ if $A$ is discrete)
satisfying the following properties,  
\begin{equation}\label{qtriangular1}D(R)=\Delta(I),\quad\quad R(R)=\Delta^{{\rm        \rm op}}(I)\end{equation}
\begin{equation}\label{qtriangular2}\Delta^{{\rm        \rm op}}(a)=R\Delta(a)R^{-1},\end{equation}
\begin{equation}\label{qtriangular3}\Delta\otimes 1(R)=\Phi_{312}R_{13}\Phi_{132}^{-1}R_{23}\Phi_{123},\end{equation}
\begin{equation}\label{qtriangular4}1\otimes\Delta(R)=\Phi_{231}^{-1}R_{13}\Phi_{213}R_{12}\Phi_{123}^{-1},\end{equation}
\end{defn}
We follow the standard notation: for a simple tensor  $a=a_1\otimes\dots \otimes a_n\in A^{\otimes n}$ and a permutation $i\in{\mathbb P}_n$, $a_{i_1\dots i_n}$ is the simple tensor having $a_j$ in the $i_j$-th component. If $a\in A^{\otimes k}$ with $k<n$ then we apply this definition to $a$ tensored on the right with $n-k$ copies of the identity operator.
  Furthermore
relations (\ref{qtriangular1})--(\ref{qtriangular4}) imply the analogue of the Yang-Baxter relation, which, taking into account (\ref{qtriangular3}) and (\ref{qtriangular4}), can be written in the  following form
\begin{equation}\label{YB}\Phi_{321}^{-1} =I\otimes R 1\otimes\Delta(R)\Phi\Delta\otimes 1(R^{-1}) R^{-1}\otimes I.\end{equation}
Relations (\ref{qtriangular1}), (\ref{qtriangular2}),   (\ref{YB}), and the following property (\ref{prop1})   express the
twist relation
\begin{equation} \label{R-twist} A^{{\rm        op}}=A_R.\end{equation}

Given a $^*$-algebra $A$  endowed with the structure of a weak quasi-bialgebra,
   we can form another weak quasi-bialgebra $\tilde{A}$, the {\it adjoint algebra}
 with the same algebra structure but counit,   coproduct, and associator given by
\begin{equation}
\tilde{\varepsilon}(a):=\overline{\varepsilon(a^*)},\quad\quad \tilde{\Delta}(a):=\Delta(a^*)^*,\quad\quad \tilde{\Phi}:={\Phi^*}^{-1}.\label{tilde}
\end{equation}
  Note that if $B$ is a $^*$-algebra, and $p$ and $q$ are idempotents of $B$ and if 
 $T\in(p, q)$ then $T^*\in (q^*, p^*)$. Hence if $T$  is partially invertible in $(p, q)$, so is $T^*$ in $(q^*, p^*)$. We understand   ${\Phi^*}^{-1}$ in this way. It will be useful to observe that  
\begin{prop}\label{R-canonicity} If $R$ is an $R$-matrix for $A$ then  
\begin{itemize}
\item[{\rm        a)}]  
  {\rm        $R^{{\rm        \rm op}}:=R_{21}$ is an $R$-matrix for $A^{{\rm        \rm op}}$}, 
\item[{\rm        b)}]  {\rm        if $A$ is a $^*$-algebra, $\tilde{R}:={R^{*}}^{-1}$ is an $R$-matrix for $\tilde{A}$}, 
\item[{\rm        c)}] {\rm        if $F\in A\otimes A$ is a twist, $R_F:=F_{21}RF^{-1}$ is an $R$-matrix for $A_F$}, 
\item[{\rm        d)}] {\rm        $R_{21}^{-1}$ is  another $R$-matrix for $A$}.
 \end{itemize}\end{prop}\medskip

 \begin{defn} By a {\it quasitriangular weak  bialgebra} we understand a weak  bialgebra endowed with a quasitriangular structure   as a weak quasi-bialgebra.  \end{defn}

Note that any    $R$-matrix of a weak  Hopf algebra is a $2$-cocycle  by (\ref{R-twist}).
An important property for representation theory of quasitriangular Hopf algebras is that the square of the antipode is an inner automorphism. This was   shown by Drinfeld who explicitly constructed an implementing invertible element $u\in A$ for Hopf algebras \cite{Drinfeld_cocommutative}. Furthermore, Reshetikhin and Turaev introduced the notion of ribbon Hopf algebra
  \cite{RT}. We next show that these developments have   extensions  to weak  Hopf algebras, although the computations in the proofs become more involved. We start with the following remark giving a simplification of the axioms in the weak  Hopf algebra case.

\begin{prop}\label{qtriangular_simplified}
Equations (\ref{qtriangular3}) and (\ref{qtriangular4}) for a weak  Hopf algebra are equivalent to
\begin{equation}\label{qtriangular5}\Delta\otimes 1(R)=\Phi_{312}R_{13} R_{23}\Phi_{123},\end{equation}
\begin{equation}\label{qtriangular6}1\otimes\Delta(R)=\Phi_{231}^{-1}R_{13}R_{12}\Phi_{123}^{-1}.
\end{equation}

\end{prop}

\begin{proof}
We   prove (\ref{qtriangular5}). We have $\Phi_{123}=1\otimes\Delta(P)\Delta\otimes 1(P)$, $\Phi^{-1}=\Delta\otimes 1(P)1\otimes\Delta(P)$,   $\Phi_{312}=\Delta\otimes 1(P') a_2\otimes b\otimes a_1$, and $\Phi^{-1}_{132}=a_1\otimes b\otimes a_2 1\otimes\Delta^{{\rm        \rm op}}(P)$ where   $P=\Delta(I)$, $P'=\Delta^{{\rm        \rm op}}(I)$, and  we have used the notation $\Delta(b)=b_1\otimes b_2$ and $P=a\otimes b$. By (\ref{qtriangular2}) we have
$$R_{13}\Phi^{-1}_{132}R_{23}=R_{13}a_1\otimes b\otimes a_2 1\otimes\Delta^{{\rm        \rm op}}(P)R_{23}=a_2\otimes b\otimes a_1R_{13} R_{23}1\otimes\Delta(P)$$
and the conclusion follows. For (\ref{qtriangular6}) we similarly have $\Phi^{-1}_{231}=1\otimes\Delta(P')b_2\otimes a\otimes b_1$ and $\Phi_{213}=b_1\otimes a\otimes b_2\Delta^{{\rm        \rm op}}\otimes 1(P)$.
\end{proof}

We   give a definition of 
 ribbon weak quasi-Hopf algebra $A$ with a strong antipode extending the corresponding notion
for Hopf algebras due to  Reshetikhin and Turaev \cite{RT2}.

\begin{defn}\label{balanced_ribbon_wqh} Let $A$ be a (discrete) weak quasi-bialgebra  Then $A$ is called
  {\it balanced} if it  is   quasitriangular  and is endowed 
  with  an invertible central element $v\in A$ ($v\in M(A)$)  such that
  \begin{equation}\label{squared_matrix} R_{21}R=v\otimes v\Delta(v^{-1}),
  \end{equation}
  where $R$ is the $R$-matrix.
If in addition $A$ has
an   antipode $(S, \alpha, \beta)$ such that
$S(v)=v$, then $A$ is called a {\it ribbon weak quasi-Hopf algebra}, and $v$   the {\it ribbon element}.
A  {\it balanced  (ribbon) weak  bialgebra} is a weak  bialgebra (weak  Hopf algebra) is defined in the natural way.
\end{defn}

Note that the definition does not depend on the choice of the antipode by Prop. \ref{unique_antipode}.
We next introduce Drinfeld element $u$. For simplicity, we restrict to the case of a weak quasi-Hopf algebra with strong antipode. This will suffice for our applications.

\begin{defn}\label{Drinfeld_element_u} Let $A$ be a quasitriangular weak quasi-Hopf  algebra with strong antipode $S$ and $R$-matrix $R$. The element  
\begin{equation} u=\sum_i S(t_i) r_i\end{equation}
where $R=\sum_i r_i\otimes t_i$ is called {\it Drinfeld element}. We also set $R^{-1}=\sum_j\overline{r}_j\otimes\overline{t}_j$.\end{defn}

\begin{prop}\label{inner_antipode}
Let $A$ be a quasitriangular weak quasi-Hopf  algebra with strong antipode $S$ and $u$ the associated Drinfeld element. Then $u$ is invertible, 
$u^{-1}=\sum_j S^{-1}(\overline{t}_j) \overline{r}_j$ and 
\begin{equation}\label{inner_antipode2}S^2(x)=uxu^{-1}, \quad\quad x\in A.\end{equation}

\end{prop}

\begin{proof}
This proof is a  
generalisation of the corresponding proof for quasitriangular Hopf algebras, see e.g. \cite{Kassel}. In the following computations we use the notation $\Delta(x)=x_1\otimes x_2$ for $x\in A$, 
 $\Delta(I)=a\otimes b$, $R=r\otimes t$, $\Phi^{-1}=x'\otimes y'\otimes z'$.
We have $$\Delta^{{\rm        \rm op}}\otimes 1(\Delta(x))R\otimes I\Phi^{-1}=R\otimes I\Phi^{-1}1\otimes\Delta(\Delta(x))$$ 
that  accordingly may be written as
 $$x_{1,2}rx'\otimes x_{1,1}ty'\otimes x_2z'=rx'x_{1}\otimes ty'x_{2,1}\otimes z'x_{2,2}.$$
  Applying  $1\otimes S\otimes S^2$ and multiplying from right to left gives by (\ref{eqn:antip1}), (\ref{eqn:intro1}),
\begin{equation}
\label{square} S^2(x)w= wx,\quad\quad w:=S^2(z')S(y')ux'.
 \end{equation}
The $3$-cocycle relation $\Phi^{-1}\otimes I=\Delta\otimes 1\otimes 1(\Phi^{-1})1\otimes 1\otimes\Delta(\Phi^{-1})I\otimes\Phi1\otimes\Delta\otimes 1(\Phi)$ leads to $w=u$. The last argument extends
in a straightforward way the case of quasi-Hopf algebras, see the proof of Lemma 2.4 in \cite{Bulacu2003_qtqh}.
 The      formula for $u^{-1}$ follows from Cor. \ref{strong_antipode_twist} b).
 \end{proof}
 
 Note that this proposition does not depend on the $R$-matrix properties (\ref{qtriangular3}), (\ref{qtriangular4}).
 But when they do hold, we obtain stronger relations for $u$ in a way that extends the corresponding relations for quasitriangular Hopf algebras.
 The following extends    Lemma 2.1.1, Ch. XI, of \cite{Turaev}, or Theorem VIII.2.4 of \cite{Kassel} to weak  Hopf algebras.

\begin{prop}\label{squared_antipode}
If $A$ is a quasitriangular weak quasi-bialgebra algebra defined by $R$ then
\begin{equation}\label{prop1} \varepsilon\otimes 1(R)=I, \quad\quad 1\otimes \varepsilon(R)=I.\end{equation}
If $A$ is a weak  Hopf algebra,
\begin{equation}\label{prop2} S\otimes S(R)=f_{21}Rf^{-1},\end{equation}
where $f$ is the element defined in Prop. \ref{anticom}.
\end{prop}

\begin{proof}
The proof of (\ref{prop1}) goes as in the bialgebra case, it suffices to   apply $\varepsilon\otimes1 \otimes 1$ 
and $1\otimes 1\otimes \varepsilon$ to  (\ref{qtriangular3})  and (\ref{qtriangular4}) respectively. To show (\ref{prop2}) we tensor   both sides of   (\ref{qtriangular5}) by the identity operator $I$ on the left and multiply  by  $I\otimes 1\otimes\Delta(P)\Delta\otimes\Delta(P)$ on the right and obtain
\begin{equation}\label{basic}1\otimes\Delta\otimes 1(I\otimes R) 1\otimes 1\otimes\Delta(I\otimes P)\Delta\otimes\Delta(P)=XR_{34}\end{equation}
 where
$$X=I\otimes \Phi_{312} R_{24} 1\otimes1\otimes\sigma[1\otimes 1\otimes\Delta (I\otimes P)\Delta\otimes\Delta(P)],$$
   $\sigma: A\otimes A\to A\otimes A$ is the flip automorphism
and we have used the intertwining relations (\ref{qtriangular1}), (\ref{qtriangular2}).
We next recall from the first section the map
$V(a\otimes b\otimes c\otimes d)= S(b)c\otimes S(a)d$
  that we wish to apply to both sides of (\ref{basic}) and we   obtain
  \begin{equation}\label{first} f=V(X)R.\end{equation} To show the claim we
  perform computations taking into account the following facts: 
  a)  one of the two ways the element $f$ is defined for a weak quasi-Hopf algebra
with strong antipode is
$f=V(I\otimes\Phi^{-1}1\otimes1\otimes\Delta(\Phi))$.  For a weak  Hopf algebra we have $$I\otimes\Phi^{-1}1\otimes1\otimes\Delta(\Phi)=I\otimes \Delta\otimes 1(P) 1\otimes 1\otimes\Delta(1\otimes\Delta(P))\Delta\otimes\Delta(P)=$$
$$1\otimes\Delta\otimes 1(1\otimes\Delta(P))I\otimes 1\otimes\Delta(P)\Delta\otimes\Delta(P).$$
b) We have  $V(1\otimes\Delta\otimes 1(Z)Y)=V(Y)$ as soon as $m\circ S\otimes\varepsilon\otimes 1(Z)=I$, where $m: A\otimes A\to A$ is the multiplication map. This holds in particular for $Z= I\otimes R$ and $Z=1\otimes\Delta(P)$, by (\ref{prop1}) and (\ref{eqn:antip1}). Hence the image of the left hand side of (\ref{basic}) under $V$ is $f$.  c) $V(XR_{34})=V(X)R$.  We next apply a similar procedure to relation (\ref{basic}) for the opposite weak  Hopf algebra
 getting the relation
\begin{equation}\label{opposite}1\otimes\Delta^{{\rm        \rm op}}\otimes 1(I\otimes R_{21}) 1\otimes 1\otimes\Delta^{{\rm        \rm op}}(I\otimes P_{21})\Delta^{{\rm        \rm op}}\otimes\Delta^{{\rm        \rm op}}(P_{21})=X^{{\rm        \rm op}}R_{43}\end{equation}
where  
$X^{{\rm        \rm op}}=I\otimes \Phi^{-1}_{213} R_{42}1\otimes 1\otimes\Delta(I\otimes P_{21})\Delta^{{\rm        \rm op}}\otimes\Delta (P_{21})$
but now we apply the map $W:=\sigma\circ S\otimes S\circ V^{{\rm        \rm op}}$ to both sides of (\ref{opposite}), where
   $V^{{\rm        \rm op}}$ acts as $V$ but with  $S^{-1}$ in place of $S$. To perform these computations we remark that: d) for the left hand side we use the identity $S\otimes S\circ V^{{\rm        \rm op}}=V\circ\tau$, where $\tau$ is the automorphism of $A^{\otimes 4}$ taking $a_1\otimes a_2\otimes a_3\otimes a_4\to a_4\otimes a_3\otimes a_2\otimes a_1$. e) the image of the left hand side of (\ref{opposite}) under $\tau$ is 
   $$1\otimes\Delta\otimes 1(R\otimes I)\Delta\otimes 1\otimes 1(P\otimes I)\Delta\otimes\Delta(P)$$
   f)   the second way in which $f$ can be computed is $f=V(\Phi\otimes I\Delta\otimes 1\otimes1 (\Phi^{-1}))$,
   and recall that this was due to the $3$-cocycle relation of $\Phi$ and the previous remark b). For a weak  Hopf algebra, computations similar to those in a) give
   $$\Phi\otimes I\Delta\otimes 1\otimes1 (\Phi^{-1})=1\otimes\Delta\otimes 1(\Delta\otimes 1(P))\Delta\otimes 1\otimes 1(P\otimes I) \Delta\otimes\Delta(P).$$
   Hence using b) again, the image of the left hand side of (\ref{opposite}) under $W$ is $f_{21}$. For the right hand side, we write $W$ in the form $W=V\circ \sigma\otimes\sigma\circ\tau$. Simple computations show that if $\alpha=\sigma\otimes\sigma\circ\tau$ then $\alpha(R_{43})=R_{21}$ and that $V(YR_{21})=S\otimes S(R)V(Y)$ for $Y\in A^{\otimes 4}$.
   Summarizing, the image of (\ref{opposite}) under $W$ is
\begin{equation}\label{second} f_{21}=S\otimes S(R)V(\alpha(X^{{\rm        \rm op}})).\end{equation} Comparing
(\ref{first}) and (\ref{second}), the proof of (\ref{prop2}) will be complete provided $V(X)=V(\alpha(X^{{\rm        \rm op}}))$.
To show this,  a computation relying on
by (\ref{qtriangular1}), (\ref{qtriangular2}), (\ref{eqn:intro4}) gives
$$I\otimes\Phi_{312} R_{24} 1\otimes 1\otimes\sigma[1\otimes\Delta\otimes 1(1\otimes\Delta(P))]=1\otimes\Delta\otimes 1(Z) I\otimes\Phi_{312}R_{24}$$
where $Z=1\otimes\Delta^{{\rm        \rm op}}(P)$.
It follows,  by a), and the $3$-cocycle relation, and (\ref{qtriangular1}), (\ref{qtriangular2}) again,
$$1\otimes\Delta\otimes 1(Z) X=I\otimes\Phi_{312} R_{24} 1\otimes1\otimes\sigma[I\otimes\Phi^{-1}1\otimes 1\otimes\Delta(\Phi)]=$$
$$ I\otimes\Phi_{312} R_{24} 1\otimes1\otimes\sigma[1\otimes\Delta\otimes 1(\Phi)\Phi\otimes I\Delta\otimes1\otimes1(\Phi^{-1})]=$$
$$(I\otimes\Phi1\otimes\Delta\otimes 1(\Phi))_{1423}R_{24}1\otimes1\otimes\sigma[\Phi\otimes I\Delta\otimes1\otimes1(\Phi^{-1})].$$
 On the other hand,
$$\alpha(X^{{\rm        \rm op}})=\Phi^{-1}_{142}R_{24} 1\otimes1\otimes\sigma[\Delta\otimes 1(P)\otimes I\Delta\otimes\Delta(P)]$$
and similar computations give
$$a_{1,1}\otimes a_2\otimes b\otimes a_{1,2}\alpha(X^{{\rm        \rm op}})=\Phi^{-1}_{142}R_{24}1\otimes1\otimes\sigma[\Phi\otimes I\Delta\otimes 1\otimes1 (\Phi^{-1})].$$
 Hence 
$$1\otimes\Delta\otimes 1(Z) X=(I\otimes\Phi1\otimes\Delta\otimes 1(\Phi)\Phi\otimes I)_{1423}a_{1,1}\otimes a_2\otimes b\otimes a_{1,2}\alpha(X^{{\rm        \rm op}})=$$
    $$(1\otimes1\otimes\Delta(\Phi)\Delta\otimes1\otimes1(\Phi))_{1423} \alpha(X^{{\rm        \rm op}})$$
by the $3$-cocycle relation again. It now suffices to apply $V$ on both sides of this identity.
 \end{proof}
 
 \begin{prop}\label{ribbon_element} Drinfeld element $u$ of a quasitriangular weak  Hopf algebra satisfies
 $$R_{21}R\Delta(u)=\Delta(u)R_{21}R=f^{-1}S\otimes S(f_{21}) u\otimes u.$$
 \end{prop}
 
 \begin{proof}  
 The first equality follows easily from (\ref{qtriangular2}). 
 We show the second equality. The left hand side
 equals, by  Prop. \ref{anticom},
 $$\Delta(u)R_{21}R=\Delta(S(t))R_{21}R\Delta(r)=f^{-1}S\otimes S(\Delta^{{\rm        \rm op}}(t))fR_{21}R\Delta(r).$$
where notation  is as  before: $R=r\otimes t$,  $\Delta(I)=P=a\otimes b$,   
 $\Delta(x)=x_1\otimes x_2$. We are thus reduced to show the equality
\begin{equation}\label{toshow}S\otimes S(\Delta^{{\rm        \rm op}}(t))fR_{21}R\Delta(r)=S\otimes S(f_{21}) u\otimes u.\end{equation}
We denote by $\lambda$ and $\rho$ the left and right hand sides of (\ref{toshow}), respectively.
We use again the   map $V: A^{\otimes 4}\to A^{\otimes 2}$,  $V(a\otimes b\otimes c\otimes d)=S(b)c\otimes S(a)d$,  and recall that  $f=V(A)=V(X)$, where we have set $A=I\otimes\Phi^{-1} 1\otimes 1\otimes\Delta(\Phi)$ and
$X=\Phi\otimes I\Delta\otimes1\otimes 1(\Phi^{-1})$. We shall also need the property 
\begin{equation}V(La\otimes b\otimes c\otimes d)=S(b)\otimes S(a)V(L)c\otimes d.\label{property}\end{equation} For example, it shows that  
  $$\lambda=V(A\cdot{}\Delta(t)\otimes[R_{21}R\Delta(r)]).$$
  Furthermore, assuming that $Y\in A^{\otimes 4}$ satisfies $V(Y)=u\otimes u$, and writing $X=x\otimes y\otimes w\otimes z$,
  we have 
  $$\rho=S\otimes S((S(y)w\otimes S(x)z)_{21})u\otimes u=S(z)S^2(x)u\otimes S(w)S^2(y)u=$$
  $$S(z)ux\otimes S(w)uy=V(Yw\otimes z\otimes x\otimes y)=V(YX_{3412}).$$
  We start computing $\lambda$. By (\ref{qtriangular5}), and $(\ref{qtriangular6})$,
  $$R_{21}R\Delta(r)\otimes\Delta (t)=R_{21}R\Delta\otimes 1\otimes 1(\Phi^{-1}_{231}R_{13}R_{12}\Phi^{-1}_{123})=$$ 
 $$\Delta\otimes1\otimes 1(\Phi^{-1}_{231})R_{21}\Delta^{{\rm        \rm op}}\otimes 1\otimes 1(R_{13})R\Delta\otimes 1\otimes1(R_{12})\Delta\otimes1\otimes1(\Phi^{-1}_{123})=$$
 $$\Delta\otimes1\otimes 1(\Phi^{-1}_{231})R_{21}\Phi_{421}R_{24}R_{14}\Phi_{214}R\Phi_{312}R_{13}R_{23}\Phi_{123}\Delta\otimes 1\otimes 1(\Phi^{-1}_{123}).$$ After applying the   permutation of $(13)(24)\in{\mathbb P}_4$, and taking into account
 $$A  1\otimes1\otimes\Delta(\Phi^{-1})=I\otimes\Phi^{-1} 1\otimes1\otimes\Delta(1\otimes\Delta(P))=1\otimes\Delta\otimes 1(1\otimes\Delta(P))I\otimes\Phi^{-1}$$
 we see that $\lambda$ equals
 $$V(1\otimes\Delta\otimes 1(1\otimes\Delta(P))[I\otimes\Phi^{-1} R_{43}\Phi_{243}R_{42}R_{32}\Phi_{432}][R_{34}\Phi_{134}R_{31}R_{41}\Phi_{341}]1\otimes 1\otimes\Delta(\Phi^{-1}_{312})).$$
 The first bracketed element is the shift to the right of 
 $\Phi^{-1} R_{32}\Phi_{132}R_{31}R_{21}\Phi_{321}$, and computations similar to those of Prop. \ref{qtriangular_simplified} show that the latter equals $(1\otimes\Delta(R))_{312}R_{21} b\otimes a_2\otimes a_1$.
 Similarly, the second bracketed element acts as identity on the second factor, and as $a\otimes b_2\otimes b_1R_{23}R_{21}\Phi^{-1}_{213}R_{31}\Phi_{231}$ in the remaining factors.  This in turn equals
 $a\otimes b_2\otimes b_1(1\otimes\Delta(R))_{213}R_{31}\Phi_{231}$ by Prop. \ref{qtriangular_simplified} again and property (\ref{basic0}).
 Hence $\lambda$ equals
 $$V(1\otimes\Delta\otimes 1(1\otimes\Delta(P))(1\otimes\Delta(R))_{423}R_{32} I\otimes b\otimes a_2\otimes a_1\cdot a\otimes I\otimes  b_2\otimes b_1(1\otimes\Delta(R))_{314}R_{41}X_{3412}).$$
 Taking the range of $X$ into account, we are finally left to show that  
 $$Y:=1\otimes\Delta\otimes 1(1\otimes\Delta(P))(1\otimes\Delta(R))_{423}R_{32} I\otimes b\otimes a_2\otimes a_1\cdot
 $$
 $$ a\otimes I\otimes  b_2\otimes b_1(1\otimes\Delta(R))_{314}R_{41}(1\otimes\Delta\otimes 1(\Delta\otimes 1(P))_{3412}.$$
 indeed satisfies $V(Y)=u\otimes u$. To this aim, we move the two idempotents at both extremes towards the center
 using the commutation relations (\ref{qtriangular2}) and  (\ref{int2}) and the domain relations $(\ref{int1})$, we find 
 $$Y=(1\otimes\Delta(R))_{423}R_{32} I\otimes b_2\otimes b_1\otimes a(1\otimes\Delta\otimes 1(1\otimes\Delta(P))_{1432}\cdot$$
$$(1\otimes\Delta\otimes 1(\Delta\otimes 1(P))_{1432} a_1\otimes I\otimes b\otimes a_2(1\otimes\Delta(R))_{314}R_{41}
=$$
$$(1\otimes\Delta(R))_{423}R_{32}(1\otimes\Delta\otimes 1(\Phi))_{1432}(1\otimes\Delta(R))_{314}R_{41}.$$
 Now $1\otimes\Delta\otimes 1(1\otimes\Delta(P))(1\otimes\Delta(R))_{423}R_{32}$ and 
$(1\otimes\Delta(R))_{314}R_{41} (1\otimes\Delta\otimes 1(\Delta\otimes 1(P))_{3412}$
have $(Q_3)_{1432}$ and $(P_3)_{1432}$
as domain and range respectively, and furthermore 
$Q_31\otimes\Delta\otimes 1(\Phi) P_3=Q_3\Delta\otimes\Delta(P)P_3$
by the cocycle relation, see Remark \ref{cocycle_variant}.
Hence we can also write
$$Y=(1\otimes\Delta(R))_{423}R_{32}(\Delta\otimes\Delta(P))_{1432}(1\otimes\Delta(R))_{314}R_{41}.$$
We are now able to compute $V(Y)$ by means of an iterative use of (\ref{property}): $$V((1\otimes\Delta(R))_{423})=I\otimes I,$$ 
$$V((1\otimes\Delta(R))_{423}R_{32})=u\otimes I,$$
$$V((1\otimes\Delta(R))_{423}R_{32}(\Delta\otimes\Delta(P))_{1432})=S(b_2)\otimes S(a_1)\cdot u\otimes I\cdot b_1\otimes a_2=$$
$$S(b_2)ub_1\otimes S(a_1)a_2=S(S(b_1)b_2)u\otimes\varepsilon(a)=$$
$$\varepsilon(b)\otimes\varepsilon(a) u\otimes I=u\otimes I,$$
$$V((1\otimes\Delta(R))_{423}R_{32}(\Delta\otimes\Delta(P))_{1432}(1\otimes\Delta(R))_{314})=I\otimes S(t_1)\cdot u\otimes I\cdot r\otimes t_2=$$
$$u\otimes I 1\otimes\varepsilon(R)=u\otimes I,
$$
$$V(Y)=u\otimes u,$$
and the proof is complete.
 
     \end{proof}
 
 At the level of representation theory, the previous proposition establishes commutativity of the following diagram.
 \[
\begin{tikzcd}
\rho\otimes\sigma  \arrow{rr}{u_\rho\otimes u_\sigma} \arrow{d}{\varepsilon(\rho,\sigma)}  &&  \rho^{\vee\vee}\otimes\sigma^{\vee\vee} \arrow{d} \\
\sigma\otimes\rho  \arrow{r}{\varepsilon(\sigma,\rho)} &  \rho\otimes\sigma \arrow{r}{u_{\rho\otimes\sigma}} &   {(\rho\otimes\sigma)^{\vee\vee} }
\end{tikzcd}
\]

We shall go back to this in Sect. \ref{17}.

\begin{rem}\label{central_square_root_of_squared_Drinfeld_element}
It follows from (\ref{squared_matrix}) and (\ref{prop1}) 
that if $v$ makes $A$ balanced then $\varepsilon(v)=1$. Furthermore when $v$ is a ribbon element, applying $m\circ S\otimes 1$ to both sides of  (\ref{squared_matrix}) 
and taking into account
  (\ref{inner_antipode2})  also gives   $v^2=uS(u)$, with $u$ Drinfeld element, as in Def. \ref{Drinfeld_element_u}.
  \end{rem}

\begin{cor}\label{spherical_structure_wqh}
The elements $u$ and $v$ of  a ribbon weak  Hopf algebra satisfy
$$\Delta(uv^{-1})=f^{-1} S\otimes S(f_{21}) uv^{-1}\otimes uv^{-1}.$$
\end{cor}

Thus $\omega=uv^{-1}$ is a pivot making $A$ into a pivotal weak  Hopf algebra in the sense of Def. \ref{pivotal_wqh}.

\begin{rem}\label{charmed_element} It is shown in \cite{RT2, Panaite2} that 
when $A$ is a quasi-triangular Hopf algebra  then $v\to\omega=uv^{-1}$ is a bijective
correspondence between   ribbon elements
  and invertible elements $\omega\in A$ satisfying a) $S^2(a)=\omega a\omega^{-1}$, b) $S(\omega)=\omega^{-1}$, and c) $S(u)=\omega^{-1} u \omega^{-1}$ and d) $\Delta(\omega)=\omega\otimes\omega$ ($\omega$ is group-like). Such elements $\omega$   are called {\it charmed}.
  \end{rem}

\begin{rem}\label{AC}
  Altschuler and Coste extended ribbon structures to    quasi-Hopf algebras  \cite{Altschuler_Coste}, stated analogues of the  lemmas of this section and outlined some of the proofs.   Complete proofs   have been given in \cite{Hausser_Nill_qtqh, Bulacu2003_qtqh}.
In this passage,
   the construction of Drinfeld element $u$ 
   and the notion of ribbon   quasi-Hopf algebra needs to be suitaby modified.  
Moreover,   the proof of the analogue of Prop. \ref{squared_antipode}, Prop. \ref{ribbon_element}  become  more involved.      These works together with the results of this section lead to extensions of the main properties of ribbon structures to the more general setting of weak quasi-Hopf algebras. Furthermore, the correspondence described in Remark \ref{charmed_element} extends as well to quasi-triangular weak quasi-Hopf algebras replacing the group-like condition d) with the pivot condition $\Delta(\omega) = f^{-1} S \otimes S(f_{21})\omega \otimes \omega$.

However,  we shall refrain from doing this, and rather take an alternative categorical approach. More in detail,
 motivated also by the study of quantum dimension,  in Sect. \ref{17}
  we shall   revisit Drinfeld isomorphism and ribbon structures  in the framework of tensor categories. Moreover, we shall study more general structures (coboundary symmetries). In particular, it will follow from the results of that section  that via Tannaka-Krein duality when $A$ is a discrete weak quasi-Hopf algebra with a strong antipode then Drinfeld element is still defined as in Def. \ref{Drinfeld_element_u}. It will also follow that 
   all the special results    of this section concerning weak  Hopf algebras extend to this setting with the same statements, and this will suffice for the forthcoming developments of our paper considered in Sect. \ref{18}, and for   our applications of Sects. \ref{20}, \ref{KW}.\end{rem}

\medskip

\section{$\Omega$-involution and $C^*$-structure}\label{8}
 In this section we introduce a $^*$-involution $^*: A\to A$ to a weak quasi-Hopf algebra. In the usual approach, among the compatibility conditions with the weak quasi-Hopf algebra stucture, one requires for example that the involution
 and the coproduct   commute. We shall relax these compatibility conditions via the introduction of a twist $\Omega$ which is part of the axioms of the involutive structure.
There are several reasons to study such structures. On one hand, unlike the ordinary approach, the more general notion   is   invariant under Drinfeld  twist operation $A\to A_F$.  Another motivation for us  arises   from considering natural examples, which include the   Drinfeld-Jimbo quantum groups $U_q({\mathfrak g})$ for the values of the deformation parameter $q$  with $|q|=1$. 
 Finally, 
as we shall see more precisely in Sect. \ref{10}, $\Omega$-involutions of weak quasi-Hopf algebras describe unitary structures in fusion categories and intervene in the study of tensor $^*$-equivalences.

 \begin{defn}\label{twisted_involution} A weak quasi bialgebra $A$ will be called {\it $\Omega$-involutive} if it is endowed with a $^*$-involution $^*: A\to A$ making it into a $^*$-algebra and a selfadjoint twist $\Omega\in A\otimes A$ such that 
 $\tilde{A}=A_\Omega$, with $\tilde{A}$ the adjoint weak quasi bialgebra defined in (\ref{tilde}). Explicitly, this means that 
 $\Omega\in A\otimes A$ is a  partially invertible element satisfying
     \begin{equation}
 \Omega^*=\Omega,\label{sa}
 \end{equation}
 \begin{equation}
D(\Omega)=\Delta(I), \quad\quad\quad R(\Omega)=\Delta(I)^*,\label{domain} 
\end{equation}
 \begin{equation}
\Delta(a^*)=\Omega^{-1} \Delta(a)^*\Omega,\quad\quad\quad a\in A,\label{intertwiner}
\end{equation}
 \begin{equation}
\varepsilon\otimes1(\Omega)=I=1\otimes\varepsilon(\Omega)\label{normalized}
\end{equation}
 \begin{equation}
{\Phi^{*}}^{-1}=I\otimes \Omega 1\otimes\Delta(\Omega)\Phi\Delta\otimes1(\Omega^{-1})\Omega^{-1}\otimes I 
 \label{eqn:omega0}
\end{equation}

  A {\it unitary weak quasi bialgebra} is an $\Omega$-involutive weak quasi bialgebra such that      $A$ is a $C^*$-algebra and $\Omega$
is   positive in $A\otimes A$. Note that $\Omega^{-1}$ is also automatically positive.
Corresponding Hopf versions assume the existence of an antipode $S$. Note that in general we require     no   compatibility assumption with the involution.\end{defn}

The most important relations are the intertwining property with the coproduct (\ref{intertwiner}) and the compatibility
relation (\ref{eqn:omega0}) between $(^*, \Omega)$ and
  the associator.

The notion of   $\Omega$-involution for a semisimple weak quasi-Hopf  algebra is the most general involutive structure that gives rise to a tensor $^*$-category structure on the category of finite dimensional representations   of $A$.
For example, we shall see that every fusion tensor $^*$-category ($C^*$-category) arises from a semisimple $\Omega$-involutive
(unitary) weak quasi-Hopf algebra. We next recall   several well known and important special notions.

\begin{rem}\label{trivial_unitary_structure}
{\bf  a)} {\it $A$ is a Hopf $^*$-bialgebra} precisely when $\Delta(I)=I\otimes I$ and $\Omega=I\otimes I$, $\Phi=I\otimes I\otimes I$. These structures are widely studied when $A$ is a $C^*$-algebra in the operator algebraic approach to quantum groups see e.g. \cite{CQGRC}, \cite{Timmermann}.
{\bf  b)} When $A$ is a  bialgebra ($\Delta(I)=I\otimes I$, $\Phi=I\otimes I\otimes I$) (\ref{sa}) and (\ref{domain}) say that $\Omega$ is a selfadjoint invertible element.
Note that in this case  (\ref{eqn:omega0}) says that $\Omega$ is a $2$-cocycle in the usual sense for Hopf algebras.
In the next proposition we   discuss an extension of this property to weak bialgebras.
{\bf c)} If $A$ is as in b)  and $A$ is a $C^*$-algebra with $\Omega$  positive
 then the twisted algebra $A_F$, 
  with $F=\Omega^{1/2}$ is a quasi $C^*$-bialgebra in the   sense of a).   We shall shortly consider an extension of  the notion of triviality of $\Omega$ in the weak quasi bialgebras which is the algebraic counterpart of  the notion of unitary weak quasi-tensor functor of Def. \ref{unitarity}.
 {\bf d)} When $A$ is a quasi-bialgebra (that is $\Delta(I)=I\otimes I$ and $\Phi$ non-trivial) we recover the notion  introduced by Gould and Lekatsas   \cite{quasi_star}.
 \end{rem}
 
 \begin{ex}
 The Hopf algebras $U_q({\mathfrak g})$ for $|q|=1$ considered by Wenzl in \cite{Wenzl} are for us important examples of   $\Omega$-involutive Hopf algebras
 with a non-trivial  selfadjoint  $2$-cocycle $\Omega$  in the sense of part b) of the previous remark.
 We shall
 discuss these examples in Section \ref{20}. 
In this case, $\Omega$ is canonically induced by the $R$-matrix\footnote{To be precise, the $\Omega$-involution of $U_q({\mathfrak g})$ is not comprised in Def. \ref{twisted_involution}. This is due to the fact that
  the $R$-matrix lies in a suitable topological completion of $U_q({\mathfrak g})\otimes U_q({\mathfrak g})$. However, when we consider the category of
 finite dimensional representations of $U_q({\mathfrak g})$, this inconvenience is not source of complications
 in that it gives rise to a braided tensor category, as   explained in  
 \cite{Sawin}.
It follows that  the associated $\Omega$ also lies in the   completed algebra.
Similarly to the $R$-matrix case, in this paper we will consider applications of the notion of $\Omega$-involution to categories of
 finite dimensional representations, see Sect. \ref{9}--\ref{VOAnets2}, and  we shall refrain from giving a more general definition of $\Omega$-involution.}. Furthermore, in Sections \ref{18}, \ref{19}, \ref{20} we shall construct new examples
 of semisimple $\Omega$-involutive or unitary  weak  Hopf algebras associated to $U_q({\mathfrak g})$ for $q$ a suitable root of unity,
 corresponding  to the associated unitary fusion categories.
 
 \end{ex}
 
 We next extend the $2$-cocycle property of $\Omega$ from bialgebras to weak  bialgebras.

 \begin{prop}\label{2-cocycle2} If $(*, \Omega)$ makes a weak  bialgebra $(A,\Delta, \varepsilon, \Phi=Q_3P_3)$  $\Omega$-involutive       then
$\Omega$ is a $2$-cocycle.
\end{prop}
\begin{proof}
By definition $\Phi=Q_3P_3$ is an associator with $\Phi^{-1}=P_3Q_3$, see
  Sect. \ref{6}. Then $(A,\tilde{\Delta}, \varepsilon, \tilde{\Phi})$ is a weak  Hopf algebra as well since
 $$\tilde{\Phi}={(Q_3P_3)^{*}}^{-1}={(Q_3P_3)^{-1}}^{*}=(P_3Q_3)^*=Q_3^*P_3^*=1\otimes\tilde{\Delta}(\tilde{\Delta}(I))\tilde{\Delta}\otimes 1(\tilde{\Delta}(I))$$ and similarly $\tilde{\Phi}^{-1}=\tilde{\Delta}\otimes 1(\tilde{\Delta}(I))1\otimes\tilde{\Delta}(\tilde{\Delta}(I))$. By (\ref{eqn:omega0}) and Prop. \ref{2-cocycle}
 we see that $\Omega$ is a $2$-cocycle.
\end{proof}

\begin{defn}\label{discrete_involutive} 

Let $A$ be a discrete algebra in the sense of Def. \ref{discrete}. A positive $^*$-involution on $A$ is a
  $^*$-involution such that 
 $A$ can be completed to a $C^*$-algebra. We may then identify $A$ with an algebraic direct sum of matrix subalgebras with the usual $^*$-involution.  An $\Omega$-involutive structure on $A$
is defined as in the unital case but $\Omega$ is here allowed to be a (selfadjoint) element  in $M(A\otimes A)$. 
   A {\it unitary discrete  weak quasi-Hopf algebra} is defined by further requiring that $\Omega$ has positive  components in the matrix subalgebras. It follows that the same holds for $\Omega^{-1}$. In the particular case where $A$ is a weak  Hopf algebra, we shall refer to $A$ as a {\it unitary discrete   weak  Hopf algebra}. 
    \end{defn}

Unless otherwise stated, involutions of discrete algebras will be assumed positive. This  will hold for   most part of this paper.   
We next describe the $\Omega$-involutions on a simple class of discrete algebras.

  \begin{ex}\label{pointed_Hermitian_qh}
  We consider the   bialgebra  $C_\omega(G)$ of complex valued functions of a finite group $G$ with the usual coproduct
  and associator given by a ${\mathbb T}$-valued $3$-cocycle $\omega$,     see Ex. \ref{pointed}. A natural unitary structure is given by the $C^*$-structure of $C_\omega(G)$ and $\Omega=I$. More generally,   a general
 $\Omega$-involution
  for   $C_\omega(G)$
  over the same   $C^*$-algebra is given by a  normalized $2$-cocycle $\Omega(g, h)$ with values in ${\mathbb R}^\times$, that is a function satisfying $\Omega(1, g)=\Omega(g, 1)=1$ and $\Omega(g, h)\Omega(gh, k)=\Omega(h, k)\Omega(g, hk)$ for all $g$, $h$, $k\in G$. The corresponding quasi-Hopf algebra is unitary if and only of
  $\Omega(g, h)>0$ for all $g$, $h\in G$. This is not always the case, 
  an example is given by   $G={\mathbb Z}_2$  
  $\Omega(g, g)=-1$ with  $g$   the group generator.   \end{ex}

  }\medskip
  
  In the next sections we shall see   examples of unitary discrete weak quasi-Hopf algebras    arising from unitary tensor categories,  Sect. \ref{10}, and quantum groups Sects. \ref{18}, \ref{19}, \ref{20}. Moreover we  shall discuss  
  conditions which guarantee  unitarity, see Theorem \ref{Positivity}.
 In the following proposition we show that   the fact that $\tilde{A}$ and $A_\Omega$ have the same counit is a redundant assumption.

\begin{prop}\label{counit}
The counit $\varepsilon$ of a weak quasi-bialgebra $A$ is unique. 
If $A$ is a weak quasi-Hopf algebra with antipode $S$ the counit satisfies $\varepsilon\circ S=\varepsilon$. 
If $A$ is an $\Omega$-involutive weak quasi-bialgebra  then  $\varepsilon(a^*)=\overline{\varepsilon(a)}$, for every $a\in A$. \end{prop}
\begin{proof}
The first two statements can  be proved in the same way as for quasi-bialgebras, namely the first  follows from (\ref{eqn:intro1}) while
the second  from applying the counit to one of the equations (\ref{eqn:antip1}).
For the last statement it suffices to  show that
$\widetilde{\varepsilon}(a):=\overline{\varepsilon(a^*)}$ is a counit.  For example,
$$
(1\otimes\widetilde{\varepsilon})(\Delta(a))=a_{(1)}\widetilde{\varepsilon}(a_{(2)})= $$
$$
(a_{(1)}^*\varepsilon(a_{(2)}^*))^*=(1\otimes\varepsilon(\Delta(a)^*))^*= $$
$$(1\otimes\varepsilon(\Omega\Delta(a^*)\Omega^{-1}))^*=(1\otimes\varepsilon(\Delta(a^*)))^*=a.$$

 \end{proof}

\begin{prop}\label{twisted_omega} \begin{itemize}
\item[{\rm        a)}] 
Let $A$ be an $\Omega$-involutive weak quasi-bialgebra and  $F\in A\otimes A$  a twist
(Def. \ref{definition_of_twist}). Then $A_F$ is an $\Omega_F$-involutive weak quasi-bialgebra with the same involution as $A$ and 
\begin{equation}\label{Omega_F}
\Omega_F:={F^{-1}}^*\Omega F^{-1}, \quad\quad(\Omega_F)^{-1}:=F\Omega^{-1}F^*
\end{equation}
\item[{\rm        b)}]  If $A$ is a discrete pre-$C^*$-algebra and $\Omega$ is positive in $M(A\otimes A)$ then  $\Omega_F$ is positive as well.
\end{itemize}
\end{prop}

We discuss a useful application of the twist of the unitary structure.

\begin{defn}\label{trivial_involution} Let $A$ be a weak quasi bialgebra with a $^*$-involution. An $\Omega$-involution 
compatible with $^*$ on $A$ is called {\it trivial} if it is given by
$\Omega=\Delta(I)^*\Delta(I)$  and  $\Omega^{-1}=\Delta(I)\Delta(I)^*$. Thus $\Omega$ is a trivial twist.
We shall call it {\it strongly trivial} if in addition $\Delta(I)$ is selfadjoint, that is equivalent to require that
 commutes $\Delta$ commutes with  the $^*$-involution as in the usual   $^*$-bialgebra theory.
  In this case,    $\Delta(I)$ is a selfadjoint projection.
 \end{defn}

With a strongly trivial involution, $\Delta$ commutes with $^*$ and the associator $\Phi$ satisfies $\Phi^*=\Phi^{-1}$.  The above   notions
 of (strong) triviality has the same motivation as that of and are related to those of (strongly) unitary weak quasi tensor functor discussed before Def. \ref{unitarity}.
 
 \begin{rem}\label{twist_strongly_trivial}
a) As in the case of weak quasitensor structures,   when $A$  is a weak quasi bialgebra with a $^*$-involution and
a trivial $\Omega$-involution 
compatible with $^*$ then
 $T=\Delta(I)$ is a twist with left inverse $T^{-1}=\Delta(I)\Delta(I)^*$ 
 (or $T'=\Delta(I)^*\Delta(I)$ with ${T'}^{-1}=\Delta(I)$)
 giving a new wqh $A_T$ ($A_{T'}$) with strongly trivial involution.
 b) When $A$ is a discrete unitary weak quasi-bialgebra with a trivial $\Omega$-involution then this involution is automatically strongly trivial. This follows
   from the fact that we are in a $C^*$-setting, Prop. \ref{strongly_unitary_prop} and the following Tannaka-Krein duality, Theorem \ref{TheoremTannakaStar}.
 \end{rem}

 \begin{ex} We have the following generalization 
of  the construction in Remark c) in \ref{trivial_unitary_structure}. Let $A$ be a unitary discrete weak quasi bialgebra with an   
$\Omega$-involution given by $\Omega\in M(A\otimes A)$.
We may consider   $T=\Omega^{1/2}$ defined via continuous functional calculus in each full matrix subalgebra of $M(A\otimes A)$. This element
satisfies the properties $T\Delta(I)=T$, $\Delta(I)^*T=T$, and $\omega\otimes 1(T)=1\otimes\omega(T)=1$, so we may  regard $T$ as an element   of $A$  with the same domain $\Delta(I)$  as $\Omega$.
Applying   the same construction to $\Omega^{-1}$,
we construct   $T'=(\Omega^{-1})^{1/2}\in M(A\otimes A)$ with range $\Delta(I)$.
\end{ex}

 \medskip

\begin{cor}\label{making_omega_trivial}
Let $A$ be a discrete unitary   weak quasi bialgebra 
defined by $\Omega$ and assume that $(\Omega^{-1})^{1/2}\Omega^{1/2}=\Delta(I)$.
 Let us regard  $T=\Omega^{1/2}$ as a twist with left inverse $T^{-1}=(\Omega^{-1})^{1/2}$.   Then  the twisted $\Omega$-involution of $A_{T}$ is  trivial, and therefore strongly trivial.
\end{cor}

\begin{proof} By
 part b) of Prop. \ref{twisted_omega},
$\Omega_T=\Delta_T(I)^*\Delta_T(I)$ and  $\Omega_T^{-1}=\Delta_T(I)\Delta_T(I)^*$. Strong triviality
follows again from the fact that we are in a $C^*$-setting, Prop. \ref{strongly_unitary_prop} and Tannaka-Krein duality
Theorem \ref{TheoremTannakaStar}.
 \end{proof}

 We shall refer to $A_{\Omega^{1/2}}$ as the {\it unitarization} of $A$.
 We next introduce a deformation of an $\Omega$-involution on a given weak  bialgebra   that may be thought of as analogous to the  twist operation for the weak quasi bialgebra structure.

\begin{defn}\label{star_twist} Let $A$ be an
  $\Omega$-involutive weak quasi bialgebra $A$ defined by ($^*$, $\Omega$). A twist for the involutive structure is an invertible selfadjoint
  $t\in A$ such that $\varepsilon(t)=1$. If $A$ is discrete in the sense of Def. \ref{discrete_involutive} we allow $t\in M(A)$.    \end{defn}

  \begin{prop}\label{1_twist}
  A twist $t$ of an  involution $(^*, \Omega)$ gives rise to another involutive structure on the same weak quasi  bialgebra  
   via
  $$a^\dag:=t^{-1}a^* t,\quad\quad \Omega_t:=t^{-1}\otimes t^{-1}\Omega\Delta(t).$$ If $A$  is a $C^*$-algebra under
  $^*$,  or else if $A$ is discrete,   and $(^*, \Omega)$ is a positive involution,    then the same holds for $A$ with respect to $(^\dag, \Omega_t)$ for any positive twist $t$.
  \end{prop}

\begin{proof}
The proof of the first statement follows from routine computations. For example,  
$\Omega_t^\dag=\Omega_t$ follows from (\ref{intertwiner}). We show the second statement.
 If $\|a\|$ denotes a  $C^*$-norm of $A$ 
compatible with $^*$ then $\|a\|_t:=\|t^{1/2} a t^{-1/2}\|$ is another   $C^*$-norm on $A$ compatible with 
$^\dag$. (Note that the original and the deformed norms are equivalent, hence completeness
of one is equivalent to completeness of the other.) Furthermore if $\Omega$ is positive with respect to the original involution,  the element 
$\Xi:=t^{-1/2}\otimes t^{-1/2}\Omega^{1/2}\Delta(t^{1/2})$ satisfies $\Xi^\dag \Xi=\Omega_t$, so $\Omega_t$ is positive with the $^\dag$-involution of $A\otimes A$.

\end{proof}

In the discrete case, any other involution making $A$ into a pre-$C^*$-algebra is of the kind $a^\dag=t^{-1} a^* t$, with $t$   determined up to a normalized central positive element of $M(A)$. This implies e the following useful result. 
 
 \begin{cor}
If a discrete weak quasi bialgebra $A$ can be made unitary   with respect to an assigned pre-$C^*$-algebra involution   of $A$, the same is true for any other  such involution. 
 \end{cor}

As for twists of bialgebra structures, twists of involutive structures admit a categorical interpretation, that will be discussed
in Prop. \ref{1_twist_equivalence}.
The next   results exploit the relations between antipode and  $\Omega$-involution. 
 
 \begin{prop}\label{involution_antipode} Let $(S,\alpha,\beta)$ be an antipode of an $\Omega$-involutive weak quasi-Hopf algebra.
 There is an invertible $\omega\in A$ such that
\begin{equation}\label{first_antipode}
S(a)=\omega S^{-1}(a^*)^*\omega^{-1},\quad a\in A,\end{equation}
\begin{equation}\label{second_antipode}S^{-1}(\beta)^*=\omega^{-1}\alpha_\Omega, \quad S^{-1}(\alpha)^*=\beta_\Omega \omega\end{equation}
  uniquely determined by (\ref{first_antipode}) and one of (\ref{second_antipode}).
In particular when $S$ is a strong antipode then
 \begin{equation}\label{little_omega}
 \omega=m(S\otimes 1(\Omega^{-1})), \quad\quad\omega^{-1}=m(1\otimes S(\Omega)).
 \end{equation}
 \end{prop}
\begin{proof}
The adjoint weak quasi-bialgebra $\tilde{A}$ defined in (\ref{tilde}) has antipode $(\tilde{S}, \tilde{\alpha},\tilde{\beta})$ with 
$\tilde{S}(a):=S^{-1}(a^*)^*$, $\tilde{\alpha}:=S^{-1}(\beta)^*$, $\tilde{\beta}:=S^{-1}(\alpha)^*$. On the other hand, $\tilde{A}=A_\Omega$, and therefore it also admits $(S_\Omega,\alpha_\Omega,\beta_\Omega)$ as an antipode by Prop. \ref{twisted_wqh}. The first statement follows from   Prop.  \ref{unique_antipode} and the last  
from a computation and (\ref{twisted_antipode}). \end{proof}

\begin{cor}\label{Kac}
The following are equivalent for an antipode $(S, \alpha, \beta)$,
\begin{itemize}
\item[{\rm        a)}] $S$ commutes with $^*$,
\item[{\rm        b)}] $S^{-1}$ commutes with $^*$,
\item[{\rm        c)}] $S^2(a)=\omega a\omega^{-1}$, $a\in A$. 
 \end{itemize}
If these conditions hold then $\omega^*\omega$ and $S(\omega)\omega$ are central.
 \end{cor}

 We   study the dependence of the element $\omega$ introduced in Prop.
 \ref{involution_antipode} on twisting.  
 
 \begin{prop}\label{invariance_little_omega_under_twisting}
 Let $A$ be an $\Omega$-involutive weak quasi-Hopf algebra with antipode $(S, \alpha, \beta)$ and involutive structure
 $(*, \Omega)$ and associated element $\omega$ as in Prop.
 \ref{involution_antipode}.
 \begin{itemize}
\item[{\rm        a)}] Let $({\rm Ad}(u)S, u\alpha, \beta u^{-1})$ be another antipode of $A$. The corresponding element is given by $\omega_u=u\omega S^{-1}(u)^*$.
\item[{\rm        b)}]
 Let $F\in A\otimes A$ be  a twist and consider the   weak quasi-Hopf algebra $A_F$ with antipode
 $(S, \alpha_F, \beta_F)$ and involutive structure $(*, \Omega_F)$.
 Then the corresponding element   is given by $\omega_F=\omega.$
 \end{itemize}

 \end{prop}

 \begin{proof} a) follows from a computation. b) By the uniqueness stated in Prop. \ref{involution_antipode} we only need to verify that $S^{-1}(\beta_F)^*=\omega^{-1}(\alpha_F)_{\Omega_F}$. 
The claim follows in a straightforward way from a computation based  on (\ref{first_antipode}) and the first relation in
(\ref{second_antipode}) which
takes into account the definition of $\alpha_F$, $\beta_F$ in
(\ref{twisted_antipode}) and of $\Omega_F$ in (\ref{Omega_F}).

 \end{proof}

\begin{defn}\label{Kac_type_definition}
An $\Omega$-involutive weak quasi-Hopf  algebra is called of {\it Kac type} if it admits a (unique) strong antipode satisfying  one of the equivalent conditions stated in Cor. \ref{Kac}. We shall also refer to the antipode as being of Kac type.
\end{defn}

The definition is motivated by the fact that  if $A$ is in turn  a Hopf $^*$--algebra in the usual sense ($\Omega=I$) then $\omega=I$, and Cor. \ref{Kac} reduces to the well known characterisation of Hopf $^*$--algebras of Kac type.

\begin{prop}\label{Kac_type_sufficiency}
Let $A$ be a  Hopf algebra   such that
$$\Delta^{{\rm        \rm op}}(a)^*=\Delta(a^*), \quad a\in A.$$ 
(e.g. $A$ is $\Omega$-involutive and satisfies $\Delta^{{\rm        \rm op}}(a)=\Omega\Delta(a)\Omega^{-1}$ for $a\in A$).
Then $A$ is of Kac type.
\end{prop}

\begin{proof}  
Since $A$ is a Hopf algebra, it admits a unique strong antipode, denoted $S$. Furthermore,      our assumptions  imply  $\Delta(a^*)=\Delta^{{\rm        \rm op}}(a)^*$ for $a\in A$. It follows that the antiautomorphism
 $\tilde{S}(a):=S(a^*)^*$  is another Hopf algebra antipode of $A$, as
 
\begin{equation}
(m\circ(1\otimes \tilde{S})\circ\Delta)(a)=a_{(1)}\tilde{S}(a_{(2)})= (S(a_{(2)}^*)a_{(1)}^*)^*=
\end{equation}
\begin{equation}
[m\circ(S\otimes1)(\Delta^{\text op}(a)^*)]^*=[m\circ(S\otimes1)(\Delta(a^*))]^*=
\end{equation}
\begin{equation}
(\varepsilon(a^*)\unit)^*=\varepsilon(a)\unit.
\end{equation}
 
Hence $\tilde{S}=S$ by uniqueness.
 
\end{proof}

Wenzl shows in  \cite{Wenzl} that  the assumptions of Prop. 
\ref{Kac_type_sufficiency} are  satisfied by the quantum groups $U_q({\mathfrak g})$ for $|q|=1$, cf. also Sect. \ref{20}.
We shall extend Prop. \ref{Kac_type_sufficiency} to  weak  Hopf algebras endowed with a $^*$-involution
and a strong antipode in Sect. \ref{18}, see Prop. \ref{antipode_commuting_with_involution_wqh_case}.

 \section{The categories ${\rm        Rep}_h(A)$ and ${\rm        Rep}^+(A)$ }\label{9}
 
 Let $A$ be a complex $^*$-algebra. In this section we associate with $A$
 the linear  $^*$-category ${\rm        Rep}_h(A)$ with objects representations on non-degenerate Hermitian spaces.
 
 If $A$ has further the structure of an $\Omega$-involutive weak quasi-Hopf algebra, then we introduce in
 ${\rm        Rep}_h(A)$ the structure of a rigid tensor $^*$-category.   
 
 Most importantly, the subclass of unitary weak quasi-Hopf algebras
 leads to    rigid tensor $C^*$-categories ${\rm        Rep}^+(A)$.\bigskip

 The basic notion is that of   Hermitian space, that is a finite dimensional vector space $V$    equipped with a non degenerate Hermitian form $(\xi, \eta)$. If $W$ is another such space, any linear map $T:V\to W$ admits an adjoint 
    $T^*:W\to V$
    defined as in the more familiar case of Hilbert spaces: $(T\xi, \eta)=(\xi, T^*\eta)$. The category ${\rm        Herm}$ with objects finite dimensional Hermitian spaces and morphisms linear maps between them is   the simplest example of a   $^*$-category.

\begin{defn} Let $A$ be a unital complex $^*$-algebra (a discrete complex $^*$-algebra resp.)
with involution $^*: A\to A$. 

  \begin{itemize}
\item[{\rm        a)}] 
A $^*$--representation is a unital
(nondegenerate resp.) representation $\rho$ of $A$ on a nondegenerate Hermitian space   $V_\rho$    satisfying
  $\rho(a^*)=\rho(a)^*$ for $a\in A$. 
\item[{\rm        b)}] A $C^*$-representation of $A$ is a $^*$-representation on a Hilbert space.
  \end{itemize}
  \end{defn}

 The study of $^*$-representations on Hermitian spaces is   motivated by   $U_q({\mathfrak g})$, for $|q|=1$  \cite{Wenzl}. In this case, Wenzl showed that for generic values of $q$, or for certain roots of unity of sufficiently high order there is a  finite set of irreducible $C^*$-representations  \cite{Wenzl}. In the latter case representation theory is not semisimple. A brief review and connections with the theory of representations of weak quasi-Hopf algebras will be studied in later sections.

      \begin{prop} Let $A$ be   a  complex $^*$-algebra either unital or  a discrete.
  Let  ${\rm        Rep}_h(A)$ be the category  with objects nondegenerate  $^*$--representations of $A$ on nondegenerate  Hermitian spaces. 
If $T\in(\rho, \sigma)$ is a morphism of ${\rm        Rep}_h(A)$,   the adjoint map $T^*: V_\sigma\to V_\rho$   is still a morphism
  of ${\rm        Rep}_h(A)$. In this way ${\rm        Rep}_h(A)$ becomes a linear $^*$-category.
  \end{prop}

An  isometric   morphism $S\in(\rho, \sigma)$  between  two $^*$-representations is a morphism satisfying $S^*S=1$. Similarly, a unitary
 is an invertible isometry, that is a morphism $U\in(\rho, \sigma)$  satisfying $U^*U=1$, $UU^*=1$. Therefore there is  a natural notion of unitary equivalence between $^*$--representations $\rho$ and $\sigma$.  Unitary equivalence implies equivalence, but, unlike the case of Hilbert space $^*$--representations, the converse does not hold. In other words, a representation can be made into a $^*$-representation in more than one way, up to unitary equivalence.  
 This can be seen with the following simple construction.
  
     Given a $^*$-representation $\rho$, let $\rho_-$ denote the  $^*$-representation with the same space and action as $\rho$ but with  
     with the opposite Hermitian form:
   $(\xi, \eta)_{V_{\rho_-}}=-(\xi, \eta)_{V_\rho}$. We shall refer to $\rho_-$ as the {\it opposite $^*$-representation}. 
       Note that  $ \rho$ and ${\rho_-}$ are equivalent 
   as representations but     they are not unitarily equivalent in the following two cases, either $\rho$ is irreducible, or it may reduce but it is a $C^*$-representation. Indeed, 
   given another $^*$-representation $\sigma$  
   and a linear map $T:V_\rho\to V_\sigma$ with adjoint $T^*$ with respect to the original forms, the adjoint of $T$ as a map   $V_{\rho_-}\to V_\sigma$ or 
   $V_\rho\to V_{\sigma_-}$ is $-T^*$. 
Thus the unitarity condition for an intertwiner $T: V_\rho\to V_{\rho_-}$ becomes $T^*T=-I$, with $T^*$ the adjoint of $T$ as a map $V_\rho\to V_\rho$,  and this is incompatible   with either irreducibility ($T$ acts as a scalar) or the $C^*$-assumption on $\rho$.     
     
  A $^*$-representation $\sigma$ is called an orthogonal direct sum of  $\rho$ and $\tau$ if there are isometries $S\in(\rho, \sigma)$, $T\in(\tau, \sigma)$ such that $SS^*+TT^*=1$. This implies that $SV_\rho$ and $TV_\tau$ are spanning, orthogonal  subspaces of $V_\sigma$: $(SV_\rho, TV_\tau)=0$, and hence are    complementary by nondegeneracy of the form.  We write $\sigma=\rho\oplus\tau$ and refer to $\rho$ and $\tau$  as orthogonal summands of $\sigma$.
 If $\rho$ and $\tau$ are $^*$--representations, the direct sum Hermitian form on   $V_\rho\oplus V_\tau$ makes this space into  a $^*$-representation $\sigma$ in the natural way and we have  $\sigma=\rho\oplus\tau$ via the   inclusions $S: V_\rho\to V_\sigma$, $T: V_\tau\to V_\sigma$. Any other realisation of $\sigma$ as a direct sum of $\rho$ and $\sigma$ will be unitarily equivalent to this.

 If $A$ is not semisimple as an algebra,  representations  may admit invariant submodules which are not summands. The following  proposition shows that the   $^*$--structure   is useful to distinguish between summands and submodules.

\begin{prop} \label{orthogonal} Let $A$ be   a  complex $^*$-algebra either unital or  a discrete.
 If $S\in(\rho, \sigma)$ is an isometry in ${\rm        Rep}_h(A)$, then $E=SS^*$ is a selfadjoint idempotent with     range $SV_\rho$, defining an orthogonal summand of $\sigma$.   Conversely, every submodule $W$ of $V_\sigma$ (i.e.  a subspace of $V_\sigma$ invariant under all the $\sigma(a)$, $a\in A$) for which the restricted Hermitian form is nondegenerate, is a $^*$--representation and an orthogonal summand.

\end{prop}  

\begin{proof}   
In general,  if the restriction of the Hermitian form of $V_\sigma$ is nondegenerate on a submodule $W$ then the adjoint of the restriction of an element $\sigma(a)$ with respect to the restricted form equals the restriction of $\sigma(a^*)$ by $^*$-invariance of $\sigma$ and nondegeneracy. 
Hence $W$ defines a $^*$-representation and 
  the inclusion map $S: W\to V_\sigma$ is an isometry.

Given an isometry $S\in(\rho, \sigma)$ in ${\rm        Rep}_h(A)$, $E=SS^*$ obviously defines an algebraic summand of $\sigma$. The ranges of $E$ and $1-E$ are orthogonal subspaces of $V_\sigma$.
This implies that 
  the Hermitian form of $V_\sigma$ is nondegenerate on either subspace and therefore these are $^*$--representations 
  $\rho$ and $\tau$ such that $\sigma=\rho\oplus\tau$.

  \end{proof}
  
  We next give a criterion for nondegeneracy of Hermitian forms.
  
  \begin{prop}\label{nondegeneracy} Let $A$ be   a  complex $^*$-algebra either unital or  a discrete.
  A nonzero Hermitian form on the vector space of an irreducible representation $\rho$ of $A$ making it $^*$-invariant is nondegenerate. Any other $^*$-representation structure on $\rho$ is unitarily equivalent to $\rho$ or $\rho_-$.  \end{prop}

  \begin{proof}
The subspace $V_\rho^\perp=\{v\in V_\rho, (v, V_\rho)=0\}$ is a   submodule by $^*$-invariance of $\rho$, and it must be proper, hence $V_\rho^\perp=0$ by irreducibility, and this shows nondegeneracy. Every other nondegenerate Hermitian form on $V_\rho$ is defined by an invertible   map
 $B: V_\rho\to V_\rho$ via $(\xi, \eta)_B=(\xi, B\eta)$, with $B$  selfadjoint with respect to the given Hermitian form.
 The adjoint of a map $T: V_\rho\to V_\rho$ with respect to the new form as compared to the old changes to $B^{-1}T^*B$. 
 The $^*$-invariance condition for $\rho$ with respect to the new form reads as $B^{-1}\rho(a^*)B=\rho(a^*)$  for $a\in A$ by $^*$-invariance of $\rho$. Thus $B$ is a nonzero real  scalar.

  \end{proof}

             A tensor product of Hermitian spaces becomes an Hermitian space in the natural way:
 $(\xi\otimes\xi', \eta\otimes\eta')_p:=(\xi,\eta)(\xi', \eta')$. In this way     ${\rm        Herm}$ becomes    a tensor $^*$-category, and it is the unique $^*$-structure on ${\rm        Herm}$ compatible with the tensor structure.
 
   We next describe how to obtain a tensor $^*$-category from an $\Omega$-involutive weak quasi-bialgebra.
Note that the $^*$-structure obtained restricting that of ${\rm        Herm}$ to ${\rm        Rep}_h(A)$
  is not the correct one, as it does not make a tensor product   of two $^*$--representations into a $^*$-representation. This is due to the fact that the    coproduct and $^*$-involution do not commute. 
 On the other hand, 
 because of the twisted commutation relation they satisfy, 
 one can   consider    a twist of the product     form  by the action of $\Omega$, 
  $$(\zeta, \zeta')_\Omega:=(\zeta,\Omega\zeta')_p, \quad \zeta,\zeta'\in V_{\rho\underline{\otimes}\rho',}$$
which is   indeed a non degenerate and Hermitian form.

 \begin{thm}\label{Hermitian_category} Let $A$ be an $\Omega$-involutive  weak quasi bialgebra. For every pair of $^*$--representations $\rho$, $\rho'$, the form $(\cdot,\cdot)_\Omega$ on $V_\rho\underline{\otimes} V_\rho'$ makes $\rho\underline{\otimes}\rho'$ into a $^*$--representation. 
 In this way ${\rm        Rep}_h(A)$ becomes a tensor $^*$--category with unitary associativity morphisms. This category is strict if $A$ is a bialgebra. 
 \end{thm}
\begin{proof}
Let $V$   be  a   Hermitian space, and consider a new  Hermitian form of $V$ defined by a given  selfadjoint invertible $A\in{\mathcal L}(V)$.  Denote by $V_A$ the associated Hermitian space.   Let $W$, $B$, be another such pair.
Given $T\in{\mathcal L}(V, W)$, we denote by $T^*$ and $T^\dag$  the adjoint of $T$ with respect
  to the new  forms (that is as a map $T: V_A\to W_B$) and the original form respectively.
They are related by
  $T^*=A^{-1} T^\dag B.$
Therefore given $T\in{\mathcal L}(V_\rho\underline{\otimes} V_{\rho'}, V_\sigma\underline{\otimes} V_{\sigma'})$,  we have 
 $T^*=\Omega^{-1} T^\dag\Omega$  with adjoints referred to  the twisted  form and the restricted product form respectively. 
Thus    $T^*=T^\dag$ if $T^\dag$ commutes with the action of $\Omega$. For example, this always holds for $T=S\underline{\otimes} S'$, with   $S\in(\rho, {\sigma})$, $S'\in({\rho'},  {\sigma'})$. Indeed,   $T^\dag=S^*\underline{\otimes} S'^*$,  and $S^*$ and $S'^*$ are intertwiners. We at once find
  $(S\underline{\otimes} S')^*=S^*\underline{\otimes} S'^*$.
   Notice   that the product form is related to the involution of the tensor product $^*$--algebra $A\otimes A$:
      $$\rho\otimes\rho'(b^*)=\rho\otimes\rho'(b)^\dag,\quad\quad b\in A\otimes A.$$
Therefore for $a\in A$,
  $$\rho\underline{\otimes}\rho'(a)^*=\rho\otimes\rho'(\Delta(a))^*=$$
  $$\rho{\otimes}\rho'(\Omega^{-1})\rho\otimes\rho'(\Delta(a))^\dag\rho{\otimes}\rho'(\Omega)=\rho\otimes\rho'(\Omega^{-1}\Delta(a)^*\Omega)=$$
  $$\rho\otimes\rho'(\Delta(a^*))=\rho\underline{\otimes}\rho'(a^*).$$
   Given   $^*$-representations $\rho$, $\sigma$, $\tau$, the   $^*$--representations $(\rho\underline{\otimes}\sigma)\underline{\otimes}\tau$ and $\rho\underline{\otimes}(\sigma\underline{\otimes}\tau)$
  act  via the morphisms $\Delta\otimes 1\circ\Delta$ and $1\otimes \Delta\circ\Delta$, respectively,
 on the       subspaces of $V_\rho\otimes\ V_\sigma\otimes V_\tau$ determined  by the image of $I$ under those morphisms. With respect to the triple product form, the associated Hermitian forms are induced by $\Omega\otimes I\Delta\otimes 1(\Omega)$ and $I\otimes\Omega 1\otimes \Delta(\Omega)$, respectively. 
 To show that the associativity morphisms $\alpha_{\rho,\sigma,\tau}$ are unitary arrows of ${\rm        Rep}_h(A)$, we  compute their adjoints taking into account the remark at the beginning of the proof,
$$\alpha_{\rho,\sigma,\tau}^*=(\Omega\otimes I\Delta\otimes 1(\Omega))^{-1}\alpha_{\rho,\sigma,\tau}^\dag I\otimes\Omega 1\otimes \Delta(\Omega)=$$
$$\rho\otimes\sigma\otimes\tau(\Delta\otimes 1(\Omega^{-1})\Omega^{-1}\otimes I \Phi^* I\otimes\Omega 1\otimes \Delta(\Omega))=\rho\otimes\sigma\otimes\tau(\Phi^{-1})=$$
$$\alpha_{\rho,\sigma,\tau}^{-1}.$$
  If in addition $A$ is a   bialgebra then $\Phi$ is the trivial associator, hence   $\Omega$ is a $2$-cocycle by Prop. \ref{2-cocycle2}. This means that 
  $(\rho\underline{\otimes}\sigma)\underline{\otimes}\tau$ and $\rho\underline{\otimes}(\sigma\underline{\otimes}\tau)$  also coincide as $^*$--representations. Since   the associativity morphisms are trivial, ${\rm        Rep}_h(A)$ is strict.

 \end{proof}

\begin{cor}\label{CorollaryRep^+(A)} Suppose that $A$ is a unitary     weak quasi bialgebra $A$. Then the full subcategory ${\rm        Rep}^+(A)$ of ${\rm        Rep}_h(A)$ with objects $C^*$-representations is a tensor $C^*$-category.
\end{cor}

\begin{proof}
The $\Omega$-twisted inner product of  a tensor product of two   $C^*$-representations
is still a positive inner product by positivity of $\Omega$.  

\end{proof}

\begin{prop} Let $A$ be an $\Omega$-involutive  weak quasi bialgebra.
The forgetful functor ${\mathcal F}: {\rm        Rep}_h(A)\to{\rm        Herm}$ (or ${\mathcal F}: {\rm        Rep}^+(A)\to{\rm        Hilb}$ in the $C^*$-case)   is a   $^*$-functor. 
The natural transformations satisfy
\begin{equation}\label{adjoints} F_{\rho,\sigma}^*=\rho\otimes\sigma(\Omega)\circ G_{\rho,\sigma},\quad\quad G_{\rho,\sigma}^*= F_{\rho,\sigma}\circ \rho\otimes\sigma(\Omega^{-1}).\end{equation}

\end{prop}

\begin{proof}
$^*$-invariance of ${\mathcal F}$ is clear. Relations (\ref{adjoints}) follow from computations as in the proof of Theorem
\ref{Hermitian_category}.\end{proof}

We   observe that thanks to $G_{\rho,\sigma}\circ F_{\rho,\sigma}=\rho\otimes\sigma(\Delta(I))$, relations (\ref{adjoints}) can also be written in the form
 \begin{equation}\label{adjoints2}
 F_{\rho,\sigma}^*\circ F_{\rho, {\sigma}}=\rho\otimes\sigma(\Omega),\quad\quad G_{\rho, \sigma}\circ G^*_{\rho, \sigma}=\rho\otimes\sigma(\Omega^{-1}).
  \end{equation}

  \begin{prop}\label{unitary_equivalence_under_twisting}
  Let $A$ be an $\Omega$-involutive  weak quasi-bialgebra with involution $(^*, \Omega)$ and $F\in A\otimes A$ a twist. Consider the twisted algebra $A_F$ with involution $(^*, \Omega_F)$ as in Prop. \ref{twisted_omega}. Then the tensor equivalence ${\mathcal E}$ defined in Prop. \ref{class} restricts to a unitary tensor equivalence between ${\rm Rep}_h(A)\to{\rm Rep}_h(A_F)$ (${\rm Rep}^+(A)\to{\rm Rep}^+(A_F)$ in the unitary case).
  \end{prop}
  
  \begin{proof}
  The two algebras have the same $^*$-involution, hence the equivalence is a $^*$-functor. We show unitarity
  of the associated natural transformation, which is given by the action of $E_{\rho, \sigma}=\rho\otimes\sigma(F^{-1})$ from ${\mathcal E}(\rho)\otimes{\mathcal E}(\sigma)$ to ${\mathcal E}(\rho\otimes\sigma)$.
  We have $E_{\rho, \sigma}^*=\Omega_F^{-1}\rho\otimes\sigma({F^{-1}}^*)\Omega=\rho\otimes\sigma(F\Omega^{-1}F^*{F^{-1}}^*\Omega)=\rho\otimes\sigma(F)=E_{\rho, \sigma}^{-1}.$
  
  \end{proof}
  
  We next note   that while  at the algebraic level, the element $\Omega$ defining a unitary  involution of a weak quasi-Hopf $C^*$-algebra may be non-unique, passing to another such operator gives rise to a unitarily equivalent tensor
  $C^*$-category.
  
\begin{prop}\label{uniqueness}  Let $A$ be a weak quasi-bialgebra endowed with the structure of a $C^*$-algebra
(or a discrete weak quasi-bialgebra with positive involution). Let $\Omega$ and $\Omega'$ define
  two unitary $\Omega$-involutive structures. Let us  upgrade the category of $C^*$-representations of $A$ to corresponding   tensor $C^*$-categories   ${\rm        Rep}^+_\Omega(A)$ and ${\rm        Rep}^+_{\Omega'}(A)$. Then the functor
  ${\mathcal F}: {\rm        Rep}^+_\Omega(A)\to {\rm        Rep}^+_{\Omega'}(A)$ acting as identity on objects and morphisms 
     admits the structure of a   unitary tensor equivalence.
   \end{prop}
   
   \begin{proof} It is easy to check that
the functor ${\mathcal F}$ becomes a tensor $^*$-equivalence with the natural transformations ${\mathcal F}(\rho)\underline{\otimes}{\mathcal F}(\sigma)\to{\mathcal F}(\rho\underline{\otimes}\sigma)$ acting as identity. 
The unitary part of the polar decomposition equips ${\mathcal F}$ with the structure of a unitary tensor equivalence by Prop. \ref{polar_decomposition} b).
 \end{proof}

  We next discuss classification of $^*$-representations for the important class of  discrete $\Omega$-involutive weak quasi-bialgebras in the sense of Def. \ref{discrete_involutive}. So we may write, up to $^*$-isomorphism,  $A=\bigoplus_r M_{n_r}({\mathbb C})$.
The   projections $\rho_r: A\to M_{n_r}({\mathbb C})$ are irreducible   $C^*$-representations.

\begin{prop}\label{complete} The   $^*$-representations  $\rho_r$  
together with their opposites $\rho_{-r}$,  exhaust  the irreducible $^*$--representations    of $A$  up to unitary equivalence. Furthermore   any $^*$-representation of $A$ decomposes as an orthogonal direct sum of copies of them.
Finally, $\Omega$ is positive if and only of for all $s$, $t$, $\rho_s\underline{\otimes}\rho_t$ is an orthogonal direct sum of 
$\rho_r$ only.
\end{prop}
  
 \begin{proof} When we forget about  the $^*$-structure, an irreducible representation $\rho$ of $A$ is equivalent to some  $\rho_r$.
 Therefore to classify irreducible $^*$--representations, we need to classify up to unitary equivalence the Hermitian forms on ${\mathbb C}^{n_r}$ making  $\rho_r$ into a $^*$-representation. By    Prop. \ref{nondegeneracy} these are $\rho_r$ and $\rho_{-r}$.
We have already noticed that $\rho_r$ and $\rho_{-r}$ are not unitarily equivalent, hence    altogether they form a complete list of irreducible $^*$--representations, up to unitary equivalence. Let now $\sigma$ be a reducible $^*$-representation of $A$ and let us decompose it, as a representation, as a direct sum of  certain $\sigma_r$, where $\sigma_r$ is a multiple of $\rho_r$. Each $\sigma_r$ acts on $V_{\sigma_r}=\sigma(e_r)V_\sigma$, with $e_r$ a minimal central projection of $A$. Hence these subspaces are pairwise orthogonal by $^*$-invariance of $\sigma$. In particular, the form of $V_{\sigma}$ is nondegenerate on all the $V_{\sigma_r}$.
  In turn, the pairwise equivalent  irreducible summands $\tau_i$ of a fixed $\sigma_r$ act on the linear span $V_i$ of 
 $\{\sigma_r(e_{11})v_i, \sigma_r(e_{21})v_i,\dots, \sigma_r(e_{n_r 1})v_i\}$ respectively, where $v_i$ form a linear basis of $\sigma_r(e_{11})V_{\sigma_r}$ and we claim that it is possible to choose $v_i$ pairwise orthogonal.  The claim shows that these   copies of $\rho_r$ act on pairwise orthogonal subspaces. To show the claim, notice that the map 
 $v\in \sigma_r(e_{11})V_{\sigma_r}\to \sigma_r(e_{i1})v\in \sigma_r(e_{ii})V_{\sigma_r}$ is unitary between pairwise orthogonal subspaces of $V_{\sigma_r}$, hence the form of $V_{\sigma_r}$ must be nondegenerate on each of them, and the claim follows.  To show the last assertion, we use an orthogonal decomposition into irreducibles in the general case, given by isometries $S^{\pm}_{r, j}\in (\rho_{\pm r}, \rho_s\underline{\otimes}\rho_t)$. These determine the components 
 $ \rho_s\otimes\rho_t(\Omega)$ in the full matrix $C^*$-subalgebras of $A\otimes A$ by the formula
 $(\xi, \rho_s\otimes\rho_t(\Omega)\eta)_p=\sum ({S^{\pm}_{r, j}}^*\xi, {S^{\pm}_{r, j}}^*\eta)$, where $\xi$, $\eta$ vary in the vector space of $\rho_s\underline{\otimes}\rho_t$ and the inner products at the right hand side refer to $\rho_{\pm r}$.
 The claim easily follows from this equation.
 
   \end{proof}

 \begin{rem}\label{non-unitarizability}
Examples
have been found by  Fr\"ohlich and Kerler  \cite{FK} and    Rowell \cite{Rowell1, Rowell2, Rowell3}
of braided fusion categories which are not unitarizable.

 \end{rem}

 \section{Unitary braided symmetry and involutive Tannaka-Krein duality}\label{10}

In this section we discuss properties of the involutive structure in a weak quasi-Hopf algebra concerning the  twisting operation, quasitriangular structure and Tannaka-Krein duality. We start with categorical  interpretation of a twist of the 
$^*$-structure of a weak quasi bialgebra, in analogy with Prop. \ref{class} for a twisted bialgebra structure.

Let $A$ be a (discrete) weak quasi bialgebra and $(^*,\Omega)$ an $\Omega$-involution in the sense
of Def.  \ref{twisted_involution}. Let $t\in A$ (or $t\in M(A)$ if $A$ is discrete) be a selfadjoint twist, and consider the corresponding twisted involution
$(^\dag, \Omega_t)$, see Prop. \ref{1_twist}. We thus have two structures $(A, \varepsilon, \Delta, \Phi, ^*, \Omega)$
and $(A, \varepsilon, \Delta, \Phi, ^\dag, \Omega_t)$ which differ only for their involution. For brevity, we denote them
respectively as $A$ and $A_t$, in analogy with a twist of the bialgebra structure. Consider the functor 
 ${\mathcal E}: {\rm Rep}_h(A)\to {\rm Rep}_h(A_t)$ defined as follows. If $\rho$ is a $^*$-representation of $A$ then
 we modify the Hermitian form  $(\xi, \eta)_{V_\rho}$   of $V_\rho$ as $(\xi, \eta)_t:=(\xi, \rho(t)\eta)_{V_\rho}$, and consider
 the representation $\rho_t$ of $A_t$ on the Hermitian space $V_{\rho_t}$ so obtained acting as $\rho$. By construction,
 $\rho_t$ is a $^\dag$-representation of $A_t$. 
  
 \begin{prop}\label{1_twist_equivalence} Let $A$ be a unitary (discrete) weak quasi bialgebra and $t$ a positive twist for the involutive structure.
 Then the functor ${\mathcal E}: {\rm Rep}^+(A)\to {\rm Rep}^+(A_t)$
 taking   $\rho$    to $\rho_t$, acting identically on morphisms  and with identity natural transformations is a unitary tensor equivalence of tensor $C^*$-categories.

 \end{prop}
 
 \begin{proof} Pick $\rho$, $\sigma\in {\rm Rep}^+(A)$. For any linear map $T: V_\rho\to V_\sigma$,
 the adjoint of $T$ with respect to the original and modified Hermitian forms are related by 
 $T^\dag=\sigma(t^{-1})T^*\rho(t)$.  Thus if   $T\in(\rho, \sigma)$ then $T^\dag = T^*$, and this shows that ${\mathcal E}$ is a $^*$-functor, which is clearly full, faithful and essentially surjective, hence a $^*$-equivalence. On the other hand, the tensor structures of ${\rm Rep}^+(A)$ and ${\rm Rep}^+(A_t)$ are identical, hence ${\mathcal E}$
  is a tensor equivalence
 under the identity natural transformations. To show unitarity we are left to verify that the inner products of 
 $\rho_t\otimes \sigma_t$ and $(\rho\otimes\sigma)_t$ coincide, but this follows from a straightforward computation.

 \end{proof}

  It is well known that if $A$ is a  quasitriangular quasi-Hopf algebra with $R$-matrix $R$, the category ${\rm        Rep}(A)$
has a braided symmetry $\varepsilon$, where the action of  $\varepsilon(\rho, \sigma)$ on the representation space $V_\rho\otimes V_\sigma$ is given by $\Sigma R$, with $\Sigma: V_\rho\otimes V_\sigma\to  V_\sigma\otimes V_\rho$ the permutation operator. This construction extends to the weak case.  Similarly, if $A$ has an $\Omega$-involution, ${\rm        Rep}_h(A)$ is a braided tensor category as well. We next observe a condition on $R$ assuring unitarity
of $\varepsilon$ in ${\rm  Rep}_h(A)$.

 \begin{prop}\label{unitary_braided_symmetry2}
 Let $A$ be an $\Omega$-involutive   weak quasi-bialgebra with quasitriangular structure defined by $R$ and satisfying
 $\tilde{R}=R_\Omega$. Then the associated braided symmetry of ${\rm        Rep}_h(A)$ is unitary. If $A$ is discrete the converse holds.
 \end{prop}
 
 \begin{proof} Our assumption on the $R$-matrix means ${R^*}^{-1}=\Omega_{21} R\Omega^{-1}$.
The relation between the adjoint morphism $\varepsilon(\rho, \sigma)^*$ with respect to the $^*$-structure of ${\rm        Rep}_h(A)$ and the adjoint $\varepsilon(\rho, \sigma)^\dag$ with respect to the product form is
 $\varepsilon(\rho, \sigma)^*=\Omega^{-1}\varepsilon(\rho, \sigma)^\dag\Omega$. 
 Therefore
  $$\varepsilon(\rho, \sigma)^*=\Omega^{-1}(\Sigma \rho\otimes\sigma(R))^\dag\Omega=\Omega^{-1}( \rho\otimes\sigma(R^*)) \Omega_{21}\Sigma=R^{-1}\Sigma=\varepsilon(\rho,\sigma)^{-1}.$$

 \end{proof}

\begin{rem}\label{invariance_unitarity_of_braiding_under_twisting} The assumptions in   Prop. \ref{unitary_braided_symmetry2} may be read as saying  that the twist relation $\tilde{A}=A_\Omega$ holds not only at the level of   weak quasi-bialgebras, but also for their natural quasitriangular structures, cf. Prop. \ref{R-canonicity}. Furthermore if  $\tilde{R}=R_\Omega$ holds for a given $\Omega$-involutive quasitriangular weak quasi-bialgebra with $R$-matrix $R$ and involution $\Omega$ then they hold
 for any twisted algebra with twisted $R$-matrix $R_T$ and twisted involution $\Omega_T$, $\Omega_T^{-1}$
 as defined in c) of Prop. \ref{R-canonicity} and Prop. \ref{twisted_omega} respectively.
 \end{rem}
 
 \begin{cor}
 Let $A$ be a finite dimensional  discrete weak quasi-Hopf algebra with a quasitriangular structure $R$. Then any
 involution $(*, \Omega)$ making $A$ into a unitary weak quasi  bialgebra satisfies $\tilde{R}=R_\Omega$.
 \end{cor}
 
 \begin{proof}
The tensor $C^*$-category ${\rm Rep}^+(A)$ is  braided and fusion, hence by Theorem 3.2 in \cite{Gal} the braided symmetry is unitary. We may then apply Prop.  \ref{unitary_braided_symmetry2}.
   
 \end{proof}

  We next discuss a version of Tannaka-Krein duality for $\Omega$-involutive weak quasi  bialgebras.
  Recall that unitarity of a weak quasi tensor $^*$-functor was defined in Def. \ref{unitarity}, and that
  triviality of an $\Omega$-involution is introduced in Def. \ref{trivial_involution}.

  \begin{thm}  
  \label{TheoremTannakaStar}
  Let ${\mathcal C}$ be a semisimple  tensor $^*$-category, with finite dimensional morphism spaces ${\mathcal F}: {\mathcal C}\to{\rm        Herm}$ a faithful weak quasi tensor $^*$-functor defined by $(F, G)$ and   
  $A={\rm        Nat}_0({\mathcal F})$ be the discrete weak quasi bialgebra associated to ${\mathcal F}$ as in Th. \ref{TK_algebraic_quasi} and Th. \ref{TK_algebraic} endowed with its natural involution $^*$. Then  
  \begin{itemize}
\item[{\rm        a)}] 
   the element $\Omega\in A\otimes A$ defined by $\Omega_{\rho, \sigma}=F_{\rho, \sigma}^*\circ F_{\rho, \sigma}$
  makes $A$  into an $\Omega$-involutive   weak quasi  bialgebra, 
  \item[{\rm        b)}]   there is a canonical unitary tensor $^*$-functor  ${\mathcal E}: {\mathcal C}\to{\rm        Rep}_h(A)$   and is   an equivalence.  Furthermore,  the composite of ${\mathcal E}$ with the forgetful functor ${\rm        Rep}_h(A)\to{\rm        Herm}$ is unitarily monoidally isomorphic to ${\mathcal F}$,
   \item[{\rm        c)}] $({\mathcal F}, F, G)$ is (strongly) unitary if and only if $A$   the $\Omega$-involution of $A$ as in a) is (strongly) trivial,
    \item[{\rm        d)}] when ${\mathcal C}$ is unitary and  ${\mathcal F}: {\mathcal C}\to{\rm        Hilb}$ then 
    $A$ is a unitary weak quasi-bialgebra and ${\mathcal E}$  is a unitary tensor equivalence  between ${\mathcal E}: {\mathcal C}\to{\rm        Rep}^+(A)$.

  \end{itemize}
  
   \end{thm}
   
   \begin{proof} a)
  For simplicity in the following computations we drop the indices of the natural transformations. Note that $\Omega$ is selfadjoint, and in particular positive when ${\mathcal F}$ takes values in ${\rm Hilb}$.
  Furtherore, $\Omega$
  has $\Delta(I)=GF$ as   domain   and $\Delta(I)^*=(GF)^*$ as range. We set 
  $\Omega^{-1}:=GG^*$. We have: $\Omega^{-1}\Omega=GG^*F^*F=G(FG)^*F=GF=\Delta(I)$ and similarly
  $\Omega\Omega^{-1}=\Delta(I)^*$. Furthermore, for $\eta\in A$,
  $$\Omega\Delta(\eta^*)=F^*F\Delta(\eta^*)=F^*FG\eta_{\rho\otimes\sigma}^*F=$$
  $$F^*\eta_{\rho\otimes\sigma}^*F=
  F^*\eta_{\rho\otimes\sigma}^*G^*F^*F=\Delta(\eta)^*\Omega.$$ We have thus verified axioms (\ref{sa}), (\ref{domain}),
  (\ref{intertwiner}), while
   (\ref{normalized}) follows easily from (\ref{normalized_iota}) and   (\ref{eqn:omega0}) can be checked with  computations
  similar to those above. b) By assumption, ${\mathcal F}(\rho)$ is an Hermitian space and by Theorem 
   \ref{TK_algebraic_quasi} ${\mathcal E}(\rho)$ is a representation of $A$ on ${\mathcal F}(\rho)$ and ${\mathcal E}$ is a tensor equivalence with ${\rm Rep}(A)$ and therefore also with ${\rm Rep}_h(A)$.
  It is easy to check that ${\mathcal E}$ is $^*$-preserving, it follows that ${\mathcal E}$ takes values in ${\rm Rep}_h(A)$. To show unitarity of ${\mathcal E}$
  recall that the tensor structure  of ${\mathcal E}$ regarded   as a morphism in  ${\rm        Rep}_h(A)$ is   $F_{\rho, \sigma}$.
We compute the adjoint  
  $F_{\rho, \sigma}^*$ in ${\rm        Rep}_h(A)$.
  As before,
  we momentarily denote by $^\dag$ the usual adjoint of the tensor category of Hilbert spaces.
  We have $$F_{\rho, \sigma}^*=\Omega^{-1}_{\rho, \sigma} F_{\rho, \sigma}^\dag=G_{\rho, \sigma}G_{\rho, \sigma}^\dag F_{\rho, \sigma}^\dag=G_{\rho, \sigma}(F_{\rho, \sigma}G_{\rho, \sigma})^\dag=G_{\rho, \sigma}.$$
  c) By definition of unitarity of $({\mathcal F}, F, G)$,   $F^*F=P^*P$ and $GG^*=PP^*$, with $P=GF$
  and this by construction corresponds to triviality of the $\Omega$-involution of $A$, and similarly for the relation between strong unitarity of the weak quasi-tensor structure and strong triviality of the $\Omega$-involution.

   \end{proof}

    \begin{rem} 
Theorem  \ref{TheoremTannakaStar} for unitary weak quasi-bialgebras  has origin  in \cite{HO} where the author assumes 
that 
$F_{\rho, \sigma}^*=G_{\rho, \sigma}$ and are  isometries, that is a strongly unitary structure in our terminology.  In this case he proves that   the $^*$-involution of $A$ commutes with the coproduct.
We note that  the   examples that we discuss in Sect. \ref{20} arising from quantum groups at roots of unity    do not satisfy this property, and this motivated us to consider the more general case.    \end{rem}

\begin{ex}
Consider the  $^*$-category   ${\mathcal C}={\rm Herm}^\omega_G$ of $G$-graded Hermitian spaces. It becomes a tensor $^*$-category with natural tensor product 
and associator given by
a ${\mathbb T}$-valued $3$-cocycle $\omega$ over $G$.  For every $g\in G$, denote by ${\mathbb C}_g^+$  (${\mathbb C}_g^-$) the   one-dimensional Hermitian space of degree $g$ with positive (negative) scalar product. Then  
${\mathbb C}_g^+$ and ${\mathbb C}_g^-$
 are two irreducible equivalent but not unitarily equivalent 
  objects of ${\rm Herm}^\omega_G$,  and   ${\mathbb C}_g^{\pm}$ and ${\mathbb C}_h^{\pm}$ are inequivalent
   for $g\neq h$.
  The category ${\rm Herm}^\omega_G$ contains ${\rm Hilb}^\omega_G$ as a full tensor $C^*$-subcategory with restricted
  $^*$-structure.
  Consider  ${\mathcal F}: {\rm Herm}^\omega_G\to{\rm Herm}$ the forgetful functor. Note that ${\mathcal F}$ preserves the Hermitian forms, thus it takes a definite sign on the unitarily inequivalent simple objects. It follows that
  ${\rm Nat}_0({\mathcal F})$ is a pre-$C^*$-algebra that may be identified with the $C^*$-algebra of complex-valued functions on $G$. Note that $F(g)\otimes F(h)$ and $F(gh)$ are unitarily equivalent Hermitian spaces with definite forms, thus
every quasitensor structure ${F}_{g, h}$ on ${\mathcal F}$ satisfies $\Omega(g, h):={F}_{g, h}^*F_{g, h}>0$.
It follows from Theorem \ref{TheoremTannakaStar} that $A={\rm Nat}_0({\mathcal F})$ is a unitary pointed quasi-bialgebra
which identifies with $C_\omega(G)$ with unitary structure defined by $\Omega$
as in Example \ref{pointed_Hermitian_qh}. Note that by the last part of Example \ref{pointed_Hermitian_qh}
 there exist examples of  pointed   tensor $^*$-categories which are not unitarily equivalent to some 
 ${\rm Herm}^\omega_G$.

\end{ex}

\begin{rem}\label{polar_decomposition_and_twist}  
 In Sect. \ref{8} we have constructed the unitarization $A_{\Omega^{1/2}}$ associated to a unitary discrete weak quasi-bialgebra $A$ in the case where $({\Omega^{-1}})^{1/2}$ is a left inverse of $\Omega^{1/2}$. This construction may be described categorically as follows. 
  Let (${\mathcal F}$, $F$, $G$) be a faithful weak quasi tensor $^*$-functor of a semisimple unitary tensor category ${\mathcal C}$ and $A$ the associated unitary discrete
weak quasi bialgebra with involution denoted $(^*, \Omega)$ following  Theorem \ref{TheoremTannakaStar}. If this functor is non-unitary and for example we know that satisfies the left inverse property
(\ref{left_inverse_condition}) then
we may consider the   unitarized functor
(${\mathcal F}$, $F'$, $G'$) as in part a) of Prop. \ref{polar_decomposition}, see also
Def.    \ref{unitarization_of_a_functor}. This new structure    in turn gives rise to a new unitary weak quasi bialgebra $B$ corresponding to the unitarization $A_{\Omega^{1/2}}$ of $A$,   by the proof
of Theorem \ref{propweakdim} with trivial unitary structure   by Cor. \ref{making_omega_trivial}. This structure is also strongly trivial  by Prop. \ref{strongly_unitary_prop}. 
\end{rem}

The notion of unitarization will have a useful  extension  in Sect. \ref{19}  in that will be applied to more useful situations in subsequent sections.

We ask how to construct and parameterise faithful   $^*$-functors ${\mathcal G}: {\mathcal C}\to{\rm Hilb}$ from a $C^*$-category.  If   ${\mathcal G}$ is given, we may construct new   $^*$-functors  to ${\rm Hilb}$ via a categorical analogue of the twist deformation of the involution of an algebra of Prop. \ref{1_twist} in the following way. Let $t\in{\rm Nat}_0({\mathcal G})$ be a positive invertible 
natural transformation and let ${\mathcal G}_t(\rho)$ be ${\mathcal G}(\rho)$ as a vector space, but with modified inner product  $(\xi, \eta)_t:=(\xi, t_\rho\eta)_{{\mathcal G}(\rho)}$. The action of ${\mathcal G}_t$ on  morphisms  is the same as that of ${\mathcal G}$.  The fact that ${\mathcal G}$ is $^*$-preserving   together with naturality of $t$ easily  imply that 
${\mathcal G}_t$ is $^*$-preserving as well, hence a $^*$-functor. 
 The $^*$-algebras $A={\rm Nat}_0({\mathcal G})$ and $B={\rm Nat}_0({\mathcal G}_t)$ are related by $B=A_t$.

Faithful functors ${\mathcal F}: {\mathcal C}\to{\rm Vec}$ are described, up to isomorphism, by 
functions $D:{\rm Irr}({\mathcal C})\to{\mathbb N}$ thanks to Theorem  \ref{propweakdim} (a). We thus need to parameterize  the ways how  ${\mathcal F}$ can be written as ${\mathcal  F}={\mathcal H}{\mathcal G}$ with  ${\mathcal G}: {\mathcal C}\to {\rm Hilb}$ is a $^*$-functor and $${\mathcal H}:{\rm Hilb}\to{\rm Vec}$$   the forgetful functor.

 \begin{prop}\label{Haring_functor_$C^*$} Let ${\mathcal C}$ be a $C^*$-category with finite dimensional morphism spaces and ${\mathcal F}:{\mathcal C}\to{\rm Vec}$ a faithful  functor. Then ${\mathcal F}$ 
 factors through 
${\mathcal  F}={\mathcal H}{\mathcal G}$ where
  ${\mathcal G}: {\mathcal C}\to {\rm Hilb}$ is a faithful  $^*$-functor and ${\mathcal H}:{\rm Hilb}\to{\rm Vec}$  the forgetful functor. Any other  $^*$-functor ${\mathcal G}'$  with the same properties is of the form ${\mathcal G}_t$ for a unique positive invertible $t\in{\rm Nat}_0({\mathcal G})$.
  \end{prop}
  
  \begin{proof} We choose, for each $\rho\in {\rm Irr}({\mathcal C})$, a positive inner product on ${\mathcal F}(\rho)$, and let ${\mathcal G}(\rho)$ the corresponding Hilbert space. 
  Note that ${\mathcal F}(T^*)={\mathcal F}(T)^*$  holds for $T\in(\rho, \rho)$ for any choice of inner product when $\rho$ is irreducibile, since these morphisms are scalars and ${\mathcal F}$ is linear.
  We use orthogonal complete reducibility of   $\mu\in {\mathcal C}$ via isometries $S_{\rho, i}\in(\rho, \mu)$ with $\rho$ irreducibles,
  to extend the construction of a Hilbert space ${\mathcal G}(\mu)$ to all objects $\mu$ via 
  $(\xi, \eta)_{{\mathcal G}(\mu)}:=
  \sum_{\rho, i}({\mathcal F}(S_{\rho, i}^*)\xi, {\mathcal F}(S_{\rho, i}^*)\eta)_{{\mathcal G}(\rho)}$. It follows that the inner product is independent of the choice of the  isometries $S_{\rho, i}$. Letting ${\mathcal G}$ act as ${\mathcal F}$ on morphisms, one sees that ${\mathcal G}(S_{\rho, i}^*)={\mathcal G}(S_{\rho, i})^*$ and this implies ${\mathcal G}$ is $^*$-preserving. Another decomposition ${\mathcal F}={\mathcal H}{\mathcal G}'$ gives a new  Hilbert space structure 
  ${\mathcal G}'(\rho)$ on the same vector space as ${\mathcal G}(\rho)$, hence we may find a unique positive invertible operator $t_\rho$ on ${\mathcal G}(\rho)$ such that $(\xi, \eta)_{{\mathcal G}'(\rho)}=(\xi, t_\rho\eta)_{{\mathcal G}(\rho)}$.
  Since ${\mathcal G}'$ is a $^*$-functor, this   implies that $t\in{\rm Nat}_0({\mathcal G})$.
    
  \end{proof}
  
  We summarise the main results of this and previous  sections.
  
  \begin{cor}\label{HO_star}
  Let ${\mathcal C}$ be a tensor  $C^*$-category with finite dimensional morphism spaces and $D$ a weak dimension function on ${\mathcal C}$. Then there is a faithful weak quasi-tensor $^*$-functor ${\mathcal G}:{\mathcal C}\to{\rm Hilb}$ such that $D(\rho)=
  {\rm dim}({\mathcal G}(\rho))$. If $A$ is the discrete unitary   weak quasi bialgebra corresponding to ${\mathcal G}$ via duality then all the others corresponding to different weak quasi-tensor $^*$-functors with the same dimension function
  are isomorphic to $A_{F, t}$ for some twist $F\in A\otimes A$ and $t\in A^+$ of the bialgebra and $\Omega$-involutive structure of $A$ respectively.
  \end{cor}

\section{Unitarizability of representations and  rigidity  
}\label{11}

Let $A$   be a complex $^*$-algebra  either unital or discrete. In this section we are motivated by    studies of direct constructions of  unitary structures on some objects of ${\rm Rep}(A)$ where $A$ is a suitable algebra arising in an application. In categorical terms we want to construct objects of ${\rm        Rep}^+(A)$ and obtain a nontrivial linear $C^*$-category. 
 
   In order to construct objects of ${\rm        Rep}^+(A)$,  we need to know which representations of $A$ are equivalent to  $C^*$--representations.   In this section we abstract to $A$ procedures used  in the setting of group theory,
     quantum groups at generic parameters or at roots of unity. This last case is of main interest in our paper, because
 we see a   similar procedure in the setting of vertex operator algebras, and wish to look for a unified description between the two settings.
  
 In the case of representations of quantum groups at roots of unity $U_q({\mathfrak g})$ and rational vertex operator algebras $V$,  the question splits into a first step which reduces the problem to the search of invariant
 Hermitian forms on the representation space. 
 We discuss this first step   in abstract terms here looking for objects of ${\rm        Rep}_h(A)$.  Both cases are covered
 by the assumption that $A$ is discrete, but this can not be our starting assumption as $U_q({\mathfrak g})$
 is not semi-simple. Later on   we shall discuss Wenzl functor defined on the semi-simple fusion category associated to $U_q({\mathfrak g})$ which gives
 rise to a discrete quotient $A$
 of $U_q({\mathfrak g})$ that keeps all the information of the linear structure of the category. This discrete algebra plays in fact
 an important role in showing that all its irreducible representations   become $C^*$-representations compatibly with a natural $^*$-involution of
 $U_q({\mathfrak g})$ in  \cite{Wenzl}. We shortly comment some more on this, for details on the construction of this algebra, see   Sect. \ref{unitary_structure_of_fusion_category}.
 
 The question of making representations of the discrete algebra $A$ into $C^*$-representations is equivalent
 to making $A$ into a $C^*$-algebra. Until discussing more in detail the construction of the discrete algebra, we may think of $A$ to be the non-semisimple unital $^*$-algebra
 $A=U_q({\mathfrak g})$.
 
Regarding the case of vertex operator algebras, representations are infinite dimensional.
 To unify with the case of quantum groups, the main important point is  to look at the Zhu algebra 
 $$A=A(V)$$  associated to a vertex operator algebra $V$, which is associative, finite dimensional and discrete under some rationality assumptions on $V$  recalled in Sect.  \ref{VOAnets2}. We develop a general correspondence
 between invariant Hermitian forms of modules of a vertex operator algebra under suitable conditions and invariant Hermitian forms on modules of the Zhu algebra, see Prop.\ref{propUnitaryZhu}.

The second step regards
 verification of positivity of Hermitian forms, and
   is studied in a direct way in applications. We do not develop an abstract theory for this, although it is apparent that methods are of a similar
   nature.
We discuss  in detail in this paper the case of affine vertex operator algebras at a positive integer level and quantum groups
 at certain roots of unity via the mentioned unified algebraic method. The complete result is of course  not new, and may be regarded as going back to the work
 of Kac on the vertex operator algebra or affine Lie algebra side (where one has  unitary representations of the compact real form of ${\mathfrak g}$) and of Wenzl on the quantum group side. 
However,  the unified   
 viewpoint via discrete algebras $A$ may have been overlooked in the literature. This viewpoint is perhaps of help
 to discuss how  positivity of the Hermitian form in the quantum group side may be seen as derived from
 positivity on the vertex operator algebra side following Wenzl  continuous path argument that connects the two setting on each relevant
 irreducible representation, of the Zhu algebra and of the discrete algebra associated to the quantum group. The unitary structure gives a connection from the module category of the VOA to the module category of a semisimple quotient of the quantum group.
 Then after getting to the setting of quantum groups, over the discrete algebra we construct a unitary weak  Hopf algebra structure compatible
 with the fusion category structure, and extend Wenzl connection to the module categories.
 
 This suggests to study more connections between unitarizable  quantum group fusion categories and  vertex operator algebra fusion categories.

 \medskip

 After these motivating remarks, we start with $A$, and pose a question more precisely. Given a  (nondegenerate f.d.) representation $\rho$ of $A$, is there and 
how  do we construct a $^*$-representation $\sigma$ of $A$ on a nondegenerate Hermitian space equivalent to $\rho$? 
Do Hermitian representations naturally  arise or are they just a more general convenient setting to study unitary representations?

 This first step is related to the study of the {\it contragredient}  representation of a representation $\rho$ of $A$
 and its relation with the  {\it conjugate} representation. This is related to the study of rigidity of  ${\rm        Rep}(A)$ 
 if $A$ has the structure of a weak quasi-Hopf algebra.

Recall that the  contragredient  representation $\rho^c$ (or $^c\rho$) has been defined 
in Def. \ref{duals_wqh} on the dual vector space $V_\rho^*$ by means of an invertible antimultiplicative map $S:A\to A$.
 The  conjugate representation acts on the conjugate 
vector space $\overline{V_\rho}$,  and is defined here below in (\ref{conjugate_representation}) and depends on $S$ as well and the $^*$-operation.

 In order to have a nondegenerate Hermitian form on the space $V_\rho$,
by Riesz representation theorem
we at least need to find a canonical linear invertible   map
$$\Phi: \overline{\xi}\in\overline{V_\rho}\to\Phi_{\overline{\xi}}\in V_\rho^*.$$ 

Our first result, Prop. \ref{conjugates}, shows that the existence of such a $\Phi$ which in addition is
an intertwiner between the   conjugate representation $\rho_c$  on $\overline{V_\rho}$
and the contragredient representation $\rho^c$   on $V_\rho^*$
   is a necessary condition for the existence of a nondegenerate Hermitian form on $V_\rho$ making $\rho$ into a $^*$-representation.
   This intertwining condition corresponds to invariance in the setting of vertex operator algebras.
   Moreover, this condition is also sufficient if $\rho$ is irreducible, and the Hermitian form is unique up to scalar multiples. From this point one can further study the case where $\rho$ is a $C^*$-representation in applications. 

For example, this is the case if $A$ is a discrete
$C^*$-algebra, by Cor. \ref{cor_conjugates}.
This Corollary   is well known in some applications,  for example in the theory of compact groups.

  Recall the definition of contragredient representations $\rho^c$ and ${}^c\rho$ given in \ref{duals_wqh}.
 Fur possible future applications, it is important to note that by Prop. \ref{rigidity} the contragredient representations have a natural description in the setting of
 weak quasi-Hopf algebras, and therefore of tensor categories, because they solve the duality equations of the tensor category ${\rm        Rep}(A)$ making it into a rigid tensor category. On the other hand, one may work with the contragredient representations without any 
 tensor structure. Indeed 
the contragredient representations  may be defined more generally in the case where $A$ is only a complex algebra, always assumed either unital or discrete to avoid degeneracies,
endowed with a linear invertible antimultiplicative
map $S:A\to A$. Thus we only need a  pair $(A, S)$.   It follows that  
the contragredient representations $\rho^c$ and ${}^c\rho$ are nondegenerate by invertibility of $S$ and clearly with the same dimension as
$\rho$. 

This more general starting point is important for applications
to a vertex operator algebra $V$, where one may start with the Zhu algebra $A(V)$ associated to $V$, which is known to be associative, unital
and endowed an invertible map $S$, known as Zhu's antipode.

Let now $(A, ^*, S)$   be a triple with   $(A, ^*)$ a complex $^*$-algebra as at the beginning of the section, and $S:A\to A$ an invertible linear antimultiplicative map.
 Then  to any nondegenerate finite dimensional representation $\rho$ on a vector space we may associate two more representations, $\rho_c$ and ${}_c\rho$ both acting on the conjugate vector space $\overline{V_\rho}$ via
\begin{equation}\label{conjugate_representation}
\rho_c(a)\overline{\xi}=\overline{\rho(S(a)^*)\xi},\quad\quad  {}_c\rho(a)\overline{\xi}=\overline{\rho(S^{-1}(a)^*)\xi}.
\end{equation}
   
 \begin{rem}
 Alternatively, we may consider the    representations acting on $\overline{V_\rho}$ via $\overline{\xi}\to \overline{\rho(S(a^*))\xi}$
  and $\overline{\xi}\to \overline{\rho(S^{-1}(a^*))\xi}$.
However, if  $(A, S, ^*)$ can be completed to
the structure of    an $\Omega$-involutive  a weak quasi-Hopf  algebra $(A, \Delta, \Phi, S, \alpha, \beta, ^*, \Omega)$,  then the two further
representations defined above are
         equivalent to ${}_c\rho$ and $\rho_c$ respectively. Moreover,
there are equivalences $\rho_{cc}\simeq\rho\simeq{}_{cc}\rho$. These properties easily follow from     Prop.   \ref{involution_antipode}.
  \end{rem}

  The following proposition  is an abstraction of the construction of Kashiwara inner product, and of Lemma. 2.2 and Prop. 2.3 of \cite{Wenzl}.

\begin{prop}\label{conjugates}
Let $(A, ^*, S)$ be  be given as above, and let $\rho$ be a finite dimensional nondegenerate vector space representation of  $A$. 
\begin{itemize}
\item[{\rm        a)}] If $\rho$ is equivalent to a $^*$-representation  on a nondegenerate Hermitian space  then there is a linear invertible map
$$\overline{\xi}\in\overline{V_\rho}\to\Phi_{\overline{\xi}}\in V_\rho^*, \quad \quad (\text{resp. }
 \Phi': \overline{\xi}\in\overline{V_\rho}\to\Phi'_{\overline{\xi}}\in V_\rho^*)$$ which is an 
equivalence between the representations $$\Phi:\rho_c\to\rho^c, \quad \quad (\text{resp. }
 \Phi':{}_c\rho\to{}^c\rho)$$ such that the Hermitian form of $V_\rho$ is given by  $$(\xi, \eta)=\Phi_{\overline{\xi}}(\eta),  \quad\quad (\text{resp. }
 (\xi, \eta)=\Phi'_{\overline{\xi}}(\eta))$$
 
\item[{\rm        b)}] if $\rho$ is irreducible and if $\rho^c\simeq\rho_c$ (or ${}^c\rho\simeq{}_c\rho$)   then $\rho$ is  equivalent to a $^*$-representation on a nondegenerate Hermitian space and the associated Hermitian form is unique up to a nonzero real scalar.
 
\end{itemize}

\end{prop}

\begin{proof}  a) If $\rho$ is equivalent to the $^*$-representation $\sigma$ via the invertible $T\in(\rho, \sigma)$ we may endow the space of $\rho$ with the nondegenerate Hermitian form making $T$ unitary, and in this way  $\rho$ 
   becomes a $^*$-representation. It follows that we may    canonically  identify the conjugate space $\overline{V_\rho}$ with $V_\rho^*$,
       via the invertible map $$\overline{\xi}\in\overline{V_\rho}\to\Phi_{\overline{\xi}}\in V_\rho^*,$$ where
       $\Phi_{\overline{\xi}}$  is the functional $\eta\to (\xi,\eta)$. A computation shows that $\Phi$ is an equivalence between  $\rho^c$
       acting on  $V_\rho^*$ and the representation
   acting on $\overline{V_\rho}$ as $\overline{\xi}\in\overline{V_\rho}\to \overline{\rho(S(a))^*\xi}=\overline{\rho(S(a)^*)\xi}=\rho_c(a)\overline{\xi}$.
    (Similarly, ${}^c\rho$  turns into   ${}_c\rho$.)
  b) Let $\Phi\in(\rho_c, \rho^c)$ be an invertible morphism, and introduce a sesquilinear form on $V_\rho$ by
  $(\xi, \eta)=\Phi_{\overline{\xi}}(\eta)$,   clearly nondegenerate. 
 Let us define the right and left adjoint of
  a linear map $T: V_\rho\to V_\rho$ respectively by $(T^*\xi, \eta)=(\xi, T\eta)$ and $(\xi, {}^*T\eta)=(T\xi, \eta)$.
    A   computation using the intertwining property of $\Phi$ shows that for $a\in A$, $\rho(a)^*=\rho(a^*)={}^*\rho(a)$.
Let us introduce an inner product $(\xi, \eta)_{{\rm        pos}}$  in $V_\rho$ making some basis orthonormal,   let
   $T\to T^\dag$ be the corresponding adjoint map and
   $B: V_\rho\to V_\rho$ be the unique invertible map such that $(\xi, \eta)=(\xi, B\eta)_{{\rm        pos}}$. Then 
   $T^*={B^{-1}}^\dag T^\dag B^\dag$ and ${}^*T=B^{-1}T^\dag B$. Equating $\rho(a)^*={}^*\rho(a)$ gives  ${B^{-1}}^\dag B\rho(a)=\rho(a) {B^{-1}}^\dag B$, hence $B^\dag$ is a scalar multiple of $B$ by irreducibility
   of $\rho$. But $\|B\|=\|B^\dag\|$ (norm associated to $(\xi, \eta)_{{\rm        pos}}$) and it follows that this scalar lies in ${\mathbb T}$.
   Hence after rescaling $B$ we get $B^\dag=B$, and finally derive that $(\xi, \eta)$ is Hermitian. 
Finally, with a similar argument, if $\rho$ is irreducible and unitary on a Hermitian space with Hermitian form $(\xi,\eta)$ then any other nondegenerate Hermitian form on the same space making $\rho$ $^*$-invariant, when written as  $(\xi, A\eta)$ with $A$ invertible and selfadjoint, implies that $A$ is a real scalar.
  \end{proof}

Let us now discuss the tensorial aspects.
If   in particular $(A, S)$   upgrades to  a weak quasi-Hopf algebra
$(A, \Delta, \Phi, S, \alpha, \beta)$   then  ${\rm        Rep}(A)$ is a rigid tensor category.  Indeed, the contragredient $\rho^c$ and 
$^c\rho$ are a right and left dual
of $\rho$ in  ${\rm  Rep}(A)$ respectively,  by Prop. \ref{rigidity}.

   If   $(A, ^*)$ upgrades to the structure of an $\Omega$-involutive weak quasi-bialgebra $(A, \Delta, \Phi, ^*, \Omega)$, then   ${\rm        Rep}_h(A)$ is a tensor $^*$-category
by Thoerem \ref{Hermitian_category}. When is ${\rm        Rep}_h(A)$ rigid?

Given $(A, ^*, S)$, and a representation $\rho$ of $A$ equivalent to a $^*$-representation, e.g. it satisfies     condition b) of Prop. \ref{conjugates}, then $\rho$ becomes a $^*$-representation   of $A$. On the other hand,
the contragredient representations $\rho^c$ and ${}^c\rho$ are not necessarily
 $^*$-representations if $S$ does not commute with $^*$.
 
 This   is well known to occur in the setting of compact quantum groups,
 although it does not occur in the setting of quantum groups at roots of unity, or in the setting of unitary vertex operator algebras because
 $S$ commutes with $^*$.
In the case where one is able to make any nondegenerate
 f.d. representation of $A$ equivalent to a $^*$-representation, then one can apply this procedure   to the contragredient representations, and this second step provides with equivalent $^*$-representations on nondegenerate Hermitian spaces.
 This connects with the work done in the setting of compact quantum groups. In this section we discuss this second step as well in the setting of weak quasi-Hopf algebras.

Recall    that the definition  of conjugate object  in a tensor $^*$-category is usually given as  in Def. (\ref{conjugates_def}). But
is equivalent to that of the underlying tensor category. A more precise discussion is before formula (\ref{conjugates_def}).
On the other hand, the underlying category  of ${\rm        Rep}_h(A)$ is not   ${\rm        Rep}(A)$, but ${\rm        Rep}(A)$
and ${\rm        Rep}_h(A)$ are equivalent iff every object of ${\rm        Rep}(A)$ is equivalent to a $^*$-representation.

A $^*$--representation $\rho\in{\rm        Rep}_h(A)$ has a conjugate   in ${\rm        Rep}_h(A)$ if and only if    the canonical right dual
  $\rho^c$ introduced in Def. \ref{duals_wqh}  is equivalent to a $^*$-representation, the latter is a conjugate $\overline{\rho}$ of $\rho$. If this is the case  the canonical left dual ${}^c\rho$ will be automatically
 equivalent to $\rho^c$ and $\overline{\rho}$     as well, and the double dual $\rho^{cc}$   to $\rho$.

Summarizing, if  $(A, ^*, S)$ has the structure of an $\Omega$-involutive weak quasi-Hopf algebra
 then   ${\rm        Rep}_h(A)$ a tensor $^*$-category and  is rigid as a tensor category. If
 every representation is equivalent
 to a $^*$-representation then  
  ${\rm        Rep}_h(A)$ is  a rigid tensor $^*$-category.

In the setting of this section, this amounts to verify that
the contragredient representations, which solve the right and left duality equations  in ${\rm        Rep}(A)$ are equivalent to $^*$-representations.

  \begin{cor}\label{cor_conjugates}
  Let $(A, ^*, S)$ be a discrete $C^*$-algebra.
  Then every representation is equivalent to a $C^*$-representation. In particular if $A$ 
  has the structure of a $\Omega$-involutive (unitary)  weak quasi-Hopf algebra then
  ${\rm Rep}_h(A)$ (${\rm Rep}^+(A)$) is rigid and 
  the forgetful functor ${\rm Rep}_h(A)\to {\rm Rep}(A)$ (${\rm Rep}^+(A)\to{\rm Rep}(A)$) is a 
  tensor equivalence.
  
    \end{cor}

  \begin{proof} Let $\rho$ be a representation of $A$ that we may assume irreducible by   complete reducibility. Note that the antipode $S$ permutes the minimal central idempotents of $A$ and that these idempotents are selfadjoint  since the involution of $A$ is positive
  by assumption, see  Def. \ref{discrete_involutive}.
  This implies that $\rho_c$ and $\rho^c$ have the same central support, and therefore they are related by an isomorphism $T$. We may  then apply   Prop. \ref{conjugates}. Note also that  a nonzero scalar multiple of $T$   induces  a  positive inner product on the space of $\rho$
  by the classification of Hermitian forms associated to involutive discrete weak quasi-Hopf algebras, Prop. \ref{complete}, hence the conclusion follows also in the case where $\Omega$ is positive.

  \end{proof}

  By the end of the section we shall identify the conjugates in the tensor $C^*$-category ${\rm Rep}^+(A)$ in the discrete weak Hopf case.
     We next discuss some results guaranteeing rigidity in possibly non-semisimple tensor categories
 motivated by the work of Kashiwara, Kirillov, Wenzl for $U_q({\mathfrak g})$ at  roots of unity \cite{Kashiwara, Kirillov, Wenzl}.
   Recall that the element $\omega$ was defined in Prop. \ref{involution_antipode}.

\begin{prop} \label{conjugates2}
Let $A$ be an $\Omega$-involutive weak quasi-Hopf algebra and $\rho$ a  $^*$-representation   equivalent to $\rho^{cc}$. Then 
\begin{itemize}
\item[{\rm        a)}] 
$\rho_c$ is equivalent to a $^*$--representation if and only if there is an invertible $K_\rho\in(\rho, \rho^{cc})$ such that $F_\rho:=K_\rho\rho(\omega^*)$ is selfadjoint with respect to the Hermitian form of $\rho$. In this case,    the   forms making  
$\rho_c$ into a $^*$-representation are parametrised by $K_\rho$ via
$(\overline{\xi}, \overline{\eta})=(\eta, F_\rho\xi)$,
\item[{\rm        b)}] 
  if $\rho$ is a $C^*$-representation
then $\rho_c$ is equivalent to a $C^*$-representation if and only if $F_\rho$ can be chosen positive,
\item[{\rm        c)}] if $\rho$ is irreducible then $\rho_c$ is equivalent to a $^*$-representation. The associated $K_\rho\in(\rho, \rho^{cc})$ is  unique up to a real scalar multiple. 
  \end{itemize}
  \end{prop}
  
  \begin{proof}
a)  If $\rho_c$ is equivalent to a $^*$--representation then $\rho_{cc}$ and $\rho_c^c$ are equivalent by the previous proposition, and let $\Psi$ be this equivalence. We may write $\Psi$ as the composite of $\rho(\omega^*): \rho_{cc}\to\rho$
with an equivalence $K_\rho: \rho\to\rho^{cc}$ in turned followed by $\Phi^t: \rho^{cc}\to(\rho_c)^c$, where $\Phi:\rho_c\to\rho^c$ is defined as in the proof of a) of Prop. \ref{conjugates}, and $\Phi^t$ is the transposed of $\Phi$. The Hermitian form making 
$\rho_c$ into a $^*$-representation is given by $(\overline{\xi},  \overline{\eta})=\Psi_{\xi}(\overline{\eta})$.  An explicit computation shows that this is precisely the form in the statement. Conversely, for any $K_\rho\in(\rho, \rho^{cc})$,  
the sesquilinear form defined by $F_\rho=K_\rho \rho(\omega^*)$ is Hermitian (positive) precisely when $F_\rho$ is selfadjoint (positive). A computation shows that $(\rho_c(a)\overline{\xi}, \overline{\eta})=(\overline{\xi}, \rho_c(a^*)\overline{\eta})$, in other words  $\rho_c$ is a $^*$-representation. The proof of b) is now clear. c)   By irreducibility and b) of  Prop. \ref{conjugates},  it suffices to show that
${\rho_{cc}}\simeq (\rho_c)^c$. Now such an equivalence can be obtained as in the proof of a) starting from the choice of an invertible $K_\rho\in(\rho, \rho^{cc})$.

  \end{proof}

  \begin{cor} Let $A$ be an $\Omega$-involutive weak quasi-Hopf algebra   with an antipode $(S, \alpha, \beta)$ such that $S$ commutes with $^*$. Then \begin{itemize}
\item[{\rm        a)}]      every $^*$-representation   $\rho$ has $\rho_c$ as  a conjugate in ${\rm        Rep}_h(A)$ with respect to  the form conjugate to that of $\rho$:
  $(\overline{\xi}, \overline{\eta})=(\eta, \xi)$. Hence ${\rm        Rep}_h(A)$ is rigid. 
\item[{\rm        b)}]   
   If $\rho$ is a $C^*$-representation, so is $\rho_c$. Hence if $A$ is a unitary   weak quasi-Hopf algebra,   ${\rm        Rep}^+(A)$ is rigid as well.
   \end{itemize}
\end{cor}

\begin{proof}
We may take $K_\rho=\rho({\omega^*}^{-1})$ by  Prop. \ref{Kac}, hence $F_\rho=I$ for all $^*$-representations $\rho$. 
\end{proof}

\begin{rem}\label{antipode_commuting_with_*} Let   $A={\rm Nat}_0({\mathcal F})$ be the discrete weak quasi-bialgebra associated to a semisimple tensor category ${\mathcal C}$ 
endowed a weak quasi-tensor functor ${\mathcal F}: {\mathcal C}\to {\rm Vec}$ as in Theorem
\ref{TK_algebraic_quasi}. When ${\mathcal C}$ is also a $C^*$-category and ${\mathcal F}$ factors through
a $^*$-functor
${\mathcal F}: {\mathcal C}\to {\rm Hilb}$ then $A$ has a natural
  pre-$C^*$-algebra involution. 
 If ${\mathcal C}$ is rigid and the dimension assumption of Prop. \ref{TK_algebraic_quasi} (d) hold (e.g. ${\mathcal C}$ is a fusion category) then  $A$ has an antipode $(S, \alpha, \beta)$. We note that $S$ may always be chosen commuting  with $^*$. Indeed, following the proof 
of Theorem \ref{TK_algebraic_quasi} (d),   for each $\rho$, ${\mathcal F}(\rho)^*$ identifies naturally with the conjugate vector space $\overline{{\mathcal F}(\rho)}$, which we endow with the unique Hilbert space structure 
making the conjugation map $J:{\mathcal F}(\rho)\to \overline{{\mathcal F}(\rho)}$ antiunitary. It also follows that a  transposed linear map
$L^t$ identifies with $JL^*J^{-1}$. On the other hand, we may choose the natural transformation $U$ unitary. It follows from the antipode formula given in the proof   that    $S(\eta^*)=S(\eta)^*$.
For example if ${\mathcal C}={\rm Rep}(A)$ with $A$ a discrete weak quasi-Hopf
algebra which is also a pre-$C^*$-algebra then  the procedure reconstructs the original antipode of $A$ when this commutes with $^*$ by 
 Remark \ref{faithfulness_and_antipode_reconstruction} c), but it gives a    new one otherwise.  
\end{rem}

\begin{ex}\label{computing_conjugates} 
 We next describe   the conjugate equations in  ${\rm        Rep}_h(A)$ (${\rm        Rep}^+(A)$) under the assumption that $S$ commutes with $^*$.
Given a $^*$-representation $\rho$, we may use the canonical identification of $\rho^c$ with  $\rho_c$ and
obtain from  Prop. \ref{rigidity}
  the following solution for the pair $\rho$, $\rho_c\in{\rm Rep}_h(A)$,
\begin{equation}\label{computing_conjugates2}
r_\rho=d_\rho^*=\Omega^{-1}\sum_{i=1}^n\mu_i\overline{e_i}\otimes\alpha^*e_i, \quad\quad
\overline{r}_\rho=b_\rho=\sum_{i=1}^n\beta e_i\otimes\mu_i\overline{e_i}
\end{equation}
with $e_i$ a basis of the space of $\rho$ satisfying $(e_i, e_j)=\delta_{i,j}\mu_i$ and $\mu_i=\pm1$.
Let us consider the case of ${\rm Rep}^+(A)$, so  
$\mu_i=1$. Then it follows by a straightforward computation that
$r_\rho=\Omega^{-1}\sum_i \overline{e_i}\otimes\alpha^*e_i=\sum_i\overline{e_i}\otimes(\alpha_\Omega)^*e_i$,
$\overline{r}_\rho^*\xi\otimes\overline{\eta}=(\beta_\Omega\eta, \xi)$ and this implies
\begin{equation}
\label{intrinsic_dimension_wqh}
r_\rho^*r_\rho=d_\rho r_\rho={\rm Tr}(\alpha(\alpha_\Omega)^*), \quad\quad \overline{r}_\rho^*\overline{r}_\rho=
\overline{r}_\rho^*b_\rho={\rm Tr}((\beta_\Omega)^*\beta),\end{equation}
where $\alpha_\Omega$ and $\beta_\Omega$ are defined in (\ref{twisted_antipode}).
When  $\Omega=\Delta(I)$ is  trivial
then a   computation shows that $\alpha_\Omega=\alpha$ and $\beta_\Omega=\beta$.  
If $\alpha$, $\beta$  are in addition unitary   then  the intrinsic dimensions    coincide with   the vector space dimensions. 
  In Sect.  \ref{20} we shall discuss   examples of $\Omega$-involutive weak quasi-Hopf algebras $A= A({\mathfrak g}, q, \ell)$  arising from a certain semisimplified quotient category associated to quantum groups at roots of unity $U_q({\mathfrak g})$. In this case the antipode is of Kac type  but $\Omega$ is non-trivial,   compatibly with non-integrality of intrinsic dimensions.

  \end{ex}

 We next construct a natural solution of the conjugate equations
 for objects of ${\rm Rep}^+(A)$, with $A$ a unitary discrete  weak  Hopf algebra 
 not necessarily of Kac type. Our methods extend those of 
  \cite{Wor, ER, Van_Daele_multiplier} for the case of discrete or compact quantum groups.
We first establish existence of a Haar  element.

  \begin{prop} Let $A$ be a discrete $\Omega$-involutive weak quasi-bialgebra.
  There is a unique nonzero selfadjoint idempotent $h\in A$ such that $ah=ha=\varepsilon(a)h$ for all $a\in A$.
  \end{prop}
  
  \begin{proof} The proof is as in Prop. 3.1 in \cite{Van_Daele_multiplier}.
  The counit is an irreducible $^*$-representation of $A$. As such, it coincides with the projection onto one of its one dimensional matrix subalgebras. The idempotent defining  this component is the desired element $h$.  
  \end{proof}
  
\begin{defn}  The element $h$ is called the {\it Haar element}.\end{defn}

The following lemma extends a known  idea in the framework of   coassociative quantum groups which, to our knowledge, 
dates back  to \cite{Wor}. Here we consider a modification due to non-triviality of the associator, where the need 
of  the special form that the associator takes for weak  Hopf algebras is apparent. We are not aware of validity of an analogous lemma in a general quasi-coassociative framework.
   
  \begin{lemma}\label{nondegeneracy2}
  Let $A$ be a discrete weak  Hopf algebra. Then $$\Delta(I)A\otimes A=\Delta(A)I\otimes A,\quad\quad A\otimes A\Delta(I)=A\otimes I\Delta(A).$$
  \end{lemma}
  
 \begin{proof} We write $\Delta(I)=a\otimes b$ and for a generic $x\in A$, $\Delta(x)=x_1\otimes x_2$.
Consider the linear map $T: A\otimes A\to A\otimes A$ defined by $T(x\otimes y)=xy_1\otimes y_2$. We show that $T$ is surjective, and this gives the second stated relation. It is straightforward to see that $T$ coincides with the map $\tilde{T}: A\otimes A\to A\otimes A$ defined by $\tilde{T}(x\otimes y)=xS(a_1)a_2y_1\otimes by_2$. Consider also the map $R$ given by $R(x\otimes y)=xS(y_1)\otimes y_2$.
We have
$$TR(x\otimes y)=T(xS(y_1)\otimes y_2)=\tilde{T}(xS(y_1)\otimes y_2)=xS(a_1y_1)a_2y_{2,1}\otimes by_{2,2}.$$
We use the associativity relation $a_1y_1\otimes a_2y_{2,1}\otimes by_{2,2}=y_{1,1}a\otimes y_{1,2}b_1\otimes y_2b_2$
and get after a brief computation
$T(xS(y_1)\otimes y_2)=[x\otimes y][S(a)b_1\otimes b_2]$. A slight modification of this idea gives $T(xS(y_1)\otimes y_2b')=\tilde{T}(xS(y_1)\otimes y_2b')=[x\otimes y][S(a)b_1b_1'\otimes y b_2b_2']$. We replace $x$ by $\tilde{x}=xS(a_2')$ and $y$ by $\tilde{y}=yS(a_1')$ and obtain 
$TR(\tilde{x}\otimes \tilde{y})=[x\otimes y] f$, where the element $f$ was defined in Prop.  \ref{anticom} for general weak quasi-Hopf algebras  and considered again in Prop. \ref{squared_antipode} for weak  Hopf algebras. Since $f$ is partially invertible with domain $\Delta(I)$, the proof is complete. The first relation can be  proved in a similar way  with the maps $T'(x\otimes y)=x_1\otimes x_2y$ as $R'(x\otimes y)=x_1\otimes S(x_2)y$.
 
 \end{proof}    
 
 The following relations extend Prop. 4.1  of \cite{Van_Daele_multiplier} to our setting.
 
 \begin{prop}\label{Haar_element}
  Let $A$ be a discrete weak  Hopf algebra. For all $x$, $y\in A$ we have
  $$\Delta(h) x\otimes y=\Delta(h)I\otimes S(x)y, \quad\quad x\otimes y\Delta(h)=xS(y)\otimes I\Delta(h).$$
 \end{prop}
 
 \begin{proof} We only show the first relation.
 We write $\Delta(I)x\otimes y$ as a finite sum of elements of the form $\Delta(p) I\otimes q$, thanks to the first relation of Lemma \ref{nondegeneracy2}.
 Evaluating $m\circ S\otimes 1$ on this element  gives $S(x)y=\varepsilon(p)q$. On the other hand 
 $$\Delta(h)x\otimes y=\Delta(hp)I\otimes q=\Delta(h\varepsilon(p))I\otimes q=\Delta(h)I\otimes\varepsilon(p)q=\Delta(h)I\otimes S(x)y,$$ and the   relation follows.
 \end{proof}
 
 The following result gives a canonical implementing element for the squared antipode. We omit the proof as it equals that of Prop. 4.3 in \cite{Van_Daele_multiplier}.
 For every full matrix subalgebra $M_{r}({\mathbb C})$ we let  $e_r$ denote its identity, regarded as a central projection of $A$, ${\rm        Tr}_r$   the trace map   which takes value $1$ on the minimal idempotents, and  $r'$   the unique index such that $S(M_{r}({\mathbb C}))=M_{r'}({\mathbb C})$, which is the same as $S(M_{r'}({\mathbb C}))=M_{r}({\mathbb C})$.
 
 \begin{prop} Let $A$ be a discrete weak  Hopf algebra. 
 Then
  $$S^2(x)=Kx{K}^{-1}$$ for all $x\in A$, where $K=(K_r)\in M(A)$ is given by $K_r=[{\rm        Tr}_{r'}\otimes 1(\Delta(h))]^{-1}\in M_{r}({\mathbb C})$.
 \end{prop}

  \begin{thm}
  Let $A$ be an unitary discrete weak  Hopf algebra. Then for every $C^*$-representation $\rho$, the invertible operator $F_\rho:= \rho(K\omega^*)$  is positive. Therefore $\rho_c$ becomes a conjugate of $\rho$ in ${\rm        Rep}^+(A)$ with   inner product $(\overline{\xi}, \overline{\eta})=(\eta, F_\rho\xi)$.
  \end{thm}
  
  \begin{proof}
    It suffices to show positivity of $F_\rho$ for   the   $C^*$-representations $\rho_r$ which project onto the matrix algebras
  $M_{r}({\mathbb C})$, since any other $\rho$ is unitarily equivalent to a direct sum of them.
We   note that $\Delta(h)\Omega^{-1}$ is positive in $M(A\otimes A)$, as $$\Delta(h)\Omega^{-1}=\Delta(h)^2\Omega^{-1}=\Delta(h)\Omega^{-1}\Delta(h)^*.$$
Hence $\Delta(h)\Omega^{-1}e_{r'}\otimes e_r$ is positive as well. Using the notation $\Omega^{-1}=x\otimes y$,
we have, thanks to Prop. \ref{Haar_element},  $$\Delta(h)\Omega^{-1}e_{r'}\otimes e_r=\Delta(h)I\otimes S(xe_{r'})ye_r=\Delta(h)I\otimes S(x)ye_r=\Delta(h)I\otimes \omega_r,$$ with $\omega_r$ the component of $\omega$ along $M_{r}({\mathbb C})$. Evaluating the positive map
${\rm        Tr}_{r'}\otimes 1$ on this element we see that $K_r^{-1}\omega_r=\rho_r(K^{-1}\omega)$ is positive. Hence 
$\rho_r(K\omega^*)=\omega_r\rho_r(K^{-1}\omega)^{-1}\omega_r^*$ is positive as well.
  \end{proof}

\section{Turning $C^*$-categories into tensor $C^*$-categories, I}\label{12}

The problem of constructing unitary tensor categories is of great importance   in  connection with the study of fusion categories  from quantum groups at roots of unity or conformal field theory. In the former setting,  a natural $^*$-structure was introduced by Kirillov \cite{Kirillov} for certain even roots of unity, and unitarity was shown by Wenzl and Xu \cite{Wenzl, Xu_star}.
A tensor category is called {\it unitarizable} if it is tensor equivalent to a tensor $C^*$-category. We have   observed in Remark \ref{non-unitarizability} that
examples of non-unitarizable fusion categories from quantum groups and certain roots of unity are known.

We start with the following setting, which will be called condition a).

\medskip

a) Let     ${\mathcal C}$ be a tensor category (resp. a braided tensor category) and ${\mathcal C}^+$ a $C^*$-category, and assume that we have an equivalence of linear categories $$ {\mathcal F}: {\mathcal C}^+\to{\mathcal C}.$$ 
  Assuming that every object of ${\mathcal C}^+$ is a finite direct sum of irreducible objects, we wish to use ${\mathcal F}$ to upgrade ${\mathcal C}^+$ to a tensor $C^*$-category
   (resp. a braided tensor $C^*$-category with unitary braided symmetry). 

In this section we discuss   a result which characterizes when a solution exists and is unique. 

We shall derive three main applications. The first one is the proof that Beilinson-Feigin-Mazur category $\tilde{\mathcal O}_\ell$
is a   braided tensor $C^*$-category with unitary braided symmetry, and is presented in Theorem \ref{BFM}. The second
is the explicit construction of a similar structure on the module category of the affine vertex operator algebra 
$V_{{\mathfrak g}_k}$ endowed with Huang-Lepowsky braided tensor category structure, see  Theorem \ref{TheoremUnitaryBraidRepAffine}.
A third application a positive answer to a question posed by Galindo in \cite{Gal} on uniqueness of
unitary tensor structures on tensor categories, Theorem \ref{Galindo}.

A  variant of the abstract results of this section will be useful in Sect. \ref{18}, \ref{19}, \ref{20} were unitary weak quasi-Hopf algebras will be constructed with a direct method from
the braiding for certain general ribbon categories and in particular for those arising from the quantum groups at roots of unity.
\bigskip

\begin{defn} Let ${\mathcal F}: {\mathcal C}^+\to{\mathcal C}$ satisfy a). We shall say that   the tensor structure (resp. the braided tensor structure) of 
${\mathcal C}$ is {\it transportable compatibly with the $C^*$-structure}, or simply {\it $C^*$-transportable} to ${\mathcal C}^+$ if 
   ${\mathcal C}^+$ can be upgraded  to a   tensor $C^*$-category (resp. a braided tensor $C^*$-category with unitary braided symmetry) in such a way that ${\mathcal F}:  {\mathcal C}^+\to{\mathcal C}$ becomes a tensor equivalence (resp. a braided tensor equivalence).
\end{defn}

We note that $C^*$-transportability will be possible only in certain circumstances. For example, 
 if ${\mathcal C}$
  is a finite semisimple tensor category then we know that ${\mathcal C}$ is tensor equivalent to some ${\rm Rep}(A)$, with $A$ a semisimple weak quasi bialgebra. Since $A$ admits the structure of a $C^*$-algebra, the $C^*$-category ${\mathcal C}^+$ of $C^*$-representations of $A$ satisfies a). On the other hand, if ${\mathcal C}$ is not tensor equivalent to   a tensor $C^*$-category, see Remark \ref{non-unitarizability}, then  ${\mathcal C}^+$ does not admit any tensor $C^*$-structure that makes it tensor equivalent to ${\mathcal C}$. 

We shall describe two main classes of tensor categories for which the  tensor structure (resp. the braided tensor structure)  are transportable compatibly with the $C^*$-structure, and  two   upgrading of  ${\mathcal C}^+$ corresponding to a $C^*$-transportable 
tensor structure of ${\mathcal C}$ provide unitarily (braided) tensor  equivalent tensor $C^*$-categories. The notion of weak quasi-Hopf algebra will play a prominent role. 

In the mentioned applications, ${\mathcal C}$   plays the role of a   category of infinite dimensional representations of interest of some algebraic structure  endowed with a `fusion' tensor product, 
and ${\mathcal C}^+$  the category of unitary   representations on Hilbert spaces.   The functor ${\mathcal F}: {\mathcal C}^+\to{\mathcal C}$   is the functor which forgets the unitary structure. The assumption that it be an equivalence   means that every object of ${\mathcal C}$ can be made into  a unitary representation, an assumption which is known to hold in a variety of circumstances.

 We thus see   that the abstract problem in our formulation admits as applications
  that of making representation categories of VOAs or of affine Lie which are known to have the structure of
  linear $C^*$-categories, and the structure of (braided) tensor categories,   into   tensor $C^*$-categories (with unitary braided symmetry).

\begin{ex}
An  example illuminating the general methodology, and useful for our application,  is obtained with ${\mathcal C}$ the Andersen fusion category of a quantum group $U_q({\mathfrak g})$ associated to a complex simple Lie algebra, with   $q=e^{i\pi/\ell}$ with $\ell$ divisible by $d$,  the ratio between the squared lengths of the longest to the shortest root, and $\ell/d>\check{h}$, the dual Coxeter number of ${\mathfrak g}$.

 In this case, Wenzl shows in the first part of \cite{Wenzl}   that ${\mathcal C}$ is a linear $C^*$-category in a natural way. 
 In the second part, that ${\mathcal C}$ is in fact a tensor $C^*$-category with unitary braided symmetry. The second part is based on the construction of a positive inner product of a fusion tensor product of two representations satisfying the extra axioms, and compatible with the starting linear $C^*$-structure of the category.
  \end{ex}

A note on notation. Since we shall deal at the same time     with  semisimple linear or $C^*$ or  tensor categories, or braided tensor categories,
and sometimes we shall use only part of the structure,
 we shall use 
a suffix $^+$ on a category to denote that it is a $C^*$-category and on a functor   if it is $^*$-preserving.
  The continuous arrow in a diagram denotes a (braided) tensor equivalence, while the     dashed arrows are linear equivalences.
  Thus, if only some of the categories or equivalences carry a (braided) tensor structure, the commutativity of the diagram is understood to hold at the level of linear functors. 
  
\begin{defn}\label{CompatibleTriple}
Let  ${\mathcal F}: {\mathcal C}^+\to{\mathcal C}$ be as in a). Let $A$ be a   discrete (quasitriangular) weak quasi bialgebra    endowed with an involution of    pre-$C^*$-algebra, and consider, accordingly, the (braided) tensor category ${\rm        Rep}(A)$
and the $C^*$-category ${\rm        Rep}^+(A)$. A triple $(A, {\mathcal E}^+, {\mathcal E})$ constituted by
   a $^*$-equivalence  
  ${\mathcal E}^+: {\mathcal C}^+\to{\rm        Rep}^+(A)$ and a (braided) tensor equivalence
 ${\mathcal E}: {\mathcal C}\to{\rm        Rep}(A)$ will be called {\it compatible} with ${\mathcal F}$ if the following diagram commutes
 
 \[
\begin{tikzcd}
{\mathcal C}^+\arrow[dashrightarrow] {rr}{{\mathcal E}^+} \arrow[dashrightarrow]{d} {{\mathcal F}}  &&   {{\rm        Rep}^+(A) } \arrow[dashrightarrow]{d}{{\mathcal F}_A}  \\
{\mathcal C}  \arrow{rr}{{{\mathcal E}}} &&   {{{\rm        Rep}(A) } }
\end{tikzcd}
\] where ${\mathcal F}_A: {\rm Rep}^+(A)\to{\rm Rep}(A)$ is the forgetful functor.
  \end{defn}
   
    A compatible triple defines a weak dimension function on ${\mathcal C}$ via $D(\rho):={\rm dim}({\mathcal E}'(\rho))$,
    where ${\mathcal E}'$ is the composite of ${\mathcal E}$ with the forgetful functor ${\rm Rep}(A)\to{\rm Vec}$.
 We next see that   compatible triples may be constructed and classified under   mild assumptions.

 \begin{prop}\label{twisted_compatible_triples} If $(A, {\mathcal E}, {\mathcal E}^+)$ is a compatible triple for ${\mathcal F}:{\mathcal C}^+\to{\mathcal C}$ then for any twist $F\in A\otimes A$ of the (quasitriangular) weak quasi bialgebra structure and any positive twist $t\in A$ of the $^*$-involution, the twisted (quasitriangular) algebra $A_{F, t}$ is part of another compatible triple with the same weak dimension function and they are all of this form.\end{prop}
 
 \begin{proof}
The proof follows from Prop.
  \ref{class}, Theorem \ref{propweakdim}  and part of  Prop. \ref{1_twist_equivalence}. The braided case follows from
  considering part c) of Prop. \ref{R-canonicity} in addition.
  \end{proof}

 \begin{rem}\label{equivalent_compatible_triple} As we shall see,  natural constructions in conformal field theory  give rise
 to canonically associated associative   algebras $A$, the Zhu algebras,   and also to  {\it linear} functors ${\mathcal E}$,
 which are already known to play  an important role in the theory of VOAs.
The construction of compatible triples for these remarkable examples  is our main motivation 
  in the definition.     \end{rem}

  \begin{prop} \label{equivalent_compatible_triple2}
  ${\mathcal F}:{\mathcal C}^+\to{\mathcal C}$ admits a compatible triple if and only if ${\mathcal C}$ admits an integral weak dimension function.
  \end{prop}
  
  \begin{proof}
The notion of a compatible triple   $(A, {\mathcal E}, {\mathcal E}^+)$ may equivalently be given via an abstract construction as follows.
There is a canonical isomorphism of algebras $\phi: A\to {\rm Nat}_0({\mathcal E}')$ 
which induces an isomorphism of categories $\phi_*:{\rm Rep}({\rm Nat}_0({\mathcal E}'))\to{\rm Rep}(A)$
  such that $\phi_*\widetilde{{\mathcal E}'}={\mathcal E}$, where $\widetilde{{\mathcal E}'}:{\mathcal C}\to{\rm Rep}({\rm Nat}_0({\mathcal E}'))$ is the   equivalence arising from Tannaka-Krein reconstruction of ${\mathcal E}'$.
 There is also an isomorphism of $^*$-algebras  $A\to {\rm Nat}_0({{\mathcal E}^+}')$. The compatibility condition implies
  ${\mathcal H}({\mathcal E}^+)'={\mathcal E}'{\mathcal F}$, with ${\mathcal H}:{\rm Hilb}\to{\rm Vec}$ the forgetful functor. These remarks together with Tannaka-Krein duality results  imply that   giving a compatible triple is the same thing as giving a faithful $^*$-functor ${{\mathcal E}^+}':{\mathcal C}^+\to{\rm Hilb}$ and a faithful weak quasi-tensor functor ${\mathcal E}':{\mathcal C}\to{\rm Vec}$ such that  ${\mathcal E}'{\mathcal F}={\mathcal H}{{\mathcal E}^+}'$.
Now it suffices to apply Theorem \ref{propweakdim}, \ref{Haring_functor_$C^*$}.

  \end{proof}

\begin{thm}\label{transportability} Let  ${\mathcal F}: {\mathcal C}^+\to{\mathcal C}$ satisfy a).  Assume that the (braided) tensor category ${\mathcal C}$ admits a weak dimension function $D$, and let   $(A, {\mathcal E}^+, {\mathcal E})$ be a compatible triple with dimension $D$. Then
the tensor structure of ${\mathcal C}$ is $C^*$-transportable to ${\mathcal C}^+$ via ${\mathcal F}$ if and only if
  $A$ can be upgraded to a unitary (resp. quasitriangular)  weak quasi bialgebra (resp. such that the $R$-matrix of $A$ satisfies the unitarity condition stated in Prop. \ref{unitary_braided_symmetry2}). If this is the case, the diagram defining the triple becomes a commuting diagram of (braided) tensor equivalences and  ${\mathcal E}^+$  can be chosen unitary. Furthermore, any two tensor  $C^*$-completions of ${\mathcal C}^+$ obtained from a  $C^*$-transportable ${\mathcal F}$   yield  unitary     tensor equivalent tensor $C^*$-categories.
\end{thm} 

\begin{proof} Following the proof of Prop. \ref{equivalent_compatible_triple2}, and adopting the same notation, we shall identify $A$ with 
${\rm Nat}_0({{\mathcal E}^+}')$ as a $^*$-algebra.
 If ${\mathcal C}^+$ admits the structure of a tensor $C^*$-category over the underlying $C^*$-category such that ${\mathcal F}: {\mathcal C}^+\to{\mathcal C}$ becomes a tensor equivalence then  the composite of the left with the bottom equivalences in the diagram is a tensor equivalence
hence,  by commutativity of the diagram, the composite of  top with the right equivalences ${\mathcal C}^+\to{\rm        Rep}^+(A)\to{\rm        Rep}(A)$ is a tensor equivalence as well. On the other hand,  ${\rm        Rep}(A)\to{\rm        Vec}$ is a weak quasi-tensor functor,
hence so is the composite ${\mathcal C}^+\to{\rm        Rep}^+(A)\to{\rm        Rep}(A)\to{\rm        Vec}$. 
But this functor factors through ${\mathcal C}^+\to{\rm        Rep}^+(A)\to{\rm        Hilb}\to{\rm        Vec}$ and ${\rm        Hilb}\to{\rm        Vec}$
is both a forgetful functor and a tensor equivalence, and this implies that ${{{\mathcal E}^+}': \mathcal C}^+\to{\rm        Rep}^+(A)\to{\rm        Hilb}$ is a ($^*$-preserving) weak quasi-tensor functor. It follows that 
${\rm        Nat}_0({{\mathcal E}^+}')$ can be made into a unitary weak quasi  bialgebra and ${\mathcal E}^+$ into a unitary tensor equivalence by Theorem \ref{TheoremTannakaStar}. This structure can be transferred to $A$, and therefore is compatible with the given $^*$-involution of $A$.  It  is now easy to see that it extends the given weak quasi-bialgebra structure on $A$.

Conversely, if  $A$ admits  the structure of a unitary weak quasi  bialgebra with the given $^*$-structure then, by Corollary \ref{CorollaryRep^+(A)},  ${\rm        Rep}^+(A)$ is a tensor $C^*$-category tensor equivalent to ${\rm Rep}(A)$ and hence to ${\mathcal C}$.   The top equivalence of the  diagram defining a compatible triple acts from the linear category ${\mathcal C}^+$ to the 
tensor category ${\rm Rep}^+(A)$. It is a general fact that under this circumstance, ${\mathcal C}^+$ can be made into a tensor category
in such a way that ${\mathcal E}^+$ is becomes a tensor equivalence. Indeed, given objects $\rho$, $\sigma\in{\mathcal C}^+$, we define a tensor product object  $\rho\otimes\sigma$ in ${\mathcal C}^+$ 
$$ \rho \otimes \sigma :=  \mathcal{S}^+ \left( \mathcal{E}^+(\rho) \otimes \mathcal{E}^+(\sigma) \right) \,,$$
and a tensor product morphism by a similar formula,
$$S\otimes T:=  \mathcal{S}^+ \left( \mathcal{E}^+(S) \otimes \mathcal{E}^+(T) \right) .$$
  Here ${\mathcal S}^+:{\rm Rep}^+(A)\to{\mathcal C}^+$is an inverse equivalence of ${\mathcal E}^+$,
Moreover, if $\alpha$ denotes the unitary associator in ${\rm        Rep}^+(A)$  we define 
the morphisms 
$$\alpha'_{\rho,\sigma,\tau } : ( \rho \otimes \sigma) \otimes \tau \to   
\rho \otimes ( \sigma \otimes \tau ) $$
by 
\begin{eqnarray*}
&&\alpha'_{\rho,\sigma,\tau} := \\
&&\mathcal{S^+}(1_{\mathcal{E}^+(\rho)}\otimes\eta^{-1}_{\mathcal{E}^+(\sigma)\otimes \mathcal{E}^+ (\tau)}  \circ \alpha_{\mathcal{E}^+(\rho), \mathcal{E}^+ (\sigma), \mathcal{E}^+   (\tau)}   \circ 
\eta_{\mathcal{E}^+(\rho)\otimes \mathcal{E}^+(\sigma)}\otimes 1_{\mathcal{E}^+(\tau)}) \, . 
\end{eqnarray*}
where $\eta:{\mathcal E}^+{\mathcal S}^+\to 1$ is a   natural transformation.
 Then, thanks to the fact that ${\mathcal E}^+$ is $^*$-preserving, ${\mathcal S}^+$ may be chosen $^*$-preserving, 
  $\eta$ unitary
 by Prop. \ref{MacLane},
   and ${\rm        Rep}^+(A)$ is a tensor $C^*$-category,  it is immediate to check that  that the  relation $(S\otimes T)^*=S^*\otimes T^*$ holds on morphisms and $\alpha'$ is unitary.   This gives the  $C^*$- tensor structure on ${\mathcal C}^+$.

   Moreover,
${\mathcal E}^+$ becomes a tensor equivalence with  unitary tensor structure 
$E_{\rho, \sigma}:=\eta^{-1}_{{\mathcal E}^+(\rho)\otimes {\mathcal E}^+(\sigma)}$.
Since the forgetful ${\rm        Rep}^+(A)\to{\rm        Rep}(A)$  is a tensor equivalence as well with the trivial tensor structure, it follows that ${\mathcal F}: {\mathcal C}^+\to{\mathcal C}$ has a unique tensor 
structure such that ${\mathcal E}{\mathcal F}={\mathcal F}_A{\mathcal E}^+$ as tensor functors.  

If in addition $A$ is quasitriangular with $R$.matrix satisfying the stated unitarity condition then ${\rm Rep}^+(A)$ is a braided tensor category with unitary braided symmetry, let it be $c$.
We define a unitary   natural transformation $c'$ in ${\mathcal C}^+$  
$$c'(\rho, \sigma)={\mathcal S}^+(c({\mathcal{E}}^+(\rho), {\mathcal{E}}^+(\sigma))).$$
It is straightforward to verify that this natural transformation satisfies the two hexagon equations with respect to $\alpha'$.

Uniqueness. Let us next   consider a new      tensor $C^*$-category ${\mathcal C}'$
coinciding with ${\mathcal C}^+$ as a $C^*$-category and   making ${\mathcal F}$ into a new tensor equivalence ${\mathcal F}'$.
Applying the above construction in the opposite direction, that is with ${\mathcal S}^+$ in place of ${\mathcal E}^+$, gives a new tensor $C^*$-category structure to ${\rm Rep}^+(A)$, denoted ${\rm Rep}'(A)$
and  new unitary tensor equivalences ${\mathcal S}': {\rm Rep}'(A)\to{\mathcal C}'$ and ${\mathcal E}': {\mathcal C}'\to{\rm Rep}'(A)$ coinciding with ${\mathcal S}^+$ and ${\mathcal E}^+$ as functors, respectively.
We obtain    a new tensor structure on
the identity functor ${\mathcal F}_A: {\rm Rep}'(A)\to{\rm Rep}(A)$   solving now the equation 
for the tensor structures obtained from  
 ${\mathcal E}{\mathcal F}'={\mathcal F}_A{\mathcal E}'$. This gives a weak quasi-tensor structure to the forgetful functor
${\rm Rep}'(A)\to{\rm Hilb}$, and therefore   ${\rm Rep}'(A)$ becomes unitarily tensorially equivalent to ${\rm Rep}^+(A')$
where $A'$ is a new unitary   weak quasi-bialgebra compatible with the original $C^*$-algebra $A$,
thanks to Theorem \ref{TheoremTannakaStar}.
 It follows from  Prop. \ref{class} that $A'$ as a weak quasi-bialgebra  is only varying by a twist of $A$. Therefore
 ${\rm Rep}^+(A')$ is unitarily tensor equivalent to ${\rm Rep}^+(A)$   by Prop. \ref{unitary_equivalence_under_twisting}
 and \ref{uniqueness}, and finally to ${\mathcal C}^+$.

\end{proof}

It follows in particular from the previous characterization that if the tensor structure of ${\mathcal C}$ is $C^*$-transportable to $ {\mathcal C}^+$ then ${\mathcal C}$ is tensor equivalent to a tensor $C^*$-category, namely ${\rm Rep}^+(A)$.
We next show more interestingly that the converse  implication holds. The following result will find important applications
in  the categories arising from affine vertex operator algebras,  Sect. \ref{VOAnets2}.

\begin{thm}\label{unitarizability} Let ${\mathcal F}:{\mathcal C}^+\to{\mathcal C}$ satisfy a) and assume that ${\mathcal C}$ admits a weak dimension function (e.g. ${\mathcal C}$ is a finite semisimple  (braided)  tensor category). If
  ${\mathcal C}$ is (braided) tensor equivalent to a (unitarily braided) tensor $C^*$-category ${\mathcal D}^+$,   then the (braided) tensor structure of ${\mathcal C}$ is $C^*$-transportable to $ {\mathcal C}^+$ in a unique way up to  unitary (braided) tensor equivalence.
  Moreover in this way ${\mathcal C}^+$ becomes unitarily (braided) tensor equivalent to ${\mathcal D}^+$.
  \end{thm}
  
\begin{proof}
Let $D$ be a weak dimension function on ${\mathcal C}$,   
and ${\mathcal G}:  {\mathcal D}^+\to{\mathcal C}$ a tensor equivalence. Then $D'(\rho):=D({\mathcal G}(\rho))$ is a weak dimension function on ${\mathcal D}^+$ since   ${\mathcal G}(\rho\otimes\sigma)$ is isomorphic to ${\mathcal G}(\rho)\otimes{\mathcal G}(\sigma)$ and $D$ is isomorphism invariant. We may then construct a faithful $^*$-functor of $C^*$-categories ${\mathcal D}^+\to{\rm Hilb}$ corresponding to $D'$ and a weak quasi-tensor structure on the composite
${\mathcal D}^+\to{\rm Hilb}\to{\rm Vec}$. By the $C^*$-version of Tannaka-Krein duality, see Theorem \ref{TheoremTannakaStar},
the algebra $A$ of natural transformations of this functor becomes a unitary   weak quasi  bialgebra, with a  corresponding involutive structure  $(^*, \Omega)$ and
such that ${\rm Rep}^+(A)$ is unitarily tensor equivalent to ${\mathcal D}^+$.
Let ${\mathcal G}':{\mathcal C}\to{\mathcal D}^+$ be an inverse tensor equivalence of ${\mathcal G}$ and let
${\mathcal E}$ be the composed tensor equivalence ${\mathcal C}\to{\mathcal D}^+\to {\rm Rep}(A)$
where the latter functor is obtained from the   duality theorem in the tensor linear case, see Theorem \ref{TK_algebraic_quasi} (or equivalently, forgetting the $C^*$-structure of $A$). We may then pick a factorisation of ${\mathcal E}{\mathcal F}$ through a $^*$-equivalence
${\mathcal E}^+: {\mathcal C}^+\to{\rm Rep}^+(A)$ and the forgetful functor ${\rm Rep}^+(A)\to{\rm Rep}(A)$ by Prop. \ref{Haring_functor_$C^*$}. Let $^\dag$ denote the corresponding involution on $A$. Since all pre-$C^*$-algebra involutions of $A$ are twisted from one another, we may find a twist
$t\in A$, positive with respect to $^*$, such that $t^\dag=t^{-1}a^*t$. We may   endow $A$ with the twisted involutive 
structure $(^\dag, \Omega_t)$ by Prop. \ref{1_twist} and obtain the complete structure and an associated tensor $C^*$-category ${\rm Rep}^+_t(A)$. We have thus shown that $(A, {\mathcal E}^+, {\mathcal E})$ is a compatible triple for ${\mathcal F}$ satisfying the necessary and sufficient condition 
 of Theorem \ref{transportability}
of $C^*$-transportability. Thus   ${\mathcal C}^+$ becomes a tensor $C^*$-category unitarily tensor equivalent to 
${\rm Rep}^+_t(A)$ and therefore to ${\mathcal D}^+$  by Prop. \ref{1_twist_equivalence}.
In the special case that ${\mathcal C}$ is a finite semisimple tensor category, it always admits a weak dimension function 
by Remark   \ref{rational}. In the braided case, the same $R$-matrix satisfies the unitarity condition stated in Prop. \ref{unitary_braided_symmetry2} with respect to the twisted involutive structure $(^\dag, \Omega_t)$, and therefore it
makes ${\rm Rep}^+_t(A)$ into a unitarily braided tensor category. We may apply the braided case of Theorem \ref{transportability}. 
\end{proof}

 \bigskip

 \section{Positive weak dimension and amenability}\label{13}

The Grothendieck ring   ${\rm Gr}({\mathcal C})$ of a  rigid semisimple tensor category ${\mathcal C}$ is called amenable if it admits a dimension function satisfying a certain analytic property. 
Such a function, called amenable,  is unique and    bounds below any other dimension function, see e.g. \cite{CQGRC}.
 In this section 
  we extend  the framework to  weak dimension functions.  
We show that the
   amenable    dimension function   is already unique among weak dimension functions and minimizes them. 
This gives a weaker criterion for amenability.
 It follows in particular that  if ${\mathcal C}$ is a fusion category the lower bound of weak dimension functions is given by the Frobenius-Perron dimension, and this   was our original motivation for the study of amenability. 
    \medskip

 Let ${\mathcal C}$ be a rigid semisimple tensor category and $D$ a weak dimension
 function on the Grothendieck ring   ${\rm Gr}({\mathcal C})$, see Def. \ref{weak_dimension_function}, that will always be assumed positive and symmetric in this section. As already mentioned, we  first aim to introduce a notion of amenability for $D$ extending  the usual amenability for a genuine dimension.
To do this, we closely follow the treatment in Sect. 2.7 in \cite{CQGRC}, dropping the unitarity 
assumption on ${\mathcal C}$.
Therefore for $\rho\in{\rm Irr}({\mathcal C})$ let $\Lambda_\rho$ be the operator of left multiplication by $\rho$
on the complexified algebra ${\rm Gr}_{\mathbb C}({\mathcal C}):={\rm Gr}({\mathcal C})\otimes_{\mathbb Z}{\mathbb C}$.
 It follows from associativity of ${\rm Gr}({\mathcal C})$ that 
  \begin{equation}\label{representation}
 \Lambda_\rho\Lambda_\sigma=\sum_\tau m^\tau_{\rho, \sigma}\Lambda_\tau,\end{equation}
 with $m^\tau_{\rho,\sigma}={\rm dim}(\tau, \rho\otimes\sigma)$ and therefore $\Lambda$ linearly extends to a representation of 
 ${\rm Gr}({\mathcal C})$. 

\begin{prop}\label{bound_below} Let ${\mathcal C}$ be a rigid semisimple tensor category admitting a weak dimension function.
The operator $\Lambda_\rho$ extends to a bounded linear operator on $\ell^2({\rm Irr}({\mathcal C}))$. We have
$\|\Lambda_\rho\|\leq D(\rho)$ for $\rho\in {\rm Irr}({\mathcal C})$ and for every
 weak   dimension function $D$.
\end{prop}

\begin{proof}
The proof extends the corresponding proof for dimension functions, see Prop. 2.7.4 in \cite{CQGRC}, with the   modification that $u_\sigma=v_\sigma=D(\sigma)$ is replaced by $u_\sigma=D(\sigma)$ and
$v_\sigma=(\Gamma u)_\sigma=\frac{D(\overline{\rho}\sigma)}{D(\rho)}\leq D(\sigma)$ which implies 
$\Gamma^t(v)_\sigma\leq\frac{D(\rho\sigma)}{D(\rho)}\leq D(\sigma)=u_\sigma$ and in turn replaces 
$\Gamma^t(v)=u$. Note indeed that these modifications are still compatible with Lemma 2.7.3 in  \cite{CQGRC}
and the proof may be completed.
\end{proof}

Given a dimension function $D$ we consider operators $\lambda_\mu=\sum_{\rho\in {\rm Irr}({\mathcal C})}\frac{\mu(\rho)}{D(\rho)}\Lambda_\rho$ associated to probability measures $\mu$
on ${\rm Irr}({\mathcal C})$ and then we find that a composition $\lambda_\mu\lambda_\nu=\lambda_{\mu*\nu}$,
with $\mu*\nu$ the convolution measure defined as at page 71 in \cite{CQGRC},
$$\mu*\nu(\tau)=\sum_{\rho,\sigma\in {\rm Irr}({\mathcal C})}m^\tau_{\rho,\sigma}\frac{D(\tau)}{D(\rho)D(\sigma)}\mu(\rho)\nu(\sigma),$$ with $m^\tau_{\rho,\sigma}$ the multiplicity of $\tau$ in $\rho\otimes\sigma$.
  For a weak dimension function a similar formula holds but $\mu*\nu$ may not be a probability measure.
  Indeed    $\|\mu*\nu\|=\sum_{\tau\in {\rm Irr}({\mathcal C})}\mu*\nu(\tau)=\sum_{\sigma,\tau\in  {\rm Irr}({\mathcal C})}\frac{\mu(\sigma)\nu(\tau)D(\sigma\tau)}{D(\sigma)D(\tau)}\leq 1$. Thus if 
    ${\rm Irr}({\mathcal C})$ is countable and if $\mu$ and $\nu$ have support 
  ${\rm Irr}({\mathcal C})$ then
   $\mu*\nu$ is a probability measure precisely when $D$ is a genuine dimension function.
  Therefore we more generally consider the operators $\lambda_\mu$ for any positive measure $\mu$ with $\|\mu\|\leq1$.
One has $\|\lambda_\mu\|\leq\|\mu\|$, so $\|\lambda_\mu\|=1$ is possible only if $\mu$ is a probability measure. 

\begin{prop}\label{amenability}
Let $D$ be  a   weak dimension function on ${\rm Gr}({\mathcal C})$. Then the following properties are equivalent.
 \begin{itemize}
\item[(a)] $1\in{\rm Sp}\lambda_\mu$ for every probability measure $\mu$,
\item[(b)] $\|\lambda_\mu\|=1$  for every probability measure $\mu$,
\item[(c)] $(\check{\mu}*\mu)^n(\iota)^{1/n}\to1$ for every probability measure $\mu$, with $\check{\mu}(\rho)=\mu(\overline{\rho})$,
\item[(d)] there is a net $\xi_\alpha\in\ell^2({\rm Irr}({\mathcal C}))$ of positive unit vectors such that
$\|\Lambda_\rho\xi_\alpha-D(\rho)\xi_\alpha\|\to 0$ for all $\rho\in {\rm Irr}({\mathcal C})$.
\end{itemize}
If they hold then $D$ is a dimension function.
\end{prop}

\begin{proof}
 The equivalence of properties (a)--(d)  may be proven just as in the case of ordinary dimension functions, cf. Lemma 2.7.5 in \cite{CQGRC}, 
 taking into account the slight modifications mentioned before the statement. The last statement   follows
 from  the observation that $\Lambda$ is a representation of ${\rm Gr}({\mathcal C})$ in the sense of 
(\ref{representation}),  and
 a $3\varepsilon$-argument applied to the vanishing net  $(\Lambda_{\rho}(\Lambda_{\sigma}-D(\sigma)))\xi_\alpha$
 with $\xi_\alpha$ as in  (d).
\end{proof}

We recall the definition of amenability.

\begin{defn}
A  dimension function on ${\rm Gr}({\mathcal C})$ satisfying the equivalent properties of Prop. \ref{amenability}  is called amenable.
The category ${\mathcal C}$ is called amenable if the intrinsic dimension function on ${\rm Gr}({\mathcal C})$ is amenable.\end{defn}

 The following    result   
   extends to weak dimension functions the uniqueness result known for an amenable dimension function, see Prop. 2.7.7 in \cite{CQGRC}.

\begin{thm}\label{amenability2} An amenable dimension function on ${\rm Gr}({\mathcal  C})$ is
unique among   weak dimension functions satisfying the equivalent properties of Prop. \ref{amenability} and  is
  given by $D(\rho)=\|\Lambda_\rho\|$ for $\rho\in{\rm Irr}({\mathcal C})$. Any  other weak dimension $D'$ satisfies
$D'(\rho)\geq D(\rho)$ for all $\rho$. 
\end{thm}

\begin{proof}
The first statement   follows from Prop. \ref{amenability} and property (b) applied to the probability measures with support a single irreducible. The second part   follows
 from the first and  Prop.  \ref{bound_below}.

\end{proof}

Existence of an amenable dimension function   is characterized by the following property of
the left regular representation,
for all $\rho\in{\rm Irr}({\mathcal C})$, $\Lambda_\rho$ is bounded and $\|\sum_i\mu_i \Lambda_{\rho_i}\|=\sum_i\mu_i\|\Lambda_{\rho_i}\|$   for finite linear combinations of basis elements with positive coefficients.
These conditions are clearly necessary as by the previous theorem the amenable dimension function is unique and
explicitly given by $\|\Lambda_\rho\|$. Conversely,
when  $\Lambda_\rho$ is bounded, we define the operators $\lambda_\mu$ as before with
 $ \|\Lambda_\rho\|$ in place of $D(\rho)$, $\rho\in{\rm Irr}({\mathcal C})$. Then it is easy to see using
continuity
of  $\mu\in\ell^1({\rm Irr}({\mathcal C}))\to\lambda_\mu\in{\mathcal B}(\ell^2({\rm Irr}({\mathcal C})))$
that the positive linearity of $\|\Lambda_\rho\|$ is equivalent to
property (b) of Prop. \ref{amenability}. It follows that   the linear extension of $\|\Lambda_\rho\|$
is automatically   an amenable weak dimension function by submultiplicativity of the norm.

It follows that    every  fusion category  ${\mathcal C}$ admits the unique amenable dimension function, and moreover has a unique positive dimension function, the Frobenius-Perron 
dimension  determined by  ${\rm FPdim}(\rho)=\|\Lambda_\rho\|$. Indeed, (d) of  Prop.  \ref{amenability},  has a solution given by  the vector with coordinates the dimensions of the simple objects, and by
 Sect. 8 in \cite{ENO} or Chapter 4 in \cite{EGNO}, $\|\Lambda_\rho\|$ is indeed 
  a   dimension function on  ${\rm Gr}({\mathcal C})$.

\begin{cor}
If ${\mathcal C}$ is a fusion category then $D(\rho)\geq {\rm FPdim}(\rho)$ for every weak dimension function $D$
on ${\rm Gr}({\mathcal C})$.
\end{cor}

 Another important class of examples   is that for which  ${\rm Gr}({\mathcal C})$ is commutative.
Yamagami showed that in this case that for the existence of the amenable dimension function  is necessary and sufficient to verify that  $\Lambda_\rho$ is bounded,   see Theorem 3.5 in \cite{Yamagami}.

\begin{rem} The examples that we have studied in the paper show that there may be more than a   natural choice of integral weak dimension functions associated to a fusion category.
For example, 
for the pointed fusion categories arising from quantum groups at roots of unity (or vertex operator algebras) at the minimal root (level),
  ${\rm FPdim}(g)=1$  on every irreducible object $g$, so ${\rm FPdim}$  is already an integral dimension function. Another natural choice is associated  to Wenzl functor or to Zhu's functor. Consider for each level $k$,   ${\rm Gr}({\mathcal C}({\mathfrak g}, q, \ell)))$ for $q=e^{i\pi/\ell}$, $\ell=d(k+\check{h})$
  and regard it as a quotient of the classical representation ring $R({\mathfrak g})$ associated to ${\mathfrak g}$.
  Then  the sequence $D_k$ of weak dimension
  functions on   ${\rm Gr}({\mathcal C}({\mathfrak g}, q, \ell))$ defined by Wenzl's functor defines in the pointwise limit
  the classical dimension function of $R({\mathfrak g})$, which is also the unique amenable dimension function of this based ring.
  \end{rem}
  
  We next apply Theorem \ref{amenability2} to a weak tensor functor between tensor $C^*$--categories studied in Sect. \ref{3} and we find a  useful upper and lower bound for the associated weak dimension function.

\begin{cor}\label{upper_lower_bound}
Let ${\mathcal C}$ and ${\mathcal C}'$ be  rigid tensor $C^*$-categories such that ${\rm Gr}({\mathcal C})$ admits the amenable dimension function $D$. Then every  weak tensor
$^*$-functor ${\mathcal F}: {\mathcal C}\to{\mathcal C}'$ defined by $F$ and $G$  satisfies 
$$D(\rho)\leq d'({\mathcal F}(\rho))\leq \|F_{\overline{\rho},\rho}\|\|G_{\rho, \overline{\rho}}\|d(\rho), \quad \rho\in{\mathcal C},$$ where  $d$, $d'$ are   the intrinsic dimensions of
${\mathcal C}$ and ${\mathcal C}'$ respectively.
 
\end{cor}

\begin{proof}
Note that the  weak dimension function
$\rho\to {\rm dim}({\mathcal F}(\rho))$  is symmetric as
${\mathcal F}(\overline{\rho})$ is a conjugate of ${\mathcal F}(\rho)$
by Prop. \ref{weak_duality}. The lower bound then follows from Theorem \ref{amenability2}. For the upper bound see
Cor. \ref{intrinsic_dimension_bound}.
\end{proof}

  We conclude  with a result concerning a dimension preserving property of unitary weak tensor functor
  between rigid $C^*$-tensor categories in the amenable case. This result extends 
   a known property for unitary tensor functors, see Cor.  2.7.9 in \cite{CQGRC} and references therein.

\begin{cor}\label{unitarity_obstruction}
Let ${\mathcal C}$ and ${\mathcal C}'$ be   rigid tensor $C^*$-categories with intrinsic dimensions $d$ and $d'$ respectively, and let $({\mathcal F}, F, G):{\mathcal C}\to{\mathcal C}'$ be a  unitary weak tensor $^*$-functor.
 If    ${\mathcal C}$ is amenable (e.g. ${\mathcal C}$ is a fusion category) then $d(\rho)=d'({\mathcal F}(\rho))$  for all $\rho$.
In particular, when ${\mathcal C}'={\rm Hilb}$ then $d(\rho)={\rm dim}({\mathcal F}(\rho))$ and therefore
${\mathcal F}$ is already tensor.
 \end{cor}

\begin{proof}  
By assumption $F^*$ and $G$  are isometric, so $\|F_{\rho, \sigma}\|=\|G_{\rho, \sigma}\|=1$. 
 It follows from Cor. 
 \ref{upper_lower_bound} that  $d'({\mathcal F}(\rho))=d(\rho)$
 as  $d(\rho)$ is the unique amenable dimension function.
  In particular when ${\mathcal C}'={\rm Hilb}$
then  
$\rho\to {\rm dim}({\mathcal F}(\rho))$
is a genuine dimension function and this implies that ${\mathcal F}$ is tensorial.
\end{proof}

\begin{rem} Theorem 5.31  of Longo and Roberts  \cite{LR} shows that a  strict tensor category with a unitary braided symmetry satisfying some extra conditions is  amenable. Note that amenability in their sense is defined in a   different way than the notion that we are using
(the intrinsic dimension function is amenable), which in turn is
closely related to amenability in the sense of Popa \cite{Popa}. They discuss the relation between the two notions  in  their Sect. 5.
\end{rem}

 By the previous corollary   when the range category for a weak quasitensor functor is ${\rm Hilb}$ then the properties of unitarity and weak tensoriality may coexist only when 
the functor is automatically tensorial and the intrinsic dimension takes integral values.
Thus when a specific $^*$-functor ${\mathcal F}:{\mathcal C}\to{\rm Hilb}$ on a fusion category is given such that the intrinsic dimension differs
from the associated vector space dimension then ${\mathcal F}$ admits no  unitary weak tensor structure $(F, G)$.
 On the other hand by the results of Sect. \ref{20}, preceded by \cite{CP} for the type $A$ case,   non-unitary weak tensor structures exist.
 For example, this applies to the functor $W$ on
${\mathcal C}({\mathfrak g}, q, \ell)$ at level $k\geq1$.   
  
    \bigskip

\section{Constructing integral  wdf; solution to Galindo's problem}\label{14}
 
In  \cite{Gal} Galindo asks whether a fusion category may admit more than
a unitary structure making it into a unitary tensor category. In \cite{GHR} the authors solve the problem
in some special cases, e.g. pointed and weakly group theoretical categories, and show in these cases a stronger property called
complete unitarity. A  proof   has been given by Reutter in \cite{Reutter} with different methods.
The following consequence of Theorem \ref{unitarizability} gives a positive answer to Galindo's question for a wide class of tensor categories 
with possibly infinitely many simple objects. Note that we do not assume   rigidity.
We prove the following Theorem \ref{Galindo2} stated in the introduction.

\begin{thm}\label{Galindo2}
Let ${\mathcal C}_1$ and ${\mathcal C}_2$ be tensor equivalent $C^*$-tensor categories endowed with an integral weak dimension function (e.g. they are finite semisimple tensor categories). Then ${\mathcal C}_1$ and ${\mathcal C}_2$ are also unitarily tensor equivalent.
\end{thm}

\begin{proof}
It follows from Theorem \ref{unitarizability} with ${\mathcal C}_1={\mathcal C}^+={\mathcal C}$, ${\mathcal F}$ identity, and ${\mathcal C}_2={\mathcal D}^+$.
\end{proof}

  In Sect. 5 we have remarked about the role of integral weak dimension functions for semisimple tensor categories
in relation to Tannaka-Krein duality and weak quasi-Hopf algebras. Moreover in Sect. \ref{12} we have used them to
turn $C^*$-categories into tensor $C^*$-categories.
We next show how to construct these functions for  a wide classes of   categories.

\begin{prop}\label{existence_integral_dimension_functions1}
Let ${\mathcal C}$  be a semisimple tensor category and $d$ be a positive (symmetric) weak dimension function on ${\rm Irr}({\mathcal C})$ taking values $\geq 1$. Then for any integer $M\geq4$,
$D(\rho)=M\lfloor d(\rho)\rfloor$ $\rho\in{\rm Irr}({\mathcal C})$, $\rho\neq\iota$, defines an integral (symmetric) weak dimension function.

\end{prop}

\begin{proof}
We need to show (\ref{weak_dim_funct_inequality}) for any pair of non-trivial objects $\rho$, $\sigma\in{\rm Irr}({\mathcal C})$.
We have 
 $$\sum_{\tau\in{\rm Irr}({\mathcal C})}D(\tau) {\rm dim}(\tau, \rho\otimes\sigma) \leq   
 \sum M\lfloor d(\tau)\rfloor {\rm dim}(\tau, \rho\otimes\sigma)  \leq $$
 $$\sum M  \lfloor d(\tau) {\rm dim}(\tau, \rho\otimes\sigma) \rfloor \leq
M\lfloor d(\rho)d(\sigma)\rfloor\leq M(\lfloor d(\rho)\rfloor + 1)(\lfloor d(\sigma)\rfloor + 1) \leq $$
$$4 M \lfloor d(\rho)\rfloor\lfloor d(\sigma))\rfloor
= \frac{4}{M} D(\rho) D(\sigma) \leq D(\rho)D(\sigma).$$   \end{proof}
Thus all we need to construct integral weak dimension functions is a positive weak dimension function, 
and we   then
ask when such a function exists and how to construct it. 
By  Prop. \ref{bound_below} a necessary condition is that the operators
$\Lambda_\rho$ of left regular
representation of ${\rm Gr}_{\mathbb C}({\mathcal C})$ on $\ell^2({\rm Irr}({\mathcal C}))$
are bounded. 
This is also a sufficient condition when 
${\rm Gr}_{\mathbb C}({\mathcal C})$ is commutative by Theorem 3.5 in \cite{Yamagami}. In the general case,
 we describe two more   classes of examples.

\begin{thm}\label{existence_integral_dimension_functions2}
Any semisimple rigid $C^*$-tensor category or any semisimple rigid  tensor category with amenable fusion ring 
(e.g. any  fusion category) admits a natural positive symmetric dimension function, and therefore infinitely many integral symmetric weak dimension functions. 
\end{thm}

\begin{proof}
The categories in the statement are all known to admit positive symmetric dimension functions, they are respectively given by the intrinsic dimension \cite{LR}, the norm of the left regular representation, see  \cite{CQGRC}  
 and also Sect \ref{13}. Fusion categories are amenable and  the Frobenius-Perron dimension is the unique positive dimension of the representation ring, cf. Cor. 2.7.8 in \cite{CQGRC} and \cite{ENO}.
\end{proof}

   \begin{rem}\label{rational}
The previous result for fusion categories was   observed in \cite{MS, Schomerus, HO}. More precisely,   a    semisimple tensor category ${\mathcal C}$ with finitely many inequivalent simple objects always admits    positive integral weak dimension functions and
when ${\mathcal C}$ is  a fusion category then $D$ may be chosen symmetric. An example is given by the function
taking constant value  ${\rm Max}_{\rho, \sigma}\sum_{\tau\in{\rm Irr}({\mathcal C})}{\rm dim}(\tau, \rho\otimes\sigma)$ 
  for non-trivial $\xi\in {\rm Irr}({\mathcal C})$
 \cite{Schomerus}.  Note that   any other   integer     larger than the constant value   of the previous remark defines another weak dimension function and this immediately shows that a fusion category admits infinitely many weak dimension functions. \end{rem}
 
 It follows from Theorem \ref{amenability2}  that when ${\mathcal C}$ is a semisimple rigid  tensor category with amenable fusion ring then every symmetric
  positive  integral weak dimension function bounds from above   the amenable dimension. 
  This interesting bound together with the results of this section
  shows the great  flexibility of weak quasi-Hopf algebras for this class of categories.

\section{Examples of fusion categories with different natural  integral  wdf}\label{15}

  Motivated by Remark \ref{rational}, it is natural to ask 
    whether     a given fusion category ${\mathcal C}$ may admit more than one weak integral dimension function corresponding to a weak  Hopf algebra. In this subsection we construct   examples
    indicating that this eventuality occurs. The first class of examples   is associated to   pointed fusion categories over the cyclic group ${\mathbb Z}_N$
      and relies on the basic example 
    $A_W({\mathfrak sl}_N, q, \ell)$ for the minimal value of $\ell$. The second example  shows that already for ${\mathbb Z}_2$ there are infinitely many weak dimension functions of this kind, and are obtained using the general constructions of Sect. \ref{6}.
     We shall need
    the ribbon structure naturally associated to the $R$-matrix   of ${\mathcal C}({\mathfrak sl}_N, q, \ell)$. These formulas   will be  recalled   in the next section. 
      
    \begin{ex}\label{pointed_case}      Let $G$ be a finite group and $\omega\in H^3(G, {\mathbb C}^\times)$. The pointed fusion category ${\rm Vec}^\omega_G$ admits the natural dimension function taking value $1$ on every irreducible
    and the associated quasi-Hopf algebra is ${\rm Fun}(G, {\mathbb C})_{\omega}$, see Example 5.13.6 in \cite{EGNO}.
    In particular, we obtain a Hopf algebra if and only if $\omega$ is trivial in $H^3(G, {\mathbb C}^\times)$.
    We next see that for $G={\mathbb Z}_N$ and $\omega=1$ for $N$ odd ($\omega=-1$ for $N$ even) this fusion category may also be described as the representation category of $A_W=A_W({\mathfrak sl}_N, q, N+1)$.
    In other words, if $g$ denotes  the natural generator of ${\mathbb Z}_N$, 
      $D(g)=N$ corresponds to a   weak  Hopf algebra.

Consider the fusion category   ${\mathcal C}({\mathfrak sl}_N, q, \ell)$ for $q= e^{i\pi/N+1}$  and let $X$ denote the object corresponding to the vector representation of $U_q({\mathfrak sl}_N)$. We have
 $d(X)=1$ and   the Grothendieck ring  ${\rm Gr}({\mathcal C}({\mathfrak sl}_N, q))$ is ${\mathbb Z}{\mathbb Z}_N$
 with basis given by
 the objects $X=X_{\Lambda_1}, \dots, X_{\Lambda_{N-1}}$ corresponding to the fundamental weights.
 The fusion rules are given by $X^k=X_{\Lambda_k}$ for $k\leq N-1$ and $X^N=1$   \cite{KW}.
 It follows that ${\mathcal C}({\mathfrak sl}_N, q)$ is tensor equivalent to
   ${\rm Vec}^\omega_{{\mathbb Z}_N}$ for some $\omega\in H^3({\mathbb Z_N}, {\mathbb T})$, cf. Ex. \ref{pointed}. 
   Hence in particular ${\rm Vec}^\omega_{{\mathbb Z}_N}$  admits a weak dimension function as required, and we are left to determine $\omega$.
 The group   $H^3({\mathbb Z_N}, {\mathbb T})$ is isomorphic  to the cyclic group ${\mathbb Z}_N$, that we write in multiplicative notation.
   An explicit  isomorphism associates 
  the $N$-th root of unity $w$ to the 
  $3$-cocycle $\omega$  is given by   (\ref{root}). For the category ${\mathcal C}({\mathfrak sl}_N, q, \ell)$ the corresponding $w$ may be determined following
the procedure   at the end of page  126 in \cite{KW}. In this case, the middle map is identity since the category is strict.
Taking into account the equation appearing in
   Prop. A.5 in \cite{CP} with the additional information that $S$ is an isomorphism for the minimum value of the level,
    we find   that   $\omega=1$   for $N$ odd and $\omega=-1$ for $N$ even. 
    
    Alternatively, we may determine $\omega$ in a more direct way as follows.
    On one hand it is not difficult to see that the only possible values are $\omega=\pm 1$.  (We shall 
see a  more general statement for higher levels in   Prop. \ref{braiding_constraint}.) On the other, by the
        general criterion in Exercise 8.4.11 (iii) pag. 206 in \cite{EGNO}, if a pointed fusion category  ${\rm Vec}^\omega_G$ is braided with braiding $c$ then
  $\omega=1$ if and only if for any element $\gamma\in G$ of order some power of $2$, say $2^r$, the associated quadratic form $q(\gamma)=c(Y, Y)$, with $Y$ simple of class $\gamma$, is of order $\leq 2^r$. This immediately leads to triviality of $\omega$  if   $N$ is odd. For $N$ even 
 we use the fact that $q(X_{\Lambda_k})$ equals the ribbon structure $\theta_{X_{\Lambda_k}}$,
  see Subsect. \ref{22.1}, and that   $\theta_{X_{\Lambda_k}}=q^{\frac{k(N-k)(N+1)}{N}}=e^{\frac{i\pi k(N-k)}{N}}$, by   the proof of Prop.
  \ref{twist_general_case}. Writing $N=2^r h$ with $h$ an odd integer, it follows  that 
  $\Lambda_h$ has order $\frac{N}{h}=2^r$ but $q(X_{\Lambda_h})^{p}\neq 1$ for all $1\leq p\leq 2^r$.

  \end{ex}

  \begin{ex}\label{infinitely_many_weak_dim_funct}
  We give  examples of infinitely many weak dimension functions corresponding to weak  Hopf algebras on the fusion categories ${\rm Vec}^\omega_{{\mathbb Z}_2}$. They are given by $D(g)=2h+1$ for  ${\rm Vec}_{{\mathbb Z}_2}$ and $D(g)=2h$ for   ${\rm Vec}^{-1}_{{\mathbb Z}_2}$, for $h\geq1$, with $g$ the group generator and  
$\omega\in H^3({\mathbb Z}_2, {\mathbb T})\simeq  {\mathbb Z}_2$.

 Consider  the fusion category ${\mathcal C}({\mathfrak sl}_2, q, \ell)$ with $q=e^{i\pi/\ell}$ and $\ell\geq 3$, and the associated Grothendieck ring (the Verlinde ring)
 $R_{2,\ell}$   with  basis  given  by the equivalence classes of selfconjugate irreducible objects 
$X_0=I$, $X_1, \dots, X_k$. Fusion rules are given by 
  $X_iX_j=\sum_{\max\{i+j-k, 0\}}^{\min\{i,j\}} X_{i+j-2r}$, see \cite{Chari_Pressley, EGNO}. The element $X=X_k$ 
   satisfies   $X^2=I$, so it generates a pointed full fusion subcategory   ${\mathcal C}_k\simeq{\rm Vec}^\omega_{{\mathbb Z}_2}$. We 
  determine   $\omega\in\{\pm1\}$  
   by means of Ex. 8.4.11 iii) in \cite{EGNO} again, so in this case  $\omega=1$  precisely when the quadratic form
 $q(g)=c(X, X)$ associated to   restricted braiding of ${\mathcal C}_k$   satisfies  $q(g)=1$ or $q(g)^2=1$.  Arguments similar to those of the previous example give
 $q(g)=\theta_X$, with $\theta$ the usual ribbon structure of $c$, whose value on $X=X_k$ is  $\theta_X=q^{k(k+2)/2}=e^{i\pi k/2}$  cf. Prop. \ref{twist_general_case}. It follows that  $\omega=1$ if and only if $k$ is even.
 On the other hand,  ${\mathcal C}({\mathfrak sl}_2, q, \ell)$ is tensor equivalent to
  the representation category of 
  $A_W({\mathfrak sl}_2, q, \ell)$  so  ${\mathcal C}_k$
  is tensor equivalent to a quotient weak  Hopf algebra
  $A\to B_k$ by Cor. \ref{quotients}. Since
   $X$ corresponds to a representation of   $A$ of dimension $k+1$,   we have $B_k={\mathbb C}\oplus M_{k+1}({\mathbb C})$ and a weak dimension function $D$ on ${\mathcal C}_k$, and therefore on ${\rm Vec}^\omega_{{\mathbb Z}_2}$ as required.
   \end{ex}
   
 \begin{ex}  The methods of the above examples may be combined to construct more examples of weak  Hopf algebras. a) For example, if $g\in{\mathbb Z}_N$
   is the natural generator, for $k\leq N-1$, $g^k$ generates a cyclic subgroup of order $M=\frac{N}{\gcd\{k, N\}}$.
   Therefore the  full subcategory of ${\mathcal C}({\mathfrak sl}_N, q, \ell)$ for $q=e^{i\pi/N+1}$ generated by $X_{\Lambda_k}$, which is pointed  over ${\mathbb Z}_M$,    corresponds to a quotient  of 
   $A_W({\mathfrak sl}_N, q, N+1)$
  (with   dimension    of the natural generator  of ${\mathbb Z}_M$ given by
   $D(h)={{N}\choose{k}}$) and also  to 
   $A_W({\mathfrak sl}_M, q, M+1)$
 (with   dimension $D'(h)=M$) with a possibly twisted associator. b) The even subcategory of 
 ${\mathcal C}({\mathfrak sl}_2, q, \ell)$ for $q=e^{i\pi/\ell}$ is an example of non-pointed full fusion subcategory, and therefore it gives rise to a quotient weak  Hopf algebra $B={\mathbb C}\oplus M_2\oplus M_4\dots$. c) More information on full fusion subcategories of 
${\mathcal C}({\mathfrak sl}_N, q, \ell)$ for $q=e^{i\pi/\ell}$ may be found in \cite{S}.
   \end{ex}
   
   \begin{rem}  Ribbon structures    first appeared as statistics phases   for  WZW and coset models in conformal field theory. Some formulae for the statistics phases, including the automorphism case of interest in   Ex. 
   \ref{infinitely_many_weak_dim_funct},
       have been generalized  by Rehren in the   framework of  conformal nets. Most importantly, in that paper the author derives the axioms of a modular category extending previous work for certain conformal models  \cite{Rehren} and references therein.
The  ribbon structure in the conformal net approach to CFT   is given by $\theta_X=e^{2\pi i h_X}$ 
with $h_X$ the minimal eigenvalue of the conformal Hamiltonian $L_0$ in the irreducible representation $X$,
by the conformal spin and statistics theorem  \cite{GL}. In the framework of vertex operator algebras 
one has an analogous formula \cite{Huang(modularity), Huang1, Huang2}.
  
\end{rem}

  \begin{rem}     In the setting of rigid tensor $C^*$-categories with infinitely many simple objects,  
   Van Daele and Wang constructed compact quantum groups $A_o(F)$ associated to an invertible matrix $F$ with ${\rm rk}(F)\geq 2$ satisfying suitable properties \cite{W_VD}, which reduce to Woronowicz compact quantum groups ${\rm SU}_q(2)$ for ${\rm rk}(F)= 2$.
     For a given $q>1$, ${\rm Rep}(A_o(F))$ turns out to be tensor equivalent to  ${\rm Rep}({\rm SU}_q(2))$ with $q$ suitably determined by $F$    \cite{Banica1, Banica2}. It follows that ${\rm Rep}({\rm SU}_q(2))$ admits  the      (non-weak) dimension 
  function taking  the   generating representation to
  the rank of $F$.  Note that only finitely many    (non-weak)  integral dimension functions arise in this way. This follows from the fact that    ${\rm rk}(F)$
  is bounded above  by the quantum dimension \cite{Wor}.
 In this setting, it is important to   recall the remarkable work by Neshveyev and Yamashita on the classification of compact quantum groups
 that beyond the   fusion rules, share the     integral dimensions with a given compact simple simply connected Lie group $G$, see \cite{NY_towards} and references therein.

  \end{rem}

  \section{Quantum group $U_q({\mathfrak g})$ at roots of unity, $^*$-involution, classical limit}\label{73}

  Let ${\mathfrak g}$ be a complex simple Lie algebra and $q$ a primitive  complex root of unity. We denote by $\ell$ the order of $q^2$. 
Let $U_q({\mathfrak g})$ be the quantized universal enveloping algebra  in the sense of Lusztig, see below for a definition 
and references. 
It is known that the category of  
    finite dimensional representations of $U_q({\mathfrak g})$ is not semisimple,
   but it gives rise to a semisimple ribbon fusion category that we denote by ${\mathcal C}({\mathfrak g}, q, \ell)$
   following \cite{Rowell2}.
    Moreover,   the categories ${\mathcal C}({\mathfrak  g}, q, \ell)$
   are known to be  modular  for certain values of $q$ see
   \cite{Andersen, Andersen1, GK, RT,   Rowell2, Sawin} and references therein, see also Subsects. \ref{18.1}--\ref{18.3}.

   Furthermore, by   work of Kirillov, Wenzl, Xu one can construct 
   a unitary ribbon category ${\mathcal C}^+({\mathfrak g},  q, \ell)$
 equivalent to  ${\mathcal C}({\mathfrak g}, q, \ell)$   for    certain   primitive roots of unity, that we call the {\it minimal roots} and precisely define   in the following Def. \ref{minimal_root}. \cite{Kirillov, Wenzl, Xu}. 
 
    In this section we construct   semisimple   weak quasi-Hopf algebras  
     associated to ${\mathcal C}({\mathfrak  g}, q, \ell)$ and unitary weak quasi-Hopf algebras associated to
     ${\mathcal C}^+({\mathfrak g},  q, \ell)$ when $q$ is a minimal root.
    Our approach may broadly be summarized as follows. 
 
From the categories  ${\mathcal C}({\mathfrak  g},  q, \ell)$, we construct  weak quasi-tensor functors to the category of vector spaces and then we use Tannaka-Krein reconstruction to obtain our examples.  

We shall do this following  two alternative approaches, and both turn out useful for us in the study
 of unitary tensor categories. The first approach goes back
 to \cite{MS, S, HO}. It consists in identifying a certain integral valued weak dimension function $D$ on ${\mathcal C}({\mathfrak g}, q, \ell)$,
 and then we apply the abstract reconstruction result, Theorem \ref{propweakdim}.
 This leads to the construction of a  ribbon weak quasi-Hopf algebra
 $A({\mathfrak g}, q, \ell)$ corresponding to
 ${\mathcal C}({\mathfrak g}, q, \ell)$
  which is defined up to twist and isomorphism. Moreover, when $q$ is a minimal root, we apply
  Theorem  \ref{TheoremTannakaStar} and we obtain a unitary  structure  $A^+({\mathfrak g}, q, \ell)$ on $A({\mathfrak g}, q, \ell)$.
   By the results of Sect. \ref{12}  this general approach 
 addresses  the study
 of unitary structures via the associator. It follows that this
 viewpoint will turn out fruitful for the construction of unitary ribbon
 structures for representation categories of affine VOA in Sect. \ref{VOAnets2}. It perhaps conveys the idea of the amount 
 of   information    needed to obtain these unitary structures from other sources for which they are known to exist.

A second  approach consists in identifying a natural functor $W: {\mathcal C}({\mathfrak  g}, q, \ell)\to{\rm Vec}$  associated to the same dimension function $D$ as before,  and thus it is a particular case of the former, and will be studied in Sect. \ref{20}.

When $q$ is a minimal root, the work of \cite{Wenzl} shows that $U_q({\mathfrak g})$
is a  Hermitian   coboundary Hopf algebra with compatible involution (in a topological sense).
We shall introduce this notion in Sect. 
\ref{20} and summarize   this result
in Theorem \ref{U_q_as_a_Hermitian_ribbon_h}. However, in this section we shall not need to go into
these details. \medskip
  
 \subsection{Quantum groups at roots of unity, the real form ${\mathcal U}^\dag_{{\mathcal A}'}({\mathfrak g})$ and   ribbon structure}\label{18.1}

 We briefly recall the basic results on quantum groups at roots of unity  that  we shall need. For a complete presentation
 we refer to  \cite{Chari_Pressley,  Lusztig_book}.

  Let  ${\mathfrak g}$ be a f.d. complex simple Lie algebra, and $q$ a complex root of unity whose order we denote by $\ell'$.   (Thus the order $\ell$ of $q^2$ is given by $\ell=\ell'$  if $\ell'$ is odd and $\ell=\ell'/2$  if $\ell'$ is even.)
    Note that our     
 $\ell$ has the same meaning  in \cite{Rowell2}, our $\ell$ corresponds to $\ell'$ (and conversely) in \cite{Sawin}, \cite{Turaev}.
 We anticipate an important lower bound condition that we require on $\ell$, that will be briefly explained immediately after the introduction
 of background material on Drinfeld-Jimbo quantum group, their specialization
at a root of unity and their representation theory.

   \begin{defn}\label{large_enough} We shall say that the order $\ell$ of $q^2$ is large enough if
 $\ell>{h}$ when $\ell$ is not divisible by $d$ and $\ell>dh^\vee$ otherwise, with $h$ the Coxeter number and $h^\vee$ the dual Coxeter number  of ${\mathfrak g}$. We define the positive integer $k$ by $\ell=d(h^\vee+k)$ if $d|\ell$
 and $\ell=h+k$ otherwise, and refer to $k$  as   the {\it    level} associated to $\ell$. We define
  $\ell=k=\infty$ if $q$ is not a root of unity.

   \end{defn}

 The
   level $k$ introduced in Def. \ref{large_enough}   is an important parameter for its relation with affine Lie algebras or affine vertex operator algebras.
   We are mainly interested in the case 
   $\ell$   divisible by $d$, but also the case of  other type of larger orders of roots of unity (possibly also of infinite order)
   will be important for us to connect modules of $U_q({\mathfrak g})$  with modules of a corresponding vertex operator algebra at fixed level $k$.

    Provided $q^2$ has finite order  $\ell$  large enough, we shall recall in Theorem \ref{irreducibles_of_fusion_category} that  the
    corresponding  Weyl alcove $\Lambda^+(q)$
    labels the irreducible objects of the   fusion category ${\mathcal C}({\mathfrak g}, q, \ell)$. The dominant weights
belonging to the alcove will correspond to simple modules with positive quantum dimension.

 In particular, quantised Weyl modules $V_\lambda(q)$ with $\lambda$ in  $\Lambda^+(q)$   are simple, with a positive definite hermitian form, by Prop.
   2.4 in \cite{Wenzl}.  
  
   The category of tilting modules
  was first introduced by Andersen \cite{Andersen}.
  This category admits an alternative definition as noted by Wenzl
  in   \cite{Wenzl},  
  as the category of direct sums of direct summands of full tensor powers of  a specific generating representation $V(q)$ for each Lie type, for convenience we recall this fact in Theorem \ref{tilting_as_fundamental}. This  
  is the definition that we will also adopt. It   gives the category ${\mathcal C}({\mathfrak g}, q, \ell)$ a perhaps  non-commutative geometric aspect, that will be crucial
  to connect with the fusion categories of the corresponding vertex operator algebra.

   As noted in \cite{Wenzl}, the basic requirement that makes the mentioned alternative definition work,
   is that $\ell$ be sufficiently large  so that $V(q)$ belongs to the Weyl alcove $\Lambda^+(q)$ corresponding to $q$. So we need
   to determine the minimal order $\ell$ of $q^2$ divisible by $d$ so that $V(q)$ belongs to the Weyl alcove $\Lambda^+(q)$ for each Lie type.
   We shall see in Prop. \ref{fundamental_alcove} that this requirement on $V(q)$ puts essentially no   lower bound
   on   the level $k$ (and that $\ell$ will automatically satisfy  Def. \ref{large_enough})
    for all the Lie types  (with the only exception of $E_8$, for which the minimum allowed
    level is $2$).
   
  \medskip

Let  ${\mathfrak h}$ a Cartan subalgebra, $\alpha_1,\dots,\alpha_r$ a set
of simple roots, and $A=(a_{ij})$ the associated Cartan matrix. Consider the unique invariant symmetric and bilinear form on ${\mathfrak h}^*$ such that 
\begin{equation}\langle\alpha,\alpha\rangle=2, \quad\text{for a short root $\alpha$}
\end{equation}\label{normalization_inner_product_qg}
and let
   $\theta$ denote the highest root.
 Let
$E$ be the real vector space generated by the roots endowed with its euclidean structure $\langle x, y\rangle$. Let
 $\Lambda$ be the weight lattice of $E$ and $\Lambda^+$ the cone of dominant weights.

Consider  the complex $^*$--algebra 
${\mathcal A}={\mathbb C}[x,x^{-1}]$ 
of   Laurent polynomials with involution
$x^*=x^{-1}$, and let ${\mathbb C}(x)$ be the associated quotient field, endowed with the involution naturally induced from ${\mathbb C}[x,x^{-1}]$. 
 We consider Drinfeld-Jimbo quantum group     $U_x({\mathfrak g})$, i.e. the algebra over ${\mathbb C}(x)$ defined by generators
$E_i$, $F_i$, $K_i$, $K_i^{-1}$, $i=1,\dots, r$, and relations  
$$K_iK_j=K_jK_i,\quad K_iK_i^{-1}=K_i^{-1}K_i=1,$$
$$K_iE_jK_i^{-1}=x^{\langle\alpha_i,\alpha_j\rangle}E_j,\quad K_iF_jK_i^{-1}=x^{-\langle\alpha_i,\alpha_j\rangle}F_j,$$
$$E_iF_j-F_jE_i=\delta_{ij}\frac{K_i-K_i^{-1}}{x^{d_i}-x^{-d_i}},$$
$$\sum_0^{1-a_{ij}}(-1)^kE_i^{(1-a_{ij}-k)}E_j E_i^{(k)}=0,\quad\sum_0^{1-a_{ij}}(-1)^kF_i^{(1-a_{ij}-k)}F_j F_i^{(k)}=0,\quad i\neq j,$$
where $d_i=\langle\alpha_i,\alpha_i\rangle/2$, and, for $k\geq0$, $E_i^{(k)}=E_i^k/[k]_{d_i}!$,  $F_i^{(k)}=F_i^k/[k]_{d_i}!$. Note that $d_i$ is an integer, hence so is
every inner product $\langle\alpha_i,\alpha_j\rangle$.  Quantum integers and factorials are defined in the usual way, $[k]_x=\frac{x^k-x^{-k}}{x-x^{-1}}$; $[k]_x!=[k]_x\dots[2]_x$, $[k]_{d_i}=[k]_{x^{d_i}}$, and result selfadjoint scalars of ${\mathbb C}(x)$. There is a unique $^*$--involution on $U_x({\mathfrak g})$ making it into a $^*$--algebra over ${\mathbb C}(x)$ such that     $$K_i^*=K_i^{-1},\quad E_i^*=F_i.$$  

This algebra becomes a Hopf algebra, with coproduct, counit, and antipode  defined, as follows, see e  in \cite{Wenzl}, where his $\tilde{K}_i$ corresponds to our $K_i$,
see also \cite{Chari_Pressley, Sawin},
$$\Delta(K_i)=K_i\otimes K_i,$$
$$\Delta(E_i)=E_i\otimes K_i+1\otimes E_i, \quad \Delta(F_i)=F_i\otimes 1+K_i^{-1}\otimes F_i,$$
$$S(K_i)=K_i^{-1}, \quad S(E_i)=-E_iK_i^{-1}, \quad S(F_i)=-K_iF_i,$$
$$\varepsilon(K_i)=1,\quad \varepsilon(E_i)=\varepsilon(F_i)=0.$$

One has the following relations between coproduct, antipode and involution for $a\in U_x({\mathfrak g})$,
\begin{equation}\label{Delta} \Delta(a^*)= \Delta^{\text{op}}(a)^*\end{equation}
\begin{equation}\label{counit_antipode} \varepsilon(a^*)=\overline{\varepsilon(a)}, \quad S(a^*)=S(a)^*, \quad S^2(a)=K_{2\rho}^{-1}aK_{2\rho},\end{equation}
where   $2\rho$  the sum of the positive roots, and, for an element $\alpha=\sum_i k_i\alpha_i$ of the root lattice, $K_\alpha:=K_1^{k_1}\dots K_r^{k_r}$.\medskip

\begin{rem}\label{compact_real_form}
Let us regard the universal enveloping algebra $U({\mathfrak g})$ as the classical limit of $U_x({\mathfrak g})$, that is with generators $E_i$, $F_i$, $H_i$ such that $K_i=x^{d_iH_i}$, for a rigorous explanation of this equality see
e.g. \cite{Chari_Pressley} pag. 304. Since  $(K_i)^*={K_i}^{-1}$
($x^*=x^{-1}$), the classical limit of the
$^*$-involution of $U_x({\mathfrak g})$ is given by 
$$E_i^*=F_i, \quad\quad F_i^*=E_i \quad\quad H_i^*=H_i.$$
This is the same as the classical limit arising from the usual $^\dag$--involution
of $U_x({\mathfrak g})$ making it into a Hopf $^*$-algebra ($E_i^\dag=F_i$,  $F_i^\dag=E_i$, $K_i^*=K_i$, $x^*=x$)
we refer to Sect. 13.1 \cite{Chari_Pressley} and
 Sect. 2.4 \cite{CQGRC} for more details.
In both cases the    real form ${\mathfrak g}_{\mathbb R}$ of ${\mathfrak g}$ corresponding to the classical
limit is the usual compact real form, the real Lie subalgebra  of skeweak  adjoint elements $a^*=-a$.

 Let $G$ be the connected simply connected complex Lie group with Lie algebra ${\mathfrak g}$. Then the Lie subgroup
 $K$ of $G$ with Lie algebra ${\mathfrak g}_{\mathbb R}$ is compact.
 Moreover for any $\lambda\in\Lambda^+$, the Weyl module $V_\lambda$ of $U({\mathfrak g})$ is unitarizable, thus
 for a given highest weight vector $v_\lambda$, $V_\lambda$
  admits a unique positive Hermitian form  $(\xi, \eta)$ such that $(v_\lambda, v_\lambda)=1$, $(a\xi, \eta)=(\xi, a^*\eta)$ for $\xi$, $\eta\in V_{\lambda}$, $a\in U({\mathfrak g})$,   and it follows that $V_\lambda$ defines a 
  unitary representation of    $K$. At the level of specialization to a complex number $q$, we recall that $U({\mathfrak g})$ is isomorphic to 
  the quotient of $U_1({\mathfrak g})$ by the ideal generated by the elements $K_i-1$.
 \end{rem}
 
\begin{rem} Note that the $^*$-involution of 
$U({\mathfrak g})$ is related to the order $2$-antilinear automorphism $T$ of  ${\mathfrak g}$ defining the compact real form
 by $a^*=T(S(a))$, $a\in U({\mathfrak g})$, where $S$ is the antipode. The presence of the antipode $S$ is necessary because the conjugation
$T$ extends to an antilinear multiplicative map on $U({\mathfrak g})$, but the $^*$-involution needs to be antimultiplicative.
 \end{rem}

Lusztig   form  ${\mathcal U}^{\rm res}_{\mathcal A}({\mathfrak g})$ is   the algebra over ${\mathcal A}$
generated by $E_i^{(k)}$, $F_i^{(k)}$, $K_i$. 
It is an integral form of $U_x({\mathfrak g})$ in that
$U_x({\mathfrak g})={\mathcal U}^{\rm res}_{\mathcal A}({\mathfrak g}){\mathbb C}(x)$ and ${\mathcal U}^{\rm res}_{\mathcal A}({\mathfrak g})$ is
 free over ${\mathcal A}$.
The algebraic relations among the generators only involve coefficients in ${\mathcal A}$, thus this form may be specialized at roots of unity. 
    To construct an $R$-matrix,  that we briefly recall in the next subsection, we   need to embed ${\mathcal U}^{\rm res}_{\mathcal A}({\mathfrak g})$ into a larger algebra. By \cite{Sawin}, 
  we   need to extend the ring of scalars from ${\mathcal A}$ to
 $${\mathcal A}':={\mathbb C}[x^{1/L}, x^{-1/L}],$$ with $L$ the smallest positive integer such that
$L\langle\lambda, \mu\rangle\in{\mathbb Z}$ for all dominant weights $\lambda$, $\mu$.
The   values of $L$  for all Lie types are listed in \cite{Sawin}. For example,   $L=N$ for ${\mathfrak g}={\mathfrak sl}_N$.

 We define the integral form   ${\mathcal U}_{{\mathcal A}'}({\mathfrak g})$ as
  the ${\mathcal A}'$--subalgebra generated by the elements $E_i^{(k)}$, $F_i^{(k)}$ and $K_i$.
This is known to be a $^*$--invariant Hopf ${\mathcal A}'$--algebra with the structure inherited from $U_x({\mathfrak g})$. 
Applying the construction in Sect. 1 of \cite{Sawin} to the modified polynomial ring, we associate to 
${\mathcal U}_{{\mathcal A}'}({\mathfrak g})$    an  extended   Hopf 
algebra ${\mathcal U}^\dag_{{\mathcal A}'}({\mathfrak g})$
in a topological sense, that is the coproduct takes values in a topological completion of 
${\mathcal U}^\dag_{{\mathcal A}'}({\mathfrak g})\otimes {\mathcal U}^\dag_{{\mathcal A}'}({\mathfrak g})$.

Sawin gives a detailed construction of the $R$-matrix of ${\mathcal U}^\dag_{{\mathcal A}'}({\mathfrak g})$. To develop this, considers among other things a definition of ${\mathcal U}^\dag_{{\mathcal A}'}({\mathfrak g})$ containing the function
algebra ${\rm Map}(\Lambda, {\mathcal A}')$ as a Hopf subalgebra. Moreover
 ${\mathcal U}^\dag_{{\mathcal A}'}({\mathfrak g})$
 has the structure of a
topological ribbon Hopf algebra. Another relevant aspect of this construction that also plays an important tole
for us  is that ${\mathcal U}^\dag_{{\mathcal A}'}({\mathfrak g})$ embeds faithfully in the discrete algebra given by the direct product  
\begin{equation}\label{direct_product_integral_form}{\mathcal U}^\dag_{{\mathcal A}'}({\mathfrak g})\to\Pi_V {\rm End}(V)\end{equation}
of linear endomorphism algebras on weight modules $V$. A weight module is defined as 
an ${\mathcal A}'$-free  ${\mathcal U}^\dag_{{\mathcal A}'}({\mathfrak g})$-module which is a finite direct sum of finitely generated $\lambda$-weight ${\mathcal A}'$-modules $V_\lambda$ associated to an integral weight $\lambda\in\Lambda$. In turn, the $\lambda$-weight ${\mathcal A}'$-submodule is defined by
$V_\lambda=\{v\in V:  K_i v=x^{\langle\lambda,\alpha_i\rangle} v, fv=f(\lambda)v\}$.
Thus elements of ${\mathcal U}^\dag_{{\mathcal A}'}({\mathfrak g})$ act as matrices with entries in ${\mathcal A}'$
on basis elements.
The coproduct \begin{equation}\label{the_coproduct_range_integral_form}
\Delta: {\mathcal U}^\dag_{{\mathcal A}'}({\mathfrak g}) \to \Pi_{V, W} {\rm End}(V)\otimes {\rm End}(W)\end{equation}
has range the closure of ${\mathcal U}^\dag_{{\mathcal A}'}({\mathfrak g})\otimes {\mathcal U}^\dag_{{\mathcal A}'}({\mathfrak g})$ in 
$\Pi_{V, W} {\rm End}(V)\otimes {\rm End}(W)$.

\bigskip

\subsection{The $R$-matrix of ${\mathcal U}^\dag_{{\mathcal A}'}({\mathfrak g})$ and ribbon structure.}\label{18.2}

The formula for the $R$-matrix $R$ of ${\mathcal U}^\dag_{{\mathcal A}'}({\mathfrak g})$
may be found by the end of Sect. 1 in \cite{Sawin}. 
 The ribbon element $v$ is constructed from $\omega=K_{2\rho}$ which is a pivot (or charmed)   element
and Drinfeld element $u$ via $\omega=uv^{-1}$ as described in the remark following Cor. \ref{spherical_structure_wqh}.

Our $R$-matrix corresponds to $R_{21}^{-1}$ in \cite{Sawin}
and agrees with \cite{Wenzl}.
We shall need to recall a characterization of $R$ following \cite{Wenzl} and going back to   to Lusztig, Chapter 32 \cite{Lusztig}.
Consider $r$-ples of non negative integers $(\nu_1, \dots, \nu_r)$ and $U_\nu^{+}$ the span
of products of $E_i^{(k)}$ where $\nu_i$ is the sum of all the $k$ occurring in a monomial for a fixed $i$.
One similarly defines $U_\nu^{-}$. Consider
$$X=\bigoplus_{\nu} U_{\nu}^{-}\oplus U_{\nu}^{+}.$$

Moreover, consider the unique antilinear (involutive) automorphism 
$T\to\overline{T}$
of ${\mathcal U}^\dag_{{\mathcal A}'}({\mathfrak g})$
 that fixes the generators $E_i^{(k)}$, $F_i^{(k)}$ and such that $\overline{x}=x^{-1}$, 
$\overline{K_i}=(K_i)^{-1}$. This automorphism induces a new coproduct $\overline{\Delta}$ on
 ${\mathcal U}^\dag_{{\mathcal A}'}({\mathfrak g})$ by conjugation.
 
 Then there is a unique element $\Theta\in X$ such that $\Theta_0=I\otimes I$
and $\overline{\Delta}(a)=\Theta\Delta(a)\Theta^{-1}$, $a\in {\mathcal U}^\dag_{{\mathcal A}'}({\mathfrak g})$.
The following result summarizes properties of the $R$-matrix, which go back to ideas of Drinfeld.
Our approach follows Ch. 32 in \cite{Lusztig}, the computation of the action $(c)$ of the ribbon element may be found in
\cite{Drinfeld_cocommutative},
compatibility properties $(a)$ and $(b)$ with the $^*$-involution are in
were proved by Wenzl in Prop. 1.4.1 in \cite{Wenzl}.

 \begin{thm}\label{Properties_of_R} The $R$-matrix of ${\mathcal U}^\dag_{{\mathcal A}'}({\mathfrak g})$ is given by
 $R=\Pi\overline{\Theta}$, where $\Pi$   acts on a tensor product $V\otimes W$ of ${\mathcal A}'$-free weight modules by
 $$\Pi v\otimes w=x^{\langle \mu, \nu\rangle}v\otimes w$$ where $\mu$ and $\nu$ are weights of the weight vectors $v$ and $w$ respectively.
 Moreover we have
 
 \begin{itemize}
\item[(a)] 
$\Delta(a^*)=\Delta^{(\rm op)}(a)^*$,
\item[(b)] $R^*=(R^{21})^{-1}$,
\item[(c)] $v$ acts as the scalar $x^{-{\langle\lambda, \lambda+2\rho\rangle}}$ on an irreducible highest weight module 
with highest weight $\lambda$.

\end{itemize}
 
 \end{thm}

Note that our $\Pi$ denotes $\Pi^{-1}$ in \cite{Wenzl}.
Moreover $\Pi$ belongs to a certain topological completion of 
 ${\mathcal U}^\dag_{{\mathcal A}'}({\mathfrak g})\otimes{\mathcal U}^\dag_{{\mathcal A}'}({\mathfrak g})$, see e.g.
 \cite{Sawin}, hence the same holds for $R$.\bigskip
 
 \begin{defn}\label{quantum_casimir}
 The ribbon element $v\in {\mathcal U}^\dag_{{\mathcal A}'}({\mathfrak g})$ is called the (formal) quantum Casimir operator.
 \end{defn}

\subsection{Weyl modules for ${\mathcal U}^\dag_{{\mathcal A}'}({\mathfrak g})$ and Lusztig and Kashiwara bases.}\label{18.3}

Let $\lambda$ be a dominant integral weight and $V_\lambda(x)$ an
the  irreducible representation of $U_x({\mathfrak g})$ with highest weight $\lambda$,
 and let $v_\lambda$ be the highest weight  vector of $V_\lambda(x)$ (unique up to a nonzero scalar multiple),
 that is $K_i v_\lambda=x^{\langle\alpha_i, \lambda\rangle} v_\lambda$ and $\lambda$ is maximal with this property.
 This module is called a Weyl module, it is direct sum of weight spaces. Tensor product of Weyl modules
 decomposes into a direct sum of Weyl modules with the same multiplicities as in the classical case,
 $$ V_{\lambda}(x)\otimes V_{\mu}(x)=\oplus_\gamma {m_{\gamma}^{\lambda, \mu}} V_{\gamma}(x),$$
and multiplicities $m_{\gamma}^{\lambda, \mu}$  
involve the Weyl group and rely on the classical Racah formula, 
We shall not need these multiplicities  in detail, and we refer the interested reader to
\cite{Humphreys},
Sect. 10 in \cite{Chari_Pressley}.

 A canonical basis $B$ for the ${\mathbb C}(x^{1/L})$-subalgebra $U_x({\mathfrak g})^{-}$ of $U_x({\mathfrak g})$ generated
 by the $F_i$ has been constructed by Lusztig and Kashiwara \cite{Lusztig_canonical_bases}, \cite{Kashiwara}. 
 
 Set $B_\lambda=\{b\in B: bv_\lambda\neq 0\}$. Then $b\in B_\lambda\to bv_\lambda\in V_\lambda(x)$
is a bijection of $B_\lambda$ onto a finite basis of $V_\lambda(x)$, and this is the canonical basis of 
$V_\lambda(x)$. Moreover, $V_\lambda(x)$ has a free ${\mathcal A}'$-submodule $V_{\lambda, {\mathcal A}'}$
such that ${\mathbb C}(x)V_{\lambda, {\mathcal A}'}=V_\lambda(x)$ with basis
$B_\lambda v_\lambda$ for a suitable choice of $v_\lambda$. This is a basis of weight vectors.  (in particular
$V_{\lambda, {\mathcal A}'}$ is a weight module for   ${\mathcal U}^\dag_{{\mathcal A}'}({\mathfrak g})$
in the sense of the previous subsection.) In other words the generators of
${\mathcal U}^\dag_{{\mathcal A}'}({\mathfrak g})$ act on  $B_\lambda v_\lambda$ as matrices with coefficients in ${\mathcal A}'$. Canonical bases plays an important role in Wenzl construction of the 
unitary structure of the fusion category ${\mathcal C}({\mathfrak g}, q, \ell)$, (see   also the following Subsect.
\ref{18.6}.). Loosely speaking,    a canonical basis $v_i$ of $V_{\lambda, {\mathcal A}'}$ specialized in a complex number 
$q$ which varies continuously on an arch, 
is a   basis on $V_\lambda(q)$, and   certain analytic properties that  hold at one extreme of the arc  may be
carried to the other extreme by continuity. We shall recall this more precisely in Theorem \ref{Wenzl_positivity} and Cor.
\ref{corollary_of_positivity}.

For the same   motivations canonical bases  will play a role in the motivation of our abstract Drinfeld-Kohno theorem in Sect. \ref{19} and  our main result in Sect. \ref{21}.

If $\mu$ is another dominant integral weight then following \cite{Wenzl} we endow $ V_{\lambda, {\mathcal A}'}\otimes_{{\mathcal A}'}V_{\mu, {\mathcal A}'}$ with the tensor product of the canonical bases $b v_\lambda\otimes b'v_\mu$, $b$, $b'\in B$. This tensor product decomposes into 
a direct sum of Weyl modules, with their canonical bases,
$$ V_{\lambda, {\mathcal A}'}\otimes_{{\mathcal A}'}V_{\mu, {\mathcal A}'}=\oplus_\gamma {m_{\gamma}^{\lambda, \mu}} V_{\gamma, {\mathcal A}'},$$
and multiplicities $m_{\gamma}^{\lambda, \mu}$  with the same multiplicities  as in the classical case.

\bigskip

\subsection{Specialization $U_q({\mathfrak g})$ at a complex primitive root of unity $q$ of order $\ell'$.}\label{18.4}

We next describe a specialization of ${\mathcal U}^\dag_{{\mathcal A}'}({\mathfrak g})$
which differs slightly from the specialization given in Sect. 2 of \cite{Sawin}, in that we work with   a   complex root of unity, combining with the presentation of  \cite{Wenzl}, and references therein among them  Drinfeld and Lusztig work. 
See also  the restricted specialization in Sect. 11.2 in \cite{Chari_Pressley}.
 
We fix $q\in{\mathbb T}$ a primitive
root of unity of order $\ell'$ (see the beginning of the section for our notation as compared to \cite{Sawin}) and we set $\ell'=\infty$ if $q^n\neq 1$ for all $n\in{\mathbb N}$.
We consider the $^*$--homomorphism ${\mathcal A}'\to {\mathbb C}$ which evaluates 
 $x^{1/L}$ to a specified complex $L$-root $q^{1/L}$ of $q$.
We form the tensor product
  $^*$--algebra,
$$U_q({\mathfrak g}):={\mathcal U}_{{\mathcal A}'}^\dag({\mathfrak g})\otimes_{{\mathcal A}'}{\mathbb C}.$$ 
The algebra $U_q({\mathfrak g})$   becomes
a ribbon complex Hopf algebra with a $^*$--involution, and is topological in the sense of 
\cite{Sawin}. 
Let $V$ be a weight module (as in the last part of Subsect. \ref{18.1}). We set
$$V(q)=V\otimes_{{\mathcal A}'}{\mathbb C}.$$
By (\ref{direct_product_integral_form}), 
we have a faithful embedding into a direct product of full matrix algebras over ${\mathbb C}$
\begin{equation}\label{direct_product_specialized_form}
{\mathcal U}_q({\mathfrak g}) \to\Pi_V {\rm End}(V(q))
\end{equation}
and by (\ref{the_coproduct_range_integral_form})
the coproduct has range 
\begin{equation}\label{the_coproduct_range_specialized_form}
\Delta: {\mathcal U}_q({\mathfrak g})\to \Pi_{V(q), W(q)} {\rm End}(V(q))\otimes {\rm End}(W(q))\end{equation}
in a direct product of full matrix algebras.

Note that the $R$-matrix $R$ and the ribbon element $v\in U_q({\mathfrak g})$ depend (only) on the choice of $q^{1/L}$, see Sect. 1 in \cite{Sawin},
Sect. 1.4 in \cite{Wenzl}. This variability of $R$-matrices and ribbon elements will be useful in Sect. \ref{KW}.

On the other hand, in Sect. \ref{20} in connection with the study of unitary structure and construction of a Drinfeld
twist, it will be important to follow Wenzl \cite{Wenzl} and   specify the value of $q$ (of possibly infinite order) written in the form
$q=e^{i\pi t}$ by $t\in(-1, 1]$
and the $L$-root   $q^{1/L}=e^{it\pi/L}$.\bigskip

\subsection{Specialization of Weyl modules $V_\lambda(q)$ and simplicity in the closed Weyl alcove $\overline{\Lambda^+(q)}$.}\label{18.5}

\noindent Let $V_\lambda(x)$ be
the  irreducible representation of $U_x({\mathfrak g})$ with highest weight $\lambda$ and let $v_\lambda$ be the highest weight  vector of $V_\lambda(x)$ and $V_{\lambda, {\mathcal A}'}$   the associated free ${\mathcal A}'$-module
with   basis $B_\lambda$ introduced   in Sect. \ref{18.3}. Then we have a
specialized complex  ${ U}_q({\mathfrak g})$-modules   at a complex number $q$
$$V_\lambda(q):=V_{\lambda, {\mathcal A}'}\otimes_{{\mathcal A}'}{\mathbb C}.$$ 
This is   a cyclic module for $U_q({\mathfrak g})$ generated by $v_\lambda(q)=v_\lambda\otimes 1$, simple
if $q$ is not a root of unity, but it is not always
  so if $q$ is a root of unity.

 Let as before $\ell'$ be the order of $q$ and $\ell$ the order of $q^2$.  The   {\it linkage principle} gives information on irreducibility  of $V_\lambda(q)$ at primitive roots of unity. We
   consider the      affine Weyl group $W_{\ell'}$,  and   its  {\it translated } action on the real vector space $E$ spanned by the roots, defined by $w.x=w(x+\rho)-\rho$.    Let $d$ denote the ratio between the squared lengths of the longest to the shortest root, so $d=1$ for Lie types $ADE$, $d=2$ for $BCF$ and $d=3$ for $G_2$. 
    The structure of  $W_{\ell'}$ depends on the parity and divisibility by $d$ of the order  $\ell'$ of $q$.

We recall that affine Weyl group   $W_{\ell'}$ is the group of isometries of ${\mathfrak h}^*$ generated by reflections in the hyperplanes $\{x\in E: \langle x, \alpha_i\rangle=k\ell_id_i\}$, where $k\in{\mathbb Z}$, $\ell'_i$ is the order of $q^{d_i}$, and $\ell_i$ is $\ell'_i$ or $\ell'_i/2$ according to whether $\ell'_i$ is odd  or even.
  The translated action   admits a fundamental domain, called the principal Weyl alcove.
  The open Weyl alcove is defined by 
 $${\Lambda^+(q)}:=\{ \lambda\in\Lambda^+: \langle   \lambda+\rho, \theta\rangle < \ell\},$$ 
 if $d$ divides $\ell$ or
 $${\Lambda^+(q)}:=\{ \lambda\in\Lambda^+: \langle   \lambda+\rho, \theta_s\rangle < \ell\},$$
 if $d$ does not divide $\ell$, with $\theta$ the highest root and $\theta_s$ the highest short root in the root system.
 The element $\rho$ denotes the Weyl vector, defined as the sum of the fundamental weights.
 The linkage principle then    implies that $V_\lambda(q)$ is  simple   for  
  $\lambda\in \overline{\Lambda^+(q)}$. The following result is well known see \cite{Chari_Pressley, Sawin} and references therein.

\begin{thm}\label{irreducibles_of_fusion_category}
The quantized Weyl modules $V_\lambda(q)$  with   $\lambda\in \Lambda^+(q)$
  have positive quantum dimension provided $q$ satisfies Def. \ref{large_enough}.
The indecomposable tilting modules $T_\lambda$ 
 with $\lambda$ not in  $\Lambda^+(q)$ are negligible. 
 The open Weyl alcove $\Lambda^+(q)$ labels the simple objects of ${\mathcal C}({\mathfrak g}, q, \ell)$.
 \end{thm} 
 
 The quantum dimension can be defined using the rigidity equations. The explicit computation
 for quantized Weyl modules $V_\lambda(q)$ has roots in the classical Weyl dimension formula
 and is given by
 $${\rm dim}_q(V_\lambda(q))=\Pi_{\lambda\in\Phi^+}\frac{q^{\langle \lambda+\rho, \alpha\rangle}-q^{-\langle\lambda+\rho, \alpha\rangle}}{q^{\langle\rho, \alpha\rangle}-q^{-\langle\rho, \alpha\rangle}}.$$
 The lower bound on the order $\ell$ of $q^2$ required in Def. \ref{large_enough} is exactly that each of the denominators
 in the above formula does not vanish.
 We   recall the notion of tilting modules and negligible modules
  in the following subsection \ref{18.6}, and the construction of ${\mathcal C}({\mathfrak g}, q, \ell)$ in the following
  section \ref{74},
 see  also Prop. 2.4 in \cite{Wenzl}, see also   \cite{Andersen, Sawin, Chari_Pressley}.

\medskip

   \subsection{The tilting  category ${\mathcal T}({\mathfrak g}, q, \ell)$.}\label{18.6}
   
    In this subsection we assume $\ell'<\infty$.
 Constructions due to \cite{Andersen, Andersen1, GK, RT}, give rise to a semisimple,   ribbon, fusion category,     ${\mathcal T}({\mathfrak g}, q, \ell)$ that we briefly outline. Notice that the   constructions impose no restriction on the order  $\ell'$
   of $q$, and depend on the order $\ell$ of $q^2$.

  Since a Weyl module $V_\lambda(q)$   may fail to be irreducible,
Andersen developed   the
  notion of tilting module,  see Sect. 1, 3 \cite{Sawin}, Ch. 11.2 \cite{Chari_Pressley} and references to the original papers.
  
 A tilting module is a finite dimensional representation $W$ of $U_q({\mathfrak g})$ admitting together with its dual,
 a Weyl filtration, i.e. a sequence of modules $\{0\}\subset W_1\subset\dots \subset W$ such that
$W_{i+1}/W_i$ is isomorphic to a Weyl module $V_{\lambda_i}(q)$ with $\lambda_i\in\Lambda^+$.
  Weyl filtrations are non-unique, but for all filtrations of $W$ the number of  factors isomorphic to a given 
$V_\lambda(q)$  is unique, and it is in fact given by the multiplicity of $V_\lambda(x)$ in $W(x)$ if $W$ is obtained from a specialisation  $x\to q$ of a module 
$W(x)$ of $U_x({\mathfrak g})$, see Prop. 3 and Remark 2 in  \cite{Sawin} for a  precise statement.

 By Sect. 11.3 \cite{Chari_Pressley} or Cor. 5 in \cite{Sawin}, every tilting module decomposes into a direct sum of indecomposable tilting modules, and every indecomposable tilting module is isomorphic to a unique indecomposable tilting module $T_\lambda$ with maximal vector of weight $\lambda$,
  with $\lambda\in\Lambda^+$. 
Thus $T_\lambda$  has a    filtration by submodules $0\subset V_\lambda\subset V_2\subset
 V_3\subset
 \dots \subset T_\lambda$ such that $V_2/V_\lambda\simeq V_\mu(q)$, $V_3/V_2\simeq V_\nu(q), \dots$ with $\lambda>\mu>\nu\dots$, \cite{Chari_Pressley} p. 363, and   the dual $T_\lambda^*$ has a similar filtration and is isomorphic to  $T_{{-w_0}\lambda}$.
  Moreover   tensor products of tilting modules is tilting, thus  the category ${\mathcal T}({\mathfrak g}, q, \ell)$ of tilting modules becomes a tensor category  with duals \cite{Chari_Pressley, Sawin}.
 
It is important for us that  every Weyl module $T_\lambda=V_\lambda(q)$ for $\lambda\in\overline{\Lambda^+(q)}$,
  by   e.g. \cite{Andersen_Stroppel} Subsect. 1.1.
   This in particular implies  that   the multiplicities of the dominant weights of the factors in the Weyl filtrations of   tensor products $V_{\lambda_1}(q)\otimes\dots\otimes V_{\lambda_n}(q)$   with   $\lambda_i\in \overline{\Lambda^+(q)}
$  are the same as those determined by decomposition into irreducibles of the corresponding tensor product in the classical (or generic) case.

It follows from Sebsect. \ref{18.1}, see also
  Theorems 3, 4 in \cite{Sawin}, that the category of tilting modules over $U_q({\mathfrak g})$ is a ribbon category.
  For a fixed choice of $q^{1/L}$, 
  the corresponding  $R$-matrices define corresponding braided symmetries for the representation category, for more details on the classification in the type $A$ case, and references see Sect. \ref{KW}.

\begin{prop}\label{ribbon_specialized_case} Properties $(a)$, $(b)$,   of Theorem \ref{Properties_of_R} hold for the coproduct and $R$-matrix 
of $U_q({\mathfrak g})$.
The ribbon element $v$ of  $U_q({\mathfrak g})$ acts as the constant $v_\lambda=q^{-\langle \lambda, \lambda+2\rho\rangle}$ for $\lambda\in\overline{\Lambda^+(q)}$, with
  $\overline{\Lambda^+(q)}:=\{\lambda\in \Lambda^+: \langle \lambda+\rho, \theta\rangle \leq\ell \}$ the closed
   Weyl alcove.  
   \end{prop}

      \bigskip

\section{Fusion categories ${\mathcal C}({\mathfrak g}, q, \ell)$ and   unitary ribbon wqh   via integral wdf}\label{74}
 
 \subsection{The open Weyl alcove $\Lambda^+(q)$ and negligible tilting modules. }\label{19.1}   We follow  Gelfand and Kazhdan for the construction of the quotient category  \cite{GK}.
   Every object of ${\mathcal T}({\mathfrak g}, q, \ell)$  
  decomposes into a direct sum of indecomposable submodules, and this decomposition is unique up to isomorphism.
 One can   form two full linear  subcategories, ${\mathcal T}^0$, and ${\mathcal T}^\perp$ of 
 ${\mathcal T}({\mathfrak g}, q, \ell)$, with objects, respectively, 
  those representations which can be written as direct sums of $V_\lambda$, with $\lambda\in\Lambda^+(q)$ only,
   and those which have no such $V_\lambda$ as a direct summand.\medskip
   
  The objects of ${\mathcal T}^\perp$ and ${\mathcal T}^0$ are called   negligible and non-negligible, respectively. A morphism $T:W\to W'$ of ${\mathcal T}({\mathfrak g}, q, \ell)$ is called negligible 
if it is a sum of morphisms that factor through $W\to N\to W'$ with $N$ negligible. \medskip

The category ${\mathcal T}^\perp$ of negligible modules satisfies the following properties, 
 \cite{Andersen, GK}, that we call Gelfand-Kazhdan properties,  \medskip
 
 \noindent (1) Any object $W\in{\mathcal T}({\mathfrak g}, q, \ell)$ is isomorphic to a direct sum $W\simeq W_0\oplus N$ with $W_0\in{\mathcal T}^0$ and
  $N\in{\mathcal T}^\perp$. 
  
  \noindent(2) For any pair of morphisms $T:W_1\to N$, $S:N\to W_2$ of ${\mathcal T}({\mathfrak g}, q, \ell)$, with $N\in{\mathcal T}^\perp$, $W_i\in{\mathcal T}^0$, then $S T=0.$
  
  \noindent (3) For any pair of objects   $W\in{\mathcal T}_\ell({\mathfrak g})$, $N\in{\mathcal T}^\perp$, then   
  $W\otimes N$ and $N\otimes W\in{\mathcal T}^\perp$.\medskip
  
 We shall extensively use the previous properties to construct the weak  Hopf algebras $A_W({\mathfrak g}, q, \ell)$ in Sect.
 \ref{20}.

Property $(1)$ follows easily from the mentioned decomposition of objects of ${\mathcal T}({\mathfrak g}, q, \ell)$, while property $(2)$ means that no non-negligible module can be a summand of a negligible one (however, it can be a factor of a Weyl filtration of a negligible). 

We recall  that negligible indecomposable tilting modules are characterized by the property of having zero quantum dimension.  
 A morphism  $T: W\to W'$ is   negligible if and only if  ${\rm Tr}_W(ST) = 0$ for all morphisms $S: W'\to W$.

\begin{rem} In particular, for $\lambda$, $\mu\in\Lambda^+(q)$,
$$V_\lambda(q)\otimes V_\mu(q)=\oplus_{\gamma\in \Lambda^+(q)}\tilde{m}_\gamma^{\lambda, \mu} V_\gamma(q)\oplus N$$
with $N\in{\mathcal T}^\perp$. The multiplicities $\tilde{m}_\gamma^{\lambda, \mu}$ are given by the quantum Racah
formula, Sect. 5 \cite{Sawin}. In particular, $\tilde{m}_\gamma^{\lambda, \mu}$ depends on the affine Weyl group $W_{{\ell}'}$. It follows from the description of $W_{{\ell}'}$ given in  Lemma 1 of \cite{Sawin} that this group
depends only on ${\mathfrak g}$ and the order $\ell$ of $q^2$, thus the same holds for $\tilde{m}_\gamma^{\lambda, \mu}$.
\end{rem}

 \subsection{The quotient category ${\mathcal C}({\mathfrak g}, q, \ell)$.}\label{19.8}

Let $\text{Neg}(W, W')$ be the subspace of negligible morphisms of   $(W, W')$.
Then the quotient category, ${\mathcal C}({\mathfrak g}, q, \ell)$, is the category with the same objects as 
${\mathcal T}({\mathfrak g}, q, \ell)$ and morphisms between the objects $W$ and $W'$ the quotient space,
$$(W, W')_{{\mathcal C}({\mathfrak g}, q, \ell)}:=(W, W')/\text{Neg}(W, W').$$
 Gelfand and Kazhdan endow ${\mathcal C}({\mathfrak g}, q, \ell)$ with the  unique structure of a tensor category   such that the quotient map ${\mathcal T}({\mathfrak g}, q, \ell)\to{\mathcal C}({\mathfrak g}, q, \ell)$ is a tensor functor.  The   tensor product  of  objects and  morphisms of ${\mathcal C}({\mathfrak g}, q, \ell)$ is usually denoted   by 
  $W\underline{\otimes} W'$
 and
 $S\underline{\otimes} T$ respectively,
 and referred to  as the truncated tensor product  in the physics literature. This is now a
  semisimple tensor category and 
  $\{V_\lambda, \, \lambda\in\Lambda^+(q)\}$ is a complete set of irreducible objects.
  
In ${\mathcal C}({\mathfrak g}, q, \ell)$ we have
 $$V_\lambda\underline{\otimes } V_\mu\simeq \oplus_{\nu\in\Lambda^+(q)} \tilde{m}^\nu_{\lambda, \mu} V_\nu.$$
 Notice that  this decomposition of $V_\lambda\otimes V_\nu$  is unique up to isomorphism but   not canonical
 (cf. \cite{Wenzl}, and also Sect. 11.3C in \cite{Chari_Pressley}  and references therein.)
 \medskip

   The ribbon structure of ${\mathcal C}({\mathfrak g}, q, \ell)$   is induced by that of the tilting category.
      Also the formulas for the fusion
coefficients and quantum dimensions are well-known, and regulated by the affine Weyl group in the sense mentioned in Sect. \ref{18.2}, see  Sect. 2, 5 of \cite{Sawin}, but we shall only need them in some special cases later on, so we refrain from recalling
them in full generality. 
Those fusion rules indeed pass to the ${\mathcal C}({\mathfrak g}, q, \ell)$ and give the fusion rules for this category,
However, it will be important for us to recall that ${\mathcal C}({\mathfrak g}, q, \ell)$
depends on $q$
but  the Grothendieck semiring $R({\mathcal C}({\mathfrak g},  q, \ell))$ depends only on $\ell$. 
We shall refer to $R({\mathcal C}({\mathfrak g},  q, \ell))$ as the Verlinde fusion ring.
 \medskip

 \subsection{Modularity and unitarizability.}\label{19.9}
 Further properties of modularity   ${\mathcal C}({\mathfrak g}, q, \ell)$ depend on on the choice of $q^{1/L}$ as a primitive root of unity of order $\ell'L$ and on
 the order $\ell'$ of $q$.   We refer to the     papers by Rowell and Sawin
 \cite{Rowell2, Sawin}   for a detailed treatment. For example the cases where $2d|\ell'$ give modular categories
 and this is the case of most physical interest, and also that meeting
   the purpose of our paper.
   
  More in particular, 
   we shall
 mostly be interested in the``minimal roots'' .
 
 \begin{defn}\label{minimal_root}
 Let $q$ be a complex root of unity and let $\ell$ be the order of $q^2$.
  We shall say that $q$ is a {\it minimal root} if
 $q$ is of the form
 $$q=e^{i\pi/\ell},  \quad q^{1/L}=e^{i\pi/\ell L}, \quad d|\ell.$$ 
 \end{defn}
 
Indeed if $q$ is a minimal root and if in addition $\ell$ is large enough as in Def. \ref{large_enough}, then ${\mathcal C}({\mathfrak g}, q, \ell)$   becomes a unitary ribbon fusion category 
(that we sometimes denote by ${\mathcal C}^+({\mathfrak g},  q, \ell)$ to stress the unitary structure) for all levels $k=\ell/d-h^\vee\geq 1$ for ${\mathfrak g}\neq E_8$ 
and $k\geq 2$ for ${\mathfrak g}=E_8$
by  \cite{Wenzl}.

\bigskip

\subsection{A general construction, the wqh algebras $A({\mathfrak g}, q, \ell)$.}\label{19.10} We   introduce the function $D$ on the Grothendieck ring of ${\mathcal C}({\mathfrak g}, q, \ell)$, which assigns 
the vector space dimension of the corresponding representation of ${\mathfrak g}$ to each irreducible $\lambda\in\Lambda^+(q)$.
It follows easily  from the quotient construction and from the fact that every tilting module decomposes uniquely up to isomorphism into  a direct sum of indecomposable tilting modules,  that $D$ is indeed a weak dimension function on 
${\mathcal C}({\mathfrak g}, q, \ell)$. We shall refer to it as the {\it classical dimension function}.
We may then   apply Theorem \ref{propweakdim}  and we have, up to isomorphism and twist, a finite dimensional   weak quasi-Hopf 
$C^*$-algebra $A({\mathfrak g}, q, \ell)$.
We next fix a root of unity of the form $q=e^{i\pi/\ell}$ with $d | \ell$. Then by
     \cite{Wenzl, Xu_star},   and  
Theorem \ref{HO_star} $A({\mathfrak g}, q, \ell)$ becomes a unitary weak quasi-Hopf algebra.
\bigskip

We shall return to specific constructions of weak Hopf algebras associated to ${\mathcal C}({\mathfrak g}, q, \ell)$ in Sects.
\ref{unitary_structure_of_fusion_category} and \ref{20}.

  \section{Kazhdan-Lusztig-Finkelberg theorem for $\tilde{\mathcal O}_\ell$, unitarity of Kirillov inner product}\label{23}

In this section we state the   Kazhdan-Lusztig-Finkelberg theorem for the category $\tilde{\mathcal O}_\ell$ of certain modules
of affine Lie algebra $\hat{{\mathfrak g}}_k$ associated to a simple Lie algebra ${\mathfrak g}$ at positive integer level $k$, proved
by Finkelberg and prove positivity of Kirillov hermitian form on the Beilinson-Feigin-Mazur category  in Theorem \ref{BFM}, solving Problem \ref{problem_kirillov}.

 We consider the linear category  $\tilde{\mathcal O}_\ell$ of $\hat{\mathfrak g}$-modules of finite length, with central charge the positive integer level $k$, $\ell=d(k+h^\vee)$,   which are integrable in the sense of Kac \cite{Kac2}. It is known that 
  $\tilde{\mathcal O}_\ell$ is   semisimple and the simple object are the level $k$ modules $L_{\lambda, k}$ defined in Subsect. \ref{30.1}.  
  Let  $\tilde{\mathcal O}_\ell$ be endowed with the braided tensor structure by Beilinson, Feigin, and Mazur \cite{BFM}. Let $\tilde{\mathcal O}_{-\ell}$ the subquotient fusion category of Kazhdan-Lusztig category   ${\mathcal O}_{-\ell}$ described in Subsect. \ref{1.1}.
  Then the following result has been shown in \cite{Finkelberg}, \cite{Finkelberg_erratum}.

   \begin{thm}  \label{Finkelberg} Let ${\mathfrak g}_k$ be an affine Lie algebra at positive integer level $k$ such that
   for ${\mathfrak g}=AD$, $k\geq1$,  for ${\mathfrak g}=BCFGE_8$, $k\geq3$, for ${\mathfrak g}=E_6, E_7$, $k\geq2$.
  Then category $\tilde{\mathcal O}_\ell$ endowed with Beilinson-Feigin-Mazur   braided tensor structure
is equivalent to $\tilde{\mathcal O}_{-\ell}$ as a ribbon fusion category and is rigid. As a consequence $\tilde{\mathcal O}_\ell$ is equivalent to ${\mathcal C}({\mathfrak g}, q, \ell)$.
  
  \end{thm}

More historical information may be found in \cite{Huang2018}, \cite{HL},  Appendix B of \cite{On_a_problem_posed_by_Huang}. For a historical perspective unifying Finkelberg equivalence with the braided Doplicher-Roberts program as considered in this paper, we refer to \cite{pinzari_survey}.

      Kirillov \cite{Kirillov}, \cite{Kirillov3} defined a hermitian modular tensor category  $\tilde{\mathcal O}_\ell$   starting from
   a modular tensor structure  of $\tilde{\mathcal O}_\ell$,  conjectured unitary. To the authors' knowledge, the following result is new.
     
   \begin{thm}\label{BFM}
Under the same assumptions as in Finkelberg's theorem \ref{Finkelberg},  $\tilde{\mathcal O}_\ell$ becomes a unitary modular fusion category with Kirillov inner product.
   \end{thm}
   
   \begin{proof}
   We apply Theorem \ref{unitarizability} to ${\mathcal C}=\tilde{\mathcal O}_\ell$, ${\mathcal C}^+$ the linear $C^*$-category with objects the affine Lie algebra modules with Kac inner product, ${\mathcal D}^+={\mathcal C}({\mathfrak g}, q, \ell)$ endowed with Wenzl's unitary structure \cite{Wenzl}, and to Finkelberg's braided tensor equivalence between these categories.
   \end{proof}
 
   In  Corollary \ref{cor_Zhu_as_a_compatible_unitary_wqh}  we    give an independent constructive modular tensor structure 
 ${\rm Rep}_{\rm T}(V_{{\mathfrak g}_k})$  on the module category of the affine vertex operator algebra for all Lie types and all positive integer levels with a compatible unitary coboundary structure on the Zhu algebra $A(V_{{\mathfrak g}_k})$, and a proof   of the equivalence with the quantum group fusion category ${\mathcal C}({\mathfrak g}, q, \ell)$.
Moreover Theorem
 \ref{Finkelberg_HL} gives a constructive proof of the equivalence with respect to Huang-Lepowsky ribbon braided tensor structure
 on ${\rm Rep}(V_{{\mathfrak g}_k})$ for the classical Lie types and $G_2$.

We may thus apply the abstract unitarizability Theorem \ref{unitarizability} to these structure and establish the following constructive result.

 \begin{thm}\label{Kirillov} 
  Let ${\mathfrak g}_k$ be an affine Lie algebra at positive integer level $k\geq1$ of type $ABCDG_2$. Then
Kirillov hermitian structure on ${\rm Rep}_{\rm HL}(V_{{\mathfrak g}_k})$ is a unitary modular tensor structure
 for Huang-Lepowsky structure.
 Let $F: u\otimes v\to\overline{u\boxtimes v}$ define Huang-Lepowsky tensor product module. Let $F_0: u[0]\otimes v[0]\to(\overline{u\boxtimes_{\rm HL} v})[0]$, $G_0$ be defined as in Theorem \ref{Zhu_from_qg_if_of_CFT_type0}, via inclusion  and projection   of $F$ of the top level ${\mathfrak g}$-modules.
 Then $(F_0)_{V_\lambda, V}(F_0)_{V_\lambda, V}^*=1$, that is the pair $(F_0, G_0)$ is a strongly unitary wqh structure for Zhu functor $Z$ on pairs of the form $(V_\lambda, V)$.
 
 For the Lie types $EF$ the same conclusions hold for ${\rm Rep}_{\rm T}(V_{{\mathfrak g}_k})$ and $k\geq1$.
 \end{thm}

   To show the above theorem we need to study the inner structures and consider naturally associated integral dimension functions. This is the content of the rest of the paper.

  \section{Vertex operator algebras,  the linear category ${\rm Rep}(V)$   and    Zhu algebra $A(V)$}\label{VOAnets}
  
   In this section we   recall the basic theory of vertex operator algebras, their modules.
   We give a reformulation in a Tannakian framework of a result that has roots in Zhu's work on the correspondence between $V$-modules and
  modules of the associated Zhu algebra $A(V)$. 
   
We also recall the important rationality conditions on $V$ under which
 Huang   proved modularity of the category ${\rm Rep}(V)$ of $V$-modules.

    In the next section we recall the basic theory of intertwining operators
  and Huang and Lepowsky's tensor product theory.
 We shall discuss applications in 
 Sects.    \ref{21},  \ref{23}, \ref{22}.
  We refer the reader to textbooks  \cite{FLM}, \cite{FHL}, \cite{Kac} for the general theory of vertex operator algebras.

 \begin{defn}
 A {\it vertex operator algebra} is a ${\mathbb Z}$-graded vector space $V$
 $$V=\bigoplus_{n\in{\mathbb Z}}V_{(n)},$$
with finite dimensional homogeneous spaces $V_{(n)}$ and $V_{(n)}=0$ for $n$ sufficiently small,
together with a triple $(Y, 1, \nu)$, where $Y$ is the {\it state-field correspondence}, 
a linear map
$$Y: V\to{\rm End}(V)[[x, x^{-1}]], \quad\quad a\to Y(a, x)=\sum_{n\in{\mathbb Z}} a_{(n)}x^{-n-1},$$
with $x$ is a formal variable. The element $Y(a, x)$ is called the {\it vertex operator} associated with $a$.
The element $1$ lies in $V_{(0)}$ and  is called the {\it vacuum}, and   $\nu\in V_{(2)}$ is
called  the {\it Virasoro vector}.
 The following conditions form the definition for all $a$, $b\in V$,\medskip

  \noindent {\bf a)}  {\it (lower truncation condition)}:   $a_{(n)}b=0$ for $n$ sufficiently large.\smallskip
  
    \noindent {\bf b)} {\it (vertex operator associated to the vacuum)} $Y(1, x)={\rm id}_V$\smallskip
    
        \noindent {\bf c)} {\it (creation property)} $Y(a, x)1\in V[[x]]$ (power series with non-negative integral powers)
        and $$\lim_{x\to0} Y(a, x)1=a,$$ \smallskip
        
        \noindent {\bf d)} {\it (Virasoro algebra relations and spectrum condition for $L_0$)} Set $$L_n=\nu_{(n+1)}$$, i.e. 
        $Y(\nu, x)=\sum_{n\in{\mathbb Z}} L_n x^{-n-2}$, then 
       \begin{equation}\label{Virasoro_relations} [L_m, L_n]=(m-n)L_{m+n} +\frac{1}{12}(m^3-m)\delta_{m+n, 0} c,
       \end{equation} with $c\in{\mathbb C}$ the 
        {\it central charge} of $V$.
        Moreover
        \begin{equation}\label{L_0}
       L_0a=n a, \quad\quad a\in V_{(n)},\end{equation}
  $n$ is called the {\it conformal weight}, {\it energy}, or {\it degree} of the vector $a\in V_{(n)}$, and denoted ${\rm deg}(a)$;
   $L_0$
  is the {\it energy operator}, or {\it conformal Hamiltonian} on $V$. \smallskip
        
        \noindent {\bf e)}  {\it $L_{-1}$-derivative (or translation) property}  
        \begin{equation}\label{translation}\frac{d}{dx}Y(a, x)=Y(L_{-1}a, x),\end{equation}
        \smallskip

       \noindent{\bf f)} {\it (Jacobi identity)} 
       for $l$, $m$, $n\in{\mathbb Z}$,
       $$\sum_{i=0}^\infty \binom{m}{i}
       Y(a_{(l+i)}b, x)x^{m+n-i}=$$
       $$\sum_{i=0}^\infty (-1)^i\binom{l}{i}a_{(m+l-i)}Y(b, x)x^{n+i}-
       \sum_{i=0}^\infty(-1)^{l+i}\binom{l}{i}Y(b, x)a_{(m+i)}x^{n+l-i},$$
       
        \end{defn}
             This completes the definition.  The formal limit    $\lim_{x\to0}$ and  the formal derivative  $\frac{d}{dx}$ are defined in the natural way.
      On a fixed vector of $V$, only finitely many addenda appearing in the series of the Jacobi identity f) are possibly non-zero, thus both sides of the identity are well-defined.  
           The axioms of vertex operator algebra have a number of consequences.
  A special case of the Jacobi identity for $l=n=0$ gives 
     \begin{equation}
     [a(m), Y(b, x)]=\sum_{i=0}^\infty \binom{m}{i}Y(a(i)b, x)x^{m-i},\end{equation}
     thus   the operators $a_{(m)}$ are in particular closed under the Lie bracket.
     Moreover, together with  
     the translation property e) of $L_{-1}$ this equation
     implies for $L_{-1}=\nu_{(0)}$ and $L_0=\nu_{(1)}$,
     $$[L_{-1}, Y(b, x)]=\frac{d}{dx}Y(b, x), \quad\quad [L_0, Y(b, x)]=\frac{d}{dx}Y(b, x)x+Y(L_0b, x).$$
  The right hand side  together with the eigenvalue property (\ref{L_0}) of $L_0$ and  implies the useful grading relation
     \begin{equation}\label{grading_equation} b(n): V_{(m)}\to V_{(m+{\rm deg}(b)-n-1)}.
     \end{equation}
     It also follows that  $a(n)1=0$ for $n\geq0$.
             
The Jacobi  identity can be written in readily equivalent forms, and we have chosen one of them.
A compact form is the following
\begin{equation}\label{Jacobi_VOA_delta}
x_0^{-1}\delta(\frac{x_1-x_2}{x_0})Y(a, x_1)Y(b, x_2)-x_0^{-1}\delta(\frac{x_2-x_1}{-x_0}) Y(b, x_2)Y(a, x_1)=
\end{equation}
 $$x_2^{-1}\delta(\frac{x_1-x_0}{x_2})Y(Y(a, x_0)b, x_2)),$$
where $\delta(x)=\sum_{n\in{\mathbb Z}} x^n$, we refer the reader to the references at the beginning of the section
for complete explanation of the use of the $\delta$.
 There is yet another set of properties equivalent to the Jacobi identity, in presence of the other axioms of vertex operator algebra.
       This set properties are   the {\it rationality of products, iterates,  commutativity of products,
       and  associativity}, the latter being a property involving products and iterates. We refer the reader
       to Sect. 8.10 in \cite{FLM} and also
       Sect. 2 in \cite{HL_tensor_products_of_modules}. This formulation of the Jacobi identity turns out useful
       for the development of tensor product theory by Huang and Lepowsky and we encourage a non-expert reader to 
       consider this viewpoint for insight. We shall come back to this in Sect. \ref{32}.
       We next give the definition of a $V$-module.

   \begin{defn}\label{V_module} ({\it $V$-module})   Given a vertex operator algebra $(V, Y, 1, \nu)$, a $V$-module  is
      a pair $(M, Y_M)$, with $M$ an ${\mathbb R}$-graded vector space
       $$M=\bigoplus_{r\in{\mathbb R}}M_{r},$$
by finite dimensional subspaces $M_{r}$ that eventually vanish for $r$ small and 
a linear map
$$Y_M: V\to{\rm End}(M)[[x, x^{-1}]], \quad\quad a\to Y_M(a, x)=\sum_{n\in{\mathbb Z}} a_{(n)}^Mx^{-n-1},$$
defining a vertex operator $Y_M(a, x)$ associated with $M$, satisfying   properties analogous to the defining
properties of the vertex operators
$Y(a, x)$ associated to $V$, except for the creation property of the vacuum vector, that is omitted.
In particular,  the endomorphism $a_{(l+i)}$ appearing at the left hand side of the Jacobi identity f)   and the operator $L_{-1}$ defining the  derivative property for the $V$-module $M$ are associated to $V$,
while the remaining equations in the definition involve vertex operators $Y_M(a, x)$ associated to $M$.
The   constant $c$ in the Virasoro algebra relations for $M$ equals the central charge of $V$.
\end{defn}

Let $M$ be a $V$-module. 
 If $\nu \in V$ is the conformal vector we write 
\begin{equation}
Y_M(\nu,x) = \sum_{n \in \mathbb{Z}} L^M_{n}x^{-n-2}\,, \quad a\in V\,.  
\end{equation}
In particular $L^M_0$ denotes the conformal Hamiltonian on $M$. 
For $m\in M_r$, the spectral property of $L^M_0$ holds,
\begin{equation}\label{conformal_weight_of_vector} L_0^M(m)=\Delta_m m, \end{equation} 
where $\Delta_m:=r$   is called the {\it conformal weight} or {\it energy} of the homogeneous vector $m$.
For $M=V$ and  $a\in V_{(n)}$, $\Delta_a={\rm deg}(a)=n$   in the notation of \cite{Zhu}.

\begin{defn}
\label{Rep(V)} 
({\it The category} ${\rm Rep}(V)$) Given   $V$-modules $M_1$ and $M_2$, a morphism $T\in {\rm Hom}_V(M_1, M_2)$ is 
is a linear map satisfying $TY_{M_1}(a, x)=Y_{M_2}(a, x)T$ for all $a\in V$. In particular,
for $a=\nu$ this condition shows that $T$ preserves the gradings of the modules.
This defines the category  ${\rm        Rep}(V)$ of $V$-modules.
\end{defn}

We wish to define a canonical functor $\mathcal{F}_V:  {\rm        Rep}(V) \to  {\rm        Vec}$. When the assumption  in 
Theorem \ref{propweakdim} are satisfied then, thanks to the Tannaka-Krein duality result in Theorem \ref{TK_algebraic_quasi} 
we will be able to associate a weak quasi-Hopf algebra to ${\rm        Rep}(V)$. 

The defining properties of a $V$-module imply, similarly to the case of a vertex operator algebra (\ref{grading_equation}).
The Jacobi identity for a $V$-module for $l=n=0$ gives 
     \begin{equation}\label{consequence_of_Jacobi_for_modules}
     [a^M(m), Y_M(b, x)]=\sum_{i=0}^\infty \binom{m}{i}Y_M(a_{(i)}b, x)x^{m-i},\end{equation}
     and implies
      \begin{equation}\label{grading_equation_for_modules2}
      [L_{-1}^M, Y_M(b, x)]=\frac{d}{dx}Y_M(b, x), \quad\quad [L_0^M, Y_M(b, x)]=\frac{d}{dx}Y_M(b, x)x+Y_M(L_0b, x).
      \end{equation}
and we similarly   derive from the right hand side and the spectral property of $L_0^M$,
the following grading relations for the {\it modes}, or {\it coefficients} $b_{(n)}^M$ of the associated vertex operators, 
for  $b\in V$  homogeneous

 \begin{equation}\label{grading_equation_for_modules}
 b^M_{(n)}: M_{r}\to M_{r+{\rm deg}(b)-n-1} \quad\text{for}\quad  r\in{\mathbb R}, \quad n\in{\mathbb Z}.\end{equation}
 It follows particular (or from the Virasoro algebra relations)
 $$L^M_n=\nu_{(n+1)}: M_r\to M_{r-n}.$$
 Every $V$-invariant subspace of $M$ is a $V$-module, see    Remark 1.5 in \cite{GuiI}. It follows from (\ref{grading_equation_for_modules}) that every $V$-module $M$ can be canonically decomposed compatibly with 
the original grading of $M$ as a direct sum of $V$-modules
$$M=\bigoplus_{s\in S}M^{(s)},$$
where $M^{(s)}=\bigoplus_{n \in \mathbb{Z}_{\geq 0}} M_{s+n}$, and $S$ is the collection of 
$s\in{\mathbb R}$ such that $M_{s-n}=0$ for all $n\in{\mathbb N}$ and $M_s\neq0$
 \cite{Zhu}.

Moreover, $M$ is irreducible if and only if 
it has no $V$-invariant subspace.
If $M$ is irreducible then $S$ consists of a single element
\begin{equation}\label{conformal_weight_of_module}
 \Delta_M:=\inf\{r\in{\mathbb R}: M_{r}\neq0\}\in S,
\end{equation}
 called the {\it conformal weight} of the irreducible module $M$. 
We set, for $n\in{\mathbb Z}_{\geq 0}$,
\begin{equation}\label{grades}
M_{(n)}:=M_{n+\Delta_M}={\rm        Ker}(L^{M}_0 - (\Delta_M + n) 1_{M}).
\end{equation} If $M$ is irreducible
then  
$$M=\bigoplus_{n\in{\mathbb Z}_{\geq 0}}M_{(n)},$$
We shell refer to $M_{(0)}$ the {\it lowest energy subspace}, or {\it top space} of $M$.

The   Virasoro algebra relations for $L^M_n$
(or the grading equation (\ref{grading_equation_for_modules}))
implies that for an integer $n\neq 0$ and $b\in M$ homogeneous
$L^M_0L^M_nb=(\Delta_b-n)L^M_nb$.
Let $M$ be irreducible.  We then have the following relations 
\begin{equation} \label{action_of_positive_Virasoro_on_lowest_energy_space}
L^M_0b=\Delta_Mb; \quad\quad L^M_nb=0, \quad\quad n>0; \quad\quad b\in M_{(0)}.
\end{equation}
This property will turn out useful in Sect. \ref{33} to discuss {\it primary fields}.
Thanks to  (\ref{consequence_of_Jacobi_for_modules}) and the first equation in
(\ref{grading_equation_for_modules2}), the conditions (\ref{action_of_positive_Virasoro_on_lowest_energy_space}) are equivalently written as
 \begin{equation}\label{differential_equation_for_primary_field}
 [L_m^M, Y_M(b, x)]=x^{m}(x\frac{d}{dx}Y_M(b, x)+(m+1)\Delta_MY_M(b, x)), \quad\quad m\in{\mathbb Z}.
 \end{equation}

\begin{rem}\label{Rationality_DLM}

\noindent 1) The study of rationality in conformal field theory
has roots in the work by Anderson and Moore \cite{Anderson_Moore}. A definition of rationality for a vertex operator algebra
was introduced by Zhu \cite{Zhu}. The definition was later simplified by Dong, Li and Mason, and we try to summarize some of their main results.
 
\noindent 2) The definition of a {\it weak module} is given in \cite{DLM1}. The weakness as compared to the Def. \ref{V_module} of $V$-module means that the  grading assumption
for $M$ and the spectral condition (\ref{L_0}) for $L^M_0$ are omitted for a weak module. Moreover the Virasoro and translation   relations  (\ref{Virasoro_relations}),  (\ref{translation}) are redundant for a weak module, and thus for a $V$-module, by Lemma 2.2
in \cite{DLM1}.

\noindent 3) Dong, Li and Mason call  a weak module   {\it ordinary} if it is endowed with a ${\mathbb C}$-grading   by finite dimensional subspaces $M_\lambda$ such that for each $\lambda\in{\mathbb C}$, $M_{\lambda+n}=0$ for $n$ small enough, and   the usual
spectral property for $L^M_0$ holds. They also call a weak module   {\it admissible} if it is graded by ${\mathbb Z}_{\geq0}$ and  the grading relation (\ref{grading_equation_for_modules}) holds for
modes of the associated vertex operators, for all $r\in{\mathbb Z}_{\geq0}$ and $n\in{\mathbb Z}$.  
Any ordinary module is admissible. If every admissible is a direct sum of irreducible admissible modules then every irreducible admissible is ordinary by \cite{DLM2}.

\noindent 4) Note that a more general  definition of $V$-module starts with a ${\mathbb C}$-graded vector space.
If each admissible $V$-module is completely reducible then $V$ has only finitely many inequivalent admissible modules, and every such module is ordinary \cite{DLM2}.

\noindent 5) By \cite{Anderson_Moore} and Theorem 11.3 in \cite{DongLiMason2}  the grading of a $V$-module
 is automatically in ${\mathbb R}$.
More in detail,   
 the conformal weight of an irreducible ordinary module is rational and the central charge is rational.
 
 \end{rem}

The rationality assumptions for $V$ in the sense of Dong, Li, Mason as in 3) of Remark \ref{Rationality_DLM}
 imply \cite{DLM2} in particular that 
there are only finitely many inequivalent irreducible $V$-modules, and that a $V$-module
$M$ can be written as a finite direct sum 
\begin{equation}
\label{moduleirrdecomposition}
M = \bigoplus_{i} M^i  
\end{equation} 
of irreducible V-modules $M^i$ compatibly with the grading,
by Theorem 8.1
\cite{DLM2}, see also Theorem 3.2 in \cite{DongLiMason2}. 
 For each $M^i$ we write
\begin{equation}
M^i = \bigoplus_{n \in \mathbb{Z}_{\geq 0}} M^i_{(n)}
\end{equation}
with $M^i_{(n)} = {\rm        Ker}(L^{M^i}_0 - (\Delta_{M_i}+ n) 1_{M^i})$ and $M^i_{(0)} \neq \{ 0\}$. Note that every $M^i_{(n)}$ is finite dimensional. We now define a finite dimensional subspace $M_{(0)} \subset M$ by 
\begin{equation}
M_{(0)} := \bigoplus_i M^i_{(0)} \,. 
\end{equation}
It is easy to see that $M_{(0)}$ is independent from the choice of the direct sum decomposition in Eq. (\ref{moduleirrdecomposition}). 
Moreover,
it can be shown that $\mathcal{U}(M,V)M_{(0)}=M$ where $\mathcal{U}(M,V)$ is the subalgebra of ${\rm        End}(M)$ generated by the vertex operator coefficients $a^M_{(n)}$, $a \in V$, $n \in \mathbb{Z}$. 

Now let $M^\alpha$ and $M^\beta$ be $V$-modules $M^\alpha$ and $M^\beta$ and 
$T: M^\alpha \to M^\beta$ a $V$-module homomorphism. Recall that from the equality $T L^{M^\alpha}_0 = L^{M^\beta}_0 T $  it follows that 
$TM^\alpha_{(0)} \subset M^\beta_{(0)}$.

We now define a linear functor $\mathcal{F}_V:  {\rm        Rep}(V) \to  {\rm        Vec}$ in the following way. If $M$ is an object in ${\rm        Rep}(V)$, i.e. a $V$-module, then $\mathcal{F}_V(M) = M_{(0)}$. If $T: M^\alpha \to M^\beta$ is a morphism in ${\rm        Rep}(V)$, i.e. a $V$-module homomorphism, then  $\mathcal{F}_V(T) = T\restriction_{M^\alpha_{(0)}}$.

If $\mathcal{F}_V(T) = 0$ then, 
$TM^\alpha = T\mathcal{U}(M^\alpha,V) M^\alpha_{(0)}  = \mathcal{U}(M^\beta,V) TM^\alpha_{(0)} = \{0\} $
so that $T=0$ and hence $\mathcal{F}_V$ is faithful. We are now in the position to apply  Theorem \ref{TK_algebraic_quasi}.  
Let $A(V):={\rm        Nat}_0(\mathcal{F}_V)$.

\begin{thm}
Let $V$ be a VOA such that every admissible module is a direct sum of simple admissible modules. Then
$A(V)$ is a semisimple associative algebra   that  can be identified with the Zhu's algebra of $V$.
Moreover, there is an equivalence $\mathcal{E}_V: {\rm        Rep}(V) \to {\rm        Rep}(A(V))$ 
which, after composition with the forgetful functor  $:{\rm        Rep}(A(V)) \to {\rm        Vec}$ is isomorphic to  $\mathcal{F}_V$. 
\end{thm}

\begin{defn}\label{Zhu_functor}  We refer to $\mathcal{F}_V:  {\rm        Rep}(V) \to  {\rm        Vec}$ as {\it Zhu's functor} and to  $\mathcal{E}_V: {\rm        Rep}(V) \to {\rm        Rep}(A(V))$ as {\it Zhu's equivalence}.
\end{defn}

\begin{rem}\label{Zhu_inverse} Zhu' functor ${\mathcal F}_V$  and equivalence ${\mathcal E}_V$ recalled at the end
of  Sect. \ref{VOAnets} play an important role in our paper.
 Moreover,  Zhu constructs a canonical linear equivalence  $${\mathcal S}_V: {\rm        Rep}(A(V))\to
 {\rm        Rep}(V), $$
that is   a right inverse of ${\mathcal E}_V$, 
$${\mathcal E}_V{\mathcal S}_V=1.$$
\end{rem}

The previous theorem is a reformulation of part of the work by Zhu, and Dong, Li and Mason
  in a Tannakian setting. By Theorem 8.1 in \cite{DLM2}, see also Theorem 3.2 in \cite{DongLiMason2},  \cite{Li}, every simple admissible is ordinary and $V$ admits
  only finitely many inequivalent simple ordinary modules. Thus $V$ is a
 rational vertex operator algebra in the sense of \cite{Zhu}. We then apply Theorems 2.1.2, 2.2.1, 2.2.2, and 2.2.3 in \cite{Zhu}. These results   say that $A(V)$ is associative and finite dimensional semi-simple under the assumptions of the previous theorem, and that ${\mathcal F}$ induces a linear equivalence of semisimple categories
 between the category of finite-dimensional $A(V)$ modules and the category of ordinary $V$-modules.
 Moreover an inverse equivalence is explicitly described in \cite{Zhu} and \cite{DLM2}. It should also be noted 
 that the functor ${\mathcal F}$ may be described intrinsically, i.e. without resorting to an irreducible decomposition, by Prop. 5.4 in \cite{DLM2}.

 \section{Vertex operator algebras and weak quasi-Hopf algebras, unitarizing
   ${\rm Rep}(V_{{\mathfrak g}_k})$}\label{VOAnets2}
 
In this section we   describe some general Tannakian constructions of weak quasi-Hopf algebras   from the theory of vertex operator algebras (VOAs).  This leads to some interesting applications to unitarizability
of important examples of module categories of vertex operator algebras including the affine examples and  to new questions. 

 We will restrict to VOAs and conformal nets whose representation category are known to be modular tensor categories. These are the rational VOAs satisfying the assumptions in \cite{Huang2} and the completely rational conformal nets  defined and studied in \cite{KLM}.
 
   The  main result of this section are Theorems \ref{TheoremUnitaryRepVOA} (a general result) and its most important application, 
Theorem \ref{TheoremUnitaryBraidRepAffine} which give  the construction of unitary structures on   ribbon braided tensor categories of   modules of  some vertex operator algebras, including the affine vertex operator algebras ${\rm Rep}(V_{{\mathfrak g}_k})$ at positive integer levels $k$ of classical Lie types and $G_2$. 

Our result is based on Wenzl work about the construction of unitary structures on fusion categories
 of quantum groups at certain roots of unity \cite{Wenzl} and on our proof of an analogue of Kazhdan-Lusztig-Finkelberg equivalence theorem in the setting of vertex operator algebras, with Huang-Lepowsky ribbon braided tensor category structure,
 Theorem \ref{Finkelberg_HL}.
To obtain these unitary structures,  weak quasi-Hopf algebras will play a role in transporting the unitary structure
from the setting of quantum groups to the setting of vertex operator algebras following the general construction of Sect. \ref{12}.
In this section we shall adopt general 
Tannakian constructions described in the previous sections.
Therefore at the end of this section, we shall conclude with a positive result on unitarization of those ${\rm Rep}(V_{{\mathfrak g}_k})$
for the Lie types stated in Theorem \ref{Finkelberg_HL} (the classical types and $G_2$), but we shall not give a full account on
the unitary structure that we obtain.
To this aim, we shall need more work on {\it canonical} Tannakian constructions of weak Hopf algebras associated to the quantum group fusion categories
${\mathcal C}({\mathfrak g}, q, \ell)$, and the complete  proof of Theorem \ref{Finkelberg_HL}. This will be done
in Sect. \ref{18}--\ref{33}. Therefore the unitary structure
of ${\rm Rep}(V_{{\mathfrak g}_k})$ will be clarified in a more complete way by the end of the paper.
\medskip

In the setting of rational conformal field theory, Moore and Seiberg observing an analogy between certain polynomial
equations that they had obtained and the   structure of a braided tensor category, discovered these structures in conformal field theory \cite{Moore-Seiberg1}, \cite{Moore-Seiberg2}. It took a long time and the work of many authors to construct braided tensor category structures on the WZW models, and show that 
they they satisfy all the axioms of modular tensor categories.
The reader may find some history in the introduction of \cite{On_a_problem_posed_by_Huang}.

 The first constructions were accomplished by the monumental work by
Kazhdan and Lusztig, Finkelberg.
In a series of papers, inspired by the work by Kazhdan and Lusztig, Huang and Lepowsky introduced a new   notion 
of  {\it braided vertex tensor category} and constructed this structure on ${\rm Rep}(V)$ under certain general assumptions, one of the
most important being   {\it associativity of intertwining operators}. This notion is stronger than the notion of braided tensor category, and their
  methods for the tensor product module are quite different from the  work by Kazhdan and Lusztig and verification of the coherence properties  for the associativity morphisms and the braiding (pentagon and hexagon equations) is mostly included in their construction.

Let $V$ be a VOA satisfying the rationality assumptions in \cite{Huang2}, namely: 
\medskip
\begin{itemize}
\item[(a)] $V$ is simple and of CFT type (i.e. $V_{(n)}=0$ for $n<0$, $V_{(0)}={\mathbb C}1$) and
the contragredient module $V'$ is isomorphic to $V$ as a $V$-module; 
\item[(b)] every $\mathbb{Z}_{\geq 0}$-graded weak module is a direct sum of irreducible $V$-modules;
\item[(c)] $V$ is $C_2$-cofinite, that is $V/C_2(V)$ is finite dimensional, with $C_2(V)$ the subspace of $V$ spanned by $a_{(-2)}b$, for $a$, $b\in V$.
\end{itemize}
\medskip

 These rationality assumptions are  motivated by the affine vertex operator algebras at positive integer levels, and Virasoro vertex operator algebra.
We   refer the reader  to  \cite{HL} for  a review on the construction of vertex braided tensor category structure of ${\rm Rep}(V)$ by Huang and Lepowsky, and references to the original articles, to the introduction of \cite{Huang_differential_equations}  and
  Remark 3.8, and Theorem 3.9  therein for an explanation on the implication from the rationality conditions (a), (b), (c) to  the verification of the conditions needed for the construction of Huang and Lepowsky vertex tensor category structure and to   \cite{Huang2}
  for more complete historical information and references, where the properties of  rigidity and modularity under these rationality conditions were proved by Huang. 
  
  It is important to note  that this general braided tensor category construction
 was preceded by a direct application of Huang and Lepowsky tensor product theory to affine vertex operator algebras  at positive integer level (WZW model)  
 in \cite{Huang_Lepowski_affine}, where the authors proved the conditions to apply their theory, including
  associativity of intertwining operators, building on previous work of several authors, including
  Knizhnik-Zamolodchikov \cite{KZ}, Tsuchiya-Kanie \cite{Tsuchiya_Kanie}, Frenkel-Huang-Lepowsky \cite{FHL},
   Frenkel-Zhu \cite{Frenkel_Zhu}, Dong, Li and Mason \cite{DLM1}.
In their work, validity of Knizhnik-Zamolodchikov differential equations for products of intertwining operators
   is established at the level of formal variables, and it is used to show the needed convergence of products   of intertwining operators as functions of complex variables, by Theorem 3.2 in \cite{Huang_Lepowski_affine}.

\begin{thm} 
\label{thmZhuTensor}
Let $V$ be a VOA satisying the rationality assumptions (a), (b), (c). Assume moreover that 
$M\mapsto D(M):= {\rm        dim}(\mathcal{F}_V(M))$, $M$ irreducible, gives a weak dimension function on the modular tensor category ${\rm        Rep}(V)$. Then, the Zhu's algebra $A(V)$ admits the structure of a weak quasi-Hopf algebra with a tensor equivalence 
$\mathcal{E}_V: {\rm        Rep}(V) \to {\rm        Rep}(A(V))$ which, after composition with the forgetful functor 
$:{\rm        Rep}(A(V)) \to {\rm        Vec}$ is tensor isomorphic to  $\mathcal{F}_V$. 
\end{thm}
\begin{proof} By Theorem \ref{propweakdim} $\mathcal{F}_V$ admits a weak quasi-tensor structure and the conclusion follows from Theorem \ref{TK_algebraic_quasi}.
\end{proof}

\begin{rem} The functor $\mathcal{E}_V: {\rm        Rep}(V) \to {\rm        Rep}(A(V))$ already appeared in the literature without mention to the tensor structure, see \cite{DongLiMason,HuangYang,Zhu}. 

\end{rem} 

\begin{rem} 
\label{RemarkWeakDimension}
The condition on $M \mapsto D(M)$, which we will call the weak dimension condition, is not satisfied in general. For example if $V$ is a rational unitary Virasoro VOA then $D(M)=1$ for all irreducible V-modules $M$. Moreover,  from the known fusion rules of these models, see e.g. \cite[Sec. 2.2]{KL}, it follows that one can always find an irreducible $M$ with $D(M\otimes M) = 2 > D(M)^2$ 
and hence the weak dimension condition is not satisfied. On the other hand the class of rational VOAs satisfying the weak dimension condition include many remarkable examples such as the unitary simple affine VOAs and the lattice VOAs.   
\end{rem}

We now discuss the case of unitary affine VOAs. Let $\mathfrak{g}$ be a complex simple Lie algebra and let $k$ be a positive integer. Moreover let $\mathfrak{g}_{\mathbb{R}}\subset \mathfrak{g}$ be a real form of ${\mathfrak{g}}$ and let $G$ be the corresponding simply connected compact simple Lie group.  We denote by $V_{\mathfrak{g}_k}$ the level $k$ affine simple unitary VOA associated to the pair $(\mathfrak{g},k)$. It is known to satisfy the assumptions (a), (b), (c) so that 
${\rm        Rep}(V_{\mathfrak{g}_k})$ is a modular tensor category with Huang-Lepowsky ribbon braided tensor category structure.  
Accordingly we can consider the functor $\mathcal{F}_{V_{\mathfrak{g}_k}}$ which satisfies the weak-dimension condition so that 
the Zhu's algebra $A(V_{\mathfrak{g}_k})$ admits a weak quasi-Hopf algebra structure. 

Now, let us consider the quantum group 
$U_q(\mathfrak{g})$ with $$q=e^{\frac{i\pi}{d(k+h^\vee)}}, \quad \ell= d(k+h^\vee),$$ see Sect. \ref{73} for a brief review, where $h^\vee$ is the dual Coxeter number of $\mathfrak{g}$.  Thus $q$ is a minimal root by Def. \ref{minimal_root}.  Consider  
 the fusion category ${\mathcal C}({\mathfrak g}, q, \ell)$ obtained from the category of tilting modules recalled in Sect. \ref{74}. It is a modular braided category admitting a compatible $C^*$-structure by \cite{Wenzl,Xu}. 
Let $\mathcal{F}_{(\mathfrak{g},q)}: {\mathcal C}({\mathfrak g}, q, \ell) \to {\rm        Vec}$ be   Wenzl functor. Then 
$\mathcal{F}_{(\mathfrak{g},q)}$ satisfies the weak-dimension condition and hence it defines a weak quasi-Hopf algebra 
$A({\mathfrak g}, q, \ell)$ following the general procedure applied  in Subsect. \ref{19.10}.

We next have the following remark that plays an important role in this paper.

\begin{rem} By Theorem \ref{Finkelberg_HL}, the category ${\mathcal C}({\mathfrak g}, q, \ell)$ is tensor equivalent to ${\rm        Rep}(V_{\mathfrak{g}_k})$ for the classical Lie types and $G_2$. The weak dimension functions for the functors $\mathcal{F}_{V_{\mathfrak{g}_k}}$
and $\mathcal{F}_{(\mathfrak{g},q)}$ regarded as defined on the same family of irreducible objects, have the same range in $\mathbb{Z}_{\geq 0}$. it follows that $A(V_{\mathfrak{g}_k})$ and $A({\mathfrak g}, q, \ell)$ are, up to a twist, isomorphic weak quasi-Hopf algebras, cf. the discussion after Theorem \ref{propweakdim}. 
 (By the main result of
\cite{CP},  $A({\mathfrak sl}_N, q, \ell)$ admits a natural structure of  weak  Hopf algebra. We shall extend this result to the other Lie types later on, and it will be useful)
Here we have followed   the original ideas of \cite{MS}, \cite{HO}.

Note that  the quantum group fusion category ${\mathcal C}({\mathfrak g}, q, \ell)$ is unitary and admits a weak tensor functor to 
${\rm Hilb}$. Therefore $A({\mathfrak g}, q, \ell)$ admits the structure of an $\Omega$-involutive ribbon weak quasi-Hopf algebra.
Next we are going to use this fact to transport this unitary structure to $A(V_{\mathfrak{g}_k})$ following our general unitarizability results.

  \end{rem}

We now discuss the unitary aspects of the above constructions. To this aim, we shall not need the detailed construction of weak  Hopf algebras alluded to in the previous remark, and general weak quasi-Hopf algebra constructions will suffice.
We first need to recall some properties of the Zhu's algebra and fix some notation. 
From now on we shall mostly pass to the following convenient notation, which is well known and very convenient.\bigskip

\begin{defn} \label{new_notation_for_gradinf_for_modules} {\it (New grading notation for modes of module maps)} 
Let  $a \in V$ be a homogeneous element of conformal weight (or degree) $d \in \mathbb{Z}$, i.e. such that $L_0 a=da$.
 For every $V$-module $M$ then $a^M_n$ is defined by 
 \begin{equation}\label{new_grading} a^M_n := a^M_{(n -1)}, \quad\quad n \in \mathbb{Z}.\end{equation} 
 For a general $a \in V$ $a^M_n$ is defined by linearity. 
 With respect to the new notation (\ref{new_grading}) and  to the new
 gradation for the homogeneous subspaces of a module defined in (\ref{grades}),
 the grading relation  (\ref{grading_equation_for_modules})
 becomes
 \begin{equation}\label{grading_equation_for_modules_new}
  a^M_{n}: M_{(r)}\to M_{(r+d-n)} \quad\text{for}\quad  r\in{\mathbb R}, \quad n\in{\mathbb Z}. 
  \end{equation}
Thus after a first  increase  by $d={\rm deg}(a)$, energy further   decreases by $-n$ for $n>0$ and increases by $-n$ for $n<0$.
 \end{defn}
 
As a vector space the Zhu's algebra is a quotient $V/O(V)$ for a certain subspace $O(V) \subset V$ and we denote by $a \mapsto [a]$ 
the quotient map $:V \to A(V)$. When $V$ satisfies the assumption (a), (b), (c)   then 
\begin{equation}
O(V) = \{a\in V: a^M_d\restriction_{M_{(0)}} = 0\; \text{for all $V$-modules $M$}\}\,,
\end{equation}  
where $d={\rm deg}(a)$.
Moreover, for every $V$-module $M$ the map $[a] \mapsto a^M_d\restriction_{M_{(0)}}$ is a representation of the associative algebra $A(V)$
on $M_{(0)}$ which is the one  corresponding to $\mathcal{E}_V(M)$ in Theorem \ref{thmZhuTensor}.

Let $V$ be a unitary VOA \cite{CKLW,DongLin} satisfying the rationality assumptions (a), (b), (c). Note that if $V$ is simple and unitary then a is necessarily of CFT type and isomorphic to the  contragredient module $V'$ as a V-module so that (a) is {\it a priori} satisfied. Let $\theta$ be the  PCT operator giving the unitary structure on $V$.  By  \cite[Eq. 5.3.1]{FHL} and \cite[Prop. 2.3.]{DongLiMason} the map 
$$[a] \mapsto [e^{L_1}(-1)^{L_0}a]$$ is an involutive anti-automorphism of $A(V)$. On the other hand, being $\theta$ an anti-linear involutive automorphism of $V$, we have that $\theta(O(V)) = O(V)$ and the map $[a] \mapsto [\theta a]$ is an anti-linear involutive automorphism of the associative algebra $A(V)$.  It follows that $$[a] \mapsto [a]^* := [e^{L_1}(-1)^{L_0}\theta a]$$
is an anti-linear involutive automorphism of $A(V)$ i.e. it gives a *-algebra structure on $A(V)$ canonically associated to the unitary structure of $A(V)$.  

\begin{prop}
\label{propUnitaryZhu} 
Let $M$ be a unitary $V$-module then the restriction to $M_{(0)}$ of the invariant scalar product of $M$ makes $\mathcal{E}_V(M)$ into a *-representation of $A(V)$. Moreover, the above restriction gives a one-to-one correspondence between the invariant scalar product on $M$ and 
the scalar products making $\mathcal{E}_V(M)$ into a *-representation of $A(V)$.
\end{prop}

\begin{proof} The first claim follows in a straightforward way from the definition of invariant scalar product and the 
*-operation on $A(V)$. Now, let $\mathcal{U}(M,V)$ be the associative algebra generated by the vertex operator coefficients $a^M_{(n)}$, $a \in V$, $n \in \mathbb{Z}$ as before.  $\mathcal{U}(M,V)$ carries a $\mathbb{Z}$-grading 
$$ \mathcal{U}(M,V) = \bigoplus_{n \in \mathbb{Z}}   \mathcal{U}(M,V)_n$$ 
where 
$$\mathcal{U}(M,V)_n := \{X \in \mathcal{U}(M,V) : e^{it L^M_0}Xe^{-it L^M_0} = e^{itn}X \} \,. $$ 
Accordingly, we have $a^M_n \in \mathcal{U}(M,V)_n$. Moreover, for every $X \in M$ there is an $X^*\in  \mathcal{U}(M,V)$ such that 
$(m_1,X m_2) = (X^*m_1,m_2)$ for all $m_1, m_2 \in M_{(0)}$, where $(\cdot,\cdot)$ is the invariant scalar product on $M$. Note that 
$(a^M_n)^* = (e^{L_1}(-1)^{L_0}\theta a)^M_{-n}$ for all $a \in V$ and all $n \in \mathbb{Z}$ so that 
$\left( \mathcal{U}(M,V)_n\right)^* = \mathcal{U}(M,V)_{-n}$ for all $n \in \mathbb{Z}$. In particular $\mathcal{U}(M,V)_0$ is a *-subalgebra
of $\mathcal{U}(M,V)$. For every $X \in  \mathcal{U}(M,V)_0$ we have $XM_{0} \subset M_{0}$ and hence $X$ restricts to an endomorphism $\tilde{X}$ of $M_{0}$. Now, given $m_1, m_2 \in M_{0}$  we have $(X_k m_1, Y_n m_2) = 0$ if $k \neq n$. Accordingly we have 
$$(Xm_1,Ym_2) = (m_1, \sum_{n\in \mathbb{Z}}(X_n)^*Y_n m_2) $$ 
which shows that the invariant scalar product on $M$ is determined by its restriction to $M_{(0)}$ Now, let $(\cdot,\cdot)$ be a fixed invariant scalar product on $M$ and let $\{\cdot,\cdot\}$ any scalar product on $M_{(0)}$ making $\mathcal{E}_V(M)$ into a *-representation of $A(V)$. Then there is an $A(V)$-module isomorphism $T_0: M_{(0)} \to M_{(0)}$ such that $\{m_1,m_2\}=(m_1,T_0 m_2)$ for all 
$m_1, m_2 \in M_{(0)}$. Since $\mathcal{E}_V$ is an equivalence of categories there is a unique $V$-module map $T: M\to M$ such that 
$\mathcal{E}_V(T)=T_0$ and we can define a sesquilinear form $\{\cdot,\cdot\}_M$ on $M$ by 
$\{m_1,m_2\}_M = (m_1,Tm_2)$, $m_1, m_2 \in M$. It is now straightforward to check that $\{\cdot,\cdot\}_M$ is an invariant scalar product on 
$M$ whose restriction to $M_{(0)}$ is $\{\cdot,\cdot\}$. 
\end{proof}

\begin{rem} Let $V$ a unitary vertex operator algebra satisfying the assumptions (a), (b), (c) so that ${\rm        Rep}(V)$ is a modular tensor category. Let ${\rm        Rep}^+(V)$ be the $C^*$-category of unitary representations of $V$. Then the forgetful functor 
$:{\rm        Rep}^+(V) \to {\rm        Rep}(V)$ is linear equivalence if and only if every $V$-module is unitarizable.  In this case ${\rm        Rep}^+(V)$ is equivalent as a $C^*$-category to the representation category ${\rm        Rep}^+(A(V))$ of finite dimensional *-representations of the $C^*$-algebra 
$A(V)$. It is not clear in general if the linear equivalence ${\rm        Rep}^+(V) \simeq {\rm        Rep}(V)$ can be used to make  ${\rm        Rep}^+(V)$ into a tensor $C^*$-category 
tensor equivalent to ${\rm        Rep}(V)$. This is an important problem which has been recently solved in some special cases by B. Gui \cite{GuiI,GuiII}. We also recall a work by Kirillov on the construction of a 
tensor $^*$-category closely related to ${\rm Rep}(V)$ which preceded the work by Huang and Lepowsky  \cite{Kirillov3}.
\end{rem}

\begin{prop}
\label{propUnitaryZhuB}
Let $V$ a unitary vertex operator algebra satisfying the assumptions (a), (b), (c). Then the equivalence 
$\mathcal{E}_V: {\rm        Rep}(V) \to {\rm        Rep}(A(V))$ gives in a canonical way a faithful *-functor 
$\mathcal{E}^+_V: {\rm        Rep}^+(V) \to {\rm        Rep}^+(A(V))$. 
If the forgetful functor ${\rm        Rep}^+(V) \to {\rm        Rep}(V)$ is an equivalence of linear categories then  $A(V)$ is a $C^*$-algebra and $\mathcal{E}^+_V: {\rm        Rep}^+(V) \to {\rm        Rep}^+(A(V))$ is an equivalence of $C^*$-categories. Moreover, in the latter case, any equivalence of linear categories $\mathcal{S}_V: {\rm        Rep}(A(V)) \to  {\rm        Rep}(V)$ together with an isomorphism 
$\eta: \mathcal{E}_V\circ \mathcal{S}_V \to 1_{{\rm        Rep}(A(V))}$ gives a canonical *-equivalence  
$\mathcal{S}^+_V: {\rm        Rep}^+(A(V)) \to  {\rm        Rep}^+(V)$ with 
$\mathcal{E}^+_V\circ \mathcal{S}^+_V$  unitarily equivalent to the identity. 
\end{prop}
\begin{proof} 
Let $M$ be a unitary $V$-module. Then $\mathcal{E}^+_V(M)$ is defined to be the $A(V)$-module $\mathcal{E}_V(M)$ together with the
scalar product obtained by restricting the given invariant scalar product on $M$. Then, thanks to Prop. \ref{propUnitaryZhu} 
$\mathcal{E}^+_V$ is a faithful *-functor which become an equivalence if the forgetful functor $:{\rm        Rep}^+(V) \to {\rm        Rep}(V)$ is a linear equivalence.. In the latter case we have the linear eqivalence ${\rm        Rep}(A(V)) \simeq  {\rm        Rep}^+(A(V))$ and hence $A(V)$ is a $C^*$-algebra. Assume now the linear equivalence ${\rm        Rep}^+(V) \simeq {\rm        Rep}(V)$
and let $\mathcal{S}_V: {\rm        Rep}(A(V)) \to  {\rm        Rep}(V)$ be an equivalence with a natural isomorphism 
$\eta: \mathcal{E}_V\circ \mathcal{S}_V \to 1_{{\rm        Rep}(A(V))}$. Let $W$ be a $C^*$-module for $A(V)$ and let $(\cdot,\cdot)_W$ be the corresponding scalar product. Then $(\eta_W\, \cdot \, \eta_W\, \cdot \, )_W$ is a scalar product on $\mathcal{E}_V\circ \mathcal{S}_V(W)$ making it into a *-representation of $A(V)$. Then it follows from Prop. \ref{propUnitaryZhu} and the assumption  ${\rm        Rep}^+(V) \simeq {\rm        Rep}(V)$    that there is a unique invariant scalar product on $\mathcal{S}_V(W)$ which restricts to 
$(\eta_W  \, \cdot \, \eta_W \, \cdot \, )_W$. This scalar product defines a unitary $V$-module $\mathcal{S}^+_V(W)$ and it is not hard to see that the map $W \mapsto \mathcal{S}^+_V(W)$ defines a functor with the desired properties. 
\end{proof}

\medskip

\begin{thm} 
\label{TheoremUnitaryRepVOA}
Let $V$ be a unitary vertex operator algebra satisfying assumptions (a), (b), (c) and such that  the forgetful functor 
$:{\rm        Rep}^+(V) \to {\rm        Rep}(V)$ is a linear equivalence and assume that the functor $\mathcal{F}_V: {\rm        Rep}(V) \to {\rm        Vec}$ satisfies the weak dimension condition in Remark \ref{RemarkWeakDimension}. Then  ${\rm        Rep}^+(V)$ admits a structure of tensor $C^*$-category with unitary braided symmetry such that the forgetful functor$:{\rm        Rep}^+(V) \to {\rm        Rep}(V)$ is a braided tensor equivalence if and only if the weak quasi-Hopf algebra  on $A(V)$ obtained from a weak quasi-tensor structure on the functor $\mathcal{F}_V: {\rm        Rep}(V) \to {\rm        Vec}$ admits  the structure of a $\Omega$-involutive weak quasi-Hopf $C^*$-algebra compatible with the canonical *-structure on $A(V)$. 
\end{thm}
\begin{proof} 
The functor $\mathcal{F}^+_V: {\rm        Rep}^+(V) \to {\rm        Hilb}$ obtained by composition of the equivalence $\mathcal{E}^+_V: {\rm        Rep}^+(V) \to {\rm        Rep}^+(A(V))$ with the forgetful functor $:{\rm        Rep}^+(A(V)) \to {\rm        Hilb}$ is a *-functor as a consequence of Prop. \ref{propUnitaryZhuB}.  If  ${\rm        Rep}^+(V)$ admits a structure of tensor $C^*$-category such that the forgetful functor 
$:{\rm        Rep}^+(V) \to {\rm        Rep}(V)$ is a tensor equivalence then $\mathcal{F}^+_V$ admits a weak  quasi-tensor *-structure so that 
${\rm        Nat}_0(\mathcal{F}^+_V)$ admits  the structure of a  weak quasi-Hopf $C^*$-algebra  as a consequence of Theorem \ref{TheoremTannakaStar}. By construction the $C^*$-algebra $A(V)$ with its canonical *-operation is isomorphic to  ${\rm        Nat}_0(\mathcal{F}^+_V)$  so that it inherits from the latter the structure of a $\Omega$-involutive weak quasi-Hopf $C^*$-algebra coinciding , up to a twist,  with the  weak quasi-Hopf algebra structure on $A(V)$ obtained from a weak quasi-tensor structure on the functor 
$\mathcal{F}_V: {\rm        Rep}(V) \to {\rm        Vec}$.

Conversely, if $A(V)$ admits  the structure of a $\Omega$-involutive weak quasi-Hopf $C^*$-algebra with the canonical *-structure then, by Corollary \ref{CorollaryRep^+(A)}  ${\rm        Rep}^+(A(V))$ is a tensor $C^*$-category tensor equivalent  ${\rm        Rep}(A(V))$ and hence to
${\rm        Rep}(V)$ . 
Let $\mathcal{S}_V: {\rm        Rep}(A(V)) \to  {\rm        Rep}(V)$  be any tensor equivalence together with an isomorphism of tensor functors
$\eta: \mathcal{E}_V\circ \mathcal{S}_V \to 1_{{\rm        Rep}(A(V))}$ and let    
$\mathcal{S}^+_V: {\rm        Rep}^+(A(V)) \to  {\rm        Rep}^+(V)$  be the corresponding canonical *-equivalence as in Prop. \ref{propUnitaryZhuB}
so that $\mathcal{E}^+_V\circ \mathcal{S}^+_V$  unitarily equivalent to the identity.

Given unitary $V$-modules $M^\alpha, M^\beta \in   {\rm        Rep}^+(V) $ we define a unitary module 
$M^\alpha \otimes M^\beta$ by      
$$ M^\alpha \otimes M^\beta :=  \mathcal{S}_V^+ \left( \mathcal{E}_V^+(M^\alpha) \otimes \mathcal{E}_V^+(M^\beta) \right) \,.$$
Moreover, if $\alpha$ denotes the unitarty associator in ${\rm        Rep}^+(A(V))$  we define 
the unitaries 
$$\alpha'_{M^\alpha,M^\beta,M^\gamma } : ( M^\alpha \otimes M^\beta) \otimes M^\gamma \to   
M^\alpha \otimes ( M^\beta \otimes M^\gamma ) $$
by 
\begin{eqnarray*}
&&\alpha'_{M^\alpha,M^\beta,M^\gamma } := \\ 
&&\mathcal{S_V^+}(1_{\mathcal{E}_V^+(M^\alpha)}\otimes\eta^{-1}_{\mathcal{E}_V^+(M^\beta)\otimes \mathcal{E}_V^+ (M^\gamma)}  \circ \alpha_{\mathcal{E}_V^+(M^\alpha), \mathcal{E}_V^+ (M^\beta), \mathcal{E}_V^+   (M^\gamma)}   \circ 
\eta_{\mathcal{E}_V^+(M^\alpha)\otimes \mathcal{E}_V^+(M^\beta)}\otimes 1_{\mathcal{E}_V^+(M^\gamma)}) \, 
\end{eqnarray*}
where $\eta: \mathcal{E}_V\circ \mathcal{S}_V \to 1_{{\rm        Rep}(A(V))}$ is the isomorphism used to define the functor $\mathcal{S}_V^+$.

Then, thanks to Prop. \ref{propUnitaryZhuB}, one can   check that this gives the desired $C^*$- tensor structure on ${\rm        Rep}^+(V)$. From the tensor equivalence ${\rm        Rep}^+(V) \simeq {\rm        Rep}(V)$ we see that ${\rm        Rep}^+(V)$ admits a braiding making the equivalence a braided tensor equivalence and this braided symmetry on ${\rm        Rep}^+(V)$ is necessarily unitary by \cite{Gal}. 
\end{proof}  

\begin{thm} 
\label{TheoremUnitaryBraidRepAffine}
Let $\mathfrak{g}$ be a complex simple Lie algebra of classical Lie type or $G_2$, and let $k$ be a positive integer and let  $V_{\mathfrak{g}_k}$ be the corresponding level $k$ affine unitary vertex operator algebra. Then ${\rm        Rep}^+(V_{\mathfrak{g}_k})$ admits the structure of tensor $C^*$-category with unitary 
braided symmetry such that the forgetful functor $:{\rm        Rep}^+(V_{\mathfrak{g}_k}) \to {\rm        Rep}(V_{\mathfrak{g}_k})$ is a braided tensor equivalence. 
\end{thm}
\begin{proof} It is known that every $V_{\mathfrak{g}_k}$-module is unitarizable and hence 
${\rm        Rep}^+(V_{\mathfrak{g}_k}) \simeq {\rm        Rep}(V_{\mathfrak{g}_k})$. Let $q=e^{\frac{i\pi}{d(k+h^\vee)}}$. Then the quantum group category ${\mathcal C}({\mathfrak g}, q, \ell)$ is a tensor $C^*$-category by \cite{Wenzl,Xu_star}. It follows from Theorem \ref{Finkelberg_HL} that 
${\mathcal C}({\mathfrak g}, q, \ell) \simeq {\rm        Rep}(V_{\mathfrak{g}_k})$ that $A(V_{\mathfrak{g}_k})$ admits the structure of a $\Omega$-involutive weak quasi-Hopf $C^*$-algebra and the conclusion follows from  Theorem \ref{TheoremUnitaryRepVOA}. 
\end{proof}

 \begin{rem}
 \label{RemarkGui}
 Theorem \ref{TheoremUnitaryBraidRepAffine} has been recently proved by B. Gui  by directly working on modules \cite{GuiI}, \cite{GuiII}.
 \end{rem}

Theorem \ref{TheoremUnitaryBraidRepAffine} can be seen as a special case of the following more abstract result which we will use to give other examples of unitary VOAs such that ${\rm Rep}^+(V)$ admits a tensor $C^*$-structure. 

\begin{thm}
\label{TheoremUnitaryBraidRepVOAbis} 
Let $V$ be a unitary vertex operator algebra satisfying assumptions (a), (b), (c) and such that  the forgetful functor 
$:{\rm        Rep}^+(V) \to {\rm        Rep}(V)$ is a linear equivalence and assume that the functor $\mathcal{F}_V: {\rm        Rep}(V) \to {\rm        Vec}$ satisfies the weak dimension condition in Remark \ref{RemarkWeakDimension}. Assume that ${\rm Rep}(V)$ is tensor equivalent to a 
tensor $C^*$-category. Then  ${\rm        Rep}^+(V)$ admits a structure of tensor $C^*$-category with unitary braided symmetry such that the forgetful functor$:{\rm        Rep}^+(V) \to {\rm        Rep}(V)$ is a braided tensor equivalence.
\end{thm} 
\begin{proof}

\end{proof}

We now give some examples of applications of Theorem \ref{TheoremUnitaryBraidRepVOAbis}.

\begin{ex}
\label{exLattices}
Let $L$ be an even positive definite lattice and let $V_L$ be the corresponding unitary VOA. It satisfies assumptions (a), (b), (c). 
It follows from \cite[Th. 4.12]{DongLin}  the forgetful functor ${\rm Rep}^+(V) \to {\rm Rep}(V)$ is a linear equivalence. The fusion ring of 
${\rm Rep}(V_L)$ is isomorphic to the finite abelian group $L^*/L$, where $L^*$ is the dual lattice of $L$.  For an irreducible $V_L$-module 
$M_{[x]}$, with equivalence class corresponding to $[x] \in L^*/L$ we have $D(M_{[x]}) = N_{[x]}$, where $N_{[x]}$ is the number of 
elements of $L^*$ in the equivalence class $[x]$ having minimal norm, see e.g. \cite{GG}. In some cases, e.g. square lattices, one can easily check that $D$ is a weak dimension function. The irreducible objects of ${\rm Rep}(V_L)$ are all invertible their equivalence classes form a finite abelian group $G \simeq L^*/L$. It follows that ${\rm Rep}(V_L)$ is tensor equivalent to ${\rm Vec}^\omega_G$ for some 3-cocycle 
$\omega \in Z^3(G, \mathbb{T})$, where ${\rm Vec}^\omega_G$ is the category of $G$-graded finite dimensional vector spaces with associators twisted by $\omega$, see \cite{EGNO}. ${\rm Vec}^\omega_G$ is tensor equivalent to the 
tensor $C^*$-category ${\rm Hilb}^\omega_G$ of $G$-graded finite dimensional Hilbert spaces with associators twisted by $\omega$
and hence, if $V_L$ satisfies the weak dimension property, ${\rm Rep}^+(V_L)$ admits a  structure of a tensor $C^*$-category with unitary braided symmetry making the forgetful functor $: {\rm Rep}^+(V_L) \to {\rm Rep}(V_L)$  into a braided tensor equivalence.  
\end{ex}

\begin{ex}
\label{exTwistedQuantum}
Let $V$ be a unitary VOA satisfying assumptions (a), (b), (c) and assume that $V$ is holomorphic i.e. that  ${\rm Rep}(V)$ is  equivalent to ${\rm Vec}$.  Let $G$ be a finite subgroup of the unitary automorphism group of $V$ and let $V^G$ be the corresponding orbifold unitary sub VOA.  It is conjectured that always $V^G$ satisfies (a), (b), (c) and that ${\rm Rep}(V^G)$ is braided tensor equivalent to  
${\rm Rep}(D^\omega(G)) \simeq \mathcal{Z}({\rm Vec}^\omega_G)$, for some 3-cocycle $\omega \in Z^3(G, \mathbb{T})$. Here, 
$D^\omega(G)$ the twisted quantum double quasi-Hopf algebra introduced in \cite{DPR} and $\mathcal{Z}({\rm Vec}^\omega_G)$ is the 
center of ${\rm Vec}^\omega_G$, \cite{EGNO}.  This conjecture is known to be true in various cases, see e.g.  \cite{CM,DongRenXu,Kirillov2,MuegerICMP}. Assume now that the above conjecture is true for a given $V$ and $G$ and also assume that every irreducible $V^G$-module is unitarizable. Since ${\rm Rep}(D^\omega(G))$  is tensor equivalent to a tensor $C^*$-category then, if 
$V^G$ satisfies the weak dimension property, ${\rm Rep}^+(V^G)$ admits a  structure of a tensor $C^*$-category with unitary braided symmetry making the forgetful functor $: {\rm Rep}^+(V^G) \to {\rm Rep}(V^G)$  into a braided tensor equivalence.  Let us know consider an explict example. Let $\Lambda$ be the Leech lattice, the even unimodular lattice of rank 24 with trivial root system, and let let $V_\Lambda$ be the corresponding lattice VOA. Since $\Lambda = \Lambda^*$, $V_\Lambda$ is holomorphic. $V_\Lambda$ as special automorphism of order two which can easily seen to be unitary, see \cite[Sec. 4.4]{DongLin} where this automorphim is denoted by $\theta$. As usual we denote by 
$V_\Lambda^+$ the corresponding unitary fixed point subalgebra. $V_\Lambda^+$ satisfies (a), (b) and (c). Moreover, up to equivalence it has 
exactly four irreducible modules $V_\Lambda^+$, $V_\Lambda^-$, $(V_\Lambda^T)^+$ and $(V_\Lambda^T)^-$ which are all invertible and unitarizable \cite{DongLin,DongRenXu}.  Hence the equivalence classes of irreducibles form an abelian group of order 4 which 
in fact is isomorphic to $\mathbb{Z}_2\times \mathbb{Z}_2$, see e.g. \cite[Prop.5.6]{DongRenXu}.
Arguing as before can conclude that ${\rm Rep}(V_\Lambda^+)$ is tensor equivalent to a tensor $C^*$-category. The characters (graded dimensions) of the irreducible modules of $V_\Lambda^+$ are known, see \cite[Sec. 10.5]{FLM} and 
\cite[Prop. 2.5]{Shimakura} and from them one can easily compute the function $M \mapsto D(M) = \mathrm{dim}\mathcal{F}_{V_\Lambda^+}(M)$ and we find $D(V_\Lambda^+) =1$, $D(V_\Lambda^-)=24$, $D\left((V_\Lambda^T)^+\right) = 2^{12}$ and $D\left( (V_\Lambda^T)^-\right) = 24\cdot 2^{12}$. It follows that $V_\Lambda^+$ has the weak dimension property and hence
 ,  by Theorem \ref{TheoremUnitaryBraidRepVOAbis}, ${\rm        Rep}^+(V_\Lambda^+)$ admits a structure of tensor $C^*$-category with unitary braided symmetry such that the forgetful functor$:{\rm        Rep}^+(V_\Lambda^+) \to {\rm        Rep}(V_\Lambda^+)$ is a braided tensor equivalence. With this structure ${\rm        Rep}^+(V)$ is a modular tensor $C^*$-category because ${\rm        Rep}(V_\Lambda^+)$ is modular. 
The modular $T$ matrix of ${\rm        Rep}^+(V)$ can also be computed from the characters and it is given 
by the diagonal matrix with diagonal entries $1, 1, 1, -1$. By \cite{RSW} there is, up to equivalence, a unique unitary fusion category with fusion rules  $\mathbb{Z}_2\times \mathbb{Z}_2$, the above $T$ matrix and topological central charge 24 mod 8 and  it is realized by the representation category of the quantum double $D(\mathbb{Z}_2)$, with trivial twist $\omega \in H^3(\mathbb{Z}_2, \mathbb{T}) \simeq \mathbb{Z}_2$. Note that,  $A(V_\Lambda^+)$ and $D(\mathbb{Z}_2)$ have equivalent representation categories but are inequivalent associative algebras. $D(\mathbb{Z}_2)$ is commutative while  $A(V_\Lambda^+)$ is not. Note also that  $D(\mathbb{Z}_2)$ is a Hopf algebra while $A(V_\Lambda^+)$ is a weak quasi-Hopf algebra. 
\end{ex} 

\medskip 

\section{Conformal nets and weak quasi-Hopf algebras}\label{VOAnets3}

In this brief section   we explain how  most of the constructions and results we have discussed in the case of rational vertex operator algebras have an analogue in the case of completely rational conformal nets. These two picture are perhaps related by the correspondence between unitary vertex operator algebras, conformal nets and their representations \cite{CKLW,CWX,GuiI,GuiII}. 

Let $\mathcal{A}$ be a completely rational conformal net on $S^1$. We denote by ${\rm        Rep}(\mathcal{A})$ the category of (Hilbert space *-) representations of $\A$ with finite index. Note that every irreducible locally normal representation of $\A$ has finite index and hence is an object in  ${\rm        Rep}(\mathcal{A})$. Accordingly the finite index condition is assumed only to rule out infinite Hilbertian direct sums. Its known that ${\rm        Rep}(\mathcal{A})$ is a modular tensor $C^*$-category \cite{Kawahigashi,KLM}. Here we briefly describe how this structure of modular tensor category is defined. Let $I \subset S^1$ be a given non-empty non-dense open interval. Then one can define a full $C^*$-subcatefgory 
${\rm        Rep}_I(\A)$, ${\rm        Rep}(\mathcal{A})$ whose objects are the representations localized in $I$, see e.g. \cite[Sec. 3.2.]{Kawahigashi}. 
The objects in ${\rm        Rep}_I(\A)$ gives rise to unital endomorphisms of the type III factor $\A(I)$ and the composition of endomorphisms makes 
${\rm        Rep}_I(\A)$ into a strict tensor $C^*$-category which turns out to be modular as a cosequence of the results in \cite{KLM}. It is known that 
Every representation in ${\rm        Rep}(\mathcal{A})$ is unitary equivalent to a representation in ${\rm        Rep}_I(\A)$ so that the embedding 
$\mathcal{I}: {\rm        Rep}_I(\A) \to {\rm        Rep}(\A)$ is a unitary equivalence of $C^*$-categories.   Accordingly, given any equivalence 
$\mathcal{E}: {\rm        Rep}(\A)  \to {\rm        Rep}_I(\A)$ with a unitary isomorphism $\eta: \mathcal{E}\circ\mathcal{I} \to 1_{{\rm        Rep}_I(\A)}$  one can transport the modular tensor $C^*$-category structure on ${\rm        Rep}(\A)$ and give to $\mathcal{E}$ a tensor structure making it into a unitary 
tensor equivalence. Note that one can chose $\mathcal{E}$ such that $\mathcal{E}\circ\mathcal{I} = 1_{{\rm        Rep}_I(\A)}$  and accordingly 
$\eta$ such that $\eta_\pi = 1_\pi$ for all $\pi$ in ${\rm        Rep}_I(\A)$. With this choice ${\rm        Rep}(\A)$ turns out to be a strict tensor $C^*$-category.

Given a representation $\pi$ of $\A$ with finite index 
we denote by $L^\pi_0$ the corresponding conformal Hamiltonian.  $L^\pi_0$ is a self-adjoint operator with pure point-spectrum.  In the following we will assume that $\A$ satisfies the following 
\medskip
\begin{itemize}
\item[(d)] For every representation $\pi$ of $\A$ with finite index $L^\pi_0$ has finite dimensional eigenspaces. 
\end{itemize}
\medskip

Assumption (d) is believed to be always satisfied. It would follow e.g. from \cite[Conjecture 9.4]{CKLW} or from \cite[Conjecture 4.18]{Kawahigashi}.

We now want to define a conformal net analogue of the functor $\mathcal{F}_V$ defined at the beginning of this section. 
 Every representation $\pi$ of $\A$ with finite index on the Hilbert space $\mathcal{H}^\pi$  can be written as a finite direct 
sum of irreducibles 
\begin{equation}
\label{repirrdecomposition}
\pi = \bigoplus_{i} \pi^i  
\end{equation}   
and correspondingly a Hilbert space decomposition. 
\begin{equation}
\label{Hilbertirrdecomposition}
\mathcal{H}^\pi = \bigoplus_{i} \mathcal{H}^{\pi_i}  \, .
\end{equation}
We denote by $h_i \geq 0$ the lowest eigenvalue of  $L^{\pi_i}_0$ and by $\mathcal{H}^{\pi_i}_{(0)}$ the corresponding eigenspace which is finite dimensional by our previous assumption. 
We now define a finite dimensional closed subspace $\mathcal{H}_{(0)} \subset \mathcal{H}$ by 
\begin{equation}
\mathcal{H}^\pi_{(0)} := \bigoplus_i \mathcal{H}^{\pi_i}_{(0)} \,. 
\end{equation}
$\mathcal{H}_{(0)}$ is independent from the choice of the direct sum decomposition in Eq. (\ref{repirrdecomposition}). 
Moreover, $\pi(\A)'' \mathcal{H}^\pi_{(0)} = \mathcal{H}^\pi$ where $\pi(\A)''$ is the von Neumann algebra on $\mathcal{H}^\pi$ generated by the algebras $\pi_I(\A(I))$, with $I$ an open non-dense non-empty interval of $S^1$. 

In complete analogy with the VOA case one can define a linear functor $\mathcal{F}_\A: {\rm        Rep}(\A) \to {\rm        Hilb}$ 
by $\mathcal{F}_\A(\pi):= \mathcal{H}^\pi_{(0)}$ for any representation with finite index $\pi$ of $\A$ and 
 $\mathcal{F}_\A(T):= T\restriction \mathcal{H}^{\pi^\alpha}_{(0)}$ for any intertwiner operator $T\in (\pi^\alpha, \pi^\beta)$ and it turns out that
 $\mathcal{F}_\A$ is a faithful *-functor.  The algebra $A(\A):={\rm        Nat}_0(\mathcal{F}_\A)$ is a finite dimensional $C^*$-algebra and there is a *-equivalence of $C^*$-categories $\mathcal{E}_\A: {\rm        Rep}(\A) \to {\rm        Rep}^+(A(V))$ which, after composition with the forgetful functor  $:{\rm        Rep}(\A) \to {\rm        Hilb}$ is isomorphic to  $\mathcal{F}_\A$. The algebra $A(\A)$ is the conformal net analogue of the Zhu's algebra.
The following is the conformal net version of Theorem \ref{thmZhuTensor}

\begin{thm} 
\label{thmNetZhuTensor}
Let $\A$ be a completely rational conformal net satisfying assumption (d). Assume moreover that 
$\pi \mapsto D(\pi):= {\rm        dim}(\mathcal{F}_\A(\pi))$, $\pi$ irreducible, gives a weak dimension function on the modular tensor category 
${\rm        Rep}(\A)$. Then, the algebra $A(\A)$ admits a structure of a $\Omega$-involutive weak quasi-Hopf $C^*$-algebra with a *-tensor equivalence 
$\mathcal{E}_\A\: {\rm        Rep}(\A) \to {\rm        Rep}^+(A(\A))$ which, after composition with the forgetful functor 
$:{\rm        Rep}(\A) \to {\rm        Hilb}$ is tensor isomorphic to  $\mathcal{F}_\A$. 
\end{thm}
\begin{proof} By Theorem \ref{propweakdim} $\mathcal{F}_\A$ admits a weak quasi-tensor structure and the conclusion follows from Theorem 
\ref{TheoremTannakaStar}.
\end{proof}

\medskip

We conclude this section with a brief comparison of the VOA and the conformal net quasi-Hopf algebras discussed in this section. In \cite{CKLW} a class of unitary simple VOAs called {\it strongly local} VOAs has been introduced and a map $V \to \A_V$ form strongly local VOAs to conformal nets has been defined. It is conjectured in \cite{CKLW} that every simple unitary vertex operator algebra $V$ i strongly local 
and that the map $V\to \A_V$ gives a one-to-one correspondence between unitary simple VOAs and (irreducible) conformal nets. Moreover, it is conjectured in \cite[Conjecture 4.43]{Kawahigashi} that the unitary VOA satisfies assumptions (a), (b) and (c) if and only if $\A_V$ is completely rational and that, in this case ${\rm        Rep}(\A_V)$ and ${\rm        Rep}(V)$ are tensor equivalent, see also \cite{GuiI,GuiII,Huang3}.  
This conjecture appears to be a very hard and important problem and whose solution for even for a representative class of  examples is of great interest.  We hope that our work could give some useful hints in this directions and we hope to come back to this in future work. Here we limit 
ourselves to give some hints in the special case of the type $A$ affine vertex operator algebras $V_{{\mathfrak{sl}_N}_k}$. 

We now from \cite{CKLW} that, for all $N\geq2$ and all $k\geq1$,  $V_{{\mathfrak{sl}_N}_k}$ is a simple unitary strongly local 
VOA and that the conformal net $\A_{V_{{\mathfrak{sl}_N}_k}}$ is isomorphic to the loop group conformal net $\A_{{\rm        SU}(N)_k}$.
The latter is known to be completely rational as a consequence of Wassermann's work \cite{Wassermann} and the fusion rules 
of ${\rm        Rep}(\A_{{\rm        SU}(N)_k})$ are known to agree with those of  ${\rm        Rep}(V_{{\mathfrak{sl}_N}_k})$. Actually the two modular tensor categories are known to have the same modular data, i.e. the same modular $S$ and $T$ matrices. Moreover by \cite{CWX}, see also 
\cite{GuiII},  every unitary  $V_{{\mathfrak{sl}_N}_k}$-module $M$ ``integrates'' to a representation $\pi^M$ of  $\A_{{\rm        SU}(N)_k}$ on the Hilbert space completion $\mathcal{H}_M$ of $M$ and the map $M \mapsto \pi^M$ gives rise to a *-isomorphism of $C^*$-categories 
$$\mathcal{E}_{{\rm        SU}(N)_k}:{\rm        Rep}^+(V_{{\mathfrak{sl}_N}_k})  \to {\rm        Rep}(\A_{{\rm        SU}(N)_k}) $$
and it is straightforward to see that 
$$\mathcal{F}^+_{V_{{\mathfrak{sl}_N}_k}} =\mathcal{F}_{\A_{{\rm        SU}(N)_k}} \circ \mathcal{E}_{{\rm        SU}(N)_k}\, . $$
As a consequence we have a canonical isomorphism $A(V_{{\mathfrak{sl}_N}_k}) \simeq A(\A_{{\rm        SU}(N)_k})$ and we have a tensor equivalence ${\rm        Rep}^+(V_{{\mathfrak{sl}_N}_k})  \simeq {\rm        Rep}(\A_{{\rm        SU}(N)_k})$ if and only if the weak quasi-Hopf algebra structures on $A(V_{{\mathfrak{sl}_N}_k})$ and $A(\A_{{\rm        SU}(N)_k})$ agree up to a twist. \bigskip

 \bigskip  
 
\section{Kazhdan-Wenzl theory and equivalence of ribbon ${\mathfrak sl}_{N, q, \ell}$-categories}\label{KW}
This section can be read independently of the rest of the paper. The main result is the proof of Theorem \ref{classification_type_A}, that highlights the relevance of our weak Hopf algebras
in identifying associators of two categories in presence of a braided symmetry using classification methods by Kazhdan-Wenzl in the type $A$ case. This principle will be used for the other Lie types  in this paper.

   Let ${\mathfrak g}$ be a simple complex Lie algebra.
  We keep the notation fixed in the first paragraph of the previous section for  ${\mathcal C}({\mathfrak g}, q, \ell)$. 
    We recall that the fusion categories  
${\mathcal C}({\mathfrak g}, q, \ell)$ arising from quantum groups at roots of unity are deeply related to fusion categories arising from chiral CFT on the circle. Let
  $k$ be  a positive integer and let $V_{{\mathfrak g}_k}$ denote the affine Vertex Operator Algebra (VOA) of level   $k$ with ${\rm Rep}(V_{g_k})$ the associated
representation category.     By results of Huang \cite{Huang(modularity), Huang1, Huang2, Huang3} this is a  modular fusion category. 
The work by  Kazhdan and Lusztig and Finkelberg \cite{Finkelberg}, \cite{Finkelberg_erratum}
culminated in the    construction of a second modular tensor category $\tilde{\mathcal O}_\ell$ associated
to modules of affine Lie  algebras at positive integer levels for all the Lie types except  $E_6$, $E_7$ $k=1$ and $E_8$, $k=1$, $2$.  The combination of these works prove that 
$\tilde{\mathcal O}_k$ and  ${\mathcal C}({\mathfrak g}, q, \ell)$
  are  equivalent as ribbon categories  for the specific roots of unity $q=e^{i\pi/d\ell}$ with   $\ell=k+\check{h}$.

On the other hand, the approach to CFT via conformal nets \cite{GF} provides   examples of modular fusion categories as well \cite{KLM}. A general connection  from VOA satisfying suitable analytic conditions   to conformal nets has recently been established \cite{CKLW}. 
 
An important example is the   fusion category associated to the loop group conformal net   over ${\rm SU}(N)$ which is known to have    the same fusion rules  \cite{Wassermann} and modular data (the $S$ and $T$ matrices) as   the corresponding affine VOA or quantum group categories. 
More precisely, the associated Verlinde fusion ring  $R_{N, \ell}$  arises from positive energy representations of the level $k$ central extension of the loop group  of ${\rm SU}(N)$ and also as the Grothendieck ring of
     ${\rm Rep}(V_{g_k})$ or
    ${\mathcal C}({\mathfrak sl}_N, q, \ell)$  for any $q$ such that $q^2$ is a primitive root of unity of order $\ell$,      in this case $\ell=k+N$ see e.g. \cite{BK, Sawin, Andersen1}.
       
It is then natural to ask whether there
is a   classification of   ribbon  fusion categories  with Verlinde fusion rules of type $A$ showing in particular ribbon equivalence 
of the fusion categories arising from the three different settings. In this   section  we give a classification
result   independent of Finkelberg equivalence theorem.
 We shall   not assume that   our categories have a unitary structure, and we replace this condition with the possibly weaker assumption 
of pseudounitarity in the sense of \cite{ENO}.  In this way   our result may be useful for   the purposes of Sect. 
 \ref{VOAnets}, \ref{VOAnets2} for this special case. In that section   
 we construct unitary structures of the representation category of all the affine vertex operator algebras.

  a)  Let $R_N$   denote the  representation ring  
  of ${\rm SL}(N, {\mathbb C})$. It is freely generated   with 
     basis $e_\lambda$ parameterised by the set of dominant integral weights $\Lambda$, so every $\lambda\in\Lambda$ is a  non increasing sequence $(m_1,\dots, m_{N-1})$ of non negative integers.             
     
    b)   For a positive integer $\ell>N$, let $\Lambda^+(q)$ be the Weyl alcove  recalled at the beginning of the previous section. For ${\mathfrak g}={\mathfrak sl}_N$,  $\Lambda^+(q)$ may be described  by   weights $\lambda\in\Lambda$ satisfying  $m_1\leq\ell-N$. The Verlinde fusion ring $R_{N,\ell}$ has a natural basis $e_\lambda$ with $\lambda\in\Lambda^+(q)$. The structure constants are determined by the Verlinde formula (see \cite{Huang1} in the setting of vertex operator algebras, and references therein), or via characters 
    of the affine Weyl group, the Kac-Walton formula Exercise 13.35 in \cite{Kac2},  \cite{Walton90}.
     The fusion ring    $R_{N,\ell}$   may also described as a quotient of $R_N$ by a certain ideal, see \cite{Douglas}.

 We set $R_{N,\infty}=R_N$,   so the general notation
  $R_{N, \ell}$ will   include $N+1\leq \ell\leq\infty$ unless otherwise stated (as it will be  for example in the main theorem of the section).    Furthermore, $R_{N, \ell}$   will be regarded as a based ring   in the sense, e.g.,
 of \cite{Ostrik}.

 Note that  a semisimple rigid tensor category ${\mathcal C}$ with based Grothendieck ring isomorphic to $R_{N, \ell}$ for $\ell$ finite
  is a fusion category.
  
Frobenius-Perron dimensions   of basis elements  ${\rm FPdim}(X_i)$ of a commutative based ring were introduced in \cite{FK}, and one has
 ${\rm FPdim}(X_i)>0$.
We refer to Sect. 8 in \cite{ENO} or Chapter 4 in \cite{EGNO} for the development of the theory in   generality. We shall
be interested in the case of the based Grothendieck ring ${\rm Gr}({\mathcal C})$ of a fusion category ${\mathcal C}$     endowed with its natural basis given by the equivalence classes of irreducible objects.

The main result   is that  $X_i\to {\rm FPdim}(X_i)$ extends uniquely to a homomorphism of algebras
 $\phi: {\rm Gr}({\mathcal C})\to{\mathbb R}$, and $\phi$ is the unique homomorphism   such that $\phi(X_i)>0$ for all $i$, see   Theorem 8.2  and Lemma 8.3 in \cite{ENO}.
 The global Frobenius-Perron dimension  is defined as ${\rm FPdim}({\mathcal C})=\sum_i{\rm FPdim}(X_i)^2$.

The global categorical dimension is in turn defined as the sum of the {\it squared dimensions} $|X_i|^2$ of simple objects $X_i$.  
Squared and global categorical dimensions  were introduced and studied
 by M\"uger  for spherical fusion categories in  \cite{Mueger0}, \cite{Mueger1} and extended to general fusion categories in \cite{ENO}. It 
 is known that    $|X_i|^2>0$ and,   if ${\mathcal C}$ is  spherical,
$|X_i|^2=d(X_i)^2$, with $d$ the categorical dimension  defined via the spherical  structure, see     
  Sec. \ref{17}. 
 In particular, $d(X_i)^2$ is independent of the choice of the spherical structure. 
 
A fusion category ${\mathcal C}$  is called {\it pseudo-unitary} if the global  dimension ${\rm dim}({\mathcal C})$  equals the Frobenius-Perron dimension ${\rm FPdim}({\mathcal C})$.

The squared dimension of every simple object $X_i$ is bounded above by ${\rm FPdim}(X_i)^2$,
 hence ${\mathcal C}$ is pseudo-unitary if and only if these are all equalities, see Prop. 8.21 in \cite{ENO}.
 By Prop. 8.23 of the same paper a pseudo-unitary fusion category admits a unique pivotal structure, in fact spherical, 
 such that the categorical dimensions of simple objects $X_i$ are positive, or equivalently coincide with the 
 ${\rm FPdim}(X_i)$.

 We next specialise to  braided fusion categories. In this case, pivotal (spherical) structures are in a natural bijective correspondence with balanced (ribbon) structures for the braided symmetry, and the correspondence is recalled in Sect. \ref{17}.
  It also follows from the previous paragraph that  a pseudo-unitary braided fusion category admits a unique  ribbon structure         inducing
  positive categorical dimensions. We shall refer to it as the {\it positive } ribbon structure. 
 The aim of this section is to show Theorem \ref{classification_type_A} stated in the introduction.

   For $N=2$ Theorem  \ref{classification_type_A} has recently been shown in  \cite{Bischoff} using
 Fr\"ohlich-Kerler classification \cite{FK}.
  It follows from Ex. \ref{non_uniqueness}    that   the positivity requirement on the ribbon structures can not be removed.
  Moreover it will be clear from the proof how a   ribbon structure can be positive only for a unique braiding.
   We reformulate Theorem  \ref{classification_type_A}  in  a  form useful  for applications.

   \begin{thm}\label{corollary_of_KW}
   Let ${\mathcal C}$ and ${\mathcal C}'$ be modular fusion categories  with   positive categorical dimensions  and with Grothendieck rings isomorphic to the Verlinde fusion ring $R_{N, \ell}$ 
   via an isomorphism compatible with the corresponding  $T$-matrices.
Then ${\mathcal C}$ and ${\mathcal C}'$ are equivalent as ribbon tensor categories.
   \end{thm}
   
   \begin{proof}
 The categories are pseudo-unitary by positivity of the categorical dimensions. Compatibility of the $T$-matrices implies compatibility of the ribbon structures. The conclusion follows from Theorem \ref{classification_type_A}.  \end{proof}
 
   \begin{rem}\label{closely_related} ({\it closely related results in the literature}) We recall  that a  characterization of braided ${\mathfrak sl}_{N, \infty}$-type categories has been made explicit Neshveyev and Yamashita in \cite{NY_twisting}, see also the thesis by Jordans \cite{Jordans}, and a classification of the braided symmetries may be found in \cite{PR_rigidity}. In \cite{NY_twisting} the authors settle down the problem of reconstructing the   twisted categories obtained from Kazhdan-Wenzl classification, see Theorem \ref{Kazhdan_Wenzl},
    as representation categories of quantum groups of their own, for $q$ a positive real number.
   We note that a  result closely related to our following Prop. \ref{braiding_constraint}
has also been obtained in \cite{Makoto_et_al}  in the $C^*$-case with different methods.
We also note the   recent paper \cite{Rowell_et_al} on closely related topics.
Moreover, an analogue of some of the results in \cite{NY_twisting} have recently been studied by Giannone in  his thesis \cite{Giannone_thesis}, where the weak Hopf algebra of \cite{CP} plays the role of the universal enveloping algebra $U_q({\mathfrak g})$ for $q$ positive.
    \end{rem}
  
 \begin{defn} Following \cite{KW}, a  semisimple rigid tensor category  ${\mathcal C}$ together with      an isomorphism  of based rings $\phi_{\mathcal C}:  R_{N,\ell}\to {\rm Gr}({\mathcal C})$ 
   is called of  ${\mathfrak sl}_{N,\ell}$-type. Two  ${\mathfrak sl}_{N, \ell}$-type categories $(\phi_{\mathcal C}, {\mathcal C})$ and $(\phi_{{\mathcal C}'}, {\mathcal C}')$  are 
equivalent if there is a tensor equivalence ${\mathcal E}:{\mathcal C}\to{\mathcal C}'$ inducing an isomorphism between the Grothendieck rings compatible with $\phi_{\mathcal C}$ and $\phi_{{\mathcal C}'}$.
 \end{defn}

  The proof of Theorem \ref{classification_type_A} will occupy the rest of this section and it is based on Kazhdan-Wenzl theory \cite{KW}.  
     To summarize, Kazhdan-Wenzl   theory gives a classification of ${\mathfrak sl}_{N, \ell}$-type tensor categories in terms of the   categories arising from quantum groups both for generic or root of unity values of the deformation parameter $q$, and a $3$-cocycle
  on the dual of the center of ${\rm SU}(N)$   which modifies the natural associator. We start recalling    the main result. We shall then show that the positive ribbon structure completely determines the ribbon tensor category under our assumptions. The most delicate part of our analysis is a characterization of braided pseudo-unitary ${\mathfrak sl}_{N, \ell}$-type fusion categories among general ${\mathfrak sl}_{N, \ell}$-type categories, stated as Cor. \ref{classification_pseudo-unitary_braided_sl_N-type}, 
  and relies on the theory of quasitriangular weak  Hopf algebras developed in the paper. We also give   a parameterisation of the braided symmetries and a classification of their ribbon structures that is useful in our proof. 
\bigskip
    
\subsection{Proof of Theorem \ref{classification_type_A}, case   $\ell=N+1$.} \label{22.1}
  The based ring
    $R_{N, N+1}$ identifies with ${\mathbb Z}{\mathbb Z}_N$, with basis ${\mathbb Z}_N$ the cyclic group of order $N$.
 Hence a ${\mathfrak sl}_{N, N+1}$-type fusion category ${\mathcal C}$ is pointed over   ${\mathbb Z}_N$. By Prop 4.1 in \cite{KW},  see also Example \ref{pointed} and references therein,
   ${\rm Vec}_{{\mathbb Z}_{N}}^\omega$
exhaust  the  ${\mathfrak sl}_{N, N+1}$-categories, which are  classified by $\omega\in H^3({\mathbb Z}_{N}, {\mathbb T})$.
  A general braided pointed fusion category over the finite abelian group $G$   of equivalence classes of irreducible objects determines   a quadratic form on $G$ via $q(g)=c(\gamma, \gamma)$, where $g$ is the class of $\gamma$.
  The pair $(G, q)$ determines ${\mathcal C}$ as a braided tensor category by Theorem 8.4.9 in \cite{EGNO}. 
  By Remark 4.13 in
   \cite{Mueger2}, if $\theta$ is the ribbon structure associated to a braided symmetry $c$ and  a spherical structure in  a
   fusion category then on every object $X$, $\theta_X={\rm Tr}_X\otimes 1(c(X, X))$. In  a pointed fusion category
   $c(X, X)$ is a scalar and $d(X)=\pm1$ if $X$ is irreducible, and hence $d(X)=1$  under the  positivity requirement, and therefore
   $\theta_X=c(X, X)$. Hence 
  $q(g)=\theta_\gamma$. In other words the datum $(G, q)$ is equivalent to that of the fusion rules
and the positive ribbon structure. 
The result applies in particular to ${\mathfrak sl}_{N, N+1}$-type categories   and the proof   is complete in this case.
\medskip

As remarked in Ex. \ref{pointed_case}   these categories are   unitary in a natural way,   so the pseudo-unitarity assumption holds automatically.

\bigskip

\subsection{Kazhdan-Wenzl theory.}\label{22.2}
Examples of ${\mathfrak sl}_{N, \infty}$-type categories are the  representation categories of    quantum ${\mathfrak sl}_N$-groups for generic values of the deformation parameter.
Specifically, the quantum group of \cite{FRT} was  originally  considered in \cite{KW}. Being a quantization of a Hopf algebra of functions, the category is described by corepresentations. 
In the setting of tensor $C^*$-categories,    it is natural to consider the category of unitary corepresentations of Woronowicz ${\rm SU}_q(N)$  group, where $q$ is real, this is e.g. the starting point of \cite{P_rep, PR_rigidity, NY_twisting, Jordans}. We refer   to   \cite{CQGRC}    for details on the natural tensor  
$C^*$-structure. In a more general framework where a $C^*$-structure is not assumed,
   one may consider the category of representations of the
  Drinfeld-Jimbo quantum group $U_q({\mathfrak sl}_N)$ for   $q$ a non-zero complex number, not a nontrivial root of unity. 
By representations we  understand those which can be obtained as direct sums of subrepresentations of tensor products 
  of Weyl modules. 
   If $q$ is positive,  ${\rm SU}_q(N)$ and $U_q({\mathfrak sl}_N)$ induce equivalent tensor    categories, see \cite{CQGRC}. To unify with the examples ${\mathcal C}({\mathfrak sl}_N, q)$ at roots of unity,   we shall adopt Drinfeld-Jimbo framework.

 \medskip

In the following we assume   $\ell>N+1$. Then   ${\mathfrak sl}_{N, \ell}$-category is determined up to tensor equivalence by
 two invariants, $q_{\mathcal C}$ and $\tau_{\mathcal C}$, a pair of    nonzero complex numbers, unique up to passing to the pair with reciprocal values, which determines the tensor category, together with the fixed isomorphism $\phi$, up to   equivalence. These invariants are defined, and related to each other, 
 as follows.

Let $X\in{\mathcal C}$ be an object in the class of the image of $(1,0,\dots,0)$ under $\phi_{\mathcal C}$. The tensor product of $X$ with any irreducible is multiplicity free, and
  the   fusion rules
 can be found   in \cite{KW}.  Let
 $a\in(X^2, X^2)$ be the idempotent onto the subobject    $(1,1,0,\dots,0)$. Then there is 
  a nonzero complex   number  $q_{\mathcal C}$ (unique up to passing to the inverse) such that   $T:=q_{\mathcal C}(I-a)-a\in(X^2, X^2)$ gives rise via the usual   construction $T_i=1^{i-1}\otimes T\otimes 1^{n-i-1}$ to a representation of the braid group $\pi_n: {\mathbb B}_n\to(X^n, X^n)$. If $g_1,\dots, g_{n-1}$ are the generators of ${\mathbb B}_n$,   thus satisfying the presentation relations $g_ig_{i+1}g_i=g_{i+1}g_ig_{i+1}$,  $\pi_n$ takes $g_i\to T_i$.
  In our formulas, for simplicity, we are assuming that the category is strict. This representation factors through   the defining   relations $(g_i-q_{\mathcal C})(g_i+1)=0$, $i=1,\dots, n-1$, of the Hecke algebra $H_n(q_{\mathcal C})$ since $a$ is an idempotent.
Thus  we have representations of the Hecke algebras   denoted with the same symbol,
 $$\pi_n^+: H_n(q_{\mathcal C})\to (X^n, X^n)$$
compatible with the tensor structure. 
The ambiguity  in the choice of $q_{\mathcal C}$ also gives $\pi'_n: H_n({q_{\mathcal C}}^{-1})\to (X^n, X^n)$, which may equivalently
be thought of as another Hecke algebra representation on the same parameter $$\pi_n^-: H_n(q_{\mathcal C})\to (X^n, X^n),$$ the {\it opposite}, or {\it dual} representation
via $\pi_n^-:=\pi_n'\beta=\pi_n^=\alpha$
using the canonical isomorphism $\beta: H_n(q)\to H_n(q^{-1})$  which relates the corresponding canonical generators via 
$g_i \to -q h_i$, and $\alpha: g_i\in H_n(q)\to q-1-g_i\in H_n(q)$.

 Let ${\mathcal C}_{q, N, \ell}$ denote ${\mathcal C}({\mathfrak sl}_N, q, \ell)$ for $q^2$ a primitive root of unity of order $\ell$, for $\ell<\infty$ and ${\mathcal C}_{q, N, \infty}$
the category ${\rm Rep}(U_q({\mathfrak sl}_N))$ for $q$ not  a non-trivial root of unity. Note that
${\mathcal C}_{q, N, \ell}$   does not change, up to tensor equivalence,  under the passage from  $q$ to $q^{-1}$. This 
 may    be seen as follows. For $\ell=\infty$ there is an    isomorphism from  the quantum group $U_q({\mathfrak g})$,     to $U_{q^{-1}}({\mathfrak g})$
 given by $E_i\to K_i F_i$, $F_i\to E_i K_i^{-1}$, $K_i\to K_i$. For $\ell<\infty$ we may use an analogous isomorphism
 for $U_x({\mathfrak g})$, where $x$ is now an indeterminate, and the quantum group is regarded over ${\mathbb C}(x)$,    (see \cite{Sawin}, with our $x$ corresponding to $q$), and taking into account Lusztig's specialization of $U_x({\mathfrak g})$ to $U_q({\mathfrak g})$ for $q$ is a complex primitive root of unity. For details see   e.g. in Sect. 9.3, and 11.2 in \cite{Chari_Pressley}
 (for  $q$   of odd order) and
 \cite{Sawin}.

The category  ${\mathcal C}_{q, N, \ell}$ becomes an ${\mathfrak sl}_{N, \ell}$-type category as follows.  Set $X=X_q$,   the natural $N$-dimensional representation of  $U_q({\mathfrak sl}_N)$, and $\phi_q:R_{N, \ell}\to{\rm Gr}({\mathcal {\mathcal C}_{q, N, \ell}})$  the natural identification.  
 We realize $T$ as the element $-\sigma$ defined in $(4.13)$ of \cite{WoronowiczTK}, with $q$ in place of $\mu$ and consider 
the associated Hecke algebra representations $\pi_n$.

 For a general   ${\mathfrak sl}_{N, \ell}$-category,  it turns out that  $q_{\mathcal C}$ is a primitive root of unity of order $\ell$ for  $\ell$ finite, and is not a nontrivial root of unity for  $\ell$ infinite. In the first case, $H_n(q_{\mathcal C})$ is not semisimple for large values of $n$.   In both cases, the kernels of   $\pi_n^+$ and $\pi_n^-$ are completely determined by the fusion rules, and the two representations
   are distinguished by the value taken by a certain scalar invariant $\mu_{\mathcal C}$, see Theorem 4.1 in \cite{KW}, which corresponds to the value of a categorical  left inverse of $X$ on $T$, in the sense of  \cite{LR} in the Hecke category.
            
   The second invariant, called the twist of the category, is   given by 
 $\tau_{\mathcal C}=p\otimes 1_X\circ T_{1, N}\circ 1_X\otimes \nu\in(X, X)\simeq{\mathbb C}$, where $\nu\in(\iota, X^N)$ and $p\in(X^N, \iota)$ satisfy $p\circ\nu=1$ and $T_{1, N}= T_N\dots T_1$ is an Hecke algebra element in the representations $\pi_n$ exchanging the first 
 factor in a tensor product of $N+1$ objects with the following $N$ factors. 
More precisely, if the category is not strict,    $X^N=((X\otimes X)\otimes X)\dots$ and we need to use associativity morphisms in defining $\tau_{\mathcal C}$.

Given ${\mathcal C}$, with associativity morphisms $\alpha$, and given a $N$-th root of unity $w$, we may consider a new tensor category, ${\mathcal C}^w$ with the same representation ring, the
  same structure as ${\mathcal C}$ except for the associativity morphisms, which are modified as follows,
\begin{equation}\label{root}\alpha^w_{X_\lambda, X_\mu, X_\nu}:=w^{\gamma(|\lambda|, |\mu|)|\nu|}\alpha_{X_\lambda, X_\mu, X_\nu},\end{equation}
for $\lambda$, $\mu$, $\nu\in\Lambda$ (or in $\Lambda^+(q)$ accordingly), where $\gamma$ is the function  $\gamma(a,b)=[\frac{a+b}{N}]-[\frac{a}{N}]-[\frac{b}{N}]$ and $|\lambda|=m_1+\dots+m_{N-1}$.

It is easy to see that 
     $q_{\mathcal C}$ does not change when passing to a twisted category. This is not the case for $\tau_{\mathcal C}$,
     which does change and in fact 
 determines the root of unity $w$ defining the twist. Indeed, starting with a given ${\mathcal C}$ as before, if $\overline{X}$ is the conjugate of $X$ naturally realized as a subobject of $X^{N-1}$, we have $\nu\in (\iota, \overline{X}\otimes X)$. Taking into account the associativity morphisms, it follows that  $\tau_{\mathcal C}$ is
 the composite ($X_0=X$)
 \begin{equation}\label{firstline}X_0\stackrel{}\longrightarrow X_0(\overline{X} X)\stackrel{\alpha^{-1}}\longrightarrow(X_0\overline{X}) X\stackrel{T_{1, N-1}\otimes1 }\longrightarrow(\overline{X}X_0)X\end{equation}
\begin{equation}\label{secondline}\stackrel{\alpha}\longrightarrow\overline{X}(X_0X)\stackrel{1\otimes T_{1, 1}}\longrightarrow\overline{X}(XX_0)\stackrel{\alpha^{-1}}\longrightarrow(\overline{X}X)X_0\longrightarrow X_0.\end{equation}
Passing from ${\mathcal C}$  to   ${\mathcal C}^w$ gives rise to a modification in
the computation of the corresponding invariant
 only on the associativity morphisms. More precisely,
 the second part,   (\ref{secondline}), does not change,  by centrality of the deforming factor in $\alpha^w$,  
see (\ref{root}), while   (\ref{firstline}) changes by a factor $w^{-1}$. This   follows from a simple computation, since 
$\overline{X}$ corresponds to $(1,\dots, 1)$. Thus $\tau_{{\mathcal C}^w}=w^{-1}\tau_{\mathcal C}$.

  The following theorem is due to Kazhdan and Wenzl \cite{KW}. For completeness sake we include a proof.

\begin{thm}\label{Kazhdan_Wenzl}
Let ${\mathcal C}$ be a ${\mathfrak sl}_{N,\ell}$-type tensor category with $N+1<\ell\leq\infty$,
 $\phi_{\mathcal C}: R_{N,\ell}\to{\rm Gr}({\mathcal C})$ an isomorphism
and let $X$, $q_{\mathcal C}$ and $\tau_{\mathcal C}$ be defined as above. Then there is a $N$-th root of unity $w$ such that $\tau_{\mathcal C}=(-1)^N w^{-1} q^{N-1}$,
where $q$ is a complex square root of $q_{\mathcal C}$. The pair $(q_{\mathcal C}, \tau_{\mathcal C})$ is unique up to   the pair
with reciprocal values and determines the pair $({\mathcal C}, \phi_{\mathcal C})$ up to  equivalence. Furthermore,  there is an equivalence of
$({\mathcal C}, \phi_{\mathcal C})$ with $({\mathcal C}_{q, N, \ell}^w,$ $\phi_q$).
\end{thm}

   \begin{proof}
 Kazhdan-Wenzl left inverse
 $\mu_{\mathcal C}$   takes the value stated in Theorem 4.1 in \cite{KW} on $T$.
    It follows that the
  representation $\pi_n$ of the Hecke algebra is quasi equivalent to that arising from the quantum group  in ${\mathcal C}_{q, N, \ell}$.  
      In the generic case, a computation of $\tau_{\mathcal C}$ as in the statement may be found e.g. in    Lemma 8.1 of \cite{PR_rigidity}, version in arXiv.
  with $T_i$
    corresponding to $-g_i$ there, based on  a computation of the left inverse   on the generator $T$ for ${\mathcal C}={\mathcal C}_{q, N, \ell}$
   and the mentioned Hecke algebra representation of the quantum group, see
     Prop. 4.1 and Theorem 3.3 (a) in \cite{P_rep}, 
with $N$ and $\mu_{\mathcal C}$ in turn  corresponding to $d$ and 
  $\lambda_{-d}$ there. See also \cite{Jordans}.
 In the root of unity case, we may argue in the same way, using now Theorem 3.3 (b) \cite{P_rep} and  replacing
 $S$   with  the morphism still denoted $S$   of  the appendix of \cite{CP},
  and derive in  a similar way an $N$-th root of unity  $w$ such that $\tau_{\mathcal C}$ takes the stated value.  We then conclude following \cite{KW}:  up to  passing to   ${\mathcal C}^{w^{-1}}$,  we may     assume with no loss of generality that $\tau_{\mathcal C}=(-1)^Nq^{N-1}$, by  (\ref{firstline}), and (\ref{secondline}). We have thus reduced the values of the invariant $q_{\mathcal C}, \tau_{\mathcal C}$ to those it would take on ${\mathcal C}_{q, N, \ell}$.
It is easy to see that this  value of $\tau_{\mathcal C}$  in the twisted category   means that the element $\nu\in(\iota, X^N)$ of \cite{KW} and the Hecke algebra representations together
 satisfy the setting of section 6 in \cite{P_rep}, that is equations $(6.1)$--$(6.4)$, where 
 a (braided) tensor equivalence with  ${\mathcal C}_{q, N, \ell}$ has been  exhibited  for $q$ real taking $X$ to $X_q$, thus compatible with $\phi$ and $\phi_q$. More precisely,   braided symmetries are constructed from   certain normalizations of the Hecke algebra generator which is necessary to match $T$ with   the $R$-matrix of the quantum group in the representation $X$.
There is minimal change for other generic values of $q$. For the root of unity case, we may argue similarly, using  the information and analogous equations  in the appendix of \cite{CP} again. 
 \end{proof}

    It will be useful for us to specialize  Kazhdan-Wenzl theory to the untwisted tensor categories.
   In the following result, $\simeq$ denotes an equivalence  between pairs (${\mathcal C}_{q, N, \ell}^w$, $\phi_q$).

\begin{cor}\label{relative_class}
Let $q\in{\mathbb C}^\times$ be either not a non trivial root of unity or  such that $q^2$ is a primitive root of unity of order $\ell>N+1$, 
and let $q'\in{\mathbb C}^\times$ be another complex number with the same property.    Then: 

For $N$ even, 
\begin{itemize}
\item[{\rm        a)}] 
   ${\mathcal C}_{q, N, \ell}\simeq {\mathcal C}_{q', N, \ell}$   if  and only if
  $q'=q$ or $q'=\frac{1}{q}$; 
 \item[{\rm        b)}]  
 ${\mathcal C}_{-q, N, \ell}\simeq {{\mathcal C}^{-1}_{q, N, \ell}}$.\end{itemize}
   
  For $N$ odd,  ${\mathcal C}_{q, N, \ell}\simeq {\mathcal C}_{q', N, \ell}$ if  and only if
  $q'=\pm q$, $q'=\pm\frac{1}{q}$.

\end{cor}
\bigskip

\subsection{Braided symmetries in ${\mathfrak sl}_{N, \ell}$-type categories.}\label{22.3} Since the work of \cite{FRT, WoronowiczTK}
    and the theory of universal $R$-matrix of Drinfeld, see e.g. \cite{Chari_Pressley},
   it has been known  that    $U_q({\mathfrak g})$
      gives rise to braided tensor categories. For the  case of ${\mathcal C}({\mathfrak sl}_N, q)$ see e.g. \cite{Sawin}.
     There is a simple parameterisation of all the possible braided symmetries of ${\mathcal C}_{q, N, \ell}$.  
We start with the two canonical braided symmetries, $\varepsilon^+$ and its opposite $\varepsilon^-$ derived
 from the $R$-matrix of the quantum group and its opposite, $R_{21}^{-1}$,  respectively, see also remark \ref{roots}.

\begin{prop}\label{symmetries_list}
Let $z$ and $z'$ vary among the $N$-th roots of unity. Then for $N+1<\ell\leq\infty$ there are $2N$  braided symmetries, $\varepsilon^+_z$ and $\varepsilon_{z'}^-$ of 
${\mathcal C}_{q, N, \ell}$ uniquely determined by 
$$\varepsilon^+_z(X, X)=z\varepsilon^+(X, X), \quad\quad \varepsilon_{z'}^-(X, X)=z'\varepsilon^-(X, X).$$
Furthermore, this is a complete list.
\end{prop}
 
\begin{proof}
Since ${\mathcal C}_{q, N, \ell}$ admits  $X$ as a generating object, any braided symmetry $c$ is determined by $c(X, X)$ thanks to  (\ref{normalization_symmetry}), (\ref{braided_symmetry1}), (\ref{braided_symmetry2})
By the fusion rules of $X^2$, a suitable normalization of $c(X, X)$ will induce a representation of a Hecke algebra. By Kazhdan-Wenzl theory, the eigenvalue of the properly normalized $c(X, X)$
corresponding to $I-a$ can only be $q_{\mathcal C}^{\pm 1}$, so that  $c(X, X)$ is a scalar multiple of $\varepsilon^+(X, X)$
or $\varepsilon^-(X, X)$. By naturality of $c(X, X)$ on the morphism $\nu\in(\iota, X^N)$, the scalar
is a $N$-th root of unity. Conversely,
for any $N$-th root of unity $z$,   the modified morphisms $c_z(X^n, X^m)=z^{nm}c(X^n, X^m)$  still   satisfy 
the same relations and also
the  naturality property on the full subcategory with objects tensor powers of $X$, and hence everywhere,  as a consequence of  $(X^n, X^m)\neq0$ if and only if $n\equiv m({\rm mod} N)$. We may then apply these considerations to $\varepsilon^+$ and $\varepsilon^-$.
\end{proof}

  \begin{rem}\label{roots}
The   braided symmetries described in  the previous proposition are perhaps more clearly explained by the specialization process of the
$R$-matrix of the quantum group. More precisely, this matrix, at the level of 
the integral form $U_{{\mathcal A}'}^\dag({\mathfrak g})$ of
$U_x({\mathfrak g})$, with $x$ a formal variable as in  \cite{Sawin} where our $x$ corresponds to $q$ in that paper,
depends 
on
a  root  $s$ of order $L$ of $x$ via $s^L=x$, where $L$ is the smallest integer such that for any pair of dominant weights $\lambda$, $\mu$,  $L\langle \lambda, \mu\rangle$ is an integer. The values of $L$ are listed in table 1 in \cite{Sawin}.
We then specialize $x$ to a primitive complex root of unity $q$, and let $\ell'$ be its order
and $s$ to a fixed but arbitrary complex $L$-th root $q^{1/L}$ of $q$. 
Note that our $q^{1/L}$ is not necessarily a primitive root of unity of order $L\ell'$ as in Sect. 2 in \cite{Sawin}, thus our specialization 
needs to be slightly generalized.
One has $L=N$ for ${\mathfrak g}={\mathfrak sl}_N$.
It follows that $X\otimes X(R)$ corresponds to the operator
  computed in Sect. 8.3G of \cite{Chari_Pressley}, where  $e^h$ corresponds to $x$ and gives rise to our braided symmetry $\varepsilon^-$ through $\varepsilon^-(X, X)=\Sigma X\otimes X(R)$. 
The $N$ possible choices of $s=q^{1/N}$ give the symmetries   $z'\varepsilon^-$, and a similar relation holds between the specialization of opposite $R$-matrix $R_{21}^{-1}$ and  the symmetries $z\varepsilon^+$.
 \end{rem}
 
 The $2N$ braided symmetries of Prop. \ref{symmetries_list} give rise to  braided tensor categories 
 $({\mathcal C}_{q, N, \ell}, \varepsilon^\pm_z)$.
  We shall need the following property.

 \begin{rem}\label{inequivalence}
 Our aim is to show that the identity isomorphism between the   representation rings of any two of  $({\mathcal C}_{q, N, \ell}, \varepsilon^\pm_z)$
can not be induced by a braided tensor equivalence. 
An explicit proof of this fact     between two categories of the kind  $({\mathcal C}_{q, N, \ell}, \varepsilon^+_z)$ (or $({\mathcal C}_{q, N, \ell}, \varepsilon^-_z)$) which fixes the generating object $X$ may be found 
e.g.  at page 8 of  \cite{PR_rigidity}
 (arxiv version)
  for $q$ real.  Those arguments extend to a nonzero complex generic $q$ or to the root of unity case with the same modifications indicated in the proof of Theorem \ref{Kazhdan_Wenzl}. 
Since an isomorphism between two objects in a braided tensor category induces
 a braided tensor equivalence between the full braided tensor subcategories they generate, it
 also follows that there is no braided tensor equivalence   which takes the generating object $X$ to an equivalent object, and the conclusion follows in this case.
On the other hand for a pair of the kind $({\mathcal C}_{q, N, \ell}, \varepsilon^+_z)$ and  $({\mathcal C}_{q, N, \ell}, \varepsilon^-_{z'})$,  an argument may be found in the proof of Theorem \ref{classification_type_A} relying on the comparison of the ribbon structures.  
  \end{rem}
  
 We refer the reader also to \cite{Br1, Br2} for further studies on these braided symmetries.
     Up to a sign change of $q=(q_{\mathcal C})^{1/2}$ for $N$ even,   an ${\mathfrak sl}_{N,\infty}$-type {\it braided} tensor category  
  ${\mathcal C}$ is tensor equivalent to some   ${\rm Rep}(U_q({\mathfrak sl}_N))$. (The case $N=2$  holds without the   braided symmetry requirement, as it follows from the work of \cite{FK}, or also from 
  Theorem \ref{Kazhdan_Wenzl}, since  $H^3({\mathbb Z_2}, {\mathbb T})\simeq {\mathbb Z}_2$, see also Cor. \ref{relative_class}.) For $N>2$ a
 proof  has been given in Remark 4.4 of  \cite{NY_twisting}  based 
 on the Tannakian property of ${\rm Rep}(U_q({\mathfrak sl}_N))$ which provides
 a discrete  Hopf algebra.
We  need to extend this result    ${\mathfrak sl}_{N, \ell}$-type categories for  $\ell<\infty$.
However it is not obvious   how to modify   the methods of  \cite{NY_twisting} 
  for general   $q$ (with $q^2$ is a primitive root of unity of order $\ell$) as the categories 
${\mathcal C}({\mathfrak sl}_N, q)$  are not associated to Hopf algebras. 
Perhaps the most natural way to proceed is    to 
restrict to some    subclass large enough   to hold our applications.  
We shall thus first consider    only   the roots of unity  $q$ such that $q^2$ is of order large enough as in Def. \ref{large_enough}.
This  will   enable us to
  replace the role of the discrete Hopf algebra of  \cite{NY_twisting} 
 with  the  weak  Hopf algebra $A=A_W({\mathfrak sl}_N, q, \ell)$
  of Sect. \ref{20},
    its quasi-triangular structure developed in Sect. \ref{7}
  and  the notion of $3$-coboundary associator for weak Hopf algebras, Sect.  \ref{6}. 
  We shall include a proof since it  becomes slightly more technical due to non-triviality of the associator of $A$.

\begin{prop}\label{braiding_constraint} Let $({\mathcal C}, \phi_{\mathcal C})$ be an ${\mathfrak sl}_{N, \ell}$-type tensor category and assume that either $\ell=\infty$ or
$q=q_{\mathcal C}^{1/2}=\pm e^{\pm i \pi/\ell}$ (or, more generally $q^2=q_{\mathcal C}$ of order large enough as
in Def. \ref{large_enough}). Then ${\mathcal C}$ admits a braided symmetry  
  if and only if $w=1$ for $N$ odd and $w=\pm 1$ for $N$ even.
        \end{prop}

\begin{proof} The case $\ell=\infty$ ($q$ generic) has been   considered in \cite{NY_twisting}.
 By Kazhdan-Wenzl theory   an ${\mathfrak sl}_{N,\ell}$-type category $({\mathcal C}, \phi_{\mathcal C})$ is equivalent to  
    $({\mathcal C}({\mathfrak sl}_N, q, \ell))^w, \phi_q)$.
 For the case  $q=\pm e^{\pm i \pi/\ell}$ recall that the weak Hopf algebra 
  $A_W({\mathfrak sl}_N, q, \ell)$
 of Sect. \ref{20}  
has representation category tensor equivalent to ${\mathcal C}({\mathfrak sl}_N, q, \ell)$.
   Let, as before, $\Delta$ and $\Phi=1\otimes\Delta(P)\Delta\otimes1(P)$, $P=\Delta(I)$, be the natural coproduct and associator of $A$.    Consider the weak quasi bialgebra
   $A_w=(A, \Delta, \Phi_w)$,
with the   new associator $\Phi_w=\Phi\Upsilon_w$,  where $\Upsilon_w=\Upsilon\in A\otimes A\otimes A$
is the central invertible element  given by $\Upsilon=w^{\gamma(|\lambda|, |\mu|)|\nu|}$. 
Let us   regard ${\mathcal C}({\mathfrak sl}_N, q, \ell)$ as tensor equivalent to ${\rm Rep}(A)$   and therefore
${\mathcal C}({\mathfrak sl}_N, q, \ell)^w$ to ${\rm Rep}(A_w)$.
Let $R_q$ denote the $R$-matrix of $A$, hence by Prop. \ref{qtriangular_simplified}
   $$\Delta\otimes 1(R_q)=\Phi_{312}(R_q)_{13}(R_q)_{23}\Phi,\quad\quad 1\otimes\Delta(R_q)=\Phi_{231}^{-1}(R_q)_{13}(R_q)_{12}\Phi_{123}^{-1}.$$
If we assume that ${\mathcal C}$ is braided then so is ${\rm Rep}(A_w)$, hence by duality $A_w$ is quasi-triangular. Let $R$ be
the corresponding $R$-matrix. Thus  $R$ satisfies equations (\ref{qtriangular1})--(\ref{qtriangular4})   with respect to $\Phi_w$.
Since $\Upsilon_{123}=\Upsilon_{213}$,  
taking also into account the computations in the proof of Prop. 
\ref{qtriangular_simplified},
equations (\ref{qtriangular3})--(\ref{qtriangular4}) become 
$$\Delta\otimes 1(R)=\Phi_{312}R_{13}R_{23}\Phi\Upsilon,\quad\quad 1\otimes\Delta(R)=
\Upsilon_{231}^{-1}\Phi_{231}^{-1}R_{13}R_{12}\Phi_{123}^{-1}.$$
We consider the   twist $F=R_q^{-1}R$, cf. (\ref{prop1}), which satisfies $\Delta_F=\Delta$ and 
$$I\otimes F 1\otimes\Delta(F)=[I\otimes R_q1\otimes\Delta(R_q)]^{-1}I\otimes R1\otimes\Delta(R),$$
$$F\otimes 1\Delta\otimes 1(F)=[\Delta^{{\rm op}}\otimes 1(R_q) R_q\otimes 1]^{-1}\Delta^{{\rm op}}\otimes 1(R) R\otimes 1.$$ 
We set,   as before, $P=a\otimes b$, $\Delta(a)=a_1\otimes a_2$, $\Delta(b)=b_1\otimes b_2$ and compute
$$I\otimes R_q1\otimes\Delta(R_q)[\Delta^{{\rm op}}\otimes 1(R_q) R_q\otimes 1]^{-1}=$$
$$(R_q)_{23}\Phi_{231}^{-1}(R_q)_{13}(R_q)_{12}\Phi_{123}^{-1}(R_q)_{12}^{-1}\Phi_{213}^{-1}(R_q)_{13}^{-1}(R_q)_{23}^{-1}\Phi_{321}^{-1}=$$
$$(R_q)_{23}\Phi_{231}^{-1}(R_q)_{13}(R_q)_{12}\Delta\otimes1 (P)1\otimes\Delta(P)\Delta\otimes 1(P)(R_q)_{12}^{-1}b_1\otimes a\otimes b_2(R_q)_{13}^{-1}(R_q)_{23}^{-1}\Phi_{321}^{-1}=$$
$$(R_q)_{23}\Phi_{231}^{-1}(R_q)_{13}(a_2\otimes a_1\otimes b)( b_1\otimes a\otimes b_2)(R_q)_{13}^{-1}(R_q)_{23}^{-1}\Phi_{321}^{-1}=$$
$$(R_q)_{23}(b\otimes a_1\otimes a_2)(R_q)_{13}(1\otimes\Delta(P)\Delta\otimes 1(P)1\otimes\Delta(P))_{213}(R_q)_{13}^{-1}(R_q)_{23}^{-1}\Phi_{321}^{-1}=$$
$$(R_q)_{23}(b\otimes a_1\otimes a_2)(b_2\otimes a\otimes b_1)(R_q)_{23}^{-1}\Phi_{321}^{-1}=$$
$$(R_q)_{23}(\Delta\otimes 1(P)1\otimes\Delta(P)\Delta\otimes 1(P))_{231}
(R_q)_{23}^{-1}(b_2\otimes b_1\otimes a)=$$
$$(b\otimes a_2\otimes a_1)(b_2\otimes b_1\otimes a)=\Phi_{321}^{-1}.$$
Hence, using centrality of $\Upsilon$,
$$[I\otimes F 1\otimes\Delta(F)]^{-1}F\otimes 1\Delta\otimes 1(F)=$$
$$[I\otimes R1\otimes\Delta(R)]^{-1}I\otimes R_q1\otimes\Delta(R_q)[\Delta^{{\rm op}}\otimes 1(R_q) R_q\otimes 1]^{-1}\Delta^{{\rm op}}\otimes 1(R) R\otimes 1=$$
$$\Upsilon_{231}\Phi_{123}R_{12}^{-1}R_{13}^{-1}\Phi_{231}R_{23}^{-1}\Phi_{321}^{-1}\Phi_{321}R_{23}R_{13}\Phi_{213}R_{12}\Upsilon_{123}=$$
$$\Upsilon_{231}\Phi_{123}R_{12}^{-1}R_{13}^{-1}\Phi_{231}R_{13}\Phi_{213}R_{12}\Upsilon_{123}=\Phi\Upsilon_{231}\Upsilon,$$
we have omitted the computations leading to the last equality, as they are very similar to the previous ones.
Hence $\Phi\Upsilon_{231}\Upsilon$ satisfies (\ref{coho1}), and  one may similarly establish validity (\ref{coho2}), thus $\Phi\Upsilon_{231}\Upsilon$ is a $3$-coboundary associator which may be twisted to $\Phi$ by $F$ by Prop. \ref{2-cocycle}.
On the other hand  as observed in \cite{NY_twisting}
$({\Upsilon_w})_{231}\Upsilon_w$ is cohomologous to $\Upsilon_{w^2}$ on the dual of the center of ${\rm SU}(N)$, and therefore we find a tensor equivalence between 
${\mathcal C}({\mathfrak sl}_N, q, \ell)$ and ${\mathcal C}({\mathfrak sl}_N, q, \ell)^{w^2}$ which identifies the generating representations, and hence is compatible with the chosen isomorphisms with
$R_{N, \ell}$. From Kazhdan-Wenzl classification we derive $w^2=1$ and we finally apply Cor. \ref{relative_class}.
\end{proof}

 \bigskip
 
\subsection{Pseudo-unitarity.} \label{22.4}
 If a given ribbon structure for the braided symmetry of a fusion category ${\mathcal C}$ induces a spherical structure making   the categorical dimensions of the simple objects positive then ${\mathcal C}$ is pseudo-unitary. It follows that ${\mathcal C}({\mathfrak sl}_N, q, \ell)$ is pseudo-unitary for $q=\pm e^{\pm i\pi/\ell}$ and $N+1<\ell<\infty$
 with respect to the natural ribbon structure, by  \cite{Andersen} (and in fact unitary by \cite{Wenzl, Xu}). In this subsection we prove that these fusion categories    
may be    intrinsically characterized 
    among general  fusion ${\mathfrak sl}_{N, \ell}$-type categories   by    the property of being both braided and pseudo-unitary.

 \begin{prop}\label{pseudo-unitarity_qg} Let $q\in{\mathbb C}$ be such that $q^2$ is a non-trivial root of unity of order 
  $\ell>N+1$. Then  ${\mathcal C}({\mathfrak sl}_N, q, \ell)$ is pseudo-unitary  if and only if $q=\pm e^{\pm i\pi/\ell}$.
 \end{prop}
 
 \begin{proof} Our proof follows that of an analogous result  for the Lie type $B$   given  in Theorem 3.8 in \cite{Rowell3}, with a slight modification due to the non-uniqueness of the spherical structures for $N$ even in our case, see the following Prop. \ref{spherical_structures}.  More in detail,
  we write $q=\pm q_z$, with $q_z=e^{i\pi z/\ell}$ and $z$   an integer with $1\leq z\leq \ell-1$ and
 $\gcd(z,\ell)=1$. Let $X$ be the object  of ${\mathcal C}({\mathfrak sl}_N, q, \ell)$  corresponding to the fundamental representation and assume $N=2k$ even. Up to a sign,   the categorical dimension $d(X)$ with respect to any spherical structure   equals
 $d_{q_z}(X):={q_z}^{N-1}+{q_z}^{N-3}+\dots+{q_z}^{-(N-1)}=2\sum_{j=1}^k\cos((2j-1)\pi z/\ell)$.
 Furthermore,   ${\rm FPdim}(X)=d_{q_1}(X')$ where $X'$ is a corresponding object in the category 
${\mathcal C}({\mathfrak sl}_N, q_1)$ for $q_1=e^{i\pi/\ell}$, since 
 these two categories have isomorphic representation rings with an isomorphism identifying $X$ to $X'$ and we know  that $d_{q_1}$ takes positive values on the irreducibles. We claim that 
 $d_{q_z}(X)<d_{q_1}(X')$ for $z\neq 1$. Thus if $d_{q_z}(X)>0$ then $|d(X)|=d_{q_z}(X)$ which then can equal ${\rm FPdim}(X)$ only if $q=\pm e^{i\pi/\ell}$. If $d_{q_z}(X)<0$ then $|d(X)|=-d_{q_z}(X)=2\sum_{j=1}^k\cos((2j-1)\pi(\ell- z)/\ell)$. Since $\ell-z$ satisfies the same properties as $z$, pseudo-unitarity again implies  $\ell-z=1$ hence $q=\pm e^{-i\pi/\ell}$. To show the claim, observe that the set $S_1$ of points ${q_1}^{2j-1}=e^{i(2j-1)\pi/\ell}$, $j=1,\dots k$ all lie in the upper semicircle.
 Furthermore the conditions $\gcd(z, \ell)=1$ and $\ell\geq N+2$ imply $\ell\nmid(2j-1)z$. In particular ${q_z}^{2j-1}\neq 1$
 for all $j$. Assume that $z$ is such that the subset $S_2$ of $\{q_z^{\pm(2j-1)}\}$ contained in the upper  semicircle
 differs from   $S_1$.
 The first point  in the natural order  of the semicircle is $q_1$. Furthermore two adjacent points of  $S_2$    correspond to arcs whose distance is at least $2\pi/\ell$. Therefore there must be an element of $S_1$ in between unless they both lie after the last $q_1^{N-1}$.  Since $\cos$ is an even function, it follows that $d_{q_z}(X)$ may be computed considering elements of $S_2$, and
 we have
 $d_{q_z}(X)<d_{q_1}(X')$ by the above remarks. We are left to show that for $z\neq 1$, $S_1\neq S_2$. For this we may apply   arguments analogous to those of the  last part of the mentioned theorem of \cite{Rowell3}.
 
 In the case where $N$ is odd the proof is   simplified by the fact that   
${\mathcal C}({\mathfrak sl}_N, q, \ell)$ admits a unique spherical structure,  so      $d(X)$  is   uniquely determined. We may thus complete the proof    with argument similar to the even case, taking into account the additional information that   $d(X)=q_z^{N-1}+\dots+q_z^{-(N-1)}=-1+2\sum_{j=0}^k\cos(2j\pi z/\ell)$ where $N=2k+1$.

 \end{proof}
 
 For completeness we recall from Example \ref{pointed_case} that ${\mathcal C}({\mathfrak sl}_N, q, \ell)$ are always unitary if $q^2$ is a primitive root of 
unity of order $\ell=N+1$.  
Here below we remark about classification of spherical structures on ${\mathcal C}_{q, N, \ell}$.
 
 \begin{prop}\label{spherical_structures} For $N$ odd, ${\mathcal C}_{q, N, \ell}$ has a unique spherical structure, for $N$ even it has two.\end{prop}
 
 \begin{proof} In a fusion category ${\mathcal C}$ spherical structures are parameterised by the group of monoidal natural transformations from the identity functor to itself and taking values $\pm 1$ on the irreducibles, see Exercise 4.7.16 of \cite{EGNO} for a precise statement. In the case where ${\mathcal C}$ admits a simple generating object $X$,
 any such natural transformation $\eta$  is   determined by the value it takes on $X$ as follows.
 If $\eta_X=\lambda 1_X$ then on any tensor power, $\eta_{X^{r}}=\lambda^r 1_{X^r}$ by monoidality. 
  It follows from naturality and complete reducibility that the values that $\eta$ takes on the simple summands of $X^r$ 
 also coincide with $\lambda^r$. Hence if $\lambda=1$ then $\eta$ is the identity natural transformation, while if $\lambda=-1$ then $\eta$ takes value $1$ ($-1$) on the simple summands of even (odd) tensor powers of $X$.
 In our specific case, if $X=X_q$ we must have $\lambda^N=1$ since the tensor unit is a subobject of $X_q^N$.
 Hence for $N$ odd the conclusion follows. For $N$ even,   the specific fusion rules of an ${\mathfrak sl}_N$-type
 tensor category show that any odd  tensor power of $X_q$ is disjoint from an even tensor power.  This implies existence of 
  a monoidal natural transformation $\eta\in(1, 1)$ taking these values.   \end{proof}

 The next step is that of      characterizing  general pseudo-unitary  ${\mathfrak sl}_{N, \ell}$-type fusion categories for $\ell>N+1$. To do this, we regard
 the relationship between  
 ${\mathcal C}({\mathfrak sl}_N, q, \ell)^w$, and  
${\mathcal C}({\mathfrak sl}_N, q, \ell)$ as an example of a general construction described in \cite{BNY_grading} of a new fusion category ${\mathcal C}^\omega$
 from from a given one ${\mathcal C}$ and a ${\mathbb T}$-valued $3$-cocyle $\omega$ on the chain group 
 ${\rm Ch}({\mathcal C})$, and we study invariance of
 pseudo-unitarity  under  $\omega$ in this   framework.

 Let ${\mathcal C}$ be a semisimple monoidal category  with associativity morphisms $\alpha$. The chain group
   ${\rm Ch}({\mathcal C})$ introduced in 
 \cite{Baumgartel_Lledo, GN}
  is defined as follows. Consider a complete family ${\rm Irr}({\mathcal C})=\{\rho_\alpha, \alpha\in A\}$
  of simple objects of ${\mathcal C}$
 endowed with the smallest equivalence relation $\simeq$ making  all the irreducible subobjects $\rho_\gamma$ appearing 
 in the decomposition of $\rho_\alpha\otimes\rho_\beta$ for fixed $\alpha$, $\beta\in A$, equivalent.
 Then ${\rm Ch}({\mathcal C})={\rm Irr}({\mathcal C})/\simeq$ is a group with $[\rho_\alpha][\rho_\beta]=[\rho_\gamma]$.
 The trivial element is the class of the tensor unit, and $[\rho_\alpha]^{-1}=[\overline{\rho_\alpha}]$.
 This is an interesting group. For example, it identifies naturally with the dual of the centre of the compact group $G$ for ${\mathcal C}={\rm Rep}(G)$ \cite{Mueger_center}.  Furthermore, the group of nonzero ${\mathbb C}$-valued homomorphisms on ${\rm Ch}({\mathcal C})$ identifies with the group of natural monoidal transformations of the identity functor on ${\mathcal C}$
 \cite{GN, BNY_grading}. Finally, for modular categories, ${\rm Ch}({\mathcal C})$ identifies with the dual
 of the (abelian) group of invertible elements of ${\mathcal C}$ \cite{GN}.
 
The chain group induces a grading on ${\mathcal C}$, 
 in the sense that there are full subcategories ${\mathcal C}_g$ indexed by elements of $g\in{\rm Ch}({\mathcal C})$ such that every object $\rho\in{\mathcal C}$ decomposes uniquely up to isomorphism into a direct sum of objects $\rho_g\in{\mathcal C}_g$
 and with the property that for $g\neq h$, objects of ${\mathcal C}_g$ are disjoint from objects of  ${\mathcal C}_h$.
 The group structure of ${\rm Ch}({\mathcal C})$   implies that the grading is compatible with the tensor structure: $\iota\in{\mathcal C}_e$ and $\rho\otimes\sigma\in{\mathcal C}_{gh}$ for $\rho\in{\mathcal C}_g$, $\sigma\in{\mathcal C}_h$.

We   consider  ${\mathcal C}^\omega$,   the    monoidal category with the same structure as ${\mathcal C}$ except for 
the   associativity 
morphisms, which are given by $\alpha^\omega_{\rho, \sigma, \tau}=\omega(g, h, k)\alpha_{\rho, \sigma, \tau}$, where $[\rho]=g$, $[\sigma]=h$, $[\tau]=k$. Note that ${\mathcal C}^\omega$ may be regarded as a special case of a categorical analogue
of Prop. \ref{deforming_Phi}. In other words, ${\mathcal C}$ and ${\mathcal C}^\omega$ have isomorphic Grothendieck rings and chain groups, and, in the framework of fusion categories, they have the same Frobenius-Perron dimension function.
We denote by $d_{\mathcal C}(\rho)$ and $d_{{\mathcal C}^\omega}(\rho)$ the categorical dimensions of an object $\rho$ considered in ${\mathcal C}$ or ${\mathcal C}^\omega$ respectively with respect to preassigned spherical structures.

   \begin{prop}\label{pseudo-unitarity_invariance}
   Let ${\mathcal C}$ be a  fusion category and $\omega\in Z^3({\rm Ch}({\mathcal C}); {\mathbb T})$ normalized. If $D$ is the right duality
   functor of ${\mathcal C}$ associated to the right duality
   $(\rho^\vee, b_\rho, d_\rho)$ and $\eta\in(1, D^2)$ is a pivotal (spherical) structure then $(\rho^\vee, b^\omega, d^\omega)$
    is a right duality for ${\mathcal C}^\omega$ where $b^\omega_\rho=b_\rho$, $d^\omega_\rho=d_\rho\omega^{-1}(g, g^{-1}, g)$, with $\rho$ simple and $[\rho]=g$. Furthermore
      $\eta^\omega=\eta$
     is a    pivotal (spherical) structure for the associated right duality functor $D_\omega$.
       In particular, if $\eta$ is spherical under the correspondence $(\eta, D)\to(\eta^\omega, D_\omega)$ we have
   $d_{\mathcal C}(\rho)=d_{{\mathcal C}^\omega}(\rho)$ for every object $\rho$.
Furthermore, ${\mathcal C}$ is pseudo-unitary if and only if  so is ${\mathcal C}^\omega$.
   \end{prop}
   
   \begin{proof}
 Let $\rho$ be an object of ${\mathcal C}_g$ and let $(b, d)$, $(b', d')$ solve the right and left duality equations respectively
for $\rho$ in ${\mathcal C}$ in the sense of (\ref{first_right})--(\ref{second_left}) with $\rho^\vee={}^\vee\rho$. Then   a solution 
of the corresponding equations in ${\mathcal C}^\omega$ is given
by
$(b^\omega, d^\omega)$, $(b'^\omega, d'^\omega)$  where $b^\omega=b$, $d^\omega=d\omega^{-1}(g, g^{-1}, g)$, $b'^\omega=b'\omega(g, g^{-1}, g)$, $d'^\omega=d'$.
To verify the duality relations it is useful to recall  the equality $\omega(g, g^{-1}, g)=\omega(g^{-1}, g, g^{-1})^{-1}$ which follows from  the $3$-cocycle equation for $\omega$.

We now start with a   right duality $(\rho^\vee, b_\rho, d_\rho)$ in ${\mathcal C}$
and recall that the associated right duality  functor $D$     was defined in (\ref{duality_functor}). The right duality functor $D_\omega$
 of ${\mathcal C}^\omega$ associated with the   solution $(\rho^\vee, b^\omega_\rho, d^\omega_\rho)$ of the previous paragraph acts as $D$ on objects, while on morphisms $T\in(\rho, \sigma)$
 with $\sigma\in{\mathcal C}_h$ we have $D_\omega(T)=\omega^{-1}(h, h^{-1}, h)D(T)$. Let $\eta\in(1, D^2)$ be
 a pivotal structure. 
 Consider  
the left duality $(\rho^\vee, b'_\rho, d'_\rho)$ defined by
 (\ref{left_duality}) with $\eta$ in place of $u$. 
It follows that ${\mathcal C}^\omega$ has left duality $(\rho^\vee, b'^\omega_\rho, d'^\omega_\rho)$.
 The natural transformation $\eta^\omega$
 in ${\mathcal C}^\omega$ 
 defined by (\ref{natural}) with $d'_\omega$ and $b_\omega$ in place of $d'$ and $b$ takes the same values as $\eta$.
 Furthermore the natural transformation say $F_{\rho, \sigma}$ in ${\mathcal C}$ making $D^2$ into a tensor functor is also natural in 
 ${\mathcal C}^\omega$ and makes $D_\omega^2$ into a tensor functor. Indeed,   it is easy to see that $D_\omega^2$ acts as $D^2$ on $\alpha_{\rho, \sigma,\tau}$ if $\rho$, $\sigma$, $\tau$ are homogeneous, and therefore in general. It follows
 that  validity of Def. \ref{wtf} for $D^2$ implies validity for $D_\omega^2$ by linearity. Hence
 $\eta^\omega$ is monoidal by Def.  \ref{monoidal_transformation}, and therefore is a pivotal structure in ${\mathcal C}^\omega$ which is spherical if so was $\eta$. The formulas also show that $d_{\mathcal C}(\rho)=db'=d_\omega b'_\omega=d_{{\mathcal C}^\omega}(\rho)$ with respect to these structures. Since ${\mathcal C}$ and ${\mathcal C}^\omega$ have the same global ${\rm FPdim}$, the last assertion is also clear.

   \end{proof}

\begin{cor}\label{classification_pseudo-unitary_braided_sl_N-type} Among the  ${\mathfrak sl}_{N, \ell}$-type tensor categories    $({\mathcal C}, \phi_{\mathcal C})$ with   $N+1<\ell<\infty$ only those equivalent to some 
$({\mathcal C}({\mathfrak sl}_N, q, \ell))^w, \phi_q)$
$(({\mathcal C}({\mathfrak sl}_N, q, \ell), \phi_q) \text{ resp.})$ with $q=e^{ i\pi/\ell}$ for $N$ odd and $q=\pm e^{ i\pi/\ell}$ for $N$ even are pseudo-unitary (pseudo-unitary   and braided resp.).
  \end{cor}

\begin{proof}
This follows immediately from Propositions \ref{pseudo-unitarity_qg}, \ref{pseudo-unitarity_invariance}, \ref{braiding_constraint}, Cor. \ref{relative_class}.
\end{proof}
\bigskip

\subsection {Proof of Theorem \ref{classification_type_A}, case  $\ell>N+1$.}  \label{22.5}
\begin{proof} Let us fix an isomorphism of based rings $\phi_{\mathcal C}: R_{N, \ell}\to {\rm Gr}({\mathcal C})$. Then $({\mathcal C}, \phi_{\mathcal C})$ is equivalent to
   $({\mathcal C}({\mathfrak sl}_N, q, \ell), \phi_q)$ with $q=e^{i\pi/\ell}$ for $N$ odd and precisely to one of  
   $({\mathcal C}({\mathfrak sl}_N, q, \ell), \phi_q)$
   where $q$ takes the values
   $q= \pm e^{i\pi/\ell}$ for $N$ even, according to an equivalence ${\mathcal E}$ inducing $\phi_{\mathcal C}$, by Cor. \ref{classification_pseudo-unitary_braided_sl_N-type}. 
 A similar conclusion holds for $({\mathcal C}', \phi_{{\mathcal C}'})$ for any choice of $\phi_{{\mathcal C}'}: R_{N, \ell}\to {\rm Gr}({\mathcal C}')$.
We fix    $\phi_{{\mathcal C}'}=f\circ \phi_{\mathcal C}$, and denote by
  ${\mathcal E}'$   the corresponding equivalence 
 with $({\mathcal C}({\mathfrak sl}_N, q', \ell), \phi_{q'})$. Using the based ring  isomorphisms induced by  ${\mathcal E}$ and ${\mathcal E}'$ between the Grothendieck rings, their compatibility
 with $\phi_{\mathcal C}$ and $\phi_{{\mathcal C}'}$, we find an isomorphism $g: {\rm Gr}({\mathcal C}({\mathfrak sl}_N, q, \ell))\to{\rm Gr}({\mathcal C}({\mathfrak sl}_N, q', \ell))$ which identifies the classes of the respective generating representations $X_q$ and $X_{q'}$. Let us now take into consideration the braided symmetries, say $c$ and $c'$ of ${\mathcal C}$ and ${\mathcal C}'$ respectively,
 and their  ribbon structures, identified with analogous structures in    the quantum group categories via the   equivalences and denoted in the same way. 
 For ${\mathcal C}({\mathfrak sl}_N, q, \ell)$ we can only have $c=z\varepsilon^+$ or $c=z'\varepsilon^-$ by Prop. \ref{symmetries_list},
where $z$ and $z'$ have the    same meaning.  Taking into consideration Remark \ref{roots},   we identify each of the $2N$ possible braided symmetries with one derived from the $R$-matrix $R$ or the opposite $R_{21}^{-1}$, subject to a   choice of a complex $N$-th root 
 $q^{1/N}$. Then ${\mathcal C}({\mathfrak sl}_N, q, \ell)$ becomes a ribbon category with positive ribbon structure 
 $\tilde{\theta}_\lambda=q^{\pm {\langle \lambda, \lambda+2\rho\rangle}}$, where  $\langle\ , \ \rangle$ is a symmetric invariant bilinear form of ${\mathfrak sl}_N$ such that $\langle \alpha, \alpha\rangle=2$ for (short) roots, the plus or minus sign
are determined by the choice of $R$ or $R_{21}^{-1}$, see \cite{Chari_Pressley}. 
   On the other hand ${\mathcal C}({\mathfrak sl}_N, q, \ell)$ also has the positive ribbon structure $\theta_\lambda$, hence $\theta_\lambda=\tilde{\theta}_\lambda$ by   uniqueness of the positive ribbon structures recalled before the statement of Theorem \ref{classification_type_A}.
      Assuming that $\theta$ corresponds to the plus sign we have that   $\theta_{X_q}=q^{\frac{N^2-1}{N}}$ (more details on this formula may be found in the proof of the following proposition).
    We claim that we may assume that $\theta'$ corresponds to a plus sign as well. Hence we similarly have
   $\theta_{X_{q'}}={q'}^{\frac{N^2-1}{N}}$.  If
   $N$ is odd we have already settled $q=q'$ and our assumption $\theta_{X_q}=\theta_{X_{q'}}$ shows that we are taking the same $N$th root of $q$, and therefore we have a braided, in fact ribbon, tensor equivalence. If $N$ is even then $N^2-1$ is odd, and since $\theta_{X_q}^N=\theta_{X_{q'}}^N$ we may exclude that $q$ and $q'$ have opposite signs. It follows again that the two $N$th roots of $q$ are the same and we get the same conclusion.
   We finally show the claim. If on the contrary we had an opposite symmetry $c' $  in ${\mathcal C}'$   then   $q={q'}^{-1}$ for $N=2$ and $q^2= q'^{-2}$ for $N>2$. 
In the first case we conclude as before since by Prop. \ref{relative_class} $q$ and $q^{-1}$ gives rise to equivalent tensor categories again. In the second case we use the twist equation
$c(X_q, X_q)^2=\theta_{X_q}\otimes\theta_{X_q}\circ\theta_{X_q\otimes X_q}^{-1}$ and similarly for $X_{q'}$, $\theta'$ and $c'$, which implies $c(X_q, X_q)^2$ and $c'(X_{q'}, X_{q'})^2$ have the same eigenvalues.
This implies $q_{\mathcal C}^4=1$ and therefore $N+1<\ell\leq 4$ giving no solution.
\end{proof}

 Note that  the positivity assumption in Theorem \ref{classification_type_A} is redundant for $N$ odd by uniqueness of the
 ribbon structure of every braided symmetry of  ${\mathcal C}_{q, N, \ell}$, Prop. \ref{spherical_structures}.
 The following example shows that  this assumption can not be dropped  for $N$ even.

\begin{ex}\label{non_uniqueness}
Consider ${\mathcal C}({\mathfrak sl}_2, q, \ell)$ for $q=e^{i\pi/\ell}$ with $\ell>3$ and the braided symmetries $\varepsilon^+_1$ and $\varepsilon^+_{-1}$ described in Prop. \ref{symmetries_list}. 
By Remark \ref{inequivalence},   the identity isomorphism between the corresponding representation rings
can not be induced by a braided tensor equivalence. On the other hand, 
each of the two braided categories has its own positive ribbon structure, say $\theta_1$ and $\theta_{-1}$ respectively.
For an irreducible $\lambda=a\Lambda_1$ we have $\theta_1(\lambda)=q^{\frac{a}{2}(a+2)}$, 
$\theta_{-1}(\lambda)=(-q^{1/2})^{{a}(a+2)}$ where $q^{1/2}=e^{i\pi/2\ell}$.  If $\eta\in(1,1)$ is the natural monoidal transformation of the identity functor taking value $-1$ on the generating object $X$ then 
 it   follows from the proof of Prop. \ref{spherical_structures} that $\eta\theta_{-1}=\theta_1$.
But $\eta\theta_{-1}$ is another
ribbon structure for $\varepsilon^+_{-1}$.

\end{ex}

We conclude the section with a partial result concerning ribbon equivalence of   examples of  ${\mathfrak sl}_{N, \ell}$-type categories where    pseudo-unitarity is not assumed but the ribbon structure is fixed.
 
\begin{prop}\label{twist_general_case} Let $q$ and $q'\in{\mathbb C}$ be either not non-trivial roots of unity or else square to primitive roots of unity of order
 $\ell>N+1$ and let us endow both
${\mathcal C}({\mathfrak sl}_N, q, \ell)$ and ${\mathcal C}({\mathfrak sl}_N, q', \ell)$ with some braided symmetry. If there is an isomorphism of based rings
  $f:\text{Gr}({\mathcal C}({\mathfrak sl}_N, q, \ell))\to \text{Gr}({\mathcal C}({\mathfrak sl}_N, q', \ell))$ identifying the generating representations and compatible with
  the canonical ribbon structures then   there is a ribbon tensor equivalence 
  ${\mathcal F}: {\mathcal C}({\mathfrak sl}_N, q, \ell)\to{\mathcal C}({\mathfrak sl}_N, q', \ell)$ inducing $f$,
  $q$ and $q'$ are related as in Prop. \ref{relative_class}, and we may arrange $q=q'$.
\end{prop}

\begin{proof}
 We write the respective  braided symmetries $c$ and $c'$ as  in the proof of the previous theorem, where now $q$ and $q'$ are general. We again have that the canonical ribbon structure of ${\mathcal C}({\mathfrak sl}_N, q, \ell)$ takes the form   
 $\theta_\lambda=q^{\pm {\langle \lambda, \lambda+2\rho\rangle}}$.   We need to be a bit more explicit on the exponents, so we write $\lambda=\sum_1^{N-1} n_j\Lambda_j$, where $\Lambda_j$ are the fundamental weights, $n_j$ are non-negative integers and $\rho=\sum_1^{N-1} \Lambda_j$.
   Then $\langle \Lambda_k, \Lambda_j\rangle=d_j d_{k,j}$, where $d_{k,j}$ are such that $\Lambda_k=\sum_j d_{k,j}\alpha_j$, with $\alpha_j$ the simple roots,   $d_j=\frac{\langle \alpha_j, \alpha_j\rangle}{2}$, hence  equal to $1$ in our case. 
 This gives $\langle \lambda, \lambda+2\rho\rangle=
 \sum_{k,j}n_k(n_j+2)d_{k,j}$.    The matrix $(d_{k,j})$ is given in Table 1 at pag. 69 of \cite{Humphreys}.
 In particular one obtains $\langle \Lambda_k, \Lambda_k+2\rho\rangle=\frac{k}{N}(N-k)(N+1)$, see e.g. Sect. 6 in 
 \cite{PR_rigidity}, 
  and more generally
 $$\langle n\Lambda_k, n\Lambda_k+2\rho\rangle=n[\langle \Lambda_k, \Lambda_k+2\rho\rangle +(n-1)d_{k, k}]=$$
 $$\frac{n}{N}[{k(N-i)(N+1)}+(n-1){k(N-k)}]=\frac{nk}{N}(N-k)(N+n).$$
 Assuming  again that   $\theta$ corresponds to the     plus sign, we have  $$\theta_{\Lambda_1}=q^{N-\frac{1}{N}}, \quad\quad \theta_{2\Lambda_1}=q^{2N+2-\frac{4}{N}}$$ 
and for $N>2$ we   in addition consider
 $$\theta_{\Lambda_2}=q^{2N-2-\frac{4}{N}}.$$
It follows that $$\theta_{2\Lambda_1}\theta_{\Lambda_1}^{-2}=q^{2-\frac{2}{N}},\quad\quad\theta_{2\Lambda_1}\theta_{\Lambda_2}^{-1}=q^4.$$
      We claim that we may assume that $\theta'$ corresponds to a plus sign as well, and we show it in the same way. The first equation   gives  $q=q'$ for $N=2$. Assuming   $N>2$,     the second equation gives
  $q'=\pm q$ or $q'=\pm iq$.
a) Case  $q'=-q$.
If $N$ is odd then $q$ and $-q$ give rise to equivalent tensor categories by Prop. \ref{relative_class}. We may thus assume with no loss of 
generality that $q=q'$.  If $N$ is even then $N^2-1$ is odd and since 
$\theta_{\Lambda_1}^N={\theta'}_{\Lambda_1}^N$
we may exclude $q'=-q$. 
b) We next show that the cases $q'=\pm iq$ are  not realized.  We need to compute the ribbon structure of weights
which are sums of different fundamental weights, and for this we  use the following addition formula
which follows from bilinearity and symmetry of the inner product 
$$\langle \Lambda_1+\lambda, \Lambda_1+\lambda+2\rho \rangle=\langle\Lambda_1, \Lambda_1+2\rho\rangle+\langle\lambda, \lambda+2\rho\rangle+2\langle\lambda, \Lambda_1\rangle.$$
On the other hand, the   equation $\theta_{\Lambda_1}^N={\theta'}_{\Lambda_1}^N$   requires $N$ odd, 
we may thus consider the weights $\mu=\Lambda_{\frac{N-1}{2}}$ and $\nu=\Lambda_{\frac{N+1}{2}}$ and since
 $\langle\Lambda_k, \Lambda_1\rangle=\frac{k}{N}$ we have
$$\langle \mu+\nu, \Lambda_1\rangle=1.$$
Applying  the addition formula to   $\mu$ and $\nu$ in place of $\lambda$ and 
comparing the ribbon structure on the weights in $\{\Lambda_1, \mu, \nu,\Lambda_1+\mu, \Lambda_1+\nu\}$ leads to $q^{2\langle\mu, \Lambda_1\rangle}=(q')^{2\langle\mu, \Lambda_1\rangle}$ and
$q^{2\langle\nu, \Lambda_1\rangle}=(q')^{2\langle\nu, \Lambda_1\rangle}$, hence after term by term multiplication we get ${(q')}^2=q^2$, contradicting $q\neq0$.

Hence in all cases we may arrange $q'=q$. The relation $\theta_{\Lambda_1}=\theta'_{\Lambda_1}$ now implies that also the corresponding two $N$-roots $q^{1/N}$ and $q'^{1/N}$ are the same, and we thus have
a ribbon tensor equivalence ${\mathcal F}: {\mathcal C}\to{\mathcal C}'$ inducing $f$.
\end{proof}

    \bigskip

 \medskip

\section{Turning $C^*$-categories  with tensor structures into tensor $C^*$-categories, II}\label{16}

  Let $A$ be a discrete weak quasi bialgebra with a pre-$C^*$-algebra structure and let $\Omega\in M(A\otimes A)$
  be a given partially invertible operator with domain $\Delta(I)$. We develop a criterion that will be useful in Sect. \ref{18}, \ref{19}, \ref{20}  to verify the axioms
  of a positive $\Omega$-involution.

  Let  $\rho\in {\rm        Rep}_h(A)$ be   a $^*$-representation.
Since the coproduct is not coassociative in general,    there are different  tensor powers of $\rho$ each given order $n\geq3$, but they are all   equivalent. 
 \begin{defn}
 A representation $\rho$ is called {\it generating} if $\rho_n(a)=0$ for all $n$ implies $a=0$, where $\rho_n$ denotes
 the choice of an $n$-th tensor power of $\rho$.
 \end{defn}
 
It suffices to check the generating condition on  a choice of a $n$-th tensor power of $\rho$ for each $n$.

Let $\sigma$ and $\tau$ be f.d. $^*$-preserving representations of $A$ on Hlbert spaces.
As for the case of $\Omega$-involutive weak-quasi bialgebras, we may define the sesquilinear form
induced by $\Omega$ on the tensor product space and consider the  $\rho\otimes\sigma$ as a representation on this
space, except we do not know whether it is a Hilbert space $^*$-representation.
Let $\rho$ be a generating Hilbert space $^*$-representation. We may consider the full   subcategory ${\mathcal C}_\rho$ of ${\rm        Rep}(A)$ with objects  the various  tensor powers 
$\rho_n$ of $\rho$ on sesquilinear spaces.
 This is
a tensor category. We may determine the Hermitian form of $\rho_n$ with an inductive procedure, as follows.
Let $\Omega_n$ be the 
element of $A^{\otimes n}$ defining this form via $(\xi, \eta)=(\xi, \Omega_n\eta)_p$, where $(\xi, \eta)_p$ denotes
the untwisted $n$-th tensor power of the original Hermitian form of $\rho$ on $V_\rho^{\otimes n}$.
Let  $\Delta_n: A\to A^{\otimes n}$ denote the homomorphism defining the $A$-action on the space of $\rho_n$.
Writing
  $$\rho_n=\rho_r\underline{\otimes}  \rho_s,\quad{\rm        with}\quad r+s=n, \quad r, s< n,$$
 we have that
 $$\Omega_n=\Omega_r\otimes\Omega_s\Delta_r\otimes\Delta_s(\Omega),\quad\quad \Delta_{n}=\Delta_r\otimes\Delta_s\circ\Delta,$$
 where $\Omega_1=I$, $\Omega_2=\Omega$, $\Delta_1=1$, $\Delta_2=\Delta$.
 
 Assume for a moment that $\Omega$ is an $\Omega$-involution. Then we
   inductively get the following relations, extending (\ref{sa})--(\ref{intertwiner}).
 
 \begin{align}\label{}
&\Omega_n=\Omega_n^*,\\
&\Omega_{n}^{-1}\Omega_n=\Delta_n(I),\quad\quad\Omega_n\Omega_n^{-1}=\Delta_n(I)^*,\\
&\Delta_n(a)^*\Omega_n=\Omega_n\Delta_n(a^*),\quad\quad a\in A.  
 \end{align}   

 We next go back to the original situation, then we  only know that the above relations holds under the image
 of $\rho_n$ if we already know that $\rho_n$ is a $^*$-representation.
 
 \begin{thm}\label{Positivity} 
 Let $A$ be a discrete   pre-$C^*$-algebra  equipped with the structure of a weak quasi-bialgebra, and let $\rho$ be a generating $C^*$-representation of $A$. Let $\Omega\in M(A\otimes A)$ be a partially invertible element with domain $\Delta(I)$ and such that  for every irreducible $C^*$-representation $\sigma$,
\begin{equation}\label{p1}
\sigma\otimes\rho(\Omega), \quad\quad  \rho\otimes\sigma(\Omega),
\end{equation}
are positive on the full tensor product space,
that $\sigma\underline{\otimes}\rho$ and $\rho\underline{\otimes}\sigma$ are $C^*$-representations w.r.t. the $\Omega$-twisted inner product and that
 \begin{equation}\label{p2}
   \sigma\otimes\rho\otimes\rho(I\otimes \Omega1\otimes\Delta(\Omega)),\quad 
    \text{and} \quad \rho\otimes\rho\otimes\sigma(\Omega\otimes I\Delta\otimes1(\Omega))
    \end{equation} are positive as well.
   Moreover, assume that the associativity morphisms
 $$\sigma\otimes\rho\otimes\rho(\Phi), \quad   \quad \rho\otimes\rho\otimes\sigma(\Phi)$$
 are unitary with respect to the $\Omega$-twisted inner products.
 Then $\Omega$   is a positive element of $M(A\otimes A)$ and  in this way  $A$ becomes   a unitary discrete   weak  quasi bialgebra and $\Omega$ is uniquely determined by the operators $\sigma\otimes\rho(\Omega)$ for every irreducible $\sigma$. \end{thm}
 
\begin{proof} It follows from the first relation in (\ref{p1}) that
$\sigma\otimes\rho\otimes\rho(\Omega\otimes I\Delta\otimes 1(\Omega))$ is positive. It also follows
that $\Omega\otimes I\Delta\otimes 1(\Omega)$ and $ I\otimes\Omega1\otimes\Delta(\Omega)$ are positive on
$V_\rho\otimes V_\sigma\otimes V_\rho$ and  $V_\rho\otimes V_\rho\otimes V_\sigma$.
Every associativity morphism $\alpha_{\rho^r, \rho^s, \rho^t}=\rho^r\otimes\rho^s\otimes\rho^t(\Phi)$ of the full subcategory ${\mathcal C}_\rho$ of ${\rm Rep}(A)$ with objects parenthisized tensor powers 
 of $\rho$ can be written as a composition of tensor products with identity of morphisms of the form
 $\alpha_{\rho^r, \rho, \rho}$, $\alpha_{\rho, \rho^r, \rho}$, $\alpha_{\rho, \rho, \rho^r}$. By complete reducibility
 of representations and naturality, our assumptions imply unitarity of the first and the last, and the pentagon equation implies unitarity of the middle one. It follows that the
 associators imply that $\alpha_{\rho^r \rho^s \rho^t}$ are unitary.
We next show that every
  $\rho_n$   is a $C^*$-representation for the choice  iteratively defined by $\rho_{n+1}=\rho_n\underline{\otimes} \rho$.
  Assuming that a fixed $\rho_n$ is so, we decompose $\rho_n$ into pairwise orthogonal irreducible components $\sigma$. Since $V_\sigma\otimes V_\rho$ is invariant under $\sigma\otimes\rho(\Omega)$,   $\rho_n\otimes\rho(\Omega)$ is positive on $V_{\rho_n}\otimes V$ as well, hence it is a positive element of the $C^*$-algebra $\rho_n(A)\otimes \rho(A)$.   We may thus find an element ${\mathcal S}\in A\otimes A$ such that $\rho_n\otimes\rho({\mathcal S})$ is selfadjoint and
$\rho_n\otimes\rho(\Omega)=\rho_n\otimes\rho({\mathcal S})^2$.
On the other hand, the Hermitian form of $\rho_{n+1}$ is defined by the  action of the operator
$\rho^{\otimes n+1}[\Omega_{n+1}]$ on $V_\rho^{\otimes n+1}$ with $\Omega_{n+1}=\Omega_n\otimes I\Delta_n\otimes 1(\Omega)$. It follows that
$$\rho^{\otimes n+1}[\Omega_{n+1}]=\rho^{\otimes n}[\Omega_n]\otimes I\rho_n\otimes\rho(\Omega)=$$
$$\rho^{\otimes n}[\Omega_n]\otimes I\rho_n\otimes\rho(\S)^2=\rho^{\otimes n+1}[\Omega_n\otimes I\Delta_n\otimes 1(\S)^2]=  $$
$$\rho^{\otimes n+1}[\Delta_n\otimes 1(\S)^*\Omega_n\otimes I\Delta_n\otimes 1(\S)] = $$
$$\rho^{\otimes n+1}[\Delta_n\otimes 1(\S)]^*\rho^{\otimes n}[\Omega_n]\otimes I\rho^{\otimes n+1}[\Delta_n\otimes 1(\S)]    $$
and this is a positive operator by positivity of $\rho^{\otimes n}[\Omega_n]$. We   consider the $C^*$-representation 
$\tau=\oplus_n\rho_n$, which is faithful as $\rho$ is generating. We are left to show that  $\tau\otimes\tau[\Omega]$ is a positive operator in this representation, since it will then be    a positive element of $\tau(A)\otimes\tau(A)$, and therefore $\Omega$ positive in $A\otimes A$. To this aim, we   observe that the action of 
$\tau\otimes\tau[\Omega]$ on the subspace $V_{\rho_r}\otimes V_{\rho_s}$ is given by that of $\rho^{\otimes n}(\Omega'_n)$, where $n=r+s$ and $\Omega'_n=\Omega_r\otimes\Omega_s\Delta_r\otimes\Delta_s(\Omega)$. Thanks to  unitarity of the associativity morphisms and an inductive argument we see that $\rho^{\otimes n}(\Omega'_n)=\rho^{\otimes n}(\Phi_n\Omega_n\Phi_n^*)$ for suitable associativity morphisms $\Phi_n$.
It follows that $\tau\otimes\tau(\Omega)$ is positive, hence $\Omega$ is positive in $M(A\otimes A)$.
Therefore ${\mathcal C}_\rho$ is a unitary tensor category with unitary structure defined by $\Omega$.
Now the axioms of the $\Omega$-involution on $A$ follow.

\end{proof}

The above theorem will be useful in the construction of the main examples of  Sect. \ref{20}.

\begin{rem}  For example, if 
 $A$ is a finite dimensional $C^*$-algebra $A=\bigoplus_r M_{n_r}({\mathbb C})$ and $\rho$ is generating, every   $\rho_r\underline{\otimes}\rho$ is unitarily equivalent
 to an orthogonal direct sum of the projection $C^*$-representations $\rho_s: A\to M_{n_s}({\mathbb C})$ and their opposites
 $\rho_{-s}$ by Prop. \ref{complete}. By the previous theorem, verification of positivity of $\Omega$ reduces to the question of whether the negative forms $\rho_{-s}$ can be ruled out for this subclass of fusion tensor products.  \end{rem}

 We conclude the section with a further discussion on $C^*$-transportability.  In comparison with Sect.  \ref{12},
the following   discussion gives a direct method
to transport the tensor structure from ${\mathcal C}$ to ${\mathcal C}^+$ that will be useful in Sect. \ref{18}, \ref{19}, \ref{20}. We note however that this method is already implicit in our main results Theorem
\ref{transportability} and \ref{unitarizability}. Let us assume codition a).
It is not difficult to see,   using a quasi-inverse of ${\mathcal F}$, that when ${\mathcal C}$ has a weak dimension function
  there always is a  faithful weak quasi-tensor functor ${\mathcal G}:{\mathcal C}\to{\rm Hilb}$ such that ${\mathcal G}{\mathcal F}$ is a $^*$-functor. 
Recall   from 
 Remark \ref{non-unitarizability}   that  there are examples  for which the tensor structure
of ${\mathcal C}$ is not   transportable to ${\mathcal C}^+$ and in these cases we have a functor ${\mathcal G}$
which does not take the associativity morphisms to unitary morphisms.
On the other hand,
it follows from Theorem \ref{transportability} that
when the tensor structure of ${\mathcal C}$ is $C^*$-transportable to ${\mathcal C}^+$ then
we may find ${\mathcal G}$  taking the associativity morphisms to unitary morphisms.
The following proposition shows that the converse holds.

 \begin{prop}\label{unitarizability2} Let ${\mathcal F}:{\mathcal C}^+\to{\mathcal C}$ satisfy a) and assume that ${\mathcal C}$ admits a weak dimension function. Let ${\mathcal G}: {\mathcal C}\to{\rm Hilb}$ be a faithful functor such that ${\mathcal G}^+={\mathcal G}{\mathcal F}$ is a $^*$-functor and the morphisms
${\mathcal G}(\alpha_{\rho, \sigma, \tau})$ are unitary. Then every weak quasi-tensor structure on ${\mathcal G}$ 
induces the structure  of a tensor $C^*$-category on ${\mathcal C}^+$ s.t. ${\mathcal F}$ is a tensor equivalence ($C^*$-transportability).
 \end{prop}

 \begin{proof}
 Let  $(F, G)$  be a weak quasi-tensor structure for ${\mathcal G}$, thus $FG=1$ and also $G^*F^*=1$.   The functors ${\mathcal G}$, ${\mathcal G}^+$  correspond to the forgetful functors of a compatible triple as in Def. \ref{CompatibleTriple}. We consider the corresponding weak quasi bialgebra  $(A, \Delta, \Phi)$ with $A={\rm Nat}_0({\mathcal G})$. The linear equivalence
  ${\mathcal F}:{\mathcal C}^+\to{\mathcal C}$ induces an algebra isomorphism
   $A\to{A}^+={\rm Nat}_0({\mathcal G}^+)$, $\eta\to\eta_{{\mathcal F}(\ )}$. Since ${\mathcal G}^+$ is a $^*$-functor,     ${A}^+$ is a $C^*$-algebra. By Theorem \ref{transportability} we     only need to make ${A}^+$ into a unitary weak quasi-bialgebra. 
   We introduce  the structure  similarly to the case of the Tannakian theorems  \ref{TK_algebraic_quasi}, \ref{TheoremTannakaStar}.   We denote by
   $x$, $y$, $z,\dots$ the  irreducible representations of ${A}^+$ and define
   $$\tilde{F}_{x. y}: =F_{{\mathcal F}(x), {\mathcal F}(y)}: {\mathcal G}({\mathcal F}(x))\otimes {\mathcal G}({\mathcal F}(y))
   \to {\mathcal G}({\mathcal F}(x)\otimes {\mathcal F}(y))$$
   $$\tilde{G}_{x. y}:= G_{{\mathcal F}(x), {\mathcal F}(y)}:  {\mathcal G}({\mathcal F}(x)\otimes {\mathcal F}(y))\to{\mathcal G}({\mathcal F}(x))\otimes {\mathcal G}({\mathcal F}(y)).
 $$
 This suffices to make ${A}^+$ into a weak quasi-bialgebra
 by $$\tilde{\Delta}(\eta_{{\mathcal F}(\ )})_{x,y}=\tilde{G}_{x, y}\eta_{{\mathcal F}(x)\otimes 
 {\mathcal F}(y)}\tilde{F}_{x, y},$$
 $$\tilde{\Phi}_{x,y,z}=\Phi_{{\mathcal F}(x), {\mathcal F}(y), {\mathcal F}(z)}=$$
 $$1\otimes \tilde{G}_{y,z}\circ G_{{\mathcal F(x), {\mathcal F}(y)\otimes {\mathcal F}(z)}}\circ{\mathcal G}(\alpha_{{\mathcal F}(x), {\mathcal F}(y), {\mathcal F}(z)})\circ F_{{\mathcal F}(x)\otimes {\mathcal F}(y), {\mathcal F}(z)}\circ \tilde{F}_{x, y}\otimes 1.$$
 We introduce an $\Omega$-involution on $A^+$ by  $\Omega=\tilde{F}^*\tilde{F}$, $\Omega^{-1}=\tilde{G}\tilde{G}^*$. In this more general setting the only non-trivial
  verification is axiom (\ref{eqn:omega0}) which reduces to   $${\mathcal G}(\alpha_{{\mathcal F}(x), {\mathcal F}(y), {\mathcal F}(z)})^*={\mathcal G}(\alpha_{{\mathcal F}(x), {\mathcal F}(y), {\mathcal F}(z)}^{-1})$$ and holds by assumption.  
 \end{proof}

 \section{Coboundary categories  and Deligne's theorem}\label{17}

  By an interesting result of Deligne,
the study of   dimension in a braided tensor category can   be addressed in two equivalent ways: via    right duality   with  extra (pivotal/spherical) structure or   else via   extra structure on the braided symmetry (balancing/ribbon structure).

In this section we introduce a notion of symmetry that is more general than that of braided symmetry, and we call {\it generalised coboundary}.   It is a generalisation of both  the notion of braided symmetry and that of a coboundary due to Drinfeld that allows to study these symmetries in a
unified way.

The generalisation  is motivated by the fact that some of the structures that we study in this paper do not need the full notion of a braided symmetry, but only the more general class of symmetries,  which have   the advantage of  stability under  certain twist deformation.

A important source of coboundaries indeed arises from deformation of braided symmetries with ribbon structure and plays
a central role  in   the unitary structure of the weak quasi-Hopf algebras arising from quantum groups at roots of unity   studied   in Sects. \ref{18}, \ref{20}. We study   pivotal or spherical structures in   tensor categories with a generalised coboundary, and we extend Deligne result to this case. We start reviewing the notion of pivotal and spherical category.

 If ${\rho}^\vee$ is a two-sided dual of $\rho$ and if $(b_\rho, d_\rho)$ and $(b'_\rho, d'_\rho)$ respectively solve   the right and the left duality equations   for this pair, then
we can associate  two  functionals on the morphism space $(\rho,\rho)$, called left and right quantum traces, via  
 \begin{equation}
\label{traces1}{\rm        Tr}^L_\rho(T)=d_\rho\circ 1_{\rho^\vee}\otimes T\circ    b'_\rho 
\end{equation}
\begin{equation}\label{traces2}
{\rm        Tr}^R_\rho(T)=d'_\rho\circ T\otimes 1_{\rho^\vee}\circ   b_\rho. 
\end{equation}
If these solutions correspond to pivotal (or spherical) structures  a  well behaved theory  of dimension can be developed.   We briefly recall the main aspects,   dropping, for simplicity, the associativity morphisms in most of  our formulae in this section.

 Let $(\rho^\vee, b_\rho, d_\rho)$ be a right duality, see Sect. \ref{3}, and   $D:{\mathcal C}\to{\mathcal C}$   the associated   functor   as in (\ref{duality_functor}).
 Note that
 $D^2:{\mathcal C}\to{\mathcal C}$ is a covariant tensor functor. We assume from now on that    our category has two-sided duals,   so there is a   natural isomorphism $u$ from the identity functor $1$,  to $D^2$, which, however,     need not be monoidal. 
 An example of this occurrence arises 
 in the framework of  representations of semisimple weak quasi-Hopf algebras. The category has two-sided duals if    the square of the antipode $S$ is an inner automorphism. The natural isomorphism 
  is monoidal if the implementing element  can be chosen group-like, but this is not always the case.
On the other hand, 
 any natural  isomorphism $u\in(1, D^2)$ in a category with two-sided duals
 defines a left duality structure    coinciding with the right one on the objects via
\begin{equation}  b'_\rho= 1_{\rho^\vee}\otimes u_\rho^{-1}\circ  b_{\rho^\vee}, \quad\quad d'_\rho=d_{\rho^\vee}\circ u_\rho\otimes 1_{\rho^\vee}.\label{left_duality}\end{equation} Furthermore any pair of right and  left dualities $(\rho^\vee, b_\rho, d_\rho)$ and $({}^\vee\rho, b'_\rho, d'_\rho)$ with   $\rho^\vee={}^\vee\rho$  is of this form with $u$ uniquely determined.   Indeed,  the morphism 
\begin{equation}u_\rho:=d'_\rho\otimes 1_{\rho^{\vee\vee}}\circ 1_\rho\otimes b_{\rho^\vee}\label{natural}\end{equation} is a natural isomorphism in $(1, D^2)$ with $u_\rho^{-1}=d_{\rho^\vee}\otimes 1_\rho\circ 1_{\rho^{\vee\vee}}\otimes b'_\rho$ and the two constructions are inverse of one another. Given $u\in(1, D^2)$, any other $\omega\in(1, D^2)$ can be written in the form
$\omega=uv^{-1}$, with $v\in (1, 1)$ uniquely determined. (The use of the inverse of $v$ matches our notation   in the framework of quantum groups, and originates with the convention in  \cite{Turaev}).  Correspondingly, any other left duality is of the form
\begin{equation}\label{deformation} \tilde{b}_\rho= 1_{\rho^\vee}\otimes v_\rho\circ  b'_{\rho}, \quad\quad \tilde{d}_\rho=d'_{\rho}\circ v^{-1}_\rho\otimes 1_{\rho^\vee}. \end{equation}
   A {\it pivotal structure} on ${\mathcal C}$ is the datum of a right duality functor $D$ together with a {\it monoidal} isomorphism $\omega\in (1, D^2)$  \cite{Freyd_Yetter}.  
In a tensor category with right duality $(b_\rho, d_\rho)$    the   left duality defined by a pivotal structure $\omega$ in place of $u$
   in 
(\ref{left_duality}) will be denoted as  $(b^\ell_\rho, d^\ell_\rho)$.
A pair of dualities $(b_\rho, d_\rho)$ and  $(b^\ell_\rho, d^\ell_\rho)$  so related  induces
    ${\mathbb C}$-valued left and right  quantum traces (\ref{traces1}), (\ref{traces2}) 
   which are multiplicative on tensor product morphisms. 
  A {\it spherical structure} on ${\mathcal C}$ is a pivotal structure such that the associated left and right traces coincide.
    In this case we shall simply write   ${\rm        Tr}_\rho$.
   A {\it spherical category} is a tensor category endowed with  a spherical structure.
In a spherical   category
  ${\rm        Tr}_\rho(ST)={\rm        Tr}_\sigma(TS)$, for any pair of morphisms $T\in(\rho,\sigma)$, $S\in(\sigma, \rho)$ the   {\it categorical (or quantum) dimension} $\rho\to d(\rho)$ is defined by
 $$d(\rho)={\rm        Tr}_\rho(1_\rho).$$  
It  is  additive, multiplicative and, for categories over ${\mathbb C}$ as those of this paper, it takes real values on the objects, see \cite{BW} and Sect. 2 in   \cite{ENO} for more details. 
It is not known whether a fusion category always admits a pivotal structure, but see \cite{ENO, Mueger2,   Mueger0} for results and references.

\begin{defn}
 A {\it generalised coboundary} is a natural isomorphism $c(\rho, \sigma)\in(\rho\otimes\sigma, \sigma\otimes\rho)$ satisfying (\ref{normalization_symmetry}) and such that the following diagram commutes.
\begin{equation}\label{generalised_coboundary} \begin{CD}
 (\rho\otimes \sigma)\otimes \tau @>{{\alpha}}>>  \rho\otimes(\sigma\otimes \tau)@>{1\otimes c}>> \rho\otimes (\tau\otimes\sigma) \\
@V{c\otimes 1}VV @. @VV{c}V \\
 (\sigma\otimes\rho)\otimes\tau@>{c}>>  \tau\otimes(\sigma\otimes\rho)@<{\alpha}<<  (\tau\otimes\sigma)\otimes\rho \end{CD}\end{equation}
 \end{defn}

If $c(\rho, \sigma)$ is a generalised coboundary then $c'(\rho, \sigma):=c(\sigma, \rho)^{-1}$ is too.

\begin{ex} 
A generalised coboundary for which
  $c$ satisfies the symmetry condition $c^2=1$ is  a coboundary in the sense introduced by Drinfeld \cite{Drinfeld_quasi_hopf}. 
  \end{ex}

\begin{rem}    
 Every braided symmetry is a generalised coboundary. Indeed, if   $c$ is such a symmetry,   we may use  the hexagonal 
 equations  
 (\ref{braided_symmetry1}), (\ref{braided_symmetry2})
   and verification of commutativity of (\ref{generalised_coboundary}) reduces to   the Yang-Baxter relation, which   follows from the braided symmetry axioms, see e.g. Prop. 8.1.10 in  \cite{EGNO}.    \end{rem}
 
 The following statement explains the notion of generalised coboundary in an
 important class of tensor categories. 
 
 \begin{prop}\label{constructing_generalized_coboundaries} Let $A$ be a weak quasi bialgebra and $Q\in A\otimes A$ a twist such that $A_Q=A^{{\rm op}}$.
 Then $c(\rho, \sigma):=\Sigma\rho\otimes\sigma(Q)$ is a generalised coboundary of ${\rm Rep}(A)$.
\end{prop}

We refer to (\ref{qtriangular1}), (\ref{qtriangular2}),   (\ref{YB}),  (\ref{prop1}), with $Q$ in place of $R$, for an explicit form of the equality $A_Q=A^{{\rm op}}$. 

\begin{rem} The construction of generalised cobounderies on ${\rm Rep}(A)$ of \ref{constructing_generalized_coboundaries} extends to the case where $A$ is a discrete weak quasi bialgebra, and the twist between $A$ and $A^{\rm op}$ satisfies $Q\in M(A\otimes A)$. In this case, we also see that all generalised coboundaries of ${\rm Rep}(A)$ are of this form, via Tannaka-Krein duality.
\end{rem}

We introduce twist deformation of generalised coboundaries.

 \begin{prop} 
  Let $c$ be a generalised coboundary and $\eta\in(1, 1)$  a natural isomorphism of the identity functor with $\eta_\iota=1_\iota$.  Then 
 $c^\eta(\rho, \sigma):=c(\rho, \sigma)\circ \eta_\rho^{-1}\otimes\eta_\sigma^{-1}\circ\eta_{\rho\otimes\sigma}$ is a generalised coboundary as well.  
\end{prop}

If  $c$ is a braided symmetry, $c^\eta$ may fail to be a braided symmetry, but it is   a generalised coboundary.

\begin{prop}\label{left_duality_coboundary}
Let ${\mathcal C}$ be a tensor category with generalised coboundary $c$. Then duals are two-sided.  If  the category has right duals and $(\rho^\vee, b_\rho, d_\rho)$ denotes a  right duality then    \begin{equation}\label{left1}
b'_\rho=c(\rho^\vee, \rho)^{-1}\circ b_\rho, \quad\quad d'_\rho=d_\rho\circ c(\rho,\rho^\vee)\end{equation} 
is a left duality with $\rho^\vee={}^\vee\rho$. Conversely, given a left duality, $({}^\vee\rho, b'_\rho, d'_\rho)$ the same formula defines a right duality.
 \end{prop}
\begin{proof}
 The left duality relations for $b'_\rho$ and $d'_\rho$ follow from a computation that uses, in order,  commutativity of the diagram (\ref{generalised_coboundary}),   naturality   of $c$ and the right duality equations for $b_\rho$ and $d_\rho$.
\end{proof}

\begin{rem}
Note that we may apply the same construction to   $c'$ and get another left duality with ${}^\vee\rho=\rho^\vee$,
\begin{equation}\label{left2}b''_\rho=c(\rho, \rho^\vee)\circ b_\rho, \quad\quad d''_\rho=d_\rho\circ c(\rho^\vee,\rho)^{-1}.\end{equation} 
In the special case where $c$ is a genuine coboundary, these two left dualities coincide, thus every 
 right duality $(b_\rho, d_\rho)$ 
has a canonically associated left duality in this way.  It is also easy to see that   the associated right and left traces coincide thanks to naturality of $c$.
We do not know whether this pair of dualities   corresponds to a pivotal structure for all coboundaries, but this is known to be the case when $c$ is a permutation symmetry or for all the examples of coboundaries that may be constructed from
braided symmetries and twist deformation.  
\end{rem}

\begin{defn}\label{balancing_ribbon_tensor_category}  Let   ${\mathcal C}$ be a tensor category and let $a_{\rho, \sigma}\in(\rho\otimes\sigma, \rho\otimes\sigma)$ be a tensor structure for the identity functor.

\noindent a) A {\it balancing structure} for $a$ is a natural isomorphism $v\in(1, 1)$ making the identity functor $1$ with tensor structure $a$   monoidally isomorphic to the trivial tensor structure of $1$, so
\begin{equation}\label{balancing_structure}
a_{\rho, \sigma}=v_\rho\otimes v_\sigma\circ v_{\rho\otimes\sigma}^{-1}.\end{equation}
b)  If  ${\mathcal C}$ has a  right duality $(\rho^\vee, b_\rho, d_\rho)$, a {\it ribbon structure} for $a$ is a balancing structure   compatible with duality, see Def. \ref{compatibility_with_duality}.

\end{defn}

\begin{rem} 
 If $v$ is a balancing structure for $a$,   the relation $v_\iota=1_\iota$ easily follows from the  fact that we are assuming that $\iota$ is a strict unit, but for   general categories it needs to be part of the axioms. 
\end{rem}

 We next see that the question of whether a rigid tensor category with a  generalised coboundary admits a pivotal
or spherical structure can be reduced to the analysis of  two   tensorial structures of the identity functor, which are naturally associated to the coboundary.  In the case where $c$ is a braided symmetry, these reduce to the same structure,
but they may be distinct in general.
 We first   generate tensor structures of $1$ from $c$.

\begin{defn}
Let $c$ be a generalised coboundary and $(b_\rho, d_\rho)$ a right duality. Consider the left duality $(b'_\rho, d'_\rho)$
described in (\ref{left1}). The   natural isomorphism    $u\in(1, D^2)$ associated  to this pair as in
(\ref{left_duality}), (\ref{natural}) is called  {\it Drinfeld isomorphism}. \end{defn}

Hence Drinfeld isomorphism is the composite 
\begin{equation}\label{Drinfeld_isomorphism}u_\rho: \rho\xrightarrow{1\otimes b_{\rho^\vee}}\rho\otimes\rho^\vee\otimes\rho^{\vee\vee}\xrightarrow{c\otimes 1}\rho^\vee\otimes\rho\otimes\rho^{\vee\vee}\xrightarrow{d_\rho\otimes 1} \rho^{\vee\vee}.\end{equation}

\begin{prop} Let $c$  be  generalised coboundary. The isomorphisms 
$$c^2(\rho, \sigma):=c(\sigma, \rho)\circ c(\rho, \sigma)\in(\rho\otimes\sigma, \rho\otimes\sigma)$$
define a tensor structure on the identity functor.
\end{prop}

\begin{proof} Naturality of $c^2$ in the two variables is obvious. The tensor structure axiom
$$c^2(\rho, \sigma\tau)\circ 1_\rho\otimes c^2(\sigma,\tau)=c^2(\rho\sigma,\tau)\circ c^2(\rho, \sigma)\otimes 1_\tau$$
is indeed a simple consequence of the generalised  coboundary axioms for $c$.
\end{proof}

 We start with a condition leading to the construction of two coinciding    quantum traces.

\begin{thm}\label{spherical}
Let $c$ be a generalised coboundary, $(b_\rho, d_\rho)$ a right duality and $v\in(1,1)$ a ribbon structure for $c^2$.
Then the left and right quantum traces corresponding to a given right duality $(b_\rho, d_\rho)$ and to the associated left duality via $\omega:=uv^{-1}\in(1, D^2)$,  coincide.
\end{thm}

\begin{proof}
The left duality defined by $\omega=uv^{-1}$ is given by  
 \begin{equation}\label{ell_duality} \tilde{b}_\rho=1_{\rho^\vee}\otimes v_\rho\circ c(\rho^\vee, \rho)^{-1}\circ b_\rho,\quad\quad \tilde{d}_\rho=d_\rho\circ c(\rho, \rho^\vee)\circ v_\rho^{-1}\otimes 1_{\rho^\vee}.\end{equation}
 The corresponding right trace is given by
$${\rm        Tr}^R_\rho(T)= d_\rho\circ 1_{\rho^\vee}\otimes T\circ c(\rho,\rho^\vee)\circ v_\rho^{-1}\otimes 1_{\rho^\vee}\circ b_\rho.
$$ 
To compare it with the left trace we compute
 $${\rm        Tr}^L_\rho(T)=d_\rho\circ 1_{\rho^\vee}\otimes T\circ1_{\rho^\vee}\otimes v_{\rho}\circ c(\rho^\vee, \rho)^{-1} \circ b_\rho=$$
$$ d_\rho\circ 1_{\rho^\vee}\otimes T\circ c(\rho^\vee, \rho)^{-1}\circ v_{\rho}\otimes 1_{\rho^\vee} \circ b_\rho=$$
 $$d_\rho\circ 1_{\rho^\vee}\otimes T\circ c(\rho^\vee, \rho)^{-1}\circ 1_\rho\otimes v_{\rho^\vee}\circ b_\rho,$$
  the last equality follows from  
 $v_\rho\otimes 1_{\rho^\vee}\circ b_\rho=1_\rho\otimes v_{\rho^\vee}\circ b_\rho$ in turn  due to compatibility of $v$ with duality.  On the other hand,
 $c(\rho^\vee, \rho)^{-1}\circ 1_\rho\otimes v_{\rho^\vee}=c(\rho,\rho^\vee)\circ v_\rho^{-1}\otimes 1_{\rho^\vee}\circ v_{\rho\otimes \rho^\vee}$  thanks to the balancing condition
 $c^2(\rho, \sigma)=v_\rho\otimes v_\sigma\circ v_{\rho\otimes\sigma}^{-1}$.
   The conclusion now follows     from this and naturality of $v$.
 
\end{proof}

We have yet another tensor structure of the identity functor induced by $c$ as follows. Let 
${d_2}({\rho, \sigma}):\rho^{\vee\vee}\otimes\sigma^{\vee\vee} \to(\rho\otimes\sigma)^{\vee\vee}$ 
denote the natural tensor structure of $D^2$. In the framework of weak quasi-Hopf algebras, we have explicitly computed
the element of $A\otimes A$ inducing $d_2$, see the discussion following Prop. \ref{rigidity} and Def. \ref{pivotal_wqh}. 
       We can   equip $1$ with the  new tensor structure, denoted $c_2$,   obtained pulling back the tensorial structure of $D^2$ via Drinfeld isomorphism.
In other words, we let ${c_2}({\rho, \sigma})\in(\rho\otimes\sigma, \rho\otimes\sigma)$ denote the isomorphisms  making the following diagram commute,  \begin{equation}\begin{CD}\label{c_2}
 \rho\otimes \sigma  @>{u_\rho\otimes u_\sigma}>>  \rho^{\vee\vee}\otimes\sigma^{\vee\vee}  \\
@V{c_2}VV @VV{d_2}V \\
 \rho\otimes\sigma @>{u_{\rho\otimes\sigma}}>>  (\rho\otimes\sigma)^{\vee\vee}
 \end{CD}\end{equation}

 We next   analyse dependence of Drinfeld isomorphism and $c_2$ on the  right duality.

\begin{lemma}
Let $(\rho^\vee, b_\rho, d_\rho)$ and $(\tilde{\rho}, \tilde{b}_\rho,\tilde{d}_\rho)$ be two right dualities with associated functors $D$ and $\tilde{D}$ respectively, and let $\xi\in(\tilde{D}, D)$ the corresponding monoidal isomorphism.
Let $u$ and $\tilde{u}$ be corresponding Drinfeld isomorphisms defined by the same generalised coboundary.
Then $$\tilde{u}_\rho=\zeta_\rho\circ u_\rho$$
where
  $\zeta_\rho:=\xi_{\tilde{\rho}}^{-1}\circ \xi_\rho^\vee: D^2\to\tilde{D}^2$ is the composite monoidal isomorphism. 
\end{lemma}

\begin{proof}
The proof follows from a computation starting from $\tilde{u}_\rho$ taking into account $b_{\tilde{\rho}}=\xi_\rho^{-1}\otimes\xi_\rho^\vee\circ b_{\rho^\vee}$,   $\tilde{b}_{\tilde{\rho}}= 1_{\tilde{\rho}}\otimes\xi_{\tilde{\rho}}^{-1}\circ b_{\tilde{\rho}}$,   $\tilde{d}_\rho= d_\rho\circ\xi_\rho\otimes 1_\rho$ and naturality of $c$.
\end{proof}

\begin{prop}
Let ${\mathcal C}$ be a tensor category with generalised coboundary $c$ and right duality $(\rho^\vee, b_\rho, d_\rho)$. Then the isomorphism $$c_2(\rho, \sigma)\in (\rho\otimes\sigma, \rho\otimes\sigma)$$ is a tensor structure of the identity functor which does not depend on the choice of the right duality.  
\end{prop} 

\begin{rem} a) It is known that $c_2=c^2$  if $c(\rho, \sigma)$ is a braided symmetry, for   a proof  in a  strict tensor category see \cite{EGNO},    Prop. 8.9.3.  
b) In Prop. \ref{ribbon_element} we have explicitly shown that $c_2=c^2$ for the braided symmetry associated to the     quasitriangular structure of any weak  Hopf algebra. 

We give an example showing that  $c_2$ and $c^2$ may be different tensor structures.

\begin{ex}\label{Example} 
 Consider 
the tensor category ${\mathcal C}={\rm Vec}_G$ of finite dimensional $G$-graded vector spaces over a finite abelian group $G$,
with tensor product defined in the  standard way, for $V=(V_g)$ and $W=(W_h)$, $(V\otimes W)_{k}=\oplus_{gh=k}V_g\otimes V_h$, and natural associator, see \cite{EGNO}. Then every group element $g$ defines a $1$-dimensional space $\delta_g$ of grade $g$ and these are all the   irreducible objects up to equivalence. We have that $\delta_g^{-1}$ is both a right and left dual of $\delta_g$ and duality equations are solved by the identity maps. A generalised coboundary
$c$ is determined by the action on $\delta_g\otimes\delta_h$, and this gives a complex-valued nonzero function $c(g, h)$ on two variables. The coboundary relation corresponds to requiring that $c(g, h)$ be a two-cocycle: $c(g, h)c(gh, k)=c(h, k)c(g, hk)$ with   $c(1, g)=c(g, 1)=1$. Drinfeld isomorphism $u_g$   acts as $c(g, g^{-1})$ on $\delta_g$, while $d_2$ acts trivially. It follows that $c_2(g, h)=c(g, g^{-1})c(h, h^{-1})c(gh, (gh)^{-1})^{-1}$ while $c^2(g, h)=c(g, h)c(h, g)$.
A computation shows that $c_2=c^2$ if and only if $c(h, h^{-1})=c(h, g)c(h, (gh)^{-1})$, and it is easy to see that  this is not always the case   for a normalised cohomologically trivial $c(g, h)=\mu(gh)\mu(g)^{-1}\mu(h)^{-1}$ for $G={\mathbb Z}_3$.
 
 \end{ex}

The following extends Deligne's result to generalised  coboundaries. 

\begin{thm}\label{Deligne} Let ${\mathcal C}$ be a tensor category with generalised coboundary $c$  and right duality $(\rho^\vee, b_\rho, d_\rho)$. There is a bijective correspondence between pivotal structures $\omega\in(1, D^2)$ and balancing structures 
$z\in (1,1)$ for $c_2$
    given by $$\omega=u z^{-1},$$ where $u\in (1, D^2)$ is   Drinfeld isomorphism associated to $c$. 

\end{thm}

\begin{proof}
The map $z\to \omega=uz^{-1}$ is a bijective correspondence between isomorphisms  $\omega\in(1, D^2)$ and $z\in(1,1)$, furthermore $\omega$ is monoidal precisely when $z$ is a balancing for $c$, by  commutativity of (\ref{c_2}).
\end{proof}

We derive a sufficient condition for existence of spherical structures.

\begin{cor}\label{spherical2} Let ${\mathcal C}$ be a tensor category with right duality $(b_\rho, d_\rho)$, generalised coboundary $c$
satisfying 
\begin{equation}\label{c_2=c^2} c_2(\rho,\sigma)=c^2(\rho, \sigma).\end{equation}
and ribbon structure $v$. Then the pivotal structure 
 $\omega=uv^{-1}$ is  spherical. 
 The corresponding left duality is given by
 \begin{equation}
b^\ell_\rho=1_{\rho^\vee}\otimes v_\rho\circ c(\rho^\vee, \rho)^{-1}\circ b_\rho,\quad\quad d^\ell_\rho=d_\rho\circ c(\rho, \rho^\vee)\circ v_{\rho}^{-1}\otimes 1_{\rho^\vee}
\end{equation} 
\end{cor} 

\begin{proof}
This is a consequence of Theorem. \ref{spherical} and Prop. \ref{Deligne}. The left duality equations follow from (\ref{deformation}), (\ref{left1}).
\end{proof}

In particular,   the   quantum dimension
is given by 
\begin{equation}\label{quantum_dimension_pivotal}
d(\rho)=d_\rho\circ 1_{\rho^\vee}\otimes v_\rho\circ c(\rho^\vee, \rho)^{-1}\circ b_\rho=d_\rho\circ c(\rho, \rho^\vee)\circ v_{\rho}^{-1}\otimes 1_{\rho^\vee}\circ b_\rho.
\end{equation}

 Theorem.  \ref{spherical} recovers corresponding results known for ribbon categories \cite{Turaev}. Note that
 Cor. \ref{spherical2}
is of little use in the  case where (\ref{c_2=c^2}) does not hold. Indeed   in the Example
 \ref{Example} as we may choose for $c$ the unique permutation symmetry, so $c_2=c^2=1$,   gives that the associated Drinfeld isomorphism $u_\rho=1$ is a spherical structure.

\end{rem}

We next discuss   properties of twisted generalised coboundaries.
  Let $(b_\rho, d_\rho)$ be  a fixed right duality, $c$   a generalised coboundary,    $u$ the associated Drinfeld isomorphism and $c_2\in(\rho\otimes\sigma, \rho\otimes\sigma)$ natural isomorphism 
  as in (\ref{c_2}).
Let $\eta\in(1,1)$ be a natural isomorphism, and $c^\eta$ the twisted coboundary. The corresponding isomorphisms   will be denoted respectively by $u^\eta$ and $c^\eta_2$.

 \begin{prop} Let $\eta\in(1,1)$ be a compatible with duality. We have
 \begin{itemize}
\item[{\rm        a)}] 
 $u^\eta=u\circ\eta^{-2},$
  \item[{\rm        b)}]  $c_2^\eta(\rho,\sigma)=\eta_{\rho\otimes\sigma}^{2} \circ c_2(\rho,\sigma)\circ \eta_\rho^{-2}\otimes\eta_\sigma^{-2},$ $((c^\eta)^2(\rho,\sigma)=\eta_{\rho\otimes\sigma}^{2} \circ c^2(\rho,\sigma)\circ \eta_\rho^{-2}\otimes\eta_\sigma^{-2},)$
  \item[{\rm        c)}]  if  $v$ is a balancing (ribbon) structure for $c_2$ $(c^2)$ then $v^\eta:=v\circ \eta^{-2}$ is a balancing (ribbon) structure for $c^\eta_2$ $((c^\eta)^2)$
    \item[{\rm        d)}] $v$ and $v^\eta$ correspond to the same pivotal structure under the map described in Prop. \ref{Deligne}, and therefore to the same left duality and quantum traces,
 \item[{\rm        e)}] if $c$ satisfies (\ref{c_2=c^2}) then so does  $c^\eta$. 
 \end{itemize}
 \end{prop} 
 
\begin{proof} a) The proof follows from a   computation starting from (\ref{Drinfeld_isomorphism}), with  $c$ replaced by $c^\eta$, where we use naturality and compatibility with duality of $\eta$ and the fact that the right duality functor (\ref{duality_functor}) 
 can equivalently be defined by $d_\sigma\circ 1_{\sigma^\vee}\otimes T=d_\rho\circ T^\vee\otimes 1_\rho$. The remaining statements follow   from one another. 
 \end{proof}

We describe a   twisting   making a generalised coboundary with a balancing structure  into a genuine coboundary and Drinfeld isomorphism into a monoidal isomorphism from the identity tensor functor. This twisting first appeared in the work of Drinfeld  \cite{Drinfeld_quasi_hopf} in the framework of quantised universal Hopf algebras. As it turns out,
the associated spherical structure is the same as that arising in the framework of ribbon categories.

\begin{thm}\label{coboundary_associated_to_a_ribbon_structure} In a tensor category with right duality, let
  $c$ be a generalised coboundary satisfying (\ref{c_2=c^2}) (e.g. a braided symmetry) with balancing structure $v$, and let $w\in(1,1)$ be a natural isomorphism compatible with duality such that $w^2=v$. Then $c^w$ is a coboundary, $c_2^w=1\otimes 1$, $v^w=1$,  and $u^w$ is a spherical structure coinciding with that defined by $c$ and  $v$ as in Prop. \ref{Deligne}.
\end{thm}

The construction of   $c^w$ is the analogue of Drinfeld construction of unitarized $\overline{R}$-matrix in a ribbon Hopf algebra. In Sect.  \ref{18} we shall study the relation with $\Omega$-involution.

\section{Hermitian coboundary wqh   and relation with Hermitian ribbon wqh}\label{18}

   In this section we  introduce  the   notion of   {\it Hermitian   coboundary weak quasi-Hopf algebra}.
   Essentially, we understand these as   as ribbon weak quasi-Hopf algebras endowed with  a $^*$-algebra
   structure satisfying various compatibility relations between the $^*$-involution, the coproduct and ribbon structure. We are mainly interested in the case
   of  discrete algebras with a pre-$C^*$-algebra structure.  
   
The most relevant structural aspect of our definition is the relation between coproduct and $^*$-involution.   
   Informally, this relation may   be interpreted as an    antimultiplicativity property of the
involution     on the `dual noncommutative function algebra', that is 
$(AB)^*=B^*A^*.$ When we take the adjoint on both sides, we get an equation that dually identifies
the opposite coproduct
$\Delta^{\rm op}$ and the adjoint coproduct $\tilde{\Delta}$.
To be more precise,  we require that
$\Delta^{\rm op}$ and  $\tilde{\Delta}$ (together with all the remaining structural data)
 are related by a trivial twist. Moreover, since we have an $R$-matrix which relates $\Delta^{\rm op}$ and $\Delta$,
 we may interpret that noncommutativity  arises explicitly from the $R$-matrix as is familiar in quantum group theory.
  This property makes these algebras rather different from the ordinary Hopf $^*$-algebras, where coproduct and $^*$-involution commute.

     Among other axioms  we assume a  relation
 involving directly the unitary structure with the braiding, or more precisely 
 with the coboundary symmetry in the representation category.
We assume the existence of a square root of the ribbon structure. Thus we have an associated coboundary in the representation category in the sense of Sect. \ref{17}. It follows from the axioms that there is an   $\Omega$-involution   on the algebra  in the sense of Sect. \ref{8} associated to the braiding data.
When the $\Omega$-involution of an
 Hermitian   coboundary weak quasi-Hopf algebra
 it  is unitary, we shall talk of a  {\it unitary coboundary weak quasi-Hopf algebra}.

 In this section we study the main properties. For example, among general $\Omega$-involutions, those associated
 to a coboundary always make the braiding unitary, see Theorem \ref{unitary_ribbon2}.

  Moreover, we shall give a characterization of the case where
    an Hermitian coboundary weak quasi-Hopf algebra gives rise to an Hermitian ribbon category, Theorem
    \ref{Hermitian_category2}.

The main result of this section is a  Tannakian   characterization
    of Hermitian coboundary weak quasi-Hopf algebras, see Theorem \ref{TK_unitary_ribbon}.
    This characterization describes such algebras as categories endowed with a faithful functor to ${\rm Herm}$ with
    a weak quasi-tensor structure $(F, G)$ and 
      compatibility equations between the
    coboundary of the  category, the permutation symmetry of ${\rm Herm}$ and $(F, G)$.
The simplest case is that of symmetric tensor categories, and 
    the Tannakian characterization becomes the notion of symmetric tensor functor.  In particular
    compact groups is a natural class of examples, and   we are  in the setting
    of   the Doplicher-Roberts duality theorem \cite{DR1}.
    More generally, the permutation symmetry is replaced by the coboundary of Drinfeld in the sense of Sect.
    \ref{17}.
    
     In the next section we discuss
  a possibly proper subclass of
 Hermitian coboundary weak quasi-Hopf algebras and we shall develop a criterion to construct such algebras.

We shall show in the next section that the  unitarization of a unitary
     coboundary weak quasi-Hopf algebra in this subclass,   is again an algebra of this kind with the advantage that both the unitary structure and the $R$-matrix take a  
   simpler form, see
   Remark \ref{properties_of_trivial_unitary_ribbon2} completely determined by the square root of the ribbon structure.
It reminds  the form taken by
   Drinfeld $R$-matrix of the quasi-Hopf algebra associated to Knizhnik-Zamolodchikov
    differential equations in
    \cite{Drinfeld_quasi_hopf}.
   It seems  valuable to us that this simple $R$-matrix may be derived in a general setting by the study of unitary structures of ribbon weak quasi-Hopf algebras.
   We hope to further develop   this  
   study in   future updates of this paper

   Our interest in discrete algebras is motivated by the unitary structure
of the fusion categories  ${\mathcal C}({\mathfrak g}, q, \ell)$ associated to $U_q({\mathfrak g})$ at certain roots of unity. Kirillov 
 defined a tensor $^*$-category tensor equivalent to ${\mathcal C}({\mathfrak g}, q, \ell)$  and conjectured that these   where unitary. The conjecture was shown to be true by Wenzl and Xu      \cite{Kirillov}, \cite{Wenzl}, \cite{Xu}.
 We may   regard our notion as an abstract version of Kirillov   $^*$-structure  
    following  the approach of Wenzl in \cite{Wenzl}. 
 We shall   recall these results in     Sect. \ref{20} and we recall in particular that the main example
 of Hermitian coboundary weak quasi-Hopf algebra is $U_q({\mathfrak g})$ itself for $|q|=1$, although not a semisimple example at roots of unity.

   Furthermore   in Sect. \ref{20} we shall construct f.d. unitary coboundary
weak  Hopf algebras as suitable quotients of
$U_q({\mathfrak g})$ with representation category equivalent to ${\mathcal C}({\mathfrak g}, q, \ell)$.

   \medskip

Recall that for a general weak quasi-bialgebra   $A$ we have defined a twisted algebra $A_F$, see Prop. \ref{twisted_wqh}, the opposite algebra $A^{{\rm        op}}$, see
(\ref{opp}) and furthermore, if $A$ is also a $^*$-algebra, we have introduced the adjoint algebra $\tilde{A}$ in (\ref{tilde}).  Note that $A_F$, $A^{{\rm        op}}$, and  $\tilde{A}$   have quasitriangular structures naturally induced by one of $A$, see Prop. \ref{R-canonicity}. Moreover, $A^{{\rm        op}}$ and  $\tilde{A}$ have a strong antipode if so does $A$, and similarly for $A_F$ if (\ref{strong_antipode_twist_eq}) holds.
In particular, $A^{{\rm        op}}$ and  $\tilde{A}$ are    weak  bialgebras if so is $A$, and similarly for $A_F$ if $F$ is a $2$-cocycle.

\begin{defn}\label{Hermitian_ribbon_wqh} A {\it Hermitian coboundary} weak quasi-Hopf   algebra $A$ is defined by the following data:
\begin{itemize}
\item[{\rm        a)}]    A   weak quasi-Hopf  algebra $A$ endowed with a $^*$-algebra
involution 
  with an antipode $(S, \alpha, \beta)$ 
\item[{\rm        b)}]   a
ribbon structure $(R, v)$ for $A$ associated to $(S, \alpha, \beta)$ (see Def. \ref{balanced_ribbon_wqh})
 such that the ribbon element $v\in A$
  is   unitary,  
\item[{\rm        c)}]     a unitary central square root $w\in A$
  of $v$ such that $\varepsilon(w)=1$, 
$S(w)=w$, 
\item[{\rm        d)}] $\tilde{A}=(A^{{\rm        op}})_E$ as quasitriangular weak quasi-bialgebras,
where $E=\Delta(I)^*\Delta^{{\rm        op}}(I)$ is a trivial twist, that  is
  $E^{-1}=\Delta^{{\rm        op}}(I)\Delta(I)^*$.  
  \end{itemize}
  
\end{defn}

\begin{rem}
Our axioms are motivated by the structure of $U_q({\mathfrak g})$ for $|q|=1$ that 
 will be important to us, and we shall
 recalled it
in Sect. \ref{20},  Theorem
\ref{U_q_as_a_Hermitian_ribbon_h}.
Notice   however that, since  
 the $R$-matrix and ribbon structure lie
in a suitable topological completion of $U_q({\mathfrak g})\otimes U_q({\mathfrak g})$ \cite{Sawin},   this 
algebra can not be included as   an example, unless we weaken our axioms. 
However      we shall refrain  from doing this. To deal with examples where the ribbon structure
lies in a larger algebra,   we shall content to consider the case of discrete algebras.

\end{rem}

\begin{defn} A {\it discrete Hermitian  coboundary} weak quasi-Hopf algebra is defined by a discrete weak quasi-Hopf algebra $A$ endowed with data $(^*, R, v, w)$ such that axioms a)-d) hold as before for discrete
algebras, that is  the $^*$-involution makes $A$ into a pre-$C^*$-algebra, and the ribbon and coboundary structure satisfies
 $R\in M(A\otimes A)$, $v$,
$w\in M(A)$.
\end{defn}

A  (discrete) weak  Hopf algebra  $A$ satisfying axioms a)--d) will be called a (discrete) Hermitian coboundary   weak  Hopf algebra.

\begin{rem}
a) Note that the definitions do not depend on the choice of the antipode by Prop. 
\ref{unique_antipode}.
Furthermore, when $A$ is discrete an antipode   may always be chosen  commuting with $^*$ by Remark
\ref{antipode_commuting_with_*}.  In the rest of the section   antipodes $(S, \alpha, \beta)$
will be chosen  with $S$ commuting with $^*$ for discrete
algebras. These antipodes are of the form $xS(\ )x^{-1}, x\alpha, \beta x^{-1})$ with $x$ unitary 
and uniquely determined. 
b) The   equality required in d) between the $R$-matrices of $\tilde{A}$ and $(A^{{\rm        op}})_E$
  amounts to  ${R^*}^{-1}=E_{21}R_{21}E^{-1}$.

 \end{rem}

   We discuss a simple example.

   \begin{ex}\label{simple_example}
   Let $G$ be a compact group and $C_\infty(G)$ the Hopf $^*$-algebra of functions on $G$ which are finite linear
   combinations of matrix coefficients of unitary finite dimensional
       representations $u$ of $G$. The coproduct and antipode are defined as usual by $\Delta(f)(g, h)=f(gh)$ and $S(f)(g)=f(g^{-1})$. Then the dual $^*$-algebra is isomorphic to $\Pi_{u\in{\rm Irr G}} B(H_u)$, with $H_u$ the Hilbert space of $u$. The algebraic direct sum $A=\bigoplus _{u\in{\rm Irr G}} B(H_u)$ is a discrete Hopf $^*$-algebra with dual
   coproduct $\hat{\Delta}$ and antipode $\hat{S}$. We have $A=A^{\rm op}$ by commutativity of
   $C_\infty(G)$, and it follows that with the trivial $R$-matrix and ribbon structure, $A$ is a discrete unitary coboundary  
   Hopf algebra.    \end{ex}
   
   The example gives a natural interpretation of axiom d) when $A$ is thought of as the dual of
   the algebra of functions on a noncommutative space.

\begin{rem} The relationship between the multiplier discrete algebra associated to
 the forgetful functor of ${\rm Rep}(U_q(g))$ and 
 $U_q(g)$ has been considered in detail by Neshveyev and Tuset
in  Sect. 2 in  \cite{CQGRC} for $q>0$, and it beautifully gives a connection between two different
 approaches to quantum groups by Woronowicz and   Drinfeld. Quite remarkably to us,
 the relevance of an analogous tannakian approach for a topological description of $U_q({\mathfrak g})$ has been
 explained by Sawin  in Sect. 1 in \cite{Sawin} motivated by the construction of the $R$-matrix.
\end{rem}

The following proposition gives a characterization of the Kac-type property for an antipode, see Def.  \ref{Kac_type_definition}.   

\begin{prop}\label{antipode_commuting_with_involution_wqh_case}
Let  $A$ be a   weak  Hopf algebra with a $^*$-involution making it into a $^*$-algebra,    strong antipode $S$
  such that $\Delta^{{\rm op}}(a)^*= \Delta^{{\rm op}}(I)^*\Delta(a^*)\Delta^{{\rm op}}(I)^*$ for all $a\in A$.
  Then   $S$   commutes with $^*$ (thus is of Kac type) if and only if $\sum_i a_iS(b_i^*)^*=I$, where $\Delta(I)=\sum_ia_i\otimes b_i$. This is always the case when $\Delta^{\rm op}(I)^*=\Delta(I)$, that is
  when $A$ is a  Hermitian   coboundary weak  Hopf algebra with compatible $^*$-involution in the sense of Sect. \ref{19}. 
\end{prop}

\begin{proof}
The necessity of the condition follows from the antipode axiom (\ref{eqn:antip1}). For the sufficiency, 
note that, if $S$ is a strong antipode then $(S, 1, 1)$
satisfies (\ref{eqn:antip1}) and by   Prop. \ref{weak_strong_antipode} and its proof this equation suffices
to make a triple $(S, 1, 1)$ into an antipode, with $S$ an antiautomorphism.
Starting with our assumptions, we may slightly modify the computations in the proof of Prop. \ref{Kac_type_sufficiency}
and show that $(\tilde{S}, 1, 1)$ satisfies (\ref{eqn:antip1}), with $\tilde{S}(a)=S(a^*)^*$, thus this is another strong antipode, and the proof is completed by uniqueness of a strong antipode.

\end{proof}

 Given  any central invertible element $z\in A$ with $\varepsilon(z)=1$
 we set $$\Theta_z:=z^{-1}\otimes z^{-1}\Delta(z), \quad \quad R_z:=R\Theta_z.$$
 (Note that when  $A$ is a weak  bialgebra, $E$ is necessarily a $2$-cocycle of $A^{{\rm        op}}$ by Prop. \ref{2-cocycle}.
Similarly,
$\Theta_z$ and $R_z$ are $2$-cocycles by Prop. \ref{2-cocycle_deformation}.)   We have   $A_{{\Theta_z}}=A$ as quasitriangular weak quasi bialgebras thanks to centrality of $z$ and since
  the twisting operation can be performed in stages, $(A_F)_G=A_{GF}$, we see that $R$ and $R_z$ twist $A$ in the same way. Therefore  
$$A_{R_z}=A_R= A^{{\rm        op}}.$$ Furthermore the deformed $R$-matrix yields a
generalised coboundary on ${\rm        Rep}(A)$ via $\Sigma\rho\otimes\sigma(R_z)$.
We set
\begin{equation}\label{overline_R_defn}\overline{R}=R\Theta_w.
\end{equation}
The element $\overline{R}$  first introduced by Drinfeld in his work on quasi-Hopf algebras    \cite{Drinfeld_quasi_hopf},  is the algebraic analogue of the  coboundary  symmetry considered  in  Sect. \ref{17}.
 
\begin{prop}\label{overline_R_prop}   The twist $\overline{R}$ satisfies
$\overline{R}_{21}\overline{R}=\Delta(I)$.  Therefore $\Sigma\rho\otimes\sigma(\overline{R})\in(\rho\underline{\otimes}\sigma, \sigma\underline{\otimes}\rho)$ is a coboundary of ${\rm        Rep}(A)$.
\end{prop}
\begin{proof}
We have
$$\overline{R}_{21}\overline{R}=R_{21}w^{-1}\otimes w^{-1}\Delta^{{\rm        op}}(w)Rw^{-1}\otimes w^{-1}\Delta(w)=$$
$$R_{21}Rw^{-2}\otimes w^{-2}\Delta(w^2)=R_{21}Rv^{-1}\otimes v^{-1}\Delta(v)=\Delta(I).$$

\end{proof}

 By axiom d), the element $E$ is required to be a trivial twist  from $A^{\rm op}$ to $\tilde{A}$. It follows that
 $\tilde{A}=(A^{\rm op})_E=({A}_{\overline{R}})_E=A_{E\overline{R}}$,
 hence \begin{equation}\label{a_tilde_equals_a_omega}\tilde{A}=A_\Omega, \quad \quad \Omega=E\overline{R}\end{equation} as quasitriangular weak quasi-bialgebras.

\begin{thm}\label{unitary_ribbon2}    Let  $A$ be a  (discrete) Hermitian coboundary   weak quasi-Hopf algebra.   Then  $A$ is    $\Omega$-involutive   with $\Omega=E\overline{R}=ER\Theta_w$. Furthermore the induced
braided  symmetry
$\Sigma\rho\otimes\sigma(R)\in(\rho\underline{\otimes}\sigma, \sigma\underline{\otimes}\rho)$, 
and therefore coboundary symmetry $\Sigma\rho\otimes\sigma(\overline{R})\in(\rho\underline{\otimes}\sigma, \sigma\underline{\otimes}\rho)$
are unitary  in ${\rm Rep}_h(A)$.
\end{thm}

\begin{proof} 
We need to show that $\Omega$ is selfadjoint.   By construction, $\Delta(I)$ and $\Delta(I)^*$ are respectively domain and range of $\Omega$.
The $R$-matrices of $\tilde{A}$  and $A^{{\rm        op}}$ are respectively given by $\tilde{R}={R^*}^{-1}$ and $R^{{\rm        op}}=R_{21}$ thanks to   Prop. \ref{R-canonicity}. Equality between the $R$-matrices of $\tilde{A}$ and $(A^{{\rm        op}})_E$ gives ${R^*}^{-1}=E_{21}R_{21}E^{-1}$, hence $R^*=ER_{21}^{-1}E_{21}^{-1}$. We may write $\Omega$ in the form 
\begin{equation}\label{comp1} \Omega=ER\Theta_w=ERw^{-1}\otimes w^{-1}\Delta(w)=\end{equation}
\begin{equation}\label{comp2} Ew^{-1}\otimes w^{-1}\Delta^{{\rm        op}}(w)R=w^{-1}\otimes w^{-1}E\Delta^{{\rm        op}}(w)R=\end{equation}
\begin{equation}\label{comp3} w^{-1}\otimes w^{-1}\tilde{\Delta}(w)ER.\end{equation}
We also have 
$$ E_{21}^{-1}E^*=(E^{-1})_{21}E^*=(\Delta^{{\rm        op}}(I)\Delta(I)^*)_{21}\Delta^{{\rm        op}}(I)^*\Delta(I)=\Delta(I).$$
Hence \begin{equation}\label{comp4}\Omega^*=R^*E^*\Delta(w^{-1})w\otimes w=\end{equation}
\begin{equation}\label{comp5} ER_{21}^{-1}E_{21}^{-1}E^*\Delta(w^{-1})w\otimes w=ER_{21}^{-1}\Delta(w^{-1})w\otimes w=\end{equation}
\begin{equation}\label{comp6} ER(R_{21}R)^{-1}\Delta(w^{-1})w\otimes w=ER\Delta(w)w^{-1}\otimes w^{-1}=\end{equation}
\begin{equation}\label{comp7} ER\Theta_w=\Omega.\end{equation}
Unitarity of the braided symmetry follows from the  property that $\tilde{A}=A_\Omega$ as quasi-triangular weak quasi-bialgebras    and 
  Prop. \ref{unitary_braided_symmetry2}.
\end{proof}

 In the rest of this section we endow $A$ with the $\Omega$-involution $\Omega=E\overline{R}$.
 Note that the Hermitian form  on the tensor product of  two representations associated to $\Omega$ is given by  
\begin{equation}
(\zeta, \zeta')_\Omega=(\zeta, \overline{R}\zeta').
\end{equation}

\begin{rem}
We may interpret the trivial twist $E$ as follows. It is non-trivial precisely when $\overline{R}$ is not selfadjoint. This follows
from the equation $\Omega=\Omega^*$. The subclass of Hermitian coboundary wqh for which $\overline{R}$ is already selfadjoint will be considered more closely in the next section.  \end{rem}

We discuss how to construct examples of Hermitian coboundary weak quasi-Hopf algebras with strongly trivial
$\Omega$-involution in the sense of Defn. \ref{trivial_involution}. The following example reduces the problem to the construction of 
Hermitian coboundary weak quasi-Hopf algebras with
trivial $\Omega$-involution. The step of constructing a unitary coboundary with trivial $\Omega$-involution will be considered in the next section. 

  \begin{prop}\label{from_trivial_to_strongly_trivial_coboundary}
Let $A$ be a Hermitian coboundary weak quasi-Hopf algebra with trivial involution 
  $\Omega=\Delta(I)^*\Delta(I)$, $\Omega^{-1}=\Delta(I)\Delta(I)^*$. Then the twist
  $T$ (or $T'$) defined in Remark \ref{twist_strongly_trivial} making the $\Omega$-involution strongly trivial
  turns $A$ into another Hermitian coboundary weak quasi-Hopf algebra $A_T$.
  \end{prop}
  \begin{proof}
Notice that axioms a)--c)   are invariant under  any twist. For axiom d), we have that 
  $T^*T_{21}=E$, and it easily follows that the twist $E_T=\Delta_T(I)^*\Delta_T^{\rm op}(I)$ has inverse
  $E_{T}^{-1}=\Delta_T^{\rm op}(I)\Delta_T(I)^*$. Moreover $\widetilde{A_T}=(\tilde{A})_{T^{-1*}}$ and 
  $(A_T)^{\rm op}=(A^{\rm op})_{T_{21}}$, it follows that 
  axiom d) for $A_T$ is equivalent to $\tilde{A}=(A^{\rm op})_{T^*E_TT_{21}}$. We have 
  $T^*E_TT_{21}=E$, thus axiom d) holds for $A_T$ also.
  \end{proof}

 \begin{defn}\label{ribbon_involutive}
Let $A$ be a  (discrete),    Hermitian   coboundary   weak quasi-Hopf (weak  Hopf) algebra.   If  
  $\Omega=E{R}\Theta_w$ is  positive in $A\otimes A$ ($M(A\otimes A)$) then $A$ will be called a  {\it unitary   
(discrete),    coboundary}   weak quasi-Hopf (weak  Hopf) algebra. \end{defn}

Recall \cite{Rowell2, Turaev} that an {\it Hermitian (unitary) ribbon  category} is a *-category ($C^*$-category)   ${\mathcal C}$ equipped with a right duality
$(\rho^\vee, b_\rho, d_\rho)$, unitary braided symmetry $\varepsilon(\rho, \sigma)$ and unitary ribbon   structure $v\in(1,1)$ such that
 \begin{equation}\label{Hermitian_ribbon}b_\rho^*=d_\rho\circ \varepsilon(\rho, \rho^\vee)\circ v_\rho^{-1}\otimes 1_{\rho^\vee}\quad\quad d_\rho^*=1_{\rho^\vee}\otimes v_\rho\circ \varepsilon(\rho^\vee, \rho)^{-1}\circ b_\rho.\end{equation}
 It follows from (\ref{quantum_dimension_pivotal}) that the quantum dimension in a Hermitian ribbon category 
 may be computed as
 $$d(\rho)=d_\rho d_\rho^*=b_\rho^*b_\rho.$$

\begin{thm}\label{Hermitian_category2} Let $A$ be a  Hermitian (unitary) coboundary    weak quasi-Hopf algebra with
an antipode $(S, \alpha, \beta)$ such that $S$ commutes with $^*$. Then
  ${\rm   Rep}_h(A)$ $({\rm Rep}^+(A))$ is a   Hermitian (unitary) ribbon category
  with the canonical duality $(\rho^\vee=\rho_c, b_\rho, d_\rho)$   associated to $A$ as in Example \ref{computing_conjugates}
   if and only if $\beta=\alpha^*$.
  This equation holds if $A$ is discrete and admits an antipode  of Kac type.
\end{thm}

\begin{proof}
 We need  to give   a right duality $(\rho^\vee, b_\rho, d_\rho)$ satisfying (\ref{Hermitian_ribbon}).  We show that this holds for  the canonical duality $(\rho^\vee=\rho_c, b_\rho, d_\rho)$   associated to $A$ as in Example \ref{computing_conjugates} and a fixed antipode $(S, \alpha, \beta)$ such that
  $S$ commutes with $^*$. We only verify the equation on the right in (\ref{Hermitian_ribbon}).
 We have  $$d_\rho^*=r_\rho=\Omega^{-1}\sum\mu_i{\overline{e}_i}\otimes \alpha^* e_i, \quad\quad b_\rho=\sum_i \beta e_i\otimes \mu_i\overline{e}_i,$$ with $e_i$ an orthonormal basis.  A computation gives for $a$, $b\in A$,      
 $$a\otimes b\sum\mu_i\overline{e}_i\otimes \alpha^* e_i=\sum\mu_i\overline{e}_i\otimes b\alpha^*S(a)e_i.$$
 Taking into account $S(w)=w$, $\varepsilon(w)=1$, and the antipode property (\ref{eqn:antip1}), it follows that 
 $$\Delta(w)^* w\otimes w \sum\mu_i\overline{e}_i\otimes \alpha^*e_i= 
 \sum\mu_i\overline{e}_i\otimes \alpha^* ve_i.$$
On the other hand,  $\Omega^{-1}=R^{-1}\Delta(w)^* w\otimes w$. It follows that   
$$d_\rho^*=\Omega^{-1}\sum\mu_i\overline{e}_i\otimes \alpha^*e_i=1\otimes v_\rho R^{-1} \sum\mu_i\overline{e}_i\otimes \alpha^*e_i=$$
$$ 1\otimes v_\rho\circ \varepsilon(\rho^\vee, \rho)^{-1} \circ \sum\alpha^*e_i \otimes 
\mu_i\overline{e}_i .$$
Thus  the equation on the right in (\ref{Hermitian_ribbon}) holds if and only if
$ \beta=\alpha^*.$
    \end{proof}

We next identify the element $\omega$ defined 
  in Prop. \ref{involution_antipode} with the element  associated to the spherical structure on ${\rm Rep}(A)$,
  as in Theorem \ref{Deligne}, see also Cor. \ref{spherical_structure_wqh},
   in the important special case of antipode of Kac type, see Def. \ref{Kac_type_definition}.

\begin{prop}\label{computing_little_omega}
Let $A$ be an Hermitian coboundary weak quasi-Hopf algebra with antipode $S$ of Kac type.
 Then $\omega=uv^{-1}$ where $u$ is Drinfeld element associated to $S$ introduced in Definition \ref{Drinfeld_element_u}.
\end{prop} 

\begin{proof}
Since  $A$ has a strong antipode $S$,   $\omega=m\circ S\otimes 1(\Omega^{-1})$ by (\ref{little_omega}). With the same notation as in
Def. \ref{Drinfeld_element_u}, we have $\Omega^{-1}=R^{-1}\Delta(w)^*w\otimes w=\sum_j \overline{r}_jw_1^*w\otimes\overline{t}_jw_2^*w$.
Recall from Prop. \ref{inner_antipode} that $S^2$ is the inner automorphism induced by $u$ and that
$u^{-1}=\sum_j S^{-1}(\overline{t}_j)\overline{r}_j$. 
 It follows that 
 $\omega=\sum_j S(\overline{r}_jw_1^*w)\overline{t}_jw_2^*w=S(S^{-1}(\overline{t}_jw_2^*w)\overline{r}_jw_1^*w)=
 S(S^{-1}(w_2^*w)u^{-1}w_1^*w)=S(u^{-1}S(w_2^*w)w_1^*w)=S(u^{-1}v),$
 for the last equality we have used axiom c) of Def. \ref{Hermitian_ribbon_wqh}. On the other hand, 
 $S(u^{-1}v)=uv^{-1}$ by 
Remark \ref{central_square_root_of_squared_Drinfeld_element}.

\end{proof}

Let $A$ be an Hermitian coboundary weak quasi-Hopf algebra. Replacing the choice of  $w$ with another   square root $w'$ of $v$ satisfying the properties in  c) of Def. \ref{Hermitian_ribbon_wqh}
gives rise to another Hermitian coboundary weak quasi-Hopf algebra
with the same the same structure as $A$ and new square root of the ribbon element given by $w'$, and correspondingly 
a new $\Omega_{w'}$,
and therefore a new tensor $^*$-category, denoted ${\rm        Rep}'_h(A)$. 
We may write   $w'=wy$ with $y$ a
   (unitary) central square root $y$ of $I$ in $M(A)$ satisfying c), that is $\varepsilon(y)=1$ and $S(y)=y$.
   Conversely, any
   $y\in M(A)$ with these properties arises in this way.
   The new $\Omega_{w'}$ differs from $\Theta_w$ by the
 $2$-coboundary $\Theta_y=y^{-1}\otimes y^{-1}\Delta(y)$, that is
$$\Omega_{w'}=\Omega_{w}\Theta_y.$$
In particular, $\Omega_{w'}=\Omega_w$ if and only if $y$ is a $1$-cocycle: $\Delta(y)=y\otimes y\Delta(I)$.

\begin{prop}\label{equivalence} Assume that $A$ is discrete. 
\begin{itemize}
\item[{\rm        a)}] 
The functor ${\mathcal F}: {\rm        Rep}_h(A)\to {\rm        Rep}'_h(A)$ acting identically on objects and morphisms with identity natural transformation $F_{\rho, \sigma}$ is a    tensor $^*$-functor and an equivalence.
There is no  unitary tensor $^*$-functor between these categories unless 
 $y=w'w^{-1}$ is a $1$-cocycle.
 \item[{\rm        b)}] If  $A$ is a unitary coboundary weak quasi-Hopf algebra with respect to $w$ and  $\Omega_{w'}$ is positive
 with respect to some   other $w'$ satisfying c) in Def. \ref{Hermitian_ribbon_wqh} then $\Omega_{w}=\Omega_{w'}$
 and $\overline{R}_w=\overline{R}_{w'}$.

 \end{itemize}
\end{prop}
\begin{proof} a)
The categories ${\rm        Rep}_h(A)$ and ${\rm        Rep}'_h(A)$ have the same tensor structure and the same $^*$-category 
structure, and the   functor  ${\mathcal F}$ becomes the identity functor for these substructures, thus it  is    a tensor $^*$-functor and an equivalence when the natural transformation $F_{\rho, \sigma}: {\mathcal F}(\rho)\otimes'{\mathcal F}(\sigma)\to{\mathcal F}(\rho\otimes\sigma)$ acts as identity. Here we have used different symbols to denote the two different tensor products. The Hermitian form of ${\mathcal F}(\rho)\otimes'{\mathcal F}(\sigma)$ differs from that of  ${\mathcal F}(\rho\otimes\sigma)$ by the action of $\Theta_y$. On the other hand $y$ acts as $\varepsilon_\rho$,
where $\varepsilon_\rho=\pm1$ on an irreducible representation $\rho$, and  all the  $\varepsilon_\rho$ determine $y$. In particular, $\Theta_y$ at most changes the sign of the Hermitian form of an irreducible component of 
${\mathcal F}(\rho\otimes\sigma)$, and if this happens then ${\mathcal F}(\rho)\otimes'{\mathcal F}(\sigma)$ and  ${\mathcal F}(\rho\otimes\sigma)$ are not unitarily equivalent, by Prop. \ref{complete}. Hence all $F_{\rho, \sigma}$ 
are unitary if and only if $\Theta_y=\Delta(I)$.
b) If two choices $w$ and $w'$  both define positive operators $\Omega_w$ and $\Omega_{w'}$  then 
 ${\mathcal F}:{\rm Rep^{+}(A)\to Rep^{+'}(A)}$ is a tensor $^*$-equivalence between tensor $C^*$-categories
 hence by Prop. \ref{polar_decomposition} c) polar decomposition of the tensor structure gives   a unitary tensor equivalence. It follows from the   previous part that
 $y=w'w^*$ is a $1$-cocycle, hence $\Omega_w=\Omega_{w'}$ and also
 $\overline{R}_w=\overline{R}_{w'}$. 

\end{proof}

 We next construct     an involutive antipode for all the twists of a unitary ribbon weak quasi-Hopf algebra of Kac type   under a spectrum condition.

 \begin{prop}
 Let $A$ be a   unitary coboundary weak quasi-Hopf algebra.  If $A$ has an antipode of Kac type $S$ with associated unitary Drinfeld   element $u$ such that $-1\notin{\rm Sp}(uv^{-1})$    then 
 for any twist $F$ of $A$, $A_F$ endowed with twisted involution $(^*, \Omega_F)$ admits an antipode
 $(\tilde{S}, \tilde{\alpha}, \tilde{\beta})$ such that   $\tilde{S}$  commutes with $^*$ and the corresponding  element
 as in Prop. \ref{involution_antipode} is $\tilde{\omega}=1$. In particular,
  ${\tilde{S}}^2=1$.
 \end{prop}

 \begin{proof}  The element $\omega$ corresponding to $S$ and $A$ is given by
  $\omega=uv^{-1}$
by Prop. \ref{computing_little_omega}, which is    unitary. Let $(S, \alpha, \beta)$ be the twisted antipode 
of $A_F$ as in (\ref{twisted_antipode}), so $\omega_F=\omega$ by Prop. \ref{invariance_little_omega_under_twisting} b).
For an invertible $x$,   the antipode $({\rm Ad}(x)S, x\alpha, \beta x^{-1})$ of $A_F$ has associated element $\tilde{\omega}=x\omega S^{-1}(x)^*$ by Prop. \ref{invariance_little_omega_under_twisting} a).
  We set $x^{-1}=\omega^{1/2}$, the continuous functional calculus of
the principal branch of the square root function, so $x$ is unitary. Since $S(\omega)=\omega^{-1}$ it follows that $S(x)=x^{-1}$ and therefore  $\tilde{\omega}=1$.
  \end{proof}

We shall see that the weak quasi-Hopf algebras arising from VOAs as satisfying the rationality assumptions of  Sect. \ref{VOAnets} have a natural involutive antipode commuting with $^*$.

\section{A categorical characterization of discrete  Hermitian coboundary wqh}\label{DR}

Let $A$ be a discrete hermitian coboundary weak quasi-Hopf algebra. In this section we give a categorical characterization of the axioms defining $A$.  For simplicity, we discuss a detailed proof only in   the unitary case.
It turns out that the categorical definition simplifies considerably this structure. 

Then the $C^*$-structure of $A$
gives rise to   the linear $C^*$-category
${\mathcal C}^+={\rm Rep}^+(A)$ of Hilbert space representations of $A$. We also have 
the tensor category ${\mathcal C}={\rm Rep}(A)$ of vector space representations of $A$ which has additional structure, the braiding, the ribbon structure and coboundary symmetry.

We recall  that ribbon and coboundary structures in tensor categories have been studied in Sect. \ref{17}. In particular, 
  by Theorem
  \ref{coboundary_associated_to_a_ribbon_structure}
 a coboundary symmetry $c^w$
  may be associated to a ribbon category ${\mathcal C}$ with braided symmetry $c$, ribbon structure $v\in(1,1)$ when
  there is a natural isomorphism $w\in(1,1)$ compatible with duality which is a square root of $v$.
  In the setting of unitary categories when $c$  and $w$ are unitary then $c^w$ is unitary, and also selfadjoint as
  $(c^w)^2=1$.

\begin{thm}\label{TK_unitary_ribbon}
Let $({\mathcal C}, \otimes, \alpha, c, v)$ be a     ribbon category,     
  $w\in(1,1)$    a   square root of $v$ compatible with duality,
   ${\mathcal C}^+$ a semisimple $C^*$-category
and let ${\mathcal F}: {\mathcal C}^+\to{\mathcal C}$ be a   linear equivalence.
Let $({\mathcal G}, F, G): {\mathcal C}\to{\rm Herm}$ (${\mathcal C}\to{\rm Hilb}$) be a faithful weak quasitensor functor
with  symmetric dimension function such that
  ${\mathcal G}^+={\mathcal G}{\mathcal F} : {\mathcal C}^+\to{\rm Herm}$ (${\mathcal C}^+\to{\rm Hilb}$) is a $^*$-functor.
 Then the discrete pre-$C^*$-algebra  $A^+={\rm Nat}_0({\mathcal G}^+)$ endowed with the natural    ribbon weak quasi-Hopf algebra structure and  $\Omega$-involution    induced
  by duality becomes an Hermitian (unitary)  coboundary weak quasi-Hopf algebra
  if and only if   ${\mathcal G}(\alpha)$, ${\mathcal G}(c)$, and ${\mathcal G}(v)$   are unitary and moreover $(F, G)$ satisfies the following conditions
  \begin{equation}\label{unitary_ribbon_functor1}    F_{\sigma, \rho}\Sigma(\rho, \sigma) F_{\rho, \sigma}^* =  {\mathcal G}(c^w(\rho, \sigma))\end{equation}
\begin{equation}\label{unitary_ribbon_functor2}   G_{\sigma, \rho}^* \Sigma(\rho, \sigma) G_{\rho, \sigma}   ={\mathcal G}(c^w(\sigma, \rho)^{-1})
\end{equation}
with $\Sigma$ the permutation symmetry of ${\rm Herm}$  (${\rm Hilb}$).
 In this case, ${\mathcal C}^+$ becomes an Hermitian (unitary) ribbon tensor category and $ {\mathcal F}: {\mathcal C}^+\to{\mathcal C}$ 
 a ribbon tensor equivalence.

Moreover   there is a   unitary ribbon tensor equivalence ${\mathcal E}:{\mathcal C}^+\to{\rm Rep}_h(A)$ 
(${\mathcal C}^+\to{\rm Rep}^+(A)$) preserving the coboundary structures such that ${\mathcal F}_A{\mathcal E}\simeq {\mathcal G}$ unitarily monoidally, with ${\mathcal F}_A$ the forgetful functor of ${\rm Rep}^+(A)$.

\end{thm}

\begin{rem}\label{DR_symmtric_functor}
a) In the particular case
of a trivial $R$-matrix, trivial ribbon element, ordinary discrete (cocommutative) Hopf $C^*$-algebra,
the categorical definition of unitary coboundary weak quasi-Hopf algebra reduces to the property that
 embedding tensor functor ${\rm Rep}(A)\to{\rm Hilb}$ to the category of Hilbert spaces is a {\it symmetric   functor} , by Theorem \ref{TK_unitary_ribbon}. For this reason, the class of unitary coboundary weak Hopf $C^*$-algebras plays a natural role in the Doplicher-Roberts program for braided tensor $C^*$-categories, see Problem \ref{problem1}. 

b) In the non-cocommutative case ($R\neq I$) with   structure maps satisfying $F=G^*$ the   equations \ref{unitary_ribbon_functor1}, \ref{unitary_ribbon_functor2} say that the weak quasitensor functor $({\mathcal G}, F, G)$ preserves the coboundary symmetries of ${\mathcal C}^+$ and ${\rm Hilb}$.
In particular, under the requirements $F=G^*$,   equations \ref{unitary_ribbon_functor1}, \ref{unitary_ribbon_functor2}  uniquely determine the structure of the braided symmetry $c(\rho, \sigma)$ of  ${\mathcal C}^+$ on an irreducible component of the minimum energy subspace (and hence determine the braided symmetry completely) as a scalar multiple of the permutation symmetry by   the action
 of the balancing element $w$.
 This form of the praided symmetry occurs in categories arising from conformal field theory, with the braiding associated to the monodromy of the KZ equations, see Corollary \ref{Toledano_laredo} for the loop group conformal net and Remark \ref{HL_braiding} for the affine vertex operator algebra setting. This highlights the relevance of making module categories of vertex operator algebras into unitary categories.
\end{rem}

  \begin{proof}  
  By theorem \ref{TK_algebraic_quasi}, $A={\rm Nat}_0({\mathcal G})$ becomes a ribbon weak quasi-Hopf algebra with coproduct and associator defined by $(F, G)$ and $R$-matrix $R$ defined in the proof and ribbon structure
  ${\mathcal G}(v_\rho)$.
  We transfer this structure to $A^+$ via the isomorphism $A\to A^+$ induced by   ${\mathcal F}$.
By Prop. \ref{unitarizability2},  $A^+$  becomes naturally a unitary weak quasi-bialgebra   if and only if ${\mathcal G}(\alpha)$ is unitary. The $\Omega$-involution of $A^+$ is given by $\Omega=\tilde{F}^*\tilde{F}$ and $\Omega^{-1}=\tilde{G}\tilde{G}^*$ respectively, where $\tilde{F}$ and $\tilde{G}$ correspond to $F$ and $G$ via the isomorphism, as in the proof of
  Prop. \ref{unitarizability2}. We have $\tilde{A^+}=A^+_\Omega$ as weak quasi-bialgebras. When ${\mathcal G}(v_\rho)$ is unitary then the natural transformation ${\mathcal G}(w_{{\mathcal F}(x))})$ defines a unitary
  square root of the ribbon structure of $A^+$, and  axioms a), b), c) of Def. \ref{Hermitian_ribbon_wqh} hold.
   We   show   that with this structure axiom d) is equivalent to (\ref{unitary_ribbon_functor1}) and 
   (\ref{unitary_ribbon_functor2}) if ${\mathcal G}(c)$ is unitary.
 Note that ${A}^{\rm op}=(A)_{\overline{R}}$ as quasitriangular  weak quasi-bialgebras. It follows that d)
 may equivalently be formulated as $\tilde{A^+}=(A^+)_{E\overline{R}}$
   as quasitriangular weak quasi-bialgebras together with the requirement that $E$ is a trivial twist, that is 
   $E^{-1}=E'$, where
   $E=\Delta(I)^*\Delta^{\rm op}(I)$ and $E'=\Delta^{\rm op}(I)\Delta(I)^* $.
   On the other hand, 
   equations (\ref{unitary_ribbon_functor1}) and  (\ref{unitary_ribbon_functor2}) are respectively equivalent to
 \begin{equation}\label{unitary_ribbon_functor5}
 \Sigma G_{\sigma, \rho}F_{\sigma, \rho}\Sigma F_{\rho, \sigma}^*F_{\rho, \sigma}=\Sigma G_{\sigma, \rho}{\mathcal G}(c^w(\rho, \sigma))F_{\rho, \sigma},\end{equation}
   \begin{equation}\label{unitary_ribbon_functor6} G_{\sigma, \rho}G_{\sigma, \rho}^*\Sigma G_{\rho, \sigma}F_{\rho, \sigma}\Sigma=G_{\sigma, \rho}{\mathcal G}(c^w(\sigma, \rho)^{-1})F_{\rho, \sigma}\Sigma.\end{equation}
We know that $R$ and $R^{-1}$ correspond   to   $\Sigma G_{\sigma, \rho}{\mathcal G}(c(\rho, \sigma))F_{\rho, \sigma}$
   and   $G_{\sigma, \rho}{\mathcal G}(c(\sigma, \rho)^{-1})F_{\rho, \sigma}\Sigma$. It follows from a computation that
$\overline{R}$  and $\overline{R}^{-1}$ in turn correspond to 
   $\Sigma G_{\sigma, \rho}{\mathcal G}(c^w(\rho, \sigma))F_{\rho, \sigma}$
   and  $G_{\sigma, \rho}{\mathcal G}(c^w(\sigma, \rho)^{-1})F_{\rho, \sigma}\Sigma$.
 It follows that equations (\ref{unitary_ribbon_functor5}) and  (\ref{unitary_ribbon_functor6}) are in turn equivalent  to $E'\Omega=\overline{R}$, $\Omega^{-1}E=\overline{R}^{-1}$, in other words $E'=E^{-1}$ and $\Omega=E\overline{R}$.  
 On the other hand, the $R$-matrices of $\tilde{A}$ and
 $A_\Omega$   coincide  by
    Prop. \ref{unitary_braided_symmetry2}   as ${\mathcal G}(c)$ is a unitary braided symmetry. Thus 
   the proof of      axiom d) is complete. Conversely,
  when $A$ is a unitary coboundary weak quasi-Hopf algebra and ${\mathcal C}={\rm Rep}(A)$ then the natural
 weak quasi-tensor structure of the forgetful ${\rm Rep}(A)\to{\rm Hilb}$ satisfies
$F_{\rho, \sigma}^*=\rho\otimes\sigma(\Delta(I)^*)\Omega_{\rho, \sigma}$ and similarly 
$G_{\rho, \sigma}^*=\Omega_{\rho, \sigma}^{-1}\rho\otimes\sigma(\Delta(I)^*)$.
Moreover $c^w$ corresponds to $\Sigma\overline{R}$. It follows that the unitarity statements and
(\ref{unitary_ribbon_functor1}) and (\ref{unitary_ribbon_functor2}) are verified.
 The property that ${\mathcal C}^+$   is an Hermitian  (unitary) ribbon category   follows from Theorem \ref{Hermitian_category2}.    
 In the unitary case  it also follows  that the canonical tensor equivalence ${\mathcal E}$ described in Theorem \ref{TK_algebraic_quasi} is unitary by Prop. \ref{unitarizability2}, see also Theorem \ref{TheoremTannakaStar} and preserves the coboundary symmetries by   construction.
      \end{proof}

\begin{rem}\label{remark_to_unitary_ribbon}
a) It follows from the proof of Theorem \ref{TK_unitary_ribbon} and that of Theorem \ref{Positivity}  that when
  ${\mathcal C}$ has a generating object (i.e. its powers contain every irreducible as a subobject) then
 Theorem   \ref{TK_unitary_ribbon} holds if equations
   (\ref{unitary_ribbon_functor1}) and (\ref{unitary_ribbon_functor2}) 
   are known to hold only for   pairs  $\rho$, $\sigma$ such that one of them, say $\rho$, is the generating object and the other varies among the irreducible objects
   of ${\mathcal C}$, or alternatively among the choice of a tensor power $\rho^n$ for each integer $n$.     
 b) It follows from Theorem \ref{TK_unitary_ribbon} and Prop.  \ref{equivalence} b) that $c^w$ does not depend on the choice of $w$. 
  \end{rem}

  \begin{rem}
It follows  that Theorem \ref{Hermitian_category2} admits a categorical formulation as well.
Indeed, we may define a discrete coboundary weak quasi-Hopf algebra $A$ equivalently as a semisimple
ribbon tensor category  $({\mathcal C}, \otimes, \alpha, c, v)$ endowed   
with a square root $w$ of the ribbon element $v$ compatible with duality and 
 the structure of a $C^*$-category with a weak quasi-tensor faithful functor $({\mathcal G}, F, G)$ such that all the conditions of Theorem   \ref{TK_unitary_ribbon} hold (we are choosing ${\mathcal C}^+={\mathcal C}$ and ${\mathcal F}$ identity).
 Then any right duality $(\rho^\vee, b_\rho, d_\rho)$ is of the form described in Example \ref{computing_conjugates}
 by the proof of Theorem \ref{TK_algebraic_quasi} (d),
 as an antipode $(S, \alpha, \beta)$ may always be chosen such that $S$ commutes with $^*$ by Remark \ref{antipode_commuting_with_*}.
 Thus by Theorem  \ref{Hermitian_category2}, the condition $\beta=\alpha^*$ is equivalent to the compatibility equations  (\ref{Hermitian_ribbon}) making ${\mathcal C}$
 into an Hermitian ribbon category with respect to $(c, v, \rho^\vee, b_\rho, d_\rho)$.

 \end{rem}
 
Taking into account the historical motivation briefly discussed in the introductory part of  Section
\ref{18}, 
we are led to look for special examples with $R$-matrix given by a weak analogue of a $2$-coboundary. 
 The next remark shows  that  
the construction of these    examples is   related  to  the study of unitary structures,
having a suitable triviality property.

\begin{rem}\label{properties_of_trivial_unitary_ribbon2}
If  an Hermitian coboundary $A$  has trivial involution as  introduced in Def. \ref{trivial_involution} then 
by definition $\Omega$ is a trivial twist, thus  
we have   from relation (\ref{a_tilde_equals_a_omega})  
\begin{equation}\label{data1}\Omega=\Delta(I)^*\Delta(I), \quad\quad {R}=\Delta^{\rm op}(I)\Delta(I)^*w\otimes w\Delta(w^{-1}),\quad \quad
\overline{R}=\Delta^{\rm op}(I)\Delta(I)^*\Delta(I).\end{equation}
Conversely, if $A$ is   Hermitian coboundary and the $R$-matrix takes the previous form then necessarily the 
involution is trivial.
We have a particular case,  when the $\Omega$-involution of $A$  is strongly trivial   ($\Omega=\Omega^{-1}=\Delta(I)$)
then 
\begin{equation}\label{data2}\Omega=\Delta(I), \quad\quad {R}=\Delta^{\rm op}(I)w\otimes w\Delta(w^{-1}),\quad \quad
\overline{R}=\Delta^{\rm op}(I)\Delta(I).\end{equation}
We recall from Example \ref{from_trivial_to_strongly_trivial_coboundary}
that strongly trivial $\Omega$-involutions can be obtained from trivial $\Omega$-involutions via suitable twisting.
Moreover, when $A$ is in addition unitary discrete than any trivial $\Omega$-involution is strongly trivial
by Prop. \ref{strongly_unitary_prop}.

\end{rem}

 In the next  section we consider the question of constructing  
new examples of unitary coboundary weak quasi-Hopf algebras with such triviality properties from old ones, and Theorem \ref{TK_unitary_ribbon}
will turn out useful.
To construct such examples,  we
look for twist deformation of given examples that
respect the structure, that may perhaps be regarded as an abstract analytic analogue of part of the arguments
involved  Drinfeld-Kohno theorem following \cite{Drinfeld_quasi_hopf}.

  \section{Compatible  unitary coboundary  wqh,  an abstract Drinfeld-Kohno}\label{19}
  
  We know from a theorem of Galindo \cite{Gal} that a braiding of a unitary fusion category is always unitary.
  Now we reverse the question and ask is there a way of constructing a unitary braided tensor category with a unitary braiding, equivalent to a given   a braided semisimple
  tensor category ${\mathcal C}$?
  In   applications we may already have a linear $C^*$-category ${\mathcal C}^+$ and a linear equivalence
  ${\mathcal F}:{\mathcal C}^+\to{\mathcal C}$  and  we want to turn ${\mathcal C}^+$ into
  a unitary braided  tensor category.
  Furthermore, if we have two   braided tensor categories ${\mathcal C}_1$ and ${\mathcal C}_2$ 
  which are linearly equivalent
 to  the same $C^*$-category ${\mathcal C}^+$
 via ${\mathcal F}_i:{\mathcal C}^+\to{\mathcal C_i}$, $i=1$, $2$,
 under what circumstances the corresponding constructions give unitarily   equivalent braided tensor categories
 ${\mathcal C_1}^+$ and ${\mathcal C_2}^+$?
  If this can be achieved, it will follow in particular that ${\mathcal C}_1$ and 
 ${\mathcal C}_2$ 
   are  also equivalent as braided tensor categories.    
  In this section we set up a specific situation and we construct a unitary braided tensor quasi-equivalence
   $({\mathcal E}, E): {\mathcal C_1}^+\to {\mathcal C_2}^+$. In other words we reduce the problem
  to verification of  the equation concerning the associativity morphisms only, that is equation
  (\ref{wt1}), (\ref{wt2}) with (${\mathcal E}$, $E$, $E^{-1}$) in place of
 (${\mathcal F}$, $F$, $G$). In doing this, we follow ideas of Drinfeld \cite{Drinfeld_quasi_hopf} in his work on Drinfeld-Kohno theorem,
 except for as already said we forget the associativity morphisms, and again ideas of Wenzl \cite{Wenzl} in his work of the unitary structures
 of fusion categories ${\mathcal C}({\mathfrak g}, q, \ell)$ of quantum groups at roots of unity.

In the   introduction of Section \ref{18}  we have  interpreted axiom d) of Definition \ref{Hermitian_ribbon_wqh}
as a noncommutativity property of the  function algebra from a dual viewpoint.
    This interpretation disregards the trivial twist $E$, and therefore becomes more meaningful 
  when the trivial twist is actually trivial.
   This leads us to the following stronger definition.
   
   \begin{defn}\label{strongly_Hermitian_ribbon_wqh}
A  Hermitian   coboundary weak quasi-Hopf algebra $(A, \Delta, \Phi, R, v, w, {}^*, S, \alpha, \beta)$  is called {\it compatible with the $^*$-involution} if it satisfies one of the following equivalent conditions,
\begin{itemize}
\item[{\rm  1)}] $E=\Delta(I)^*=\Delta^{\rm op}(I)$, 
\item[{\rm  2)}]  $\Delta(a)^*=\Delta^{\rm op}(a^*)$, $a\in A$, 
\item[{\rm  3)}]  $\Omega=\overline{R}$,
\item[{\rm  4)}] $\overline{R}$ is selfadjoint.

\end{itemize}
Thus
axiom d) of Def. \ref{Hermitian_ribbon_wqh} is replaced by the stronger axiom
\begin{itemize}
\item[{\rm        d')}] 
$\tilde{A}=A^{{\rm        op}}$ as quasitriangular weak quasi-bialgebras.
\end{itemize}
In particular we have ${R^*}^{-1}=R_{21}$.
Unitary, discrete, or weak Hopf versions are naturally defined.

\end{defn}

Example \ref{simple_example} is of this kind. In Sect. \ref{20} we construct examples associated to fusion categories
${\mathcal C}({\mathfrak g}, q, \ell)$ associated to $U_q({\mathfrak g})$ at certain roots of unity with compatible $^*$-involution.

In the rest of the paper we restrict to the unitary case.

\begin{prop}\label{Positivity_coboundary} Let $A$ be a discrete unitary coboundary wqh with a generating representation $\rho$.
Then 
$A$ has compatible $^*$-involution if and only if
 $$\sigma\otimes\rho(\overline{R}), \quad \sigma\otimes\rho(\overline{R}^{-1})$$ 
$$\sigma\otimes\rho\otimes\rho(I\otimes \overline{R}1\otimes\Delta(\overline{R})), \quad \rho\otimes\rho\otimes\sigma(\overline{R}\otimes I\Delta\otimes 1(\overline{R}))$$
are positive for every irreducible representation $\sigma$. \end{prop}

\begin{proof} Necessity is clear.
We note that the associativity morphisms $\sigma\otimes\rho\otimes\rho(\Phi)$ and $\rho\otimes\rho\otimes\sigma(\Phi)$
are unitary w.r.t. the given unitary coboundary structure, which is defined by $\overline{R}$ on the involved subspaces.
By Theorem \ref{Positivity}, $A$ becomes a unitary coboundary wqh with compatible $^*$-involution. On the other hand
the original coboundary structure $\Omega=\Delta(I)^*\overline{R}$  and the new compatible coboundary structure $\Omega'=\overline{R}$  coincide of the spaces of $\sigma\otimes \rho$ and therefore coincide everywhere by the conclusion of Theorem \ref{Positivity}.

\end{proof}

The following remark is an analogue of Remark \ref{properties_of_trivial_unitary_ribbon2} for the subclass
of wqh of this section, and takes a perhaps remarkable stronger form that seems to remind of the form taken by the $R$-matrix in the specific case of 
Drinfeld category
\cite{Drinfeld_quasi_hopf} for quasi-Hopf algebras.

\begin{rem}\label{unitary_braiding_data_of_strong_unitarization}
Let
$A$ be a discrete Hermitian   coboundary wqh with compatible $^*$-involution  and let    $\Omega$ by a trivial involution
with respect to a tensor product representation $\rho\otimes\sigma$.
Thus we have on the space of $\rho\otimes\sigma$,
\begin{equation}\label{data3}
\Omega=\Delta(I)^*\Delta(I)=\overline{R}, \quad\quad {R}=\Delta(I)^*w\otimes w\Delta(w^{-1}).\end{equation}
When the $\Omega$-involution is in addition strongly trivial with respect to
$\rho\otimes\sigma$ (recall that this is automatic when $A$ is discrete unitary
by Prop. \ref{strongly_unitary_prop})
 then in the representation space of $\rho\otimes\sigma$,
\begin{equation}\label{data4}\Omega=\Delta(I)=\overline{R}, \quad\quad {R}=w\otimes w\Delta(w^{-1}).\end{equation}
In particular, if this holds for any pair of simple representations
$\rho$, $\sigma$ then $A$ has a cocommutative coproduct $(\Delta=\Delta^{({\rm op})}$) by centrality of $w$.
\end{rem}

 Note that  if $T$ is a twist of $A$ with left inverse $T^{-1}$ then by definition
 $\Delta(I)$ is the domain of $T$ and range of $T^{-1}$.
If $A$ has a compatible $^*$-involution then we also have that
 ${T^{-1}}^*_{21}$ has domain ${\Delta^{\rm op}(I)}^*=\Delta(I)$,
 $T_{21}$ has domain $\Delta^{\rm op}(I)=\Delta(I)^*$, $T^*$ has range $\Delta(I)^*$.
 \bigskip
 
 The following result is our abstract analogue of Drinfeld-Kohno theorem.

\begin{thm}\label{Drinfeld_Kohno}
Let $A=(A, \Delta, \Phi, R, v, w)$ be a discrete  unitary   coboundary weak quasi-Hopf algebra with compatible $^*$-involution $(^*, \Omega=\overline{R}\geq 0)$. Let $(T, T^{-1}, P, Q)$ be a quadruple of elements in $M(A\otimes A)$:
such that $T$ is a twist of $A$ with left inverse $T^{-1}$, 
$P$, $Q$ are selfadjoint projections in $M(A\otimes A)$ such that  
$$PQ=0,\quad\quad P+Q=I $$
$$T={(T^{-1})_{21}^*}, \quad\quad \overline{R}=T^*T_-,  \quad\quad \overline{R}^{-1}=T^{-1}_-(T^{-1})^*.$$
where $$T_-=(P-Q)T, \quad\quad T^{-1}_-=T^{-1}(P-Q).$$
  Then  \begin{itemize}
\item[{\rm  a)}] 
  $A_T$ is another discrete unitary coboundary weak quasi-Hopf algebra with compatible involution,
 and associated quadruple $(\Delta_T(I), \Delta_T(I), P, Q)$ satisfying the same properties with respect to the twisted
 structure  of $A_T$: 
$$\Omega_T=\Delta_T^*(I)\Delta_T(I)_-=\Delta_T^{\rm op}(I)\Delta_T(I)_-=\overline{R}_T\geq0$$
$$ {R}_T=\Delta_T^{\rm op}(I)\Delta_T(I)_-w\otimes w\Delta_T(w^{-1}),$$
where $\Delta_T(I)_-=(P-Q)\Delta_T(I)$,

 \item[{\rm  b)}] If $\rho$ and $\sigma$ are two Hilbert space $^*$-representations of $A$ such that 
 $\rho\otimes\sigma(QT)=0$ and $\rho\otimes\sigma(T^{-1}Q)=0$ then 
 $$\rho\otimes\sigma(T_-)=\rho\otimes\sigma(T), \quad\quad \rho\otimes\sigma((T_-)^{-1})=\rho\otimes\sigma(T^{-1}).$$
   Moreover,
  $$(F_T)_{\rho, \sigma}(F_T)_{\rho, \sigma}^*=1, \quad\quad (G_T)_{\rho, \sigma}^*(G_T)_{\rho, \sigma}=1,\quad\quad (G_T)_{\rho, \sigma}=(F_T)_{\rho, \sigma}^*.$$
 Thus on the tensor product space of $\rho\otimes\sigma$ the twisted Hermitian form  and $\overline{R}$-matrix  are trivial $\Omega_T=\Delta_T(I)=\overline{R}_T$. It follows  that  $$R_T=w\otimes w\Delta_T(w^{-1}), \quad\quad \Delta_T(a)^*=\Delta_T(a^*)$$ on this space.
  
 \item[{\rm c)}] If the assumptions in b) hold for any pair of irreducible $^*$-representations $\rho$, $\sigma$
 of $A$ then
 $$T=T_-, \quad\quad T^{-1}=(T_-)^{-1}, \quad\quad \Delta_T(I)=\Delta_T(I)_-.$$ Moreover the twisted structure $(F_T, G_T)$ is strongly unitary. Thus the $R$-matrix $R_T$ and the Hermitian form $\Omega_T$ simplify further as in (\ref{data4}). 
In particular, the coproduct $\Delta_T$ of $A_T$ is cocommutative, $\Delta_T=\Delta_T^{\rm op}$.
 \end{itemize}

\end{thm}
\bigskip

\begin{proof} a) We have
 $T^{-1}({T^{-1}})^*_{21}=\Delta(I)$,  $T^*T_{21}=\Delta(I)^*.$
 Let $(F, G)$ be the weak quasi-tensor structure  defining the forgetful functor of $A$.
 Then 
 by  Theorem \ref{TK_unitary_ribbon}   equations 
(\ref{unitary_ribbon_functor1}) and (\ref{unitary_ribbon_functor2}) hold  for $(F, G)$.
Let $(F_T, G_T)$  be the new weak quasi-tensor structure obtained from the twist $T$,
$F_T=FT^{-1}$, $G_T=TG$.
We have
$$F_T\Sigma F_T^*=FT^{-1}\Sigma (T^{-1})^*F^*=FT^{-1}({T^{-1}})^*_{21}\Sigma F^*=F(GF)\Sigma F^*=F\Sigma F^*$$
and similarly
$$G_T^*\Sigma G_T=G^*T^*T_{21}\Sigma G=G^*(GF)^*\Sigma G=G^*\Sigma G.$$
It follows that 
equations 
(\ref{unitary_ribbon_functor1}) and (\ref{unitary_ribbon_functor2}) hold    for
 $(F_T, G_T)$. The twisted $R$-matrix $R_T$ induces a unitary braided symmetry in ${\rm Rep}^+(A_T)$ by
 Remark \ref{invariance_unitarity_of_braiding_under_twisting}. Moreover the twisted associator of $A_T$ is unitary
 ${\rm Rep}^+(A_T)$ by invariance of axioms of $\Omega$-involution under twisting.
It follows from   Theorem \ref{TK_unitary_ribbon} again that $A_T$ is a unitary coboundary 
 weak quasi-Hopf algebra. 
   It follows from Prop. \ref{R-canonicity} c) that
 $$\overline{R}_T=T_{21}\overline{R}T^{-1}=T_{21}T^*(P-Q)TT^{-1}=T_{21}(T^{-1})_{21}(P-Q)TT^{-1}=
 \Delta_T^{\rm op}(I)(P-Q)\Delta_T(I).$$
 We also have $\Delta_T(I)^*=(T^{-1})^*T^*=T_{21}T_{21}^{-1}=\Delta_T^{\rm op}(I)$
    thus $A_T$ has a compatible $^*$-involution.
 This is also equivalent to $\overline{R}_T=\Omega_T$. The formula for $R_T$ follows from the definition of $\overline{R}$
in (\ref{overline_R_defn}) for a general Hermitian coboundary wqh.
b) In this case we have $\rho\otimes\sigma(T_{-})=\rho\otimes\sigma(PT)=\rho\otimes\sigma((P+Q)T)=
\rho\otimes\sigma(T)$. 
In a similar way, $\rho\otimes\sigma(T_{-}^{-1})=\rho\otimes\sigma(T^{-1})$.
It follows that
$(F_T)_{\rho, \sigma} (F_T)_{\rho, \sigma}^*=F_{\rho, \sigma}\rho\otimes\sigma(T^{-1}(T^{-1})^*)
F_{\rho, \sigma}^*=F_{\rho, \sigma}\rho\otimes\sigma(\overline{R})^{-1}F_{\rho, \sigma}^*=F_{\rho, \sigma}G_{\rho, \sigma}G_{\rho, \sigma}^*F_{\rho, \sigma}^*=1$. 
One similarly shows that  $(G_T)_{\rho, \sigma}^*(G_T)_{\rho, \sigma}=1$. 
The equality $(G_T)_{\rho, \sigma}^*=(F_T)_{\rho, \sigma}$ follows from Prop. \ref{strongly_unitary_prop}.
Since $\Delta_T$ is induced by $(F_T, G_T)$ on the tensor product space, it follows that  $\Delta_T$ is $^*$-invariant
on this space, and also that $\Omega_T$ acts trivially.
c) This follows from b) and Tannaka-Krein duality.

\end{proof}

 We have seen that general Hermitian (unitary) coboundary wqh are Hermitian (ribbon) if and only if the elements 
 $\alpha$, $\beta$ of the antipode satisfy $\beta=\alpha^*$, and that it suffices that the antipode be of Kac type,
 see Theorem  \ref{Hermitian_category2}.
 We next ask when in addition we have a (unitary) modular fusion category.

\begin{thm}\label{modularity}
Let $A$ be a discrete unitary coboundary wqh algebra with compatible $^*$-involution  satisfying all the assumptions in the statement and in part c)
of theorem \ref{Drinfeld_Kohno} (e.g. $A$ has strongly trivial involution, see Def. \ref{trivial_involution}).  Let $T$ be the twist defined in  Theorem \ref{Drinfeld_Kohno}. Then $T$ induces 
a canonical unitary equivalence of ribbon categories
 ${\rm Rep}^+(A)\to {\rm Rep}^+(A_T)$. Moreover, if     $A$ is semisimple then ${\rm Rep}(A)$ is modular if and only if ${\rm Rep}(A_T)$ is modular.
   If $A$
  has   antipode defined by elements $\alpha$ and $\beta$ such that $\beta=\alpha^*$ (e.g. of Kac type)
then   ${\rm Rep}^+(A)$ is a unitary modular fusion category if and only if so is ${\rm Rep}^+(A_T)$.
\end{thm}

\section{Wenzl's  unitary structure of ${\mathcal C}({\mathfrak g}, q, \ell)$, square root of the  quantum Casimir}\label{unitary_structure_of_fusion_category}

Recall that the algebra $U_q({\mathfrak g})$ at complex roots of unity was introduced in Sect. \ref{73}, and we assume the same setting as there.
In particular, it  becomes
a (topological) ribbon complex Hopf algebra with a $^*$--involution.

For the reader convenience we recall two main results from \cite{Wenzl} concerning the unitary structure.
They  center on a continuity  argument that plays an important role in our paper in connecting
structures from quantum groups at roots of unity to corresponding structures from affine Lie algebras,
  developed in Sects. \ref{16}--\ref{21}. To this aim, as we shall see, it is important that the root of unity
  is minimal of large enough order in the sense of Def. \ref{large_enough}, \ref{minimal_root}.
  
  \begin{rem}\label{non_minimal_roots}
  Note that for non-minimal roots, the continuity argument may not be applied, but  the construction of this section
  possibly extend to these cases and lead to Hermitian structure. This may have possible applications to the study of Hermitian ribbon categories or applications to Chern-Simons theory, see Remark 8 in \cite{Sawin}.
  \end{rem} \bigskip

\subsection{Unitary structure of irreducible modules.}\label{25.1} The first result  describes unitarizability of a fixed  specialized Weyl module
 $V_\lambda(q)$   of the Hopf algebra $U_q({\mathfrak g})$  when $q$ varies in an arc of the unit circle depending on $\lambda$.
 Recall that if $q$ is any complex root of unity, $\ell'$ denotes the order of $q$,   $\ell$ the order of $q^2$,
$\Lambda^+$ the cone of dominant integral weights, $\Lambda^+(q)$ the open Weyl alcove associated to $q$,
$\Lambda^+(q)=\{\lambda\in\Lambda^+: \langle\lambda+\rho, \theta_0\rangle<\ell\}$,
where $\theta_0$ is the highest root $\theta$ if $d$ divides $\ell$ and the highest short root $\theta_s$ otherwise.
Forthermore, $\overline{\Lambda^+(q)}$ is  contained in a fundamental domain for the action of the affine Weyl group $W_{\ell'}$.
To cover the case where $q\in{\mathbb T}$ is not a root of unity, we set $\ell'=\infty$, $\Lambda_{\ell'=\infty}=\Lambda^+$
  and $W_{\ell'}$ is the ordinary Weyl group.

\begin{thm}\label{Wenzl_positivity} (\cite{Wenzl} Proposition 2.4) Let $\lambda$ be a dominant weight and 
$$I=\{q=e^{i\pi t}: |t|<\frac{1}{m-d}\}$$
where $m=\langle\lambda+\rho, \theta\rangle$. Then $V_\lambda(q)$ is simple and the Hermitian form defined in Prop. 2.2 in \cite{Wenzl} is positive definite for $q\in I$.
\end{thm}

\begin{proof}
We claim that $\lambda\in\overline{\Lambda^+(q)}$ for $q\in  I$, that is $\langle\lambda+\rho, \theta_0\rangle\leq\ell$.
We may assume that $\ell'$ is finite. Note that $|t|<\frac{1}{m-d}$ implies $\ell>m-d$ as $2\pi r |t|<2\pi$ for $r\leq m-d$.
We have the following cases.

a) $d=1$; In this case $\theta_0=\theta$, and the claim is $m\leq\ell$, but   we already know this.

b) $1<d$ and $d$ divides $\ell$; In this case $\theta_0=\theta$ and we need to show $m\leq\ell$.
As pointed out in \cite{Wenzl}, $\langle \mu, \theta\rangle$ is divisible by $d$ for any dominant weight $\mu$, thus $m-d<\ell$ is equivalent to $m\leq\ell$,

c) $2=d$ and $\ell$ is not divisible by $d$; we have $\theta_0=\theta_s$ and we need to show $\langle \lambda+\rho, \theta_s\rangle<\ell$. We have $\langle \lambda+\rho, \theta_s\rangle=m-\langle(\lambda+\rho, \theta-\theta_s\rangle<\ell+d-\langle\lambda+\rho, \theta-\theta_s\rangle$, but $\langle\lambda+\rho, \theta_s\rangle<\langle\lambda+\rho, \theta\rangle$ are both integers thus 
$\langle\lambda+\rho, \theta_s\rangle<\ell+1$ in this case,

d) $3=d$ and $\ell$ not divisible by $d$; we have $\theta_0=\theta_s$ and we verify $\langle \lambda+\rho, \theta_s\rangle\leq\ell$. We also have ${\mathfrak g}=G_2$, $\langle\lambda+\rho, \theta-\theta_s\rangle=4$ and we conclude as in the previous case, and the claim is proved.
Assume that  $V_\lambda(q)$ is not simple.  It follows it admits a submodule of highest weight $\mu<\lambda$, thus $\lambda$ and $\mu$   lie in $ \overline{\Lambda^+(q)}$ that is a fundamental domain for the translated action
of  $W_{\ell'}$, thus these dominant weights can not be conjugate under this action, contradicting the linkage principle.
The Hermitian form $(\xi, \eta)$ of $V_{\lambda, {\mathcal A}'}$ as in the statement is Hermitian over the Laurent polynomial ring ${\mathcal A}'$ by Prop. 2.3 in
\cite{Wenzl} and thus it becomes an Hermitian form for $V_\lambda(q)$ specializing $x$ to any complex number $q$ with $|q|=1$. Moreover it is not trivial as $(v_\lambda(q), v_\lambda(q)>0$ and it makes $V_\lambda(q)$ into a $^*$-representation, and it follows that the radical of the Hermitian form is a $^*$-subrepresentation, that vanishes by simplicity, and it follows that the form is non-degenerate. Let $v_i(q)$ be the specialization of Kashiwara-Lusztig basis at $q$ and consider the matrix with entries $(v_i(q), v_j(q))$, which depend continuously on $q$, and therefore 
the same holds for the eigenvalues. The $^*$-involution of the Hopf algebra $U({\mathfrak g})$ for
$q=1$ (classical limit)
corresponds to the compact real form ${\mathfrak g}_{\mathbb R}$ of ${\mathfrak g}$, and the Hermitian form
of the Weyl module $V_\lambda$ of $U({\mathfrak g})$ is a 
  unitary representation of the compact Lie subgroup $K$  of $G$,
see Remark \ref{compact_real_form}. Hence $((v_i(q), v_j(q)))$ is positive for all $q\in I$.
\end{proof}

\begin{rem}
Note that by the proof of Lemma 2.2 in \cite{Wenzl},  the Hermitian form of the statement is the specialization at $q$ of an Hermitian form 
on the Weyl module $V_{\lambda, {{\mathcal A}'}}$ of $U_{{\mathcal A}'}^\dag({\mathfrak g})$ satisfying the assumptions in b) of Prop. \ref{conjugates}.

\end{rem}

We summarize some of the results recalled in Sect. \ref{73} and in  proof of the previous theorem  in the following corollary that plays an important role in  \cite{Wenzl} and also in our paper.
Recall that minimal roots of large enough order have been defined in Def. \ref{large_enough}, \ref{minimal_root}.

\begin{cor}\label{corollary_of_positivity}
Let  $q_0$ be a  minimal root of unity 
of   large enough order. Then for all $q$ varying on the arc ${\mathbb T}_{q_0, 1}$ of the unit circle clockwise connecting  
 $q_0$ to $1$, we have that $\overline{\Lambda^+(q_0)}\subset \overline{\Lambda^+(q)}$ and the latter equals the whole
 Weyl chamber $\Lambda^+$ for $q=1$;   for   $\lambda\in\overline{ \Lambda^+(q_0)}$, and $q$ varying continuously in 
 ${\mathbb T}_{q_0, 1}$,
 the Kashiwara-Lusztig basis of $V_{\lambda, {\mathcal A}'}$ specializes to a basis of $V_\lambda(q)$;  the  invariant Hermitian form of $V_\lambda(q)$ is positive definite and is uniquely determined by $(v_\lambda(q), v_\lambda(q))=1$ and is invariant with respect to the compact real form of $U({\mathfrak g})$ for $q=1$.

\end{cor}
\bigskip

We next recall the unitary structure on certain tensor products of specialized Weyl modules of $U_q({\mathfrak g})$.\bigskip

\subsection{Square root of the quantum Casimir operator and  compatible Hermitian coboundary structure of $U_q({\mathfrak g})$.}\label{25.2}
Loosely speaking, the strategy of \cite{Wenzl}  is to define   a Hermitian structure on tensor products $V_\lambda(q_0)\otimes V_\mu(q_0)$ for $\lambda$, $\mu\in\overline{\Lambda^+(q_0)}$
via the action of a    coboundary matrix of $U_{q_0}({\mathfrak g})$ on the space. By \cite{Sawin}, $U_{q_0}({\mathfrak g})$ is a ribbon Hopf algebra in a topological sense. However the ribbon element $v$ may not have a square root in
$U_{q_0}({\mathfrak g})$  arising  as the specialization  at $q_0$ of an element in the integral form
${\mathcal U}_{{\mathcal A}'}^\dag({\mathfrak g})$. This problem is mentioned
at the beginning of Sect. 3.6 in   \cite{Wenzl}.

  Indeed, if this were the case then following a central idea in
\cite{Wenzl},
by our   theorem \ref{unitary_ribbon2} we would have an Hermitian form on the full tensor
product obtained by specialization of a nondegenerate  ${\mathcal A}'$-valued Hermitian form.
Since it is obtained as a specialization,   the Hermitian form  would be continuous in $q$, and thus positive since it reduces to the
usual inner product for $q=1$. It would then follow that the tensor product is completely reducible, but this is not always the case.

  For our purposes we shall need an Hermitian form defined on every tensor product space $V_\lambda(q_0)\otimes V_\mu(q_0)$ for $\lambda$, $\mu\in\Lambda^+(q_0)$, thus we look for a square root operator of the action of $R_{21}R$. Using the ribbon structure given by the specialized quantum Casimir $v$ see Def. \ref{quantum_casimir}, Theorem \ref{ribbon_specialized_case},
    Sect. \ref{73}, \ref{74}, \ref{unitary_structure_of_fusion_category},
   we are reduced to construct a square root of the action of $v$ on a full tensor product.

\begin{defn} For an invertible operator $T$ on a finite dimensional vector space, we define $T^{1/2}$ via Jordan decomposition, that is
for a Jordan block $J=c(I+N)$ of $T$ with $c\in{\mathbb C}$, $c\neq 0$ and $N$ nilpotent we set $J^{1/2}=c^{1/2}(1+N)^{1/2}$, where $(1+N)^{1/2}$ is defined via Taylor expansion (that is eventually constant since $N$ is nilpotent)
and for $c=|c| e^{2\pi it}\in{\mathbb C}$, $t\in(-1, 1]$ we set $c^{1/2}=|c|^{1/2} e^{\pi it}$, $|c|^{1/2}>0$.

\end{defn}

It follows from the inclusions (\ref{direct_product_specialized_form}), (\ref{the_coproduct_range_specialized_form})
that for any invertible elements $T\in U_{q_0}({\mathfrak g})$, the square roots $T^{1/2}$ and $\Delta(T)^{1/2}$ are well
defined as operators of $\Pi_{V(q_0)} {\rm End}(V(q_0))$ and $\Pi_{V(q_0), W(q_0)} {\rm End}(V(q_0))\otimes {\rm End}(W(q_0)) $
respectively. 

\begin{defn}\label{square_root} In particular, let $v$ be the specialized quantum Casimir operator of $U_{q_0}({\mathfrak g})$
as in Def. \ref{quantum_casimir}.
Then $v^{1/2}$ and $\Delta(v)^{1/2}$   are well defined invertible operators of 
$\Pi_{V(q_0)} {\rm End}(V(q_0))$ and $\Pi_{V(q_0), W(q_0)} {\rm End}(V(q_0))\otimes {\rm End}(W(q_0)) $
respectively. 

\end{defn}

 \begin{prop}\label{Wenzl_Hermitian_form}  We have  for  $\lambda$, $\mu\in{\overline{\Lambda^+(q_0)}}$,
   \begin{itemize}
   \item[{\rm        a)}] for $T\in U_{q_0}({\mathfrak g})$ invertible, $R\Delta(T)^{1/2}=\Delta^{\rm op}(T)^{1/2}R$,
\item[{\rm        b)}] 
  $v^{1/2}$
 acts on $V_\lambda(q_0)$ as a   scalar in ${\mathbb T}$, given by $w_\lambda:={q_0}^{-\frac{\langle\lambda, \lambda+2\rho\rangle}{2}}$,
 \item[{\rm        c)}] 
  $((R_{21}R)^{-1})^{1/2}$ acts on $V_\lambda(q_0)\otimes V_\mu(q_0)$ as $\Delta(v)^{1/2}w_\lambda^{-1}\otimes w_\mu^{-1}$,
  \item[{\rm        d)}]   the operator 
 $\overline{R}_{\lambda, \mu}=R((R_{21}R)^{-1})^{1/2}$ is invertible selfadjoint on $V_\lambda(q_0)\otimes V_\mu(q_0)$
 w.r.t. the product of Wenzl inner products and satisfies $(\overline{R}_{\mu, \lambda})_{21}\overline{R}_{\lambda, \mu}=1$.
   \end{itemize}
  \end{prop}

 \begin{proof} The proof of a), b), c) follow easily from the definitions and Prop.  \ref{ribbon_specialized_case}. d) Selfadjointess of $\overline{R}_{\lambda, \mu}$
 can be shown as in the proof of Theorem \ref{unitary_ribbon2} with the following modifications. In equations (\ref{comp1})--(\ref{comp3}) we put $E=I$, replace $\Delta(w)$ with $\Delta(v)^{1/2}$, $\Delta^{\rm op}(w)$  with $\Delta^{\rm op}(v)^{1/2}$, $\tilde{\Delta}(w)$
 with $\tilde{\Delta}(v)^{1/2}$, and in (\ref{comp2}) we use a). In (\ref{comp4}) we use $(\tilde{\Delta}(v)^{1/2})^*=
 (\tilde{\Delta}(v)^*)^{1/2}$ (the coefficients of the Taylor expansion are real), in (\ref{comp4}) we replace $\Delta(w^{-1})$ with $\Delta(v^{-1})^{1/2}$. The coboundary equation
 $(\overline{R}_{\mu, \lambda})_{21}\overline{R}_{\lambda, \mu}=1$ can be shown as in Prop. \ref{overline_R_prop} with similar 
 replacements of $\Delta(w)$ and $\Delta^{\rm op}(w)$ with $\Delta(v)^{1/2}$ and $\Delta^{\rm op}(v)^{1/2}$ respectively.
 
 \end{proof}
 
  \begin{rem} Thus $\overline{R}_{\lambda, \mu}$ defines a non-degenerate Hermitian form on the full tensor product space
 $V_\lambda(q_0)\otimes V_\mu(q_0)$ making the braided symmetry unitary, cf. the
 problem mentioned  
at the beginning of Sect. 3.6 in   \cite{Wenzl}. On the other hand, not all the idempotents
$p_\gamma:   V_{\lambda, {\mathcal A}'}\otimes_{{\mathcal A}'} V_{\mu, {\mathcal A}'}\to V_{\gamma, {\mathcal A}'}$ describing classical fusion
    specialize to corresponding idempotents $p_\gamma(q)$ by non-semisimplicity. By \cite{Wenzl}, this is an obstruction in  the proof
    of positivity of the Hermitian form, in that when this is possible then the Hermitian form of the specialized
    module $V_\lambda(q_0)\otimes V_\mu(q_0)$ is positive on $V_\gamma(q_0)$ by a continuity argument that
    links $q_0$ to the classical limit $1$ that  comes from specialization (unlike the construction of Prop. \ref{Wenzl_Hermitian_form}).
    Notably, there are sufficiently many modules for which this holds, and in a nutshell from them we shall construct the wqh.
    \end{rem}

We next show that $U_{q_0}({\mathfrak g})$ satisfies the axioms of a Hermitian coboundary Hopf algebra.  With this we mean that we show that axioms
a), b), c), d) of Def. \ref{Hermitian_ribbon_wqh} hold for the topological Hopf algebra $U_{q_0}({\mathfrak g})$  defined by \cite{Sawin}.  
The   result   is motivating
for the construction of the unitary discrete coboundary weak  Hopf algebras $A_W({\mathfrak g}, q_0, \ell)$ of the next section. Moreover, in Sect. \ref{20}, see Theorem \ref{main_wh}, and more in detail Prop. \ref{$2$-cocycle property} and in Sect.
\ref{21} we shall need the $2$-cocycle property of the coboundary $\overline{R}=\overline{R}_U$   of $U_{q_0}({\mathfrak g})$  for the proof of our analogue of KL-F theorem for affine vertex operator algebras \ref{Finkelberg_HL}.

We need at least to show that the square root $v^{1/2}$ of the ribbon element $v$ lies  in $U_{q_0}({\mathfrak g})$ and that $\Delta(v^{1/2})=\Delta(v)^{1/2}$. We shall see that combining with the work in \cite{Wenzl}, this is all
is left to show.
Note that the $R$-matrix $R$ and the ribbon element $v\in U_{q_0}({\mathfrak g})$ depend only on the choice of 
$q_0^{1/L}$, see Sect. 1 in \cite{Sawin},
Sect. 1.4 in \cite{Wenzl}. 
\medskip

\begin{thm}\label{U_q_as_a_Hermitian_ribbon_h} Let $q_0^{1/L}$ be a fixed $L$-th root of $q_0$ and consider the associated ribbon structure $(R, v)$ on  $U_{q_0}({\mathfrak g})$ defined as in Sect. \ref{73}, (i.e. following \cite{Sawin} for the algebraic structure and \cite{Wenzl} for the $^*$-involution).
Let $v^{1/2}$ and $\Delta(v)^{1/2}$ be defined 
   as in Def. \ref{square_root}. 
Then 
 \begin{itemize}
   \item[{\rm        a)}] 
$v^{1/2}\in U_{q_0}({\mathfrak g})$ and  $\Delta(v)^{1/2}\in U_{q_0}({\mathfrak g})\overline{\otimes} U_{q_0}({\mathfrak g})$, the completed tensor product defined in \cite{Sawin}, and $\Delta(v^{1/2})=\Delta(v)^{1/2}$.  
   \item[{\rm        b)}] with $w=v^{1/2}$,  $U_{q_0}({\mathfrak g})$ becomes a  (topological) 
 Hermitian   coboundary Hopf algebra   with compatible involution and   antipode of Kac type.
    \item[{\rm        c)}]  In particular, $\overline{R}^U:=R^U\Delta(w)w^{-1}\otimes w^{-1}$ is a $2$-cocycle for  $U_{q_0}({\mathfrak g})$,
    $$\overline{R}^U\otimes I\Delta\otimes1(\overline{R}^U)=I\otimes \overline{R}^U1\otimes\Delta(\overline{R}^U).$$

 \end{itemize}

 \end{thm}
 
 \begin{proof} The Kac-type property of the antipode follows from 
 properties (\ref{Delta})--(\ref{counit_antipode}) that still hold for 
$U_{q_0}({\mathfrak g})$.
Axioms a) b), d) of Def. \ref{Hermitian_ribbon_wqh} are shown in
 Lemma 1.4.1 of \cite{Wenzl} in the case of the algebra over a formal variable $x$, and easily extend to the specialized algebra $U_{q_0}({\mathfrak g})$ of \cite{Sawin}.
 For axiom c), we set $w=v^{1/2}$ defined on each full matrix algebra $M$ defining $U_{q_0}({\mathfrak g})$ as in Def. \ref{square_root} via the Jordan form of $v$. If $c(I+N)$ is a Jordan block for $v$ then $c(I\otimes I+\Delta(N))$
 is a Jordan block for $\Delta(v)$, and this puts $\Delta(v)$ in Jordan form. Thus $v^{1/2}$ and $\Delta(v)^{1/2}$ are   
 eventually constant limit on the full matrix algebras defining $U_{q_0}({\mathfrak g})$ and $U_{q_0}({\mathfrak g})\overline{\otimes} U_{q_0}({\mathfrak g})$ respectively.
 This shows that $v^{1/2}$ and $\Delta(v)^{1/2}$ lie in $U_{q_0}({\mathfrak g})$ and $U_{q_0}({\mathfrak g})\overline{\otimes} U_{q_0}({\mathfrak g})$ respectively.  Since $v$ is central, the nilpotent parts occurring in $v$ are central, thus $v^{1/2}$ is central.  The equality $\Delta(v^{1/2})=\Delta(v)^{1/2}$ follows
 from the homomorphism property of $\Delta$.  \end{proof}

     The previous result was  implicitly used in \cite{CP}  but no detailed proof was given.
In that paper   two of us constructed weak  Hopf algebras associated
to Verlinde fusion categories in the type $A$ case with specific methods. In Sect. \ref{20}  we shall extend and expand the main result of \cite{CP}  to the other Lie types   with general      methods.

    \begin{prop}\label{selfadjointness_three_coboundary}
    The elements $\overline{R}^U$
    \begin{equation}\label{2-cocycle_equation}
  \overline{R}^U\otimes I\Delta\otimes1(\overline{R}^U)=I\otimes \overline{R}^U1\otimes\Delta(\overline{R}^U)
     \end{equation}
are selfadjoint in 
$U_{q_0}({\mathfrak g})\overline{\otimes} U_{q_0}({\mathfrak g})$ and
$U_{q_0}({\mathfrak g})\overline{\otimes} U_{q_0}({\mathfrak g})\overline{\otimes} U_{q_0}({\mathfrak g})$ respectively.
    
    \end{prop}
    
    \begin{proof} Selfadjointness of  $\overline{R}^U$ follows from part d) of Prop. \ref{Wenzl_Hermitian_form}.
Recall that $\overline{R}^U$ twists $\Delta$ to $\Delta^{\rm op}$ and that $\Delta(a)^*=\Delta^{\rm op}(a^*)$ for 
$a\in U_{q_0}({\mathfrak g})$.
We have
    $$(\overline{R}^U\otimes I\Delta\otimes1(\overline{R}^U))^*=\Delta^{\rm op}\otimes1(\overline{R}^U)\overline{R}^U\otimes I=\overline{R}^U\otimes I\Delta\otimes1(\overline{R}^U).$$
    \end{proof}
    
    \begin{rem}\label{remark_sect_25_on_the _2_cocycle_property}
 Note that the ribbon element of
    $U_{q_0}({\mathfrak g})$  defined in detail in
 \cite{Sawin}   acts on any irreducible Weyl module $V_\lambda(q_0)$
    as ${q_0}^{-{\langle\lambda, \lambda+2\rho\rangle}}$, but this may not suffice to determine an expression of
    the ribbon element in $U_{q_0}({\mathfrak g})$.  Moreover, there is no unitary square root of 
    the ribbon element  in $U_{q_0}({\mathfrak g})$ which is continuous in $q$ from $1$ to $q_0$ as otherwise 
    by Theorem \ref{unitary_ribbon2} this induces an Hermitian form  on any full tensor product 
    of Weyl modules and by continuity in $1$ by Wenzl argument it is positive definite.
    Thus the tensor product would be completely reducible, but this is not always the case.

    On the other hand, when a Wenzl idempotent $p_\gamma: V_\lambda(q_0)\otimes V_\mu(q_0)\to V_\gamma(q_0)$ onto an irreducible Weyl module  $V_\gamma(q_0)$ 
is well-defined, a unitary square root $w$ of the action of $v$  can be defined as in part b) of Prop.
\ref{Wenzl_Hermitian_form} on  $V_\gamma(q_0)$ and is a specialization at $q_0$ of  a corresponding action in the formal variable $x$. In particular, this holds for $\lambda$, $\mu$, $\gamma$ in the open Weyl alcove $\Lambda^+(q_0)$.

     \end{rem}

    \bigskip

\subsection{Hermitian form on tensor product modules (following Wenzl).}\label{25.3}

 We next recall a second main result in \cite{Wenzl} on the construction of an Hermitian form on certain tensor products of specialized Weyl modules.
 When we specialize $x$ to a   root of unity $q$, not all the specialized Weyl modules $V_\lambda(q)$ are irreducible, and  not all
the idempotents $p_\gamma:   V_{\lambda, {\mathcal A}'}\otimes_{{\mathcal A}'} V_{\mu, {\mathcal A}'}\to V_{\gamma, {\mathcal A}'}$ describing classical fusion
    specialize to corresponding idempotents $p_\gamma(q)$. 
    
    For example, when   $\lambda$, $\mu\in\Lambda^+(q)$
    then   $V_\lambda(q)\otimes V_\mu(q)$  is tilting and  the decomposition into indecomposable tilting modules
    depends in general on the affine Weyl group associated to $q$ and thus on the order of $q^2$.  When $q_0$ is a minimal root 
    of large enough order (Def. \ref{minimal_root}, \ref{large_enough}) there are special
      cases   where  $\lambda$, $\mu$, $\gamma$  are  fixed    dominant weights of   $\overline{\Lambda^+(q_0)}$ when
      we let  $q$ varying continuously
    from $q_0$ to $1$ as in Corollary \ref{corollary_of_positivity}, and we regard $\lambda$, $\mu$, $\gamma$  as   elements of
    $\overline{\Lambda^+(q)}$ (that contains $\overline{\Lambda^+(q_0)}$   by the same Corollary).
    
    In these cases   $p_\gamma(q)$ can be constructed as an  idempotent onto the isotypic component corresponding to $\gamma$ arising from
    specialization in $q$ of a corresponding idempotent between modules of the integral form. Assume that we are in one of these cases.
    On the range of $p_\gamma(q)$,    $\Delta(v)^{1/2} v^{-1/2}\otimes v^{-1/2}$ acts as as a square matrix  
    with respect to the product of the two Kashiwara-Lusztig bases with entries powers
    in the specified   power $q^{1/L}=e^{2\pi i t/L}$
    (the weights of the tensor product are congruent to $\lambda+\mu$, 
    and the action of $v$ on each weight space is given by the quantum Casimir, part d) of Theorem \ref{Properties_of_R},
    thus $\Delta(v)^{1/2} v^{-1/2}\otimes v^{-1/2}$ acts diagonally on each weight space with eigenvalues integral powers of $q^{1/L}$). Note in particular that this action is continuous. 
    
   Moreover, $R(q)$ acts as a rectangular matrix, as $R(q)$ may not leave the range of $p_\gamma(q)$ invariant.
    It follows that $\overline{R}_{\lambda, \mu}(q)$ defines a non-degenerate Hermitian form on the product Hilbert space making the ranges
    pairwise orthogonal if we in addition know that  $V_\lambda(q)\otimes V_\mu(q)$ is completely reducible
    into irreducible Weyl modules as in the classical case. Thus in this case  $p_\gamma(q)$ is selfadjoint with respect to
    the Hermitian form of the domain, and it follows that the form is nondegenerate on the range of $p_\gamma(q)$,
    and therefore positive by the continuity argument.
    The selfadjointness property reads as
\begin{equation}p_\gamma(q)^\dag\overline{R}_{\lambda, \mu}(q)=\overline{R}_{\lambda, \mu}(q)p_\gamma(q),\end{equation}
    where $p_\gamma(q)^\dag$ is the adjoint of $p_\gamma(q)$ with respect to the tensor product Hilbert space structure.

 More precisely, by Sect. 3.5 in \cite{Wenzl},
 we choose  a fundamental   representation $V$ of the Lie algebra   for each Lie type. Let $\kappa$ be its highest weight.
 This representation can be defined explicitly
 as follows (we follow \cite{Humphreys} for notation of the fundamental weights here below):
 
 $A_N$) the vector module of $\mathfrak{sl}_{N+1}$, $\kappa=\lambda_1$,\smallskip
 
 $B_N$) the spin module of  ${\mathfrak o}_{2N+1}$, $\kappa= \lambda_N$,\smallskip
 
 $C_N$) the vector module of $\mathfrak{sp}_{2N}$, $\kappa=\lambda_1$,\smallskip
 
 $D_N$) the two spin modules of ${\mathfrak o}_{2N}$, $\kappa=\lambda_{N-1}$, $\lambda_N$,\smallskip
 
 $E_6$) $\kappa$ one among $\lambda_1$ or $\lambda_6$ the corresponding vertex in the Dynkin diagram is
chosen between the two farthest from the branching point,\smallskip

$E_7$, $E_8$) $\kappa=\lambda_7$, $\lambda_8$ respectively,\smallskip

$F_4$) $\kappa=\lambda_4$,\smallskip

$G_2$) $\kappa=\lambda_1$.\bigskip

We shall need the following properties of $V$. 
As anticipated in Subsect. \ref{18.1} Wenzl \cite{Wenzl} notes the following useful fact.

\begin{thm}\label{tilting_as_fundamental}
Suppose that $\kappa$ lies in the open principal Weyl alcove $\Lambda^+(q)$.
 Then a $U({\mathfrak g})$-module
  is tilting if and only if it is a direct sum of direct summands of tensor
powers of $V(q)$.
\end{thm}

Recall that in Subsect. \ref{18.5} we introduced the open Weyl alcove $\Lambda^+(q)$, and recalled
that it labels the irreducible objects of ${\mathcal C}({\mathfrak g}, q, \ell)$.
We here recall that the   Coxeter number $h$ and its dual $h^\vee$
may are also computed as
\begin{equation}\label{Coxeter_numbers} h^\vee=\langle\rho, \theta^\vee\rangle+1, \quad h=\langle\rho, \theta_s^\vee\rangle+1,\end{equation}
with as before $\theta$ the highest root and $\theta_s$ the highest short root.
Recall also from Def. \ref{large_enough} the level $k$ associated to the order $\ell$ of $q^2$.
We introduce 
$$h(\lambda):=\langle\lambda, \theta^\vee\rangle \quad \text{if } d |\ell, \quad \quad h(\lambda)=\langle\lambda, \theta_s\rangle  \quad \text{otherwise}.$$
The following result (mentioned in Subsect. \ref{18.1})
is simple but useful, especially for the use of Theorem \ref{tilting_as_fundamental},   to determine
the minimum possible value of the level $k$ for which the above Theorem applies.

\begin{prop}\label{fundamental_alcove}
The dominant weight $\lambda$ lies in the Weyl alcove $\Lambda^+(q)$ if and only if $h(\lambda)\leq k$. Moreover
$h(\kappa)=1$ for the classical Lie types, $E_6$, $E_7$, and $h(\kappa)=2$ for $E_8$. For $F_4$ and $G_2$
 then $h(\kappa)=1$ if $d|\ell$, and $h(\kappa)=2$ otherwise.
\end{prop}

The fundamental representation $V$ is irreducible if ${\mathfrak g}$ is not of type $D$, and is the sum of the two half spin representations in the type $D$ case; 
  every   irreducible of ${\mathfrak g}$ is a subrepresentation of a power of $V$; the  dominant weight of $V$ (or of each summand in type $D$)
 lies in $\Lambda^+(q)$. 
 
 We consider the associated Weyl module of $U_q({\mathfrak g})$ denoted in the same way.
 For    $\lambda\in\Lambda^+(q_0)$,
 $V_\lambda(q)\otimes V(q)$ decomposes into a   direct sum of indecomposable tilting modules $T_\gamma$ 
 $$V_\lambda(q)\otimes V(q)=\oplus_\gamma 
 T_\gamma\otimes {\mathbb C}^{m_\gamma}$$ with the property that for ${\mathfrak g}\neq E_8$, the dominant weights $\gamma$ appearing in the decomposition 
  at most lie in   $\overline{\Lambda^+(q)}$, and for ${\mathfrak g}\neq F_4$, $E_8$ the decomposition is multiplicity-free
  in the classical case.
Thus  $T_\gamma=V_\gamma(q)$ for all $\gamma$ and the decomposition is completely reducibile and unique, thus $p_\gamma(q)$
  is defined.
  
 For ${\mathfrak g}=F_4$   multiplicity may arise
 for $\gamma=\lambda$. In this case, $p_\gamma(q)$ is first defined for $\gamma\neq\lambda$   then  $p_\lambda(q)$
 is defined as the complement idempotent.
 
 For ${\mathfrak g}=E_8$, the summand $T_\gamma$ may not lie in 
$ \overline{\Lambda^+(q)}$
 for $\gamma=\lambda+\kappa$  with $\kappa$ the dominant weight of $V$.  Multiplicity may arise  for $\gamma=\lambda$ in this case also. In Prop. 3.6 case 2 \cite{Wenzl} it is shown that a selfadjoint projection
 $p_\gamma(q)$ may be defined, and we   shall comment  more on this in the next section, and more precisely in the proof of Lemma \ref{k=1}.

These constructions
 hold for any primitive  root $q$ such that the order $\ell$ of $q^2$ is large enough (see Def 
 \ref{large_enough})
 and  under this assumption we shall
 construct the weak  Hopf algebra $A_W({\mathfrak g}, q, \ell)$ at an algebraic level in the first part of the following section.
 
 For the unitary structure, Wenzl shows that the restriction to the minimal roots Def. \ref{minimal_root} is essential. 
 The requirement of having large enough order includes   the fundamental representations in the open Weyl alcove
 for all Lie types and positive integer levels,  
see  Prop. 2.4 in \cite{Wenzl}, also Theorem \ref{Wenzl_positivity}.
 
\begin{defn}\label{Wenzl_idempotents} We refer to $p_\gamma$ (also denoted $ p_{\lambda, \gamma}$ in the next section)
\begin{equation}\label{Wenzl_idempotents2}   V_\lambda(q)\otimes V(q)\to V_\gamma(q)\otimes {\mathbb C}^{m_\gamma}, \quad\quad\lambda, \gamma\in\Lambda^+(q)\end{equation} as  
   Wenzl  idempotents.
   \end{defn}

   \begin{rem}
For the unitarization result of \cite{Wenzl},  one  only needs  $\lambda$,  $\gamma\in {\Lambda^+(q_0)}$   (but to
    construct   $p_\gamma$ for $\gamma=\lambda$ in the $F_4$ one also needs $p_\nu$ with $\nu\in\overline{ {\Lambda^+(q_0)}}$, 
    see the proof of the proposition at page 274 in \cite{Wenzl}). On the other hand, for our purposes (a search of strongly unitary structures), see Theorem
    \ref{Zhu_as_a_compatible_unitary_wqh}, we need the stronger positivity results of Lemma 3.6.2 (b) in \cite{Wenzl}
    that involve positivity of the Hermitian form also on the range of $p_\gamma$ for $\gamma\in\overline{\Lambda(q_0)}$
    for  ${\mathfrak g}\neq E_8$.
    \end{rem}

\bigskip

    \section{Compatible unitary coboundary weak  Hopf algebras $U_q({\mathfrak g})\to A_W({\mathfrak g}, q, \ell)$}\label{20}
    
    The aim of this section is to prove Theorem
\ref{main_wh} stated in Sect. \ref{5++}. The proof is divided in three parts. In Theorem \ref{equivalence_G_q_with_qg_fusion_category},
  we introduce a full strict tensor category ${\mathcal G}_q$ generated tensorially by the fundamental
  representation $V$ equivalent to
${\mathcal C}({\mathfrak g}, q, \ell)$. We define Wenzl functor $W:\tilde{\mathcal G}_q\to{\rm Vec}$ on an equivalent concrete realization $\tilde{\mathcal G}_q$ of ${\mathcal G}_q$, in Def. \ref{W_functor}.
 In Theorems \ref{first1}, \ref{second2} we construct a natural weak tensor structure on $W$ and study the algebraic and unitary aspects, respectively.
In particular, the algebraic Theorem  \ref{first1} considers more general roots of unity, and we hope it will be useful 
in other circumstances, discussed in \cite{Sawin}.

Our first aim is to construct  a natural functor $W: {\mathcal C}({\mathfrak  g}, q, \ell)\to{\rm Vec}$  associated to the same dimension function $D$ as in Sect. \ref{73},  and thus this functor  is a particular case of the former.

To do this, we   consider  the tensor structure of ${\mathcal C}({\mathfrak g}, q, \ell)$ ($q$ not necessarily minimal)
of \cite{Wenzl}. This gives rise to
the mentioned forgetful functor $W$,  and we shall construct a natural weak tensor structure   on $W$, and  in this way  
we have a canonical weak  Hopf algebra $A_W({\mathfrak g}, q, \ell)$.

When $q$ is a minimal root of unity,    the work of \cite{Wenzl} implies that $U_q({\mathfrak g})$
is a  Hermitian   coboundary Hopf algebra with compatible involution (in a topological sense),  see
  Theorem \ref{U_q_as_a_Hermitian_ribbon_h}.
This Hermitian structure
underlies the unitary structure of ${\mathcal C}^+({\mathfrak  g}, q, \ell)$. 

To ease notation we shall write $\overline{R}^U\Delta(I)$ in place of ${\pi}\otimes{\pi}(\overline{R}^U)\Delta(I)$, 
and we shall use a similar notation for other elements of $A_W({\mathfrak g}, q, \ell)\otimes A_W({\mathfrak g}, q, \ell)$
that come from elements of   $U_q({\mathfrak g})\overline{\otimes} U_q({\mathfrak g})$.

\begin{rem}\label{how_it_acts}
Note that by the presence of the factor $\Delta(I)$, the element $\Delta(v)^{1/2}$ involved in  
(\ref{the_formula_for_the_weak_coboundary}) is diagonal on the irreducible components, and acts
as ${q}^{-\frac{\langle\gamma, \gamma+2\rho\rangle}{2}}$ on an irreducible component of highest weight $\gamma$.
\end{rem}
 
The weak tensor structure is not unique but when we change it  then the weak  Hopf algebra changes only by a trivial twist.
The special case ${\mathfrak g}={\mathfrak sl}_N$  will be useful for the construction of tensor equivalences studied
in  Sect. \ref{KW}.  We also note that in this case we recover the example   constructed in \cite{CP} with a different method.

  \bigskip

  Recall  that the quotient category ${\mathcal C}({\mathfrak g}, q, \ell)$ was outlined in Sect. \ref{73}. We assume
  $\ell'<\infty$.

 \begin{rem}\label{G_q}
 By Lemma 1.1 in \cite{GK}, composition of  inclusion   ${\mathcal T}^0\to{\mathcal T}({\mathfrak g}, q, \ell)$ with   projection   
 ${\mathcal T}({\mathfrak g}, q, \ell)\to {\mathcal C}({\mathfrak g}, q, \ell)$ is an equivalence of linear categories.
 Hence ${\mathcal T}^0$ becomes a semisimple tensor
category   tensor equivalent to 
 ${\mathcal C}({\mathfrak g}, q, \ell)$. In the following subsection we shall construct among other things
 a specific tensor structure of a certain interesting equivalent full subcategory ${\mathcal G}_q\subset{\mathcal T}^0$ making ${\mathcal G}_q\to
 {\mathcal C}({\mathfrak g}, q, \ell)$ an equivalence of tensor categories.\end{rem}

\subsection{A specific construction, the weak  Hopf algebras  $A_W({\mathfrak g}, q, \ell)$.}\label{29.1} In this subsection
$q$ is any root of unity of order   in the sense of Def. \ref{large_enough}.
We obtain a   functor
  ${\mathcal C}({\mathfrak g}, q, \ell)\to {\rm Vec}$ together with a weak    tensor structure 
  $(F, G)$ associated to the same dimension function $D$ as in the previous subsection, and correspondingly a weak  Hopf algebra $A_W({\mathfrak g}, q, \ell)$. In the next subsection we consider the case where $q$ is a minimal root.
For this construction we mostly take into consideration ideas in \cite{Wenzl}  that we   review and extend to a general root of unity $q$ such that $\ell$ is  in the sense of Def. \ref{large_enough}. When $q$ is a minimal root, $A_W({\mathfrak g}, q, \ell)$ becomes a   unitary coboundary weak  Hopf algebra. To do this, as briefly anticipated in Remark \ref{G_q},
  we shall introduce a   linear category ${\mathcal G}_q$ of non-negligible tilting modules associated to a fundamental representation of ${\mathfrak g}$. This category appears implicitly in \cite{Wenzl}. In 
  \cite{CP} we have shown that ${\mathcal G}_q$ has a natural structure of a strict (ribbon) tensor category when
   ${\rm Vec}$ is regarded as strict and $q$ is a minimal root and is unitarily ribbon equivalent to
   ${\mathcal C}({\mathfrak g}, q, \ell)$. In this subsection we extend this to  roots of order  in the sense of Def. \ref{large_enough}
     and moreover we shall  define a functor $ {W}:{\mathcal G}_q\to{\rm Vec}$ and then introduce a weak tensor structure on $W$ that corresponds by Tannakian reconstruction to $A_W({\mathfrak g}, q, \ell)$.
     
 We keep the notation of the previous section. Let     
$$p_{\lambda, \gamma}: V_\lambda\otimes V\to V_\gamma\otimes {\mathbb C}^{m_\gamma}, \quad\quad\lambda, \gamma\in\Lambda^+(q)$$
denote Wenzl idempotents defined in (\ref{Wenzl_idempotents2}).

  \begin{defn}\label{projections_p_n}
 We define the projection $p_2:=\sum_{\gamma\in\Lambda^+(q)} p_{\kappa, \gamma}$   and   set $V\underline{\otimes}V=p_2V\otimes V$. 
 We  use $p_{\lambda, \gamma}$ to iteratively define projections $p_n: V^{\otimes n}\to V^{\underline{\otimes} n}$  onto the maximal non-negligible submodule   $V^{\underline{\otimes} n}$ induced by the decomposition of $V^{\underline{\otimes} n-1}\otimes V$.  
 \end{defn}
 
 \begin{rem}\label{canonical_iterative} By the iterative argument in the construction, every representation $V^{\underline{\otimes} n}$ has a canonical decomposition into irreducible subrepresentations $V_{\gamma, j}^{(n)}$, where $\gamma$ denotes the highest weight
     of $V_{\gamma, j}^{(n)}$ and $j$ counts the multiplicity up to isomorphism.
     \end{rem}

 \begin{defn}\label{defn_of_G_q} Let  ${\mathcal G}_q$ denote the   completion with idempotents and direct sums of
 the   full linear subcategory     of 
 ${\mathcal T}({\mathfrak g}, q, \ell)$ with objects the truncated tensor powers
$V^{\underline{\otimes} n}$. 
\end{defn}

Thus by construction ${\mathcal G}_q$ is a semisimple linear category, that we regard it as an abstract  category.
Let ${\rm Vec}$ be realized as a strict tensor category.
We regard  $V^{\underline{\otimes} n}$ as a summand on $V^{{\otimes} n}$ via $p_n$, and identify the morphism  space $(V^{\underline{\otimes} m}, V^{\underline{\otimes} n})$ in ${\mathcal G}_q$ with
 the subspace of morphisms $T\in (V^{{\otimes} m}, V^{{\otimes} n})$ in ${\mathcal T}({\mathfrak g}, q, \ell)$ satisfying $Tp_m=p_nT=T$.
We set
 \begin{equation} V^{\underline{{\otimes}} m}\underline{\otimes}V^{\underline{{\otimes}} n}:=V^{\underline{{\otimes}} m+n},\end{equation}
 \begin{equation}  S\underline{\otimes} T:=p_{m'+n'}\circ S\otimes T\circ p_{m+n},\quad\quad S\in (V^{\underline{{\otimes}} m}, V^{\underline{{\otimes}} m'}), \quad T\in (V^{\underline{{\otimes}} n}, V^{\underline{{\otimes}} n'}).\end{equation}

\begin{thm}\label{equivalence_G_q_with_qg_fusion_category}
With the above tensor product and trivial associativity morphisms
  ${\mathcal G}_q$,   becomes a ribbon strict ribbon tensor category ribbon tensor equivalent
 to ${\mathcal C}({\mathfrak g}, q, \ell)$.
 \end{thm}

The previous theorem is essentially  in \cite{Wenzl},  perhaps the strictness property was noticed in
 Theorem 5.4 in \cite{CP}. 
 This property plays a role in our proof of Theorem \ref{Finkelberg_HL}, in Sect. \ref{21}.
 Moreover, we   describe the equivalence   in this section,
 see Remark \ref{remark_on_equivalence}.
 We next introduce a   concrete version of ${\mathcal G}_q$.
  
 \begin{defn}\label{W_functor} Let $\tilde{\mathcal G}_q$ denote
 the full representation subcategory   of ${\mathcal T}({\mathfrak g}, q, \ell)$ with objects representations which are finite direct sums of summands
 of the representations $ V^{\underline{{\otimes}} n}$.  
 \end{defn}

 Then $\tilde{\mathcal G}_q$ is also a linear semisimple category.
 There is a canonical linear equivalence
 $$\tilde{\mathcal G}_q\to {\mathcal G}_q$$
 taking the summand $P^{(n)}V^{\underline{\otimes} n}$ 
 of $V^{\underline{\otimes} n}$
 defined by an idempotent
 $P^{(n)}$
 regarded as an object of $\tilde{\mathcal G}_q$ to
 $P^{(n)}$ regarded as an object of ${\mathcal G}_q$ and acting trivially on morphisms.
 We shall make $\tilde{\mathcal G}_q$ into a tensor category with  a   tensor structure  $(\tilde{\mathcal G}_q, \boxtimes, \alpha)$ in such a way that ${\mathcal E}$ becomes a tensor equivalence
 $({\mathcal E}, E)$.

Let then $${W}: \tilde{\mathcal G}_q\to {\rm Vec}$$ be the forgetful functor. 
To define a tensor structure  on $\tilde{\mathcal G}_q$ we first define   linear maps
     $(F_{\lambda, \mu}, G_{\lambda, \mu})$    on  ${W}$ defined on pairs $(\lambda, \mu)\in\Lambda^+(q)\times\Lambda^+(q)$
   that will correspond to a  tensor structure $(\boxtimes, \alpha)$ of  $\tilde{\mathcal G}_q$ and subsequently also to a 
   weak tensor structure for ${W}$.

For every $\lambda\in\Lambda^+(q)$ choose an integer $n_\lambda$ such that $\lambda$ appears as the dominant weight of a summand $V_\lambda$ of $V^{\underline{\otimes} n_\lambda}$ as observed in Remark \ref{canonical_iterative}.
Let 
 $p_\lambda: V^{\underline{\otimes} n_\lambda}\to V_\lambda$  denote the corresponding idempotent onto $V_\lambda$ for each $\lambda\in\Lambda^+(q)$. In the following formulae we extend $p_\lambda$ to $V^{\otimes n_\lambda}$ in a trivial way on $(1-p_{n_\lambda})V^{\otimes n_\lambda}$.
 
 \begin{prop}\label{maximal_non_negligible} We have that
$p_\lambda\underline{\otimes} p_\mu=p_{n_\lambda+n_\mu} p_\lambda\otimes p_\mu p_{n_\lambda+n_\mu}$
is an idempotent in the semisimple category  $\tilde{\mathcal G}_q$ onto a module   isomorphic to a maximal non-negligible submodule 
of $V_\lambda\otimes V_\mu$ in ${\mathcal T}({\mathfrak g}, q, \ell)$.
 \end{prop}
 
 \begin{proof} Notice that
$p_\lambda\underline{\otimes} p_\mu$ is   a morphism in ${\mathcal T}({\mathfrak g}, q, \ell)$ and is an idempotent by $(2)$ in Subsect. \ref{18.6} with range in the semisimple part $V^{\underline{\otimes}(n_\lambda+n_\mu)}$,
thus this range is a semisimple representation depending only on $\lambda$, $\mu$ up to isomorphism.
If $M_{\lambda, \mu}$ is a maximal idempotent onto a nonnegligible summand of $p_\lambda V^{\otimes n_\lambda}\otimes p_\mu V^{\otimes {n_\mu}}$
then $p_\lambda\underline{\otimes} p_\mu= p_{n_\lambda+n_\mu} M_{\lambda, \mu} p_{n_\lambda+n_\mu}$.
We have that $T=p_{n_\lambda+n_\mu} M_{\lambda, \mu}\in(M_{\lambda, \mu}, p_{n_\lambda+n_\mu})$
and $T^{-1}=M_{\lambda, \mu}p_{n_\lambda+n_\mu}  \in(p_{n_\lambda+n_\mu}, M_{\lambda, \mu})$
satisfy $T^{-1}T=M_{\lambda, \mu}$ and $TT^{-1}=p_\lambda\underline{\otimes} p_\mu$.
 \end{proof}

We define
$V_\lambda\boxtimes V_\mu:=p_{n_\lambda}\underline{\otimes} p_\mu V^{{\otimes}(n_\lambda+n_\mu)}$ as a module
of $\tilde{\mathcal G}_q$, 
thus $W(V_\lambda\boxtimes V_\mu)=p_{n_\lambda+n_\mu}p_\lambda\otimes p_\mu V^{\underline{\otimes}(n_\lambda+n_\mu)}$ as a linear space. 

\begin{defn}\label{weak_tensor_structure}
For $\lambda$, $\mu\in\Lambda^+(q)$, let
$$F_{\lambda, \mu}: V_\lambda\otimes V_\mu\to   V_\lambda\boxtimes V_\mu,\quad\quad G_{\lambda, \mu}:
V_\lambda\boxtimes V_\mu\to V_\lambda\otimes V_\mu$$ be the morphisms in ${\mathcal T}({\mathfrak g}, q, \ell)$
  respectively defined as the   restriction of $p_\lambda\underline{\otimes} p_\mu=p_{n_\lambda+n_\mu} p_\lambda\otimes p_\mu p_{n_\lambda+n_\mu}$
  to  $V_\lambda\otimes V_\mu$ and that of $p_\lambda\otimes p_\mu$ to $V_\lambda\boxtimes V_\mu$.
  Thus   we have linear maps
$$F_{\lambda, \mu}: W(V_\lambda)\otimes W(V_\mu)\to   W(V_\lambda\boxtimes V_\mu),\quad\quad G_{\lambda, \mu}:
W(V_\lambda\boxtimes V_\mu)\to W(V_\lambda)\otimes W(V_\mu).$$

  \end{defn}

\begin{prop}\label{properties_of_weak_tensor_structure}
We have that $F_{\lambda, \mu}G_{\lambda, \mu}=1$ and $G_{\lambda, \mu}F_{\lambda, \mu}$ is an idempotent
of ${\mathcal T}({\mathfrak g}, q, \ell)$ onto a maximal non-negligible submodule of $V_\lambda\otimes V_\mu$.
\end{prop}

\begin{proof} The first statement is again a simple consequence of $(2)$ in Subsect. \ref{18.6} the remaining part
 follows   from this and Prop. \ref{maximal_non_negligible}.
\end{proof}

We next extend $\boxtimes$ and $F_{\lambda, \mu}$, $G_{\lambda, \mu}$ to all objects of $\tilde{\mathcal G}_q$.
Let $P^{(n)}\in (V^{\underline{\otimes} n}, V^{\underline{\otimes} n})$,
$P^{(m)}\in (V^{\underline{\otimes} m}, V^{\underline{\otimes} m})$
 be   idempotents in ${\mathcal G}_q$ and consider morphisms describing decomposition into irreducibles,
 that is 
$$S_{\lambda, j}: V_\lambda \to P^{(n)}V^{\underline{\otimes} n},\quad
S'_{\lambda, j}: P^{(n)}V^{\underline{\otimes} n}\to V_\lambda,$$
 $$T_{\mu, k}: V_\mu \to P^{(m)}V^{\underline{\otimes} m}, \quad
 T'_{\mu, k}: P^{(m)}V^{\underline{\otimes} m}\to V_\mu,$$
$$S'_{\lambda,  j} S_{\lambda', j'}=\delta_{(\lambda, j), (\lambda', j')}, \quad
\sum_{\lambda, j} S_{\lambda, j}S'_{\lambda, j}=P^{(n)},$$
$$T'_{\mu,  k} T_{\mu', k'}=\delta_{(\mu, k), (\mu', k')}, \quad \sum_{\mu, k} T_{\mu, k}T'_{\mu, k}=P^{(m)}.$$
We set
\begin{equation}\label{F} F_{P^{(n)}V^{\underline{\otimes} n}, P^{(m)}V^{\underline{\otimes} m}}=\sum_{\lambda, j, \mu, k}
  S_{\lambda, j} \underline{\otimes} T_{\mu, k}  \circ F_{\lambda, \mu}\circ S'_{\lambda, j}\otimes T'_{\mu, k},\end{equation}
 \begin{equation}\label{boxtimes_for_objects} P^{(n)}V^{\underline{\otimes} n}\boxtimes P^{(m)}V^{\underline{\otimes} m}=F_{P^{(n)}V^{\underline{\otimes} n}, P^{(m)}V^{\underline{\otimes} n}}(P^{(n)}V^{\underline{\otimes} n} \otimes P^{(m)}V^{\underline{\otimes} m})\end{equation}
 and we let $G_{P^{(n)}V^{\underline{\otimes} n}, P^{(m)}V^{\underline{\otimes} m}}$ be the restriction of 
 \begin{equation}\label{G} \sum_{\lambda, j, \mu, k}
S_{\lambda, j} \otimes T_{\mu, k} \circ  G_{\lambda, \mu}\circ   S'_{\lambda, j}\underline{\otimes} T'_{\mu, k}\end{equation} to $P^{(n)}V^{\underline{\otimes} n}\boxtimes P^{(m)}V^{\underline{\otimes} m}$.
 Notice that $F_{P^{(n)}V^{\underline{\otimes} n}, P^{(m)}V^{\underline{\otimes} m}}$ and $G_{P^{(n)}V^{\underline{\otimes} n}, P^{(m)}V^{\underline{\otimes} m}}$ are independent of the choice of 
 $S_{\lambda, j}$, $S'_{\lambda, j}$, $T_{\mu, k}$, $T'_{\mu, k}$ by bilinearity of $\otimes$. In particular,
 these maps and tensor products extend the previous ones on the chosen class of irreducibles.
 Finally, we extend this structure to any object of $\tilde{\mathcal G}_q$ by bilinearity.

\begin{rem}\label{maximal_non_negligible2}
 It follows that $F_{\lambda, \mu}$, $G_{\lambda, \mu}$  of Def. \ref{weak_tensor_structure} and Prop. \ref{properties_of_weak_tensor_structure} are morphisms of the tilting category onto a maximal non-negligible
 addendum.
 This property extends to all pairs of objects of $\tilde{\mathcal G}_q$ in place of $V_\lambda$ and $V_\mu$.  Note that 
   this property is not needed in the definition of weak quasi-tensor structure, but it turns out useful to construct the
 weak  Hopf algebra of this section. Moreover,   this flexibility of the notion of weak quasi-tensor structure will turn out useful
 in Sect. \ref{21}, where we shall construct a strongly unitary weak quasi-tensor structure
 $(F_0, F_0^*)$ for the Zhu algebra $A_{V_{{\mathfrak g}_k}}$ associated to an affine VOA for a suitable positive integer $k$
that is a morphism for the action of the simple Lie algebra ${\mathfrak g}$ via a twist that is not a morphism of the tilting category.

 \end{rem}
 
We next define a tensor product between morphisms  and associativity morphisms in $\tilde{\mathcal G}_q$
 as follows. Let $\rho$, $\sigma$, $\tau$ be objects of $\tilde{\mathcal G}_q$.
 For $S: \rho\to \rho'$, $T:\sigma'\to\sigma'$, set
  \begin{equation}\label{boxtimes}S\boxtimes T=F_{\rho', \sigma'} S\otimes T G_{\rho, \sigma}.\end{equation}
We endow 
  $\tilde{\mathcal G}_q$ with associativity morphisms
  \begin{equation}\label{associator_for_boxtimes} \alpha_{\rho, \sigma,\tau}=F_{\rho, \sigma\boxtimes\tau}\circ 1_\rho\otimes F_{\sigma, \tau}\circ   G_{\rho, \sigma}\otimes 1_\tau \circ G_{{\rho\boxtimes \sigma}, \tau}\end{equation}
     
  Note that the maps $F$ and $G$ are   defined up to varying   the choice of the integers $n_\lambda$ or the   definition of 
 the idempotents $p_n$.

  \begin{thm} \label{first1} Let ${\mathfrak g}$ be a complex simple Lie algebra, $q$ a complex root of unity such that $q^2$ has order $\ell$ in the sense of Def. \ref{large_enough}. Then 
  \begin{itemize}
\item[{\rm        a)}] 
$(\tilde{\mathcal G}_q, \boxtimes, \alpha)$ is a semisimple tensor category,
\item[{\rm        b)}] the canonical linear equivalence ${\mathcal E}: \tilde{\mathcal G}_q\to{\mathcal G}_q$ admits
a unique structure  of  tensor equivalence $({\mathcal E}, E): (\tilde{\mathcal G}_q, \boxtimes, \alpha)\to ({\mathcal G}_q, \underline{\otimes }, 1)$
  such that
$$E_{\lambda, \mu}: {\mathcal E}(V_\lambda)\underline{\otimes} {\mathcal E}(V_\mu)\to{\mathcal E}(V_\lambda\boxtimes V_\mu),\quad \quad \lambda, \mu\in\Lambda^+(q)$$ acts as $F_{\lambda, \mu}$ and we have that
$E_{\lambda, \mu}^{-1}$ acts as $G_{\lambda, \mu}$, 
  \item[{\rm        c)}] 
 the pair $(F, G)$ is a weak tensor    structure for the forgetful functor ${W}: \tilde{\mathcal G}_q\to{\rm Vec}$, therefore $A_{{W}}({\mathfrak g}, q, \ell)={\rm Nat}_0({W})$ is a ribbon weak  Hopf algebra,
  \item[{\rm        d)}] 
a different   choice of $p_n$, $p_\lambda$ changes $A_{{W}}({\mathfrak g}, q, \ell)$ by a trivial twist.

 \end{itemize}
   \end{thm}
   
   \begin{proof}
   a) Note that $S\boxtimes T$ is composition of morphisms in ${\mathcal T}({\mathfrak g}, q, \ell)$ with domain and range 
   representations of $\tilde{\mathcal G}_q$, thus it is a morphism in $\tilde{\mathcal G}_q$.
   By Remark \ref{maximal_non_negligible2}, for any pair of objects
   $\rho$, $\sigma\in\tilde{\mathcal G}_q$,
   $G_{\rho, \sigma}F_{\rho, \sigma}$ is an idempotent in ${\mathcal T}({\mathfrak g}, q, \ell)$ 
   with range a maximal non-negligible summand of the tensor product  tilting module $\rho\otimes\sigma$.
   Thus $1-G_{\rho, \sigma}F_{\rho, \sigma}$ is an idempotent onto the negligible summand. This observation
   together with property $(2)$ in Subsect. \ref{18.6} implies that $\boxtimes$ is a bifunctor of $\tilde{\mathcal G}_q$.
   The pentagon equation can be shown again taking into account property (2) and we also need (3).
   For example computing the short side of the pentagon equation (\ref{pentagon_equation})
   $$\alpha_{\nu,\rho, \sigma\tau}\alpha_{\nu\rho, \sigma, \tau}=$$
   $$F_{\nu, \rho(\sigma\tau)}\circ 1_{\nu}\otimes F_{\rho, \sigma\tau}\circ G_{\nu, \rho}\otimes 1_{\sigma\tau}\circ G_{\nu\rho, \sigma\tau}F_{\nu\rho, \sigma\tau}\circ 1_{\nu\rho}\otimes F_{\sigma, \tau}\circ G_{\nu\rho, \sigma}\otimes 1_{\tau}\circ G_{(\nu\rho)\sigma, \tau}$$
   we may first eliminate the central term $G_{\nu\rho, \sigma\tau}F_{\nu\rho, \sigma\tau}$, then use the commutation relation
   $$G_{\nu, \rho}\otimes 1_{\sigma\tau}\circ 1_{\nu\rho}\otimes F_{\sigma, \tau}= 1_{\nu\rho}\otimes F_{\sigma, \tau}\circ
   G_{\nu, \rho}\otimes 1_{\sigma\tau}$$
   thus
$$\alpha_{\nu,\rho, \sigma\tau}\alpha_{\nu\rho, \sigma, \tau}=F_{\nu, \rho(\sigma\tau)}\circ 1_{\nu}\otimes F_{\rho, \sigma\tau}\circ1_{\nu\rho}\otimes F_{\sigma, \tau}\circ
   G_{\nu, \rho}\otimes 1_{\sigma\tau}\circ G_{\nu\rho, \sigma}\otimes 1_{\tau}\circ G_{(\nu\rho)\sigma, \tau}.$$
   The computation involving the long side of the pentagon equation is slightly longer because of the use of $\boxtimes$ at both sides. However it can   patiently be carried out and it leads to equating the left hand side.
   b) It is clear that $E_{\lambda, \mu}$ and $E_{\lambda, \mu}^{-1}$ are morphisms and are inverse of each other.
   Furthermore extending these morphisms by naturality to every pair of objects we see that they act as
   $F$ and $G$ respectively. Then we may verify the tensoriality equation (\ref{wt1}) for $({\mathcal E}, E)$.
   To do this, notice that the tensor product $\underline{\otimes}$ at right hand side of (\ref{wt1}) modifies $\otimes$ by inserting suitable idempotents $p_n$ which may then be disregarded thanks to (2) again.
   c) Naturality of $F$ and $G$ as transformations from $\tilde{\mathcal G}\to{\rm Vec}$ may be checked with direct computation. Notice also that by construction $F$ and $G$ are natural as transformations ${\mathcal G}_q\to{\rm Vec}$, therefore by composition $\tilde{\mathcal G}_q\to{\mathcal G}_q\to{\rm Vec}$ we find that they are also natural with respect to ${\boxtimes}$. Property d) follows again from (2).
   
   \end{proof}

\begin{rem}\label{remark_on_equivalence} It follows from part b) of the previous theorem that the composition
  $${\mathcal Q}:\tilde{\mathcal G}_q\to{\mathcal T}({\mathfrak g}, q, \ell)\to {\mathcal C}({\mathfrak g}, q, \ell)$$ of the natural  inclusion followed by quotient is an equivalence of tensor categories.
  In this way, $\tilde{\mathcal G}_q$ admits a unique structure of a ribbon category in a way that ${\mathcal Q}$ is a ribbon equivalence.  On the other hand, this can also be  seen   directly.
  \end{rem}
  \bigskip
  
  \subsection{Unitary coboundary structure of $A_W({\mathfrak g}, q, \ell)$.} \label{29.2}In this subsection we consider
  the case of unitary structures, thus we assume that
  $q=e^{i\pi/\ell}$ is a minimal root of  order  in the sense of Def. \ref{large_enough} and we study the unitarity property of 
  $A_W({\mathfrak g}, q, \ell)$. Note that by Remark \ref{non_minimal_roots}, most of the constructions hold
  for non-minimal roots of this kind of order.

  We recall from Theorem \ref{U_q_as_a_Hermitian_ribbon_h} that $U_q({\mathfrak g})$ is a (topological) Hermitian coboundary Hopf algebra with compatible involution and antipode of Kac type.
  Furthermore recall also that 
 by Prop. 2.4 in \cite{Wenzl}, 
   for $\lambda\in\Lambda^+(q)$ the natural Hermitian form of $V_\lambda$ in the sense of Sect. \ref{11} is a positive definite inner product, so $V_\lambda$ is a $C^*$-representation of $U_q({\mathfrak g})$.

   \begin{defn}
Let ${\mathcal T}_W$ denote the full subcategory of ${\mathcal T}({\mathfrak g}, q, \ell)$ with objects orthogonal direct sums of   summands defined by selfadjoint idempotents of
finite tensor products of $V_\lambda$ with $\lambda\in\Lambda^+(q)$ endowed with
   the non-degenerate Hermitian form induced by iterates of $\overline{R}$ of $U_q({\mathfrak g})$.

\end{defn}
   
   Consider a finite tensor product $W$ of $V_\lambda$ with $\lambda\in\Lambda^+(q)$ endowed with
   the non-degenerate Hermitian form induced by iterates of $\overline{R}$ of $U_q({\mathfrak g})$, or more generally a 
an orthogonal direct sum of   summands defined by selfadjoint idempotents of a module of this kind such that the form is nondegenerate on $W$.  For any morphism $T: W\to W'$ of ${\mathcal T}_W$,
  the adjoint $T^*: W'\to W$
   is well defined.   
   We next consider in particular the canonical decomposition into indecomposable tilting modules
     recalled in the previous subsection
   $V_\lambda\otimes V=\oplus_\gamma 
 T_\gamma\otimes {\mathbb C}^{m_\gamma}$. The  Hermitian form induced by  $\overline{R}$ is positive definite on the 
 on the isotypic component  $T_\gamma\otimes {\mathbb C}^{m_\gamma}=V_\gamma\otimes {\mathbb C}^{m_\gamma}$    for $\gamma\in\Lambda^+(q)$.   The idempotents $p_{\lambda, \gamma} V_\lambda\otimes V\to V_\gamma\otimes{\mathbb C}^{m_\gamma}$
 are selfadjoint  with respect to this inner product. It follows that the iterated tensor powers $V^{\underline{\otimes} n}$
 are Hilbert space representations of $U_q({\mathfrak g})$ with this iterated Hermitian form.
Let  ${\mathcal G}_\ell$ denote the completion under selfadjoint idempotents and orthogonal direct sums of the full
 subcategory of ${\mathcal T}({\mathfrak g}, q, \ell)$ with objects $V^{\underline{\otimes} n}$.
Thus ${\mathcal G}_\ell$
 has the structure of a linear semisimple $C^*$-category. Furthermore with tensor product $\underline{\otimes}$ defined
 as in the previous subsection, ${\mathcal G}_\ell$ becomes a unitary strict tensor category.
  
  \begin{prop}\label{selfadjointness}
  The idempotents $p_\lambda$ and $p_\lambda\underline{\otimes} p_\mu$ are selfadjoint in ${\mathcal G}_\ell$.
   \end{prop}
  
\begin{proof}  Notice that for all $n$, the idempotents say $p_{\gamma, j}$ onto the irreducible decomposition  $V^{(n)}$
  by the  $V^{(n)}_{\gamma, j}$  
  described in Remark \ref{canonical_iterative} have pairwise orthogonal ranges with respect to the inner product
by  orthogonality of the addenda of $V_\lambda\otimes V$ with dominant weights in $\Lambda^+(q)$   and the iterative construction of $V^{\underline{\otimes} n}$. It follows in particular that    $p_\lambda$ are selfadjoint idempotent
in ${\mathcal G}_\ell$. Let  $c^w(\lambda, \mu)$ be the coboundary operators in ${\mathcal T}({\mathfrak g}, \ell, q)$ associated to the $\overline{R}$ matrix of $U_q({\mathfrak g})$. By naturality
we have $p_\lambda\otimes p_\mu c^w(\mu, \lambda)=c^w(\mu, \lambda) p_\mu\otimes p_\lambda$, 
and thus $\overline{R}$ commutes with $p_\lambda\otimes p_\mu$. It follows that $p_\lambda\otimes p_\mu$ is selfadjoint with respect to 
 the iterated Hermitian form of   $V^{\otimes{n_\lambda+n_\mu}}$, and therefore also $p_\lambda\underline{\otimes} p_\mu$ are selfadjoint.
 \end{proof}
 
 We then similarly  introduce   the concrete category  $\tilde{\mathcal G}_\ell$ of ${\mathcal T}({\mathfrak g}, q, \ell)$
 taking into consideration summands defined by selfadjoint idempotents and orthogonal direct sums.
  By the previous proposition, for 
  $\lambda$, $\mu\in\Lambda^+(q)$, $F_{\lambda, \mu}$ and $G_{\lambda, \mu}$ introduced as in the previous subsection are morphisms in ${\mathcal T}_W$ and
    satisfy in addition the property $F_{\lambda, \mu}^*=G_{\lambda, \mu}$. Using orthogonal decompositions
    of objects of $\tilde{\mathcal G}_\ell$, we obtain natural transformations $F_{\rho, \sigma}$, $G_{\rho, \sigma}$ as in (\ref{F}), and (\ref{G}), satisfying $F_{\rho, \sigma}^*=G_{\rho,\sigma}$ in ${\mathcal T}_W$. Then we introduce in $\tilde{\mathcal G}_\ell$ the structure of a tensor category
    $(\tilde{\mathcal G}_\ell, \boxtimes, \alpha)$
   as in (\ref{boxtimes_for_objects}), (\ref{boxtimes}),  and (\ref{associator_for_boxtimes}).
   Finally, we   consider  the forgetful functor,
  $$W: \tilde{\mathcal G}_\ell\to {\rm Hilb}$$ and is a $^*$-functor endowed with the weak tensor structure $(F, G)$ regarded with values in ${\rm Hilb}$. (Notice that as natural transformations of $W$,   we do not have 
  $F_{\rho, \sigma}^*=G_{\rho, \sigma}$, more details will be discussed in the proof of the following result.)

\begin{thm}\label{second2}  Let $q$ be a minimal root of order  in the sense of Def. \ref{large_enough},  then
 \begin{itemize}
\item[{\rm        a)}] $(\tilde{\mathcal G}_\ell, \boxtimes, \alpha)$ is a unitary semisimple tensor category and the  tensor equivalence
$({\mathcal E}, E):(\tilde{\mathcal G}_\ell, \boxtimes, \alpha) \to({\mathcal G}_\ell, \underline{\otimes}, 1)$ is unitary,
\item[{\rm        b)}]    $A_{{W}}({\mathfrak g}, q, \ell)={\rm Nat}_0(W)$ becomes a  unitary coboundary weak weak  Hopf algebra with compatible involution,  weak tensor structure defined by $(F, G)$ and  antipode of Kac type
such that ${\mathcal G}_\ell\to{\rm Rep}(A_W({\mathfrak g}, q, \ell))$ is a unitary equivalence of ribbon categories.
\end{itemize}
 
\end{thm}

\begin{proof}  
a) The property $(S\boxtimes T)^*=S^*\boxtimes T^*$ follows from the relation $F^*=G$ in ${\mathcal T}_{\Lambda^+(q)}$ and arguments similar to those in the proof of Prop. \ref{selfadjointness}. Unitarity of the associator follows from
$F^*=G$ and $(2)$, and $(3)$.
b) By theorem
Theorem \ref{TK_unitary_ribbon} we need to show
 (\ref{unitary_ribbon_functor1}) and  (\ref{unitary_ribbon_functor2}). We only show the former.
  By Remark \ref{remark_to_unitary_ribbon} a), it is enough to do this for $\rho=V^{\underline{\otimes} n}$, $\sigma=V$. In this case $F_{\rho, \sigma}=p_{n+1}$ as $(1-p_n)\otimes 1_V$ is negligible. This follows by construction as
  $F_{\rho, \sigma}^*=p_{n+1}^*\overline{R}=\overline{R}p_{n+1}$ where $p_{n+1}^*$ is the adjoint with respect to the standard inner product
  of $V^{\underline{\otimes} n}\otimes V$.    We next show that $A_W({\mathfrak g}, q, \ell)$ has an   antipode of Kac type. It is shown in the proof 
  of Lemma 10.4 in \cite{CP} that a solution of the conjugate equations in ${\mathcal G}_\ell$  is of the form (\ref{computing_conjugates2}) with $\alpha=\beta=I$, $\mu_i=1$. It follows from  the proof of Theorem \ref{TK_algebraic_quasi} that the corresponding antipode is strong and therefore of Kac type. 
  The proof  of compatibility with the $^*$-involution is the
  content of the following lemmas \ref{compatibility1}, \ref{k=1}, \ref{k=2}.

  \end{proof}

 \bigskip

\begin{lemma}\label{compatibility1}
Let $\overline{c}_U(\rho, \sigma)$ be the natural coboundary symmetry associated to $U_q({\mathfrak g})$. Then   the unitary coboundary weak  Hopf algebra
$A_W$ has compatible $^*$-involution if and only if 
\begin{equation}\label{compatibility_condition1}
\overline{c}_U(V_\lambda, V^{\underline{\otimes} k}) 
G_{V_\lambda, V^{\underline{\otimes} k}}F_{V_\lambda, V^{\underline{\otimes} k}}
\overline{c}_U(V^{\underline{\otimes} k}, V_\lambda)=G_{V^{\underline{\otimes} k}, V_\lambda}
 F_{{V^{\underline{\otimes} k}, V_\lambda}} \quad
 \lambda\in\Lambda^+(q), \quad\quad k=1, 2.
 \end{equation}
It suffices that the following two equations involving the braided symmetries $c_U$
   and  $c_U^{-1}$ associated to $U_q({\mathfrak g})$ and also the braided symmetries $c$ and $c^{-1}$ associated to $A_W$, hold,
 \begin{equation}\label{compatibility_condition2}
c(V_\lambda, V^{\underline{\otimes} k}) F_{V_\lambda, V^{\underline{\otimes} k}}
c_U(V_\lambda, V^{\underline{\otimes} k})^{-1}= 
 F_{{V^{\underline{\otimes} k}, V_\lambda}}, \quad
 \lambda\in\Lambda^+(q), \quad\quad k=1, 2,
\end{equation}
\begin{equation}\label{compatibility_condition3}
c(V^{\underline{\otimes} k}, V_\lambda)^{-1} F_{V_\lambda, V^{\underline{\otimes} k}}
c_U(V^{\underline{\otimes} k}, V_\lambda)= 
 F_{{V^{\underline{\otimes} k}, V_\lambda}}, \quad
 \lambda\in\Lambda^+(q), \quad\quad k=1, 2.
\end{equation}

\end{lemma}

\begin{proof}
Taking the adjoint of equations (\ref{compatibility_condition2}), (\ref{compatibility_condition3}) and multiplying them term by term we get
\begin{equation}\label{compatibility_condition4}
c_U(V_\lambda, V^{\underline{\otimes} k}) 
G_{V_\lambda, V^{\underline{\otimes} k}}F_{V_\lambda, V^{\underline{\otimes} k}}
c_U(V_\lambda, V^{\underline{\otimes} k})^{-1}=G_{V^{\underline{\otimes} k} , V_\lambda }
 F_{{V^{\underline{\otimes} k}, V_\lambda}} \quad
 \lambda\in\Lambda^+(q), \quad\quad k=1, 2.
\end{equation}
\begin{equation}\label{compatibility_condition5}
c_U(V^{\underline{\otimes} k}, V_\lambda)^{-1}
G_{V_\lambda, V^{\underline{\otimes} k}}F_{V_\lambda, V^{\underline{\otimes} k}}
c_U(V^{\underline{\otimes} k}, V_\lambda)=G_{V^{\underline{\otimes} k} , V_\lambda }
 F_{{V^{\underline{\otimes} k}, V_\lambda}} \quad
 \lambda\in\Lambda^+(q), \quad\quad k=1, 2.
\end{equation}
In turn it follows that
 $c_U^2$ commutes with
$G_{V^{\underline{\otimes} k} , V_\lambda }
 F_{{V^{\underline{\otimes} k}, V_\lambda}}$. It follows that the principal branch square root commutes also, and this implies
(\ref{compatibility_condition1}). 

We next show the first statement.
By Prop. \ref{Positivity_coboundary},  compatibility of the $^*$-involution is equivalent to
$\Delta^{\rm op}(I)=\Delta(I)^*$ on the   spaces
of $V^{\underline{\otimes} k}\otimes V_\lambda$  and $V_\lambda\otimes V^{\underline{\otimes} k}$ for $k=1$, $2$. We have
  $\Delta(I)^*=\overline{R}^U\Delta(I)(\overline{R}^U)^{-1}$,   and it follows that   the desired equalities reduce to our assumptions. Note that  equation (\ref{compatibility_condition1}) together with 
   the
 coboundary property $\overline{c}_U^2=1$, see Prop. \ref{Wenzl_Hermitian_form}, d),
 imply
  that the symmetric equation 
 with $V_\lambda$ on the right and $V^{\underline{\otimes}k}$ on the left  at the l.h.s. of the equation holds and this completes the proof.
  \end{proof}

  \begin{lemma}\label{k=1} The natural transformation $F$ defining   $A_W$ satisfies equations (\ref{compatibility_condition2}), 
  (\ref{compatibility_condition3})  for $k=1$.

  \end{lemma}

  \begin{proof} Assume ${\mathfrak g}\neq E_8$.
For $k=1$, by \cite{Wenzl}, $V\otimes V_\lambda$ is completely reducible into irreducible
 components $\oplus_\mu m_\mu V_\mu$ (with multiplicity $0$ or $1$ except for ${\mathfrak g}=F_4$ where  $\mu_\mu>1$ only for $\mu=\lambda$) and we have that $\mu\in\overline{\Lambda^+(q)}$. Thus there is 
 a unique morphism idempotent onto a maximal non-negligible submodule $V\otimes V_\lambda\to\oplus_{\mu\in\Lambda^+(q)} V_\mu$ which then coincides with $F_{V, V_\lambda}$.
 This uniqueness property and unitarity of the braided symmetries imply that (\ref{compatibility_condition2}), 
  (\ref{compatibility_condition3}) hold for $k=1$.
The case ${\mathfrak g}=E_8$ is more delicate that the others, and is not covered by the above proof.
In this case we consider the decomposition of $V_\lambda\otimes V$ and of $V\otimes V_\lambda$
into indecomposable tilting modules
given at page 274 in \cite{Wenzl}. 
Let  $F_{V_\lambda, V}$ and $F_{V, V_\lambda}$ be
  the corresponding idempotents onto the maximal non-negligible submodules.
  Then it follows from the proof therein and unitarity of the braided symmetry, that
   (\ref{compatibility_condition2}), 
  (\ref{compatibility_condition3}) hold in this case. More precisely, we consider $W_l=V_\lambda\otimes V$ considered
  in the mentioned proposition of \cite{Wenzl} and we pair it with $W_r=V\otimes V_\lambda$, thus $W_l$ is $W$ in \cite{Wenzl}. Similarly, for the submodules of $W_l$ and $W_r$ we follow the same notation, with an addition of
  a subscript $l$ or $r$, for example $M_l=(p_\lambda+p_{\lambda+k})W_l$ and $M_r=(p_\lambda+p_{\lambda+k})W_r$.
  The braiding $c_U$ intertwines $W_l$ with $W_r$, and therefore also $M_l$ with $M_r$, as each addendum isomorphic to some $V_\mu(q)$ with $\mu\neq\lambda$ has multiplicity $1$. Moreover $M_l$ and $M_r$ contain a submodule isomorphic to $V_{\lambda+\kappa}(q)$. We note that these submodules, that we denote by 
  $V_{\lambda+\kappa}(q)_l$ and $V_{\lambda+\kappa}(q)_r$ are generated by the $(\lambda+k)$-weight space, that is one dimensional, thus $c_U$ takes $V_{\lambda+\kappa}(q)_l$ to $V_{\lambda+\kappa}(q)_r$. Moreover $c_U$ is unitary. It follows from these two facts that choosing $T(q)_r=c_UT(q)_l$,
  $c_U$ takes $(M_\lambda)_l$ to $(M_\lambda)_r$, and therefore  $c_U(p_\lambda)_l=(p_\lambda)_rc_U$
  and $c_U(p_{\lambda+\kappa})_l=(p_{\lambda+\kappa})_rc_U$. Note that $T(q)$ involve a choice. More in detail,
  $T(q)/V_{\lambda+\kappa}$ complements the submodule $N(q)^\perp/V_{\lambda+\kappa}(q)$ of $M(q)/V_{\lambda+\kappa}(q)$ (this is semisimple). By the property of $c$, we may choose $T(q)_r=cT(q)_l$.
     \end{proof}

   \begin{lemma}\label{k=2}
   The natural transformation $F$ defining   $A_W$ satisfies equations (\ref{compatibility_condition2}), 
  (\ref{compatibility_condition3})  for $k=2$ and all Lie types.
   \end{lemma} 
    
 \begin{proof}
We use the weak  Hopf property in categorical form (\ref{wt1}), (\ref{wt2}). Working with ${\rm Vec}$ strict,

\begin{equation}\label{a}
\quad {\mathcal F}(\alpha_{V_\lambda, V, V})       =F_{V_\lambda, V\underline{\otimes} V}\circ 1{\otimes} 
F_{V, V}\circ G_{V_\lambda, V}{\otimes}1\circ G_{V_\lambda\underline{\otimes} V, V},  \end{equation}

 \begin{equation}\label{b}\quad {\mathcal F}(\alpha_{V_\lambda, V, V})^{-1}=F_{V_\lambda\underline{\otimes} V, V}\circ F_{V_\lambda, V}{\otimes} 1\circ 1{\otimes} G_{V, V}
\circ G_{{V_\lambda, V\underline{\otimes} V}},\end{equation}

\begin{equation}\label{c}\quad {\mathcal F}(\alpha_{V, V, V_\lambda})=F_{V, V\underline{\otimes} V_\lambda}\circ 1{\otimes} F_{V, V_\lambda} \circ  G_{V, V} {\otimes}1
\circ G_{V\underline{\otimes} V, V_\lambda},\end{equation}

\begin{equation}\label{d}\quad {\mathcal F}((\alpha_{V, V, V_\lambda})^{-1})  =F_{V\underline{\otimes} V, V_\lambda}\circ 
F_{V, V} {\otimes} 1\circ 1 {\otimes}G_{V, V_\lambda}\circ G_{V, V\underline{\otimes} V_\lambda}.\end{equation}
We set

 $$\tilde{F}_{1,2}:=F_{V_\lambda, V\underline{\otimes} V}\circ 1{\otimes} 
F_{V, V},$$
$$G_{2,1}=G_{V_\lambda, V}{\otimes}1\circ G_{V_\lambda\underline{\otimes} V, V},$$
$$F_{2,1}=F_{V_\lambda\underline{\otimes} V, V}\circ F_{V_\lambda, V}{\otimes} 1,$$
$$\tilde{G}_{1,2}=1{\otimes} G_{V, V}
\circ G_{{V_\lambda, V\underline{\otimes} V}}.$$
$$F'_{1,2}=F_{V, V\underline{\otimes} V_\lambda}\circ 1{\otimes} F_{V, V_\lambda},$$
$$\tilde{G'}_{2,1}=G_{V, V} {\otimes}1
\circ G_{V\underline{\otimes} V, V_\lambda},$$
$$\tilde{F'}_{2,1}=F_{V\underline{\otimes} V, V_\lambda}\circ F_{V, V} {\otimes} 1,$$
$$G'_{1,2}=1 {\otimes}G_{V, V_\lambda}\circ G_{V, V\underline{\otimes} V_\lambda}.$$
Note that by Lemma \ref{k=1}, naturality of all the transformations and the braiding, and the two hexagonal 
equations (\ref{braided_symmetry1}), (\ref{braided_symmetry2}),
the map
$G_{2,1}$, ($F_{2,1}$ resp.), is conjugate to $G'_{1,2}$, ($F'_{1,2}$ resp.) via a specific braiding (that is the representative of the braid group element
$b_1b_2b_1=b_2b_1b_2$ in the category)
For example,
$$F_{V, V\underline{\otimes} V_\lambda}=c(V\underline{\otimes} V_\lambda, V)F_{V\underline{\otimes} V_\lambda, V}c_U(V\underline{\otimes} V_\lambda, V)^{-1}$$
$$1{\otimes} F_{V, V_\lambda} =1{\otimes} c(V_\lambda, V)\circ 1\otimes F_{V_\lambda, V}\circ 1\otimes c_U(V_\lambda, V)^{-1}$$
imply
$$F'_{1,2}=c(V\underline{\otimes} V_\lambda, V)c(V_\lambda, V)\otimes 1\circ F_{2,1}\circ (1\otimes c_U(V_\lambda, V) c_U(V\underline{\otimes} V_\lambda, V))^{-1}=$$
$$c(V\underline{\otimes} V_\lambda, V) c(V_\lambda, V)\otimes 1\circ F_{2,1}\circ (c_U(V\underline{\otimes} V_\lambda, V)  c_U(V_\lambda, V)\otimes 1)^{-1}.$$
 Multiplying together (\ref{a}) and (\ref{b}) and then (\ref{c}) and (\ref{d}) gives respectively
\begin{equation}\label{e}1= F_{2,1}\circ 1{\otimes} G_{V, V}
\circ P_{V_\lambda, V\underline{\otimes} V}\circ 1{\otimes} 
F_{V, V}\circ  G_{2,1},\end{equation}
\begin{equation}\label{f}1=F'_{1,2} \circ  G_{V, V} {\otimes}1
\circ P_{V\underline{\otimes} V, V_\lambda}\circ 
F_{V, V} {\otimes} 1\circ  G'_{1,2}.\end{equation}
Conjugating (\ref{f}) by the same braid group element gives
\begin{equation}\label{g}1=F_{2,1} \circ  1\otimes G_{V, V} 
\circ P^c_{V\underline{\otimes} V, V_\lambda}\circ 
1\otimes F_{V, V}  \circ  G_{2, 1},\end{equation}
where 
$$P_{ V_\lambda, V\underline{\otimes} V}=G_{ V_\lambda, V\underline{\otimes} V}F_{ V_\lambda, V\underline{\otimes} V},$$
$$P^c_{V\underline{\otimes} V, V_\lambda}=c_U(V\underline{\otimes} V, V_\lambda)\circ G_{V\underline{\otimes} V, V_\lambda}\circ F_{V\underline{\otimes} V, V_\lambda}\circ c_U(V\underline{\otimes} V, V_\lambda)^{-1}.$$
It follows from (\ref{e}) and (\ref{g}) that
\begin{equation}\label{h} 0=F_{2,1} \circ  1\otimes G_{V, V} 
\circ A \circ 
1\otimes F_{V, V}  \circ  G_{2, 1},\end{equation}
where $A=P_{V_\lambda, V\underline{\otimes} V}-P^c_{V\underline{\otimes} V, V_\lambda}$ is a morphism of the tilting category may be regarded a selfadjoint element 
of a $C^*$-algebra, hence it can be written as the difference of two orthogonal positive operators still morphisms of the tilting category
$$A=A_+-A_-, \quad A_+A_-=A_-A_+=0.$$
Being $G_{2,1}F_{2,1}$ an idempotent onto a maximal non-negligible submodule, we have by the Gelfand-Kazhdan properties (1)-(3) in Sect. \ref{74}
$$0=F_{2,1} \circ  1\otimes G_{V, V} 
\circ A_+ \circ 
1\otimes F_{V, V}  \circ  G_{2, 1}F_{2,1} \circ  1\otimes G_{V, V} 
\circ A_- \circ 
1\otimes F_{V, V}  \circ  G_{2, 1}.$$
The right hand side (\ref{h}) is selfadjoint, thus its square is positive. By the last equation, the square
takes the same form as (\ref{h}) with $A$ replaced by $A^2=(A_+)^2+(A_-)^2$. Both
$F_{2,1} \circ  1\otimes G_{V, V} 
\circ A_+^2 \circ 
1\otimes F_{V, V}  \circ  G_{2, 1}$ and $F_{2,1} \circ  1\otimes G_{V, V} 
\circ A_-^2 \circ 
1\otimes F_{V, V}  \circ  G_{2, 1}$ are positive, thus
 \begin{equation}
0=F_{2,1} \circ  1\otimes G_{V, V} 
\circ A_+^2 \circ 
1\otimes F_{V, V}  \circ  G_{2, 1}, \end{equation}
\begin{equation} 0=F_{2,1} \circ  1\otimes G_{V, V} 
\circ A_-^2\circ 
1\otimes F_{V, V}  \circ  G_{2, 1}\end{equation}
Taking the categorical trace we have
${\rm Tr}(A_+^2)={\rm Tr}(A_-^2)=0,$
hence $A_+=A_-=0$ by the $C^*$-property and it follows that $A=0$, that is 
$P_{V_\lambda, V\underline{\otimes} V}=P^c_{V\underline{\otimes} V, V_\lambda}$.

\end{proof}

The following proposition concludes the proof of Theorem \ref{main_wh}.

\begin{prop}\label{$2$-cocycle property} Let
$\overline{R}^U$ be the coboundary associated to $U_{q}({\mathfrak g})$ as   in 
Theorem \ref{U_q_as_a_Hermitian_ribbon_h}.
Then $\overline{R}^U\Delta(I)$ is a $2$-cocycle for the weak Hopf algebra $A_W({\mathfrak g}, q, \ell)$ as defined in Def.
\ref{def_2-cocycle}.
\end{prop}.

\begin{proof}
Formula (\ref{relation_between_the_two_coproducts}) in the statement  of Theorem \ref{main_wh} describes the relation between the coproducts of $U_q({\mathfrak g})$ and $A_W({\mathfrak g}, q, \ell)$.
Using this relation with 
$F=\overline{R}^U\Delta(I)$ $F^{-1}=\Delta(I)(\overline{R}^U)^{-1}$
   the left hand side of the first displayed equation of Def. \ref{def_2-cocycle},
   becomes 
   $$1\otimes\Delta(\Delta(I))[1\otimes\Delta^U((\overline{R}^U)^{-1})I\otimes (\overline{R}^U)^{-1}   \overline{R_U}\otimes I \Delta^U\otimes I(\overline{R}^U)]\Delta\otimes 1(\Delta(I)),$$
   with $\Delta^U$ the coproduct of $U_q({\mathfrak g})$ and $\Delta$ that of $A_W({\mathfrak g}, q, \ell)$.
   By Theorem \ref{U_q_as_a_Hermitian_ribbon_h}, $\overline{R}^U$ is a $2$-cocycle for $U_q({\mathfrak g})$, thus
   the central part of this equation equals $I\otimes I\otimes I$, and the proof of the first displayed formula of Def. \ref{def_2-cocycle} is complete. One similarly shows the second formula.
    (Note that by Remark \ref{how_it_acts}, see also Remark
\ref{remark_sect_25_on_the _2_cocycle_property}, the action of $\overline{R}^U\Delta(I)$
on $V_\lambda(q)\otimes V_\mu(q)$ on a simple component of highest weight $\gamma$, with $\lambda$, $\mu$, $\gamma\in\Lambda^+(q)$ is given by $R^U\Theta_w$, where $\Theta_w$
acts as
${q_0}^{\frac{\langle\lambda, \lambda+2\rho\rangle+\langle\mu, \mu+2\rho\rangle-\langle\gamma, \gamma+2\rho\rangle}{2}}$.)

\end{proof}

  \section{From $A_W({\mathfrak g}, q,\ell)$
  to a   compatible unitary coboundary wqh   structure on    $A(V_{{\mathfrak g}_k})$}\label{21}

  In this section we are interested in the module categories of affine vertex operator algebras $V_{{\mathfrak g}_k}$ at positive integer levels $k$. This is an important class
  of  vertex operator algebras associated to affine Lie algebras.
  Every vertex operator algebra $V$ has an   associated associative   algebra, called the Zhu algebra $A(V)$ \cite{Zhu}   recalled in Sect. \ref{VOAnets}. 
  We shall briefly recall  a natural identification of the Zhu algebra in the case of affine VOAs.
  
  The aim of this section is to apply the Drinfeld-Kohno theorem \ref{Drinfeld_Kohno} and then transport all the untwisted unitary coboundary structure from $A_W({\mathfrak g}, q_0, \ell)$ to the Zhu algebra $A(V_{{\mathfrak g}_k})$
via a beautiful continuous path argument that has been discovered by Wenzl in \cite{Wenzl}.
Using the methods of Sect. \ref{12}, we derive in this way a unitary modular tensor category structure on
${\rm Rep}(V_{{\mathfrak g}_k})$. In the second part of the
section, we describe this structure of  $A(V_{{\mathfrak g}_k})$ in more detail. 

Summarizing, in this section, we prove part 
(a) and (b) of the following main Theorem  \ref{Zhu_as_a_compatible_unitary_wqh}. Part (c) constitutes a statement about our
analogue of Kazhdan-Lusztig-Finkelberg theorem, Theorem \ref{Finkelberg_HL}, in the setting of affine vertex operator algebras
with Huang-Lepowsky ribbon braided tensor category structure, and will be   reformulated more precisely in the following section, see Theorem \ref{Finkelberg_HL}, where we shall give a first idea of proof connecting with work by Kirillov.  In Sect. \ref{22} we shall study the comparison between the  braided symmetry of ${\rm Rep}(V_{{\mathfrak g}_k})$ obtained in this section and the braided symmetry that arises from CFT on affine Lie algebras,
mostly following the exposition by Wassermann and Toledano-Laredo, in this way we  connect ideas of Wenzl on quantum groups with the setting of loop groups. In Sect. \ref{32} we describe the tensor category structure of 
${\rm Rep}(V)$ due to Huang and Lepowsky, in Sect. \ref{33} we shall study
the comparison between the associativity morphisms of  ${\rm Rep}(V_{{\mathfrak g}_k})$ obtained in this section with  Huang and Lepowsky associativity morphisms in the setting of vertex operator algebras. To do this, the notion 
of {\it primary field} of Knizhnik and Zamolodchikov will play a central role.  
  \bigskip
  
  \subsection{Affine Lie algebra $\hat{\mathfrak g}$,  affine vertex operator algebra   $V_{{\mathfrak g}_k}$, Zhu algebra 
  $A(V_{{\mathfrak g}_k})$.}\label{30.1}     In this subsection we recall some basic facts about affine Lie algebras and their associated affine vertex operator algebras. We are interested in the case of positive integer levels
  (see \cite{Kac2}, \cite{Frenkel_Zhu}).
  
  Let ${\mathfrak g}$ be a complex finite dimensional simple Lie algebra, ${\mathfrak h}$ a Cartan subalgebra,
  $\alpha_1,\dots,\alpha_r$ a set
of simple roots, and $A=(a_{ij})$ the associated Cartan matrix. Consider the unique invariant symmetric and bilinear form on ${\mathfrak h}^*$ such that 
\begin{equation}\label{normalization_inner_product_affine}
\langle\langle\theta,\theta\rangle\rangle=2\quad \text{for the highest root
  $\theta$}
  \end{equation}
Consider the affine Lie algebra $\hat{\mathfrak g}={\mathfrak g}\otimes{\mathbb C}[t, t^{-1}]\oplus{\mathbb C}{\bf  k}$,
with ${\bf k}$ in the center of $\hat{\mathfrak g}$ and Lie algebra structure given by
\begin{equation}\label{affine_Lie_algebra}
[a\otimes t^n, b\otimes t^m]=[a, b]\otimes t^{m+n}+{\bf k}\langle\langle a, b\rangle\rangle\delta_{m+n, 0}.\end{equation}
Let us fix $k\in{\mathbb C}$. Every ${\mathfrak g}$-module $W$ gives rise to a $\hat{\mathfrak g}$-module ${W}_k$ 
such that ${\bf k}$ acts as the scalar $k$. For a fixed irreducible ${\mathfrak g}$-module
$L(\lambda)$ with dominant weight $\lambda\in{\mathfrak h}^*$, corresponding $\hat{\mathfrak g}$-module
$L_{k, \lambda}$ is characterized up to isomorphism by the following three properties,\medskip

\noindent{\bf a)} $L_{k, \lambda}$ is irreducible, \smallskip

\noindent{\bf b)} ${\bf k}$ acts as $k$,\smallskip

\noindent{\bf c)}  $L_{k, \lambda}$ contains an isomorphic copy of
$L(\lambda)$   given by $\{a\in L_{k, \lambda}, \hat{\mathfrak g}_+a=0\}.$  \medskip

\noindent where 
$\hat{\mathfrak g}_+={\mathfrak g}\otimes{\mathbb C}[t]t$. By \cite{Frenkel_Zhu}, $V_{{\mathfrak g}_k}:=L_{k, 0}$ has the structure of a vertex operator algebra for 
$k\neq h^\vee$, the dual Coxeter number and when $k$ is a positive integer, $V_{{\mathfrak g}_k}$ is a rational VOA,
see also Sect. \ref{VOAnets}, \ref{VOAnets2} for more details and references to the original papers. 
By Theorem 3.1.2 in \cite{Frenkel_Zhu}, in this case the Zhu algebra $A(V_{{\mathfrak g}_k})$ is canonically isomorphic
to a quotient of $U({\mathfrak g})$ by the two-sided ideal generated by $e_\theta^{k+1}$, where
$e_\theta$ is an element in the root space ${\mathfrak g}_\theta$ of the maximal root $\theta$.
By Theorem 3.1.3 in \cite{Frenkel_Zhu}, the family of modules $L_{\lambda, k}$, where $\lambda$ is a dominant weight
in  $$\Lambda^+_k:=\{\lambda\in\Lambda^+: \langle\langle\lambda, \theta\rangle\rangle\leq k\}$$ is a complete list of irreducible
$V_{{\mathfrak g}_k}$-modules. 

Taking into account the different normalizations
(\ref{normalization_inner_product_qg}) and (\ref{normalization_inner_product_affine}) of the two inner products  of ${\mathfrak h}^*$,
we have
$$\langle \xi, \eta\rangle=d\langle\langle\xi, \eta\rangle\rangle.$$
Recalling the expression  of the dual Coxeter number $h^\vee$ given  in  (\ref{Coxeter_numbers}), one sees that
 the simple objects of ${\mathcal C}({\mathfrak g}, q, \ell)$
 summarized in Theorem \ref{irreducibles_of_fusion_category} and those of ${\rm Rep}(V_{{\mathfrak g}_k}$) are labelled by the same set $$\Lambda^+_k=\Lambda^+(q)$$ provided the order $\ell$ of $q^2$ and the level $k$ are related by
 $\ell=d(h^\vee+k)$ and $d|\ell$, see  Defs.   \ref{large_enough},  \ref{minimal_root}.

\bigskip

 \subsection{Compatible unitary coboundary structure on $A(V_{{\mathfrak g}_k})$, proof of     Theorem \ref{Zhu_as_a_compatible_unitary_wqh}, parts (a), (b)}\label{30.2}

In the following proof of Theorem \ref{30.2} and in the rest of the section, the root of unity $q$ in the statement   is renamed   $q_0$,  and  $q$   denotes complex numbers in a neighborhood of $q_0$.

\begin{proof} 
(a), (b) We show that the unitary coboundary weak  Hopf algebra $A_W({\mathfrak g}, q_0, \ell)$ satisfies the assumptions of the introductory part of the statement of Theorem \ref{Drinfeld_Kohno}, and those in parts a), b) of the same theorem, and  that the twisted algebra
 is canonically isomorphic to the Zhu algebra.
Let $A$ be a selfadjont operator on a finite dimensional Hilbert space, and let $A=A_+-A_-$ denote the spectral decomposition 
 of $A$ with $A_+$ and $A_-$ positive and $A_+A_-=0$. We set $$A^{1/2}=(A_+)^{1/2}+i(A_-)^{1/2},$$ where 
 $(A_+)^{1/2}$ and  $(A_-)^{1/2}$ are the positive square roots. If $A$ is invertible then 
 $[(A^{-1]})^{1/2}]^*=(A^{1/2})^{-1}.$ For $\lambda$, $\mu\in\Lambda^+(q_0)$, let $\overline{R}_{\lambda, \mu}$ be the
 selfadjoint invertible operator   on the Hilbert space $V_\lambda(q_0)\otimes V_\mu(q_0)$ defined in Prop. \ref{Wenzl_Hermitian_form} d). We set $T_{\lambda, \mu}=(\overline{R}_{\lambda, \mu})^{1/2}\Delta(I)$,  as
 an operator on the same Hilbert space with domain $\Delta(I)$ and similarly $T_{\lambda, \mu}^{-1}=\Delta(I)(\overline{R}_{\lambda, \mu}^{1/2})^{-1}$ with range $\Delta(I)$. This $T$ is a twist by Def. \ref{definition_of_twist}.
  Then 
 $$(T_{\lambda, \mu}^{-1})^*=(\overline{R}_{\lambda, \mu}^{-1})^{1/2}\Delta^{\rm op}(I)$$ by compatibility of $\Delta$ with $^*$. Thus
 $$(T_{\lambda, \mu}^{-1})^*_{21}=(\overline{R}_{\lambda, \mu})^{1/2}\Delta(I)=T_{\lambda, \mu}$$ as $(\overline{R}_{\lambda, \mu})_{21}=(\overline{R}_{\lambda, \mu})^{-1}$ by the same Proposition. 
 Let $P_{\lambda, \mu}$ and $Q_{\lambda, \mu}$ be the selfadjoint projections domain and range of ${\overline{R}_{\lambda, \mu}}_+$ and ${\overline{R}_{\lambda, \mu}}_-$ respectively, thus clearly $P_{\lambda, \mu}Q_{\lambda, \mu}=0$ and $P+Q=1$.
 Moreover,
 $$T_{\lambda, \mu}^*{T_{\lambda, \mu}}_-= T_{\lambda, \mu}^*(P-Q)T_{\lambda, \mu}=\Delta(I)^*((({\overline{R}_{\lambda, \mu}})^{-1})^{1/2})^{-1}(P-Q)(\overline{R}_{\lambda, \mu})^{1/2}\Delta(I)=$$
 $$\Delta(I)^*(\overline{R}_{\lambda, \mu})^{1/2}(\overline{R}_{\lambda, \mu})^{1/2}\Delta(I)
 =\Delta(I)^*\overline{R}_{\lambda, \mu}\Delta(I)=\overline{R}_{\lambda, \mu}\Delta(I),$$
 and this coincides with the $\overline{R}$-matrix of $A_W({\mathfrak g}, q_0, \ell)$. One  similarly shows that
 $({T_{\lambda, \mu}})_{-}^{-1}(T_{\lambda, \mu}^{-1})^*$ coincides with the left inverse of the $\overline{R}$-matrix of 
  $A_W({\mathfrak g}, q_0, \ell)$. Thus all the assumptions of Theorem \ref{Drinfeld_Kohno} a) hold, and we have the twisted   wqh structure on $(A_W({\mathfrak g}, q_0, \ell))_T$ with the properties of being unitary coboundary with compatible involution.
  For $\lambda\in\Lambda^+(q_0)$, let $V_\lambda$ be the classical representation of $U({\mathfrak g})$.
  We consider the linear isomorphism between the Hilbert spaces
  $$\phi_\lambda: V_\lambda(q_0)\to V_\lambda$$ taking an element $v(q_0)$ of the specialized Kashiwara-Lusztig basis
  of $V_\lambda(q_0)$ to $v(1)$, see Subsect. \ref{18.3} and Theorem \ref{Wenzl_positivity}.
By Remark \ref{compact_real_form}, the induced inner product of $V_\lambda$ is invariant with respect to
the compact real form of ${\mathfrak g}$.   
  Let $U_\lambda$ be the unitary part of the polar decomposition of $\phi_\lambda$.   Then $U_\lambda$ induces a canonical $^*$-isomorphism
of finite dimensional $C^*$-algebras 
  $${\rm Ad}(U_\lambda)_{\lambda\in\Lambda^+(q_0)}: (A_W({\mathfrak g}, q_0, \ell))_T\to A(V_{{\mathfrak g}_k})$$ 
that identifies the $^*$-involution 
  of $(A_W({\mathfrak g}, q_0, \ell))_T$ with the $^*$-involution of $A(V_{{\mathfrak g}_k})$ coming from
  the compact real form of ${\mathfrak g}$. Thus the $^*$-involution of
  $A(V_{{\mathfrak g}_k})$ is induced by the classical $^*$-involution of $U({\mathfrak g})$ via the quotient map.
  Via this algebra $^*$-isomorphism we transfer all the structure to $A(V_{{\mathfrak g}_k})$. 
  
  We verify the assumptions in part  b) of Theorem \ref{Drinfeld_Kohno} for $\rho=V(q_0)$ and $\sigma=V_\lambda(q_0)$,
  $\lambda\in\Lambda^+(q_0)$ and ${\mathfrak g}\neq E_8$. By the proof of Lemma 3.6.2 (b) and the proof
  of {\it Case 1} in the proposition that follows it in \cite{Wenzl}, $\overline{R}(q_0)$ is positive on the full tensor product Hilbert space $V_\lambda(q_0)\otimes V(q_0)$
  (with tensor product structure) thus $Q(q_0)=0$. Assume ${\mathfrak g}=E_8$.
   If $Q(q_0)p_\gamma(q_0)\neq0$ then $Q(q)p_{\gamma}(q)\neq0$ for $q$ in a neighborhood of $q_0$ by the continuity
   argument in the proof of {\it Case 2} in the same proposition in \cite{Wenzl}. On the other
hand in this neighborhood we may find values of $q$ for which $\overline{R}_{\lambda, \mu}(q)$ is positive by Lemma 3.6.2 a) in \cite{Wenzl}. It follows that $Q(q_0)p_{\gamma}(q_0)=0$ for $\gamma\in\Lambda^+(q_0)$.

\end{proof}
  
 \begin{rem}\label{non_uniqueness_Wenzl} In Theorems \ref{Zhu_from_qg_if_of_CFT_type0}, \ref{Zhu_from_qg_if_of_CFT_type} we describe the weak quasi-tensor structure $(Z, F_0, G_0)$ of parts (a) and (b)   theorem \ref{Zhu_as_a_compatible_unitary_wqh} in   detail. Note that $(F_0, G_0)$ is note unique as the same holds for the corresponding structure $(W, F, G)$ of Wenzl
 functor $W:{\mathcal C}({\mathfrak g}, q_0, \ell)\to{\rm Hilb}$. We have already noted that the possible variations on Wenzl functor are described by trivial twists of the associated weak  Hopf algebra.
 \end{rem}

  \begin{rem}
  We give another argument concerning the assumption in part (c)   of Theorem \ref{Drinfeld_Kohno}. 
Let $\lambda$, $\mu\in\Lambda^+(q_0)$, and  $F_{\lambda, \mu}(q_0): V_\lambda(q_0)\otimes V_\mu(q_0)\to V_\lambda\underline{\otimes} V_\mu(q_0)$ be
  projections onto canonical maximal non-negligible addenda defined in \ref{weak_tensor_structure}.

  We consider 
$F_{\lambda, \mu}(q): V_\lambda(q)\otimes V_\mu(q)\to V_\lambda\underline{\otimes} V_\mu(q)$
in the same vein as Wenzl path of idempotents $p_{\lambda}(q)$, with $q$ varying in ${\mathbb T}_{q_0, 1}$, see
Subsect. \ref{25.3}.  We have that the quantum Casimir acts by matrices with entries polynomials
in   ${\mathcal A}'$  specialized at $x^{1/L}\to q^{1/L}$ (the minimal $L$-th square root)
on the range of $F_{\lambda, \mu}(q)$ by complete reducibility, and since this holds for $R_{\lambda, \mu}(q)$ also, it follows that
$\overline{R}_{\lambda, \mu}(q)F_{\lambda, \mu}(q)$ is continuous in $q^{1/L}$. Thus for any continuous function $f:{\mathbb R}\to{\mathbb R}$ with $f(0)=0$, $f(\overline{R}_{\lambda, \mu})(q)F_{\lambda, \mu}(q)$ is continuous.
In particular, following \cite{Wenzl}, see also Subsect. \ref{18.3},  if $f(t)=-1$  for $t\in (-\infty, -a]$, for $a>0$ sufficiently small and $f(t)=0$ for $t\geq0$, then $f(\overline{R}_{\lambda, \mu})(q)$  is the spectral projection $Q(q)$
of $\overline{R}_{\lambda, \mu}$ corresponding to the negative eigenvalues for $q$
varying in a   neighborhood $U$ of $q_0$. Thus   $Q(q)F_{\lambda, \mu}(q)$ is continuous
in $U$.
If $Q(q_0)F_{\lambda, \mu}(q_0)\neq0$ then $Q(q)F_{\lambda, \mu}(q)\neq0$ in a neighborhood of $q_0$. On the other
hand in this neighborood we may find values of $q$ for which $\overline{R}_{\lambda, \mu}(q)$ is positive by Lemma 3.6.2 a) in \cite{Wenzl}. It follows that $Q(q_0)F_{\lambda, \mu}(q_0)=0$.

For ${\mathfrak g}\neq E_8$ the proof of the assumptions of c) of Theorem \ref{Drinfeld_Kohno} may be simplified.
Let $n_\mu$ be the smallest integer such that $V_\mu(q_0)$ is an addendum of a left parenthesized truncated power $V^{\underline{\otimes}n_\mu}(q_0)$ as in Def. \ref{projections_p_n}.
  We use induction on $n_\mu$. The case $n_\mu=1$, is the previous step.
  We want to verify that the possibly negative part $Q(q_0)$ of $R_{\lambda, \mu}(q_0)$ (with respect to the usual tensor product Hilbert space $V_\lambda(q_0)\otimes V_\mu(q_0)$) annihilates the range of $F_{\lambda, \mu}(q_0)$.
 Recall that $F_{\lambda, \mu}(q_0)=p_\lambda\underline{\otimes}p_\mu(q_0)p_\lambda(q_0)\otimes p_\mu(q_0)$.
By the $2$-cocycle
 property 
 of $\overline{R}(q_0)$ for $U_{q_0}({\mathfrak g})$, $\overline{R}_\lambda(q_0)\otimes \overline{R}_\mu(q_0)\overline{R}_{\lambda, \mu}(q_0)$ equals the restriction $\overline{R}_{n_\lambda+n_\mu}(q_0)$ of the form induced by $\overline{R}(q_0)$ with respect to the left parenthesized addendum $V^{\underline{\otimes}n_\lambda+n_\mu}(q_0)$ of $V^{n_\lambda+n_\mu}(q_0)$. By induction on $n_\mu$ and the previous step
($(V_\mu(q_0)=V(q_0)$), the negative part of $\overline{R}_{n_\lambda+n_\mu}(q_0)$ annihilates the range of $p_{n_\lambda+n_\mu}(q_0)$. On the other hand, this negative part equals the negative part of $\overline{R}_\lambda(q_0)\otimes \overline{R}_\mu(q_0) \overline{R}_{\lambda, \mu}(q_0)$ and this equals   
$\overline{R}_\lambda(q_0)\otimes \overline{R}_\mu(q_0) Q(q_0)$.
\end{rem}

\subsection{The coproduct of $A(V_{{\mathfrak g}_k})$. Universality property of $\boxtimes$-bifunctor
in ${\rm Rep}(A(V_{{\mathfrak g}_k}))$, derivation of   Frenkel-Zhu isomorphic images of the spaces of VOA intertwining operators}\label{30.3} We next describe more explicitly the weak bialgebra structure of the Zhu algebra $A(V_{{\mathfrak g}_k})$ of Theorem \ref{Zhu_as_a_compatible_unitary_wqh} (a), and we shall use this description to identify the structure with that arising from CFT, that is the braiding following the treatment of \cite{Kirillov3}, \cite{Wassermann}, \cite{Toledano_laredo} and the associator following Huang and Lepowsky tensor product structure,  in Sect. \ref{22}, \ref{23} respectively.

To compare with Huang-Lepowsky theory, we aim to show in more detail that the coproduct of $A(V_{{\mathfrak g}_k})$ induces
a tensor product bifunctor $\boxtimes$ in ${\rm Rep}(A(V_{{\mathfrak g}_k}))$ with the following property.
Any $A(V_{{\mathfrak g}_k}))$-module is also a $U({\mathfrak g})$-module using the realization of 
$A(V_{{\mathfrak g}_k}))$ as a quotient of $U({\mathfrak g})$ \cite{Frenkel_Zhu}.

\begin{thm}\label{FZ_condition}
For any pair $\lambda_1$, $\lambda_2$ of dominant weights in the open Weyl alcove, the  bilinear map 
$$F_0: V_{\lambda_1}\otimes V_{\lambda_2}\to V_{\lambda_1}\boxtimes V_{\lambda_2}$$ 
of Theorem \ref{Zhu_as_a_compatible_unitary_wqh}  intertwines the action of $U({\mathfrak g})$ and satisfies the following universality property.
For any $U({\mathfrak g})$-intertwining map
$$f_{FZ}: V_{\lambda_1}\otimes V_{\lambda_2}\to V_{\lambda_3}$$
satisfying Frenkel-Zhu   condition of Theorem 3.2.3 in \cite{Frenkel_Zhu}, there is a unique
$$\eta\in{\rm Hom}_{A(V_{{\mathfrak g}_k})}(V_{\lambda_1}\boxtimes V_{\lambda_2}, V_{\lambda_3})$$
such that the following diagram commutes

\[
 \begin{tikzcd}
V_{\lambda_1}\otimes V_{\lambda_2} \arrow{r}{F_0}\arrow[swap] {dr} {f_{FZ}} & {V_{\lambda_1}\boxtimes V_{\lambda_2}} \arrow{d}{{\eta}} \\
 & {V}_{\lambda_3}
 \end{tikzcd}
 \]

\end{thm}

\begin{rem}\label{FZ_condition_remark} Frenkel-Zhu condition was first proved by Tsuchiya and Kanie for the initial terms of primary fields (vertex operators in their terminology) in the case of the affine Lie algebra associated to
${\mathfrak sl}_2$ at level $k$ \cite{Tsuchiya_Kanie}. Similarly,  the relevance of Frenkel-Zhu condition in the statement of Theorem \ref{FZ_condition} is that 
the space of $U({\mathfrak g})$-intertwining maps $f_{FZ}$ that satisfies it is naturally isomorphic to
the space of primary fields for the affine Lie algebra associated to ${\mathfrak g}$ at level $k$, that we shall introduce later, see Prop. \ref{initial_term} and the comment following it. This is also 
a characterization under a natural isomorphism
of the space of {\it intertwining operators} ${{\mathcal M}^{W_3}_{W_1, W_2}}$ that
arise in the theory of vertex operator algebras for $W_i$ the $V_{{\mathfrak g}_k}$-module $L_{k, \lambda_i}$ associated to ${\lambda_i}$, see Subsect. \ref{30.1}. These spaces are basic building blocks for the construction
of Huang-Lepowsky vertex tensor category structure in ${\rm Rep}(V)$ under suitable conditions. We shall sketch their theory in Sects. \ref{32} and part of \ref{33}. The notation ${{\mathcal M}^{W_3}_{W_1, W_2}}$ is defined before Def. \ref{fusion_Rule_for_VOA_intertwiner}. The isomorphism is given by a correspondence between
three spaces, a suitable space of intertwining
operators, a corresponding space of primary fields and a corresponding space of initial terms of primary fields, that are the maps $f_{FZ}$.
Theorem \ref{FZ_condition} and these isomorphisms will be continued in Theorem \ref{tensor_product} in the setting of loop groups
and in Subsect. \ref{35.2} in the setting of affine vertex operator algebras at positive integer level. This series of theorems are all related and central in our discussion.  \end{rem}

Recall from Sect. \ref{20} that we have introduced a category ${\mathcal G}_{q_0}$ starting with a fundamental representation $V$ of
 ${\mathfrak g}$, see Def.   \ref{defn_of_G_q}, and that by  Theorem  \ref{equivalence_G_q_with_qg_fusion_category}, ${\mathcal G}_{q_0}$ is naturally
 equivalent to ${\mathcal C}({\mathfrak g},  q_0, \ell)$ for the suitable value of $q_0$. In this section, 
 we work with ${\mathcal G}_{q_0}$, and with abuse of notation we shall use the more standard notation ${\mathcal C}({\mathfrak g}, q_0, \ell)$ in place.
 With this variant,  every object $\rho$ is defined by construction as a decomposition into a direct sum of simple objects
 $V_{\lambda_\alpha}(q_0)$, with $\alpha\in F$ in a finite set  defined by isometries $S_\alpha\in (V_{\lambda_\alpha}(q_0), \rho)$.

  \begin{defn}\label{definition_of_the_equivalence}  If $\rho$ is a not necessarily simple object of the unitary fusion category ${\mathcal C}({\mathfrak g}, q, \ell)$ that
  is given as a specific orthogonal direct sum decomposition  $\rho=\oplus_{{\alpha}\in{F}} V_{\lambda_\alpha}(q_0)$
  (that is we have specific isometries $S_\alpha\in(V_{\lambda_\alpha}(q_0), \rho)$ satisfying $\sum_\alpha S_\alpha S_\alpha^*=1$) then we   extend the unitary isomorphism $U_\lambda: V_\lambda(q_0)\to V_\lambda$ defined
  in the proof of  Theorem \ref{Zhu_as_a_compatible_unitary_wqh} in a natural way to a linear unitary map between Hilbert spaces
  $$U_\rho: \rho\in{\mathcal C}({\mathfrak g}, q_0, \ell)\to \oplus_{\alpha\in F} V_{\lambda_\alpha}\in{\rm Rep}(A(V_{{\mathfrak g}_k}))$$ using the isometries $\{S_\alpha, \alpha\in F\}$, $U_\rho S_\beta=S'_\beta U_{\lambda_\beta}$, where $S'_\beta: V_{\lambda_\beta}\to  \oplus_{\alpha\in F} V_{\lambda_\alpha}$ is
  the canonical inclusion in the direct sum. If $T\in(\rho, \sigma)$ is a morphism, $\rho$ is defined by $S^\rho_\alpha$ and $\sigma$ by $S^\sigma_\beta$ then $(S^\sigma_\beta)^*TS^\rho_\alpha\in(V_{\lambda_\alpha}(q_0), V_{\lambda_\beta}(q_0))$ is nonzero only if $\lambda_\alpha=\lambda_\beta$, and is a scalar $t_{\alpha,\beta}$.
  We set ${\mathcal E}(\rho)=\oplus_{\alpha\in F} V_{\lambda_\alpha}$,
  ${\mathcal E}(T)=\sum_{\alpha,\beta} S'_\beta t_{\alpha, \beta} (S'_\alpha)^*=U_\sigma T U_\rho^*$. Then ${\mathcal E}$ is a linear $^*$-equivalence $${\mathcal E}: {\mathcal C}({\mathfrak g}, q_0, \ell)\to {\rm Rep}(A(V_{{\mathfrak g}_k})).$$
  \end{defn}

 Moreover, for any pair of simple objects $\lambda$, $\mu\in\Lambda^+(q_0)$ we have defined
 a fusion submodule $V_\lambda(q_0)\boxtimes V_\mu(q_0)$ and a weak tensor structure $(F_{\lambda, \mu}, G_{\lambda, \mu})$ of Wenzl functor $W:{\mathcal C}({\mathfrak g}, q_0, \ell)\to{\rm Hilb}$
 in Def. \ref{weak_tensor_structure}.

 The following result summarizes some properties of the weak quasi-tensor structure of Zhu's functor,
and the coproduct of the Zhu algebra transported from
the isomorphism $A(V_{{\mathfrak g}_k})\simeq^\phi (A_W({\mathfrak g}, q_0, \ell))_T$
of Theorem \ref{Zhu_as_a_compatible_unitary_wqh} (a).

\begin{thm}\label{Zhu_from_qg_if_of_CFT_type0}  The linear category ${\rm Rep}(A(V_{{\mathfrak g}_k}))$ becomes a pre-tensor category with tensor product   defined by
$$V_\lambda\boxtimes V_\mu={\mathcal E}(V_\lambda(q_0)\boxtimes V_\mu(q_0)). $$
Let $$Z: {\rm Rep}(A(V_{{\mathfrak g}_k}))\to{\rm Vec}$$ be the forgetful functor
associated to the Zhu algebra and 
let $(F_0, G_0)$ be the weak quasi-tensor structure of  $Z$
 defined by 
$$(F_0)_{\lambda, \mu}=U_{V_\lambda(q_0){\boxtimes}V_\mu(q_0)}(F_T)_{\lambda, \mu} U^{-1}_{V_\lambda(q_0)}\otimes U^{-1}_{V_\mu(q_0)}, $$
 $$(G_0)_{\lambda, \mu}=U_{V_\lambda(q_0)}\otimes U_{V_\mu(q_0)} (G_T)_{\lambda, \mu}  U_{V_\lambda(q_0){\boxtimes} V_\mu(q_0)}^{-1}$$
 Then
the linear equivalence ${\mathcal E}: {\mathcal C}({\mathfrak g}, q_0, \ell)\to {\rm Rep}(A(V_{{\mathfrak g}_k}))$ 
defined in \ref{definition_of_the_equivalence} satisfies $${\mathcal E}(V_\lambda(q_0))\boxtimes V_\mu(q_0))=
{\mathcal E}(V_\lambda(q_0))\boxtimes {\mathcal E}(V_\mu(q_0)), \quad\quad \lambda, \mu\in\Lambda^+(q_0).$$
Moreover,
for $\lambda$, $\mu\in\Lambda^+_k$,  
 \begin{itemize}
\item[(a)]
  $V_{\lambda}\boxtimes V_\mu$ is a ${\mathfrak g}$-invariant submodule and {\bf addendum} of 
  $V_\lambda\otimes V_\mu $, corresponding to the idempotent $(G_0)_{\lambda, \mu}(F_0)_{\lambda, \mu}$,
  \item[(b)] decomposition multiplicities of $V_{\lambda}\boxtimes V_\mu$ in the pre-tensor category ${\rm Rep}(A(V_{{\mathfrak g}_k}))$  equal 
those of the Grothendieck ring of ${\mathcal C}({\mathfrak g}, q_0, \ell)$,
 \item[(c)] The structure $(F_0, G_0)$ induces on $A({V_{{\mathfrak g}_k}})$ the coproduct described in the proof
 of Theorem \ref{Zhu_as_a_compatible_unitary_wqh},
 \item[(d)] In the special cases
  $$(F_0)_{\lambda, V}: V_\lambda\otimes V\to V_{\lambda}\boxtimes V, \quad\quad (G_0)_{\lambda, V}: V_{\lambda}\boxtimes V\to V_\lambda\otimes V$$ are respectively a ${\mathfrak g}$-invariant  orthogonal projection and   is a ${\mathfrak g}$-invariant inclusion
        with respect to
 the usual tensor product inner product of the classical compact real form of ${\mathfrak g}$.
\end{itemize}
\end{thm}

\begin{proof} (a) By construction,  $V_{\lambda}\boxtimes V_\mu$ is a $U({\mathfrak g})$-module that factors to a
$A(V_{{\mathfrak g}_k})$-module. The property that this is a submodule of $V_\lambda\otimes V_\mu $
follows from the fusion rules of affine vertex operator algebras, that have been shown in Theorem 3.2.3 in \cite{Frenkel_Zhu}.  (b)
We know that $F_{\lambda, \mu}$, $G_{\lambda, \mu}$ are morphisms in the tilting category relating the full tensor product
module with a canonical maximal non-negligible addendum.  The twisted weak quasi-tensor structure
$(F_T=FT^{-1}, G_T=TG)$    acts between the same spaces as $(F, G)$. By construction,  $F_0$ is an idempotent   onto a ${\mathfrak g}$-invariant addendum of
$V_\lambda\otimes V_\mu$ with the same decomposition multiplicities as $V_\lambda(q_0)\boxtimes V_\mu(q_0)$.
(c) For $\eta\in A_W({\mathfrak g}, q, \ell)={\rm Nat}(W)$, the isomorphism $\phi: {\rm Nat}(W)\to A(V_{{\mathfrak g}_k})={\rm Nat}(Z)$, with $Z:{\rm Rep}(V_{{\mathfrak g}_k})\to{\rm Vec}$ Zhu's functor, is given by
$\phi(\eta)_{V_\lambda}=U_\lambda\eta_{V_\lambda(q_0)}U_\lambda^{-1}$.  By naturality of $\eta$ and $\phi(\eta)$ and definition of $U_\rho$ on non-simple objects we also have $\phi(\eta)_{{\mathcal E}(\rho)}=U_\rho\eta_\rho U_\rho^{-1}$. The coproduct $\Delta$ of  $A(V_{{\mathfrak g}_k})={\rm Nat}(Z)$ is defined trasporting   the twisted coproduct $\Delta_T$ oft
$A_W({\mathfrak g}, q, \ell)={\rm Nat}(W)$ by $T$ via $\phi$, thus for $\rho$, $\sigma$ simple,
$$\Delta(\eta)_{{\mathcal E}(\rho), {\mathcal E}(\sigma)}=\phi\otimes\phi\circ\Delta_T(\phi^{-1}(\eta))_{{\mathcal E}(\rho), {\mathcal E}(\sigma)}=$$
$$U_\rho\otimes U_\sigma\circ (TG)_{\rho, \sigma}\circ U^{-1}_{\rho\otimes\sigma}(\eta_{{\mathcal E}(\rho)\otimes{\mathcal E}(\sigma)})\circ U_{\rho\otimes\sigma}(FT^{-1})_{\rho, \sigma}U^{-1}_\rho\otimes U^{-1}_\sigma.$$ (d) has been shown in the proof of Theorem \ref{Zhu_as_a_compatible_unitary_wqh}.

\end{proof}

\subsection{The associator of $A(V_{{\mathfrak g}_k})$}\label{30.4}
 We next consider the associator.

  \begin{thm}\label{Zhu_from_qg_if_of_CFT_type} The associator of $A(V_{{\mathfrak g}_k})$ obtained 
  transporting the associator of $A_W({\mathfrak g}, q_0, \ell)$ via
    the twist $T$ and isomorphism $\phi$ 
  $$A(V_{{\mathfrak g}_k})\simeq^\phi (A_W({\mathfrak g}, q_0, \ell))_T$$
   as in Theorem \ref{Zhu_as_a_compatible_unitary_wqh}  (a)
  is an associator satisfying the Def.
  \ref{CFT_type_associator} of     
    CFT-type ${\mathcal V}$-pre-associator $\Phi_{F_0, G_0}$, with ${\mathcal V}$
  the collection of triples
  of simple objects of $A(V_{{\mathfrak g}_k})$ of the form $(V_\lambda, V, V)$, $(V, V_\lambda, V)$,
  $(V, V, V_\lambda)$. 
  \end{thm}

   \begin{proof}
   Up to the isomorphism $\phi$, the structure of $B$ is induced by the twisted structure $(F_T=FT^{-1}, G_T=TG)$) for ${\mathcal F}$.  Let $\Delta$ and $\Phi$ be the coproduct and associator of $A_W({\mathfrak g}, q_0, \ell)$ induced by $(F, G)$. By  (\ref{twisted_associator}),  the associator 
   of $A_T$ is given by
   $\Phi_T= I\otimes T 1\otimes\Delta(T)\Phi\Delta\otimes1(T^{-1})T^{-1}\otimes I$ with $\Phi$ the associator of $A$.
   By the weak  Hopf property,
   $\Phi=1\otimes\Delta(\Delta(I))\Delta\otimes1(\Delta(I))$.
   It follows that 
  \begin{equation}\label{twisted_associator_cat}(\Phi_T)_{x,y,z}=I_{{\mathcal F}(x)}\otimes T_{y,z} (1\otimes\Delta(T))_{x,y,z}(\Delta\otimes1(T^{-1}))_{x,y,z}(T^{-1})_{x,y}\otimes I_{{\mathcal F}(z)}=\end{equation}
   $$I_{{\mathcal F}(x)}\otimes (TG)_{y,z} T_{x, y\otimes z}I_{{\mathcal F}(x)}\otimes F_{y,z} G_{x,y}\otimes I_{{\mathcal F}(z)} T^{-1}_{x\otimes y, z}(FT^{-1})_{x, y}\otimes I_{{\mathcal F}(z)} =$$
   $$I_{{\mathcal F}(x)}\otimes (G_T)_{y,z} (G_T)_{x, y\otimes z}F_{x, y\otimes z}I_{{\mathcal F}(x)}\otimes F_{y,z} G_{x,y}\otimes I_{{\mathcal F}(z)} G_{x\otimes y, z}(F_T)_{x\otimes y, z}(F_T)_{x, y}\otimes I_{{\mathcal F}(z)}$$
      (This formula may alternatively be derived from the categorical form of the twisted associator
given by Tannakian duality (\ref{associativity})  applied to $(F_T, G_T)$ and   the fact that $(F, G)$ is a weak tensor structure for ${\mathcal F}$.)
Assume   that two elements among $(x, y, z)$ are the fundamental representation $V$, and the third element is
a simple representation $V_\lambda$ in the  alcove. We claim that
$T=\overline{R}^{1/2}\Delta(I)$ on the full
tensor product Hilbert space $H_x\otimes H_y\otimes H_z$  acts as a $2$-cocycle for $A_W({\mathfrak g}, q_0, \ell)$ in the sense of Def. \ref{def_2-cocycle}. Writing down the definition of $2$-cocycle in categorical terms 
for the specific triple $(x, y, z)$ (that is, using
the Tannakain form of the coproduct induced by $(F, G)$, as in Sect. 5), and inserting the expression   in the center of the last line of (\ref{twisted_associator_cat})
  we have, following the notation of Def. \ref{CFT_type},
\begin{equation}\label{computation}(\Phi_T)_{x, y, z}=((G_T)_{1,2}(F_T)_{1,2}(G_T)_{2,1}(F_T)_{2,1})_{x,y,z}.\end{equation}
In detail, for ${\mathfrak g}\neq E_8$, by the first statement of  Lemma \ref{claim},
$\overline{R}^{1/2}$ on the full
tensor product Hilbert space $H_x\otimes H_y\otimes H_z$ satisfies the $2$-cocycle property
$$\overline{R}^{1/2}\otimes 1\Delta^U\otimes 1(\overline{R}^{1/2})=1\otimes \overline{R}^{1/2} 1\otimes\Delta^U(\overline{R}^{1/2}),$$
where $\Delta^U$ is the usual coproduct of $U_{q_0}({\mathfrak g})$.
Set $P=\Delta(I)$, with $\Delta$ the coproduct of $A_W({\mathfrak g}, q_0, \ell)$.
By Theorem \ref{main_wh} there is an epimorphism of algebras  
$$\pi: U_{q_0}({\mathfrak g})\to A_W({\mathfrak g}, q_0, \ell)$$
that has support the simple representations of $U_{q_0}({\mathfrak g})$ in the alcove
and   that satisfies
$$P\pi\otimes\pi(\Delta^U(a))=\Delta(\pi(a))=\pi\otimes\pi(\Delta^U(a))P.$$
It follows that on the full tensor product $H_x\otimes H_y\otimes H_z$ we have the $2$-cocycle property
with respect to $\Delta$:
$$1\otimes\Delta(\pi\otimes\pi(\overline{R}^{-1/2}))I\otimes \pi\otimes\pi(\overline{R}^{-1/2}) \pi\otimes\pi(\overline{R}^{1/2})\otimes I\Delta\otimes 1(\pi\otimes\pi(\overline{R}^{1/2}))=$$
$$I\otimes\Delta(I)\Delta(I)\otimes I=I\otimes (G_{y, z}F_{y,z}) (G_{x, y}F_{x, y})\otimes I.$$
We may then insert the left hand side of this equation in the center of the last line of (\ref{twisted_associator_cat})
and we get (\ref{computation}).
For ${\mathfrak g}=E_8$ we use the second statement of lemma \ref{claim}.

   \end{proof} 
   
     \begin{rem}
An associator on an arbitrary triple is uniquely determined by
its values on the family ${\mathcal V}$ of  triples $(V_\lambda, V, V)$, $(V, V_\lambda, V)$, $(V, V, V_\lambda)$ under some extra assumptions.
This claim  is stated and proved in Theorem \ref{claim1}.
The claim implies that   any other associator on $A(V_{{\mathfrak g}_k})$ that restricts to a CFT-type pre-associator associated to the same pair $(F_0, G_0)$,   it must coincide with the   associator of $A(V_{{\mathfrak g}_k})$  constructed in Theorem
  \ref{Zhu_as_a_compatible_unitary_wqh} (a) on every triple of representations. We shall find in Theorem \ref{HL_is_of_CFT_type} that also Huang-Lepowsky associativity
  morphisms restrict to the same CFT-type pre-associator $\Phi_{F_0. G_0}$ on ${\mathcal V}$.  
    \end{rem} 
   
  \begin{rem}
If one knew that this CFT-type pre-associator is an associator,
   then it must coincide with the former associator of $A(V_{{\mathfrak g}_k})$  constructed in Theorem
  \ref{Zhu_as_a_compatible_unitary_wqh}. It would also follow that $A(V_{{\mathfrak g}_k})$ becomes a weak Hopf algebra in this way.
  \end{rem} 
    
   The proof of the following lemma is   interesting. It     gives an explanation of the displacement
   of the spectral decomposition of the action of $\overline{R}^U$ on a tensor product of Weyl modules of $U_{q_0}({\mathfrak g})$ in the alcove with respect to the decomposition into irreducible components, met in the
   proof of Theorem \ref{Zhu_as_a_compatible_unitary_wqh},
   as a phenomenon
   arising from quantization, that in some sense fixes the position of the eigenspaces of $\overline{R}^U$ and displaces
   the simple components of a tensor product along Wenzl continuous arch on the circle from   $1$ to $q_0$.

   \begin{lemma}\label{claim}
Let $(x, y, z)$ be a triple of       modules of $U_{q_0}({\mathfrak g})$
 of the form $(V_\lambda, V, V)$, $(V, V_\lambda, V)$, $(V, V, V_\lambda),$
with $\lambda\in \Lambda^+(q_0)$.  Then   for ${\mathfrak g}\neq E_8$,
the positive square root   operator $T^U:=(\overline{R}^U)^{1/2}$ on $H_x\otimes H_y$ or $H_y\otimes H_z$
satisfies the $2$-cocycle property for $U_{q_0}({\mathfrak g})$
\begin{equation}\label{2-cocycle-for-the-square-root}T^U\otimes I\Delta^U\otimes 1(T^U)=I\otimes T^U 1\otimes\Delta^U(T^U)\end{equation}
on the full
tensor product Hilbert space $H_x\otimes H_y\otimes H_z$ with respect to well defined operators
 $\Delta^U\otimes 1(T^U)$ and $1\otimes\Delta^U(T^U)$.
 
For ${\mathfrak g}=E_8$, $T=T^U\Delta(I)$ satisfies the $2$-cocycle property for $A_W({\mathfrak g}, q_0, \ell)$ on the same Hilbert space as in Def. \ref{def_2-cocycle}.
    \end{lemma}
    
    \begin{proof}
We first assume ${\mathfrak g}\neq E_8$. 
By
c) of Theorem
\ref{U_q_as_a_Hermitian_ribbon_h}, $\overline{R}^U$
is a $2$-cocycle. By Prop. \ref{selfadjointness_three_coboundary} both sides of the $2$-cocycle equation
(\ref{2-cocycle_equation})
are selfadjoint. We claim that they are positive as operators on $H_x\otimes H_y\otimes H_z$.
By Lemma 3.6.2 (b) in \cite{Wenzl}, $\overline{R}^U$ is positive (and invertible) on $H_x\otimes H_y$ and $H_y\otimes H_z$ . It follows that $T^U$   is a positive operator
on each of these Hilbert spaces, with usual Hilbert space structure.  This implies  that the twisted coproduct
$\Delta_{T^U}^U(a)=T^U\Delta^U(a)(T^U)^{-1}$ of $U_{q_0}({\mathfrak g})$ commutes with the adjoint on the same spaces. On $H_x\otimes H_y\otimes H_z$, 
we have
$$\overline{R}^U\otimes I \Delta^U\otimes 1(\overline{R}^U)=T^U \otimes I\Delta^U_{T^U}\otimes 1(\overline{R}^U)T^U\otimes I=A^*A$$
$$I\otimes\overline{R}^U 1\otimes\Delta^U(\overline{R}^U)=I\otimes T^U 1\otimes\Delta^U_{T^U}(\overline{R}^U)I\otimes T^U=B^*B,$$ 
where 
$$A=\Delta^U_{T^U}\otimes 1(T^U)T^U\otimes I=T^U\otimes I\Delta^U\otimes 1(T^U),$$
$$B=1\otimes\Delta^U_{T^U}(T^U)I\otimes T^U=I\otimes T^U1\otimes\Delta^U(T^U),$$
and this   the claim follows. 

We want to show that $A=B$ on $H_x\otimes H_y\otimes H_z$.
We already know that by (\ref{2-cocycle_equation}), $A^*A=B^*B$. By uniqueness of the positive square root operator,
it suffices to show that $A$ and $B$ are positive operators. Let $S_i\in (V_i(q_0), H_x\otimes H_y)$ be   morphisms giving an orthogonal decomposition into simple objects with respect to the twisted inner product of $H_x\otimes H_y$. Positivity of $A$ is equivalent to positivity of the matrix $(S_i^\dagger\otimes I A S_j\otimes I)$,
with ${}^\dagger$ the adjoint with respect to the 
usual inner product of a tensor product Hilbert space. Recall that  if $^*$ denotes the adjoint with respect to the $\overline{R}^U$-twisted inner product then $S_i^*=S_i^\dagger\overline{R}^U$.
We have $$S_i^\dagger\otimes I A S_j\otimes I=S_i^\dagger\otimes IT^U\otimes I\Delta^U\otimes 1(T^U) S_j\otimes I=$$
\begin{equation}\label{conjugate}
(S_i^\dagger\otimes IT^U\otimes IS_j\otimes I) T^U=((S_i^*(T^U)^{-1}S_j)\otimes I)T^U.\end{equation}
Let $U_\rho$ be the unitaries defined as in Def. \ref{definition_of_the_equivalence}, here extended to a
 representation $\rho$ given as a direct sum of Weyl modules $V_{\lambda}(q_0)$ till the closure of the Weyl alcove, $\lambda\in\overline{\Lambda^+(q_0)}$ with respect to isometries $S_i$, $S'_i$. Set $\rho=x\otimes y$ and $p_i=S_iS_i^*$.
 Then $\Sigma\overline{R}^Up_i(\overline{R}^U)^{-1}\Sigma=:q_i$ is a $U_q({\mathfrak g})$-morphism in from $y\otimes x$ to itself. Let $p'_i:=U_{x\otimes y}p_i U_{x\otimes y}^{-1}$ and $q'_i:=U_{y\otimes x}q_i U_{y\otimes x}^{-1}$ be the corresponding $U({\mathfrak g})$-morphisms.
 Then $q'_i=\Sigma p'_i\Sigma$ by construction. It follows that  
 $p'_i$ commutes with $U_{x\otimes y}\overline{R}^UU_{x\otimes y}^{-1}$, and therefore also with $T^U$. Thus, unlike $\overline{R}^U$,
 the spectral decomposition of $U_{x\otimes y}\overline{R}^UU_{x\otimes y}^{-1}$ is compatible with the decomposition of the $U({\mathfrak g})$-representation corresponding to $x\otimes y$  into simple components, and the same holds for $U_{x\otimes y}T^UU_{x\otimes y}^{-1}$.
 
Let us multiply the last term of (\ref{conjugate}) on the left by $U_{i}\otimes U_z$ and on the right by $U_{j}^*\otimes U_z^*$. The relation $U_{i}S_i^*=S_i'^*U_{x\otimes y}$ and the previous observations imply that the first factor vanishes for $i\neq j$ and is a positive scalar otherwise. We are left to show that $T^U$ is positive on $V_i(q_0)\otimes H_z$ for all $i$, since then the right hand side of
(\ref{conjugate}) will be positive for all $i$. 
We repeat the computation in \ref{conjugate} with $\overline{R}^U\otimes I \Delta^U\otimes 1(\overline{R}^U)=A^*A$
in place of $A$. By positivity of $A^*A$ the computations gives that $\overline{R}^U$ is positive on 
$V_i(q_0)\otimes H_z$ for all $i$. Thus a positive solution $T^U$ on each $V_i(q_0)\otimes H_z$ may be found.

For ${\mathfrak g}=E_8$, $T=T^U\Delta(I)$   twists of the coproduct of $A_W({\mathfrak g}, q_0, \ell)$
into one that commutes with the adjoint,   as discussed in the proof of Theorem \ref{Zhu_as_a_compatible_unitary_wqh}.
We may slightly modify the previous proof to this case.
    \end{proof}
    
 \subsection{The braiding of $A(V_{{\mathfrak g}_k})$}\label{30.5}
   
   We finally determine the braiding of ${\rm Rep}(A(V_{{\mathfrak g}_k}))$.  By
      Prop.  \ref{braided_symmetry_with_generating_object},
   the braiding is completely determined by  the operators $c(V, V_\lambda)$, with $V$ the   generating representation,       provided the associator is determined. We thus restrict to compute these braiding operators.
 
            \bigskip
         
         \begin{thm}\label{braiding_Zhu_algebra} Let $ A(V_{{\mathfrak g}_k})$ be endowed with the unitary  coboundary wqh structure as in Theorem
\ref{Zhu_as_a_compatible_unitary_wqh} (a).
 Then on  the space of $V_\lambda\underline{\otimes} V$ the braiding
  of ${\rm Rep}(A(V_{{\mathfrak g}_k}))$ is given by 
  $$c(V_\lambda, V)=\Sigma e^{\frac{i\pi}{2(k+h^\vee)} A }$$ where $A$ is the selfadjoint operator with eigenvalue $\langle\langle\gamma, \gamma+2\rho  \rangle\rangle-\langle\langle\lambda, \lambda+2\rho  \rangle\rangle-\langle\langle\kappa, \kappa+2\rho  \rangle\rangle$ on
  a simple addendum $V_\gamma$ and
    $\kappa$ is the dominant weight of $V$ (the last addendum is replaced by the sum of the addenda of the dominant weights $\kappa_i$, $i=1, 2$ of $V$ in the type $D$ case). 
    
    \end{thm}

 \begin{proof}
 It follows   from 
 equation (\ref{data4}) in Remark \ref{unitary_braiding_data_of_strong_unitarization}
 that the braiding acts as $\Sigma\Delta(w^{-1})w\otimes w$ on a tensor product $\rho\otimes\sigma$ of 
 representations of a unitary coboundary wqh with compatible involution such that the Hermitian form $\Omega$
 is strongly trivial, $\Omega=\Delta(I)$ on this space. Part b) of
  Drinfeld-Kohno  theorem \ref{Drinfeld_Kohno} derives this property on the untwisted algebra derived from a general unitary coboundary wqh with compatible involution under suitable assumptions.
  The assumptions of this theorem in our application $A=A_W({\mathfrak g}, q_0, \ell)$ have been verified in 
  Theorem \ref{Zhu_as_a_compatible_unitary_wqh} for the representations $\rho=V_\lambda$, $\sigma=V$.
The value of $w$ is computed in
   part b), c) of Prop. \ref{Wenzl_Hermitian_form} and $q_0$ is the minimal root with $\ell=d(k+h^\vee)$ of Def. \ref{minimal_root}.

 \end{proof}
 
 Note that $$C_\gamma:=\langle\langle\gamma, \gamma+2\rho\rangle\rangle.$$ is the action of the classical Casimir
 operator of $U({\mathfrak g})$ on a simple representation $V_\gamma$. 
  
 \begin{rem}\label{braiding_Zhu_algebra2} We shall see in Cor. \ref{Toledano_laredo}
 Remark \ref{HL_braiding},
that the braided symmetry and ribbon structure of  ${\rm Rep}(A(V_{{\mathfrak g}_k}))$ of Theorem
 \ref{braiding_Zhu_algebra} arising from quantum groups
 coincides with the corresponding structure   induced on the Zhu algebra by the opposite braiding of ${\rm Rep}(V_{{\mathfrak g}_k})$ 
    in the setting of vertex operator algebras  \end{rem}

\section{Loop groups,  primary fields,
  tensor product  theorem for initial terms  \ref{tensor_product}}\label{22} 

 In Sect. \ref{21}, Subsect. \ref{30.1} we recalled the basis of affine Lie algebras, the connection with affine vertex operator algebras and the classification and correspondence between their irreducible representations.
 In  this section we    recall the basic elements of  affine Lie algebra CFT.

A first main result of this section is Theorem \ref{tensor_product}, where  we relate the spaces of initial terms of primary fields   with the tensor product bifunctor constructed in Theorem  \ref{Zhu_as_a_compatible_unitary_wqh} a).
 The constructions  of Verma modules for affine Lie algebras, see V.G. Kac \cite{Kac2} (10.4.6), and proofs of Prop. 2.1 and Theorem 2.3 in Tsuchiya and Kanie \cite{Tsuchiya_Kanie} should
 be regarded as helpful and preliminary material.

   Then with  Remark \ref{HL_braiding}.
 we identify     the braiding of ${\rm Rep}(V_{{\mathfrak g}_k})$  constructed
in Cor. \ref{cor_Zhu_as_a_compatible_unitary_wqh} and described explictly in Theorem \ref{braiding_Zhu_algebra}
with that arising from   the setting of CFT,   following the work of  Knizhnik-Zamolodchikov \cite{KZ}  and Tsuchiya-Kanie \cite{Tsuchiya_Kanie} and focusing on the subsequent developments by Wassermann  \cite{Wassermann} and Toledano-Laredo, \cite{Toledano_laredo}, in   the analytic functional setting. It will suffice to apply this identification to special pairs of the form $(V, V_\lambda)$ and $(V_\lambda, V)$.
These   works are in the setting of loop groups, and we start  recalling the basic terminology and some of their results. 
 (See  also the works by Gui, \cite{GuiI}, \cite{GuiII}, \cite{GuiIII}, \cite{GuiIV}, 
       \cite{GuiV} which expand these results and include  the setting of vertex operator algebras.)
       
The comparison between the   braided symmetry of ${\rm Rep}(A(V_{{\mathfrak g}_k}))$ and the braiding properties
in the setting of loop groups will be achieved using
 the following Corollary \ref{Toledano_laredo}, discussed in \cite{Toledano_laredo}, \cite{Wassermann},\ref{HL_braiding}.

The tensor product bifunctor of Huang-Lepowsky    \cite{Huang_LepowskiIII}
$(W_1, W_2)\to W_1\boxtimes_{\rm HL}W_2$ is explicitly determined by the irreducible modules and {\it intertwining operators} between irreducible modules, and in the case of interest  ${\rm Rep}(V_{{\mathfrak g}_k})$, this is in turn closely related to the study of initial terms of primary fields. We shall describe these connections
in Sect. \ref{32}, \ref{33}.

\bigskip

 \subsection{Classification of simple projective unitary representations of $LG$, connection with the VOA setting.}\label{32.1}
 A simple $V_{{\mathfrak g}_k}$-modules $L_{\lambda, k}$ corresponds to a {\it projective
   unitary representation} $\pi$ of the loop group $LG=C^\infty(G, {\mathbb T})$ on a Hilbert space ${\mathcal H}$  {\it of positive energy}. 
   This is a strongly continuous representation $\pi: LG\to {\rm PU}({\mathcal H})$ that extends to the semidirect product
   $LG\rtimes{\rm Rot}({\mathbb T})$, where ${\rm Rot}({\mathbb T})$ is the rotation group of the
   circle which acts on $LG$ by translation. It is required that the infinitesimal generator $d$ of rotations acts
   with spectrum bounded below and has finite dimensional subspaces. 
          A simple projective representation $\pi$ is uniquely determined by an integer $k$, the level, and the top level ${\mathcal H}[0]$, a simple representation  of  
   $G$,  called the {\it lowest energy subspace} in \cite{Toledano_laredo}.
For the sake of completeness we summarize the 
          classification result.
   We have a decomposition of ${\mathcal H}$ as a Hilbert space direct sum
   $${\mathcal H}=\bigoplus^{\rm HS}_{n\in{\mathbb N}}{\mathcal H}[n],$$
   with ${\mathcal H}[n]$ is the subspace of
   ${\mathcal H}$ such that the rotation $\pi(R_\theta)$ by $\theta$ acts as multiplication by $e^{in\theta}$.
   Consider the algebraic direct sum $${\mathcal H}^{\rm fin}=\bigoplus_{n\in{\mathbb N}}{\mathcal H}[n]$$ as a representation of the 
    Kac-Moody algebra, the 
  semidirect product $\hat{\mathfrak g}\rtimes{\mathbb C}d$   (it corresponds to ${\mathfrak g}_{\mathbb C}$ in Kac book \cite{Kac2} and $\hat{{\mathfrak g}_{\mathbb C}}$ in \cite{Toledano_laredo}). 
     This space is a core for the infinitesimal selfadjoint generator $d$ of the rotation group $\pi(R_\theta)=e^{i\theta d}$.
     The operator $d$ acts as multiplication by $n$ on ${\mathcal H}[n]$.
     We write $$a{(m)}:=a\otimes t^m, \quad\quad a\in{\mathfrak g}.$$ An element $X\in{\mathfrak g}\otimes{\mathbb C}[t, t^{-1}]$
     gives rise to a one-parameter projective group in $LG$  and therefore to an operator $\pi(X)$ via Stones' theorem which has
      ${\mathcal H}^{\rm fin}$ as a core  and is determined up to a constant. The constant may be fixed by a suitable convention, see    Theorem 1.2.1 of \cite{Toledano_laredo}.
      $X\to\pi(X)$, $d\to d$ is a unitarizable  representation of  the Kac-Moody algebra on
       ${\mathcal H}^{\rm fin}$. Thus relation (\ref{affine_Lie_algebra}) holds in the representation space, and the central element ${\bf k}$ acts as a positive integer $k$, the level.
       Fusion tensor product in this setting has been studied by Wassermann \cite{Wassermann} for $G={\rm SU}(N)$,
   Toledano Laredo for $G={\rm Spin}(2n)$ \cite{Toledano_laredo}. The other Lie types and
   connections between Connes fusion tensor product of the category of finite projective representations
   of $LG$ of finite energy and the tensor product of the category of $V_{{\mathfrak g}_k}$-modules
   of the corresponding affine vertex algebra with tensor product structure of Huang and Lepowsky    of affine vertex operator algebras  has been given by Gui \cite{GuiI},
   \cite{GuiII}, \cite{GuiIII}, \cite{GuiIV}, \cite{GuiV}.  Gui described
     applications to a comparison between representation categories of conformal nel and VOAs.
   In the setting of loop groups we shall need some formulas
   regarding the braiding that may be found in   \cite{Wassermann} and
   \cite{Toledano_laredo}. These formulas hold also in the other cases
      by the mentioned developments by Gui.

      \bigskip

\subsection{The conformal Hamiltonian $L_0$}\label{32.2}   Let $C\in U({\mathfrak g})$ be the   Casimir element, that is $C=\sum_iX_iX^i$, where $X_i\in{\mathfrak g}$ is a basis, $X^i$ a dual basis
   with respect to the   form $\langle\langle \xi, \eta \rangle\rangle$. Since $C$ is a central element, for any representation $\pi$ of ${\mathfrak g}$
   $$C_\pi=\pi(C)\in(\pi, \pi).$$
   On an irreducible highest weight representation $V_\lambda$, $C_{V_\lambda}$ acts as the positive scalar
   $$C_\lambda=\langle\langle \lambda, \lambda+2\rho\rangle\rangle.$$

   For $X\in{\mathfrak g}$, we define the formal Laurent series
   $$X(t)=\sum_{n\in{\mathbb Z}} X(n)t^{-n-1}.$$
   A construction by Segal-Sugawara (\cite{Segal2}, \cite{Sugawara})
 gives gives the Virasoro (or {energy-momentum) field,
     $$T(t)=\sum_{n\in{\mathbb Z}} L_n t^{-n-2}.$$
 Then the endomorphisms act on the representation space and    satisfy the Virasoro algebra relations
   $$[L_m, L_n]=(m-n)L_{m+n}+\frac{c_k}{12}(m^3-m)\delta_{-m, n},$$
   where $c_k$ is called the central charge and is a known constant \cite{Frenkel_Zhu}.
   The element $$L_0$$ is the {\it conformal Hamiltonian}, and in our case is given by a special case of the Segal-Sugawara
   construction, Sect. 9.4 in \cite{Pressley_Segal}, see also Sect. 1.2 in \cite{Toledano_laredo},
   $$L_0=\frac{1}{k+h^\vee}(\frac{1}{2}\sum_i X_i(0)X^i(0)+\sum_{m>0}\sum_i X_i(-m)X^i(m)).$$
   On the lowest energy space   ${\mathcal H}[0]$ of a simple level $k$ projective positive energy representation  of $LG$
   corresponding to the weight $\lambda$ in the alcove $\Lambda^+_k$,
   (previously denoted the top space $L_{k, \lambda}[0]=L(\lambda)$ of the semple $V_{{\mathfrak g}_k}$-module $L_{k, \lambda}$,   in the VOA setting) the conformal Hamiltonian $L_0$ acts as the scalar
\begin{equation}\label{scalar}\Delta_\lambda=\frac{C_\lambda}{2(k+h^\vee)}.
\end{equation}
 
Moreover, $L_0$ satisfies the same commutation relation as $d$, that is $[L_0, X(n)]=-nX(n)$, thus
$$L_0=d+\Delta_\lambda$$
on $L_{k, \lambda}$. It follows that $\Delta_\lambda$ is the minimal eigenvalue of $L_0$ on $L_{k, \lambda}$.

   \bigskip
   
   \subsection{Primary field, initial term.}\label{32.3} The important notion of {\it primary field} was introduced by   Knizhnik and Zamolodchikov \cite{KZ}, continued by Tsuchiya and Kanie \cite{Tsuchiya_Kanie}. Our exposition in the setting of loop groups is influenced by Toledano-Laredo thesis \cite{Toledano_laredo}.  
   
   Let $V$ be a finite dimensional unitary representation of $G$. Then $V[t, t^{-1}]$ is a representation of
   ${\mathfrak g}\otimes{\mathbb C}[t, t^{-1}]$ that extends to ${\mathfrak g}\otimes{\mathbb C}[t, t^{-1}]\rtimes{\rm Rot}({\mathbb T})$ by $R_\theta f(t)=e^{-i\theta}t$,
    $$X(m)v(n)=(Xv)(m+n)\quad\quad dv(n)=-nv(n).$$
   
   \begin{defn}\label{primary_field}
   Let ${\mathcal H}_i$, ${\mathcal H}_j$ be simple (positive energy) representations of $LG$ of level $k\in{\mathbb N}$
   and let $V_h$ be a unitary simple $G$-module.
   A {\it primary field} is a linear map
   $$\phi: V_h[t, t^{-1}]\otimes {\mathcal H}_i^{\rm fin} \to {\mathcal H}_j^{\rm fin}$$
   that intertwines the action of ${\mathfrak g}\otimes{\mathbb C}[t, t^{-1}]\rtimes{\rm Rot}({\mathbb T})$.
   The $G$-module $V_h$ is called the {\it charge} of $\phi$, the representations ${\mathcal H}_i$ and ${\mathcal H}_j$
   are called the {\it source} and {\it target} respectively. 
   \end{defn}

   We set, for $v\in V_h$,
   $$\phi(v, n)=\phi(v\otimes t^n): {\mathcal H}_i^{\rm fin}\to {\mathcal H}_j^{\rm fin}.$$
   We represent $\phi$ as a formal operator-valued distribution
   $$\phi(v, t)=\sum_{n\in{\mathbb Z}} \phi(v, n) t^{-n-\Delta_\phi},$$
   where
   $$\Delta_\phi:=\Delta_i+\Delta_h-\Delta_j,$$ is called the {\it conformal weight} (or conformal dimension, scaling dimension) of the field $\phi$,
   and    $\Delta_i$, $\Delta_h$, $\Delta_j$ are defined as in (\ref{scalar}).
   The intertwining relation of $\phi$ can be written as
   \begin{equation}\label{primary_field_eq}[X(m), \phi(v, t)]=\phi(Xv, t)t^m, \quad\quad [d, \phi(v, t)]=(t\frac{d}{dt}+\Delta_\phi)\phi(v, t).\end{equation}
   The following form for (\ref{primary_field_eq}) is also useful,
\begin{equation}\label{primary_field_eq2} [X(m), \phi(v, n)]=\phi(Xv, m+n), \quad\quad [d, \phi(v, n)]=-n\phi(v, n).\end{equation}
\begin{defn}\label{gauge_condition} 
   The   equation on the left hand side of (\ref{primary_field_eq}) 
   or (\ref{primary_field_eq2}) is called {\it gauge condition}, see \cite{Tsuchiya_Kanie}.
   \end{defn}
   The equation of the right hand side of (\ref{primary_field_eq}) 
   or (\ref{primary_field_eq2}) also has an interesting interpretation, see Remark \ref{equation_of_motion}.
It implies the following important grading relations,
\begin{equation}\label{grading_loop_group}
\phi(v, n): {\mathcal H}_i[k]\to{\mathcal H}_j[k-n], \quad\quad k\in{\mathbb Z}_{\geq0}, \quad n\in{\mathbb Z}.\end{equation}
It follows in particular that
$ \phi(v, 0) $ restricts to a $G$-intertwiner  
$$T_\phi:  V_h\otimes {\mathcal H}_i[0]\to{\mathcal H}_j[0].$$
\begin{defn} \label{initial_term_def} The $G$-intertwiner $T_\phi$  is called the {\it initial term} of $\phi$.
 We denote by $I_G(V_h\otimes{\mathcal H}_i[0], {\mathcal H}_j[0])$ the subspace of
${\rm Hom}_G(V_h\otimes{\mathcal H}_i[0], {\mathcal H}_j[0])$ of initial terms $T_\phi$ of primary fields with charge
the irreducible $G$-representation $V_h$, source ${\mathcal H}_i$ and target ${\mathcal H}_j$.
\end{defn}

\begin{prop} \label{initial_term} (Prop. 2.1  in \cite{Tsuchiya_Kanie})
The initial term $T_\phi$ of a primary field $$\phi: V_h[t, t^{-1}]\otimes{\mathcal H}_i^{\rm fin}\to {\mathcal H}_j^{\rm fin}$$  of  fixed charge $V_h$ acting between fixed irreducible representations
uniquely determines $\phi$. Thus the linear map  
  $$\phi\to T_\phi\in{\rm Hom}_G(V_h\otimes {\mathcal H}_i[0], {\mathcal H}_j[0])$$ is faithful.  
\end{prop}

\begin{rem} By Theorem 2.3 in \cite{Tsuchiya_Kanie} and Theorem 3.2.3 in \cite{Frenkel_Zhu}, a quotient relation is well known that
 describes morphisms of ${\rm Hom}_G(V_h\otimes {\mathcal H}_i[0], {\mathcal H}_j[0])$ that arise as initial terms of primary fields. This quotient relation is a property of the irreducible representation  ${\mathcal H}_i^{\rm fin}$. It is 
 derived from a presentation of the ideal of a Verma module  at level the positive integer $k$
 with quotient the irreducible ${\mathcal H}_i^{\rm fin}$ (or ${\mathcal H}_j^{\rm fin}$). 
 \end{rem}

For all ${\mathfrak g}$, the  linear space $I_G(V_h\otimes {\mathcal H}_i[0], {\mathcal H}_j[0])$  has dimension given by the dimension of the corresponding primary fields  by Prop. \ref{initial_term}. The term 
{\it fusion rules} in the setting of vertex operator algebras is used to denote this dimension,
 see Def. \ref{fusion_Rule_for_VOA_intertwiner}, but   also the linear space to which it refers.

 Following \cite{Tsuchiya_Kanie},   only for this section, we use the notation of $\phi$ as a {\it vertex}, on top we put
the charge of $\phi$.
$$\phi^h_{j, i}, \quad {\rm or } \quad {{h}\choose{ji}}$$
 to represent a primary field with initial term in ${\rm Hom}_G(V_h\otimes V_i, V_j)$.
    In the sections regarding vertex operator algebras Sect. \ref{32}, \ref{33}, we shall conform to the
   notation of intertwining operators, with the charge on the bottom left of the vertex, the source on the bottom right
   and the target on top.

   \begin{rem} \label{equation_of_motion}
The equation on the right hand side of of (\ref{primary_field_eq}) 
   (or \ref{primary_field_eq2}) appears in the literature written in terms of $L_0$, as
  \begin{equation}[L_0, \phi(v, t)]=(t\frac{d}{dt}+\Delta_h)\phi(v, t),\end{equation}
   with $\Delta_h$ the minimal eigenvalue of $L_0$ on the level $k$ representation with top space $V_h$.
   This equation   is  
   the part for $m=0$   
   of the {\it equation of motion} of \cite{KZ}, \cite{Tsuchiya_Kanie} that involves the whole Virasoro field $L_m$, $m\in{\mathbb Z}$. In the case of vertex operator algebras, the corresponding equation involving the Virasoro field
   follows from the Jacobi identity for intertwining operators.
   We shall comment more on  this in Sect. \ref{33}, see (\ref{equation_of_motion_VOA}).  

 \end{rem}

  The classification of irreducible representations of the affine Lie algebra $\hat{{\mathfrak g}}$ at level $k\in{\mathbb N}$
 via their restriction to   (irreducible) representations of the corresponding classical Lie algebra with dominant weights in the open Weyl alcove $\Lambda_k^+$,
 is discussed in \cite{Frenkel_Zhu}.

\subsection{The tensor product Theorem \ref{tensor_product}}\label{32.4}
  
  \begin{thm}\label{tensor_product} 
  {\bf a)} There is a   connection between \smallskip
  
 \noindent{\bf i)} the subspace $I_G(V_h\otimes{\mathcal H}_i[0], {\mathcal H}_j[0])$  of ${\rm Hom}_G(V_h\otimes{\mathcal H}_i[0], {\mathcal H}_j[0])$
  of initial terms
  of charged primary fields with charge an irreducible $G$-representation $V_h\in\Lambda_k^+$, associated to simple positive energy representations
  of $L(G)$ at level $k\in{\mathbb N}$ and\smallskip
  
\noindent{\bf  ii)} the tensor product bifunctor $\boxtimes$ of the representation category  
of the Zhu algebra   $A(V_{{\mathfrak g}_k})$ defined by the weak quasitensor structure $(F_0, G_0)$ of
  the forgetful functor of ${\rm Rep}(A(V_{{\mathfrak g}_k}))$ derived from the weak quasi-bialgebra structure 
  of $A(V_{{\mathfrak g}_k})$ arising from
   the  corresponding quantum group
  at a root of unity $U_{q_0}({\mathfrak g})$
 as in Theorem \ref{Zhu_as_a_compatible_unitary_wqh}.
 The connection is given as follows.

  {\bf b)}  The choice of a weak quasi-tensor structure $(F_0, G_0)$ for Zhu functor $Z$
  as in Theorem \ref{Zhu_as_a_compatible_unitary_wqh}
 induces the following description of spaces of initial terms,
 \begin{equation}\label{initial_term_Zhu} 
 \{T\in{\rm Hom}_G(V_h\otimes{\mathcal H}_i[0], {\mathcal H}_j[0]): TG_0F_0=T\}\simeq I_G(V_h\otimes{\mathcal H}_i[0], {\mathcal H}_j[0]). \end{equation}

An element $T\in I_G(V_h\otimes{\mathcal H}_i[0], {\mathcal H}_j[0])$ decomposes as
 $$T=\hat{T} F_0,$$ 
 with $$\hat{T}=TG_0\in{\rm Hom}_{A(V_{{\mathfrak g}_k})}(V_h\boxtimes{\mathcal H}_i[0], {\mathcal H}_j[0])\subset 
 {\rm Rep}(A(V_{{\mathfrak g}_k})).$$
 The linear map $$T\in I_G(V_h\otimes{\mathcal H}_i[0], {\mathcal H}_j[0])\to \hat{T}=TG_0\in {\rm Hom}_{A(V_{{\mathfrak g}_k})}(V_h\boxtimes{\mathcal H}_i[0], {\mathcal H}_j[0])$$ is an isomorphism
 with inverse 
 $$\hat{T}\to T=\hat{T}F_0.$$
 
  {\bf c)} Assume that one of the charge $V_h$ or the   eigenspace ${\mathcal H}_i[0]$ of the source 
with minimal eigenvalue of $L_0$
 is the fundamental representation $V$ of ${\mathfrak g}$, chosen as in \cite{Wenzl} and that ${\mathfrak g}\neq E_8$. Then $G_0F_0$ is the unique (selfadjoint) idempotent corresponding to Wenzl   idempotent onto the unique maximal non-negligible addendum of $V_h(q_0)\otimes V_i(q_0)$ of $U_{q_0}({\mathfrak g})$ via the twist construction and continuous path as described
 in Sect. \ref{21}. 
 We have 
 \begin{equation}\label{initial_term_Zhu2} I_G(V_h\otimes{\mathcal H}_i[0], {\mathcal H}_j[0])=
 {\rm Hom}_G(V_h\otimes {\mathcal H}_i[0], {\mathcal H}_j[0]).
 \end{equation}

  \end{thm}

\begin{proof}   {\bf a)}  A characterization of   $I_G(V_h\otimes {\mathcal H}_i[0], {\mathcal H}_j[0])$  
  is given in
Theorem 2.3 in \cite{Tsuchiya_Kanie} for ${\mathfrak g}={\mathfrak sl}_2$, \cite{Wassermann} for 
 ${\mathfrak g}={\mathfrak sl}_N$,    Prop. 4.1. in \cite{Toledano_laredo} for the other Lie types in the setting of affine Lie algebras or loop groups; 
 and Theorem 3.2.3 in \cite{Frenkel_Zhu} in the setting of affine vertex operator algebras.
 It follows automatically that charge space $V_h$   belongs to the Weyl alcove $\Lambda_k^+$.
 
 {\bf b)} The fact that $\hat{T}$ is a morphism in the category of modules of the Zhu algebra follows from the identification
of the Zhu algebra recalled in Sect. \ref{21}, Subsect. \ref{30.1}, and the isomorphism property is a consequence of $F_0G_0=1$.

 {\bf c)} It follows from the characterization result  that if one of the charge $V_h$ or the   eigenspace ${\mathcal H}_i[0]$ of the source 
with minimal eigenvalue of $L_0$
 is the fundamental representation $V$ of ${\mathfrak g}$, chosen as in \cite{Wenzl},
and if $i$ and $j$ vary in   $\Lambda_k^+$,
 then
 $I_G(V_h\otimes{\mathcal H}_i[0], {\mathcal H}_j[0])={\rm Hom}_G(V_h\otimes {\mathcal H}_i[0], {\mathcal H}_j[0])$  for ${\mathfrak g}\neq E_8$. 

Proposition at page 274 in \cite{Wenzl} gives a decomposition, for all ${\mathfrak g}$, of $V_h\otimes {\mathcal H}_i[0]$ as a direct sum of a complemented submodule that played a central role    in Theorem \ref{Zhu_as_a_compatible_unitary_wqh} and   denoted by $V_h\boxtimes{\mathcal H}_i[0]$.
This submodule   is defined by the projection $F_0$ and   inclusion $G_0$ corresponding to Wenzl construction in the setting of quantum groups via
the twist and continuous path. It follows that (\ref{initial_term_Zhu}) holds. In particular, 
(\ref{initial_term_Zhu2}) follows from the fusion rules of tensor product $V_\lambda\otimes V$ in ${\rm Rep}({\mathfrak g})$  where $V$ denotes the fundamental representation of ${\mathfrak g}$ as explained in \cite{Wenzl} and already used in this paper.

    \end{proof}

    \begin{rem}\label{Primary_following_TL}
  Primary fields   play a central role in  \cite{KZ}, \cite{Tsuchiya_Kanie},
and are  also called vertex operators. 

On the other hand, the notion of vertex operator is also used in the setting of vertex operator algebras with a different but related meaning, see Sect. \ref{VOAnets}, \ref{VOAnets2}. To avoid confusion we shall use the  term vertex operator only in the vertex operator algebra setting. (Please note also the historical Remark 5.2 in \cite{Kac} on vertex operators and reference to the original articles.)

In Sect. \ref{33}, we shall    describe the           corresponding   notion of primary field  as in
  Def. \ref{primary_field}, in the setting of      vertex operator algebras, see Def. \ref{primary_field_VOA}.
  
    In the case of affine vertex operator algebras,   primary fields and their initial terms   play a central role for our purposes,
 the study of the connection between Huang-Leposwky tensor product theory and tensor product of the representation category of the Zhu algebra
  as described in Theorem \ref{Zhu_as_a_compatible_unitary_wqh}. These developments are considered in  Sect. \ref{33}.
  
  More precisely, in   Subsect. \ref{35.1}, we shall introduce a subclass of the class of primary fields now  in the setting of vertex operator algebras as defined in
\cite{Kac},  \cite{FHL}. This is the subclass   defined by a vector in the top space of the charge module and as said it corresponds to Def. \ref{primary_field}
 in the setting of loop groups. We shall also introduce
their {\it initial term}.

As said, in Subsect. \ref{35.2} we briefly discuss how this notion reduces to Def. \ref{primary_field} and we re-interpret Theorem
\ref{tensor_product} in the setting of vertex operator algebras.

In Subsect. \ref{35.3} we complete the identification taking into account the associativity morphisms.

 \end{rem}

 \begin{rem} Note that  at this point of the paper, the notation $V$ seems confusing, in that it both denotes the fundamental representation of a complex Lie algebra ${\mathfrak g}$ as in \cite{Wenzl} and a general vertex operator algebra
as in Sect. \ref{VOAnets}.
When the two notions are used simultaneously, we shall refrain from using  $V$ for a vertex operator algebra. We hope that consistent use of the notation  
$V_{{\mathfrak g}_k}$ for the affine vertex operator algebra at level $k$ associated to ${\mathfrak g}$, clarifies the possible confusion.
\end{rem}

 \bigskip

\section{Knizhnik-Zamolodchikov differential equations, and the braiding}\label{32.5}
      The conformal Hamiltonian determines the braiding of ${\rm Rep}(V_{{\mathfrak g}_k})$ explicitly.
   In the setting of loop groups, formulae for the braiding may be found in 
   \cite{Wassermann} in the type $A$ case and see also Chapter 9   Lemma 6.1 in \cite{Toledano_laredo}.
   The fusion tensor product is defined via {\it Connes fusion}, a tensor product operation of bimodules
   over von Neumann algebras \cite{Connes}.
   The braiding gives rise to certain commutation properties satisfied by four-point functions of primary fields,
   called {\it braiding properties} in \cite{Wassermann}, \cite{Toledano_laredo} that reflect the action of the braiding operators in the category of projective finite energy representations with Connes fusion tensor product.
   We are interested to
    unravel the action of the braiding operators on the corresponding lowest energy subspaces.  In this way we shall be able
   to compare the braiding of ${\rm Rep}(V_{{\mathfrak g}_k})$ arising from the setting of loop groups
   directly with the braiding that we obtain from our construction in Theorem \ref{Zhu_as_a_compatible_unitary_wqh} (a)
    that   in turn has been unravelled in Prop. \ref{braiding_Zhu_algebra}. 
    
    We consider $\phi_4, \dots, \phi_1$  charged primary field
    of type      $$ {{m_4}\choose{0i_3}} {{m_3}\choose{i_3i_2}}{{m_2}\choose{i_2i_1}}{{m_1}\choose{i_10}}$$ respectively. 
    The four point functions of these   fields  is defined as the   formal Laurent series 
    $F=(\Gamma, \phi_4(v_4, t_4)\dots\phi_1(v_1, t_1)\Gamma)$, with
    $\Gamma\in{\mathcal H}_0[0]$ of norm $1$.    As for the case of a single primary field, it follows from the intertwining property  that $F$  takes values in ${\rm Hom}_G(V_{m_1}\otimes V_{m_2}\otimes V_{m_3}\otimes V_{m_4}, {\mathbb C})=
    (V_{m_4}^*\otimes V_{m_3}^*\otimes V_{m_2}^*\otimes V_{m_1}^*)^G
    $. Moreover $F$ is a formal solution of
      a first order partial differential equation, called the  Knizhnik-Zamolodchikov equation, 
    $$\frac{\partial F}{\partial t_i}=\frac{1}{k+h^\vee}\sum_{j\neq i}\frac{\Omega_{i,j}}{t_i-t_j}F$$
    with the convention that the denominators are expanded in series in the region $|t_4|>\dots>|t_1|$
    and $$\Omega_{i,j}=\sum_r\pi_i(X_r)\pi_j(X^r),$$ with $X_r\in{\mathfrak g}$ a basis and $X^r$ a dual basis, $\pi_i$ the action of ${\mathfrak g}$ on the $i$-th factor of $V_{m_4}^*\otimes V_{m_3}^*\otimes V_{m_2}^*\otimes V_{m_1}^*$. In particular, $\Omega_{i, i}$ acts as the Casimir $C_i$.
    
    The theory of  KZ equations implies, among other things, that   products of primary fields $\phi_{m_4, j}^{m_3}(v_3, t)\phi_{j, m_1}^{m_2}(v_2, s)$ are   single valued weakly holomorphic functions
    in $\{(t, s): t/s\notin[0, +\infty)\}$ on the finite energy vectors.
    
    A braiding is determined by an isomorphism between the subspace of ${\rm Hom}_G(V_2\otimes V_3, V_4)\simeq 
    (V_4\otimes V_3^*\otimes V_2^*)^G$ of the  initial terms of  and that of ${\rm Hom}_G(V_3\otimes V_2, V_4)\simeq (V_4\otimes V_2^*\otimes V_3^*)^G$. We want to know the action of $\Omega_{2,3}$.
    It follows from an algebraic manipulation that relates $\Omega_{2,3}$ to the Casimirs $\Omega_{i,i}=C_i$,
    as in Subsect. \ref{29.2}, $i=0, \dots, 3$,
    and the fact that $\pi_1+\dots+\pi_3=0$ that $\Omega_{2,3}$ acts on this space as multiplication by $\frac{1}{2}(C_4-C_2-C_3)$.
     It is explained in
    the proof of Lemma 3.1, Chaper 7, \cite{Toledano_laredo}, how this action explicitly determines the braiding previously
    defined using the conformal hamiltonian $L_0$ in the special case where the decomposition of
    ${\mathcal H}_1\boxtimes{\mathcal H}_2$ is multiplicity-free. 
    This is called {\it abelian braiding} in \cite{Wassermann}, \cite{Toledano_laredo}. We thus have the following corollary of this lemma
    that is implicit in \cite{Toledano_laredo}.
    
   \begin{cor}\label{Toledano_laredo}   Let ${\mathcal H}_2$ and ${\mathcal H}_3$ be 
   finite energy unitary projective representations of $LG$ such that Connes fusion tensor product is multiplicity-free.
   Then braiding operator $c({\mathcal H}_2, {\mathcal H}_3)$   acts between the corresponding 
   lowest energy spaces $V_2$ and $V_3$ respectively as  
   \begin{equation}\label{conformal_braiding}
   \Sigma e^{\frac{i\pi}{2(k+h^\vee)}\Omega}: V_2\otimes V_3\to V_3\otimes V_2
   \end{equation} where $\Omega$ acts as multiplication by $C_2+C_3-C_4$ on a simple module $V_4$ of the tensor product.
   \end{cor}

    \begin{proof}
  The proof follows   from a decomposition of $V_2\otimes V_3$
    into simple representations $V_4$ and the previous discussion.
    \end{proof}

  Note that formula (\ref{conformal_braiding}) is in \cite{Kirillov3},   without details on  the eigenvalues of $\Omega$.
   On  one hand, regardless the multiplicity-free property,
   a braiding is completely determined by  the operators $c(V, V_\lambda)$, with $V$ the   generating representation, by
      Prop.  \ref{braided_symmetry_with_generating_object},
      provided the associator is determined.

      On the other hand, we note that  in Prop. \ref{braiding_Zhu_algebra} we 
         have derived the same formula for the braiding as  (\ref{conformal_braiding}) for the special braidings $c(V, V_\lambda)$, and we recall that  the proof of this fact   relies on the analysis
        of Lemma 3.6.2 in \cite{Wenzl}.\bigskip
        
        \section{Interlude and organization of the next steps}\label{Interlude}
      
In Sect.   \ref{5+}  we outlined   our strategy for the proof Theorem \ref{Finkelberg_HL}. 
Recall Def. \ref{CFT_type_associator}  of ${\mathcal V}$-preassociator of CFT-type. In Sect. \ref{73} we  started to develop and apply the strategy to quantum group fusion categories. Our aim is to complete the proof of    Theorem \ref{Zhu_as_a_compatible_unitary_wqh}.

 So far regarding this theorem, we have worked on the side of quantum groups. We have constructed a unitary coboundary weak quasi-bialgebra structure on
 $A(V_{{\mathfrak g}_k})$, which induces a unitary ribbon modular tensor category structure on ${\rm Rep}(A(V_{{\mathfrak g}_k}))$. We have proved parts  (a) and (b)  of    Theorem \ref{Zhu_as_a_compatible_unitary_wqh}. In Theorem \ref{Zhu_from_qg_if_of_CFT_type} we have shown that the associator of $A(V_{{\mathfrak g}_k})$ induced from the quantum group, restricts to a 
${\mathcal V}$-preassociator of CFT-type,  building on the work of Wenzl  on the unitary structure of the fusion category.
 We need to prove  part (c) of the same theorem.

   On the side of vertex operator algebras, Huang and Lepowsky have introduced and constructed a vertex modular tensor category structure in module categories
  of a general class  of vertex operator algebras
  \cite{Huang_LepowskiI},
        \cite{Huang_LepowskiII}, \cite{Huang_LepowskiIII}, \cite{HuangIV}. They have shown that a vertex modular tensor category structure induces in a natural way
  a modular tensor category structure on the same category \cite{HL_tensor_products_of_modules}.
  
   In Sect. \ref{32} we sketch a description of Huang-Lepowsky vertex tensor category structure and the associated tensor category structure. 
          Let us   consider their ribbon tensor category structure in the case of
 ${\rm Rep}(V_{{\mathfrak g}_k})$, with $k$ a positive integer. This special case has been studied in \cite{Huang_Lepowski_affine}.

Part (c) of Theorem \ref{Zhu_as_a_compatible_unitary_wqh} is reformulated in  Theorem \ref{Finkelberg_HL} (c), and
diagram (\ref{31.2})   describes the methods that we intend to use to study  this part.
The proof of   tensor preserving equivalence of  the first map of the diagram  (\ref{31.2})  will be considered in
Subsects.    \ref{34.3}, \ref{34.4}.

We refer the reader to  \cite{On_a_problem_posed_by_Huang}
  for further
  explanations on this program,  and for  the developments of the abstract aspects of the original remark.

  We were not aware at the time of working on our paper,
     of the beautiful paper   \cite{McRae}\footnote{CP is grateful to  M. Yamashita for a discussion leading   to \cite{McRae}}   where   transport methods similar to these of our Sect. \ref{12} are developed. Note that McRae  
   does not   use of weak quasi-Hopf algebras, but we strongly advise this paper. For the present paper we follow     methods as discussed in the last years for this project, and we have kept the original arguments in the presentation of Sects. \ref{32}, \ref{33}.
            
The arguments are as follows. The mentioned proof    will break in two steps.
The {\it first step} transports Huang-Lepowsky vertex ribbon tensor category structure of
 ${\rm Rep}(V)$  to  the representation category of the Zhu algebra
 ${\rm Rep}(A(V))$ via Zhu's functor $Z$, and this is the content of Subsect.  \ref{34.3}, except for the braiding that is considered separately in Sect. \ref{22} and relies on Toledano-Laredo exposition in the application $V=V_{{\mathfrak g}_k}$.

 The {\it second step} compares this  transported vertex tensor category structure   of the Zhu algebra
 for $V=V_{{\mathfrak g}_k}$ with the structure arising from side of quantum groups summarized in the first two paragraphs of this section.           
 
The comparison with the transported  associativity morphisms of Huang and Lepowsky relies on the notion of primary field introduced by Knizhnik and Zamolodchikov \cite{KZ}, the work Tsuchiya and Kanie and Fenkel and Zhu
in Subsect. \ref{34.3}, which reduces infinite dimensional representations to finite dimensional representations of the Zhu algebra  in a natural (vertex) tensor categorical way. This is a specific property of the model and the identification is given by certain linear maps, part of a natural transformation, that uses 
the clarity of the work by Wenzl on the unitary structure on  triples of representations of which two terms are the fundamental representation and the remaining term is an arbitrary irreducible, for all Lie types on the side of quantum groups, which we show to be
a significative class of triples. On the side of conformal field theory, the identifying linear maps
uses the lack of logarithmic terms of certain solutions of a one-variable KZ equation, considered in \cite{Tsuchiya_Kanie}
and 
\cite{Huang_differential_equations}. This  indicates that an approach via unitarity is parallel to that vis the KZ equation.

The comparison between the Huang-Lepowsky associator on the Zhu algebra and the CFT-type pre-associator on the same algebra  from Theorem \ref{Zhu_as_a_compatible_unitary_wqh} is discussed in Subsect. \ref{34.4}.   This part relies on the fact that the two tensor product bifunctors may be identified. At the level of affine Lie algebras, an explanation of this identification was given in Theorem \ref{tensor_product}, and already at this level
   in turn relies
       on the central notion of primary field first introduced by Knizhnik and Zamolodchikov and studied by Tsuchiya and Kanie \cite{KZ}, \cite{Tsuchiya_Kanie}. In the setting of vertex operator algebras, identification of the two tensor product bifunctors on the Zhu algebra again based on the corresponding notion of primary field  will be discussed in Sect. \ref{33}.

    More precisely, to discuss tensor equivalence of the {\it second step} we first need to compare the two pre-tensor structures.     Theorem \ref{tensor_product}   
   relates the spaces of initial terms of primary fields  in the setting of loop groups with the tensor product of
${\rm Rep}(A(V_{{\mathfrak g}_k}))$ constructed in Sect. \ref{21}.  To make  this result useful, we need to relate primary fields in the setting of loop groups with primary fields in the setting of vertex operator algebras. This will be done in Sect. \ref{33}. Then the equivalence of the two associativity morphisms follows from the fact that Huang-Lepowsky associativity morphisms on the Zhu algebra is of CFT-type.

In Subsect. \ref{35.1} we discuss the notion of primary fields in the setting of vertex operator algebras. In Subsect.
        \ref{35.2} we discuss the correspondence with
     the setting of Sect. \ref{22} and re-interpret Theorem \ref{tensor_product} in this setting. This gives   a correspondence of the
     two pre-tensor category structures on the Zhu algebra.

Then
we will   identify the two tensor product bifunctors on ${\rm Rep}(A(V_{{\mathfrak g}_k}))$ arising from Huang-Lepowsky tensor product theory and from Theorem \ref{Zhu_as_a_compatible_unitary_wqh}.

            Since the second and last maps of the diagram (\ref{31.2})  are a ribbon tensor equivalences by the indicated Theorems,  
            completion of the proof of  
       Theorems \ref{Zhu_as_a_compatible_unitary_wqh} (c) and also of the main Theorem \ref{Finkelberg_HL}
     reduces to  prove the identification of  the associativity morphisms and braiding morphisms of the first map
     in  (\ref{31.2})  by comparing the braided tensor category structure  on  ${\rm Rep}(A(V_{{\mathfrak g}_k}))$ obtained in  parts (a) and (b) of the same Theorem  \ref{Zhu_as_a_compatible_unitary_wqh}
  with those of the braided tensor category structure on ${\rm Rep}(A(V_{{\mathfrak g}_k}))$ induced by Huang-Lepowsky braided tensor category structure of  ${\rm Rep}(V_{{\mathfrak g}_k})$.

Regarding the comparison of the two associativity morphisms and the two braiding morphisms on ${\rm Rep}(A(V_{{\mathfrak g}_k}))$, we recall that on the side of quantum groups,  Theorem \ref{Zhu_from_qg_if_of_CFT_type}  gives the explicit evaluation  of 
       the associator of  $A(V_{{\mathfrak g}_k})$ derived from Theorem \ref{Zhu_as_a_compatible_unitary_wqh} (a) on the special  triples of representations with two terms  occupied by the fundamental representation of the classical Lie algebra and the remaining term occupied by an arbitrary irreducible representation in the Weyl alcove. The theorem describes this associator as an associator extending the  ${\mathcal V}$-pre-associator of 
       CFT-type associated to $(F_0, G_0)$ and to the collection ${\mathcal V}$ of these special triples. For the braided symmetries,
       theorems \ref{braiding_Zhu_algebra} and Corollary \ref{Toledano_laredo} identify the braiding morphisms for all pairs with one term occupied by the fundamental representation. We  show in Subsect. \ref{34.4} that the associator on ${\rm Rep}(A(V_{{\mathfrak g}_k}))$ arising from Huang-Lepowsky tensor category structure
        restricts to the  CFT-type pre-associator associated  to $(F_0, G_0)$ on the same special triples ${\mathcal V}$.
       The conclusion will follow as an application of the uniqueness Theorem       
\ref{claim1}, provided the assumption on the generating properties of the braid group are satisfied.
These assumptions have been verified in Sect. 11 of \cite{On_a_problem_posed_by_Huang}, based on the existing literature, and leads to the
restricted Lie types as stated.

        We   note that an interesting problem is    to develop more explicitly  the connection between our unitary
        structure on the Zhu algebra,   the CPT operator and also with the papers by Gui, a statement of which may be found in \ref{Zhu_as_a_compatible_unitary_wqh} (d).
        The Hermitian form of the unitary wqh $A(V_{{\mathfrak g}_k})$ conjecturally
coincides with that of Gui in the setting of vertex operator algebras \cite{GuiI, GuiII}. We postpone the development of Theorem
\ref{Zhu_as_a_compatible_unitary_wqh} (d).

  \bigskip

 \section{Huang and Lepowsky  (vertex) tensor category structure
 }\label{32}

       \subsection{Intertwining operators for vertex operator algebras and $P(z)$-intertwining maps}   \label{34.0}
          
           In Sect. \ref{VOAnets2} we   briefly mentioned   the main results by Huang and Lepowsky on the construction
of the modular vertex tensor category structure on ${\rm Rep}(V)$. 
Every vertex tensor category has a naturally associated tensor category as described in Sect. 4 of \cite{HL_tensor_products_of_modules}. In this section we   describe this tensor category structure in   more detail.
          
          For a complex vector space $U$, let $U\{x\}$ denote the vector space of formal  power series 
             $\sum_{n\in{\mathbb R}}u_n x^n$ with generalized
          series over ${\mathbb R}$,
             and coefficients in $u_n\in U$.
       We recall the definition of intertwining operator from \cite{FHL}.
          
   \begin{defn}
   Let $(V, Y, 1, \nu)$ be a vertex operator algebra, see Sect. \ref{VOAnets},  and let $W_1$, $W_2$, $W_3$ be $V$-modules, Def. \ref{V_module}. An {\it intertwining operator} of type
             ${{W_3}\choose{W_1\ W_2}}$ is a linear map
             $$ {\mathcal Y}: W_1\to{\rm Hom}(W_2, W_3)\{x\}, \quad\quad    {\mathcal Y}(w^{(1)}, x)=
             \sum_{r\in{\mathbb R}}  w^{(1)}_{(r)} x^{-r-1} $$
             satisfying the following axioms (note that we use a suffix e.g. $w^{(1)}$ to indicate the space where the vector lies and an index $w^{(1)}_{(r)}$ to indicate a coefficient of the intertwining operator ${\mathcal Y}$)
            \medskip

  \noindent {\bf a)}  {\it (lower truncation condition)}:   $w^{(1)}_{(r)}w^{(2)}=0$ for $r$ sufficiently large,\smallskip
  
   \noindent {\bf b)}  {\it $L_{-1}$-derivative (or translation) property}  
        \begin{equation}\label{translation_intertwining}\frac{d}{dx}{\mathcal Y(w^{(1)}, x)}={\mathcal Y}(L_{-1}^{W_1}w^{(1)}, x),\end{equation}
        \smallskip
  
            \noindent{\bf c)} {\it (Jacobi identity)} for $v\in V$, $w^{(1)}\in W_1$,
              \begin{equation}\label{Jacobi_intertwining_operator_delta}
x_0^{-1}\delta(\frac{x_1-x_2}{x_0})Y_{W_3}(v, x_1){\mathcal Y}(w^{(1)}, x_2)-x_0^{-1}\delta(\frac{x_2-x_1}{-x_0}){\mathcal  Y}(w^{(1)}, x_2)Y_{W_2}(v, x_1)=
\end{equation}
 $$x_2^{-1}\delta(\frac{x_1-x_0}{x_2}){\mathcal Y}(Y_{W_1}(v, x_0)w^{(1)}, x_2)).$$

\noindent The module $W_1$ is called the {\it charge}, $W_2$ the {\it source} and $W_3$ the {\it target}.
     
     \end{defn}
     
     \begin{defn} \label{fusion_Rule_for_VOA_intertwiner}   The dimension of the vector space ${{\mathcal M}_{W_1, W_2}}^{W_3}$ of intertwining operators of a given type
             ${{W_3}\choose{W_1\ W_2}}$ is called the {\it fusion rule}.
    \end{defn}

For a $V$-module $W$, let $\overline{W}$ denote the algebraic completion of $W$,
$$\overline{W}=\Pi_{r\in{\mathbb R}}W_{r},$$
\begin{equation}\label{P_r}
P^W_r: \overline{W}\to W_{r}, 
\end{equation}
the natural projection maps  for $r\in{\mathbb R}$.
We recall the definition of $P(z)$-intertwining map from \cite{HL_tensor_products_of_modules}, \cite{Huang_LepowskiI},
 \cite{Huang_LepowskiII},  \cite{Huang_LepowskiIII}. 

\begin{defn}\label{P(z)-intertwining-map}
Fix $z\in{\mathbb C}^{\times}$, and let $W_1$, $W_2$, $W_3$ be $V$-modules.
A $P(z)$-{\it intertwining map} of type ${{W_3}\choose{W_1\ W_2}}$
is a linear map
\begin{equation}
F: W_1\otimes W_2\to \overline{W}_3
\end{equation} satisfying for $w^{(1)}\in W_1$, $w^{(2)}\in W_2$, $v\in V$,
 \medskip
 
 \noindent {\bf a)} ({\it lower truncation})  $P^{W_3}_{r-n}F(w^{(1)}\otimes w^{(2)})=0$ for $n\in{\mathbb Z}$ large,
 \medskip

 \noindent {\bf b)} ({\it intertwining relation}) \begin{equation}
x_0^{-1}\delta(\frac{x_1-z}{x_0})Y_{W_3}(v, x_1)F(w^{(1)}\otimes w^{(2)})=
\end{equation}
$$
z^{-1}\delta(\frac{x_1-x_0}{z}) F(Y_{W_1}(v, x_0)w^{(1)}\otimes w^{(2)})+
x_0^{-1}\delta(\frac{z-x_1}{-x_0})F(w^{(1)}\otimes Y_{W_2}(v, x_1)w^{(2)}).
 $$
\end{defn}

Intertwining operators and $P(z)$ intertwining maps are closely related.
Roughly speaking, a $P(z)$-intertwining map results from an intertwining operator by evaluation
  on a nonzero complex number $z$.
Along this procedure, the charge and source spaces do not change, but the range  is replaced by $\overline{W}_3$.  
The following is a more precise description. 

Let $\log(z)$ be the principal branch of the complex logarithm function with the complex plane cut along the positive real axis,
that is
$$\log(z)=\log(|z|)+ i\arg(z),\quad\quad 0\leq \arg(z)<2\pi,$$
and for $p\in{\mathbb Z}$, define the branches 
$$l_p(z)=\log(z)+2\pi p i.$$
\begin{prop}\label{intertwining_map_intertwining_operator0}
(Prop. 12.2 in \cite{Huang_LepowskiIII}) Given $z\in{\mathbb C}$,  a fixed branch $l_p$ induces   a canonical isomorphism
from the space of intertwining operator of type ${{W_3}\choose{W_1\ W_2}}$
onto the space of $P(z)$-intertwining maps of the same type.

 \begin{itemize} 
\item[{\rm  a)}]

The isomorphism from intertwining operators to $P(z)$-intertwining maps  is given by specification of the formal
variable $x$ associated to $l_p$,
\begin{equation} \label{intertwining_map_intertwining_operator} {\mathcal Y}\to F, \quad\quad F(w^{(1)}\otimes w^{(2)})={\mathcal Y}(w^{(1)}, e^{l_p(z)})w^{(2)},\end{equation}

\item[{\rm  b)}]
The inverse isomorphism  
is given by
\begin{equation}F\to{\mathcal Y}, \quad\quad  {\mathcal Y}(w^{(1)}, x)w^{(2)}=\sum_{r\in{\mathbb R}} e^{(r+1)l_p(z)}P_{\Delta_{w^{(1)}}+ \Delta_{w^{(2)}}-r-1}F(w^{(1)}\otimes w^{(2)})x^{-r-1},\end{equation}
with $w^{(i)}$   homogeneous with conformal weight $\Delta_{{w^{(i)}}}$, in the sense of (\ref{conformal_weight_of_vector}).  
\end{itemize}
In particular, the dimension of the vector space of $P(z)$-intertwining maps of type ${{W_3}\choose{W_1\ W_2}}$
equals the fusion rule associated to the same $V$-modules.

\end{prop}

Note that the series at the right hand side of (\ref{intertwining_map_intertwining_operator})converges in the product topology of $\overline{W}_3$ of the discrete
topological spaces $(\overline{W}_3)_{(m)}$ by (\ref{grading_for_intertwining_irreducible}).

Prop. \ref{intertwining_map_intertwining_operator0} plays a central role in the connection between
affine vertex operator algebras and weak  Hopf algebras.  

Similarly to Def. \ref{new_notation_for_gradinf_for_modules},  from now on we shall use the following convenient change of notation for modes of intertwining operators.

\begin{defn} {\it (New grading  notation for modes intertwining operators)} 
We set 
\begin{equation}\Delta_{\mathcal Y}=\Delta_{W_1}+\Delta_{W_2}-\Delta_{W_3},\end{equation}
with $\Delta_{W_i}$ is the conformal weight of $W_i$, as in
(\ref{conformal_weight_of_module}). 
We also set
  \begin{equation} {w^{(1)}}_n={w^{(1)}}_{(n-1+\Delta_{\mathcal Y})}\end{equation}
  Similarly to the case of modules, the translation and Jacobi properties for intertwining operators imply  for 
  $w^{{(1)}}$ homogeneous,    \begin{equation} 
  \label{grading_for_intertwining_irreducible}{w^{(1)}}_n: (W_2)_{(m)}\to (W_3)_{(m+ {\rm deg}(w^{{(1)}})-n)},
  \end{equation}
  where for $M$ irreducible and $n\in{\mathbb Z}_{\geq0}$, $M_{(n)}$ has been defined in (\ref{grades}) and ${\rm deg}(w^{{(1)}})=k\in{\mathbb Z}$
  means that $w^{{(1)}}\in (W_1)_{(k)}$
  (\cite{FHL}, see also Prop. 1.5.1 in \cite{Frenkel_Zhu} for details)    
\end{defn}

We assume for convenience that ${\rm Rep}(V)$ is semisimple as a linear category (e.g. $V$ is rational in the sense of
\cite{Huang_LepowskiI}).
Then, similarly to the case of the vertex operators $Y_M(a, x)$ cf.  Sect. \ref{VOAnets}, 
when $W_i$ are all simple, $i=1$, $2$, $3$,
the defining series of an intertwining operator ${\mathcal Y}$ of type ${{W_3}\choose{W_1\ W_2}}$ is over the countable set, 
\begin{equation}\label{intertwining_irreducible_case}
{\mathcal Y}(w^{(1)}, x)=\sum_{n\in{\mathbb Z}}{w^{(1)}}_n x^{-n} x^{-{\Delta}_{\mathcal Y}}.\end{equation}

\begin{defn}\label{HL_P(z)_tensor_product}
Let $W_1$ and $W_2$ be $V$-modules. A $P(z)$-tensor product is a $V$-module
($W_1\boxtimes_{P(z)}W_2, Y_{W_1\boxtimes_{P(z)}W_2}$) with a $P(z)$-intertwining map 
\begin{equation}
F: W_1\otimes W_2\to \overline{W_1\boxtimes_{P(z)}W_2}
\end{equation}
of type ${{W_1\boxtimes_{P(z)}W_2}\choose{W_1\ W_2}}$ such that for any other
$V$-module ($W_3, Y_{W_3})$ and a $P(z)$-intertwining map
\begin{equation}
F': W_1\otimes W_2\to \overline{W_3}
\end{equation}
of type ${{W_3}\choose{W_1\ W_2}}$ there is a unique morphism of $V$-modules 
$$\eta: W_1\boxtimes_{P(z)}W_2\to W_3$$
defined in Def. 
\ref{Rep(V)} 
such that the following diagram commutes
 
\[
 \begin{tikzcd}
W_1\otimes W_2 \arrow{r}{F}\arrow[swap] {dr} {F'} & \overline{W_1\boxtimes_{P(z)}W_2} \arrow{d}{\overline{\eta}} \\
 & \overline{W}_3
 \end{tikzcd}
 \]
 where $\overline{\eta}$ is the natural extension of $\eta$ to the algebraic completions. 

\end{defn}
The above universality property implies that a $P(z)$-tensor product $V$-module is unique up to isomorphism.
\medskip

A $P(z)$-tensor product space is constructed as follows.
Let $V$ be rational in the sense of \cite{Huang_LepowskiI}, that is
1) $V$ admits only finitely many inequivalent  irreducible $V$-modules, 2)
 every V -module is completely reducible, 3) the fusion rules associated to triples of  irreducible modules are finite.
  Then by Prop. 12.5  \cite{Huang_LepowskiI} 
 \begin{equation}\label{P(z)_tensor_product} W_1\boxtimes_{P(z)} W_2=\bigoplus_i ({\mathcal M}[P(z)]_{W_1, W_2}^{M_i})^* \otimes M_i
 \end{equation}
 is a $P(z)$-tensor product module, with $({\mathcal M}[P(z)]_{W_1, W_2}^{M_i})^*$
 the dual space of the space of $P(z)$-interwining maps of type
 ${{M_i}\choose{W_1\ W_2}}$
 and $\{M_i\}$ is a complete family of irreducible $V$-modules.
 The vertex operator $Y_{W_1\boxtimes_{P(z)} W_2}$ is the direct sum of the corresponding vertex operators.
 \medskip

  \begin{defn} ({\it Products and iterates of intertwining operators})
Intertwining operators ${\mathcal Y}_1$ and ${\mathcal Y}_2$ associated to different triples  of $V$-modules
of types
 ${{W_4}\choose{W_1\ W_5}}$ and ${{W_5}\choose{W_2\ W_3}}$ respectively,
 may be composed in a way similar to composition of elements of a category,
\begin{equation}\label{product}
{\mathcal Y}_1(w^{(1)}, x_1){\mathcal Y}_2(w^{(2)}, x_2),
\end{equation}
(i.e. with the only requirement that the target of ${\mathcal Y}_2$ equals the source of ${\mathcal Y}_1$.)
This composition is called the {\it product} of intertwining operators.

Similarly, intertwining operators ${\mathcal Y}_3$ and ${\mathcal Y}_4$ of types  
 ${{W_6}\choose{W_1\ W_2}}$ and ${{W_4}\choose{W_6\ W_3}}$ may be composed along the charge space
 \begin{equation}\label{iterate}
 {\mathcal Y}_4({\mathcal Y}_3(w^{(1)}, x_1-x_2)w^{(2)}, x_2).
 \end{equation}
 This composition is called the {\it iterate}.
 \end{defn}

The product and iterate of intertwining operators give   rise to  the triple tensor products $W_1\boxtimes_{\rm HL}(W_2\boxtimes_{\rm HL} W_3)$ and $(W_1\boxtimes_{\rm HL} W_2)\boxtimes_{\rm HL} W_3$, respectively, in the following way.

A central property for intertwining
operators, called {\it associativity} of intertwining operators, a property relating products with iterates.
Associativity is very important, and draws its origin from the {\it Operator Product Expansion}, introduced by 
Belavin, Polyakov and Zamolodchikov in 2-dimensional conformal field theory \cite{BPZ}.

 \begin{rem}\label{on_OPE} Associativity of intertwining operators  is referred     the {\it (Non-meromorphic) Operator Product Expansion (OPE) of chiral vertex operators},   in the works by Huang and Lepowsky. Note that in physics literature one often meets OPE for {\it primary fields}.
 Since primary fields    are of main interest for our purposes,
to avoid possible confusion, we shall refer to Huang-Lepowsky associativity   as
   {\it  OPE of intertwining operators}. See also Def. \ref{OPE_intertwining_operators} for a reformulation in terms of irreducible modules.\end{rem}

 \begin{defn} \label{associativity_intertwining_operators_def} ({\it Theorem 14.8 in \cite{HuangIV}})  {\it Associativity of intertwining operators}  holds when
for any pair of intertwining operators ${\mathcal Y}_1$, ${\mathcal Y}_2$  
of types ${{W_4}\choose{W_1\ W_5}}$ and ${{W_5}\choose{W_2\ W_3}}$
  respectively
there are $W_6$ and intertwining operators ${\mathcal Y}_3$, ${\mathcal Y}_4$ of types  ${{W_6}\choose{W_1\ W_2}}$ and ${{W_4}\choose{W_6\ W_3}}$ 
 respectively such that
 for any two nonzero complex numbers $z_1$ and $z_2$ satisfying 
 \begin{equation}\label{position_of_complex_numbers}
 |z_1|>|z_2|>|z_1-z_2|>0,
 \end{equation} and vectors $w^{(i)}$ in the corresponding charge spaces 
 we have 
 \begin{equation}\label{associativity_intertwining_operators} {\mathcal Y}_1(w^{(1)}, z_1){\mathcal Y}_2(w^{(2)}, z_2)={\mathcal Y}_4({\mathcal Y}_3(w^{(1)}, z_1-z_2)w^{(2)}, z_2).
 \end{equation}
 Moreover, given ${\mathcal Y}_3$ and ${\mathcal Y}_4$  as above there are  $W_5$ and
 ${\mathcal Y}_1$ and ${\mathcal Y}_2$ as above such that (\ref{associativity_intertwining_operators}) holds in the same domain.
 \end{defn}

 \begin{rem}\label{convergence_and_extension} A term of the form of the left hand side of (\ref{associativity_intertwining_operators}) 
 is a series with the summation developed over the graded subspaces of the target of ${\mathcal Y}_2$ (or the source of ${\mathcal Y}_1$),
 and
 converges absolutely in the weak topology defined by the source of ${\mathcal Y}_2$ and the restricted dual of the
 target of ${\mathcal Y}_1$ if and only if
 the series at the right hand side
 of (\ref{associativity_intertwining_operators}) converges absolutely weakly as well.
 A property implying these equivalent conditions is developed in \cite{HuangIV} and called {\it convergence and extension property}. The property is also formulated in 
 Def. 3.4 in \cite{Huang_differential_equations}. In that paper, the convergence and extension property
  has been derived  
  under certain conditions. In particular  in the same paper the author
 showed that it holds  
 under the rationality conditions a), b) c) of Sect. \ref{VOAnets2}. 
 Moreover, it was shown in \cite{HuangIV} that associativity of intertwining operators follows from the
 convergence and extension property. Thus the convergence and extension property is a main analytical
 step  to construct the structure of a braided tensor
(in fact a braided vertex tensor) category on the category of modules of a vertex operator algebra.
The detailed discussion and derivation of this property and its connection to systems of differential 
equations may be found in \cite{Huang_differential_equations}.
We shall only be concerned with some simple aspects of convergence properties
that will be useful for us to expose the connection between Huang-Lepowsky associativity morphisms and
that arising from quantum groups (cf. Theorem \ref{Zhu_as_a_compatible_unitary_wqh})
 in the case of affine vertex operator algebras.
 
 \end{rem}

 \begin{defn} The  $(W_1\otimes W_2\otimes W_3\otimes (W_4)' )'$-valued functions defined by evaluating the left hand side of (\ref{associativity_intertwining_operators}) 
 on vectors of the source space and the restricted dual of the target space are called {\it $4$-point correlation 
 functions.}
 \end{defn}

\subsection{Huang-Lepowsky (vertex) tensor category associativity morphisms  of ${\rm Rep}(V)$  (OPE  of intertwining operators and OPE constants $F_{\alpha, \beta}^{\beta', \alpha'}$)}\label{34.2}
As already mentioned,  associativity of intertwining operators has the important consequence, together with some other properties, that
  the conditions of Huang and Lepowsky
tensor product theory  are verified, and therefore their theory may be applied \cite{Huang_LepowskiI}, \cite{Huang_LepowskiII}, \cite{Huang_LepowskiIII}, \cite{HL_tensor_products_of_modules}, \cite{Huang_Lepowski_affine}, \cite{HuangIV}. It follows that
  the $P(z)$-tensor product module (\ref{P(z)_tensor_product})
extends to   vertex tensor category structure by Theorem 3.7 in \cite{Huang_differential_equations}, 
which is in fact modular by \cite{Huang(modularity)}, \cite{Huang2}. To describe this tensor structure, we
 follow the exposition given by Gui,
 Sect. 2.4 in \cite{GuiI}.

Let $W_1$, $W_2$, $W_3$ be $V$-modules and let ${\mathcal M}_{W_1, M_r}^{M_s}$ and 
${\mathcal M}_{W_2, W_3}^{M_r}$ be the vector spaces of
 intertwining operators of  type
 ${{M_s}\choose{W_1\ M_r}}$ and  ${{M_r}\choose{W_2\ W_3}}$ respectively,
 with $\{M_r\}$   a complete family of irreducible $V$-modules.
Let \begin{equation}\label{basis}{\mathcal Y}_\alpha\in {\mathcal M}_{W_1, M_r}^{M_s} \quad\quad {\mathcal Y}_\beta\in {\mathcal M}_{W_2, W_3}^{M_r}\end{equation} be bases for these spaces, respectively, with $r$ varying.

When $|z_1|>|z_2|$ are nonzero fixed complex numbers, the collection of products
  $\{{\mathcal Y}_\alpha\upharpoonright_{z_1} {\mathcal Y}_\beta\upharpoonright_{z_2}\}$
is linearly independent. A proof may be found in Prop. 2.3 in \cite{GuiI}.
When the modules $W_i$ are irreducible and the $\arg$-function is chosen continuously
then these intertwining operators are single valued   functions.
If in addition $$|z_1|>|z_2|>|z_1-z_2|>0$$ then there is a basis development of products
of intertwining operators with respect to iterates of intertwining operators independent of the choice of the points,
corresponding as said to the notion of OPE in physics.

\begin{defn}  \label{OPE_intertwining_operators} Given ${\mathcal Y}_\alpha\in {\mathcal M}_{W_1, M_r}^{M_s} $ and ${\mathcal Y}_\beta\in {\mathcal M}_{W_2, W_3}^{M_r}$ bases as in
(\ref{basis}) with $r$ varying,
 let
$${\mathcal Y}_{\alpha'}\in {\mathcal M}^{{M_i}}_{{W_1\ W_2}}, \quad\quad
{\mathcal Y}_{\beta'}\in {\mathcal M}^{{M_s}}_{{M_i\ W_3}}$$ be bases of the vector spaces
 ${\mathcal M}^{{M_i}}_{{W_1\ W_2}}$ and ${\mathcal M}^{{M_s}}_{{M_i\ W_3}}$ respectively, with $i$ varying.
The basis basis development 
\begin{equation}\label{F_coefficients}
{\mathcal Y}_\alpha(w^{(1)}, z_1){\mathcal Y}_\beta(w^{(2)}, z_2)=\sum_{\alpha',\beta'} F_{\alpha, \beta}^{\beta', \alpha'} {\mathcal Y}_{\beta'}({\mathcal Y}_{\alpha'}(w^{(1)}, z_1-z_2)w^{(2)}, z_2),
\end{equation}
 is called {\it  OPE of intertwining operators}. We shall refer to $F_{\alpha, \beta}^{\beta', \alpha'}$ as the {\it (OPE) structure constants.}
 \end{defn}

In Huang-Lepowsky theory, the vertex tensor category gives rise to a tensor category. We shall omit the
description of this passage in their work. It will be important for us to notice that the structure constants
$F_{\alpha, \beta}^{\beta', \alpha'}$ will also be the structure constants for the associativity morphisms of the
associated tensor category structure. We summarize this passage in the following result, following the exposition in \cite{GuiII}.

 The OPE of intertwining operators
 determines the tensor category structure of ${\rm Rep}(V)$ as follows.
 We define the linear map
  \begin{equation}\label{coefficients2}B=B(W_1, W_2, W_3; s)
 \end{equation} to emphasize the dependance of $B$ on the variables.
 $$B: \bigoplus_{r} {\mathcal M}_{W_1, M_r}^{M_s} \otimes {\mathcal M}_{W_2, W_3}^{M_r}\to
 \bigoplus_{i}  {\mathcal M}_{{M_i\ W_3}}^{M_s}\otimes {\mathcal M}_{{W_1\ W_2}}^{M_i}  
 $$
 $${\mathcal Y}_\alpha\otimes{\mathcal Y}_\beta\to \sum_{\alpha', \beta'} F_{\alpha, \beta}^{\beta', \alpha'}{\mathcal Y}_{\beta'}\otimes{\mathcal Y}_{\alpha'}.$$
 Given finite dimensional vector spaces $X$, $Y$, we  identify the dual space of a direct sum and a tensor product with
 $$(X\oplus Y)^*\simeq X^*\oplus Y^*, \quad\quad (X\otimes Y)^*\simeq Y^*\otimes X^*.$$
 Finally we set
 \begin{equation} \label{coefficients} A=B^t: 
 \bigoplus_{i} ({\mathcal M}_{{W_1\ W_2}}^{M_i})^*\otimes  ({\mathcal M}_{{M_i\ W_3}}^{M_s})^*\to
  \bigoplus_{r}  ({\mathcal M}_{W_2, W_3}^{M_r})^*\otimes ({\mathcal M}_{W_1, M_r}^{M_s})^* ,\end{equation}
 with $B^t$ the transpose of $B$. We shall write
 \begin{equation}\label{coefficients3}A=A(W_1, W_2, W_3; s).
 \end{equation}
 The following is  a major achievement of Huang-Lepowsky theory.  \bigskip
 
\begin{thm} \label{HL} (Huang-Lepowsky) Let $V$ be a vertex operator algebra satisfying the
the rationality conditions a), b), c) of Sect. \ref{VOAnets2},  
 let $W_1$, $W_2$, $W_3$ be $V$-modules, and
let ${\mathcal M}_{W_1, W_2}^{M_i}$ be the vector space 
 interwining operators of type
 ${{M_i}\choose{W_1\ W_2}}$
 with $\{M_i\}$   a complete family of irreducible $V$-modules.
 Then
 \begin{itemize} 
\item[{\rm  1)}]
 \begin{equation}\label{tensor_product_module}
 W_1\boxtimes_{\rm HL} W_2:=\bigoplus_i ({\mathcal M}_{W_1, W_2}^{M_i})^*\otimes  M_i,
\end{equation}
 the vertex operator $Y_{W_1\boxtimes_{\rm HL} W_2}$ is defined as the direct sum of the corresponding vertex operators.
Writing
\begin{equation}
(W_1\boxtimes_{\rm HL} W_2)\boxtimes_{\rm HL} W_3=\bigoplus_i ({\mathcal M}_{W_1, W_2}^{M_i})^* \otimes (M_i\boxtimes_{\rm HL} W_3)= 
\end{equation}
$$\bigoplus_{i , s} ({\mathcal M}_{W_1, W_2}^{M_i})^* \otimes ({\mathcal M}_{M_i, W_3}^{M_s})^*\otimes M_s
$$
and 
\begin{equation}
W_1\boxtimes_{\rm HL} (W_2\boxtimes_{\rm HL} W_3)=\bigoplus_r ({\mathcal M}_{W_2, W_3}^{M_r})^* \otimes (W_1\boxtimes_{\rm HL} M_r)= 
\end{equation}
$$\bigoplus_{r, s} ({\mathcal M}_{W_2, W_3}^{M_r})^* \otimes ({\mathcal M}_{W_1, M_r}^{M_s})^*\otimes M_s.
$$
then the associativity morphism is given by $$\alpha_{W_1, W_2, W_3}=\oplus_s A(W_1, W_2, W_3; s) \otimes 1_{M_s}.$$
\item[{\rm  2)}]
There are  natural identifications of the following spaces of intertwining operators and their duals  with corresponding morphism spaces in ${\rm Rep}(V)$:
$${\mathcal M}_{W_1, W_2}^{M_i}\simeq {\rm Hom}_V(W_1\boxtimes_{\rm HL} W_2, M_i), \quad\quad ({\mathcal M}_{W_1, W_2}^{M_i})^* 
 \simeq{\rm Hom}_V(M_i, W_1\boxtimes_{\rm HL} W_2).$$
\end{itemize}

\end{thm}

We  comment on 2). After the tensor category structure of ${\rm Rep}(V)$ has been constructed then
an element $\varphi\in ({\mathcal M}_{W_1, W_2}^{M_i})^*$ induces a morphism in
${\rm Hom}_V(M_i, W_1\boxtimes_{\rm HL} W_2)$ via left tensor product with $\varphi$; and an element  $T\in{\rm Hom}_V(M_i, W_1\boxtimes_{\rm HL} W_2)$
induces a linear functional on ${\rm Hom}_V(W_1\boxtimes_{\rm HL} W_2, M_i)$ via left composition by $T$ in ${\rm Rep}(V)$.

\subsection{Transporting Huang-Lepowsky vertex tensor category structure to the Zhu algebra}\label{34.3}
In this subsection we describe   the first step mentioned in Sect. \ref{Interlude}, useful  to
 construct   the first tensor equivalence in (\ref{31.2}). Recall that the first step aims to transport the vertex tensor category structure
 of ${\rm Rep}(V_{{\mathfrak g}_k})$ to the Zhu algebra, and preserves all the structure by construction. 
   
 To study  
 associativity morphisms of module categories of vertex operator algebras and the transported
 structure  to the Zhu algebra, we find it very clarifying the following passage from composition of intertwining operators to composition of Huang-Lepowsky  {\it intertwining maps}. 
  This passage is analogous to the case of   intertwining operators, see Prop. \ref{intertwining_map_intertwining_operator0}, and an important step to the construction
 of vertex tensor product theory.
 We next sketch Huang-Lepowsky associativity morphisms $\alpha_{\rm HL}$ of the vertex tensor category.
\begin{rem}
 In what follows, our notation $\alpha_{\rm HL}$ will correspond to
 $$\alpha_{\rm HL}:=({\mathcal A}_{P(z_1), P(z_2)}^{P(z_1-z_2), P(z_2)})^{-1}$$
 of Theorem 14.10 in \cite{HuangIV}.
 \end{rem}

 Given a branch of the complex logarithm and distinct nonzero complex numbers $z_1$, $z_2$
 satisfying  (\ref{position_of_complex_numbers}), if associativity of intertwining operators holds as in Def.
 \ref{associativity_intertwining_operators_def} then the associativity equation (\ref{associativity_intertwining_operators})
 can be written as
 \begin{equation}\label{composition_intertwining_maps}   F_1\circ 1\otimes F_2=F_4\circ F_3\otimes1,\end{equation}
 where $F_1$ and $F_2$ are $P(z_1)$- and $P(z_2)$- intertwining maps corresponding to ${\mathcal Y}_1$ and ${\mathcal Y}_2$ respectively, and $F_3$ and $F_4$ are $P(z_1-z_2)$- and $P(z_2)$- intertwining maps
corresponding to ${\mathcal Y}_3$ and ${\mathcal Y}_4$ respectively. The definition of both sides
of (\ref{associativity_intertwining_operators}) is as series 
obtained inserting in the middle projection maps onto the homogeneous subspaces, and is well defined
under convergence conditions. Note indeed  that the inner maps have targets in the algebraic completions of modules.
 In the construction of their tensor product theory, the study of products and iterates of intertwining maps
is the starting point. For an overview we refer the reader to \cite{HL}, and for complete explanations to
\cite{HuangIV}, and references to previous papers of the same series. In Subsect. \ref{35.3} we discuss in more detail significative examples of these series for the model $V=V_{{\mathfrak g}_k}$.

In terms of decomposition into irreducible modules, with the same notation as in Def.
\ref{OPE_intertwining_operators}, equation (\ref{F_coefficients}) becomes
\begin{equation}\label{associativity_in_terms_of_F}
F_\alpha\circ 1\otimes F_\beta=\sum_{\alpha', \beta'}F_{\alpha, \beta}^{\beta', \alpha'} F_{\beta'}\circ F_{\alpha'}\otimes 1,\end{equation}
where recall $F_\alpha$ and $F_\beta$ are linear bases $P(z_1)$- and $P(z_2)$- intertwining maps of type
${{M_s}\choose{W_1\ M_r}}$ and  ${{M_r}\choose{W_2\ W_3}}$  respectively
and  $F_{\alpha'}$, $F_{\beta'}$ bases of $P(z_1-z_2)$- and $P(z_2)$- intertwining maps of type 
${{M_i}\choose{W_1\ W_2}}$ and  ${{M_s}\choose{M_i\ W_3}}$   respectively.
Let 
\begin{equation}\label{F_1} F_{W_2, W_3}: W_2\otimes W_3\to \overline{W_2\boxtimes_{P(z_2)}W_3}
\end{equation} be the universal $P(z_2)$-intertwining map
of type ${{W_2\boxtimes_{P(z_2)}W_3}\choose{W_2\ W_3}}$
and \begin{equation}\label{eta_1}\eta_\beta: W_2\boxtimes_{P(z_2)}W_3\to M_r\end{equation} the unique morphism of $V$-modules
derived from $F_\beta$ by universality of the $P(z_2)$-tensor product,  see Def. \ref{HL_P(z)_tensor_product},
\begin{equation} \label{eta'_1} F_\beta=\overline{\eta_\beta}\circ F_{W_2, W_3}.\end{equation}
Similarly, we have universal maps and associated $V$-module morphisms for the various intertwining maps
composing (\ref{associativity_in_terms_of_F}) as follows;

\begin{equation}\label{additivity1} F_{W_1, M_r}: W_1\otimes M_r\to \overline{W_1\boxtimes_{P(z_1)} M_r},\end{equation}
the universal $P(z_1)$-intertwining  of type
${{W_1\boxtimes_{P(z_1)} M_r}\choose{W_1\ M_r}}$,
and the $V$-module morphism
\begin{equation}\label{eta_2}\eta_\alpha: W_1\boxtimes_{P(z_1)}M_r\to M_s\end{equation}
  such that
  \begin{equation} \label{eta'_2} F_\alpha=\overline{\eta_\alpha}\circ F_{W_1, M_r};\end{equation}
  the universal $P(z_1-z_2)$-intertwining 
  \begin{equation} \label{F_2} F_{W_1, W_2}: W_1\otimes W_2\to \overline{W_1\boxtimes_{P(z_1-z_2)} W_2}\end{equation}
of type
${W_1\boxtimes_{P(z_1-z_2)}W_2}\choose{{W_1\ W_2}}$
and the $V$-module morphism 
\begin{equation} \label{eta_3} \eta_{\alpha'}: W_1\boxtimes_{P(z_1-z_2)} W_2\to M_i\end{equation}
such that
\begin{equation}\label{eta'_3} F_{\alpha'}=\overline{\eta_{\alpha'}}\circ F_{W_1, W_2}\end{equation}
and finally the universal $P(z_2)$-intertwining

\begin{equation}\label{additivity2} F_{M_i, W_3}: M_i\otimes W_3\to\overline{M_i\boxtimes_{P(z_2)}W_3}\end{equation}
of type 
${M_i\boxtimes_{P(z_2)}W_3}\choose{{M_i\ W_3}}$,
the $V$-module morphism
\begin{equation}\label{eta_4} \eta_{\beta'}: M_i\boxtimes_{P(z_2)}W_3\to M_s\end{equation}
such that
\begin{equation}\label{eta'_4} F_{\beta'}=\overline{\eta_{\beta'}}\circ F_{M_i, W_3}.\end{equation}

The following formula gives Huang-Lepowsky associativity morphisms in terms of the OPE coeffiecients
$F_{\alpha, \beta}^{\beta', \alpha'}$.
\begin{thm} \label{HL_associator}  
Let $\eta_\alpha$, $\eta_\beta$, $\eta_{\alpha'}$, $\eta_{\beta'}$ be the $V$-module morphisms defined as in
(\ref{eta_1}), (\ref{eta'_1}), (\ref{eta_2}), (\ref{eta'_2}), (\ref{eta_3}), (\ref{eta'_3}), (\ref{eta_4}), (\ref{eta'_4}).

Let 
$\xi_\alpha$, $\xi_\beta$ be   $V$-module morphisms   with source (range) the range (source) of $\eta_\alpha$ and
$\eta_\beta$ respectively, and such that such that
\begin{equation}\label{Cuntz}\sum_\alpha \xi_\alpha\eta_\alpha=1, \quad\quad \sum_\beta \xi_\beta\eta_\beta=1.\end{equation} 
Then Huang-Lepowsky vertex tensor category associativity morphisms   are given by
$$\alpha_{\rm HL}:=\sum F_{\alpha, \beta}^{\beta', \alpha'}(1_{W_1}\boxtimes_{P(z_1)} \xi_\beta)\circ \xi_\alpha\circ \eta_{\beta'}\circ (\eta_{\alpha'}\boxtimes_{P(z_2)} 1_{W_3}).$$
\end{thm}
 
The morphisms $\xi_\alpha$, $\xi_\beta$ exist by semisimplicity.
The naturality of the universal intertwining maps gives the following form of $\alpha_{\rm HL}$  in terms of the universal
$P(z)$-intertwining maps on distinct points
(cf. Theorem 14.10 in \cite{HuangIV}).  
\begin{thm}\label{HL_associativity_vertex} The associativity equation
(\ref{associativity_in_terms_of_F}) in turn becomes
\begin{equation}\label{34.42} F_{W_1, W_2\boxtimes_{P(z_2)}W_3}\circ 1_{W_1}\otimes F_{W_2, W_3}=\overline{\alpha_{\rm HL}}\circ F_{W_1\boxtimes_{P(z_1-z_2)}W_2, W_3}\circ
F_{W_1, W_2}\otimes 1_{W_3},\end{equation}
where  $\overline{\alpha_{\rm HL}}$ is the extension of $\alpha_{\rm HL}$ to the algebraic completion, $F_{W_1, W_2}$, $F_{W_2, W_3}$  are the universal $P(z_1-z_2)$- and $P(z_2)$- intertwining maps defined in    (\ref{F_2})
and
(\ref{F_1})  respectively and
$F_{W_1\boxtimes_{P(z_1-z_2)}W_2, W_3}$ and $F_{W_1, W_2\boxtimes_{P(z_2)}W_3}$ are
 the   universal $P(z_2)$-, $P(z_1)$- intertwining maps 
\begin{equation}F_{W_1\boxtimes_{P(z_1-z_2)}W_2, W_3}: (W_1\boxtimes_{P(z_1-z_2)}W_2)\otimes W_3\to\overline{(W_1\boxtimes_{P(z_1-z_2)}W_2)\boxtimes_{P(z_2)}W_3},\end{equation}
\begin{equation}F_{W_1, W_2\boxtimes_{P(z_2)}W_3}: W_1\otimes (W_2\boxtimes_{P(z_2)}W_3)\to \overline{W_1\boxtimes_{P(z_1)}(
W_2\boxtimes_{P(z_2)}W_3)}\end{equation}
defined by       naturality in the first and second variable from (\ref{additivity2}),  (\ref{additivity1}) respectively.
\end{thm}

We next define pointed tensor products on the representation category of the Zhu algebra ${\rm Rep}(A(V))$
with irreducible objects labelled by objects of ${\rm Rep}(V)$ by Zhu theorem \cite{Zhu}, inspired 
by   the methods of Prop. 12.5 in \cite{Huang_LepowskiIII}.

\begin{defn}\label{Zhu_as_a_vertex_tensor_category} Consider the vertex tensor category structure of ${\rm Rep}(V)$ transported to $A(V)$ as follows. For a non-zero $z\in{\mathbb C}$, and irreducible $V$-modules $W_1$ and $W_2$, set
$$(W_1)_{(0)}\boxtimes_{P(z)} (W_2)_{(0)}:=(W_1\boxtimes_{P(z)}W_2)_{(0)},$$
and let $$F_0=F_0^{HL}: (W_1)_{(0)}\otimes (W_2)_{(0)}\to (W_1)_{(0)}\boxtimes_{P(z)} (W_2)_{(0)}$$
be the composition of Huang-Lepowsky   universal $P(z)$-intertwining map
$$F=F^{HL}: W_1\otimes W_2\to \overline{W_1\boxtimes_{P(z)} W_2}$$ with inclusion and projection onto top level subspaces.
Let 
$$Z:{\rm Rep}(V)\to {\rm Rep}(A(V))$$
be Zhu functor. 
Let us endow ${\rm Rep}(A(V))$ with associativity morphisms
 \begin{equation}\label{associativity_Zhu} Z(\alpha_{\rm HL})=\sum F_{\alpha, \beta}^{\beta', \alpha'}(1_{W_1}\boxtimes_{P(z_1)} Z(\xi_\beta))\circ Z(\xi_\alpha)\circ Z(\eta_{\beta'})\circ (Z(\eta_{\alpha'})\boxtimes_{P(z_2)} 1_{W_3}).
 \end{equation} 
\end{defn}

 Let us endow $A(V)$ with the transported structure as above. 
 
 \begin{rem}\label{important} Consider the natural transformation associated to
$Z$ given by the identity on the pointed tensor products of irreducible objects.
This setting has been described in Sect. 6 in an abstract form in the case of   tensor categories.
Theorem \ref{strict_equivalence} gives conditions under which this natural transformation becomes a tensor equivalence. On the other hand, these conditions are met by
an extension of the methods of Sect. \ref{12}, regarding
the transport of   tensor structure from a   tensor category ${\rm Rep}(V)$ 
to the Zhu algebra ${\rm Rep}(A(V))$ to the case of vertex tensor categories. In the following, we shall keep in mind
the form of the associativity morphisms explained in the proof of Theorem \ref{transportability}, which is conceptually more useful than (\ref{associativity_Zhu}), although the relevance of the latter resides in the possibility of defining
the OPE coeffiecients $F_{\alpha, \beta}^{\beta', \alpha'}$.
\end{rem}

 \begin{thm}\label{equivalence_of_vertex_tensor_categories}
 $Z: ({\rm Rep}(V), \boxtimes_{P(z)}, \alpha_{\rm HL}) \to ({\rm Rep}(A(V)), \boxtimes_{P(z)}, Z(\alpha_{\rm HL}))$ becomes an equivalence of vertex tensor categories.
 \end{thm}

\begin{rem} An approach to the construction of Def. \ref{Zhu_as_a_vertex_tensor_category} emphasizing Zhu functor $Z$ and its  right inverse has been considered in \cite{McRae}.  \end{rem}

\subsection{An equivalence between the tensor category structure on ${\rm Rep}(A({V_{{\mathfrak g}_k}}))$ from quantum groups  as in Theorem \ref{Zhu_as_a_compatible_unitary_wqh} and from Huang-Lepowsky (vertex) tensor category}\label{34.4}

In this subsection we construct
 an identification of (quasi-tensor equivalence between)  
   the pre-tensor structure of   the transported structure  and that arising from quantum groups following Theorem
   \ref{Zhu_as_a_compatible_unitary_wqh} and conclude with the second step mentioned in Sec. \ref{Interlude}, that is an identification of the two associativity morphisms,
   one arising from Huang-Lepowsky theory and the other from Theorem \ref{Zhu_as_a_compatible_unitary_wqh}.

For each nonzero $z\in{\mathbb C}$ there is a unique isomorphism depending on $z$ 
$$ \eta^z_{Z(W_1), Z(W_2)}\in{\rm Hom}_{A(V_{{\mathfrak g}_k})}(Z(W_1)\boxtimes Z(W_2), Z(W_1)\boxtimes_{P(z)}Z(W_2))$$ such that the following diagram commutes,

\begin{equation}\label{cd}
 \begin{tikzcd}
Z(W_1)\otimes Z(W_2) \arrow{r}{F_0}\arrow[swap] {dr} {F^{HL}_0} & {Z(W_1)\boxtimes Z(W_2)} \arrow{d}{{\eta}^z_{Z(W_1), Z(W_2)}} \\
 & Z(W_1)\boxtimes_{P(z)}Z(W_2)
 \end{tikzcd}
\end{equation}
 
 The map $\eta^z$ exists by  universality of $\boxtimes$-tensor product, see Theorem \ref{FZ_condition}. 
 We give a direct description of $\eta^z$.
 
 In the case of rational vertex operator algebras, Prop. 12.5 in \cite{Huang_LepowskiIII} describes
$P(z)$-tensor product of $V$-modules. It follows   that
 for a non-zero complex number $z$, and irreducible $V$-modules $W_i$, 
\begin{equation}\label{decomposition1}Z(W_1)\boxtimes_{P(z)} Z(W_2)=\oplus_j ({\mathcal M}(z)_{W_1, W_2}^{M_j})^*\otimes Z(M_j)
\end{equation} is a $V$-module,
where $M_j$ is a complete family of irreducible $V$-modules and ${\mathcal M}(z)_{W_1, W_2}^{M_j}$ is the space
of $P(z)$-intertwining maps of type  ${{M_j}\choose{W_1\ W_2}}$.

Consider the case $V=V_{{\mathfrak g}_k}$. Then by Theorem \ref{Zhu_as_a_compatible_unitary_wqh}, ${\rm Rep}(A(V_{{\mathfrak g}_k}))$ is a tensor category. Recall that   the tensor product bifunctor is denoted by $\boxtimes$.
Thus we have a decomposition into isotypic components,
\begin{equation}\label{decomposition2}Z(W_1)\boxtimes Z(W_2)=
\oplus_j ({\rm Hom}_{A(V_{{\mathfrak g}_k})}(Z(W_1)\boxtimes Z(W_2), Z(M_j)))^*\otimes Z(M_j).
\end{equation}

Recall from Sect. \ref{22}  the notion of {\it primary field} $\phi$
and {\it initial term} $T_\phi$ in the setting of affine Lie algebras,  first introduced in \cite{KZ} in conformal field theory.
See also Sect. \ref{33} for an exposition in the setting of vertex operator algebras. For the affine vertex operator algebra associated to an affine Lie algebra, the corresponding notions correspond in a natural way.   

The initial term corresponds  to the map denoted by $f_{FZ}$ in Theorem
\ref{FZ_condition}, is a ${\mathfrak g}$-morphism satisfying Frenkel-Zhu condition, see also Remark \ref{FZ_condition_remark}. (These two properties, that is
faithfulness of $f\to f_{FZ}$ and characterization of the image are specific for the case affine vertex operator algebras
at positive integer level.)

By a  combination of Prop. \ref{intertwining_map_intertwining_operator0} with   Theorem 3.2.3 in \cite{Frenkel_Zhu}, 
 ${\mathcal M}(z)_{W_1, W_2}^{M_j}$ is   isomorphic to the space
of ${\mathfrak g}$-morphisms $Z(W_1)\otimes Z(W_2)\to Z(M_j)$ that satisfy Frenkel-Zhu condition.

 An isomorphism is given as follows. Let $f\in {\mathcal M}(z)_{W_1, W_2}^{M_j}$ be a $P(z)$ intertwining map.
Consider the corresponding intertwining operator ${\mathcal Y}$,   and associate to it the corresponding primary field $\phi$. Consider
 the initial term $T$ of $\phi$ and the morphism $\hat{T}$ in $({\rm Rep}(A(V_{{\mathfrak g}_k})), \boxtimes)$ defined in Theorem \ref{tensor_product}.
 
 We obtain an
 isomorphism by composition
 \begin{equation}\label{isomorphism_alpha} \alpha_{W_1, W_2}^{M_j}: f\in {\mathcal M}(z)_{W_1, W_2}^{M_j}\to \hat{T}\in {\rm Hom}_{A(V_{{\mathfrak g}_k})}(Z(W_1)\boxtimes Z(W_2), Z(M_i)).
 \end{equation}
 Consider the dual isomorphism between the dual spaces, rescaled by $z^{-\Delta_{\phi}}$, with $\Delta_\phi$ the conformal weight of $\phi$,
 $$\eta^{z, Z(M_j)}_{{Z(W_1), Z(W_2)}}: \phi\in ({\rm Hom}_{A(V_{{\mathfrak g}_k})}(Z(W_1)\boxtimes Z(W_2), 
Z(M_j)))^*\to z^{-\Delta_\phi}\phi\circ \alpha_{W_1, W_2}^{M_j}\in ({\mathcal M}(z)_{W_1, W_2}^{M_j})^*.$$
With respect to (\ref{decomposition1}), (\ref{decomposition2}), set
$$\eta^{z}_{{Z(W_1), Z(W_2)}}:=\oplus_j \eta^{z, Z(M_j)}_{{Z(W_1), Z(W_2)}}\otimes 1_{Z(M_j)}.$$
Then $\eta^{z}_{{Z(W_1), Z(W_2)}}$ is an isomorphism in ${\rm Rep}(A(V_{{\mathfrak g}_k}))$ between the two tensor product representations
$$\eta^{z}_{{Z(W_1), Z(W_2)}}: Z(W_1)\boxtimes Z(W_2)\to Z(W_1)\boxtimes_{P(z)}Z(W_2)$$
and makes the diagram (\ref{cd}) commutative.

 We extend $\eta^z$ to an (invertible) natural transformation on all pairs of objects of ${\rm Rep}(A(V_{{\mathfrak g}_k}))$.
 Then $\eta^z$ is a natural transformation between the two tensor product bifunctors 
  \begin{equation}\label{identity}({\rm Rep}(A(V_{{\mathfrak g}_k})), \boxtimes_{P(z)})\to({\rm Rep}(A(V_{{\mathfrak g}_k})), \boxtimes)
 \end{equation}
  
 \begin{thm}\label{transport}
 Let us consider $({\rm Rep}(A(V_{{\mathfrak g}_k})), \boxtimes_{P(z)}, Z(\alpha_{\rm HL}))$ as a  vertex tensor category     as in Def. \ref{Zhu_as_a_vertex_tensor_category}.  Then $({\rm Rep}(A(V_{{\mathfrak g}_k})), \boxtimes)$
 becomes a tensor category with the unique associativity morphisms $({\rm Rep}(A(V_{{\mathfrak g}_k})), \boxtimes, \hat{\alpha}_{\rm HL})$ making   
 the identity functor with natural transformation $\eta$ of (\ref{identity}) into a tensor equivalence.
 \end{thm}

 Recall from Theorem \ref{equivalence_of_vertex_tensor_categories} that
 $$Z: ({\rm Rep}(V_{{\mathfrak g}_k}), \boxtimes_{P(z)}, \alpha_{\rm HL})\to({\rm Rep}(A(V_{{\mathfrak g}_k})), \boxtimes_{P(z)}, Z(\alpha_{\rm HL}))$$  is an equivalence of vertex tensor categories.
 \begin{cor}\label{tensor_equivalence} Consider $({\rm Rep}(V_{{\mathfrak g}_k}), \boxtimes_{\rm HL}, \alpha_{\rm HL})$ with   tensor category structure
 associated to the vertex tensor category $({\rm Rep}(V_{{\mathfrak g}_k}), \boxtimes_{P(z)}, \alpha_{\rm HL})$
 of \cite{HL_tensor_products_of_modules}, 
sketched  in Theorem \ref{HL}.
Then one obtains a natural equivalence of tensor categories
$$({\rm Rep}(V_{{\mathfrak g}_k}), \boxtimes_{\rm HL}, \alpha_{\rm HL})\to ({\rm Rep}(A(V_{{\mathfrak g}_k})), \boxtimes, \hat{\alpha}_{\rm HL}).$$
\end{cor}

  A   computation using the definition of tensor equivalence
 given  as a special case of (\ref{wtf}) and naturality of $\eta$ gives the following formula
 for $\hat{\alpha}_{\rm HL}$.
 
 \begin{prop}\label{HL_associativity_morphism_for_boxtimes}
 The associativity morphisms of the tensor category $({\rm Rep}(A(V_{{\mathfrak g}_k})), \boxtimes, \hat{\alpha}_{\rm HL})$ are given by
  $$(\hat{\alpha}_{\rm HL})_{Z(W_1), Z(W_2), Z(W_3)}=\sum \ F_{\alpha, \beta}^{\beta', \alpha'} \ 1_{Z(W_1)}\boxtimes\hat{\xi}_\beta\circ \hat{\xi}_\alpha\circ \hat{\eta}_{\beta'}\circ
\hat{\eta}_{\alpha'}\boxtimes 1_{Z(W_3)}
$$
where 
$$\hat{\eta}_{\alpha'}=Z(\eta_{\alpha'})\eta_{Z(W_1), Z(W_2)}\in{\rm Hom}_{A(V_{{\mathfrak g}_k})}(Z(W_1)\boxtimes Z(W_2), Z(M_i)),$$
$$\hat{\eta}_{\beta'}=Z(\eta_{\beta'})\eta_{Z(M_i), Z(W_3)}\in{\rm Hom}_{A(V_{{\mathfrak g}_k})}(Z(M_i)\boxtimes Z(W_3),  Z(M_s)),$$
$$\hat{\xi}_\alpha=  \eta^{-1}_{Z(W_1), Z(M_r)}Z(\xi_\alpha)\in{\rm Hom}_{A(V_{{\mathfrak g}_k})}(Z(M_s), Z(W_1)\boxtimes Z(M_r)),$$
$$\hat{\xi}_\beta=\eta^{-1}_{Z(W_2), Z(W_3)}Z(\xi_\beta)\in{\rm Hom}_{A(V_{{\mathfrak g}_k})}(Z(M_r), Z(W_2)\boxtimes Z(W_3)).$$
 \end{prop}
 
 \begin{proof}
Writing $M$ in place of $Z(M)$ for all modules of the Zhu algebra, 
 $$\hat{Z}(\alpha_{\rm HL})=1_{W_1}\boxtimes\eta^{-1}_{W_2, W_3}\circ\eta^{-1}_{W_1, W_2\boxtimes_{P(z_2)}W_3}\circ
 Z(\alpha_{\rm HL})\circ\eta_{W_1\boxtimes_{P(z_1-z_2)} W_2, W_3}\circ\eta_{W_1, W_2}\boxtimes 1_{W_3}=$$
$$ \sum F_{\alpha, \beta}^{\beta', \alpha'}1_{W_1}\boxtimes\eta^{-1}_{W_2, W_3}\circ\eta^{-1}_{W_1, W_2\boxtimes_{P(z_2)}W_3}\circ
(1_{W_1}\boxtimes_{P(z_1)} Z(\xi_\beta))\circ Z(\xi_\alpha)\circ Z(\eta_{\beta'})\circ$$
$$ (Z(\eta_{\alpha'})\boxtimes_{P(z_2)} 1_{W_3}))\circ\eta_{W_1\boxtimes_{P(z_1-z_2)} W_2, W_3}\circ\eta_{W_1, W_2}\boxtimes 1_{W_3}=$$
$$\sum F_{\alpha, \beta}^{\beta', \alpha'} 1_{W_1}\boxtimes(\eta^{-1}_{W_2, W_3}Z(\xi_\beta))\circ(\eta^{-1}_{W_1, M_r}Z(\xi_\alpha))\circ (Z(\eta_{\beta'})\eta_{M_i, W_3})\circ
(Z(\eta_{\alpha'})\eta_{W_1, W_2})\boxtimes 1_{W_3}=
$$
$$\sum \ F_{\alpha, \beta}^{\beta', \alpha'} \ 1_{W_1}\boxtimes\hat{\xi}_\beta\circ \hat{\xi}_\alpha\circ \hat{\eta}_{\beta'}\circ
\hat{\eta}_{\alpha'}\boxtimes 1_{W_3}.
$$
\end{proof}

To describe the associativity morphisms $\hat{\alpha}_{\rm HL}$ of $({\rm Rep}(A(V_{{\mathfrak g}_k})), \boxtimes)$
we need to describe the
tensor equivalence of Theorem \ref{tensor_equivalence}.

This tensor equivalence is built on the following important  two   constructions on the Zhu algebra $A(V_{{\mathfrak g}_k})$, among the number of constructions that we have outlined.
The representation category of this algebra is endowed with the structure of a vertex tensor category induced by Huang-Lepowsky theory, a family of tensor product bifunctors and associativity morphisms that depend on the
variable $z$, as in Def. \ref{Zhu_as_a_vertex_tensor_category}.

Huang-Lepowsky associativity morphisms $\alpha_{\rm HL}$ are morphisms in ${\rm Rep}(V_{{\mathfrak g}_k})$,
$$\alpha_{\rm HL}: (W_1\boxtimes_{P(z_1-z_2)}W_2)\boxtimes_{P(z_2)}W_3\to W_1\boxtimes_{P(z_1)}(W_2\boxtimes_{P(z_2)}W_2)$$
such that the extension $\overline{\alpha}_{\rm HL}$ is   determined by
the universal intertwining maps defined from complete tensor products of vector spaces $\otimes$
as in Theorem \ref{HL_associativity_vertex}.
By Theorem \ref{HL},
$$F_{W_1\boxtimes_{P(z_1-z_2)}W_2, W_3}\circ
F_{W_1, W_2}\otimes 1_{W_3}: $$
$$W_1\otimes W_2\otimes W_3\to \bigoplus_{i,s} ({\mathcal M(P(z_1-z_2))}_{{W_1\ W_2}}^{M_i})^*\otimes  ({\mathcal M(P(z_2))}_{{M_i\ W_3}}^{M_s})^*\otimes \overline{M_s},$$
$$F_{W_1, W_2\boxtimes_{P(z_2)}W_3}\circ 1_{W_1}\otimes F_{W_2, W_3}:$$
$$W_1\otimes W_2\otimes W_3\to \bigoplus_{r, s} ({\mathcal M((P(z_2))}_{W_2, W_3}^{M_r})^* \otimes ({\mathcal M}(P(z_1))_{W_1, M_r}^{M_s})^*\otimes \overline{M_s}.$$
  Using the indicated decomposition into irreducible $V$-modules, the top level subspace of the domain of $\alpha_{\rm HL}$
is obtained in two stages as the direct sum of the top level subspaces of the decomposition of
$W_1\boxtimes_{P(z_1-z_2)} W_2$ and $M_i\boxtimes_{P(z_2)}W_3$. 
The top level subspace of $W_1\boxtimes_{P(z_1-z_2)} W_2$ that we are interested in for the Zhu algebra only arises from the restriction of the universal
$P(z_1-z_2)$-intertwining map
$F_{W_1, W_2}$ to $Z(W_1)\otimes Z(W_2)$, and similarly for the top level subspace of $M_i\boxtimes_{P(z_2)} W_3$.

The result is
$$\bigoplus_{i,s} ({\mathcal M(P(z_1-z_2))}_{{W_1\ W_2}}^{M_i})^*\otimes  ({\mathcal M(P(z_2))}_{{M_i\ W_3}}^{M_s})^*\otimes Z(M_s).$$
By definition, $Z(\overline{\alpha}_{\rm HL})$ is the restriction of $\overline{\alpha}_{\rm HL}$ to the top level subspace
of the triple tensor product, see Remark \ref{important} and   the indicated reference therein to  the proof of Theorem \ref{transportability} for the associativity morphisms.
Thus $Z(\alpha_{\rm HL})$ acts on vectors of the kind
$$(F^0_{\rm HL})_{W_1\boxtimes_{P(z_1-z_2)}Z(W_2), Z(W_3)}\circ
(F^0_{\rm HL})_{Z(W_1), Z(W_2)}\otimes 1_{Z(W_3)}(w^{(1)}\otimes w^{(2)}\otimes w^{(3)}),$$
with $w^{(i)}\in Z(W_i)$.
As before,   $F^0_{\rm HL}$ denotes the restriction of the $P(\zeta)$-intertwining map $F$ to the tensor product of top level subspaces for the suitable $\zeta$.

We have a connection between the two structure maps the variable dependent $F_0^{HL}$ and the constant $F_0$,
given by the natural transformation $\eta^z$ in (\ref{cd}).

Using the natural transformation $\eta^\zeta$ for the suitable $\zeta$, it follows that $\hat{\alpha}_{\rm HL}$ acts 
on vectors of the kind
$$(F^0)_{W_1\boxtimes Z(W_2), Z(W_3)}\circ
(F^0)_{Z(W_1), Z(W_2)}\otimes 1_{Z(W_3)}(w^{(1)}\otimes w^{(2)}\otimes w^{(3)}),$$
where $F^0$ is part of the weak quasitensor structure of Zhu functor constructed in
Theorem \ref{Zhu_as_a_compatible_unitary_wqh}.
In terms of the irreducible decomposition, these are vectors of the space
$$\bigoplus_{i,s} ({\rm Hom}_{A(V_{{\mathfrak g}_k})}({{W_1\boxtimes W_2}}, {M_i}))^*\otimes  ({\rm Hom}_{A(V_{{\mathfrak g}_k})}({{M_i\boxtimes W_3}}, {M_s}))^*\otimes Z(M_s).$$
The action of $\hat{\alpha}_{\rm HL}$ is defined by a formula similar to (\ref{34.42}) with $F_0$, $Z(W_i)$, $\hat{\alpha}_{\rm HL}$   in place of $F$,   $W_i$   $\overline{\alpha_{\rm HL}}$ respectively.
We ask whether this is of CFT-type. By Remark \ref{restriction_of_CFT_type}, we need to show that 
 the inclusion map $G_0$ of weak-quasi-tensor structure
 constructed in Theorem   \ref{Zhu_as_a_compatible_unitary_wqh}, is compatible with $\hat{\alpha}_{\rm HL}$.
 Recall from Sect. \ref{19} that the weak tensor structure associated to Wenzl functor is non-unique, and correspondingly
 $G_0$ is non-unique. One may then try to choose accordingly $G_0$ and try to construct a weak tensor structure for Zhu functor. We shall avoid to go into these details, and we restrict to the specific situation where
 $$(Z(W_1), Z(W_2),  Z(W_3))\in{\mathcal V},$$
 where as before $ {\mathcal V}$ is constituted by the triples of the form $$(V_\lambda, V, V), \quad (V, V_\lambda, V), \quad 
 (V, V, V_\lambda)$$ with $V$ the fundamental representation of ${\mathfrak g}$ and $V_\lambda$ an irreducible representation of ${\mathfrak g}$ in the open Weyl alcove.
 By the results of Sects. \ref{5+}, \ref{5++} we shall content ourselves of showing that $\hat{\alpha}_{\rm HL}$
 arises from ${\mathcal V}$-pre-associator of CFT-type. This is a simpler situation, in that in this case $G_0$ is uniquely
 determined for all Lie type for ${\mathfrak g}\neq E_8$ by the work of Wenzl already widely used in the first part of this paper.
 For pairs of the first two kinds, we shall use the analogue (\ref{34.42}) previously described, while for
 pairs of the third kind we shall use an equation derived from (\ref{34.42}) multiplying in the left by
 $\overline{\alpha_{\rm HL}}^{-1}$.
 In the case ${\mathfrak g}=E_8$ Wenzl give a construction of $(F, G)$ in the setting of quantum groups, and accordingly also for $(F_0, G_0)$, in the proof at the Proposition at page 274 in \cite{Wenzl}.
 Then one may work in the same way in the setting of affine Lie algebras. 
  
 \begin{thm}\label{HL_is_of_CFT_type}
 We have that $\hat{\alpha}_{\rm HL}$  is an  associativity morphism   of ${\rm Rep}(A({V_{{\mathfrak g}_k}}))$ corresponding to a ${\mathcal V}$-pre-associator of CFT-type for the triple $(Z, F_0, G_0)$.
 \end{thm} 

It follows from Theorem \ref{claim1} that the proof of the two steps described in Sect. \ref{Interlude} is complete, and with them the proof
of part c) of Theorem \ref{Zhu_as_a_compatible_unitary_wqh} and \ref{Finkelberg_HL}.

\begin{rem} \label{Mooore_Seiberg} The   induced Huang-Lepowsky 
vertex tensor category construction on ${\rm Rep}(A(V_{{\mathfrak g}_k})$   with the vertex tensor category structure 
may be regarded as associated to a variable dependent family
of weak quasi-bialgebra construction
on $A(V_{{\mathfrak g}_k})$, a variable dependent family of coproducts and associators $(\Delta_z, \Phi_{(z_1-z_2), z_2}^{z_1, z_2})$. It is natural to ask whether this is related to Moore-Seiberg conjecture on a quantum group interpretation in conformal field theory in Section 9 of 
\cite{Moore-Seiberg2}. We shall come back to this in more detail in Remark \ref{on_tensor_product_theorem_for_VOA}.
\end{rem}

\section{Primary fields for affine VOAs, products, iterates, initial terms} \label{33}

   In this section, we complement the discussion of the associativity morphisms of the Huang-Lepowsky tensor product theory 
  in the case of an affine vertex operator algebra $V_{{\mathfrak g}_k}$ 
  at a positive integer level $k$. 
    Central notions     are composition of primary fields 
  operators. As for primary fields themselves, compositions are uniquely determined by restriction and compression to top level spaces. This gives rise to a single variable operator valued function that solves a linear one variable differential equation, known as the reduced KZ equation first studied in \cite{KZ}, \cite{Tsuchiya_Kanie}.
  The initial term is  understood in the setting of a single variable differential equation with regular singular points.
    
In the setting of vertex operator algebras,  we are interested in a subclass of the class of primary fields following the definitions in
 \cite{FHL}. 
 Primary fields in the sense of \cite{FHL} are intertwining operators  evaluated on a lowest weight vector $w$ of the charge module for the Virasoro algebra (i.e. $w$ is an eigenvector of $L_0$ and satisfies $L_nw=0$, $n>0$).
The subclass we are interested in satisfies in addition the requirement that $w$ has
 the lowest  eigenvalue of the action of $L_0$ on the module.
  
   In other words,
  if $W$ is an irreducible charge module for an intertwining operator ${\mathcal Y}$, we are interested in the primary field ${\mathcal Y}_r$ obtained restricting ${\mathcal Y}$ to $W_{(0)}$. These correspond to the primary fields  in conformal field theory \cite{Schellekens}, and also correspond to those discussed in several  papers in mathematics cited in the introduction of \cite{Huang_Lepowski_affine}, to which we refer the reader.
    
  The primary field ${\mathcal Y}_r$ in turn gives rise to an  {\it initial term} ${\mathcal Y}_0$ that we define as a linear map from a tensor product of two lowest weight spaces of irreducible modules
  to a lowest weight space of another irreducible module, that is our main object of study.
  
 We may summarize the main ideas as follows. The space of intertwining operators between three irreducible
 $V_{{\mathfrak g}_k}$ modules is naturally isomorphic to the space of primary fields as above of the same type.
We transport Huang-Lepowsky   associativity morphism to  primary fields, and from this 
 to the operator algebra of products of primary fields associated to the OPE of primary fields, and finally to the spaces of their initial terms. In this way, we are working on the Zhu algebra.

  \bigskip

\subsection{Primary fields   for vertex operator algebras associated to the lowest eigenvalue of $L_0$, initial term.}\label{35.1} Let $W$ be a $V$-module and $S\subset W$ a subset of homogeneous elements.
The linear span $W_S$ of elements image of $S$ under finite products of 
coefficients of the vertex operators $Y_W(a, x)$, with $a$ varying in $V$
(i.e. elements of the form $(a^{(1)})^W_{(n_1)}\dots (a^{(p)})^W_{(n_p)} s$ with $a^{(j)}\in V$, $s\in S$) is a $V$-submodule  
of $W$ with vertex operators $Y_{W_S}(a, x)$ given by the restriction of $Y_{W}(a, x)$ to $W_S$, for $a\in V$.

\begin{defn}
The module $W_S$ is called the submodule of $W$ generated by $S$  \cite{FHL}.
\end{defn}

 \begin{rem}\label{generation_irreducible_module}
If $W$ is irreducible then it is generated by a nonzero vector $w\in W_{(0)}$, the eigenspace of $W$ with the lowest eigenvalue of $L_0^W$ (\ref{grades}),
the conformal weight $\Delta_{W}$ (\ref{conformal_weight_of_module}).
\end{rem}

Let $W_1$, $W_2$, $W_3$ be   $V$-modules and let ${\mathcal Y}$ be
an intertwining operator of type
 ${{W_3}\choose{W_1\ W_2}}$. The following explanation reports Remark 5.4.6 in \cite{FHL}.
 The Jacobi identity for intertwining operators implies for $a\in V$,
 $$
 {\mathcal Y}(Y_{W_1}(a, x_0)w^{(1)}, x_2)=$$
 $$Y_{W_3}(a, x_0+x_2){\mathcal Y}(w^{(1)}, x_2)+{\mathcal Y}(w^{(1)}, x_2)(Y_{W_2}(a, x_2+x_0)-
 Y_{W_2}(a, x_0+x_2)).
 $$
 It follows that if $W_1$ is generated by a subset $S$, then the intertwining operator  ${\mathcal Y}(w^{(1)}, x)$
 on any vector $w^{(1)}\in W_1$ is uniquely determined by the intertwining operators ${\mathcal Y}(s, x)$
 with $s$ varying in $S$. 
 
\begin{rem}\label{determining_intertwining_1} Combining with Remark \ref{generation_irreducible_module},
 it follows that
if $W_1$ is irreducible then an intertwining operator ${\mathcal Y}(w^{(1)}, x)$ of type
 ${{W_3}\choose{W_1\ W_2}}$ is determined
by the {\it restriction} to a nonzero vector of the top space of $W_1$.
$${\mathcal Y}_r: (W_1)_{(0)}\to {\rm Hom}(W_2, W_3)\{x\}.$$
\end{rem}
\begin{defn}\label{primary_field_VOA}
We shall refer to a field of the form ${\mathcal Y}_r$ as a {\it primary field}, occasionally denoted by 
$$\phi:={\mathcal Y}_r.$$ We shall use the same notation
in the case where $x$ is substituted by a non-zero complex number.
\end{defn}

Any vector $w^{(1)}$ in the top space $(W_1)_{(0)}$
is a lowest weight vector for the Virasoro algebra $L^{W_1}_n$, that is 
$$L^{W_1}_nw^{(1)}=0, \quad \quad n>0,\quad w^{(1)}\in (W_1)_{(0)}.$$

 Using arguments similar to those for the vertex operators of a given module $M$, it follows from the
 Jacobi identity for intertwining operators that ${\mathcal Y}_r(w^{(1)}, x)$ satisfies a generalization of the differential equation
 (\ref{differential_equation_for_primary_field}), (the {\it equation of motion} in the case of affine
 Lie algebras \cite{KZ}, \cite{Tsuchiya_Kanie}, see also Remark \ref{equation_of_motion})
 \begin{equation}\label{equation_of_motion_VOA} [L_m, {\mathcal Y}_r(w^{(1)}, x)]=x^{m}(x\frac{d}{dx}{\mathcal Y}_r(w^{(1)}, x)+(m+1)\Delta_{W_1}{\mathcal Y}_r(w^{(1)}, x)), \quad m\in{\mathbb Z}, \ w^{(1)}\in(W_1)_{(0)}\end{equation}
 Thus the restricted intertwining operator ${\mathcal Y}_r(w^{(1)}, x)$ satisfies the definition of  {\it primary field} given in \cite{FHL}.
Indeed   as in the previous case,  this differential equation follows from the translation property of intertwining operators,
the lowest weight property for the Virasoro algebra of elements of $(W_1)_{(0)}$, and the following
 version of (\ref{consequence_of_Jacobi_for_modules}) extended to intertwining operators,

\begin{equation}\label{consequence_of_Jacobi_for_intertwining}
     [a_{(m)}, {\mathcal Y}(w^{(1)}, x)]=\sum_{i=0}^\infty \binom{m}{i}{\mathcal Y}(a^{W_1}_{(i)}w^{(1)}, x)x^{m-i},\quad a\in V, \ w^{(1)}\in W_1.\end{equation}
 The brackets at the left hand side of (\ref{equation_of_motion_VOA}) and (\ref{consequence_of_Jacobi_for_intertwining}) are defined  using the coefficients 
     of the vertex operators $Y(a, x_1)$ associated to $W_2$ and $W_3$, for $a\in V$
     (see also (5.4.17) in \cite{FHL}.)
     Note that (\ref{consequence_of_Jacobi_for_intertwining}) implies that
when $W_2$ is irreducible, an intertwining operator ${\mathcal Y}$ of type
 ${{W_3}\choose{W_1\ W_2}}$ is also determined by its action on a nonzero vector of $W_2$.
 
  Assume that $W_1$, $W_2$, $W_3$ are irreducible let us define the coefficients $w^{(1)}_n$ of an intertwining operator ${\mathcal Y}(w^{(1)}, x)$
   of type
 ${{W_3}\choose{W_1\ W_2}}$   as in (\ref{intertwining_irreducible_case}). By  
Remark   \ref{determining_intertwining_1}, and (\ref{grading_for_intertwining_irreducible}), the full intertwining operator 
${\mathcal Y}$ is uniquely determined by the   family of homomorphisms
\begin{equation}\label{grading_for_coefficients_of_intertwining_convenient_topspace}
 {w^{(1)}}_n: (W_2)_{(k)}\to (W_3)_{(k -n)}, \quad\quad w^{(1)}\in (W_1)_{(0)}, \ k\in{\mathbb Z}_{\geq0}, n\in{\mathbb Z}.
  \end{equation}
(These relations correspons to relations (\ref{grading_loop_group}) in the setting of primary fields for loop groups.)
 
In particular consider the {\it constant term} of the series factor in (\ref{intertwining_irreducible_case}),
$${w^{(1)}}_0: (W_2)_{(0)}\to (W_3)_{(0)},\quad\quad w^{(1)}\in(W_1)_{(0)}.$$
This term defines a linear map between finite dimensional spaces associated to
${\mathcal Y}$,
 \begin{equation}{\mathcal Y}_0: (W_1)_{(0)}\otimes (W_2)_{(0)}\to (W_3)_{(0)}\end{equation}
 \begin{equation}\label{product3} {\mathcal Y}_0(w^{(1)}\otimes w^{(2)})={w^{(1)}}_0(w^{(2)}).\end{equation}
 Borrowing terminology from the case of affine Lie algebras, see Subsect. \ref{32.3} we give the following definition.
 
 \begin{defn}\label{initial_term_def_VOA} Let ${\mathcal Y}$ be an intertwining operator
 of type ${{W_3}\choose{W_1\ W_2}}$, and ${\mathcal Y}_r$ the associated primary field.
 We refer to
${\mathcal Y}_0$ as the {\it initial term}
 of ${\mathcal Y}_r$.
 \end{defn}
 
 The following simple remark shifts the viewpoint and will play in important role.
 
 \begin{prop}\label{shift} Let $W_i$ be irreducible $V$-modules, $i=1$, $2$, $3$.
  The initial term ${\mathcal Y}_0$ of a primary field ${\mathcal Y}_r$
  associated to an intertwining operator ${\mathcal Y}$ of type ${{W_3}\choose{W_1\ W_2}}$
  can equivalently be associated to   Huang-Lepowsky $P(z)$-intertwining map $F$ of the same type   corresponding to ${\mathcal Y}$ as in Prop. \ref{intertwining_map_intertwining_operator0} 
  (for any fixed nonzero $z\in{\mathbb C}$
  and any fixed branch of the complex logarithm)
 as follows,
  $${\mathcal Y}_0=z^{\Delta_{W_1}+\Delta_{W_2}-\Delta_{W_3}}P_0^{W_3}FI_0^{W_1}\otimes I_0^{W_2}: (W_1)_{(0)}\otimes (W_2)_{(0)}\to (W_3)_{(0)}$$
  with $I_0^{W_i}:(W_i)_{(0)}\to W_i$  the canonical inclusion of the top level subspaces into the corresponding modules
  and $P_0^{W_i}:\overline{W_i}\to(W_i)_{(0)}$ the canonical projection from the completed module to the top level subspace. 
   \end{prop}

 In the next subsection, as anticipated in  in Remark \ref{Primary_following_TL}, we discuss  the fact that in the case of the affine vertex operator algebra $V_{{\mathfrak g}_k}$ associated to a loop group $LG$, the two definitions
 of primary field Def. \ref{primary_field}, \ref{primary_field_VOA} and initial terms Def. \ref{initial_term_def}, \ref{initial_term_def_VOA} are in correspondence. Moreover, in Remark
  \ref{on_tensor_product_theorem_for_VOA} we give a re-interpretation of the tensor product Theorem   \ref{tensor_product} in the setting of vertex operator algebras. 
 \bigskip

 \subsection{Primary fields   of the affine vertex operator algebra $V_{{\mathfrak g}_k}$. On a conjecture by Moore and Seiberg}\label{35.2}
  In Sect. \ref{21}, Subsect. \ref{30.1} we have recalled the basic facts of an affine Lie algebra $\hat{\mathfrak g}$ and the associated affine vertex operator algebra $V_{{\mathfrak g}_k}$ with level $k$ a positive integer, and their irreducible $V$-modules. 
  
  In Sect. \ref{22} we have recalled the basis of  loop group CFT, and we have introduced the notion of primary field following  \cite{Wassermann} \cite{Toledano_laredo}. 
  
  In  Subsect. \ref{35.1} we have discussed primary fields for a vertex operator algebra associated to a 
  vector in the charge space of lowest eigenvalue for the conformal Hamiltonian $L_0$.
  The following result is well known. Since it will be useful to us,  we  sketch of proof.
  
 \begin{thm} \label{correspondence_intertwining_initial_affine_case} If $V=V_{{\mathfrak g}_k}$, under Frenkel-Zhu correspondence of affine Lie algebras and affine vertex operator algebras, the two notions of primary field given in Sect. \ref{22} and Subsect.  \ref{35.1}
 coincide. In this correspondence $\phi$ becomes ${\mathcal Y}_r$ and the initial term $T_\phi$ becomes ${\mathcal Y}_0$. Thus  ${\mathcal Y}_0$ determines ${\mathcal Y}_r$ and ${\mathcal Y}$.
  The space of initial terms ${\mathcal Y}_0$ has been characterized in Theorem 2.3 in \cite{Tsuchiya_Kanie} for
  ${\mathfrak sl}_2$ and Theorem 3.2.3 in \cite{Frenkel_Zhu} for a general ${\mathfrak g}$, and described  in terms of a universal tensor product property of modules of $A(V_{{\mathfrak g}_k})$ in Theorem \ref{tensor_product}. Moreover, the maps ${\mathcal Y}\to{\mathcal Y}_r$, ${\mathcal Y}_r\to{\mathcal Y}_0$ are linear isomorphisms.
 \end{thm}
 
 \begin{proof} We need to take into account the correspondence between affine Lie algebras and affine vertex operator algebras 
 established by Frenkel and Zhu \cite{Frenkel_Zhu}.
 We have already discussed injectivity of ${\mathcal Y}\to{\mathcal Y}_r$ in a general setting in Remark \ref{determining_intertwining_1}
 and of $\phi\to T_\phi$ for a primary field in the setting of affine Lie algebras in Prop. \ref{initial_term}.
 One needs to verify 
  that a primary field in the sense of Sect. \ref{22}, Subsect. \ref{32.3}, is the restriction of an intertwining operator,
  for a discussion see Sect. 8 in \cite{GuiII} with reference to the original works. The space of initial terms ${\mathcal Y}_0$ has been characterized
   by Frenkel and Zhu \cite{Frenkel_Zhu}, and this gives surjectivity of ${\mathcal Y}\to{\mathcal Y}_0$.
   Validity of
  (\ref{primary_field_eq2}) follows from Jacobi identity for intertwining operators. We have discussed derivation of the equation of motion in Remark \ref{equation_of_motion} and Subsect. \ref{35.1}.
In theorem \ref{tensor_product} we have discussed a description of the space of initial terms using the coproduct
of the Zhu algebra.
\end{proof}

The following table clarifies the correspondence of notation of primary fields from Sect. \ref{22} in the loop group setting  and of Sect.
\ref{32} in the vertex operator algebra setting   for $w^{(1)}\in(W_1)_{(0)}$,

$$ \text{ charge } V_h\to \text{top level } (W_1)_{(0)} \text{ of  charge } W_1  $$

 $$\text{type} \quad {{h}\choose{ji}}\to \text{type} \quad {{W_3}\choose{W_1\ W_2}} $$
 
 $$ \phi(v, t)=\sum_{n\in{\mathbb Z}}\phi(v, n)t^{-n-\Delta_\phi}\to
   {\mathcal Y}_r(w^{(1)}, z)=\sum_{n\in{\mathbb Z}}w^{(1)}_n z^{-n-\Delta_{W_1}-\Delta_{W_2}+\Delta_{W_3}}.$$

   $$\phi(v, n)\to w^{(1)}_n,$$
   
    $$ \text{source } {\mathcal H}_i\to \text{source } W_2, \quad\text{target } {\mathcal H}_j\to\text{target } W_3$$

   $$\phi(v, n): {\mathcal H}_i[k]\to{\mathcal H}_j[k-n]$$
   
\begin{equation}\label{grading_condition}
w^{(1)}_n: (W_2)_{(k)}\to (W_3)_{(k-n)}
\end{equation}
  Recall that by definition,   
   \begin{equation}\label{minimum_energy_subspace}
  W_{(n)}=0\quad \text{for } n\in{\mathbb Z}_{<0}.\end{equation}
  
   \begin{rem} \label{on_tensor_product_theorem_for_VOA}
We remark here   on a version  of Theorem  \ref{tensor_product} in the setting of
 vertex operator algebras.   We conclude with a structure that seems closely related to
 Moore-Seiberg conjecture on quantum groups arising from conformal field theory, in Section 9 in \cite{Moore-Seiberg2}.
 
The initial term  ${\mathcal Y}_0$ of an intertwining operator ${\mathcal Y}$ defined in Subsect. \ref{35.1}, may have a non-trivial
 kernel. Motivated by Theorem \ref{tensor_product} in the case of affine Lie algebras,  one may want to use this fact to
extend the theorem  to a more general situation of vertex operator algebras.
 A necessary condition for the vertex operator algebra is the existence of  a complemented subspace of the kernel of ${\mathcal Y}_0$ when $W_3$ varies that is compatible with the vertex operator algebra action. This would lead to
a description of ${\mathcal Y}_0$ with $W_3$ fixed as
 we have obtained in Theorem \ref{tensor_product} in the case of affine vertex operator algebras, and the   construction of an associated pre-tensor category structure on ${\rm Rep}(A(V))$. A necessary condition for the existence
 of this complemented subspace is
 the assumption that
 ${\rm Rep}(V)$ satisfies the weak dimension function property defined in
 (\ref{weak_dim_funct_inequality}) in Def. \ref{weak_dimension_function} with respect to Huang-Lepowsky tensor category structure. Let us assume for a moment that  the weak dimension function property holds and that a choice of a complemented
 subspace for the kernel of  ${\mathcal Y}_0$ when $W_3$ varies compatible with the vertex operator algebra  is done for each pair of irreducible modules.
 Consider the construction of Prop. \ref{shift}  extended  to $W_1\boxtimes_{P(z)}W_2$ by additivity on the irreducible modules $W_3$. That is, let as in Subsect. \ref{34.4}
  $$F^{HL}_0: (W_1)_{(0)}\otimes (W_2)_{(0)}\to (W_1\boxtimes_{P(z)} W_2)_{(0)}, $$
  be the composition of the tensor product of the two canonical inclusions of the top level subspaces    
  $$ I^{W_1}_0\otimes I^{W_2}_0: (W_1)_0\otimes (W_2)_0\to W_1\otimes W_2$$
  with a $P(z)$-intertwining map
  defining the tensor product module
  as in Def. \ref{HL_P(z)_tensor_product},
  $$F: W_1\otimes W_2\to \overline{W_1\boxtimes_{P(z)}W_2}$$
  and composed   with the extension to the algebraic completion of the canonical projection onto the top space of the tensor product module,
  $$ P_0^{{W_1\boxtimes_{P(z)}W_2}}: \overline{W_1\boxtimes_{P(z)}W_2}\to  {(W_1\boxtimes_{P(z)} W_2)_{(0)}}.$$
 Then  we may complete $F^{HL}_0$ to   a pair $(F^{HL}_0, G^{HL}_0)$, with
  $$G^{HL}_0:
(W_1\boxtimes_{P(z)} W_2)_{(0)}\to (W_1)_{(0)}\otimes (W_2)_{(0)}$$
 such that $$F^{HL}_0G^{HL}_0=1$$ in a way that both $F^{HL}_0$ and $G^{HL}_0$ factor through the chosen complemented 
 submodule of the kernel of ${\mathcal Y}_0$.
 We may then define a coproduct on
   the Zhu algebra $A(V)$ defined by $(F^{HL}_0, G^{HL}_0)$, and by construction ${\rm Rep}(A(V))$ identifies with ${\rm Rep}(V)$ as a pre-tensor category.
   
   We shall not  study   a notion of   complemented subspace of the
   kernel of ${\mathcal Y}_0$ compatible with the action of the vertex operator algebra.
 We content ourselves with the examples  of affine vertex operator algebras at a positive integer level. For them in Subsect. \ref{34.4} we have  related the map $F^{HL}_0$ arising from Huang-Lepowsky tensor product module   with  the map $F_0$ constructed from quantum group theory via the natural transformation $\eta^z$ depending on $z$.  We have also constructed $G_0$ in Theorem \ref{Zhu_as_a_compatible_unitary_wqh}      a), b).
  Accordingly, one may set
 $$G_0^{HL}=G_0(\eta^z)^{-1}.$$
 This pair induces a coproduct $\Delta_z$ on $A(V_{{\mathfrak g}_k})$
 inducing the vertex tensor category structure of ${\rm Rep}(A(V_{{\mathfrak g}_k}))$ introduced
 in Def. \ref{Zhu_as_a_vertex_tensor_category}, see also Remark \ref{important} for more information.
 The associativity morphisms $\hat{\alpha}_{\rm HL}$ lift to an associator
 $\Phi_{(z_1-z_2), z_2}^{z_1, z_2}$ for $A(V_{{\mathfrak g}_k})$ making
$(A(V_{{\mathfrak g}_k}), \Delta_z, \Phi_{(z_1-z_2), z_2}^{z_1, z_2})$
into a structure extending the notion of weak quasi-Hopf algebra in the direction of a vertex tensor category.
  \end{rem}

  \bigskip
  
\subsection{Top level series   for products and iterates of  intertwining operators of vertex operator algebras and initial terms (solutions of the reduced ODE KZ equation). Connection with the pre-tensor structure of ${\rm Rep}(A(V_{{\mathfrak g}_k}))$.}\label{35.3}  In this subsection we   define an analogue of  the construction of
the initial terms from an intertwining operator to a primary field in Subsect. \ref{35.1}, for a product or iterate of intertwininig operators, see the definitions in
 (\ref{product}), (\ref{iterate}). We   consider the case where the formal variables $x_1$ and $x_2$ are replaced by nonzero complex numbers $z_1$ and $z_2$ respectively.
 The main difference   the construction of ${\mathcal Y}_0$ for a single primary field are certain convergent Laurent series of a single variable for a composition primary fields or associated iterates. Of these series we consider the leading terms, and refer to them again as the
{\it initial terms} of products (and iterates) following the terminology of \cite{Tsuchiya_Kanie} in the setting of affine Lie algebras.  They turn out to be initial terms of linear differential equations with regular singular points.
We shall indicate further interesting reading.

In some more detail, we
define two important  one variable convergent power series. The first was originally   studied by Tsuchiya and Kanie \cite{Tsuchiya_Kanie} and draws its origin in the compositions of primary fields and KZ equation in 2-dimensional
conformal field theory \cite{KZ}. In their work appear as $4$-point correlation functions associated  to primary fields
and vectors in the charge spaces. We write them down in an equivalent operatorial form.
These series also play an important role in the work by
Huang and Lepowsky   to establish  {\it associativity of intertwining operators}
as in Def. \ref{associativity_intertwining_operators_def} (the Operator Product Expansion in the setting of vertex operator algebras) \cite{HuangIV}.
These series are important because similarly to the case of the initial term of a single primary field,
they determine the composition (and iterate) uniquely for affine Lie (vertex operator) algebras, see Remark \ref{connection_with_KZ_equation} for further  information.

 These computations have a twofold purpose. On one hand as said they are of expository nature on a very small part of the analytic work   on convergence properties of the series arising from the OPE.     
 The main result of this subsection is the application to the case of affine vertex operator algebras discussed in Corollary \ref{leading_for_affine_VOA}, where we discuss an application to the weak quasi-bi-algebra structure
 of the Zhu algebra constructed in Theorem \ref{Zhu_as_a_compatible_unitary_wqh}. The reader may compare this result with   Prop. 9.5 of \cite{Kirillov3}.

 \begin{defn}\label{top_level_series_of_products}  Let $W_1, \dots, W_5$ be irreducible $V$-modules.
 Let ${\mathcal Y}_1(w^{(1)}, z_1)$ and ${\mathcal Y}_2(w^{(2)}, z_2)$ be intertwining operators of type  ${{W_4}\choose{W_1\ W_5}}$
 and $ {{W_5}\choose{W_2\ W_3}}$ respectively. Then the {\it top level series   of the product} 
 ${\mathcal Y}_1(w^{(1)}, z_1){\mathcal Y}_2(w^{(2)}, z_2)$ is defined as the restriction of the product to vectors
   $$w^{(1)}\in(W_1)_{0}, \quad \quad w^{(2)}\in (W_2)_{(0)}$$ and as the composition of this restriction
   with inclusion and projection maps to top level spaces, that is  
 \begin{equation}\label{initial_term_for_product} ({\mathcal Y}_1{\mathcal Y_2})_{0}(w^{(1)}, w^{(2)}, z_1, z_2):=P^{W_4}_0({\mathcal Y}_1)_r(w^{(1)}, z_1)({\mathcal Y}_2)_r(w^{(2)}, z_2)I^{W_3}_0,
 \end{equation}
 where $P_0^{W_4}$ denotes the canonical projections onto the top level space of $W_4$ defined in (\ref{P_r})
 and $$I_0^{W_3}: (W_3)_{(0)}\to\overline{W_3}$$ the canonical inclusion map of the top level space of $W_3$.
 Set $\Delta_i:=\Delta_{W_i}$. The series development is given by
\begin{equation}\label{development_of_product}({\mathcal Y}_1{\mathcal Y_2})_{0}(w^{(1)}, w^{(2)}, z_1, z_2)=\left[\sum_{n=0}^{-\infty} w^{(1)}_{-n}w^{(2)}_{n} z_1^n z_2^{-n}\right]\  z_1^{-\Delta_{1}-\Delta_{5}+\Delta_{4}}z_2^{-\Delta_{2}-\Delta_{3}+\Delta_{5}}=
\end{equation}
\begin{equation}\label{development_of_product2}
\left[\sum_{m=0}^{+\infty} w^{(1)}_{m}w^{(2)}_{-m}({z_2}/{z_1})^{m+\Delta_{1}+\Delta_{5}-\Delta_{4}}\right]\ z_2^{-\Delta_{1}-\Delta_{2}-\Delta_{3}+\Delta_{4}}
\end{equation}
where 
$$w^{(1)}_{-n}w^{(2)}_{n}: (W_3)_{(0)}\to(W_5)_{(-n)}\to (W_4)_{(0)}.$$
The summation over non positive integers $n$ (non negative integers $m$)   follows from (\ref{grading_condition}) and (\ref{minimum_energy_subspace}). 
The equality between (\ref{development_of_product}) and (\ref{development_of_product2}) follows from a direct computation.
The constant term $w^{(1)}_{0}w^{(2)}_{0}$ of the power series is called the {\it initial term} at $0$.

 \end{defn}\bigskip

 \begin{rem}\label{connection_with_KZ_equation}
The top level series of the product of primary fields has been studied analytically in  \cite{Tsuchiya_Kanie} in the setting of affine Lie algebras for  ${\mathfrak g}={\mathfrak sl}_2$ (called vertex operators in their work). They are an important tool  to study products.  They have been   extended by Huang and Lepowsky 
in the setting of vertex operator algebras in greater generality, see Remark \ref{convergence_and_extension} for some references. The series converges for
 $|z_1|>|z_2|>0$ and the limit defines a weakly holomorphic multivalued operator valued function, which for 
 $V=V_{{\mathfrak g}_k}$ determines the product $({\mathcal Y}_1)_r({\mathcal Y}_2)_r$ uniquely \cite{Tsuchiya_Kanie}, \cite{Frenkel_Zhu}, \cite{Huang_Lepowski_affine}, and hence also $({\mathcal Y}_1{\mathcal Y}_2)$ \cite{FHL}. Moreover, the one-variable function of $\zeta$ defined between square brackets of (\ref{development_of_product2}) setting $\zeta=z_2/z_1$
 satisfies   a {\it reduced}  KZ  differential equation,
 see Prop. 4.3 in \cite{Tsuchiya_Kanie} which is a first order linear ordinary differential equation with regular singular poins at $\zeta=0$, $\zeta=1$. The initial term  of the top level series in the sense of Def. \ref{top_level_series_of_products} is also the initial term as a 
solution of the reduced KZ equation at $\zeta=0$, and determines the solution uniquely.
Thus the initial term of the top level series determines the product  the corresponding    primary fields $({\mathcal Y}_1)_r({\mathcal Y}_2)_r$ (or intertwining operators) uniquely for  $V=V_{{\mathfrak g}_k}$.
 The more symmetric notion  of intertwining operator in the setting of vertex operator algebras  allows     to study iterates  on a symmetric basis and they are an important tool in the work by Huang and Lepowsky
 to study convergence and extension properties.  
 We would like to advertise the paper by
 McRae \cite{McRae}. This paper, having among other things  several points in common with the approach of our paper,
 is a good complement and expands parts that we do not consider in detail.  
 \end{rem}
 
 We   next focus on a corresponding symmetric notion of top level series and initial term for an iterate.

 \begin{defn} Let $W_1$, $W_2$, $W_3$, $W_4$, $W_6$ be irreducible $V$-modules and let
 ${\mathcal Y}_3$ and ${\mathcal Y}_4$ be intertwining operators of types  
 ${{W_6}\choose{W_1\ W_2}}$ and ${{W_4}\choose{W_6\ W_3}}$ respectively.
 Then the {\it top level series   of the iterate}  is defined restricting the iterate to vectors $$w^{(1)}\in(W_1)_{(0)}, \quad\quad w^{(2)}\in (W_2)_{(0)}$$ by
  \begin{equation}\label{initial_iterate}
 ({\mathcal Y}_4({\mathcal Y_3}))_{0}(w^{(1)}, w^{(2)}, z_1, z_2):=P_0^{W_4}{\mathcal Y}_4({\mathcal Y}_3(w^{(1)}, z_1-z_2)w^{(2)}, z_2)I_0^{W_3}.
 \end{equation}
  The series development is given by
 \begin{equation}\label{development_of_iterate} ({\mathcal Y}_4({\mathcal Y}_3))_{0}(w^{(1)}, w^{(2)}, z_1, z_2)=\left[\sum_{p=0}^{-\infty} ({w^{(1)}}_p{w^{(2)})}_{-p}(z_1-z_2)^{-p}z_2^{p}\right]\ (z_1-z_2)^{-\Delta_{1}-\Delta_{2}+\Delta_{6}} z_2^{-\Delta_{6}-\Delta_{3}+\Delta_{4}}=\end{equation}
 \begin{equation}\label{development_of_iterate2}
 \left[\sum_{q=0}^{+\infty} ({w^{(1)}}_{-q}{w^{(2)})}_{q}\left(\frac{1-z_2/z_1}{z_2/z_1}\right)^{q-\Delta_1-\Delta_2+\Delta_6}\right] \ z_2^{-\Delta_{1}-\Delta_{2}-\Delta_{3}+\Delta_{4}}
 \end{equation}
 Similarly, the summation over the non positive integers $p$ (non negative integers $q$) follows from (\ref{grading_condition}) and (\ref{minimum_energy_subspace}) applied to ${w^{(1)}}_pw^{(2)}$, as  this vector lies in $ (W_6)_{-p}$.
It follows that  the grading relations satisfied by the coefficients of the series development of 
${\mathcal Y}_4({w^{(1)}}_pw^{(2)}, z_2)$ follow the more general rule given in  (\ref{grading_for_intertwining_irreducible}).
The equality between (\ref{development_of_iterate}) and \ref{development_of_iterate2}) follows from a direct computation.
The term $({w^{(1)}}_{0}{w^{(2)})}_{0}$ is called the {\it initial term } of the top level series of the iterate.

 \end{defn}
 
 As for intertwining operators, the top level series of products and iterates are understood as multivalued functions when the series is convergent, unless a branch of the complex logarithm is specified.

We identify the constant terms of the factor series with integer powers in terms of the initial terms of the composing intertwining operators.

\begin{prop}\label{Leading_terms} ({\it  Initial terms of the top level power series   of products and iterates of intertwining operators})

Let $W_1,\dots, W_6$ be irreducible $V$-modules.
Let
${\mathcal Y}_1$ and ${\mathcal Y}_2$ be intertwining operators of type  ${{W_4}\choose{W_1\ W_5}}$
 and $ {{W_5}\choose{W_2\ W_3}}$ respectively.
  Let moreover  ${\mathcal Y}_3$ and ${\mathcal Y}_4$ be intertwining operators of types  
 ${{W_6}\choose{W_1\ W_2}}$ and ${{W_4}\choose{W_6\ W_3}}$ respectively.  
Under the correspondence
taking an intertwining operator ${\mathcal Y}$ to the initial term ${\mathcal Y}_0$ described in Subsect. \ref{35.1},
we have that

 \begin{itemize}
\item[{\rm  a)}]
the initial term ${w^{(1)}}_0{w^{(2)}}_0$ of the   series factor for $n=0$ appearing in (\ref{development_of_product}) 
 corresponds to
 \begin{equation}\label{initial_term_of _product3} ({\mathcal Y}_1)_0\circ 1_{(W_2)_{(0)}}\otimes({\mathcal Y}_2)_0: 
 \end{equation}
 \begin{equation}\label{initial_term_of _product2}(W_1)_{(0)}\otimes[(W_2)_{(0)}\otimes (W_3)_{(0)}]\to (W_1)_{(0)}\otimes  (W_5)_{(0)} \to (W_4)_{(0)},
\end{equation}
where $\otimes$ denotes the usual tensor product of finite dimensional complex vector spaces;

\item[{\rm  b)}] the initial term $({w^{(1)}}_0{w^{(2)}})_0$ of the   series factor appearing in
(\ref{development_of_iterate}) for $p=0$ 
 corresponds to
 \begin{equation} \label{initial_term_for_iterate3}  ({\mathcal Y}_4)_0\circ({\mathcal Y}_3)_0\otimes 1:
 \end{equation}
\begin{equation}\label{initial_term_for_iterate2} [(W_1)_{(0)}\otimes(W_2)_{(0)}]\otimes (W_3)_{(0)}\to (W_6)_{(0)}\otimes (W_3)_{(0)}\to(W_4)_{(0)}. 
\end{equation}

\end{itemize}

\end{prop}
\bigskip

 \begin{cor} \label{leading_for_affine_VOA} Set $V=V_{{\mathfrak g}_k},$
 $$ \rho=(W_1)_{(0)}, \quad \sigma=(W_2)_{(0)},\quad
  \tau=(W_3)_{(0)}.$$ Let $\boxtimes$ be  the bifunctor of ${\rm Rep}(A({V_{{\mathfrak g}_k}}))$ and $(F_0, G_0)$ 
   the weak quasi-tensor structure of the forgetful functor $Z$ constructed in Theorem \ref{Zhu_as_a_compatible_unitary_wqh} 
 a), b).   Then 
  \begin{itemize}
 \item[{\rm  a')}] (\ref{initial_term_of _product3}) takes the form
  $$\hat{S}\circ 1\boxtimes \hat{T}\circ (F_0)_{\rho, \sigma\boxtimes \tau} \circ 1_{\rho}\otimes (F_0)_{\sigma, \tau},$$
  $$\hat{T}\in {\rm Hom}_{A({V_{{\mathfrak g}_k}})}(\sigma\boxtimes \tau, (W_5)_{(0)}), \quad\quad
   \hat{S}\in  {\rm Hom}_{A({V_{{\mathfrak g}_k}})}(\rho\boxtimes  (W_5)_{(0)}, (W_4)_{(0)}).
  $$
  
\item[{\rm  b')}]  (\ref{initial_term_for_iterate3}) takes the form 
 $$\hat{S'}\circ  \hat{T'}\boxtimes 1\circ (F_0)_{\rho\boxtimes\sigma,   \tau} \circ (F_0)_{\rho, \sigma}\otimes 1_{\tau},$$
 $$\hat{T'}\in {\rm Hom}_{A({V_{{\mathfrak g}_k}})}(\rho\boxtimes \sigma, (W_6)_{(0)}), \quad\quad
   \hat{S'}\in  {\rm Hom}_{A({V_{{\mathfrak g}_k}})}((W_6)_{(0)}\boxtimes\tau  , (W_4)_{(0)}).
  $$

\end{itemize}

 \end{cor}

\begin{proof}
The proof follows from Proposition \ref{Leading_terms} and Theorem \ref{tensor_product}, taking into account the passage of notation from affine
Lie algebras of Sect. \ref{22} to affine vertex operator algebras recalled   in Subsect. \ref{35.2}.
\end{proof}

  \subsection{OPE of intertwining operators and OPE of primary fields}\label{35.4}

 We assume that associativity of intertwining operators holds, see Def. \ref{associativity_intertwining_operators_def}.
 Let   $W_1$, $W_2$, $W_3$ be fixed irreducible   $V$-modules,
  and $\{M_s\}$ a complete family of irreducible $V$-modules.
  Let ${\mathcal M}_{W_1, W_2, W_3}^{M_s}$ be the linear span of (linearly independent) products of intertwining operators
  ${\mathcal Y}_\alpha{\mathcal Y}_\beta$, let ${\mathcal Y}_{\beta'}({\mathcal Y}_\alpha')$ denote the linearly idependent iterates   of Def   \ref{OPE_intertwining_operators} and
  $F$ be the change-of-basis matrix corresponding to the identity map
   $${\rm id}:
  {\mathcal M}_{W_1, W_2, W_3}^{M_s}\to {\mathcal M}_{W_1, W_2, W_3}^{M_s}$$
  from products ${\mathcal Y}_\alpha{\mathcal Y}_\beta$ to iterates here briefly denoted by ${\mathcal Y}_{\beta'}({\mathcal Y}_{\alpha'})$, whose   coefficients $F_{\alpha, \beta}^{\beta', \alpha'}$ were
  defined in $(\ref{F_coefficients})$.
  
  The {\it OPE of primary fields} is the expression of the product 
  $$\phi_\alpha(w^{(1)}, z_1)\phi_\beta(w^{(2)}, z_2), \quad\quad w^{(i)}\in(W_i)_{(0)},$$
  of primary fields $$\phi_\alpha:=({\mathcal Y}_\alpha)_r, \quad\quad \phi_\beta=({\mathcal Y}_\beta)_r$$ 
  with respect to $${\mathcal Y}_{\beta'}(\phi_{\alpha'}(w^{(1)}, z_1-z_2 )w^{(2)}, z_2)$$
  derived at the bottom right corner of the following diagram

  \begin{equation}\label{OPE_intertwining_primary}
\xymatrix@C=1em{
  \{{\mathcal Y}_\alpha{\mathcal Y}_\beta\}\subset{\mathcal M}_{W_1, W_2, W_3}^{M_s}\ar[d]_{{\rm id}}\ar[r]^{}
  & \phi_\alpha\phi_\beta \ar[d]^{{\rm id}}\\
 \{{\mathcal Y}_{\beta'}({\mathcal Y}_{\alpha'})\}\subset{\mathcal M}_{W_1, W_2, W_3}^{M_s}\ar[r] 
 & \text{OPE of prim. fields}
}
\end{equation}
which gives
$$\phi_\alpha(w^{(1)}, z_1)\phi_\beta(w^{(2)}, z_2)=\sum_{\beta', \alpha'} \frac{F_{\alpha, \beta}^{\beta', \alpha'}}{(z_1-z_2)^{\Delta_{W_1}+\Delta_{W_2}-\Delta_{M_i}}}\phi_{\beta'}(\phi_{\alpha'}(w^{(1)}, 0)w^{(2)}, z_2)+\text{``less sing. terms''},$$

where the disregarded  part is given by

$$\text{``less sing. terms''}=\sum_{k=1}^{+\infty}\sum_{\beta', \alpha'} \frac{F_{\alpha, \beta}^{\beta', \alpha'}}{(z_1-z_2)^{\Delta_{W_1}+\Delta_{W_2}-\Delta_{M_i}-k}}{\mathcal Y}_{\beta'}(\phi_{\alpha'}(w^{(1)}, -k)w^{(2)}, z_2),$$
and $\phi_{\alpha'}(w^{(1)}, -k)$ denotes the mode $\phi_{\alpha'}(w^{(1)}, -k):=(w^{(1)})_{-k}$ of the restricted intertwining operator ${\mathcal Y}_{\alpha'}$ to the top level subspace $(W_1)_{(0)}$ of the charge module.

The conformal dimensions $-\Delta_{W_1}-\Delta_{W_2}+\Delta_{M_i}$ are non negative when ${\rm Rep}(V)$ is a unitary tensor category.
 For   introductory  reading to conformal field theory and primary fields, the reader may consult \cite{Schellekens}, \cite{Walton}.

\subsection{Computing OPE constants for triples in ${\mathcal V}$  via Kirillov-Wenzl unitary structure of $A(V_{{\mathfrak g}_k})$} \label{35.5}
 
 To proceed further, we specialize to the case $V=V_{{\mathfrak g}_k}$.
In this case we have a natural linear isomorphism from the space of intertwining operators with a given irreducible charge module, and irreducible source and target modules (or intertwining maps of the same type) and the space of initial terms.
The correspondence is described in Theorem \ref{correspondence_intertwining_initial_affine_case}.  
Initial terms are morphisms in the category of modules of the Zhu algebra, ${\rm Rep}(A(V_{{\mathfrak g}_k}))$, and this has been described in theorem \ref{tensor_product}.
Recall that by Theorem \ref{Zhu_as_a_compatible_unitary_wqh}, ${\rm Rep}(A(V_{{\mathfrak g}_k}))$ is a unitary tensor category with tensor bifunctor $\boxtimes$.
Let us consider orthonormal bases
\begin{equation}\label{orth_basis1}\{v_\alpha\}\subset{\rm Hom}_{A(V_{{\mathfrak g}_k})}((M_s)_{(0)}, (W_1)_{(0)}\boxtimes (M_r)_{(0)}),\end{equation}
with respect to the inner product of ${\rm Hom}_{A(V_{{\mathfrak g}_k})}((M_s)_{(0)}, (W_1)_{(0)}\boxtimes (M_r)_{(0)},$ given by
$$S^*T=(S, T)1_{M_s}.$$
Similarly, 
let 
\begin{equation}\label{orth_basis2}\{v_\beta\}\subset{\rm Hom}_{A(V_{{\mathfrak g}_k})}((M_r)_{(0)}, (W_2)_{(0)}\boxtimes (W_3)_{(0)}),\end{equation}
\begin{equation}\label{orth_basis3}\{v_{\alpha'}\}\subset{\rm Hom}_{A(V_{{\mathfrak g}_k})}((M_i)_{(0)}, (W_1)_{(0)}\boxtimes (W_2)_{(0)}),\end{equation}
\begin{equation}\label{orth_basis4}\{v_{\beta'}\}\subset{\rm Hom}_{A(V_{{\mathfrak g}_k})}((M_s)_{(0)}, (M_i)_{(0)}\boxtimes (W_3)_{(0)})\end{equation}
be orthonormal bases with respect to analogous inner products.
Correspondingly we choose bases $F_\alpha$, $F_\beta$, $F_{\alpha'}$, $F_{\beta'}$ of Huang-Lepowsky intertwing 
maps of the same types as before, with initial terms indicated as follows by the arrow $\to$
\begin{equation}\label{choice_of_intertwining_map1} F_\alpha\to v_\alpha^*\circ(F_0)_{(W_1)_{(0)}, (M_r)_{(0)}}\quad F_\beta\to v_\beta^*\circ(F_0)_{(W_2)_{(0)},(W_3)_{(0)}},\end{equation}
\begin{equation}\label{choice_of_intertwining_map2} 
F_{\alpha'}\to v_{\alpha'}^*\circ(F_0)_{(W_1)_{(0)}, (W_2)_{(0)}}, \quad F_{\beta'}\to v_{\beta'}^*\circ(F_0)_{{(M_i)_{(0)}, (W_3)_{(0)}}},\end{equation} respectively, by Theorem \ref{FZ_condition}.
Following the notation of Sect. \ref{5+}, we write for brevity
$$(F_0)_{1,2}:=(F_0)_{(W_1)_{(0)}, (W_2)_{(0)}\boxtimes (W_3)_{(0)}}\circ 1\otimes (F_0)_{(W_2)_{(0)}, (W_3)_{(0)}},$$
$$(F_0)_{2,1}:=(F_0)_{{(W_1)_{(0)}\boxtimes(W_2)_{(0)}},   (W_3)_{(0)}}\circ  (F_0)_{(W_1)_{(0)}, (W_2)_{(0)}} \otimes 1.$$
Consider linear isomorphisms of spaces of products and iterates of Huang-Lepowsky intertwining maps with corresponding compositions
of initial terms. By naturality of $F_0$, these isomorphisms are given by linear extensions of the maps taking
products and iterates to their initial term,
\begin{equation} \label{product_Zhu} F_\alpha\circ 1\otimes F_\beta\to v_\alpha^*\circ 1\boxtimes v_\beta^*\circ (F_0)_{1,2}\end{equation}
 \begin{equation}\label{iterate_Zhu} F_{\beta'}\circ F_{\alpha'}\otimes 1\to v_{\beta'}^*\circ v_{\alpha'}^*\boxtimes 1\circ (F_0)_{2,1}\end{equation}
 These maps are isomorphisms as the terms at the left hand side are linear bases by the work
of Huang and Lepowsky previously referred to. Moreover the factors at the
right hand side containing the $v$'s are orthonormal bases, and the maps $F_0$  are right invertible.
The vector spaces on both sides have the same dimension.
In Huang-Lepowsky theory, the structure constants of the associativity morphisms $F_{\alpha, \beta}^{\beta', \alpha'}$
are given by linear expansion of products in terms of iterates, see (\ref{associativity_in_terms_of_F}).
At the level of the Zhu algebra, this corresponds to linearly expand the right hand side of (\ref{product_Zhu})
in terms of the right hand side of (\ref{iterate_Zhu}).

\begin{prop}\label{linear_expansions}
Let $(W_1)_{(0)}$, $(W_2)_{(0)}$, $(W_3)_{(0)}$ be a given triple of irreducible representations of the Zhu algebra
$A(V_{{\mathfrak g}_k})$, and let us fix orthonormal bases $v_\alpha$, $v_\beta$, $v_{\alpha'}$, $v_{\beta'}$ 
of the morphism spaces of the Zhu algebra as in (\ref{orth_basis1}), (\ref{orth_basis2}), (\ref{orth_basis3}), (\ref{orth_basis4}),  resp. and correspondingly intertwining maps $F_\alpha$, $F_\beta$, $F_{\alpha'}$, $F_{\beta'}$ as in (\ref{choice_of_intertwining_map1}), (\ref{choice_of_intertwining_map2}).

 \begin{itemize}
\item[{\rm  a)}]
If    the initial term of  the top level series of Huang-Lepowsky products $F_{\alpha}\circ 1\otimes F_{\beta}$ expands into a 
linear combination of initial terms of the top level series of Huang-Lepowsky iterates $F_{\beta'}\circ F_{\alpha'}\otimes 1$ then    this expansion is given by
associativity morphisms $$\alpha_{(W_1)_{(0)}, (W_2)_{(0)}, (W_3)_{(0)}}: ((W_1)_{(0)}\boxtimes(W_2)_{(0})\boxtimes(W_3)_{(0)}\to (W_1)_{(0)}\boxtimes((W_2)_{(0}\boxtimes(W_3))_{(0)}$$ 
 given by the evaluation on the given triple of the CFT-type pre-associator 
 $$\alpha_{(W_1)_{(0)}, (W_2)_{(0)}, (W_3)_{(0)}}=(F_0)_{1,2}(G_0)_{2,1}.$$
 
 \item[{\rm  b)}]
Similarly, if   the initial term of    iterates $F_{\beta'}\circ F_{\alpha'}\otimes 1$
 expands into a 
linear combination of initial terms of   products $F_{\alpha}\circ 1\otimes F_{\beta}$ then this expansion is given
by associativity morphisms
$$\beta_{(W_1)_{(0)}, (W_2)_{(0)}, (W_3)_{(0)}}: (W_1)_{(0)}\boxtimes((W_2)_{(0)}\boxtimes(W_3)_{(0)})\to ((W_1)_{(0)}\boxtimes(W_2)_{(0)})\boxtimes(W_3)_{(0)}$$  given by
$$\beta_{(W_1)_{(0)}, (W_2)_{(0)}, (W_3)_{(0)}}=(F_0)_{2,1}(G_0)_{1,2}.$$
\item[{\rm  c)}]
Conversely, if the embedding of the initial term of the top level series of a  product $F_\alpha\circ 1\otimes F_\beta$
(or an iterate $F_{\beta'}\circ F_{\alpha'}\otimes 1$ resp.)
in ${\rm Vec}$ is compatible with that of an iterate on a specific triple
$((W_1)_{(0)}, (W_2)_{(0)}, (W_3)_{(0)})$ then the CFT-type pre-associator $(F_0)_{1,2}(G_0)_{2,1}$
(or $(F_0)_{2,1}(G_0)_{1,2}$ resp.)  induces  linear expansions as above.
\end{itemize}
 
\end{prop}

\begin{proof} a)
The elements at the right hand side of (\ref{iterate_Zhu}) are invariant under right multiplication by $(G_0)_{2,1}(F_0)_{2,1}$ by right invertibility of $F_0$.
Thus if the desired expansion exists for a given triple $((W_1)_{(0)}, (W_2)_{(0)}, (W_3)_{(0)})$ then we necessarily have
the compatibility condition described in the statement, 
\begin{equation}\label{emergence_of_CFT_type} v_\alpha^*\circ 1\boxtimes v_\beta^*\circ (F_0)_{1,2}=v_\alpha^*\circ 1\boxtimes v_\beta^*\circ (F_0)_{1,2}
(G_0)_{2,1}(F_0)_{2,1}.\end{equation}
It follows that the linear expansion   is necessarily given by
\begin{equation}\label{Fourier}\sum_{\alpha', \beta'}[v_\alpha^*\circ 1\boxtimes v_\beta^*\circ (F_0)_{1,2}\circ
(G_0)_{2,1}\circ v_{\alpha'}\boxtimes 1\circ  v_{\beta'}]  \circ v_{\beta'}^*\circ v_{\alpha'}^*\boxtimes 1\circ (F_0)_{2,1}.
\end{equation}
The scalars of the linear expansions are given by the terms in square bracket 
\begin{equation}\label{little_eff} f_{\alpha, \beta}^{\beta', \alpha'}:=v_\alpha^*\circ 1\boxtimes v_\beta^*\circ (F_0)_{1,2}\circ
(G_0)_{2,1}\circ v_{\alpha'}\boxtimes 1\circ  v_{\beta'},\end{equation}
this shows that they define the value on the given triple of the CFT-type pre-associator $(F_0)_{1,2}(G_0)_{2,1}$.
b) follows in a similar way. 

c) Conversely, the identification of the initial term of top level series of products   of Huang-Lepowsky intertwining 
operators with linear maps in the non-strict tensor category of vector spaces requires right multiplication on the right of by
the associativity morphism of ${\rm Vec}$. Since this is a morphism in the category of representations of a classical Lie algebra satisfying in addition an $k$-constrained truncated fusion rules, under the realization of the tensor product
bifuntor $\boxtimes$ of the representation category of the Zhu algebra, the initial term of a product must be multiplied
by $(G_0)_{2,1}(F_0)_{2,1}$, cf. the equality (\ref{emergence_of_CFT_type}).
We may now apply the computations (\ref{Fourier}).
\end{proof}

\begin{rem} Note that the compatibility condition refers to a specific triple, and   in that case we  write down the associativity 
morphism of the Zhu algebra as pre-associator of CFT-type. It may not be easy to construct a pre-associator of CFT-type which satisfies all the axioms of a tensor category (i.e. to solve the problem of whether   the Zhu algebra can be made into a weak  Hopf algebra). On the other hand, the quantum groups constructions imply that the Zhu algebra  admits a ${\mathcal V}$-pre-associator of CFT-type and that the corresponding unique associator is equivalent to Huang-Lepowsky theory, since it is so on the collection ${\mathcal V}$. \end{rem}

We  conclude with sufficient conditions under which the assumptions of a) in Prop. \ref{linear_expansions} hold.
These are the conditions under which elements at the right hand side of (\ref{product_Zhu}) admit a linear 
expansion in terms of elements at the right hand side of  (\ref{iterate_Zhu}).
This will give an application
  of Theorem \ref{HL_is_of_CFT_type} on the determination of the OPE constants
$F_{\alpha, \beta}^{\beta', \alpha'}$ from the unitary structure of the weak quasi-Hopf algebra $A(V_{{\mathfrak g}_k})$
on certain triples $((W_1)_{(0)}, (W_2)_{(0)}, (W_3)_{(0)})$ of representations of $A(V_{{\mathfrak g}_k})$ as follows. 

Recall that $V$ here denotes the fundamental of the classical Lie algebra ${\mathfrak g}$ chosen as in \cite{Wenzl}, and not the vertex operator algebra $V_{{\mathfrak g}_k}$.

Recall also that $((W_1)_{(0)}, (W_2)_{(0)}, (W_3)_{(0)})\in{\mathcal V}$  means that
  two entries correspond to 
the fundamental representation $V$, and another to an arbitrary irreducible $(L_{k, \lambda})_{(0)}$ with dominant weight $\lambda$ in the open Weyl alcove at positive integer level $k$ ($k\geq 1$ for ${\mathfrak g}\neq E_8$ and $k\geq 2$ ${\mathfrak g}=E_8$).

We assume that  
$(L_{k, \lambda})_{(0)}$ is either  $(W_1)_{(0)}$, or $(W_2)_{(0)}$.
Then in this case $((W_1)_{(0)}\otimes (W_2)_{(0)})\otimes (W_3)_{(0)}$ admits a unique inner product such that
$(W_1)_{(0)}\otimes (W_2)_{(0)}$ has the invariant inner product under the action of ${\mathfrak g}$ with a fixed normalization on the highest weight vectors, and such that for any dominant weight $\gamma$ that
results from the decomposition into irreducible under the action of ${\mathfrak g}$ such that
$\gamma$ is in the open Weyl alcove, we similarly endow $(L_{k, \gamma})_{(0)}\otimes (W_3)_{(0)}$ with 
the unique invariant inner product under the action of ${\mathfrak g}$.
In this way we have a well defined decomposition of $((W_1)_{(0)}\otimes (W_2)_{(0)})\otimes (W_3)_{(0)}$ 
into orthogonal irreducible subspaces, which is canonical for ${\mathfrak g}\neq E_8$. Selecting only the irreducible components
of $(L_{k, \gamma})_{(0)}\otimes (W_3)_{(0)}$ with dominant weights in the open Weyl alcove, gives
the projection $(F_0)_{2,1}$ from
$((W_1)_{(0)}\otimes (W_2)_{(0)})\otimes (W_3)_{(0)})$ to that addendum
and an   the inclusion $(G_0)_{2,1}$ from that addendum to $((W_1)_{(0)}\otimes (W_2)_{(0)})\otimes (W_3)_{(0)})$.
Then apply Theorem \ref{HL_is_of_CFT_type} which implies the assumptions of a).
Then we may proceed with the computations in the proof of Prop. \ref{linear_expansions}.

From (\ref{little_eff}), 
we have
\begin{equation}\label{little_eff_adjoint} (F_0)_{1,2}(G_0)_{2,1}=\sum_{\alpha, \beta, \alpha', \beta'} f_{\alpha, \beta}^{\beta', \alpha'} 1\boxtimes v_\beta\circ v_\alpha\circ v_{\beta'}^*\circ v_{\alpha'}^*\boxtimes 1.\end{equation}
 Recall the definition of $\eta_\alpha$, $\eta_\beta$, $\eta_{\alpha'}$, $\eta_{\beta'}$, and corresponding
$\xi_\alpha$, $\xi_\beta$, $\xi_{\alpha'}$, $\xi_{\beta'}$   in Subsect. \ref{34.3}, of the two variable natural transformation  $\eta^z$ in Subsect. \ref{34.4}. Recall also $\hat{\eta}_{\alpha'}$, $\hat{\eta}_{\beta'}$,
$\hat{\xi}_\alpha$, $\hat{\xi}_\beta$   defined in Prop. \ref{HL_associativity_morphism_for_boxtimes}.
By construction we have
$$\hat{\eta}_{\alpha'}=Z(\eta_{\alpha'})\eta^{z_1-z_2}_{Z(W_1), Z(W_2)}=v_{\alpha'}^*, \quad \quad \hat{\eta}_{\beta'}=Z(\eta_{\beta'})\eta^{z_2}_{Z(M_i), Z(W_3)}=v_{\beta'}^*.$$
Since $$v_{\alpha_1}^*v_{\alpha_2}=\delta_{\alpha_1, \alpha_2},$$ we automatically have from the image of  (\ref{Cuntz}) under Zhu functor,
$$\hat{\xi}_\alpha=(\eta^{z_1}_{Z(W_1), Z(M_r)})^{-1}Z(\xi_\alpha)=v_\alpha, \quad\quad \hat{\xi}_\beta=(\eta^{z_2}_{Z(W_2), Z(W_3)})^{-1}Z(\xi_\beta)=v_\beta.$$
Inserting these data in (\ref{little_eff_adjoint}), by Prop. \ref{HL_associativity_morphism_for_boxtimes} and Theorem 
\ref{HL_is_of_CFT_type}
we have 
$$f_{\alpha, \beta}^{\beta', \alpha'}=F_{\alpha, \beta}^{\beta', \alpha'}.$$

 The case of conditions guaranteeing the assumptions in b) will be obtained in a similar way and will apply to determine the inverse matrix of the OPE constant matrix $(F_{\alpha, \beta}^{\beta', \alpha'})$ on triples of the form $(V, V, (L_{k, \lambda})_{(0)})$.
 
\subsection{Identification  between the Huang-Lepowsky braided symmetry for vertex operator algebra modules and Wassermann and Toledano Laredo braiding operators}

\begin{rem}
\label{HL_braiding}
We emphasize that while our analysis in Sections \ref{32} and \ref{33} focuses primarily on 
the Huang--Lepowsky tensor product and the comparison of associativity 
morphisms, the corresponding braided ribbon structure is completely determined and connected to the braided ribbon structure described in Section \ref{32.5} in the setting of loop group conformal net.
Specifically, in the VOA setting the category-theoretic braiding operator $\mathcal{R}_{W_1, W_2}$ 
is defined via the analytic continuation and geometric rotation of vertex operator 
intertwining fields, satisfying the formula $\mathcal{R}(Y(v, z)) = e^{z L(-1)} 
Y(v, e^{\pi i} z)$ \cite{HuangIV}. 
 For the category of modules of affine 
Lie algebras at positive integer levels, the existence of this braided vertex 
tensor category structure was   established by Huang and Lepowsky in 
\cite{Huang_Lepowski_affine}. 
The eigenvalues of the resulting 
double-braiding (monodromy) operator on an irreducible component $W_3 \subseteq 
W_1 \boxtimes W_2$ are  governed by the ribbon balancing relation 
$\theta_{W_1 \boxtimes W_2} = \mathcal{R}_{W_2, W_1} \circ \mathcal{R}_{W_1, W_2} 
\circ (\theta_{W_1} \otimes \theta_{W_2})$, where the twist is given by the 
Virasoro exponent $\theta = e^{2\pi i L(0)}$ by Theorem 4.1 in \cite{Huang2}. 
On the highest-weight spaces, this yields the explicit monodromy eigenvalue 
$\lambda = e^{2\pi i (h_3 - h_1 - h_2)}$, where $h_i$ are 
the conformal weights of the respective VOA modules. Via the Sugawara 
construction $h_\lambda = \frac{C_\lambda}{2(k + h^\vee)}$, these spectral 
phases translate precisely into the quadratic Casimir expressions yielding the 
quantum group parameter values $\pm q^{\frac{1}{2}(C_3 - C_1 - C_2)}$. This 
coincides exactly with the explicit monodromy formulas for the 
Knizhnik--Zamolodchikov equations obtained via the operator-algebraic loop group 
framework of Toledano Laredo \cite{Toledano_laredo}, which we detailed in 
Section \ref{32.5}, Corollary \ref{Toledano_laredo}. This numerical and structural agreement completely reconciles the 
monodromy data of both settings while bypassing the need for an independent 
geometric verification of the braiding axioms.
\end{rem}

 \bigskip

\section{Acknowledgments}    
C.P. wishes to express her heartfelt   gratitude   to her daughters for their support during the preparation of this work.

She is deeply grateful to 
 Sergio Doplicher for suggesting to study the relationship between quantum groups and low dimensional quantum field theories in 1996, for inspiration and constant encouragement during all the years she had the good fortune to share scientific discussions; to Vaughan F. Jones
for encouragement on the project  in Frascati in 2016; to Alessandro D'Andrea for first letting her
know of the Zhu algebra; to  Hans Wenzl for encouragement and  kind email correspondence   on the centraliser algebras
for the spinor representations,
 and for  references on  duality for tensor powers of the fundamental representation of Drinfeld-Jimbo quantum groups for the Lie types $B$, $D$, $G_2$ in August 2024, helpful to complete the proof of Lemma 7.10 in v7 of this paper, now Theorems \ref{claim0} and \ref{claim1};
to  Yi-Zhi Huang for several email exchanges in 2022 and 2025 helpful to  clarify the history  on tensoriality, rigidity, and modularity for the categories of representations of  vertex operator algebras, and on  KLF theorem.
She is also grateful to
 Makoto Yamashita for  helpful discussions and references on associators in 2021 and for informing her
 of \cite{Wenzl_E_8} in late 2023; to 
 Bin Gui for discussions  on his work and for letting us know about the paper by Kirillov \cite{Kirillov3} in 2019; to
 Eric Rowell for correspondence  on the status of classification of low rank modular categories in 2020; to
Ludmil K. Hadjiivanov for an email discussion in 2021;
to Paolo Aschieri, John Barrett,  Marcel Bischoff,   Cesar Galindo, Simone Del Vecchio,   Luca Giorgetti,  Andre Henriques,  Yi-Zhi Huang, Masahiko Miyamoto, Sergey Neshveyev,  Julia Y. Plavnik, Paolo Papi, David Reutter, 
  Christoph Schweigert, Lars Tuset,   Paola Zurlo for their   interest and discussions.

The results of the paper have been presented at several conferences. We are very grateful to the organizers for the invitation and giving us the opportunity to talk.

The authors are  grateful to an anonymous referee of an earlier stage of this project for constructive suggestions that greatly helped highlight the uniform advantages of this approach regarding non-simply laced Lie types.

We thank Sapienza University of Rome for financial support   Finanziamento di Ateneo per la Ricerca Scientifica, years 2018; 2019; 2020; 2021, 2023.      S.C. acknowledges support by Sapienza, 2016 Avvio alla ricerca `Gruppi quantistici alle radici dell'unit\`a e $C^*$-categorie tensoriali'; M.V.G. acknowledges support by   Sapienza, I.C.E. 16/22 based on Ateneo 2018.

An early version of Section 22 was developed in collaboration with S. Carpi in 2018; he subsequently withdrew from the project in January 2023.

\end{document}